\documentclass[11pt, reqno]{amsart}
\usepackage{appendix}
\usepackage[a4paper, centering]{geometry}
\usepackage{xcolor}
\usepackage[utf8]{inputenc}
\usepackage{epigraph}

\usepackage{musicography}
\usepackage{multicol}
\usepackage{caption}
\usepackage{bbm}

\usepackage[backend=biber, bibencoding=utf8,style=alphabetic,maxnames=99,maxalphanames=99,backref=true]{biblatex}
\usepackage[colorlinks=true,hyperindex,linktocpage=true]{hyperref}

\usepackage[capitalize]{cleveref}
\usepackage{mathtools}
\usepackage{amsfonts}
\usepackage{bm}
\usepackage{tikz-cd}
\usepackage[T1]{fontenc}

\usepackage{amssymb}    %
\usepackage{graphicx}   %
\usepackage{scalerel}   %

\usepackage{amsmath}

\hypersetup{
	colorlinks,
	linkcolor={violet!60!black},
	citecolor={cyan!60!black},
	urlcolor={orange!60!black}
}
\usepackage{stmaryrd}
\usepackage{mathrsfs}  
\usepackage{varwidth}
\usepackage{enumerate}

\theoremstyle{plain}
\newtheorem{theorem}{Theorem}[subsection]
\newtheorem{proposition}[theorem]{Proposition}
\newtheorem{construction}[theorem]{Construction}
\Crefname{construction}{Construction}{Constructions}
\newtheorem{lemma}[theorem]{Lemma}
\newtheorem{corollary}[theorem]{Corollary}
\newtheorem{conjecture}[theorem]{Conjecture}
\newtheorem{question}[theorem]{Question}

\theoremstyle{definition}
\newtheorem{definition}[theorem]{Definition}

\newtheorem{example}[theorem]{Example}
\newtheorem{remark}[theorem]{Remark}
\newtheorem{convention}[theorem]{Convention}
\newtheorem{Assumption/Conjecture}[theorem]{Assumption/Conjecture}

\newtheorem{assumption}[theorem]{Assumption}
\newtheorem{assumption/question}[theorem]{Assumption/Question}

\newcommand{\nc}{\newcommand}
\nc{\on}{\operatorname}

\nc{\Q}{\mathbb{Q}}
\nc{\Z}{\mathbb{Z}}
\nc{\cl}{\mathrm{cl}}

\newcommand{\pitch}{\mathbin{\rotatebox[origin=c]{180}{\bm{$\pitchfork$}}}}

\DeclareMathOperator{\RHom}{{{\it R} \mathscr{H}\mkern-1.8mu o\mkern-0.7mu m\mkern-0.5mu}}
\nc{\fraka}{{\mathfrak a}} \nc{\bba}{{\mathbf a}}
\nc{\frakb}{{\mathfrak b}}
\nc{\frakc}{{\mathfrak c}}
\nc{\frakd}{{\mathfrak d}}
\nc{\frake}{{\mathfrak e}}
\nc{\frakf}{{\mathfrak f}}

\nc{\frakg}{{\mathfrak g}}
\nc{\frakh}{{\mathfrak h}}
\nc{\fraki}{{\mathfrak i}}
\nc{\frakj}{{\mathfrak j}}
\nc{\frakk}{{\mathfrak k}}
\nc{\frakl}{{\mathfrak l}}
\nc{\frakm}{{\mathfrak m}}
\nc{\frakn}{{\mathfrak n}}
\nc{\frako}{{\mathfrak o}}
\nc{\frakp}{{\mathfrak p}}
\nc{\frakq}{{\mathfrak q}}
\nc{\frakr}{{\mathfrak r}}
\nc{\fraks}{{\mathfrak s}}
\nc{\frakt}{{\mathfrak t}}
\nc{\fraku}{{\mathfrak u}}
\nc{\frakv}{{\mathfrak v}}
\nc{\frakw}{{\mathfrak w}}
\nc{\frakx}{{\mathfrak x}}
\nc{\fraky}{{\mathfrak y}}
\nc{\frakz}{{\mathfrak z}}
\nc{\frakA}{{\mathfrak A}}
\nc{\frakB}{{\mathfrak B}}
\nc{\frakC}{{\mathfrak C}}
\nc{\frakD}{{\mathfrak D}}
\nc{\frakE}{{\mathfrak E}}
\nc{\frakF}{{\mathfrak F}}
\nc{\frakG}{{\mathfrak G}}
\nc{\frakH}{{\mathfrak H}}
\nc{\frakI}{{\mathfrak I}}
\nc{\frakJ}{{\mathfrak J}}
\nc{\frakK}{{\mathfrak K}}
\nc{\frakL}{{\mathfrak L}}
\nc{\frakM}{{\mathfrak M}}
\nc{\frakN}{{\mathfrak N}}
\nc{\frakO}{{\mathfrak O}}
\nc{\frakP}{{\mathfrak P}}
\nc{\frakQ}{{\mathfrak Q}}
\nc{\frakR}{{\mathfrak R}}
\nc{\frakS}{{\mathfrak S}}
\nc{\frakT}{{\mathfrak T}}
\nc{\frakU}{{\mathfrak U}}
\nc{\frakV}{{\mathfrak V}}
\nc{\frakW}{{\mathfrak W}}
\nc{\frakX}{{\mathfrak X}}
\nc{\frakY}{{\mathfrak Y}}
\nc{\Ani}{{\on{Ani}}}
\nc{\frakZ}{{\mathfrak Z}}
\nc{\bbA}{{\mathbb A}}
\nc{\bbB}{{\mathbb B}}
\nc{\bbC}{{\mathbb C}}
\nc{\bbD}{{\mathbb D}}
\nc{\bbE}{{\mathbb E}}
\nc{\bbF}{{\mathbb F}} \nc{\bbf}{{\mathbf f}}
\nc{\bbG}{{\mathbb G}}
\nc{\bbH}{{\mathbb H}}
\nc{\bbI}{{\mathbb I}}
\nc{\bbJ}{{\mathbb J}}
\nc{\bbK}{{\mathbb K}}
\nc{\bbL}{{\mathbb L}}
\nc{\bbM}{{\mathbb M}}
\nc{\bbN}{{\mathbb N}}
\nc{\bbO}{{\mathbb O}}
\nc{\bbP}{{\mathbb P}}
\nc{\bbQ}{{\mathbb Q}}
\nc{\bbR}{{\mathbb R}}
\nc{\bbS}{{\mathbb S}}
\nc{\bbT}{{\mathbb T}}
\nc{\bbU}{{\mathbb U}}
\nc{\bbV}{{\mathbb V}}
\nc{\bbW}{{\mathbb W}}
\nc{\bbX}{{\mathbb X}}
\nc{\bbY}{{\mathbb Y}}
\nc{\bbZ}{{\mathbb Z}}
\nc{\calA}{{\mathcal A}}
\nc{\calB}{{\mathcal B}}
\nc{\calC}{{\mathcal C}}
\nc{\calD}{{\mathcal D}}
\nc{\calE}{{\mathcal E}}
\nc{\calF}{{\mathcal F}}
\nc{\calG}{{\mathcal G}}
\nc{\calH}{{\mathcal H}}
\nc{\calI}{{\mathcal I}}
\nc{\calJ}{{\mathcal J}}
\nc{\calK}{{\mathcal K}}
\nc{\calL}{{\mathcal L}}
\nc{\LKB}{\calL_{(b,K)}}
\nc{\FS}{{\on{FS}}}
\nc{\RKB}{\calR_{(b,K)}}
\nc{\calM}{{\mathcal M}}
\nc{\calN}{{\mathcal N}}
\nc{\calO}{{\mathcal O}}
\nc{\calP}{{\mathcal P}}
\nc{\calQ}{{\mathcal Q}}
\nc{\calR}{{\mathcal R}}
\nc{\calS}{{\mathcal S}}
\nc{\calT}{{\mathcal T}}
\nc{\calU}{{\mathcal U}}
\nc{\calV}{{\mathcal V}}
\nc{\calW}{{\mathcal W}}
\nc{\calX}{{\mathcal X}}
\nc{\calY}{{\mathcal Y}}
\nc{\calZ}{{\mathcal Z}}

\nc{\scrA}{{\mathscr A}}
\nc{\scrB}{{\mathscr B}}
\nc{\scrC}{{\mathscr C}}
\nc{\scrD}{{\mathscr D}}
\nc{\scrE}{{\mathscr E}}
\nc{\scrF}{{\mathscr F}}
\nc{\scrG}{{\mathscr G}}
\nc{\scrH}{{\mathscr H}}
\nc{\scrI}{{\mathscr J}}
\nc{\scrJ}{{\mathscr I}}
\nc{\scrK}{{\mathscr K}}
\nc{\scrL}{{\mathscr L}}
\nc{\scrM}{{\mathscr M}}
\nc{\scrN}{{\mathscr N}}
\nc{\scrO}{{\mathscr O}}
\nc{\scrP}{{\mathscr P}}
\nc{\scrQ}{{\mathscr Q}}
\nc{\scrR}{{\mathscr R}}
\nc{\D}{{\on{D}}}
\nc{\Div}{{\on{Div}}}
\nc{\Perv}{{\on{Perv}}}

\nc{\bnu}{{\bar{ \nu}}}
\nc{\olO}{\bar{\calO}}

\nc{\al}{{\alpha}} 
\nc{\be}{{\beta}}
\nc{\ga}{{\gamma}} \nc{\Ga}{{\Gamma}}
\nc{\hGa}{\hat{\Gamma}}
\nc{\ve}{{\varepsilon}} 
\nc{\la}{{\lambda}} \nc{\La}{{\Lambda}}
\nc{\om}{\omega} \nc{\Om}{\Omega} 
\nc{\sig}{{\sigma}} \nc{\Sig}{{\Sigma}}
\nc{\dR}{{\mathrm{dR}}}
\nc{\Perf}{{\mathrm{Perf}}}
\nc{\Perff}{{\mathrm{Perf}^\aff}}
\nc{\Perfff}{{(\widetilde{\Perff})_{\on{v}}}}
\nc{\perf}{{\mathrm{perf}}}
\nc{\qcqs}{{\mathrm{qcqs}}}
\nc{\pfp}{{\mathrm{pfp}}}
\nc{\rep}{{\mathrm{rep}}}
\nc{\cons}{{\mathrm{cons}}}
\nc{\PSch}{{\mathrm{PSch}}}
\nc{\PSchf}{{\mathrm{PSch}^\aff}}
\nc{\PSchff}{{(\widetilde{\PSchf})}_{\on{v}}}
\nc{\PShv}{{\mathrm{PShv}}}
\nc{\Fun}{{\mathrm{Fun}}}
\nc{\AlgSp}{{\mathrm{AlgSp}}}
\nc{\Sch}{{\mathrm{Sch}}}
\nc{\PreStk}{{\mathrm{PreStk}}}

\nc{\Sptl}{{\mathrm{Sptl}}}
\nc{\StdPerf}{{\mathrm{StdPerf}}}
\nc{\LocSptl}{{\mathrm{LocSptl}}}
\nc{\AnArt}{{\mathrm{AnArt}}}
\nc{\AnStk}{{\mathrm{AnStk}}}
\nc{\AnPreStk}{{\mathrm{AnPreStk}}}
\nc{\SchStk}{{\mathrm{SchStk}}}
\nc{\resSchStk}{{\mathrm{{}_{res}SchStk}}}
\nc{\CAlg}{{\mathrm{CAlg}}}
\nc{\Frob}{{\mathrm{Frob}}}
\nc{\Gm}{{\mathbb{G}_m}}
\DeclareMathOperator{\colim}{colim}
\nc{\et}{\mathrm{\acute{e}t}}
\nc{\proet}{\mathrm{pro\acute{e}t}}
\nc{\uop}{\mathrm{pro\acute{e}t}}
\nc{\mer}{{\on{mer}}}
\nc{\FF}{{\on{FF}}}
\nc{\sht}{{\on{sht}}}
\nc{\sss}{{\on{ss}}}
\nc{\Sets}{{\on{Sets}}}
\nc{\Grps}{{\on{Grps}}}

\nc{\AlgSt}{{\on{AlgSt_\ell}}}

\nc{\vvb}{{\on{Vect}^{\calO^\sharp}_{\on{v}}}}
\nc{\Vect}{{\on{Vect}}}
\nc{\onv}{{\on{v}}}
\nc{\fg}{{\on{f.g.}}}
\nc{\renn}{{\on{ren}}}
\nc{\Kan}{{\on{Kan}}}
\nc{\Zhu}{{\on{Zhu}}}

\nc{\T}[2]{\calT_{#1}^{#2}}
\nc{\Hkk}[1]{{\on{Hk}}^\bullet ( {#1} )}
\nc{\DM}{{\frakD\frakM}}
\nc{\DMG}{{\frakD\frakM_{\calG}}}
\nc{\IC}{{\mathfrak{B}}}
\nc{\ICG}{{\frakB(G)}}
\nc{\SHT}{{\calS\calH\calT}}
\nc{\VEC}{{\on{Bun}}_{\on{FF}}^{\on{mer}}}
\nc{\VECT}{{\on{Bun}}_{\on{FF}}}
\nc{\Catex}{{\on{Cat}^{\otimes, \on{ex}}_1}}
\nc{\Cat}{{\on{Cat}}}
\nc{\LinCat}{{\on{LinCat}}}
\nc{\Fil}{\calF{il}}
\nc{\Dm}{{\on{DM}}}
\nc{\Sht}{\on{Sht}}
\nc{\Isoc}{\on{Isoc}}
\nc{\Shv}{\mathrm{Shv}}
\nc{\IndShv}{\on{IndShv}_{\on{f.g.}}}
\nc{\fShv}{\calS\hspace{-0.1em}\textit{hv}}

\def\PrL{\mathrm{Pr}^{L}}
\def\ICoh{\mathrm{IC}}
\nc{\BC}{\calB\calC}
\nc{\PFS}{\bbP_{{\on{FS}}}}
\def\preceqdot{\mathrel{\preceq\kern-.5em\raise.22ex\hbox{$\cdot$}}}

\makeatletter
\DeclareFontEncoding{LS1}{}{}
\DeclareFontSubstitution{LS1}{stix}{m}{n}
\DeclareMathAlphabet{\rhomalpha}{LS1}{stixscr}{m}{n}
\makeatother
\nc{\Spa}{\on{{Spa}}}
\nc{\Spd}{\on{{Spd}}}
\nc{\tnb}{\psi_{\rm tame}}
\nc{\oM}{\overline{{M}}}
\nc{\op}{{\on{op}}}
\nc{\ad}{{\on{ad}}}
\nc{\alg}{{\on{alg}}}
\nc{\Art}{{\on{Art}}}
\nc{\Alg}{{\on{Alg}}}
\nc{\Ad}{{\on{Ad}}}
\nc{\Adm}{{\on{Adm}}}
\nc{\aff}{{\on{aff}}}
\nc{\Aut}{{\on{Aut}}}
\nc{\Bun}{{\on{Bun}}}
\nc{\cha}{{\on{char}}}
\nc{\der}{{\on{der}}}
\nc{\Der}{{\on{Der}}}
\nc{\diag}{{\on{diag}}}
\nc{\End}{{\on{End}}}
\nc{\Fl}{{\calF\!\ell}}
\nc{\Tr}{{\on{Transp}}}
\nc{\TR}{{\calT\!\calR}}
\nc{\Gal}{{\on{Gal}}}
\nc{\Gr}{{\on{Gr}}}
\nc{\Hk}{{\on{Hk}}}
\nc{\rH}{{\on{H}}}
\nc{\uHom}{{\ul{\on{Hom}}}}
\nc{\Hom}{{{\on{Hom}}}}
\nc{\id}{{\on{id}}}
\nc{\Id}{{\on{Id}}}
\nc{\ind}{{\on{ind}}}
\nc{\Ind}{{\on{Ind}}}
\nc{\Lie}{{\on{Lie}}}
\nc{\Pic}{{\on{Pic}}}
\nc{\pr}{{\on{pr}}}
\nc{\Res}{{\on{Res}}}
\nc{\res}{{\on{res}}} \nc{\Sat}{{\on{Sat}}}
\nc{\spc}{{\on{sc}}}
\nc{\drv}{{\on{der}}}
\nc{\sgn}{{\on{sgn}}}
\DeclareMathOperator{\Spec}{Spec}
\nc{\Spf}{\on{Spf}} 
\nc{\Sph}{\on{Sph}}
\nc{\St}{{\on{St}}}
\nc{\tr}{{\on{tr}}}
\nc{\Mod}{{\mathrm{-Mod}}}
\nc{\Hilb}{{\on{Hilb}}} 
\nc{\Ext}{{\on{Ext}}} 
\nc{\vs}{{\on{Vec}}}
\nc{\ev}{{\on{ev}}}
\nc{\nO}{{\breve{\calO}}}
\nc{\tS}{{\tilde{S}}}
\nc{\spe}{{\on{sp}}}
\nc{\loc}{{\on{loc}}}
\nc{\pre}{{\on{pre}}}
\nc{\bb}{\mathbb}

\nc{\dimt}{{\on{dim.trg}}}

\let\ol\overline

\nc{\co}{\colon}
\nc{\dia}{{\diamondsuit}}

\nc{\nscrR}{{\mathscr{R}^{\on{nr}}}}

\nc{\GL}{{\on{GL}}}

\nc{\Gl}{\on{Gl}} 
\nc{\GSp}{{\on{GSp}}}
\nc{\gl}{{\frakg\frakl}}
\nc{\SL}{{\on{SL}}} 
\nc{\SU}{{\on{SU}}} 
\nc{\SO}{{\on{SO}}}
\nc{\PGL}{{\on{PGL}}}

\nc{\Conv}{{\on{Conv}}}
\nc{\Rep}{{\on{Rep}}}
\nc{\Dom}{{\on{Dom}}}
\nc{\red}{{\on{red}}} 
\nc{\Red}{{\on{Red}}} 
\nc{\act}{{\on{act}}}
\nc{\nr}{{\on{nr}}}
\nc{\ctf}{{\on{ctf}}}

\nc{\str}{{\on{-}}} 
\nc{\os}{{\bar{s}}}
\nc{\oeta}{{\bar{\eta}}}

\newcommand{\Eis}{\mathrm{Eis}}
\newcommand{\IW}{\mathrm{IW}}
\newcommand{\un}{\mathrm{un}}
\newcommand{\spec}{\mathrm{spec}}
\newcommand{\basic}{\mathrm{basic}}

\newcommand{\QCoh}{\mathrm{QCoh}}
\newcommand{\mf}{\mathfrak}
\newcommand{\lisse}{\mathrm{lis}}
\newcommand{\can}{\mathrm{can}}
\newcommand{\ExtRHom}{\mathrm{Hom}}
\newcommand{\profet}{\mathrm{prof\acute{e}t}}
\newcommand{\coarse}{\mathrm{coarse}}

\newcommand{\IndCoh}{\mathrm{IndCoh}}
\newcommand{\Coh}{\mathrm{Coh}}
\newcommand{\IndPerf}{\mathrm{IndPerf}}
\newcommand{\unip}{\mathrm{unip}}

\newcommand{\Closed}{\mathrm{Closed}}
\newcommand{\Detale}{\mathcal{D}_{\mathrm{\acute{e}t}}}
\newcommand{\ra}{\rightarrow}
\newcommand{\an}{\mathrm{an}}
\newcommand{\sch}{\mathrm{sch}}
\newcommand{\sstd}{\mathrm{std}}
\newcommand{\onepinfty}{\frac{1}{p^{\infty}}}
\newcommand{\fine}{\mathrm{fine}}

\newcommand{\ULA}{\mathrm{ULA}}
\newcommand{\Dlis}{\mathcal{D}_{\mathrm{lis}}}
\newcommand{\ul}{\underline}
\newcommand{\cInd}{\on{c-Ind}}
\newcommand{\Nt}{\mathrm{Nt}}
\newcommand{\BZ}{\mathrm{BZ}}
\newcommand{\coh}{\mathrm{coh}}
\newcommand{\ren}{\mathrm{ren}}
\newcommand{\lis}{\mathrm{lis}}
\newcommand{\spl}{\mathrm{spl}}
\newcommand{\pl}{{\mathrm{qc},\mathrm{pl}}}
\newcommand{\sm}{\mathrm{sm}}
\newcommand{\cg}{\mathrm{cg}}
\newcommand{\qpl}{{\mathrm{qc},\mathrm{qpl}}}
\newcommand{\qome}{{\mathrm{qc},\mathrm{qpl},\mathrm{cg}}}
\newcommand{\IndStk}{\mathrm{IndStk}}
\newcommand{\sStk}{\mathrm{sIndStk}}
\newcommand{\vpl}{{\mathrm{qc},\mathrm{vpl}}}

\newcommand{\oc}{{\on{o.c.}}}

\def\Par{\mathrm{Par}}

\def\Corr{\mathrm{Corr}}
\def\ins{\mathrm{ins}}
\def\Suave{\mathrm{Suave}}
\def\SD{\mathrm{SD}}
\def\Dcons{\mathcal{D}_{\mathrm{cons}}}
\def\vShv{\mathrm{vShv}}
\def\vStk{\mathrm{vStk}}
\def\fdcs{\mathrm{fdcs}}
\def\fdcss{\mathrm{fdcss}}
\def\tame{\mathrm{tame}}

\def\LiuZheng{\calL\calZ}
\def\finexp{\mathrm{fin.exp.}}

\nc{\hookto}{\hookrightarrow}
\nc{\longto}{\longrightarrow}
\nc{\leftto}{\leftarrow}
\nc{\onto}{\twoheadrightarrow}
\nc{\lonto}{\twoheadleftarrow}

\nc{\pot}[1]{ [\hspace{-0,5mm}[ {#1} ]\hspace{-0,5mm}] }
\nc{\rpot}[1]{ (\hspace{-0,7mm}( {#1} )\hspace{-0,7mm}) }
\nc{\smallpot}{{ <\hspace{-1,0mm}<}}

\numberwithin{equation}{section}

\newcommand{\rar}{\rightarrow}
\newcommand\prolim{\mathop{\underleftarrow{\lim} }}

\DeclareMathOperator{\corr}{Corr}

\DeclareMathOperator{\Modd}{Mod}
\DeclareMathOperator{\all}{All}
\nc{\from}{\colon}
\nc{\defined}{\coloneqq}
\nc{\PrLst}{\PrL_{\mathrm{st}}}
\nc{\PrLD}{\PrL_{\D}}
\nc{\injto}{\hookrightarrow}
\newcommand{\xto}[1]{\xrightarrow{#1}}

\newcommand{\Stk}{{\mathrm{Stk}}}

\DeclareMathOperator{\Convex}{Convex}

\newcommand{\maybe}[1]{\ignorespaces}

\begin{document}
	\title{On the schematic and analytic constructions of the local Langlands category}
	\author[I. Gleason, L. Hamann, A. Ivanov, J. Louren\c{c}o, K. Zou]{Ian Gleason, Linus Hamann, Alexander B. Ivanov, Jo\~ao Louren\c{c}o, Konrad Zou}

	\address{National University of Singapore,  Department of Mathematics, Faculty of Science, National University of Singapore, Block S17, 10 Lower Kent Ridge Road, Singapore 119076}
	\email{ianandreigf@nus.edu.sg}

    \address{Harvard University,  Department of Mathematics, 1 Oxford Street Cambridge, MA, 02138}
	\email{hamannn@math.harvard.edu}

	\address{Fakult\"at f\"ur Mathematik, Ruhr-Universit\"at Bochum, Universit\"atsstrasse 150, D-44780 Bochum, Germany}
	\email{a.ivanov@rub.de}

    \address{LAGA, UMR 7539, CNRS, Universit\'e Sorbonne Paris Nord (Paris 13), 99 avenue Jean Baptiste Cl\'ement, 93430 Villetaneuse, France}
    \email{lourenco@math.univ-paris13.fr}

    \address{Institut de Math\'ematiques de Jussieu - Paris Rive Gauche, Sorbonne Universit\'e – Campus Pierre et Marie Curie, place Jussieu 4, 75252 Paris Cedex 05, France}
    \email{konrad.zou@imj-prg.fr}
	
	\begin{abstract}
    We prove a folklore conjecture identifying two categorical enhancements of the automorphic side of the local Langlands correspondence. 
    Concretely, we construct an equivalence for torsion coefficients between the category considered by Zhu and the one considered by Fargues--Scholze. 
    To achieve this, we revisit Scholze's analytification functor and apply the first author's theory of kimberlites. 
   We discuss unconditional applications to the splitting of the semi-orthogonal decomposition on $\Bun_{G}$, and the compatibility with Eisenstein functors. 
   Finally, we formulate a linearity conjecture for our functor with which we can show new vanishing statements for the cohomology of local Shimura varieties, and perverse exactness statements for Hecke operators. 
	\end{abstract}

	\maketitle
	\setcounter{tocdepth}{2}
	\tableofcontents	
    \setcounter{tocdepth}{1}
	\newpage

	\section{Introduction}
	We divide our introduction into a soft part and a technical part. 
	The reader who enjoys mathematical story-telling is invited to read the soft introduction (\S \ref{s: SoftIntroduction}). 
	Otherwise, we invite the pragmatic reader to skip directly to the technical part (\S \ref{s: Technical Introduction}). 
	Having divided our audience by taste, we will shamelessly be vague during the soft portion of this introduction, to the extent that we cannot even be wrong.
	\subsection{Soft introduction}{\label{s: SoftIntroduction}}
	Categorification, in modern representation theory, has already transitioned from being treated as a technical tool to an object of independent mathematical interest.  
	Following a key insight of Fargues in 2014, number theorists and arithmetic geometers working in the Langlands program have realized throughout the last decade that the classical local Langlands correspondence could and should be upgraded to a categorical statement.
We refer the reader to \cite{fargues_geometrization_of_the_local_langlands_correspondence_overview,gaitsgory2016geometric,Zhu2018} for informal early accounts on the subject. 

	Following the hard work of numerous experts coming from different backgrounds, the community has reached the so called ``categorical local Langlands conjecture'' (or CLLC).
This is a beautiful conjecture that interweaves classical representation theory with arithmetic geometry.
Surprisingly, however, this conjecture came in two different formulations which address different pieces of structure.
One formulation of the conjecture was made precise by Fargues--Scholze (\cite[Conjecture X.3.5]{FS21}) and the other was made precise in the work of Zhu (\cite[Conjecture 2.1.3]{Zhu20}).
We will refer to the first theory as the $\Bun_G$ perspective, and to the second one as the $\ICG$ perspective.
	The purpose of this article is to construct an equivalence of categories, $\pitch$, whose role is to bridge the two perspectives on how to formulate the CLLC.
This work builds on \cite{GIZ25,Gle23}, where some of us laid the  geometric foundations necessary for the construction of the functor $\pitch$.
    
	As we will try to convince the reader later in this introduction, this work is not a linguistic translation.
	Instead, each of the perspectives is sensitive to different phenomena, and connecting these perspectives enriches our understanding of the CLLC, and its connection to other mathematical objects of interest.   
	Although we are forced to discuss a substantial amount of formalism, the key to $\pitch$ ultimately comes from interesting geometric observations. 
\\

	From an impressionistic lens, the CLLC resembles the classical formulation of unramified Langlands.   
	To explain this, let us first fix some notation.
	We let $E$ be a non-Archimedean local field of residue characteristic $p$, with ring of integers $O_E\subseteq E$ and residue field $\bbF_q$ with $q$ elements.
We set $W_E$ to be the Weil group of $E$ and set $I_E\subseteq W_E$ to be the inertia subgroup. We let $G/E$ be a connected reductive group, which we assume to be split for simplicity in this introduction, and we fix as usual $T\subseteq B\subseteq G$ a Borel and a maximal torus.
We let $W$ denote the Weyl group of $G$.
We let $\Lambda$ be a suitable coefficient ring in which $p$ is invertible.

	Recall that the classical unramified local Langlands can be obtained by using the following recipe.
If $V$ denotes a $\bbC$-valued smooth irreducible and spherical representation, then $V^{G(O_E)}$ becomes an irreducible $\calH$-module where
	\[\calH:=(C_{c}^{\infty}(G(O_E)\backslash G(E)/ G(O_E),\mathbb{C}),\ast)\]
	denotes the spherical Hecke algebra equipped with convolution.
A consequence of the Satake isomorphism
	\begin{equation}
		\label{eq: classical satake}	
	\calH\simeq \bbC[X_*(T)^W], 
	\end{equation}
	is that every irreducible $\calH$-module is $1$-dimensional.
	In particular, $V$ is completely determined by the eigenvalues of the simultaneously diagonalizable Hecke operators acting on the only $G(O_E)$-fixed line of $V$.
 
	These eigenvalues are the so-called Satake parameters of $V$, which can be encoded by a semi-simple conjugacy class in $\hat{G}$.
	In this way, classical unramified Langlands is a spectral classification of spherical irreducible representations, i.e., it classifies spherical irreducible representations by their Hecke eigenvalues.
	It also encompasses a multiplicity-one statement saying that for each family of eigenvalues, there is exactly one irreducible spherical representation corresponding to it.  

	The first observation that the CLLC makes is that one should enlarge the (derived) category of $G(E)$-representations, to a so-called local Langlands category.
	In the $\Bun_G$ setup we will denote it as $\calD^\an_\Lambda(\Bun_G)$ (see \cite[\S 5]{FS21}), and in the $\ICG$ setup we will denote it by $\Shv^!(\ICG)$ (see \cite[\S 3.4]{Zhu25}). 
	Each of these categories come equipped with fully faithful functors 
	\[j_{1,!}:\on{Rep}(G(E))\to \calD^\an_\Lambda(\Bun_G)\]
	and 
	\[i_{1,!}:\on{Rep}(G(E))\to \Shv^!(\ICG)\]
	respectively, where $\Rep(G(E))$ denotes the derived category of smooth representations of $G(E)$ on $\Lambda$-modules.
When we want to speak of a platonic local Langlands category, without a specific construction of it, we will denote it as $\calD_{\on{LLC}}$.
	That we can take this platonic attitude is partially justified by our main result, which is the construction of an equivalence
	\[\pitch: \Shv^!(\ICG) \xrightarrow{\simeq} \calD^\an_\Lambda(\Bun_G)\]
	between LLC categories satisfying several desiderata.

	In the CLLC, instead of considering the spherical Hecke algebra $\calH$ one considers the Satake category $(\calS_G,\ast)$, which is a category with a commutative algebra structure under convolution a.k.a.~a symmetric monoidal category. 
	Just as the spherical Hecke algebra is the space of compactly supported functions on the double coset space $G(O_E)\backslash G(E)/G(O_E)$ endowed with convolution, the Satake category is obtained by considering perverse sheaves on the double coset stack $L^+G\backslash LG/L^+G$, the so-called local Hecke stack, and the multiplication is also through convolution. 
	Just as in the classical picture, the Satake category has a spectral expression through the so called geometric Satake
	\[(\calS_G,\ast)\simeq (\on{Rep} \hat{G},\otimes),\]
	which is a categorification of \eqref{eq: classical satake}.
	Indeed, highest weight theory identifies $X_*(T)^W$ with the set of isomorphism classes of irreducible algebraic representations of $\hat{G}$ for $\Lambda = \ol{\bb{Q}}_{\ell}$.

	Analogous to the action of the spherical Hecke algebra on the spherical vectors of a smooth representation one should ``treat $\calD_{\on{LLC}}$ as a vector space'', and endow $\calD_{\on{LLC}}$ with Hecke operators.
In other words, one should realize $\calD_{\on{LLC}}$ as an $\calS_G$-module.
	Since we are in a categorical context, the Hecke operators themselves have automorphisms.
	To wit, in the CLLC one should consider an arithmetic version of $(\calS_G,\ast)$ that varies with the arithmetic geometry of $E$, this is customarily done through the Beilinson--Drinfeld Grassmannian in its $B_{\on{dR}}^+$-incarnation.
	This furnishes Hecke operators with a plethora of automorphisms coming from the Weil group action: in very rough terms, a ``Galois loop'' gives rise to coherently organized automorphism of the Hecke operators (this is known as a factorization structure), which in turn gives rise to the so-called excursion operators, as introduced by V. Lafforgue \cite{Laf18}.
	Following ideas of geometric Langlands for Riemann surfaces (see \cite{nadler2019spectral}), one can package excursion operators, i.e., the data of Hecke operators with coherently organized Galois automorphisms, through the so-called spectral action. 
	That is, one forms the stack $\Par_{G,E}$ of $L$-parameters, and the $(\on{Rep}(\hat{G}),\otimes)$-action promotes to an action 
	\[(\on{QCoh}(\Par_{G,E}),\otimes) \circlearrowright \calD_{\on{LLC}}.\]
    Here $\on{QCoh}(\Par_{G,E})$ denotes a suitable derived category of quasi-coherent sheaves on the parameter stack.
The naive multiplicity one statement in this setup would posit that $\calD_{\on{LLC}}$ is a free $1$-dimensional $(\on{QCoh}(\Par_{G,E}),\otimes)$-module, but this is provably false outside the case in which $G$ is a torus. 
	Instead, after fixing a Whittaker datum $\psi$, the CLLC predicts that there is a unique $(\on{QCoh}(\Par_{G,E}),\otimes)$-linear equivalence
	\begin{equation}
		\label{CLLC conjecture}
		\bbL_\psi: \on{IndCoh}_{\mathrm{Nilp}}(\Par_{G,E})\simeq \calD_{\on{LLC}} 
	\end{equation}
	with the property that 
	\[\bbL_\psi(\calO_{\Par_{G,E}})= \calW_\psi\]
	where $\calW_\psi\in \calD_{\on{LLC}}$ is the Whittaker sheaf attached to $\psi$.
Here $ \on{IndCoh}_{\mathrm{Nilp}}(\Par_{G,E})$ is a version of $\on{QCoh}(\Par_{G,E})$ that accounts for the singularities of $\Par_{G,E}$ in a precise categorical sense, and it comes equipped with a canonical $(\on{QCoh}(\Par_{G,E}),\otimes)$-action.
In this way, the CLLC is also a spectral parametrization where the Hecke operators, and the narrative, has been categorified.
Strictly speaking, $\on{IndCoh}(\Par_{G,E})$ is not a free $1$-dimensional $(\on{QCoh}(\Par_{G,E}),\otimes)$-module, but it is so in a smooth and dense open locus of $\Par_{G,E}$. 
	  In this way, the CLLC also predicts how to precisely correct the multiplicity one statement in this context. \\

	  In the $\Bun_G$ setup, Fargues--Scholze have constructed all of the structures necessary to carry out this narrative, completely formulating the CLLC conjecture for $\calD^\an_\Lambda(\Bun_G)$, i.e., Fargues--Scholze construct\footnote{Fargues--Scholze's construction assumes some mild assumptions on $\ell$ the characteristic of $\Lambda$.} a spectral action 
\[(\on{QCoh}(\Par_{G,E}),\otimes) \circlearrowright \calD^\an_\Lambda(\Bun_{G}),\]
	incarnated geometrically as a global Hecke stack on the Fargues--Fontaine curve, and conjecture that there is a unique $\on{QCoh}(\Par_{G,E})$-linear equivalence 
	\[\bbL^\an_\psi: \on{IndCoh}(\Par_{G,E})\simeq \calD^\an_\Lambda(\Bun_G) \]
	with the property that 
	\[\bbL^\an_\psi(\calO_{\Par_{G,E}})= \calW^\an_\psi.\]

	The fundamental motto guiding Fargues--Scholze's work is that the CLLC should be interpreted as ``global unramified Geometric Langlands for the Fargues--Fontaine curve''. 
	A major part of Fargues--Scholze's work goes into providing geometric foundations to show that in this exotic arithmetic context, the key geometric players are well behaved, and one can still find all of the structures that the geometric Langlands program for Riemann surfaces had recognized.
	Fargues--Scholze's framework invites the possibility to consider several other important constructions in global Geometric Langlands and its consequences to the CLLC.
A notable example are the so-called geometric Eisenstein series which have been explored and shown to be well-behaved in \cite{HamGeomES,HIDualCpx,HHS}.  \\

	In the $\ICG$ setup, one does not yet have a full precise formulation of the CLLC as we have discussed above, at least for mixed-characteristic local fields $E$.  
	Indeed, although one already has the language and technology to pose that there should be an equivalence
	\[\bbL^\sch_\psi: \on{IndCoh}(\Par_{G,E})\simeq \Shv^!(\ICG), \]
	with the property that 
	\[\bbL^\sch_\psi(\calO_{\Par_{G,E}})= \calW^\sch_\psi,\]
	one does not have a spectral action on $\Shv^!(\ICG)$ in full generality\footnote{There is work in progress by \'Eteve--Gaitsgory--Genestier--Lafforgue constructing a spectral action on $\Shv^!(\ICG)$ when $E$ is an equicharacteristic local field.}.
	In particular, one cannot ask if this equivalence is linear, and it becomes much more subtle to formulate in what sense this equivalence is unique.
	Nevertheless, one can already formulate and prove a tame-level version of \eqref{CLLC conjecture}, as Zhu does \cite[Theorem 1.6, Theorem 1.7]{Zhu25}.
	More precisely, there is an open and closed substack 
	\[\Par^{\on{tame}}_{G,E}\subseteq \Par_{G,E},\]
	that parametrizes the tamely ramified $L$-parameters (see \S \ref{ss: TameUnipotentStackofLanglandsParameters}).
	Moreover, Zhu defines a tamely ramified local Langlands subcategory 
	\[\Shv^{!,\on{tame}}(\ICG) \subseteq \Shv^!(\ICG), \]
	and shows, after fixing adjectives appropriately, that there is an equivalence 
	\[ \bbL^{\sch,\on{tame}}_\psi: \on{IndCoh}(\Par^{\on{tame}}_{G,E}) \simeq \Shv^{!,\on{tame}}(\ICG)\]
    (see Theorem \ref{thm: tamecategoricalLanglandsEquivalence}).
From Zhu's work, it follows that one can endow $\Shv^{!,{\on{tame}}}(\ICG)$ with a $\on{QCoh}(\Par^{\on{tame}}_{G,E})$-action, but this is rather indirect (see \cite[\S 3.1.5]{YangZhuTorsion}).    
	That this definition is indirect, is partially due to the fact that the $\ICG$ setup interprets the CLLC as ``the Frobenius trace of the local geometric Langlands conjecture''.
	Indeed, the starting point is Bezrukavnikov's equivalence \cite{Bez16} (interpreted appropriately as a $2$-categorical statement in local geometric Langlands). 
	In this $2$-categorical context, one has a monoidal category which we will denote $(\calB,\ast)$ in this introduction.
Hecke operators in this context are incarnated through Gaitsgory's central functor, which is a monoidal functor
	\[\calZ: (\on{Rep}(\hat{G}),\otimes)\to (\calB,\ast)\]
	factoring through the center.
	This construction is compatible with inertial Galois loops which gives an action 
	\[\on{QCoh}(\Par_{G,\breve E}^{\on{tame}}) \circlearrowright (\calB,\ast), \]
    where $\Breve{E}$ is the completion of the maximal unramified extension of $E$.
After taking trace of Frobenius, this yields a partial spectral action
	\[\on{QCoh}(\Par_{G,E}^{\on{tame}})\simeq \on{Tr}(\on{QCoh}(\Par_{G,\breve E}^{\on{tame}})) \circlearrowright  \on{Tr}((\calB,\ast))\]
    restricted exclusively to the tame parts.
	The work of Zhu discusses patiently some of the formal foundations necessary to carry this narrative, and crucially constructs an identification 
	\[\Shv^{!,{\on{tame}}}(\ICG)\simeq \on{Tr}((\calB,\ast)).\]
It is through this identification, and through transfer of structure, that Zhu endows $\Shv^{!,{\on{tame}}}(\ICG)$ with a spectral action. \\ 

Although we believe it is profitable to have two perspectives on the CLLC, from a purist lens in which there should be only one local Langlands category $\calD_{\on{LLC}}$.
What our work explains is a precise sense in which Fargues--Scholze's geometrization relates to local geometric Langlands.  
Namely, one of the many things that the work of Zhu explains is in what sense $\Shv^!(\ICG)$ is a natural recipient for the Frobenius trace construction, and why one obtains a fully-faithful functor 
\[\on{Tr}(\calB,\ast)\hookrightarrow \Shv^!(\ICG).\]
In essence, our work explains why $\calD^\an_\Lambda(\Bun_G)$ is also a natural recipient for the trace construction.
Independently of Zhu's work, one can use our methods to directly construct a functor  
\[\on{Tr}(\calB,\ast)\to \calD^\an_\Lambda(\Bun_G),\]
but it would only be through Zhu's work and our equivalence $\pitch$, that one can clearly see why this functor is fully faithful. 

From our perspective, the clear advantage of working with $\Shv^!(\ICG)$ is its immediate connection to local geometric Langlands, whereas the clear advantage of working with $\calD^\an_\Lambda(\Bun_G)$ is that, even in mixed-characteristic, it has a geometrically defined spectral action. 
By reconciling these two perspectives, we obtained a platonic $\calD_{\on{LLC}}$ that has a geometrically defined spectral action and can be studied through its connection to the geometric local Langlands program.
As a small application of this perspective, consider the following.

Relying on techniques from local geometric Langlands, Deligne--Lusztig theory and the geometry of affine Deligne--Lusztig varieties (ADLV), Zhu constructs a direct sum decomposition 
\[\Shv^{!,\on{tame}}(\ICG)\simeq \bigoplus_{\zeta}\Shv^{!,\hat{\zeta},}(\ICG), \]
as $\zeta$ ranges over the tame inertial types. 
As a proof of concept, one of the immediate corollaries of our work (see \Cref{cor: variantsoftheequivalence} for details) is that we also have a similar decomposition 
\[\calD^{\an,\tame}_\Lambda(\Bun_G) \simeq \bigoplus_{\zeta}\calD_\Lambda^{\an,\hat{\zeta}}(\Bun_G).\]
That this direct sum decomposition holds for $\calD^{\an,\tame}_\Lambda(\Bun_G)$ is predicted by the CLLC in the $\Bun_G$ setup, but it is quite unclear otherwise (see Remark \ref{rem: DiscussionofSplitting}). 
Using $\pitch$, we can already obtain this result unconditionally.

We expect many more applications to be obtained in the near future, not only in the CLLC but also in more classical arithmetic Langlands, by using our functor to combine the methods from the two setups.
This will become even more likely once the connection between $\pitch$ and Gaitsgory's central functors has been properly elucidated, which some of us are planning to discuss in future work.
For a precise statement of the kind of desired result, see Conjecture \ref{conj: IndPerfLinearity}.\\

Leaving the purist point of view behind, the two perspectives on the CLLC interact with different geometric objects of mathematical interest. 
For this reason, we believe it is interesting to understand the CLLC through both perspectives and be able to transfer results from one perspective to the other.
We make a table highlighting important examples, we have left a question mark ``?'' when the analogue in the table is not understood. 

\begin{table}[h]
\centering
\begin{tabular}{c c}
\hline
$\Bun_G$ & $\ICG$ \\
\hline
Generic fiber of Shimura varieties & Special fiber of Shimura varieties \\
$p$-adic period domains \& local Shimura varieties (LSV) & affine Deligne--Lusztig varieties (ADLV) \\
Fargues--Scholze parameters & Genestier--Lafforgue parameters \\
? & $p$-adic Deligne--Lusztig stacks \\
Geometric Eisenstein series & ?\\
\hline
\end{tabular}
\end{table}

For example, Zhu uses the dimension formula of ADLV to show that Hecke operators are perverse t-exact on $\Shv^!(\ICG)$.
Since our functor is also perverse t-exact (see \Cref{intro:t-exact}), under the expected property that $\pitch$ intertwines Hecke operators (which will be verified in a different work), one can conclude a similar perversity statement of Hecke operators on $\calD^\an_\Lambda(\Bun_G)$, as established in certain cases by \cite[Corollary~4.27]{HamannLeeTorsion}.
One can then use this perversity statement to deduce vanishing results in the cohomology of local Shimura varieties (see \Cref{shtukas corollary}), analogous to the results in \cite[\S~10.2]{HamGeomES}.\footnote{The standard way to connect the cohomology of LSV to the geometry of ADLV is through the nearby cycles construction.
One can think that the functor $\pitch$ encapsulates a plethora of nearby cycle computations, and it is precisely by studying specific nearby cycles (see Theorems \ref{nearbycyclestheorem} and \ref{hyperbolic localization for some quotients}) that the equivalence is obtained.
The magic of $\pitch$ is that it gives for free more nearby cycle computations than those required to show that it is an equivalence.} For this last part, it was critical to use the CLLC formulated with $\calD_{\on{LLC}}=\calD^\an_\Lambda(\Bun_G)$, since local Shimura varieties do not directly interact with $\ICG$.
In addition, using our functor $\pitch$ one can use results on the splitting of the semi-orthogonal decomposition on $\Shv^{!}(\ICG,\Lambda)$ around a generic localization, as proven in \cite{YangZhuTorsion}, in order to obtain unconditional analogous results for a generic localization of $\calD_{\Lambda}^{\an}(\Bun_{G})$ (see Corollaries \ref{cor: SplitSemiOrthogonalDecomposition} and \ref{cor: semiorthogonaldecompositionsplits}), which proves results along similar lines to \cite{HamGeomES,HamannLeeTorsion} that were originally established by a careful study of geometric Eisenstein functors in the Fargues-Scholze context (see Remark \ref{rem: splittingofSemiorthogonalDecomposition} for a discussion of the relation).
Similarly, we are able to show unconditionally in certain special cases that the tame categorical Langlands equivalence of Zhu intertwines the geometric Eisenstein functors studied in \cite{HamGeomES,HHS,HIDualCpx} with spectral Eisenstein functors, by using our functor to compare the two (see Proposition \ref{prop: EisensteinCompatability}).

\subsection{Technical introduction}{\label{s: Technical Introduction}}
As above, we fix $E$ a non-Archimedean local field of residue characteristic $p$.
We let $O_E\subseteq E$ denote the ring of integers.
We fix $\pi\in O_E$ a uniformizer, an isomorphism $\bbF_q\simeq O_E/\pi$ where $q$ is a power of $p$, and an algebraic closure $k=\overline{\bbF}_q$ of $\bbF_q$. 
We let $\breve{E}$ denote the $\pi$-adically completed maximal unramified extension of $E$.
We let $W_E$ denote the Weil group of $E$, $I_E\subseteq W_E$ the inertia group and $\Gamma_E\supseteq W_E$ the Galois group.
We fix $G/E$ a quasi-split connected reductive group, we fix $T\subseteq B\subseteq G$ a maximally split rational torus and a rational Borel. 
We let $\calI$ denote the parahoric group scheme associated with an Iwahori subgroup $\calI(O_E)\subseteq G(E)$. 
Since Lurie's formalism of $\infty$-categories has become widespread in the CLLC, we adapt our terminology accordingly:  by category we mean an $\infty$-category, whereas $1$-truncated categories are called ordinary categories. 
We fix a prime number $\ell\neq p$, and we let $\Lambda$ denote a fixed torsion $\bbZ_\ell$-algebra (i.e., $\ell^n\cdot \Lambda=0$ for some $n$). 
We let $\on{Mod}_\Lambda$ denote the (derived) category of $\Lambda$-modules. 
We let $\LinCat_\Lambda$ denote the category of presentable stable $\Lambda$-linear categories.

We let $\PSch^\aff$ denote the category of perfect affine $k$-schemes and let $\PSch^{\aff,\pfp}\subseteq \PSch^{\aff}$ denote the full subcategory of those affine schemes that are perfectly finitely presented over $k$.
We let $\PreStk$ denote the category of small accessible presheaves on $\PSch^{\aff}$ with values in anima (or spaces). 
We let $\Perf^\aff$ denote the category of affinoid perfectoid spaces over $k$, and we let $\on{Std}\subseteq \Perf^\aff$ denote the full subcategory of those affinoid perfectoid spaces that are strictly totally disconnected.
We consider $\Perf^\aff$ as an ordinary site endowed with the v-topology, and we let $\AnStk_v$ denote the category of small v-sheaves with values in anima.
We now review the precise definition of the local Langlands categories of Fargues--Scholze and Zhu.
\subsubsection{Zhu's construction}
For any $X\in \PreStk$, Zhu considers a category $\Shv^!(X)\in\LinCat_\Lambda$ of \'etale sheaves that is often featured in geometric Langlands, and goes back at least to the work of Raskin \cite{raskin-placid} and Bouthier--Kazhdan--Varshavsky \cite{perverse-sheaves-infinite-dimensional-stacks} (see \cite[Remark 10.70]{Zhu25} for further references).
 When specialized to $\ICG$ and \'etale sheaves, this category will give rise to Zhu's version of the automorphic local Langlands category.
We now review this construction.

For every $X\in \PSch^{\aff,\pfp}$, we have an ordinary \'etale site $X_\et$, and we may consider the category of sheaves with values in $\Mod_\Lambda$, which we denote by $\calD(X_\et,\Lambda)$.
This coincides with the derived category of its abelian heart; in other words, it is given by the $\infty$-categorical enhancement of the ordinary triangulated derived category attached to the abelian category of $\Lambda$-modules on $X_\et$.
This category is compactly generated, and for every morphism $[f:X\to Y]\in \PSch^{\aff,\pfp}$ one has a functor 
\[f^!:\calD(Y_\et,\Lambda)\to \calD(X_\et,\Lambda).\]
Moreover, this functor preserves compact objects (i.e., bounded complexes with constructible cohomology) and restricts to a functor 
\[f^{!,\omega}:\calD(Y_\et,\Lambda)^\omega\to \calD(X_\et,\Lambda)^\omega.\]
We may interpret the data discussed above as a functor with values in small $\Lambda$-linear categories 
\[\calD_\cons^!:(\PSch^{\aff,\pfp})^\op\to \LinCat^\sm_\Lambda.\] 
We can let $\Shv^!$ denote the ind-extension of $\calD_\cons^!$ which we interpret as a functor 
\[\Shv^!:(\PSch^{\aff,\pfp})^\op\to \LinCat_\Lambda.\]
For a map $f:X\to Y$, we get a functor 
\[\Shv^!([f:X\to Y])=:\ul{f^!}:\Shv^!(Y)\to \Shv^!(X_\et)\]
defined as the Ind-extension of $f^{!,\omega}$. 
This is naturally equivalent to the mapping 
\[{f^!}:\calD(Y_\et)\to \calD(X_\et),\]
but we will keep this change of notation in what follows since we will extend the domain of definition of $\ul{f^!}$ in a way that it will no longer necessarily agree with $f^!$. 
We perform this extension in two steps:
\begin{enumerate}
	\item We let the value of $\Shv^!$ on $\PSch^{\aff}$ be given by the right Kan extension along the inclusion $\PSch^{\aff,\pfp}\subseteq \PSch^\aff$.
		In other words, if $\Spec A=\varprojlim \Spec A_i$, with each $A_i$ pfp over $k$, then 
		\[\Shv^!(\Spec A)\simeq \varinjlim_{\ul{f^!}} \Shv^!(\Spec A_i).\]
	\item We let the value of $\Shv^!$ on $\PreStk$ be given by the left Kan extension along the Yoneda embedding $\PSch^\aff\subseteq \PreStk$.
		In other words, if $X\in \PreStk$, then 
		\[\Shv^!(X)\simeq \varprojlim_{\substack{\Spec A \to X\\ \ul{f^!}}} \Shv^!(\Spec A).\]
\end{enumerate}

Recall the loop group $LG\in \PreStk$,
\[LG:(\PSch^\aff)^\op\to \Ani\]
with formula
\[LG(\Spec R):=G(\bbW(R)[\frac{1}{\pi}]), \]
where $\bbW(R)$ denotes the $O_E$-Witt vectors. 
It comes with an automorphism 
\[\varphi:LG\to LG\]
extracted from $q$-Frobenius on $R$.
We consider the Kottwitz stack parametrizing isocrystals with $G$-structure, $\ICG\in \PreStk$, which admits a description as the \'etale stack quotient by the $\varphi$-conjugation action
\[\ICG:=\frac{LG}{\on{Ad}_\varphi LG}.\]

\begin{definition}
	Zhu's version of the local Langlands category is given by
	\[\calD_{\on{LLC}}:=\Shv^!(\ICG).\]	
\end{definition}

\subsubsection{Fargues--Scholze's construction.}
For any $X\in \AnStk_v$ one attaches a presentable stable $\Lambda$-linear category $\calD^\an_\Lambda(X)\in \LinCat_\Lambda$ of $\Lambda$-\'etale sheaves whose definition we shall recall.
For every $X\in \on{Std}$ a strictly totally disconnected space in the sense of \cite[Definition~1.14]{Sch17}, we have an ordinary \'etale site $X_\et$, and, as above, one attaches the category of sheaves with values in $\Lambda$-modules that we denote by $\calD(X_\et,\Lambda)$, which is essentially sheaves on the underlying topological space of $X$ (which at the level of connected components is just a profinite set). 
From pure site theoretic considerations, for every map $[f:X\to Y]\in \on{Std}$ one has a pullback functor 
\[f^*:\calD(Y_\et,\Lambda)\to \calD(X_\et,\Lambda), \]
which we may consider as a functor 
\[\calD^\an_\Lambda:(\on{Std})^\op\to \LinCat_\Lambda.\]
As above, we extend the domain of definition of $\calD^\an_\Lambda$.
More precisely, 
\[\calD^\an_\Lambda: \AnStk_v\to \LinCat_\Lambda\]
is the left Kan extension along the Yoneda embedding $\on{Std}\subseteq \AnStk_v$.

Recall that, to any $\Spa(R,R^+)\in \Perf^\aff$, one attaches a relative adic Fargues--Fontaine curve $X_{\on{FF},R}$.
This is a sousperfectoid adic space with formula
\[X_{\FF,R}:=(\Spa(\bbW(R^+))\setminus V(\pi [\varpi]))/\varphi.\]
Here $\varpi\in R^+$ is a choice of pseudo-uniformizer, $[\varpi]\in \bbW(R^+)$ is a Teichm\"uller lift, and $\varphi$ is the action induced from $q$-Frobenius on $R^+$.
One considers the moduli stack of $G$-bundles on the Fargues--Fontaine curve, $\Bun_G\in \AnStk_v$. 
Namely, 
\[\Bun_G:\Perf^\aff\to \Ani\]
with formula
\[\Bun_G(\Spa(R,R^+)):= \{G\text{-bundles on } X_{\FF,R}\}.\]

\begin{definition}
	Fargues--Scholze's version of the local Langlands category is given by
	\[\calD_{\on{LLC}}:=\calD^\an_\Lambda(\Bun_G).\]	
\end{definition}

\subsubsection{The main theorem.}
The following is our main theorem
\begin{theorem}[{\cref{thm: MainTheoremPartial}(1)}]
	\label{intro: main theorem}
	Let $G$ be a quasi-split reductive group over $E$.
	If $\Lambda$ is an $\ell$-torsion $\bbZ_\ell$-algebra, then the categories $\Shv^!(\ICG)$ and $\calD^\an_\Lambda(\Bun_G)$ are equivalent.	
	More precisely, we construct an equivalence 
	\[\pitch:\Shv^!(\ICG)\xrightarrow{\simeq} \calD^\an_\Lambda(\Bun_G).\]
\end{theorem}
In the rest of this introduction, we will do the following.
\begin{itemize}
	\item[a)] Discuss the finer properties of our equivalence that are established in this article.
	\item[b)] Sketch a construction of $\pitch$.
	\item[c)] Summarize the strategy and discuss the methods we employ to show that $\pitch$ is an equivalence.
\end{itemize}
 
\subsubsection{Finer properties of $\pitch$.}
Recall the Kottwitz set 
\[B(G):=\frac{G(\breve{E})}{\on{Ad}_\varphi G(\breve{E})} \]
from \cite[\S3]{KottwitzII}.
It identifies with the set of isomorphism classes of isocrystals over $k$.
The set $B(G)$ comes equipped with two invariants, the so-called Newton map and the so-called Kottwitz map, which gives an injective map of sets 
\[(\nu,\kappa):B(G)\to (X_*(T)\otimes \bbQ)^{\Gamma_E,+}\times \pi_1(G)_{\Gamma_E},\]
see \cite[(4.13.1)]{KottwitzII}. The right hand side has a partial order, which $B(G)$ inherits. 
Here $(\nu_1,\kappa_1)\leq (\nu_2,\kappa_2)$ whenever $\kappa_1=\kappa_2$ and $\nu_2- \nu_1$ is a non-negative rational linear combination of positive coroots.
Moreover, this partial order is down-finite in the sense that the set $Z_{\leq b}:=\{b'\in B(G)\mid b'\leq b\}$ is a finite set, see \cite[Proposition 2.4.(iii)]{RR_96}.

One can endow $B(G)$ with its partial order topology where the up-closed sets $U_{\geq b}=\{b'\in B(G)\mid b'\geq b\}$ are a basis of open sets generating the topology.
We will denote this topological space also by $B(G)$. 
Similarly, one can consider the opposite order topology in which the sets $Z_{\leq b}:=\{b'\in B(G)\mid b'\leq b\}$ are declared to be a basis of open subsets generating the topology, we denote this topological space as $B(G)^\op$.
Recall that we have homeomorphisms
\[ |\ICG| \simeq B(G)\]
and 
\[ |\Bun_G| \simeq B(G)^\op.\]

The first homeomorphism is due to He \cite[Theorem 2.12]{He_hecke_padic}, and the second due to Viehmann \cite[Theorem 1.1]{Viehamann_Newton_BdR}.
This reversal of topologies explains plenty of phenomena concerning our functor $\pitch$. 
Using this topological input, for every finite closed subset $Z\subseteq B(G)$ we obtain a closed substack 
\[i_Z:\ICG_Z\subseteq \ICG\]
and an open substack 
\[j_Z:\Bun_G^Z\subseteq \Bun_G.\]
Overall, we get fully-faithful functors
\[\ul{i_{Z,*}}:\Shv^!(\ICG_Z)\to \Shv^!(\ICG) \text{ and } j_{Z,!}:\calD^\an_\Lambda(\Bun^Z_G)\to \calD^\an_\Lambda(\Bun_G).\]
We let $\Shv^!(\ICG_Z)_{*}\subseteq \Shv^!(\ICG)$ and $\calD^\an_\Lambda(\Bun_G^{Z})_{!} \subseteq \calD^\an_\Lambda(\Bun_G)$ denote the essential images of these functors.
\begin{proposition}[{\cref{left orthogonality is shown}(3)}]
	\label{intro:The B(G)-filtration}
	For every finite closed subset $Z\subseteq B(G)$, we have that
	\[\pitch(\Shv^!(\ICG_Z)_{*}) \subseteq \calD^\an_\Lambda(\Bun_G^{Z})_{!},\]
	and consequently there exists a unique (up to contractible choice) functor  
	\[\pitch_Z:\Shv^!(\ICG_Z)\to \calD^\an_\Lambda(\Bun_G^{Z}),\]
    making the following diagram 
    \[
    \begin{tikzcd}
        \Shv^!(\ICG_Z)  \ar{r}{\ul{i_{Z*}}}\ar{d}{\pitch_Z} & \Shv^!(\ICG) \ar{d}{\pitch} \\
        \calD^\an_\Lambda(\Bun_G^Z) \ar{r}{j_{Z!}} & \calD^\an_\Lambda(\Bun_G)
    \end{tikzcd}
    \]
    commute.
\end{proposition}

In \S \ref{sec: B(G)-filt}, we recall the theory of recollements and semi-orthogonal decompositions following \cite[Appendix~A.8]{LurieHigherAlgebra} and \cite{ayala2023stratifiednoncommutativegeometry}, relating it to the theory of presentable Nagata 2-functor formalisms in \cite{dauser2024uniquenesssixfunctorformalisms}.
We call a finite subset $S\subseteq B(G)$ \textit{convex} if it can be written (non-uniquely) in the form $S=Z_1\setminus Z_2$, with $Z_2\subseteq Z_1\subseteq B(G)$ finite closed subsets.
It follows formally from this theory that for every finite convex subset $S\subseteq B(G)$ one obtains a functor 
\[\pitch_S:\Shv^!(\ICG_S)\to \calD^\an_\Lambda(\Bun_G^S),\]
irrespective of the choice of presentation $Z_1\setminus Z_2$, that fits in the following commutative diagram

\[
\begin{tikzcd}
	\Shv^!(\ICG_{Z_2}) \arrow[r, "\pitch_{Z_2}"] \arrow[d]   & \calD^\an_\Lambda(\Bun_G^{Z_2}) \arrow[d] \\
	\Shv^!(\ICG_{Z_1}) \arrow[r, "\pitch_{Z_1}"] \arrow[d] & \calD^\an_\Lambda(\Bun_G^{Z_1}) \arrow[d] \\
	\Shv^!(\ICG_{S}) \arrow[r, "\pitch_S"] &  \calD^\an_\Lambda(\Bun_G^{S}).
\end{tikzcd}
\]
where the vertical arrows yield a Verdier quotient presentation for the bottom line.
Here $\ICG_S$ (respectively $\Bun_G^S$) denotes the locally closed substack of isocrystals (respectively $G$-bundles) whose value on each geometric point is of type of $b$ with $b\in S$. 

Since the categories participate in various recollements, the formalism of semi-orthogonal decompositions provides additional functors which we shall clarify below.
We explain the $\ICG$ setup first since it is more intuitive.
For a closed subset $Z\subseteq B(G)$, the fully faithful functor $\ul{i_{Z,*}}=\ul{i_{Z,!}}$ (associated with the closed immersion $i_Z:\ICG_Z\subseteq \ICG$) admits left and right adjoints i.e., we have
\[\ul{i_Z^*}\dashv \ul{i_{Z,*}} \dashv \ul{i_Z^!}.\]
When $S=Z_1\setminus Z_2$, we have an open immersion $i_{S,Z_1}:\ICG_S\to \ICG_{Z_1}$ and an identification 
\[ \ul{i^!_{S,Z_1}}=\ul{i^*_{S,Z_1}}:\Shv^!(\ICG)_{S,*} \xrightarrow{\simeq} \Shv^!(\ICG_S).\]
Similarly, this functor fits into a sequence of adjoint functors
\[\ul{i_{S,Z_1,!}}\dashv \ul{i^!_{S,Z_1}} \dashv \ul{i_{S,Z_1,*}}.\]
Composing the functors considered above, we obtain for all finite convex $S\subseteq B(G)$ functors
\[\ul{i_{S,!}},\ul{i_{S,*}}:\Shv^!(\ICG_S)\to \Shv^!(\ICG) \text{ and } \ul{i^!_{S}},\ul{i^*_{S}}: \Shv^!(\ICG)\to  \Shv^!(\ICG_S)  \] 
participating in adjunctions
\[\ul{i_{S,!}} \dashv \ul{i^!_{S}} \text{ and } \ul{i^*_{S}} \dashv \ul{i_{S,*}}.\] 
Moreover, these functors do not depend on the presentation $S=Z_1\setminus Z_2$.
Actually, something else is true, namely if $S_1\subseteq S_2$ is an inclusion of finite convex subsets we get analogous functors

\[\ul{i_{S_1,S_2,!}},\ul{i_{S_1,S_2,*}}:\Shv^!(\ICG_{S_1})\to \Shv^!(\ICG_{S_2}) \text{ and } \ul{i^!_{S_1,S_2}},\ul{i^*_{S_1,S_2}}: \Shv^!(\ICG_{S_2})\to  \Shv^!(\ICG_{S_1})  \] 
participating in adjunctions
\[\ul{i_{S_1,S_2,!}} \dashv \ul{i^!_{S_1,S_2}} \text{ and } \ul{i^*_{S_1,S_2}} \dashv \ul{i_{S_1,S_2,*}},\] 
and compatible with composition in the expected manner along inclusions $S_1\subseteq S_2\subseteq S_3$.
One can capture all of these compatibilities and coherences into a so called semi-orthogonal decomposition (see \Cref{def: semi-orthogonal}).

We let $\Convex_{B(G)}$ denote the category whose objects are the finite convex subsets of $B(G)$ and the morphisms are inclusions.
\begin{proposition}[{See \S\ref{subsection semi on isoc}}]
	There is a semi-orthogonal decomposition $\bbS_{\ICG}$ of $\ICG$ with respect to the partially ordered set $B(G)$ (i.e., a presentable 2-functor formalism satisfying the axioms of \Cref{def: semi-orthogonal})
	of the form
\[\bbS_{\ICG}:\Corr(\Convex_{B(G)},\on{All})\to \LinCat_\Lambda\]
with
\[\bbS_{\ICG}(V):=\Shv^!(\ICG_V)\]
and 
\[\bbS_{\ICG}([V\leftarrow W \rightarrow U]):= \ul{i_{WU!}} \circ \ul{i^*_{WV}}.\]
\end{proposition}

Understanding the structure that should come from this semi-orthogonal decomposition via the equivalence $\pitch: \Shv^{!}(\ICG) \simeq \calD^\an_\Lambda(\Bun_G)$ is more subtle, but this is to be expected by the reversal of topologies.
Indeed, $|\Bun_G|\simeq B(G)^\op$, and there is a natural semi-orthogonal decomposition with respect to $B(G)^{\op}$.
However, this is not the one we are looking for. In particular, we are trying to describe a decomposition with respect to $B(G)$ and not with respect to $B(G)^\op$.
 This will be explained by the existence of certain exceptional adjoints, which had already featured in the work of Fargues--Scholze \cite{FS21}. 
More precisely, for a closed subset $Z\subseteq B(G)$, the reversal of topologies will give us an open immersion $j_Z:\Bun_G^Z\to \Bun_G$ and the functor $j_{Z!}$ admits the usual right adjoint $j_Z^!=j_Z^*$.  
Exceptionally, it also admits a left-adjoint which is not guaranteed by the usual $6$-functor formalism (see \S\ref{ss: theanalyticcategory}). 
We will denote the left adjoint to $j_{Z!}$ by $j_Z^{\flat}$. 
In other words, we have a sequence of adjoint functors
\[j_Z^{\flat} \dashv j_{Z!} \dashv j^*_Z.\]
To obtain the functor $j^\flat_Z$, recall that Fargues--Scholze construct the so-called Bernstein--Zelevinsky duality
\[\bbD_{\BZ}:\calD^\an_\Lambda(\Bun_G)^\omega \xrightarrow{\simeq} \calD^\an_\Lambda(\Bun_G)^{\omega, \op}.\]
One can show that $\bbD_{\BZ}$ preserves the subcategory $\calD^\an_\Lambda(\Bun^Z_G)^\omega\subseteq \calD^\an_\Lambda(\Bun_G)^\omega$.
This gives rise to a unique functor fitting in the following commutative diagram
\[
	\begin{tikzcd}
		\calD^\an_\Lambda(\Bun^Z_G)^\omega	 \arrow{r}{\bbD_{\BZ,Z}} \arrow{d}{j_{Z,!}^{\omega}}  & \calD^\an_\Lambda(\Bun^Z_G)^{\omega,\op}  \arrow{d}{j_{Z,!}^{\omega, \op}} \\
		\calD^\an_\Lambda(\Bun_G)^\omega	 \arrow{r}{\bbD_{\BZ}} & \calD^\an_\Lambda(\Bun_G)^{\omega,\op}.
	\end{tikzcd}
\]
Then $j^\flat_Z$ can be computed by the formula 
\[j^\flat_Z:= \on{Ind}(\bbD^Z_{\BZ} j^{*,\omega}_Z \bbD_{\BZ}).\]

Now, when we consider a finite convex subset with presentation $S=Z_1\setminus Z_2$, the reversal of topologies gives us a closed immersion $j_{SZ_1}:\Bun^S_G\to \Bun^{Z_1}_G$, and we get a canonical map 
\[j^*_{SZ_1}:\calD^\an_\Lambda(\Bun^{Z_1}_G)\rightarrow \calD^\an_\Lambda(\Bun^S_G).\]
This functor has an evident right adjoint since $j_{SZ_1*}=j_{SZ_1!}$ and it also has an exceptional left adjoint $j_{S,Z_1,\sharp}$.
Hence, we have a sequence of adjoint functors 
\[j_{SZ_1\sharp}\dashv j^*_{SZ_1} \dashv j_{SZ_1*}.\]
Again, the exceptional left adjoint functor $j_{S,Z_1,\sharp}$ is given by the formula
\[j_{SZ_1\sharp}:= \on{Ind}(\bbD_{\BZ} j^{\omega}_{S,Z_1,!} \bbD_{\BZ,S}),\]
where $\bbD_{\BZ,S}$ is the unique functor filling in the following commutative diagram
\[
	\begin{tikzcd}
		\calD^\an_\Lambda(\Bun^{Z_1}_G)^\omega	 \arrow{r}{\bbD_{\BZ,Z_1}} \arrow{d}{j^{*,\omega}_{SZ_1}}  & \calD^\an_\Lambda(\Bun^{Z_1}_G)^{\omega,\op}  \arrow{d}{j_{S,Z}^{*,\omega,\op}} \\
		\calD^\an_\Lambda(\Bun^S_G)^\omega	 \arrow{r}{\bbD_{\BZ,S}} & \calD^\an_\Lambda(\Bun^S_G)^{\omega,\op},
	\end{tikzcd}
\]
whose existence comes from the functoriality of Verdier quotients.

As before, one can use composition to define for any finite convex subset $S\subseteq B(G)$ functors 
\[j_{S\sharp},j_{S!}:\calD^\an_\Lambda(\Bun^S_G)\to \calD^\an_\Lambda(\Bun_G) \text{ and } j_S^*, j^\flat_S: \calD^\an_\Lambda(\Bun_G)\to  \calD^\an_\Lambda(\Bun^S_G)  \] 
participating in adjunctions
\[j_{S\sharp} \dashv j_S^* \text{ and } j_S^\flat  \dashv j_{S!}.\] 
As above, we also have a version for inclusions of finite convex subsets $S_1\subseteq S_2$, as follows
\[j_{S_1S_2\sharp},j_{S_1S_2!}:\calD^\an_\Lambda(\Bun^{S_1}_G)\to \calD^\an_\Lambda(\Bun^{S_2}_G) \text{ and } j_{S_1S_2}^*, j^\flat_{S_1S_2}: \calD^\an_\Lambda(\Bun^{S_2}_G)\to  \calD^\an_\Lambda(\Bun^{S_1}_G)  \] 
participating in adjunctions
\[j_{S_1S_2,\sharp} \dashv j_{S_1S_2}^* \text{ and } j_{S_1S_2}^\flat  \dashv j_{S_1S_2,!},\] 
and compatible with composition in the expected manner along inclusions $S_1\subseteq S_2\subseteq S_3$.
Again, we package all of these compatibilities and coherences as follows.

\begin{proposition}[See {\S\ref{sec:semiorthogonalDecOfDBunG}}]{\label{sec: semiorthogonalDecompositionIntro}}
	There is a semi-orthogonal decomposition $\bbS^{\on{ex}}_{\Bun_G}$ of $\Bun_G$ by $B(G)$ (see \Cref{def: semi-orthogonal} for a precise definition)
	of the form
\[\bbS^{\on{ex}}_{\Bun_G}:\Corr(\Convex_{B(G)},\on{All})\to \LinCat_\Lambda\]
with
\[\bbS^{\on{ex}}_{\Bun_G}(V):=\calD^\an_\Lambda(\Bun^V_G)\]
and 
\[\bbS^{\on{ex}}_{\Bun_G}([V\leftarrow W \rightarrow U]):= j_{WU\sharp} \circ j^\flat_{WV}.\]
\end{proposition}

The following theorem upgrades \Cref{intro: main theorem} to explain that $\pitch$ is compatible with a variety of operations. 
\begin{theorem}[{\Cref{thm: MainTheoremPartial}}]\label{correct commutations}
	There is an equivalence of 2-functor formalisms 
	\[\pitch_\bbS:\bbS_{\ICG}\xrightarrow{\simeq}\bbS^{\on{ex}}_{\Bun_G},\]
	such that 
	\[\pitch\simeq \lim_{Z\subseteq B(G)} \pitch_Z\]
	where $Z\subseteq B(G)$ ranges over finite closed subsets of $B(G)$.
\end{theorem}

The non-technical translation of the theorem above is the following.

\begin{corollary}[{\Cref{thm: MainTheoremPartial}}]{\label{cor: MainTheoremPartialIntro}}
	For every finite convex subset $S\subseteq B(G)$, we have an equivalence 	
	\[\pitch_S:\Shv^!(\ICG_S)\xrightarrow{\simeq} \calD^\an_\Lambda(\Bun^S_G),\]
together with commutation formulas
	\begin{multicols}{2}
	\begin{enumerate}
		\item $\pitch\circ \ul{i_{S,!}}\simeq j_{S,\sharp}\circ \pitch_S$
		\item $\pitch_S\circ \ul{i_S^*}\simeq j_S^\flat  \circ \pitch $
		\item $\pitch_S\circ \ul{i_S^!}\simeq j_S^*  \circ \pitch $
		\item $\pitch\circ \ul{i_{S,*}}\simeq j_{S,!}\circ \pitch_S$.
	\end{enumerate}
	\end{multicols}
	Moreover, for a pair of convex subsets $S_1\subseteq S_2$ we also have commutation formulas

	\begin{multicols}{2}
	\begin{enumerate}
		\item $\pitch_{S_2}\circ \ul{i_{S_1,S_2,!}}\simeq j_{S_1,S_2,\sharp}\circ \pitch_{S_1}$
		\item $\pitch_{S_1}\circ \ul{i_{S_1,S_2}^*}\simeq j_{S_1,S_2}^\flat  \circ \pitch_{S_2}$
		\item $\pitch_{S_1}\circ \ul{i_{S_1,S_2}^!}\simeq j_{S_1,S_2}^*  \circ \pitch_{S_2}$
		\item $\pitch_{S_2}\circ \ul{i_{S_1,S_2,*}}\simeq j_{S_1,S_2,!}\circ \pitch_{S_1}$.
	\end{enumerate}
	\end{multicols}
\end{corollary}

Amusingly, the equivalence also provides exceptional right adjoint functors in the $\ICG$ setup.
Namely, in the $\Bun_G$ setup we naturally have adjoint functors $j^*_S\dashv j_{S,*}$ and $j_{S,!}\dashv j_S^!$ induced by the $6$-functor formalism $\calD_{\Lambda}^{\an}$. 
From here, we can deduce that we also have adjunctions $\ul{i^!_S}\dashv \ul{i_{S,\flat}}$ and $\ul{i_{S,*}}\dashv \ul{i_S^\sharp}$.\footnote{In the work of Zhu, the functors $\ul{i_{S,*}}$ and $\ul{i^!_S}$, when they exist, are automatically forced to commute with colimits.
Consequently, one can use the adjoint functor theorem to deduce the existence of these adjoints.}

	\begin{remark}
	For the convenience of the reader, we make the following table of adjunctions that are intertwined by $\pitch$
\begin{table}[h]
\centering
\begin{minipage}{0.45\textwidth}
\centering
\begin{tabular}{ccccc}
$j^\flat$ & $\dashv$ &$j_!$ & $\dashv$ &$j^!$ \\
$\ul{i^*}$ & $\dashv$ &$\ul{i_*}$ &$\dashv$ & $\ul{i^\sharp}$ \\
\end{tabular}
\end{minipage}
\hfill
\begin{minipage}{0.45\textwidth}
\centering
\begin{tabular}{ccccc}
$j_\sharp $& $\dashv$ &$j^* $ & $\dashv$ &$j_* $\\
$\ul{i_!}$ & $\dashv$ &$\ul{i^!} $& $\dashv$ &$\ul{i_\flat} $\\
\end{tabular}
\end{minipage}
\end{table}

	\end{remark}

	A consequence of the intertwining discussed in \Cref{cor: MainTheoremPartialIntro} is the following.

\begin{corollary}[{\Cref{cor: variantsoftheequivalence}}]
The equivalence 
\[\pitch:\Shv^!(\ICG)\xrightarrow{\simeq}\calD^\an_\Lambda(\Bun_G)\]
restricts to equivalences 
\[\pitch^{\on{Adm}}:\Shv^!(\ICG)^\Adm\xrightarrow{\simeq} \calD^\an_\Lambda(\Bun_G)^{\on{ULA}}\]
and 
\[\pitch^{\omega}:\Shv^!(\ICG)^\omega\xrightarrow{\simeq} \calD^\an_\Lambda(\Bun_G)^\omega.\]
\end{corollary}

From the geometric nature of the constructions, both setups come with their own version of perverse t-structure. 
Let us recall them. 
For $b\in B(G)$ we consider $\on{dim}(\Bun_G^b)$, the $\ell$-cohomological dimension of $\Bun_{G}^{b}$. 
We recall that this is always a non-positive number that can be computed by the formula $-\langle 2 \rho_G, \nu_{b} \rangle$ where $2\rho_G$ is the half-sum of positive roots.
Adhering to convention, we let $d_b:=\langle 2 \rho_G, \nu_{b} \rangle$, so that $\on{dim}(\Bun_G^b)=-d_b$.
As usual, the perverse t-structure can be described as 
\[{}^p\calD^{\an,\leq 0}_\Lambda(\Bun_G):=\{A\in \calD^\an_\Lambda(\Bun_G)\mid j_b^* A\in {}^s \calD^{\an,\leq d_b}_\Lambda(\Bun^b_G) \text{ for all } b\in B(G)\}\]
and 
\[{}^p\calD^{\an,\geq 0}_\Lambda(\Bun_G):=\{A\in \calD^\an_\Lambda(\Bun_G)\mid j_b^! A\in {}^s \calD^{\an,\geq d_b}_\Lambda(\Bun^b_G) \text{ for all } b\in B(G)\}, \]
where we use ${}^s\calD^\an_\Lambda$ to denote the standard t-structure (see \Cref{defn: BunGPerverseTStructure}).

Similarly, one would want to define the perverse t-structure in the $\ICG$ setup as 
\[{}^p\Shv^!(\ICG)^{\leq 0}:=\{A\in \Shv^!(\ICG)\mid \ul{i_b^*} A\in {}^s\mathrm{Shv}^!(\ICG_b)^{\leq d_b} \text{ for all } b\in B(G)\}\]
and 
\[{}^p\Shv^!(\ICG)^{\geq 0}:=\{A\in \Shv^!(\ICG)\mid \ul{i_b^!} A\in {}^s\mathrm{Shv}^!(\ICG_b)^{\geq d_b} \text{ for all } b\in B(G)\}.\]
Nevertheless, there is a subtlety with this formula (see \cite[Remark 3.107]{Zhu25} for a discussion of this subtlety). 
Instead, one defines ${}^p\Shv^!(\ICG)^{\leq 0}$ as the subcategory generated under small colimits by the family of objects
\[\{\ul{i_{b,!}}\cInd_K^{G_b(E)}\Lambda[n-d_b] \mid b\in B(G), n\geq 0, K\subseteq G_b(E) \text{ is a pro-p-subgroup}\},\]
but the definition of ${}^p\Shv^!(\ICG)^{\geq 0}$ stays the same (see Definition \ref{defn: IsocPerverseTStructure}).

	Our functor allow us to define exceptional t-structures on either of the perspectives.  

\begin{definition}
	\label{exceptional t-structures}
	We define the exceptional t-structures
	\[({}^e\calD^{\an,\leq 0}_\Lambda(\Bun_G),  {}^e\calD^{\an,\geq 0}_\Lambda(\Bun_G)):= (\pitch ({}^p\Shv^!(\ICG)^{\leq 0}), \pitch( {}^p\Shv^!(\ICG)^{\geq 0}))\]
	and 
	\[( {}^e\Shv^!(\ICG)^{\leq 0},  {}^e\Shv^!(\ICG)^{\geq 0}):= (\pitch^{-1}({}^p\calD^{\an,\leq 0}_\Lambda(\Bun_G)), \pitch^{-1}( {}^p\calD^{\an,\geq 0}_\Lambda(\Bun_G))).\]
\end{definition}

The exceptional t-structures on $\Shv^!(\ICG)$ and $\calD^\an_\Lambda$ of \Cref{exceptional t-structures} agree with the exceptional t-structures considered by Zhu (see \cite[Proposition 3.110]{Zhu25}) and the ($\ell$-torsion version of the) hadal t-structure considered by Hansen \cite[Theorem 1.2.3]{HansenBeijingNotes}\footnote{We note that, with $\ol{\bb{Q}}_{\ell}$-coefficients, Hansen actually proves a stronger claim. 
In particular, he shows that this gives rise to a $t$-structure on the full subcategory of compact objects. 
This is not at all obvious, since showing that the compact generators are preserved under standard truncation requires a non-trivial input.}. 
More precisely, we have the following.

\begin{corollary}[{\Cref{thm: texactnessstatement}}]\label{intro:t-exact}
	We have the following formulas: 
	\begin{enumerate}
		\item ${}^e\calD^{\an,\leq 0}_\Lambda(\Bun_G)$ is the category generated under small colimits by the family 
			\[\{j_{b,\sharp}\cInd_K^{G_b(E)}\Lambda[n-d_b] \mid b\in B(G), n\geq 0, K\subseteq G_b(E) \text{ is a pro-p-subgroup}\},\]
		\item ${}^e\calD^{\an,\geq 0}_\Lambda(\Bun_G)=\{A\in \calD^\an_\Lambda(\Bun_G)\mid j_b^* A\in {}^s \calD^{\an,\geq d_b}_\Lambda(\Bun^b_G) \text{ for all } b\in B(G)\} $
		\item ${}^e\Shv^!(\ICG)^{\leq 0}:=\{A\in \Shv^!(\ICG)\mid \ul{i_b^!} A\in {}^s\Shv^!(\ICG_b)^{\leq d_b} \text{ for all } b\in B(G)\}$
		\item \({}^e\Shv^!(\ICG)^{\geq 0}:=\{A\in \Shv^!(\ICG)\mid \ul{i_b^\sharp} A\in {}^s \Shv^!(\ICG_b)^{\geq d_b} \text{ for all } b\in B(G)\}.\)
	\end{enumerate}
\end{corollary}

\Cref{intro:t-exact} essentially follows from \Cref{correct commutations}.
Indeed, the only additional input to this statement is that the equivalence 
\[\pitch_b:\Shv^!(\ICG_b)\to \calD^\an_\Lambda(\Bun_G^b)\]
is t-exact for the standard t-structure which follows from a direct computation (see \Cref{concrete statement}), which we now explain.

A priori, the categories $\Shv^!(\ICG_b)$ and $\calD^{\an}_\Lambda(\Bun^b_G)$ are abstractly defined.    
For computations, it is sometimes useful to give more concrete descriptions of these categories.
Recall that, after fixing $\dot{b}\in G(\breve{E})$ a representative of $b\in B(G)$ one obtains morphisms $\dot{b}:\ast\to \ICG_b$ and $\dot{b}:\ast\to \Bun_G^b$.
These morphisms produce, through $\ul{\dot{b}^!}$ and $\dot{b}^*$ respectively, t-exact identifications  
\[\Rep (G_{\dot{b}}(E))\simeq_{\dot{b}} \Shv^!(\ICG_b) \quad \text{ and } \quad \Rep (G_{\dot{b}}(E))\simeq_{\dot{b}}\calD^\an_\Lambda(\Bun^b_G).\] 
Here $G_{\dot{b}}(E)$ denotes the $\varphi$-centralizer $G_{\dot{b}}(E)=\{g\in G(\breve{E})\mid g^{-1}\dot{b} \varphi(g)=\dot{b}\}$, and $\Rep (G_{\dot{b}}(E))$ is the (derived) category of smooth representations, where $\varphi$ denotes the Frobenius on $\Breve{E}$. 
This leads to an auto-equivalence $\alpha_{\dot{b}}$ characterized by fitting in a commutative diagram
\[
	\begin{tikzcd}
		\Shv^!(\ICG_b) \arrow{r}{\pitch_b} \arrow{d}{\simeq_{\dot{b}}}  & \calD^\an_\Lambda(\Bun_G^b) \arrow{d}{\simeq_{\dot{b}}} \\
				\Rep (G_{\dot{b}}(E)) \arrow{r}{\alpha_{\dot{b}}} & \Rep (G_{\dot{b}}(E)).
	\end{tikzcd}
	\]

	To make $\alpha_{\dot{b}}$ more explicit, consider the following.
	Both the $\calD^\an_\Lambda(\Bun_G)$ setup and the $\Shv^!(\ICG)$ come with canonically constructed Frobenius algebra structures in the sense of \cite[Definition~4.6.5.1]{LurieHigherAlgebra}). These give rise to identifications 
	\[\bbD_\BZ:\calD^\an_\Lambda(\Bun_G)\xrightarrow{\simeq}\calD^\an_\Lambda(\Bun_G)^\vee\]
	and 
	\[\id_\BZ:\Shv^!(\ICG)\xrightarrow{\simeq}\Shv^!(\ICG)^\vee,\]
    where $(-)^{\vee}$ denotes the dual category.
Moreover, both of these dualities are suitably compatible with semi-orthogonal decompositions, which lead to naturally defined dualities
	\[\bbD_{\BZ,b}:\calD^\an_\Lambda(\Bun^b_G)\xrightarrow{\simeq}\calD^\an_\Lambda(\Bun^b_G)^\vee\]
	and 
	\[\id_{\BZ,b}:\Shv^!(\ICG_b)\xrightarrow{\simeq}\Shv^!(\ICG_b)^\vee.\]

	Recall that $\Rep (G_{\dot{b}}(E))$ is self-dual\footnote{It is even compactly generated.} and that there is a canonical duality introduced by Bernstein 
	\[\bbD_{\coh,G_{\dot{b}}(E)}:\Rep (G_{\dot{b}}(E))\xrightarrow{\simeq} \Rep (G_{\dot{b}}(E))^\vee,\]
    which we refer to as cohomological duality. 
	As it turns out, these dualities fit in commutative diagrams

\[
	\begin{tikzcd}
		\Shv^!(\ICG_b) \arrow{r}{\id_{\BZ,b}} \arrow{d}{\simeq_{\dot{b}}}  & \Shv^!(\ICG_b)^\vee \arrow{d}{\simeq_{\dot{b}}} 
									   & \calD^\an_\Lambda(\Bun_G^b) \arrow{r}{\bbD_{\BZ,b}} \arrow{d}{\simeq_{\dot{b}}}  & \calD^\an_\Lambda(\Bun_G^b)^\vee  \arrow{d}{\simeq_{\dot{b}}} \\
		\Rep (G_{\dot{b}}(E)) \arrow{r}{-\otimes \xi_{\dot{b}}} & \Rep (G_{\dot{b}}(E)) & 
				\Rep (G_{\dot{b}}(E)) \arrow{r}{-\otimes \tau_{\dot{b}}} & \Rep (G_{\dot{b}}(E))
	\end{tikzcd}
	\]
	where $\tau_{\dot{b}}$ and $\xi_{\dot{b}}$ are invertible objects in $\Rep (G_{\dot{b}}(E))$. 
	One can show that $\xi_{\dot{b}}\simeq \delta^{\on{Zhu}}_{\dot{b}}[-2d_b]$ for a certain smooth character $\delta^{\on{Zhu}}_{\dot{b}}:G_{\dot{b}}(E)\to \Lambda^\times$.
Similarly, one can show that $\tau_{\dot{b}}\simeq\delta^{\on{FS}}_{\dot{b}}[2d_b]$ for a smooth characters $G_{\dot{b}}(E)\to \Lambda^\times$, which was explicitly computed in \cite{HIDualCpx}.
This is what we can say about the functor $\alpha_{\dot{b}}(-)$ at the moment.
Recall from \cite{MR3444233} the notion of \textit{weakly unramified character} of $p$-adic reductive group. 
These are those characters that factor through the Kottwitz homomorphism, or alternatively, those whose restriction to parahoric subgroups is trivial.

\begin{proposition}[{\cref{diagram spell pitch out}, \cref{lemma: RestrictiontoIwahoriTrivial}}]
	\label{concrete statement}
	We have the following formula
	\[\alpha_{\dot{b}}(-) \simeq (-\otimes \delta^{\on{Zhu}}_{\dot{b}}\otimes \delta^{\on{FS},-1}_{\dot{b}}).\]
	In particular, $\alpha_{\dot{b}}(-)$ and $\pitch_b$ are t-exact.
Moreover, $\chi_{\dot{b}} := \delta^{\on{Zhu}}_{\dot{b}}\otimes \delta^{\on{FS},-1}_{\dot{b}}$ is a weakly unramified character. 
\end{proposition}

\begin{remark}
	Our expectation is that $\chi_{\dot{b}}$ is always trivial i.e., that $\chi_{\dot{b}}\simeq \mathbbm{1}$.
	One of us (H.) together with Imai computed explicit formulas for the character $\delta^{\on{FS}}_{\dot{b}}$ (see \cite{HIDualCpx}). 
	One can hope that $\delta^{\on{Zhu}}_{\dot{b}}$ would also follow from explicit methods, but we have not attempted to compute it ourselves. 
	Alternatively, we expect Fargues--Scholze's duality $\bbD_\BZ$ to be intertwined with Zhu's duality $\id_\BZ$ under the equivalence $\pitch$.
If this was true, one could use abstract methods to show that $\chi_{\dot{b}}$ must be trivial. 
\end{remark}
\begin{remark}
The fact that the character $\chi_{\dot{b}}$ has trivial restriction to the Iwahori allows us to see that our functor restricts to an equivalence
\[ \pitch^{\tame}: \Shv^{!,\tame}(\ICG,\Lambda) \xrightarrow{\simeq} \calD_{\Lambda}^{\an,\tame}(\Bun_{G}) \]
for the tame part $\calD_{\Lambda}^{\an,\tame}(\Bun_{G}) \subset \calD_{\Lambda}^{\an}(\Bun_{G})$ (see \ref{defn: analyticTameUnipotentCategory}).
Indeed, this will follow essentially from \Cref{correct commutations} and the fact that tensoring by $\chi_{\dot{b}}$ will preserve the tame subcategories since it is weakly unramified.
Similarly, we can deduce variants for the unipotent subcategories (see \Cref{cor: variantsoftheequivalence}). 
\end{remark}

\subsection{Construction of $\pitch$}
To first approximation, the functor $\pitch$ may be understood as being geometrically realized via a correspondence defined by an object $\Bun_{G}^{\mer}$ sitting in a diagram of $v$-stacks
\[
	\begin{tikzcd}
		\Bun_G^\mer \arrow{r}{\sigma} \arrow{d}{\gamma}  & \Bun_G  \\
	\ICG^\diamondsuit,  
	\end{tikzcd}
\]
as introduced in \cite{GIZ25} by two of the authors together with Zillinger. 
Here, for any $X\in \PreStk$, we are letting $X^\diamondsuit\in \AnStk_v$ denote the sheafification of the rule
\[X^{\diamondsuit_\pre}(\Spa(R,R^+)):=X(\Spec R), \]
(see \S \ref{sec: the Diamond Functor}). 
Indeed, as described already in \cite{GIZ25}, one expected to be able to define a functor $\pitch^{\on{naive}}$ along the lines of the formula 
\[\pitch^{\on{naive}}:=\sigma_!\gamma^*c^*, \]
where
\[c^*:\calD_\Lambda^{\sch}(\ICG)\to \calD^\an_\Lambda(\ICG^\diamondsuit)\]
is Scholze's analytification functor (\S \ref{subsec: DiamondFunctor}). 
Here $\calD_\Lambda^{\sch}(\ICG)$ denotes the category of \'etale sheaves obtained from v-descent (see \Cref{defn: DLambdaSch}).
This approach is however too naive for the following two technical reasons:
\begin{enumerate}
	\item[a)] It is not at all clear if the map $\sigma:\Bun_G^\mer\to \Bun_G$ of $v$-stacks is nice enough to support a functor $\sigma_!$ in the $6$-functor formalism $\calD_{\Lambda}^{\an}$, i.e., we do not know if $\sigma$ is $!$-able.
	\item[b)] The precise relationship between $\Shv^!(\ICG)$ the category that Zhu studies and $\calD_\Lambda^{\sch}(\ICG)$ is not at all clear.
		Indeed, both categories are, at heart, constructed from the classical derived category of \'etale sheaves on perfectly finitely presented $k$-schemes.
However, the categorical procedures one uses to define them are of a different nature. 
\end{enumerate}

Let us discuss how we tackle $b)$.
One of the main reasons that the category $\Shv^!(\ICG)$ can be studied is because the stack $\ICG$ is sind-very placid in Zhu's terminology (see \cite[\S 10.6]{Zhu25}).
Informally, this means that it is constructed from perfectly finitely presented $k$-schemes by straightforward categorical operations.
Indeed,
the prefix ``sind-'' is the contraction of the word ``sifted'' and the prefix ``ind-'', and it is designed to remind the reader that $\ICG$ as an object in $\PreStk$ is a sifted colimit of ind-objects each of which is very-placid\footnote{According to Zhu, the suffix ``sind-'' was introduced by Hemo, in the joint unpublished work of Hemo--Zhu that gave rise to \cite{Zhu25}.}. 
More precisely, recall the Newton map
\[\on{Nt}:\Sht^\sch_\calI\to \ICG,\]
where the source space is the stack of $\calI$-shtukas (see \cite[\S 3.1.4]{Zhu25} and \S \ref{sec:isoc_and_shtukas}).
We can form the \v{C}ech nerve with respect to the Newton map to obtain a resolution by so-called Hecke correspondences
\[\Hk^\bullet(\Sht^\sch_\calI)\to \ICG,\]
with 
\[\Hk^n(\Sht^\sch_\calI):= \underbrace{\Sht^\sch_\calI \times_{\ICG} \Sht^\sch_\calI \times_{\ICG} \cdots \times_{\ICG} \Sht^\sch_\calI}_{\text{$n+1$ times}}.\]

Then $\ICG$ is the geometric realization (a sifted colimit) of $\Hk^\bullet(\Sht^\sch_\calI)$, and each $\Hk^n(\Sht^\sch_\calI)$ is an ind-very placid stack.
More precisely, each $\Hk^n(\Sht^\sch_\calI)$ has a presentation as a filtered colimit $\Hk^n(\Sht^\sch_\calI)=\colim_{\beta \in B_n} Z_{n,\beta}$ where each $Z_{n,\beta}$ is a very-placid stack, and the transition maps are closed immersions. 
We refer to the $Z_{n,\beta}$ as \textit{the bounded pieces} of $\ICG$.
When $n=0$, the bounded pieces $Z_{n,\beta}$ can be taken to be of the form $\Sht_{\calI,\leq \mu}^{\sch}$ as we range over dominant cocharacters $\mu\in X_*(T)^+$. 
Overall we get (up to \'etale sheafification) a colimit formula in $\PreStk$ of the form
\[\ICG\simeq \colim_{n\in \Delta^\op} \colim_{\beta\in B_n} Z_{n,\beta}.\]
This also gives a colimit formula at the level of sheaf theories
\[\Shv^!(\ICG)\simeq \colim_{n\in \Delta^{\op}} \colim_{\beta\in B_n} \Shv^!(Z_{n,\beta}),\]
and a colimit formula for the dual categories
\[\Shv^!(\ICG)^\vee\simeq \colim_{n\in \Delta^\op} \colim_{\beta\in B_n} \Shv^!(Z_{n,\beta})^\vee.\]
At this point, we observe that for the very-placid stacks $Z_{n,\beta}$ we have canonical comparison maps
\[\Shv^!(Z_{n,\beta})^\vee \to \calD_\et(Z_{n,\beta})\xrightarrow{c^*} \calD_\et(Z^\diamondsuit_{n,\beta}) ,\]
where the first map compares the $\Shv^!$-theory with the (left-completed) \'etale sheaves on $Z_{n,\beta}$, and the second is Scholze's analytification of \'etale sheaves.
In particular, the dual category $\Shv^!(\ICG)^\vee$ can be written as a colimit of pieces $\Shv^!(Z_{n,\beta})^\vee$ each of which can be analytified. 

Let us discuss how we tackle $a)$. 
In rough terms, what we realized is that even if $\sigma$ is not $!$-able, $\sigma_!$ can still be constructed on sheaves that come from analytification. 
We achieve this by treating $\Bun_G^\mer$ as a sind-object, as it turns out $\Bun_G^\mer$ can be treated as a sind-Artin v-stack.
More precisely, we have a presentation 
\[\Bun_G^\mer\simeq \colim_{n\in \bbN} \colim_{\beta\in B_n} (Z^\diamondsuit_{n,\beta}\times_{\ICG^\diamondsuit} \Bun_G^\mer).\] 
This presentation allows us to consider the map 
\[\sigma_{n,\beta}:Z^\diamondsuit_{n,\beta}\times_{\ICG^\diamondsuit} \Bun_G^\mer\to \Bun_G,\]
by composing the natural projection $Z^\diamondsuit_{n,\beta}\times_{\ICG^\diamondsuit} \Bun_G^\mer\to \Bun_G^\mer$ with $\sigma$.
This suggests that we may define the appropriate version of $\sigma_{!}$ by first defining $\sigma_{n,\beta!}$ and then taking the colimit.

Indeed, in the process of writing \cite{GIZ25} one of us (G.) noticed that the analogue of the map $\sigma$ for the moduli of analytic shtukas was $!$-able.
The following theorem was proved by further developing the theory of kimberlites.
\begin{theorem}[{\cite[Theorem 1.1]{Gle23}}]
The following statements hold.
\begin{enumerate}
	\item Each of $Z^\diamondsuit_{n,\beta}\times_{\ICG^\diamondsuit} \Bun_G^\mer$ is an Artin v-stack.
	\item The maps $\sigma_{n,\beta}$ are representable in locally spatial diamonds.
More precisely, they are fdcs and in particular $!$-able.
\end{enumerate}
\end{theorem}

This theorem allows us to define a functor 
\[\pitch^\vee_{n,\beta}:\Shv^!(Z_{n,\beta})^\vee\to \calD^\an_\Lambda(\Bun_G)\]
by composing the maps
\[\Shv^!(Z_{n,\beta})^\vee\to \calD_\et(Z_{n,\beta})\xrightarrow{c^*}\calD_\et(Z_{n,\beta}^\diamondsuit) \xrightarrow{\gamma^*_{n,\beta}} \calD_\et(Z^\diamondsuit_{n,\beta}\times_{\ICG^\diamondsuit}\Bun_G^\mer)\xrightarrow{\sigma_{n,\beta,!}}\calD^\an_\Lambda(\Bun_G), \]
where $\gamma_{n,\beta}$ denotes the natural base-change of $\Bun_{G}^{\mer} \ra \ICG$.
The next step is to form the colimit 
\begin{equation}
	\label{intro: colimit formula}	
\pitch^\vee:=\colim_{n\in \Delta^\op} \colim_{\beta \in B_n}\pitch^\vee_{n,\beta},
\end{equation}
this is a map 
\[\pitch^\vee:\Shv^!(\ICG)^\vee\to \calD^\an_\Lambda(\Bun_G).\]
Finally, we can precompose with Zhu's duality $\id_{\BZ}:\Shv^!(\ICG)\xrightarrow{\simeq} \Shv^!(\ICG)^\vee$ to obtain our functor
\[\pitch:=\pitch^\vee\circ \id_{\BZ}:\Shv^!(\ICG)\to \calD^\an_\Lambda(\Bun_G).\]

The necessary background to understand the construction of the functor is discussed in \S \ref{sec: analytificiation and 6-functor formalisms}, \S \ref{s:Analyt for stacks section} and \S \ref{Dictionary-with-Zhu}, while the construction is carried in \S \ref{sec:construction_functor}. 
Since we are manipulating ($\infty$)-categories and functors, one has to be careful when it comes to taking colimits of such gadgets, particularly when we deal with stacks and descent.
In particular, one has to be able to verify that the constructions employed are appropriately functorial in an $\infty$-categorical sense.
To this aim, in \S \ref{sec: analytificiation and 6-functor formalisms}, we discuss Scholze's analytification functor in functorial language by encoding analytification through a $6$-functor formalism, as in \cite{HeyerMann}.
To our knowledge, our \Cref{thm: DiamondAnalytificationUltimate} is the first entry in the literature that carefully discusses the compatibility of Scholze's analytification functor $c^*$ with the $!$-lower operation when maps of stacks are involved.
Since we ran into more than a few surprises when thinking about this, we have decided to include a very thorough discussion that spans \S \ref{sec: analytificiation and 6-functor formalisms} and \S \ref{s:Analyt for stacks section}.

Similarly, the purpose of \S \ref{Dictionary-with-Zhu} is to recall Zhu's formalism of cosheaves on ind-sifted-placid stacks, and to explain how it interacts with analytification.
This turns out to also be a subtle topic. 
Indeed, as mentioned above, one does not expect to have an analytification functor of the form 
\[\Shv^!(X)\to \calD^\an_\Lambda(X^\diamondsuit)\]
for $X\in \PreStk$.
Instead, it is more natural to expect a functor of the form
\[\Shv^!(X)^\vee\to \calD^\an_\Lambda(X^\diamondsuit),\]
but for general prestacks even this can fail to exist.
This forces us to carefully discuss conditions under which one can find such a functor. 
On the other hand, in the discussion above, we have used the Newton map 
\[\Nt:\Sht^\sch_\calI\to \ICG\]
to construct $\pitch$, which naturally leads to the question: How does $\pitch$ depend on the choice of the Iwahori $\calI$? The formalism described in \S \ref{Dictionary-with-Zhu} is precisely designed to give us a presentation of the functor $\pitch$ that does not depend on any choices.

\subsection{Why is $\pitch$ an equivalence?}
In the end, the equivalence boils down to specific nearby cycles computations, but we first do some abstract reasoning to isolate which of these computations are key.
We first exploit that our functor $\pitch$ is $B(G)$-filtered (i.e., \Cref{intro:The B(G)-filtration} holds). 
This is a consequence of the following geometric reasoning. 
Recall the following 
	\begin{theorem}
		[{\cite[Theorem 1.1]{GIZ25}}]	\label{intro: cartesian isoc bungmer}
		We have a commutative diagram 
	\begin{center}
	\begin{tikzcd}
		\calM_b \arrow{r} \arrow{d}{q_b} \arrow[bend left]{rr}{\pi_b} & \Bun_G^\mer \arrow{d}{\gamma} \ar{r}{\sigma} & \Bun_G\,.  \\
		\ICG_b^\diamondsuit \arrow{r}{i_b^\diamondsuit} & \ICG^\diamondsuit,
	\end{tikzcd}
	\end{center}
    where the square is Cartesian.
	\end{theorem}
	Here $\calM_b$ are the Fargues--Scholze charts parametrizing filtered bundles with graded piece isomorphic to $\calE_b$ (as a graded bundle).
We call the map $\sigma$ the special Newton polygon map.
The category $\Shv^!(\ICG_Z)$ is generated under colimits by the subcategories 
	\[\{\ul{i_{b,!}}(\Shv^!(\ICG_b))\subseteq \Shv^!(\ICG_Z)\}_{b\in Z},\]
	so it suffices to show that $\pitch(\ul{i_{b,!}}A)\in \calD^\an_\Lambda(\Bun^Z_G)_!$.
	This can be achieved by using \Cref{intro: cartesian isoc bungmer}, proper base change and the fact that $\pi_b(\calM_b)\subseteq \Bun^Z_G$ whenever $b\in Z$ and $Z\subseteq B(G)$ is closed.

As we mentioned above, this already provides purely abstractly functors 
\[\pitch_S:\Shv^!(\ICG_S)\to \calD^\an_\Lambda(\Bun^S_G),\]
for every finite convex subset $S\subseteq B(G)$.
It is also not hard to see that 
\[\pitch=\colim_{\{Z\subseteq B(G) \mid Z \text{ is finite and closed}\}}\pitch_{Z}.\]
In particular, we are easily reduced to showing that $\pitch_Z$ is an equivalence for $Z\subseteq B(G)$ a closed finite subset.
Intuition from a decategorified setting, like filtered vector spaces, would suggest that a filtered functor is an equivalence if and only if it is an equivalence on graded pieces, but this is not true and it is easy to construct counterexamples.
This does not mean that we cannot try to reduce to graded pieces, but it suggests that we should tread carefully.
Concretely, some conditions must be imposed on the left and right orthogonals of the graded pieces, with the following condition already being sufficient.
\begin{proposition}[{\Cref{prop: altexplicatedsemiorthogonaldecompositioncriterion2}}]
	\label{intro: recollement B(G)-filtered}
	Suppose that $\pitch$ is $B(G)$-filtered (i.e., \Cref{intro:The B(G)-filtration} holds) and that the following conditions hold.
	\begin{enumerate}
		\item For every convex subset $S\subseteq B(G)$, $\pitch(\on{Im}(\ul{i_{S, !}}))\subseteq \on{Im}(j_{S,\sharp})$.
		\item For every convex subset $S\subseteq B(G)$ and for all $b\in S$ that are minimal in $S$, i.e., a closed point in $S$, we have that $\pitch_S(\on{Im}(\ul{i_{S\setminus \{b\},S,*}}))\subseteq \on{Im}(j_{S\setminus \{b\},S,!})$.
		\item For every $b\in B(G)$, $\pitch_b:\Shv^!(\ICG_b)\to \calD^\an_\Lambda(\Bun_G^b)$ is an equivalence.
	\end{enumerate}
	Then $\pitch$ is an equivalence.
\end{proposition}
\begin{remark}
We point out that at this stage of the argument we do not have all of the commutation identities of \Cref{correct commutations} since some of them only follow after it is shown that $\pitch$ is an equivalence.
\end{remark}

The next step is to give a more concrete construction of $\pitch_S$ for $S\subseteq B(G)$ finite and convex.
The stack $\ICG_S$ is also sind-very placid and we have a correspondence 
\[
	\begin{tikzcd}
		\Bun_G^{\mer,\{\gamma,\sigma\in S\}} \arrow{r}{\sigma_S} \arrow{d}{\gamma_S}  & \Bun^S_G  \\
	\ICG_S^\diamondsuit.  
	\end{tikzcd}
\]
Here $\Bun_G^{\mer,{\gamma,\sigma}\in S}$ is the locally closed locus in $\Bun_G^\mer$ where both the special and generic Newton polygons lie in $S$.
More precisely, it is defined through the following Cartesian diagram
\[
	\begin{tikzcd}
		\Bun_G^{\mer,\{\gamma,\sigma\in S\}} \arrow{r} \arrow{d}  & \Bun^S_G\times \ICG_S^\diamondsuit  \ar{d}{(j_S,i_S)} \\
		\Bun^\mer_G \ar{r}{(\sigma,\gamma)} &  \Bun_G\times \ICG^\diamondsuit .
	\end{tikzcd}
\]

We can use the same blueprint of using a sind-resolution of $\ICG_S$ and the correspondence above to construct a functor 
\[\pitch'_S:\Shv^!(\ICG_S)\to \calD^\an_\Lambda(\Bun_G^S),\]
and we show that the abstractly defined $\pitch_S$ agrees with $\pitch'_S$.
This already allows us to show item $(3)$ in \Cref{intro: recollement B(G)-filtered} holds. 
Indeed, when $S=\{b\}$ the map 
\[\gamma_b:\Bun_G^{\mer,(b,b)}\to \ICG_b^\diamondsuit\] is an isomorphism of stacks.
Moreover, after fixing $\dot{b}$, we get a commutative diagram
\[
	\begin{tikzcd}
		\Bun_G^{\mer,(b,b)} \arrow{r} \arrow{d}{\simeq}  & \Bun^b_G \ar{d}{\simeq} \\
		  \ast/{\ul{G_{{\dot{b}}}(E)}}  \ar{r} &   \ast/\tilde{G}_{\dot{b}}, 
	\end{tikzcd}
\]
where $\tilde{G}_{\dot{b}} \simeq \ul{G_{{\dot{b}}}(E)} \ltimes \tilde{G}_{\dot{b}}^{> 0}$, is the group of automorphism of the $G$-bundle corresponding to $\dot{b}$ (see \Cref{thm: BunGGeometricFacts} (4)).
From here, it follows that 
\[\sigma_{b,!}: \calD^\an_\Lambda(\Bun_G^{\mer,(b,b)})\to \calD^\an_\Lambda(\Bun_G^b)\]
is an equivalence, and this will imply that $\pitch_b$ is also an equivalence. \\

Let us discuss condition (1) of \Cref{intro: recollement B(G)-filtered} holds.
Given a finite convex set $S\subseteq B(G)$ we may consider its closure $S\subseteq Z$ which is also a finite closed subset $Z\subseteq B(G)$. 
It suffices to show that for all $\delta\in S$ and $b\in Z\setminus S$ 
\[\on{Hom}_{\calD^\an_\Lambda(\Bun_G)}(\pitch(\ul{i_{\delta !}}A), j_{b !}B)=0\]  
for all $A\in \Shv^!(\ICG_{\delta})$, and $B\in \calD^\an_\Lambda(\Bun_G^b)$.
	Indeed, $L\in \calD^\an_\Lambda(\Bun^Z_G)_!$ lies in $\on{Im}(j_{S,Z,\sharp})$ if and only if $L$ is left-orthogonal to all objects in $\on{Im}(j_{b !})$ for all $b\in Z\setminus S$.

	One can use \Cref{intro: cartesian isoc bungmer} and the fact that the map 
	\[\pi_\delta:\calM_\delta\to \Bun_G\]
	is smooth to show that
	\[\on{Hom}_{\calD^\an_\Lambda(\Bun_G)}(\pitch(\ul{i_{\delta !}}\cInd_K^{G_\delta(E)}\Lambda), j_{b !}B)=\Gamma(\widetilde{\calM}_\delta,j^\delta_{b !}B')^K,\]  
	for some $B'$ and where we consider the locally closed immersion 
	\[j^\delta_{b}:\widetilde{\calM}^{\sigma = b}_\delta\to \widetilde{\calM}_\delta, \text{ and the $\ul{G_\delta(E)}$-torsor } \widetilde{\calM}_\delta\to {\calM}_\delta. \] 
	It follows from the results of Fargues--Scholze (see the proof of \cite[Proposition V.4.2]{FS21}) that $\Gamma(\widetilde{\calM}_\delta,j^\delta_{b !}B')$ always vanishes regardless of what $B'$ is. 
	This is because the spaces $\widetilde{\calM}_\delta$ are ``henselian'' towards their unique non-analytic closed point $\widetilde{\calM}^{\sigma=\delta}_\delta\subseteq \widetilde{\calM}_\delta$ (cf. \cite[Proposition~V.4.2]{FS21} and its proof).\footnote{More conceptually, this is a version of geometric second adjointness in the sense of Theorem \cite[Theorem~1.2.1 (3)]{HHS}, via the relationship between $\mathcal{M}_{\delta}$ and moduli spaces of parabolic structures on $G$-bundles (see \cite[Example~V.3.4]{FS21}). On a similar note, this second adjointness is also precisely the result guaranteeing the existence of the exceptional left adjoints $j_{b\sharp} \dashv j_{b}^{*}$ required to construct the semi-orthogonal decomposition in \S \ref{sec: semiorthogonalDecompositionIntro}, which we are essentially comparing in this argument.}

The proof strategy to show that condition (2) of \Cref{intro: recollement B(G)-filtered} holds follows a similar reasoning.
In this case, one wants to show that if $S$ is a finite convex subset, $b\in S$ is minimal and $\delta\in S\setminus \{b\}$, then 
\[\on{Hom}_{\calD^\an_\Lambda(\Bun^S_G)}(j_{b,S !}B, \pitch_S(\ul{i_{\delta, S*}}A))=0,\]  
for all $A$ and $B$. Since $b\in S$ is minimal, the inclusion $\Bun^b_G\to \Bun^S_G$ is an open immersion, so it suffices to show that for all $A$
\[j_{b,S}^*\pitch_S(\ul{i_{\delta,S,*}}A))=0.\]  

\begin{figure}[htbp]
  \centering
  \includegraphics[
    page=1,
    clip,
    width=0.6\textwidth
    ]{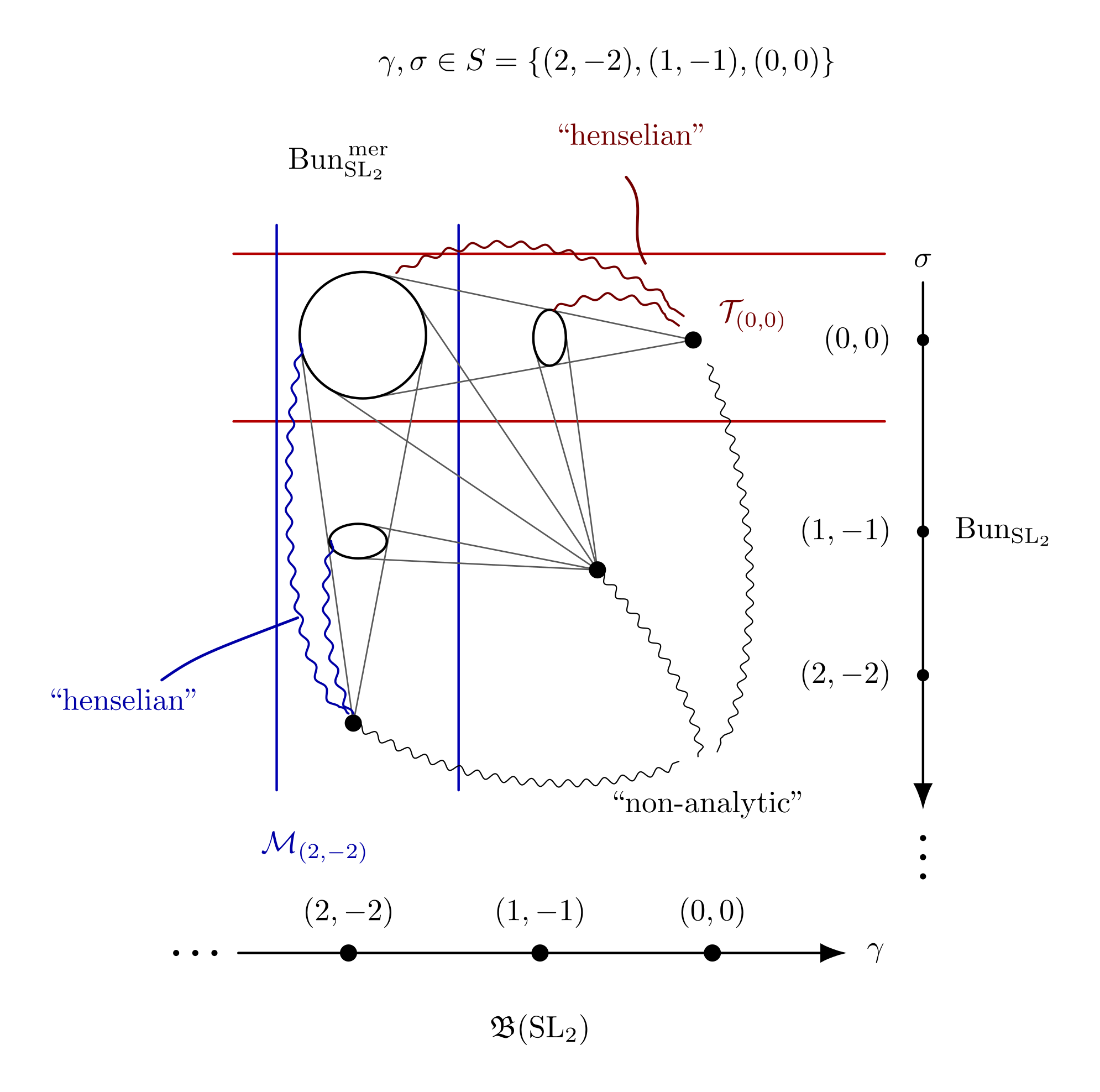}
    \caption{This picture depicts the $(\gamma,\sigma)$-strata of $\Bun_{\on{SL}_2}^{\mer}$, where $\gamma,\sigma$ run over the set $S=\{(0,0),(1,-1),(2,-2)\}$.
There are a total of 6 strata indexed by $(\sigma,\gamma)$ subject to the condition $\gamma \geq \sigma$. 
	   If $\gamma> \sigma$, then the strata is analytic: they are depicted by a circle and two smaller ellipses in the figure above. The remaining strata are given by $\gamma=\sigma$ and are thus non-analytic: we depict them by bullet points in the above picture.
    The blue parallel lines surround the locus $\calM_{(2,-2)}$ inside $\Bun^{\mer}_{\on{SL}_2}$ where $\sigma$ varies through $S$ and $\gamma=(2,-2)$, and whose unique non-analytic point is indexed by $\sigma=\gamma=(2,-2)$, corresponding to the bottom left vertex of the vertical cone. 
    Similarly, the red parallel lines surround the locus $\calT_{(0,0)}$ of $\Bun^{\mer}_{\on{SL}_2}$ where $\sigma=(0,0)$ and $\gamma$ runs through $S$, and whose unique non-analytic point is indexed by $\sigma=\gamma=(0,0)$, corresponding to the top right vertex of the horizontal cone. 
    The square enclosed by the red and blue lines contains a circle depicting the stratum given by $\sigma=(0,0)$ and $\gamma=(2,-2)$.
     }
  \label{fig:pdf-cropped}
\end{figure}

There were two key facts that allowed us to show that condition (1) of \Cref{intro: recollement B(G)-filtered} was satisfied. 
Namely, that the map $\Bun_G^{\mer,\gamma=\delta}=\calM_\delta \to \Bun_G$ is smooth and that the source of this map is henselian towards its unique non-analytic closed point.
In other words, the $\bar{\bbF}_p$-fibers of the map
\[\gamma:\Bun_G^\mer\to \ICG^\diamondsuit\]
satisfy a henselianity property.

For condition (2) of \Cref{intro: recollement B(G)-filtered}, the relevant computation goes through the correspondence
\[
	\begin{tikzcd}
		\calT_b^{\gamma\in S}:=\Bun_G^{\mer,\{\gamma\in S,\sigma=b\}} \arrow{r}{\sigma_b} \arrow{d}{\gamma_S}  & \Bun^b_G  \\
		\ICG^\diamondsuit_S    .
	\end{tikzcd}
\]
Here $\calT_b^{\gamma\in S}$ mimics $\calM_\delta$, where the role of the Newton polygons has been reversed in the definition.  
By analogy, one may hope that the map $\gamma_S$ is also smooth, and that the $\bar{\bbF}_p$-fibers of 
\[\sigma:\Bun_G^\mer\to \Bun_G\]
also satisfy a henselianity property.
Unfortunately, $\gamma_S$ is not smooth, since different fibers have different dimensions.
Nevertheless, the map $\gamma_S:\calT_b^{\gamma\in S}\to \ICG^\diamondsuit_S$ still satisfies a shadow of smooth base change, at least for sheaves coming from analytification (see \Cref{base change formulas}).
That is, for $\delta\in S\setminus \{b\}$ and a compact $K\subseteq G_\delta(E)$, we have a formula
\[j_{b,S}^*\pitch_S(\ul{i_{\delta,S,*}}\cInd_K^{G_\delta(E)}\Lambda)=\sigma_{b,!}k_{\delta,*}\gamma_\delta^*\cInd_K^{G_\delta(E)} \Lambda,\]
where $k_\delta$ fits in the following Cartesian diagram
\[
\begin{tikzcd}
	\calT^{\gamma=\delta}_b\arrow{r}{k_\delta} \arrow{d}{\gamma_\delta}  & \calT^{\gamma\in S}_b \arrow{d}{\gamma_S} \\
	\ICG^\diamondsuit_\delta \arrow{r}{i^\diamondsuit_{\delta,S}} & \ICG^\diamondsuit_S.
\end{tikzcd}
\]
To finish showing 
\[j_{b,S}^*\pitch_S(\ul{i_{\delta,S,*}}\cInd_K^{G_\delta(E)}\Lambda)=\sigma_{b,!}k_{\delta,*}\gamma_\delta^*\cInd_K^{G_\delta(E)} \Lambda\simeq 0\]
we observe that $\sigma_{b,!}k_{\delta,*}C \simeq 0$ regardless of what $C$ is. 
This is precisely the henselianity property for the spaces $\calT_b^{\gamma\in S}$. 
Indeed, the $\calT_b^{\gamma\in S}$ are henselian towards their unique non-analytic closed point $\calT_b^{\gamma=b}\subseteq \calT_b^{\gamma\in S}$.
 
\begin{remark}
Although, in essence, our proof that $\pitch$ is an equivalence follows the sketch provided above, we warn the reader that the actual details look quite different due to various technical subtleties.
In particular, one of the main reasons for this is that $\pitch$ is constructed by taking a colimit, and thus every appearance of $\ICG$, $\Bun_G^\mer$ or any of their corresponding locally closed strata that appear in the above argument must be replaced by a bounded piece of it.
For this reason, in the actual technical argument the stack of shtukas $\Sht^{\sch}_{\calI,\leq \mu}$ and various loci within it appears more often than $\ICG$.
\end{remark}
The proof that $\pitch$ is an equivalence is spread through \S \ref{sec: compute analytif}, \S \ref{sec: kimberlitecalculations} and \S \ref{sec: TheEquivalence}. 
In \S \ref{sec: compute analytif} we introduce techniques that we call pro-unipotent base change and overconvergent replacement with the intention of showing that the map $\gamma_S:\calT_b^{\gamma\in S}\to \ICG^\diamondsuit_S$ still satisfies a weak version of smooth base change.
We believe that these techniques will have further applicability beyond this work, since they should help develop an analytic version of placid geometry, which should show up in various contexts. 
In \S \ref{sec: kimberlitecalculations} we recall the theory of kimberlites developed in \cite{Gle24} and expanded in \cite{Gle23}.  
One of the achievements of this theory is to show that moduli spaces of analytic shtukas are locally spatial kimberlites, and that locally spatial kimberlites are automatically henselian towards their non-analytic locus. 
As it turns out, the spaces $\calT_b^{\gamma\in S}$ admit a sind-resolution by locally spatial kimberlites, this is the key input to show that the spaces $\calT_b^{\gamma\in S}$ are henselian towards their non-analytic locus (see \Cref{nearbycyclestheorem}).
Finally, \S \ref{sec: TheEquivalence} is the section which ties all of the previous sections together in order to show that $\pitch$ is an equivalence.
In \S \ref{sec:properties_and_applications}, we explain the applications to the splitting of the semi-orthogonal decomposition on $\calD_{\Lambda}^{\an}(\Bun_{G})$ alluded to in \S \ref{s: SoftIntroduction} and the compatibility of Zhu's equivalence with geometric Eisenstein functors, by combining with results of \cite{Zhu25} on the tame categorical Langlands correspondence. We also spell out the conjectural linearity property of our functor under the spectral action and some applications to the cohomology of local Shimura varieties alluded to above.

\section*{Acknowledgments}
I.~G.~was supported by the Deutsche Forschungsgemeinschaft (DFG, German Research Foundation) through the Leibniz Preis of Peter Scholze and by the National University of Singapore through the PYP funding scheme. 
A.~I.~was supported by the DFG through the Heisenberg program (grant nr. 462505253).
J.~L.~was supported by Germany’s Excellence Strategy EXC 2044/2 – 390685587, Mathematics Münster : Dynamics–Geometry–Structure, by the European Research Council (ERC) through the Consolidator Grant 770936 of Eva Viehmann, and by the DFG through the Leibniz Preis of Eva Viehmann, and the CRC 1442 Geometry: Deformations and Rigidity.
K.~Z.~was supported by the DFG through the Leibniz Preis of Peter Scholze, and by the Simons Collaboration on
Perfection in Algebra, Geometry, and Topology as postdoctoral researcher at CNRS. He also thanks the Max Planck Institute for Mathematics for their hospitality during part of this work.

We are very grateful to Xinwen Zhu for sharing an early draft of \cite{Zhu25}.
We thank the following mathematicians for conversations related to this work: 
Johannes Anschütz, Alexander Bertoloni-Meli, Raphäel Beuzart-Plessis, Alexis Bouthier, Ana Caraiani, Jean-François Dat, Gabriel Dospinescu, Arnaud \'Eteve, Laurent Fargues, Dennis Gaitsgory, Anton G\"uthge, David Hansen, Eugen Hellmann, David Helm, Pol van Hoften, Thibaud van den Hove, Dongryul Kim, Arthur-César Le Bras, Lucas Mann, Stefano Morra, Vincent Pilloni, Timo Richarz, Simon Riche, Marco Robalo, Juan Esteban Rodríguez Camargo, Nick Rozenblyum, Peter Scholze, Eva Viehmann, Sug Woo Shin, Zhiyou Wu, Yifei Zhao, Mingjia Zhang, and Felix Zillinger.

	\section*{Notation.}\label{sec:notation}
	\begin{enumerate}
		\item Let $E$ be a non-Archimedean local field of residue characteristic $p>0$, $O_E$ its ring of integers, and $\bbF_q$ its residue field of cardinality $q$.
		Let $\pi\in O_E$ be a uniformizer. 
		We also let $\breve{E}$ denote its completed maximal unramified extension and $\Gamma_E := \Gal(\ol{E}/E)$ be its absolute Galois group.
		\item Let $G$ be a connected reductive group defined over $E$. We often write $\calG$ for an auxiliary parahoric model of $G$ defined over $O_E$.
If $\calG$ is an Iwahori, we frequently use the letter $\calI$ instead.
		\item Fix $k$ an algebraically closed field of characteristic $p$ together with an embedding $\bbF_q\subseteq k$.
		\item Let $\CAlg$ and $\CAlg^{\perf}$ denote the category of $k$-algebras and perfect $k$-algebras respectively.
Let $\Sch$, $\PSch$, $\Sch^\qcqs$, $\PSch^\qcqs$ and $\PSchf$ denote the categories of schemes, perfect schemes, qcqs schemes, perfect qcqs schemes and perfect affine schemes over $\Spec k$ respectively. 
		\item For any $R\in \CAlg^{\perf}$, we let $\on{Frob}$ denote the arithmetic $q$-Frobenius of $R$.
		\item  
		Given $R\in \CAlg^{\perf}$, we let $\bbW(R)$ denote the ring of $O_E$-ramified Witt vectors.
		We let $\varphi:\bbW(R)\to \bbW(R)$ denote the automorphism induced from $\on{Frob}$.
		\item We let $\Perf$ and $\Perff$ denote the categories of perfectoid and affinoid perfectoid spaces over $k$. 
		\item  Given a Huber pair $(A,A^+)$ over $\bbZ_p$ we let $\Spd(A,A^+)$ denote its associated v-sheaf as in \cite[Lemma 15.1]{Sch17}. 
		Whenever $A^+=A^\circ$ we also denote $\Spd A =\Spd(A,A^\circ)$.
		In particular, if $R$ is a ring endowed with its discrete topology, then $\Spd R$ denotes $\Spd(R,R)$.
		\item In our higher categorical setup, all $2$-morphisms are invertible, so we work exclusively with $(\infty,1)$-categories (or simply $\infty$-categories).
We also follow the choice of terminology employed in \cite[Subsection 1.1, Page 5]{HeyerMann}: a $(\infty,1)$-category will be called a category from now on and a $1$-truncated category will be referred to as an ordinary category.
        \item Given a category $\calC$, we let $\calC^\omega$ denote the subcategory of compact objects. 
        \item If $\mathcal{C}$ is a category, and $E$ is a subset of arrows in the homotopy category that contains all isomorphisms, we let $\mathcal{C}_E\subseteq \mathcal{C}$ denote the wide subcategory spanned by $E$.
        \item Let $\Cat$ denote the category of small categories.
We let $\widehat{\Cat}$ denote the category of not necessarily small categories. 
			We let $\Pr^L$ denote the category of presentable categories whose morphisms are colimit-preserving functors. 
			We let $\Ani\in \Pr^L$ denote the category of anima or spaces.  
		We let $\Delta\in \Cat$ denote the simplex category.
We let $\LinCat$ denote the category of presentable stable categories whose morphisms are colimit preserving functors, and we write \(\LinCat^R\) for the category of presentable stable categories whose morphisms are right adjoint functors.
	\item Given a ring $\Lambda$, we let $\mathrm{Mod}_\Lambda \in \LinCat$ denote the derived category of $\Lambda$-modules which is an algebra object.
We also let $\LinCat_\Lambda$ denote the category of presentable stable $\Lambda$-linear categories (i.e., the category of $\Lambda$-module objects in $\LinCat$).
This has a symmetric monoidal structure $\otimes_\Lambda$ given by the relative Lurie tensor product over $\Modd_\Lambda \in \LinCat_\Lambda$, which is its $\otimes$-unit.
		We let $\LinCat_\Lambda^\cg\subseteq \LinCat_\Lambda$ denote the (non-full) subcategory of compactly generated $\Lambda$-linear categories whose functors preserve compact objects. 
			We also let $\LinCat_\Lambda^\sm$ denote the category of small stable $\Lambda$-linear categories. 
		\item We fix $\ell \neq p$. 
			All the coefficient rings $\Lambda$ we consider are classical and $\bbZ_\ell$-algebras.
For the majority of the text they will further be assumed to be $\ell$-torsion, i.e., $\bbZ/\ell^n\bbZ$-algebras for some $n$. 
		\item Given a category \(\calC\), we will need to consider \(\PreStk=\Fun(\calC^\op,\Ani)\).
This poses issues if \(\calC\) is not small.
In regards to set theory, we consider two possibilites.
The first one is to fix a sufficiently big regular cardinal \(\kappa\).
In this case, by ``small'' we mean \(\kappa\)-small and by ``accessible'' we mean the empty condition, and restrict the categories in Point (4) and (7) to \(\kappa\)-small versions.
There are no issues defining prestacks on these categories.
		The second possiblity is to follow the approach of \cite[\S 4]{Sch17} that only considers those prestacks that are accessible, equivalently small.
        By definition, this is equivalent to taking the union over all \(\kappa\) from the first approach.
We refer to \cite[Definition A.3]{hesselholt2024diracgeometryiicoherent} for more exposition on accessible prestacks.
By \cite[Remark A.4]{hesselholt2024diracgeometryiicoherent}, the Yoneda lemma and the universal characterization of the category of prestacks still works as expected.
        \item 	Given a topology $\tau$ on a category $\calC$ in the sense of \cite[Definition 6.2.2.1]{HTT}, we let $\Stk_\tau\subset \PreStk$ denote the full subcategory of accessible prestacks that satisfy descent with respect to $\tau$, i.e., they form sheaves in the sense of \cite[Definition 6.2.2.6]{HTT}.
Let us mention that there exists a stronger notion of hypercomplete sheaves, see \cite[Definition 6.5.3.2]{HTT}, that we do not make use of.
        \item The inclusion $\Stk_\tau\subset \PreStk$ has left adjoint, which we denote by $(-)_{\tau}$.
        If \(\calC\) is small, this is \cite[Lemma~A.4.15]{HeyerMann}. 
        If instead \(\calC\) is big but \(\calC^\op\) is accessible, we have fully faithful functors \(\PreStk^{\leq\kappa}\injto\PreStk^{\leq\kappa'}\) for \(\kappa\leq\kappa'\) given by left Kan extension, one checks that \((-)_\tau\) for different \(\kappa\) assemble together.
        Here the superscript means we take prestacks on \(\kappa\)-small objects.
        In any of these cases this is known as the sheafification functor, and exhibits $\Stk_{\tau}$ as a reflective subcategory of $\PreStk$. 
        \item If $\mathcal{C}_{1}$ and $\mathcal{C}_{2}$ are two compactly generated categories in $\LinCat^\cg$, we let $\mathcal{C}_{1}^{\omega}$ and $\mathcal{C}_{2}^{\omega}$ denote the full subcategories spanned by their compact objects. 
		Given $F: \mathcal{C}_{1} \ra \mathcal{C}_{2}$ a functor preserving compact objects, we write $F^{\omega}: \mathcal{C}_{1}^{\omega} \ra \mathcal{C}_{2}^{\omega}$ for the induced functor on compact objects. 
		Finally, we define $F^{o} := \Ind(F^{\omega,\op})$ regarded as a functor $\mathcal{C}_{1}^{\vee} \ra \mathcal{C}_{2}^{\vee}$ between the dual categories. 
		We refer to this as the \emph{conjugate functor}.\label{cons: ConjugateFunctor}
	\end{enumerate}

	\section*{\textbf{Part I. Formalism and generalities}}
    \section{Analytification and 6-functor formalisms}\label{sec: analytificiation and 6-functor formalisms}
    \noindent
    \textbf{Soft preamble to the section:}
        Throughout this paper, we are interested in comparing objects coming from algebraic and analytic sheaf theories. 
	We will need to understand how to pass from the algebraic context to the analytic context and how to relate the six operations that appear in each of the setups. 
	In order to do this, we will use the language of 6-functor formalisms and natural transformations between them. \\

	\noindent
    \textbf{Technical preamble to the section:}
    This section is mostly expository, but we discuss some subtle technical points that can easily be missed when thinking about analytification for the first time.
    \begin{enumerate}
	    \item[\S \ref{ss: basic setup for analytification}] Recalls the basic setup that we use to define and study analytification. 
	    \item[\S \ref{section on six functors}] Recalls the theory of 6-functor formalisms following \cite{Mann2022padic6functors,scholze6ff,HeyerMann}.
	    \item[\S \ref{section-natural-transforms}] Studies analytification as a formal construction. 
		    We explain the perspective that analytification functors should be treated as morphisms in their own 6-functor formalism. 
		    We define and study the so-called standard $!$-able maps (\Cref{define-standard-Shriek}) and standard fine maps (\Cref{lets define fine finally}), which are subsets of the usual $!$-able maps and usual fine maps that behave better with respect to analytification. 
		    The basic observation is that maps of universal $\calD^!$-descent do not behave appropriately with respect to analytification, and must be replaced by a stronger condition. 
    \end{enumerate}

    \subsection{Perfect and Perfectoid Prestacks}{\label{ss: basic setup for analytification}}
        \subsubsection{Perfect Prestacks}{\label{ss: perfectPrestacks}}
        We first setup some notation following \cite[\S 10]{Zhu25}.
	
	\begin{definition}{\label{defn: perfectprestack}}
		A \textit{perfect prestack} $X$ is an accessible functor $X:\PSch^{\aff,\op}\to \mathrm{Ani}$ from the category of perfect affine schemes over \(k\) to the category of anima.
We write $\PreStk$ for the category of perfect prestacks. 
	\end{definition}

	In the schematic setup, all of our geometric objects belong to $\PreStk$ after a suitable embedding.
As we will consider multiple different $6$-functor formalisms on schemes satisfying slightly different descent properties, it is often convenient to formulate things in this level of generality.
    
	We consider several full subcategories of $\PreStk$.
	The Yoneda embedding described in \cite[Proposition 5.1.3.1]{HTT} yields fully faithful functors 
	\[\CAlg^{\perf,\op} \simeq \PSchf \hookrightarrow \PSch\hookrightarrow \PreStk,\] and we invariably identify these categories with their essential image. 
	We let $\CAlg^{\pfp}\subset \CAlg^{\perf}$ denote the full subcategory of $k$-algebras that are perfectly finitely presented over $k$, compare with \cite[Definition 10.4]{Zhu25}. 
	We let $\PSch^\pfp\subset \PSch^\qcqs\subset  \PSch$ denote the subcategory of perfectly finitely presented $k$-schemes and qcqs perfect $k$-schemes, respectively. 
	We let $\AlgSp$ denote the category of perfect algebraic spaces over $k$ and $\AlgSp^\pfp\subset \AlgSp$ the full subcategory of those of perfect finite presentation, see \cite[Definition 10.4]{Zhu25}.
	
	Given a topology $\tau$ on $\CAlg$ in the sense of \cite[Definition 6.2.2.1]{HTT}, we let $\SchStk_\tau\subset \PreStk$ denote the full subcategory of accessible prestacks that satisfy descent with respect to $\tau$, i.e., they form sheaves in the sense of \cite[Definition 6.2.2.6]{HTT}. 
	Let us mention that there exists a stronger notion of hypercomplete sheaves, see \cite[\S 6.5.2, \S 6.5.3]{HTT}, that we do not make use of.
	Throughout, we will mostly let $\tau$ be one of the following.
	\begin{enumerate}
		\item The \'etale topology, with corresponding sheaf category $\SchStk_\et$.
		\item The schematic v-topology as in \cite[Definitions 2.1, 2.13]{bhatt_scholze_projectivity_of_the_witt_vector_affine_grassmannian}, with corresponding sheaf category $\SchStk_v$.
		\item The arc-topology as in \cite[Definition 1.2]{BhattMathew_21}, with corresponding sheaf category $\SchStk_{\on{arc}}$.
		\item The $\diamondsuit$-topology, which will be defined in \Cref{diamond topology}, with corresponding sheaf category $\SchStk_{\diamondsuit}$. 
	\end{enumerate}
        For such a Grothendieck topology $\tau$, we recall that the inclusion $\SchStk_{\tau} \ra \PreStk$ admits a left adjoint by \cite[Lemma~A.4.15]{HeyerMann}, which we denote by $(-)_{\tau}$. 
        This is known as the sheafification functor, and exhibits $\SchStk_{\tau}$ as a reflective subcategory of $\PreStk$.
        
    \subsubsection{Perfectoid Prestacks}
    \label{more subsections here}
     In the analytic context, we have the following analogue to \Cref{defn: perfectprestack}.
	\begin{definition}
		\label{accessible condition}
		A \textit{perfectoid prestack} $X$ is an accessible functor 
		\[X:\Perf^{\aff,\op} \to \Ani,\]
        from the category of affinoid perfectoid spaces over $k$ to the category of anima.
		We write $\AnPreStk$ for the category of perfectoid prestacks. 
	\end{definition}

	\begin{remark}
		The condition that a v-sheaf is small in \cite[Definition 12.1]{Sch17} agrees with the accessibility condition in \Cref{accessible condition}. 
	\end{remark}
	Just like in the schematic setup, all geometric objects take values in $\AnPreStk$. 
	While the category $\Perff$ does carry a wide variety of topologies, we will mostly consider the v-topology, as defined in \cite[Definition~8.1 (iii)]{Sch17}. 
	Occasionally, we also consider the pro\'etale topology \cite[Definition~8.1 (i)] {Sch17}.
	As in the previous section, we obtain full subcategories
        \begin{equation}{\label{eqn: SubcategoriesofAnalyticStacks}}
         \AnStk_{v} \subseteq \AnStk_\uop \subseteq  \AnPreStk, 
         \end{equation}
	and the inclusion functors admit left adjoint functors 
	\[(-)_{v}: \AnPreStk \ra \AnStk_{v} \text{ and } (-)_\uop: \AnPreStk \to \AnStk_\uop\] given by sheafification. 
	We have full subcategories $\vShv \subset \vStk \subset \AnStk_{v}$, of small v-sheaves and small v-stacks as defined in \cite[Definition~12.1]{Sch17} and \cite[Definition~11.4]{Sch17}, respectively. 

	\begin{remark}{\label{rem: smallvstacks}}
	A small v-stack $X \in \vStk$ is the same as a $1$-truncated (i.e., take values in the subcategory $\tau_{\leq 1}: \Ani_{\leq 1} \hookrightarrow \Ani$ of spaces whose homotopy groups vanish for $i > 1$) perfectoid prestack that lies in the full subcategory $\AnStk_{v}$. 
	Although the 6-functor formalism of \'etale cohomology of diamonds \cite{Sch17} was initially only constructed for $\vStk$, with modern tools (see, e.g., \cite{Mann2022NuclearSheaves}), one can adapt the construction to treat the category $\AnStk_{v}$ of v-stacks with values in anima.
	For our purposes, the difference between $\AnStk_v$ and $\vStk$ is immaterial, since all of our geometric objects take values in $\vStk$. 
	\end{remark}
        
	We recall that $\vShv$ contains the full subcategory of diamonds (see \cite[Proposition 11.9]{Sch17}), which are defined as a special kind of set-valued pro\'etale sheaves on $\Perf$ that resemble algebraic spaces (see \cite[Definition~11.1]{Sch17}). 
	Moreover, \cite[Definition~11.17]{Sch17} isolates the better-behaved subcategories of spatial (resp. locally spatial) diamonds. 
	These particular type of diamonds are very useful for cohomological considerations. 
	Indeed, one of the key features is that, for these diamonds, their underlying topological spaces are automatically spectral (resp. locally spectral) (see \cite[Proposition~11.19]{Sch17}).

        \subsubsection{The $\diamondsuit$-Functor}{\label{sec: the Diamond Functor}}
	We follow \cite[\S 27]{Sch17}.
        \begin{definition}{\label{defn: diamondfunctor}}
        We define
        \[ (-)^{\diamondsuit_{\pre}}: \PreStk \ra \AnPreStk  \]
        \[ X \mapsto \{(R,R^{+}) \mapsto X(\Spec(R)) \}, \]
	we write $(-)^{\diamondsuit}:=(-)_v\circ (-)^{\diamondsuit_{\pre}}$ for the composition of $(-)^{\diamondsuit_{\pre}}$ with the v-sheafification functor and we refer to either of these functors as the \emph{big diamond functor}. 
        \end{definition} 
        More conceptually, note that $\Spec(A)^{\diamondsuit} = \Spec(A)^{\diamondsuit_{\pre}}$ is a v-sheaf in light of \cite[Theorem~8.7]{Sch17}.
 Therefore, the functor $\diamondsuit$ is a left Kan-extension of the functor 
        \[ \diamondsuit: \PSchf \ra \AnStk_v \subset \AnPreStk \]
        \[ \Spec(A) \mapsto \Spec(A)^{\diamondsuit}, \]
        along $\PSchf \subset \PreStk$. 

	In \S \ref{subsec: threeanalytificationfunctors}, we will describe two additional constructions related to analytification $\diamond$ and $\dagger$. 
	However, among these three functors ($\diamond$, $\dagger$ and $\diamondsuit$), $\diamondsuit$ is the only functor that upgrades to a natural transformation of 6-functor formalisms as we will explain below. 
	This is why we study $\diamondsuit$ first. 
	We recall some generalities on $6$-functor formalisms, following closely \cite{Mann2022padic6functors,scholze6ff,HeyerMann} with special emphasis on the latter reference.

	\subsection{6-Functor Formalisms}
	\label{section on six functors}
        \subsubsection{Basic Notions}
	The basic context on which a 6-functor formalism exists is the following, compare also with \cite[Definition 2.1.1, Remark 2.1.2]{HeyerMann}.
    
	\begin{definition}{\label{defn: geometricsetup}}
		A \emph{geometric setup} is a pair \((\mathcal{C},E)\), where \(\mathcal{C}\) is an \(\infty\)-category admitting finite limits\footnote{This assumption is not made in \cite[Definition 2.1.1]{HeyerMann}, but it simplifies the discussion.} and \(E\) is a homotopy class of edges such that
		\begin{enumerate}
			\item \(E\) contains all isomorphisms,
			\item $E$ is stable under composition,
			\item $E$ is stable under base change in $\calC$,
			\item \(E\) is stable under forming diagonals, i.e., if \(f\in E\), then \(\Delta_f\in E\).
Equivalently, $E$ is right-cancellative i.e., if $g\circ f\in E$ and $g\in E$ then $f\in E$ (see \cite[Lemma 2.1.5]{HeyerMann}).
		\end{enumerate}
	\end{definition}
	The morphisms in $E$ are colloquially referred to as the $!$-able maps of the geometric setup.
		Geometric setups form a category, see \cite[Definition 2.1.1]{HeyerMann}, whose morphisms $(\mathcal{C}_{1},E_{1}) \ra (\mathcal{C}_{2},E_{2})$ are functors $f:\mathcal{C}_{1} \ra \mathcal{C}_{2}$ that restricts to a functor $\calC_{E_{1}}\to \calC_{E_{2}}$ satisfying that $f(X\times_Y Z)\simeq f(X)\times_{f(Y)}f(Z)$ whenever $[X \ra Y]\in \calC_{E_1}$.
    
	Given a geometric setup, one can construct a category $\corr(\mathcal{C},E)$ of correspondences, see \cite[Definition 2.2.10]{HeyerMann}.
	Objects in this category are simply the objects of $\calC$.
	Morphisms can be informally described as diagrams of the form $\{X \leftarrow Z \rightarrow Y\}$ such that $Z \rightarrow Y$ is in $E$. 
	Composition of correspondences is given by the obvious diagram involving fiber products, see also \cite[Definition 1.2.2]{HeyerMann}. 
	The mapping $(\mathcal{C},E) \mapsto \corr(\mathcal{C},E)$ is functorial with respect to morphisms of geometric setups and preserves limits, see \cite[Lemma~2.2.11]{HeyerMann}.

	By \cite[Proposition~2.3.3, Example 2.3.4]{HeyerMann}, the category $\Corr(\mathcal{C},E)$ carries a symmetric monoidal structure $\Corr(\mathcal{C},E)^{\otimes}$, which on objects is simply the Cartesian product on $\calC$ (recalling that we always assume that $\calC$ has finite limits). 
	We will reference this symmetric monoidal structure throughout. 
	Moreover, given any map of geometric setups $(\calC,E) \ra (\calC',E')$ the induced map on correspondence categories $\Corr(\calC,E)^{\otimes} \ra \Corr(\calC',E')^{\otimes}$ is a lax symmetric monoidal functor, by \cite[Proposition~2.3.5]{HeyerMann}.
	\begin{definition}[{\cite[Definition~1.2.3]{HeyerMann}}]{\label{defn: 3functorformalism}}
		Consider $\widehat{\Cat}$ equipped with the Cartesian symmetric monoidal structure $\widehat{\Cat}^\times$ in the sense of \cite[Definition 2.4.1.1, Proposition 2.4.1.5]{LurieHigherAlgebra}. 
		A \emph{3-functor formalism} on $(\calC,E)$ is a lax symmetric monoidal functor 
		\[ \mathcal{D} : \corr(\calC,E)^\otimes \ra \widehat{\Cat}^\times  \]
		\[ X \mapsto \mathcal{D}(X). \]
	\end{definition}
	\begin{remark}[{\cite[Definition~3.1.4]{HeyerMann}}]
		\label{basic-obs}
		We have a non-full embedding $\mathcal{C}_{E} \hookrightarrow \corr(\mathcal{C},E)$ which on objects is given by the mapping $[X\mapsto X]$
		and on arrows is given by the mapping
		\[ [f: X \ra Y] \mapsto \{X \xleftarrow{\id_{X}} X \xrightarrow{f} Y \}.\]
		Given a $3$-functor formalism, we write $f_{!}: \mathcal{D}(X) \ra \mathcal{D}(Y)$ for the induced functor. 
		Similarly, we have a non-full embedding $\mathcal{C}^{\op} \rightarrow \corr(\mathcal{C},E)$ which on objects is given by the mapping $[X\mapsto X]$ and on arrows can be described by the mapping
		\[ [f: X \ra Y] \mapsto \{Y \xleftarrow{f} X \xrightarrow{\id_{X}} X \}.\]
		Given a $3$-functor formalism, we write $f^{*}: \mathcal{D}(Y) \ra \mathcal{D}(X)$ for the induced functor. 
Finally, we note that the lax symmetric monoidal structure defines a map 
		\[ \calD(X) \times \calD(X) \ra \calD(X \times X), \]
		which, when composed with the pullback $\Delta^{*}: \calD(X \times X) \ra \calD(X)$ along the diagonal map, defines a symmetric monoidal structure $\otimes_{\calD,X}$ on $\calD(X)$. 
	\end{remark}
	\begin{remark}
		The above shows how to extract the three operations $f^*$, $f_!$, and $\otimes_{\calD}$ from an abstract 3-functor formalism.	
		The point of encoding this in the correspondence category is that it captures coherently all the compatibilities that one expects (see \cite[Proposition~3.1.8]{HeyerMann}).
		For example, it follows formally that on any $3$-functor formalism the functors $f^*$ are symmetric monoidal with respect to $\otimes_\calD$ (see \cite[Proposition~3.1.8.(ii)]{HeyerMann}).
		We will often use $\mathbf{1}_{X}$ to denote the tensor unit in this symmetric monoidal structure. 
		If the context is clear, we remove the decorations and simply write $\mathbf{1}$ and $\otimes$ instead.
	\end{remark}

	Recall that $\PrL$ is equipped with the so-called \emph{Lurie tensor product} denoted $\otimes$, as defined in \cite[Proposition 4.8.1.15]{LurieHigherAlgebra}.
Moreover, the forgetful functor $\PrL\to \widehat{\Cat}$ is lax symmetric monoidal.
		\begin{definition}[{\cite[Definition 3.2.1]{HeyerMann}}]
			{\label{defn: presentableSixFunctorFormalisms}} 
			Let $(\calC,E)$ be a geometric setup.
		\begin{enumerate} 
			\item A \emph{6-functor formalism} on $(\calC,E)$ is a 3-functor formalism on $(\calC,E)$ such that all the functors $f_{!} \colon \calD(X) \to \calD(Y)$, $f^{*} \colon \calD(Y) \to \calD(X)$, and $\otimes$ admit right adjoint functors.
These are denoted $f^{!}$, $f_{*}$, and $\uHom(-.-)$ respectively.
			\item A \emph{presentable 3-functor formalism} $\calD$ is a lax symmetric monoidal functor
				\[\calD:\Corr(\calC,E)^{\otimes}\to (\PrL)^\otimes.\]
				Equivalently, it is a 3-functor formalism which factors through the subcategory $\PrL$ of $\widehat{\Cat}$ and for which the map $\calD(X)\times \calD(Y)\to \calD(X\times Y)$ factors through $\calD(X)\otimes \calD(Y)$ (see \cite[Lemma 3.2.5]{HeyerMann}).
		\end{enumerate}
	\end{definition}
	\begin{remark}{\label{rem: presentable3functorsare6functors}}
		By the $\infty$-categorical adjoint functor theorem \cite[Corollary~5.5.2.9]{HTT}, any presentable $3$-functor formalism is automatically a $6$-functor formalism. 
	\end{remark}

\begin{definition}
			For a $6$-functor formalism $\calD$ and $[f:X \ra S] \in \calC_{E}$, we define the \emph{Verdier duality functor} 
			\[ \bbD_{X/S,\calD}: \calD(X) \ra \calD(X)^{\op} \]
			\[ A \mapsto \uHom(A,f^{!}(\mathbf{1})). \]
			If $\ast\in \calC$ denotes the final object, and $X \ra \ast$ is in $\calC_{E}$, we will denote this functor by $\bbD_{X,\calD}$.
We will suppress $\calD$ when it is clear from context.
\end{definition}
\begin{remark}{\label{rem: VerdierdualitySends!to*s}}
We warn the reader that although this functor is called Verdier duality, it is not a duality in many of the situations that we consider.
Indeed, in the analytic theory this functor is almost never a duality, it already fails for the unit ball.
\end{remark}

\begin{remark}
	It follows from the Yoneda lemma and the projection formula that, for a diagram in $\calC_E$,
\[ \begin{tikzcd} 
			X \arrow{rr}{f} \arrow[dr] & & Y \arrow[dl] \\
			& S &
		\end{tikzcd}, \]
		the following identities hold
\[ \bbD_{Y/S,\calD}f_{!} \simeq f_{*}\bbD_{X/S,\calD}  \]
\[ \bbD_{X/S,\calD}f^{*} \simeq f^{!}\bbD_{Y/S,\calD}.  \]
However, the other familiar identities can fail whenever $\bbD_{X/S}$ or $\bbD_{Y/S}$ are not equivalences.
Whenever $\bbD_{X/S}^{2} = \id$ and $\bb{D}_{Y/S}^{2} = \id$ hold, then we also have natural equivalences 
\[ \bbD_{Y/S,\calD}f_{*} \simeq f_{!}\bbD_{X/S,\calD}  \]
\[ \bbD_{X/S,\calD}f^{!} \simeq f^{*}\bbD_{Y/S,\calD}.  \]
\end{remark}

We will repeatedly use \Cref{prop: BiDuality} below. 
Recall that an object $A\in \calD(Y)$ is called reflexive if the natural map $A\to \bb{D}_{Y/S}^{2}A$
is an equivalence.
\begin{proposition}{\label{prop: BiDuality}}
Fix a presentable 3-functor formalism $\calD$, and let $A\in \calD(Y)$. 
The following statements hold.
\begin{enumerate}
	\item If both $A$ and $f^*\bbD_{Y/S}A$ are reflexive, then $\bbD_{X/S,\calD}f^{!} A \simeq f^{*}\bbD_{Y/S,\calD} A$.	
	\item If $A$ and $f_!\bbD_{X/S,\calD}A$ are reflexive, then $\bbD_{Y/S,\calD}f_{*} A \simeq f_{!}\bbD_{X/S,\calD} A$.	
\end{enumerate}
\end{proposition}

\subsubsection{More elaborate notions}
As \Cref{prop: BiDuality} exemplifies, it is important to isolate a set of objects on which relative Verdier duality is an equivalence. 
A key example is the set of the so-called \emph{suave objects} which we recall below. 

Let $(\mathcal{C},E)$ be a geometric setup, $\calD$ a 3-functor formalism over $(\calC,E)$ and fix an object $S \in \mathcal{C}$. 
To these data, we can attach the $2$-category of kernels, see \cite[Definition~4.1.3]{HeyerMann}, that is denoted $\calK_{\calD,S}$. 
	Objects are given by $[X \ra S ]\in \calC_{E}$ and morphism categories are given by $\calD(X \times_{S} Y)$. %
	We have the following abstract notions defined in terms of the kernel category.
\begin{definition}[{\cite[Definition~4.4.1, Definition~4.5.1]{HeyerMann}}]
	{\label{defn: Suavesmoothunipotent}}
		Fix a 3-functor formalism $\calD$ on a geometric setup $(\calC,E)$, and a map $[f: X \ra S ]\in \mathcal{C}_{E}$.
		\begin{enumerate} 
			\item We say that $A \in \calD(X) = \Hom_{\mathcal{K}_{S,\calD}}(X,S)$ is $f$-\emph{suave} if it admits a right adjoint $\mathrm{SD}_{f,\calD}(A) \in \calD(X) = \Hom_{\mathcal{K}_{S,\calD}}(S,X)$, which we refer to as the $f$-\emph{suave dual} of $A$. 
			\item We say that $A \in \calD(X) = \Hom_{\mathcal{K}_{S,\calD}}(X,S)$ is $f$-\emph{prim} if it admits a left adjoint $\mathrm{PD}_{f,\calD}(A) \in \calD(X) = \Hom_{\mathcal{K}_{S,\calD}}(S,X)$, which we refer to as the $f$-\emph{prim dual} of $A$. 
			\item We let $\mathrm{Suave}_{f,\calD}(X) \subset \calD(X)$ denote the full subcategory of suave objects inside $\calD(X)$. 
			\item We let $\mathrm{Prim}_{f,\calD}(X) \subset \calD(X)$ denote the full subcategory of prim objects inside $\calD(X)$. 
			\item We say that $f$ is $\calD$\emph{-suave} if the tensor unit $\mathbf{1} \in \calD(X)$ is suave. 
			\item We say that $f$ is $\calD$\emph{-prim} if the tensor unit $\mathbf{1} \in \calD(X)$ is prim. 
			\item We say that $f$ is $\calD$\emph{-smooth} if it is $\calD$-suave and $\mathrm{SD}_{f,\calD}(\mathbf{1})$ is invertible for $\otimes_\calD$. 
            \item We say that $f$ is $\calD$\emph{-unipotent} if it is $\calD$-smooth and the pullback functor $f^{*}$ is fully faithful.
		\end{enumerate}
	\end{definition}

	The category of suave objects is equipped with an involutive equivalence, which is a manifestation of the fact that that $\mathcal{K}_{S,\calD}^{\op} \simeq \mathcal{K}_{S,\calD}$, since $\calD(X \times_{S} Y) \simeq \calD(Y \times_{S} X)$.
	\begin{proposition}[{\cite[Lemma~4.4.4 (i)]{HeyerMann}}]{\label{cor: suavedualityinducesequivalence}}
		The suave duality functor induces an involutive equivalence
		\[ \Suave_{f,\calD}(X) \xrightarrow{\simeq} \Suave_{f,\calD}(X)^{\op}. \]
	\end{proposition}
	We note that suave duality has the following concrete interpretation in a 6-functor formalism.
	\begin{proposition}[{\cite[Lemma~4.4.5]{HeyerMann}}]{\label{prop: VerdierDualityandSuaveDuality}}
		For a 6-functor formalism $\calD$ and a map $[f: X \ra S] \in \calC_{E}$, we have a natural identification
		\[ \bbD_{X/S,\calD} \simeq \SD_{f,\calD} \]
		when we restrict the left-hand side functor to $\Suave_{f,\calD}(X)$.
	\end{proposition}
		In light of Proposition \ref{prop: VerdierDualityandSuaveDuality}, when speaking about suave duality in the setting of $6$-functor formalisms we will just refer to it as Verdier duality, and implicitly use this identification throughout.
Note that, in particular, Verdier duality becomes involutive when restricted to suave objects. 
        
	The next group of definitions concerns descent in a 3-functor formalism.
	We will take the approach of \cite{HeyerMann} and use covering sieves. 
	We refer the reader to \cite[Definition 6.2.2.1]{HTT} and \cite[Definition A.4.1]{HeyerMann} for the definition of a sieve.

	\begin{definition}[{\cite[Definition~3.4.6]{HeyerMann}}]
		{\label{defn: universal*and!covers}}
		Let $(\calC,E)$ be a geometric setup and let $\calD: \Corr(\calC,E) \ra \widehat{\Cat}$ be a 3-functor formalism.
Let $X \in \calC$.
		\begin{enumerate}
			\item We say a sieve $\mathcal{U}$ over $X$ is a $\calD^{*}$\emph{-cover} of $X$ if the natural map
				\[ \calD(X) \ra \lim_{V \in \mathcal{U}^{\op}} \calD(V) \]
			induced by $*$-pullback is an equivalence. 
			We say that it is a \emph{universal $\calD^{*}$-cover} if this is true after pulling back along any $X' \ra X$ with respect to the pullback sieve $f^{*}\calU$. 
		\item Given a $\calC_E$-sieve $\calU$ over $X$ (i.e., $\calU \subset (\calC_{E})_{/X}$), we say that $\calU$ is a $\calD^{!}$\emph{-cover} of $X$ if it is generated by a small family of maps $\{U_i\to X\}_{i\in I}$, and the natural map
				\[ \calD(X) \ra \lim_{V \in \mathcal{U}^{\op}}  \calD(V) \]
			induced by $!$-pullback is an equivalence. 
			We say that this is a \emph{universal $\calD^{!}$-cover} if this is true after pulling back along any $X' \ra X$ with respect to the pullback sieve $f^{*}\calU$.
		\item We say a map $f: U \ra X$ is a (universal) $\calD^*$-cover (resp. (universal) $\calD^!$-cover), if the sieve that it generates is a (universal) $\calD^*$-cover (resp. (universal) $\calD^!$-cover). 
		\end{enumerate}
	\end{definition}

	\begin{remark}
		In practice, most sieves that one cares about are generated by small family of maps, and very often generated by single map $U\to X$.
		In this case, the covering properties state that the map
				\[ \calD(X) \ra \lim_{n \in \Delta}  \calD(U^{n/X}),\]
				where $U^{n/X}$ denotes the $n$-fold fiber product of $U$ over $X$, is an equivalence (see \cite[Lemma A.4.6]{HeyerMann}).
	\end{remark}

 The following statements provide the main two classes of morphisms that are universal $\calD^*$-covers and universal $\calD^!$-covers.
    
    \begin{proposition}[{\cite[Lemma 4.7.1]{HeyerMann}}]
		\label{suave-implies-fine}
		Let $\calD$ be a 6-functor formalism on some geometric setup $(\calC,E)$.
		Let $[f\colon X\to Y]\in \calC_E$ and suppose that $\calD(X)$ has all countable limits and colimits.
		If $f$ is $\calD$-suave and $f^*:\calD(Y)\to \calD(X)$ is conservative, then $f$ is a universal $\calD^*$-cover and a universal $\calD^!$-cover.
	\end{proposition}

	There is a similar, yet more involved statement for prim maps.

	\begin{definition}[{\cite[Definition 3.18]{MatthewGaloisGroup}, \cite[Definition 4.7.2]{HeyerMann}}]
    		\label{matthews-defi}
		Let $\calC$ be a stable monoidal category. 
		For every object $A \in \calC$, we denote by $\langle A \rangle \subseteq  \calC$ the smallest full subcategory containing $A$ which is stable under finite limits, finite colimits, retracts and tensor products. 
		We say that $A$ is descendable if $\langle A \rangle$ contains the tensor unit.
	\end{definition}

	Recall that a $6$-functor formalism is said to be \textit{stable} if for all $X\in \calC$, the category $\calD(X)$ is a stable category.

	\begin{proposition}[{\cite[Lemma 4.7.4]{HeyerMann}}]\label{prim standard !-able}
		Let $(\calC,E)$ be a geometric setup.
Let $\calD$ be a stable $6$-functor formalism on $(\calC,E)$.
Let $[f:X\to Y]\in \calC_E$ be a map.
Suppose that $f$ is $\calD$-prim and that $f_* \mathbf{1}$ is descendable.
Then, $f$ is a universal $\calD^*$-cover and a universal $\calD^!$-cover.
	\end{proposition}

    \subsubsection{The Liu--Zheng construction}
    We now turn our attention to the construction of 3-functor formalisms. 
    There are two main types of tools available for this purpose. 
    One is the so-called Liu--Zheng construction, which takes two functors for granted ($f^*$ and $\otimes$) and produces the third ($f_!$) by passing to a left or a right adjoint of $f^*$. 
    The second type of tools are the extension procedures which start with a given 3-functor formalism and enlarge it.

    Let us give an informal account of the Liu--Zheng construction.
    In nice geometric situations, a morphism $f$ will admit a presentation as a composition $f=g\circ j$ where $j$ is an open immersion and where $g$ is a proper map.
    To define compactly supported cohomology, one would let $f_!:=g_*\circ j_\natural$ (extension by $0$ composed with usual cohomology). 
    For example, if we took $\calC = \mathrm{Sch}^{\mathrm{qcqs}}$ to be the category of qcqs schemes, $E$ the set of separated finite type morphisms, $I$ to be the set of open immersions and $P$ to be the set of proper morphisms, then, by the Nagata compactification theorem \cite[Tags 0F41, 0ATT]{StaProj}, every morphism $f\in E$ would decompose as $f=g\circ j$ with $g\in P$ and $j\in I$.

The Liu--Zheng construction takes as data subsets of morphisms $I\subseteq E$ that behave as if they were open immersions and $P\subseteq E$ that behave as if they were proper morphisms such that all the morphisms in $E$ can be obtained as compositions $f=g\circ j$ with $g\in P$ and $j\in I$.
    Given $I$ and $P$ subject to such axioms, one can define the third operation by letting $g_!:= g_*$ (the right adjoint to $g^*$) when $g\in P$ and $j_!:=j_\natural$ when $j\in I$ (the left adjoint to $j^*$).
    For general $f$, to define the third operation, one would find a decomposition $f=g\circ j$ and let $f_!:=g_!\circ j_!$.  
    This leads to the evident question: \emph{How does this depend on the presentation of $f$?}
    Intuitively speaking, the Liu--Zheng construction provides coherently organized isomorphisms between all the possible definitions of $f_!$ obtained from the various decompositions.

    Let us move to a more precise account on the subject.
    Recall the definition of a suitable decomposition.
    \begin{definition}[{\cite[Definition 3.3.2]{HeyerMann}}]{\label{defn: SuitableDecomposition}}
	    Let $(\calC,E)$ be a geometric setup. 
	    A \textit{suitable decomposition} of $E$ is a pair $I,P\subseteq E$ satisfying that:
	    \begin{enumerate}
    \item 		The classes \(I\) and \(P\) are stable under pullback, composition, and contain all isomorphisms.
    \item 		Any morphism \(f\in E\) factors as \(f=\overline{f}j\) for some \(j\in I\) and some \(\overline{f}\in P\).
    \item 		If \(f\from X\to Y\) and \(f'\from X'\to Y\) are both in \(I\) (respectively \(P\)) and \(g\from X'\to X\) is a map over \(Y\), then \(g\in I\) (respectively \(g\in P\)).
    \item 		If \(f\in I\cap P\), then \(f\) is \(n\)-truncated (\cite[Definition~A.1.20]{HeyerMann}) for some \(n \geq -2\).
	    \end{enumerate}	
    \end{definition}

    The key construction in this situation, introduced in \cite{liu2024enhancedoperationsbasechange} and streamlined in \cite{Mann2022padic6functors}, is as follows.
	\begin{construction}[{\cite[Proposition~3.3.3]{HeyerMann}}]
    \label{construction: 3 functors}
		Fix a geometric setup \((\calC,E)\) and a functor \(\calD_0\from\calC^\op\to\CAlg(\widehat{\Cat})\).
		Assume that there is a suitable decomposition \(I\) and \(P\) such that the following holds.
		\begin{enumerate}
			\item For any \(f\in I\), the functor \(f^*\) has a left adjoint \(f_\natural\) satisfying base change against any map in $\calC$ and satisfying the projection formula.
			\item For any \(f\in P\), the functor \(f^*\) has a right adjoint \(f_*\) satisfying base change against any map in $\calC$ and satisfying the projection formula.
			\item For any Cartesian diagram 
			\begin{equation*}
				\begin{tikzcd}
					X' \arrow[d, "g'"] \arrow[r, "j'"] & X \arrow[d, "g"] \\
					Y' \arrow[r, "j"]                  & Y               
				\end{tikzcd}
			\end{equation*}
			with \(j\in I\) (hence \(j'\in I\)) and \(g\in P\) (hence \(g'\in P\)), the natural map \(j_\natural g_*'\to g'_*j_\natural\) is an isomorphism.
		\end{enumerate}
		Then we obtain a 3-functor formalism \(\calD\from\corr(\calC,E)\to\widehat{\Cat}\) that extends $\mathcal{D}_{0}$ along the (non-full) embedding $\mathcal{C}^{\op} \hookrightarrow \corr(\calC,E)$. 
		Moreover, $\calD$ satisfies that for \(f\in I\), \(f_!=f_\natural\) is the left adjoint of \(f^*\) and for \(f\in P\), \(f_!=f_*\) is the right adjoint of \(f^*\). 
	\end{construction}

	\Cref{construction: 3 functors} is called the Liu--Zheng construction, and we use $\calL\calZ(\calC,E,P,I,\calD_{0})$ to denote the 3-functor formalism that it produces.
 
	It will be useful for us to have a characterization of the 6-functor formalisms constructed in this way. 
	That such a characterization exists was conjectured by Scholze and recently shown by Dauser and Kuijper \cite{dauser2024uniquenesssixfunctorformalisms}.
	To explain this, we introduce the following definition, which gives us an intrinsic characterization of the morphisms $P$ and $I$ appearing in the construction.

        \begin{definition}[{\cite[Definition~6.10, Definition 6.12]{scholze6ff}}]
		{\label{defn: cohomologicallyproperandetalemaps}}
		Given a 6-functor formalism $\calD$ on a geometric setup $(\calC,E)$, assume that $[f: X \ra Y ]\in \calC_{E}$ is $n$-truncated for some $n \geq -2$. 
		We consider the following conditions.
        \begin{enumerate}
        \item 
		We say that $f$ is $\calD$\emph{-proper} if:
		\begin{enumerate}[a)]
			\item $f$ is an isomorphism, or
			\item if $\Delta_{f}$ is $\calD$-proper and the natural map (see \cite[Proposition~6.11]{scholze6ff})
        \[ f_{!} \ra f_{*} \]
	is an isomorphism.
		\end{enumerate}

        \item 
		We say that $f$ is $\calD$\emph{-\'etale} if:
		\begin{enumerate}[a)]
			\item $f$ is an isomorphism, or
			\item if $\Delta_{f}$ is $\calD$-\'etale and the natural map (see \cite[Proposition~6.13]{scholze6ff})
        \[ f^{!} \ra f^{*} \]
	is an isomorphism.
		\end{enumerate}
        \end{enumerate}
        \end{definition}

\begin{remark}
        We note that \Cref{defn: cohomologicallyproperandetalemaps} is not self-referential. 
	Indeed, if $f$ is $n$-truncated then $\Delta_{f}$ is $n - 1$-truncated and $(-2)$-truncated maps are isomorphisms.
\end{remark}

        \begin{remark}{\label{rem: cohomologicallyetalemapsaresuave}}
		We note that $f: X \ra Y$ being $\calD$-\'etale (resp. $\calD$-proper) in the above sense immediately implies that $f$ is $\calD$-suave (resp. $\calD$-prim). 
		Indeed, this is precisely the content of \cite[Proposition~6.11, Proposition~6.13]{scholze6ff}.
        \end{remark}

	With \Cref{defn: cohomologicallyproperandetalemaps} in hand, we can recall the following. 

        \begin{definition}[{\cite[Definition~2.1, Definition 2.15]{dauser2024uniquenesssixfunctorformalisms}}]{\label{defn: nagatasetups}}\leavevmode
	\begin{enumerate}
        \item  A \emph{Nagata setup} is a tuple $(\calC,E,P,I)$ where $(\calC,E)$ is a geometric setup, $I,P\subseteq E$ is a suitable decomposition of $E$ (as in \Cref{defn: SuitableDecomposition}) subject to the stronger axiom that if $f\in I\cup P$ then $f$ is $n$-truncated for some  \(n \geq -2\).
        \item For $(\calC,E,P,I)$ a Nagata setup and $\calD$ a 6-functor formalism on the geometric setup $(\calC,E)$, we say that $\calD$ is \emph{Nagata} with respect to this Nagata setup if morphisms in $P$ (resp. $I$) are $\calD$-proper (resp. $\calD$-\'etale).
        \end{enumerate}
        \end{definition}
        The main result characterizing the essential image of the Liu--Zheng construction is the following.

        \begin{theorem}[{\cite[Theorem 3.3]{dauser2024uniquenesssixfunctorformalisms}}]{\label{thm: LiuZhengUniquenesss}}
		Let $(\calC,E,P,I)$ be a Nagata setup. 
		Suppose that $\calD$ is a 6-functor formalism that is Nagata with respect to $(\calC,E,P,I)$.
		Then we have a unique (up to contractible choice) equivalence  
        \[ \calD \simeq \LiuZheng(\calC,E,P,I,\calD|_{\calC^{\op}}), \]
	of 6-functor formalisms extending the natural identification 
	$\calD|_{\calC^{\op}}\simeq\LiuZheng(\calC,E,P,I,\calD|_{\calC^{\op}})\mid_{\calC^\op}$.
        \end{theorem}

\begin{remark}
	The reference \cite{dauser2024uniquenesssixfunctorformalisms} also characterizes 3-functor formalisms coming from the Liu--Zheng construction, but one has to be more careful defining the notions of $\calD$-proper and $\calD$-\'etale maps in the more general setup (see \cite[Definition 2.6, Proposition 2.13]{dauser2024uniquenesssixfunctorformalisms}). 
\end{remark}

        We use \Cref{thm: LiuZhengUniquenesss} throughout the paper to show that various different constructions give rise to the same 6-functor formalism. 
	We will also repeatedly use the following. 
        \begin{corollary}{\label{cor: LiuZhengRestriction}}
        Let $(\calC,E,P,I)$ be a Nagata setup, and let $E' \subset E$ be another set of edges such that $(\calC,E') \ra (\calC,E)$ defines a morphism of geometric setups.
We let $P' \subset P$ and $I' \subset I$ and suppose that $(\calC,E',P',I')$ is a Nagata setup for $(\calC,E')$.
Given a functor 
        \[ \calD_{0}: \calC^{\op} \ra \CAlg(\widehat{\Cat}) \]
        we then have a natural equivalence
        \[ \LiuZheng(\calC,E,P,I,\calD_{0})|_{\Corr(\calC,E')} \simeq \LiuZheng(\calC,E',P',I',\calD_{0}). \]
        \end{corollary}

	\begin{remark}
		\label{conventions on calD0}
		From now on, if $\calD$ is a $6$-functor formalism on a geometric setup $(\calC,E)$, we let $\calD_0$ denote the functor obtained from precomposing along the map $\calC^\op\to \corr(\calC,E)$.
	\end{remark}

    \subsubsection{Extending 6-functor formalisms}
    For our applications, we need to work with 6-functor formalisms defined for stacks. 
    When one works with stacks the Nagata compactification theorem rarely holds.
    A key insight (and the main motivation for the development of the theory) is that 6-functor formalisms can (and should) be extended using descent.

    Suppose that $(\calC,E)$ is a geometric setup such that $\calC$ is a site endowed with a Grothendieck topology $\tau$.
Let $\calX$ be the topos of $\tau$-sheaves of anima on $\calC$. 
We write $E_{\mathrm{rep}}$ (or also $E^{\mathrm{rep}}$) for the set of edges in $\calX$ that are relatively representable in $E$. 
		In other words, the set of edges $[f: X \ra Y] \in \calX$ such that, after pulling back to any $g: Y' \ra Y$ with $Y' \in \calC$, it lies in $E$. 
    Recall that a 3-functor formalism is sheafy in the sense of \cite[Definition 3.4.1]{HeyerMann} if the induced map $\calD: \calC^{\op} \ra \widehat{\Cat}$ is a $\tau$-sheaf. 
    The first extension procedure is as follows. 
	\begin{construction}[{\cite[Proposition~3.4.2]{HeyerMann}}]
	{\label{cons: ExtensiontoRepresentableMaps}}	Let $\calD$ be a sheafy 3-functor formalism on a geometric setup $(\calC,E)$, suppose that $\tau$ is a subcanonical Grothendieck topology on $\calC$.
		The following statements hold.
		\begin{enumerate}
			\item The natural map $(\calC,E) \ra (\calX,E_{\mathrm{rep}})$ is a morphism of geometric setups.
			\item There is a unique sheafy 3-functor formalism $\calD'$ on $(\calX,E_{\mathrm{rep}})$ extending $\calD$.
			\item The 3-functor formalism $\calD'$ is obtained as a right Kan-extension along the map 
				\[\Corr(\calC,E)^\otimes\to \Corr(\calX,E_{\on{rep}})^\otimes.\]
		\end{enumerate}
A similar result is true for the topos of hypersheaves if one assumes that $\calC$ is hyper-subcanonical. 
	\end{construction}

    Intuitively, this construction enlarges the set of objects over which the 3-functor formalism is defined, but one does not obtain truly new $!$-able edges. 
	It is useful to explicitly describe what this construction does on objects. 
	We record this with the following statement.
	\begin{lemma}[{\cite[Proposition~A.5.16]{Mann2022padic6functors}}]{\label{lemma: rightKanextensionDoesWhatisExpected}}
		In the situation of \Cref{cons: ExtensiontoRepresentableMaps}, and with the notational convention of \Cref{conventions on calD0}, the right Kan extension of $\calD_{0}$ along the embedding $\calC^{\op} \hookrightarrow \calX^{\op}$ is naturally equivalent to $\calD_{0}'$. 
	\end{lemma}

	\begin{remark}{\label{rem: extensionofPresentableisPresentable}}
		We note that the procedure of right Kan extension will preserve the property of being a presentable 3-functor formalism (see \cite[Remark~3.4.3]{HeyerMann}).
	\end{remark}

        The second type of extension procedure goes as follows.
	Intuitively, these extension procedures enlarge the set of $!$-able maps (i.e., $\calC_E$). 

	\begin{construction}[{\cite[Proposition~3.4.8.(i)]{HeyerMann}}]{\label{cons: ExtensiontoFineMapsI}}
		Let $\calD$ be a presentable 3-functor formalism on a geometric setup $(\calC,E)$. 
		We let $E \subseteq E'$ be another collection of morphisms so that $(\calC,E')$ is a geometric setup, and such that for every $[f:Y\to X]\in \calC_{E'}$ there is a universal $\calD^*$-cover $g:X'\to X$ such that the pullback of $f$ along $g$ lies in $\calC_E$. 
		In this situation, $\calD$ extends uniquely (up to contractible choice) along the map $\Corr(\calC,E)^\otimes \ra \Corr(\calC,E')^\otimes$ to a presentable 3-functor formalism $\calD'$ on $(\calC,E')$.
	\end{construction}

	There is also an extension procedure to access the set of so-called \textit{fine maps}. 
	We will discuss in \Cref{lets define fine finally}.
	\begin{construction}[{\cite[Proposition~3.4.8.(ii)]{HeyerMann}}]{\label{cons: ExtensiontoFineMaps}} 
		Let $\calD$ be a presentable 3-functor formalism on a geometric setup $(\calC,E)$. 
		We let $E \subseteq E'$ be another collection of morphisms so that $(\calC,E')$ is a geometric setup and such that for every $[f:X\to Y]\in E'$ there is a map $g:X'\to X$ such that 
		\begin{itemize}
			\item[a)] $g$ and $f\circ g$ are in $E$. 
			\item[b)] $g$ is a universal $\calD^!$-cover.
		\end{itemize}
		In this situation, $\calD$ extends uniquely (up to contractible choice) along the map $\Corr(\calC,E)^\otimes \ra \Corr(\calC,E')^\otimes$ to a presentable 3-functor formalism $\calD'$ on $(\calC,E')$.
		
	\end{construction}

    Finally, one can use \Cref{cons: ExtensiontoRepresentableMaps} and iterate uses of \Cref{cons: ExtensiontoFineMapsI} and \Cref{cons: ExtensiontoFineMaps} to obtain the following extension procedure. 
	\begin{construction}[{\cite[Theorem~3.4.11]{HeyerMann}}]
		\label{shriek hull}
		Fix a geometric setup $(\calC,E)$, a subcanonical topology $\tau$ on $\calC$ and a sheafy presentable 3-functor formalism $\calD$ on $(\calC,E)$. 
		Let $\calX$ denote the topos of $\tau$-sheaves on $\calC$.
        Then there is a smallest collection of edges \(E'\) in \(\calX\) with the following properties:
        \begin{enumerate}
            \item The inclusion \(\calC\injto\calX\) defines a morphism of geometric setups \((\calC,E)\to(\calX,E')\) and \(\calD\) extends uniquely to a sheafy 6-functor formalism on \((\calX,E')\).
            \item \(E'\) is \(*\)-local on the target: Let \(f\from Y\to X\) be a map in \(\calX\) whose pullback to every object in \(\calC\) lies in \(E'\); then \(f\) lies in \(E'\)
            \item \(E'\) is !-local: Let \(f\from Y\to X\) be a map in \(\calX\) that is !-locally on the source or target in \(E'\); then \(f\in E'\)
            \item \(E'\) is tame: Every map \(f\from Y\to X\) in \(E'\) with \(X\in\calC\) is !-locally on the source in \(E\).
        \end{enumerate}
		We denote this smallest collection $E'$ by $E_{\calD,!}$ and call it the \emph{$!$-able hull of $E$ with respect to $\calD$}.  
		If the context is clear, we omit $\calD$ from the notation. 
		Moreover, by abuse of notation we still refer to morphisms in $E_{\calD,!}$ as $!$-able.
	\end{construction}

	\begin{remark}
		We warn the reader that one can easily be mislead by \Cref{cons: ExtensiontoFineMapsI} into thinking that in order to verify whether a map $f:X\to Y$ is $!$-able (i.e., $f\in E_{\calD,!}$) or not it suffices to do so on a universal $\calD^*$-cover $Y'\to Y$ of $Y$.
		This is not the case since this is not what \Cref{shriek hull}.(2) is saying.
		One may also hope that if we let $X'=X\times_Y Y'$ and the resulting map $[f':X'\to Y']\in E_{\calD,!}$, then one can at the very least apply \Cref{cons: ExtensiontoFineMapsI} once again in order to ensure $f$ becomes $!$-able. 
		This might not always be possible.
		Indeed, if we let $E_*$ denote all maps that lie on $E_{\calD,!}$ after pullback by a universal $\calD^*$-cover, then it might not be true that $(\calX,E_*)$ is a geometric setup. 
		The challenge becomes finding $E'$ satisfying that
		\begin{enumerate}
			\item $(\calX,E')$ is a geometric setup,
			\item $E'\subseteq E_*$, and
			\item $f\in E'$.
		\end{enumerate}
	The first two condition are in tension against each other.
	\end{remark}

        \subsection{Analytification as a 6-Functor Formalism}
	\label{section-natural-transforms}
        Throughout this document, we compare various 3-functor formalisms. 
	 Of primary interest will be an algebraic 3-functor formalism $\calD^{\alg}$ and an analytic one denoted $\calD^{\an}$. 
	Since we work with stacks, it is important to make sure that analytification is functorial and respects the various higher coherences of the six operations. 
	This boils down to constructing a natural transformation of 3-functor formalisms $\calD^{\alg}\Rightarrow \calD^{\an}$.
	An efficient solution to construct such natural transformations is to organize them into their own 3-functor formalism, as we explain below. 
	In this section, we put ourselves in a formal context which we will later apply in \S \ref{subsec: DiamondFunctor} to the case of interest. 

	We start by recalling a very simple situation.
	Fix a category $\calC$ admitting finite limits and endowed with two functors 
	\[\calD^\alg_{0},\calD^\an_{0}\from\calC^\op\to\CAlg(\widehat{\Cat}).\] 
	For $[f: X \ra Y]\in \calC$, we write the pullback maps (i.e., $\calD_0^?(f)$) as 
	\[ f^{*\an}: \calD_0^{\an}(Y) \ra \calD_0^{\an}(X) \]
	\[ f^{*\alg}: \calD_0^{\alg}(Y) \ra \calD_0^{\alg}(X), \]
	respectively. 
	Suppose that we have a natural transformation \(\eta^*\from\calD^\alg_{0} \to\calD^\an_{0}\).
	Observe that this can be encoded in a single functor \(\calD_0\from(\calC\times\Delta^1)^\op\to\CAlg(\widehat{\Cat})\). 
	Indeed, \((\calC\times\Delta^1)^\op\simeq\calC^\op\times\Delta^1\).
	Let us unravel what \(\calD_0\) does.
	\begin{example}
		\label{remark-factoring-morphisms}
	Let $\calC^\Delta=\calC\times \Delta^1$ and let $\pi_\calC$ and $\pi_\Delta$ denote the natural projections onto $\calC$ and $\Delta^1$.
	Let $\iota_i:\calC\times\{i\}\to \calC^\Delta$ for $i\in \{\an,\alg\}$ denote the natural inclusions. 
	Objects in \(\calC\times\Delta^1\) are all of the form \((c,\an)=\iota_\an(c)\) or \((c,\alg)=\iota_\alg(c)\) for \(c\in\mathcal{C}\).
	
	We have two wide subcategories $\calC^\Delta_H,\calC^\Delta_V\subseteq \calC^\Delta$ (horizontal and vertical) spanned by morphisms of the form $(g,\id)$ and $(\id,\eta)$, respectively.
Here $\eta:\an\to \alg$ is the unique non-identity morphism in $\Delta^1$.
	Moreover, any morphism $h:(c,i)\to (d,j)$ factors canonically as $h=\eta_{h,d}\circ f$ and $g\circ \eta_{h,c}$ where $\eta_{h,c},\eta_{h,d}\in \calC_V$ and $f,g\in \calC_H$.
 
	We can explicitly write $f:=\iota_i(\pi_\calC(h))$ and $g=\iota_j(\pi_\calC(h))$.
Further, $\eta_{h,?}=(\id_?,\id_i)$ when $i=j$ or $(\id_?,\eta)$ if $i=\an$ and $j=\alg$ for $?\in\{c,d\}$.
	We note that since $\calC$ has finite limits $\calC\times \Delta^1$ also does, and the following identity holds 
	\[(d,\alg)\times_{(c,\alg)} (c,\an)=(d,\an).\]
	For an object $c\in \calC$, we write $\eta_c\defined (\id_c,\eta):(c,\an)\to (c,\alg)$.
	Then we have \[\calD_0(c,\alg)=\calD_0^\alg(c) \text{ and } \calD_0(c,\an)=\calD_0^\an(c).\]
	Moreover, $\calD_0(\eta_{c})=\eta_{c}^*\from\calD_0^\alg(c)\to\calD_0^\an(c)$ encodes the desired natural transformation.
	\end{example}

We move onto the more abstract setting.
Since we wish compare several 3-functor formalisms, it will be convenient to work in slightly more generality than the simplicial set $\Delta^{1}$.

	Recall that the category of partially ordered sets (posets) can be regarded as a full subcategory of ${\Cat}$. 
	In our convention, if $K$ is a poset,  $x\leq y$ and $\calP$ denotes the (ordinary) category attached to $K$, then there is a unique morphism $x\leftarrow y$.
	We call any category in the essential image a \emph{posetal category}.

\begin{definition}{\label{defn: nicepathcategory}}
	A \emph{lex posetal category} $\calP$ is a \emph{posetal category} that has all finite limits.
\end{definition}

If $\calP$ is a lex posetal category we equip it with its Cartesian monoidal structure. 
In the language of posets, the Cartesian monoidal structure simply takes two elements to their supremum.
If $K$ is a small poset and $\calP$ is the attached posetal category, then we abuse notation and write 
$\{k\}_{k \in K} \subset \mathcal{P}$ for the wide subcategory that only has the identity morphisms.

	\begin{definition}{\label{defn: deltaextendedgeometricsetup}}
		Given a geometric setup $(\calC,E)$ and a lex posetal category $\mathcal{P}$ attached to a small poset $K$ as in \Cref{defn: nicepathcategory}, we define a new geometric setup that we denote $(\calC^{\calP},E^{\calP})$.
		We let $\calC^{\calP} \defined \calC\times \calP$ and we let $\calC_{E^{\calP}}\defined \calC_E \times \{k\}_{k \in K}$ (we note that this gives a well-defined geometric setup in the sense of \Cref{defn: geometricsetup} by our assumption that $\calP$ has finite limits).
		More generally, given a class of arrows $J$ defining a wide subcategory $\calC_J\subseteq \calC$ we let $J^{\calP}$ be the class of arrows defining the wide subcategory $\calC_{J^{\calP}}\defined \calC_J\times\{k\}_{k \in K}$ of $\calC^{\calP}$.
	\end{definition}

	One can check that for any lex posetal category $\calP$ and any geometric setup $(\calC,E)$, the pair \((\calC^{\calP},E^{\calP})\) forms a geometric setup. 
Moreover, if $E$ admits a suitable decomposition into classes of morphisms $I$ and $\calP$, then $I^{\calP}$ and $P^{\calP}$ form a suitable decomposition for $E^{\calP}$. 
We have insertion morphisms of geometric setups for all $k \in K$; namely we have inclusions 
	\begin{equation}{\label{eqn: morrphismofgeometricsetups0}}
		\iota_k:	(\mathcal{C},E) \rightarrow (\mathcal{C}^{\calP},E^{\calP}) 
	\end{equation}
	\[ c \mapsto (c,k), \]
	and, applying the functor $\Corr(-)$, we obtain symmetric monoidal functors for $k\in K$ 
\[\ins^{(\mathcal{C},E)}_{k}: \corr(\mathcal{C},E) \ra \corr(\mathcal{C}^{\calP},E^\calP).\]
Given $[f: X \ra Y ]\in \mathcal{C}$, we write $\phantom{}^{k}f:=\iota_k(f)$.
 We will omit the superscript $(-)^{(\mathcal{C},E)}$ when the context is clear.
	
	Given a 3-functor formalism \(\calD\from\corr(\calC^\calP,E^{\calP})\to \widehat{\Cat},\)
	by precomposing with $\ins_{k}$ for $k \in K$, we can extract 3-functor formalisms 
	\[ \phantom{}^{k}\calD\from\corr(\calC\times\{k\},E)\to \widehat{\Cat}.  \]

	 This gives rise to functors $\phantom{}^{k}f_{!}: \phantom{}^{k}\calD(X) \ra \phantom{}^{k}\calD(Y)$ for all $k\in K$. 
	 Finally, given a map $[j \leftarrow i] \in \calP$, we can relate the 3-functor formalism attatched to $i$ and $j$ by a natural transformation, which we denote by $\eta_{i \ra j}^{*}: \phantom{}^{j}\calD \ra \phantom{}^{i}\calD$.
To capture this, we record the following lemma.

	\begin{lemma}{\label{lemma embeddingLemma5.7}}
		For $\calP$ a lex posetal category, there is a symmetric monoidal equivalence
		\[ \mathrm{Corr}(\calC,E) \times \calP^{\op} \xrightarrow{\simeq} \mathrm{Corr}(\mathcal{C}^{\calP},E^{\calP}), \]
		which sends $(X,j)$ to $(X,j)$ and a morphism $(X \xleftarrow{f} Z \xrightarrow{g} Y, j \leftarrow i=i)$ to the diagram 
		\[ (X,j) \xleftarrow{(f,j \leftarrow i)} (Z,i) \xrightarrow{(g,\id_{j})} (Y,i). \]
        Here \(\calP^\op\) is given the coCartesian monoidal structure.
	\end{lemma}

	\begin{proof}
        Consider the geometric setup $(\calP,\{\id_{k}\}_{k \in K})$.
	Then \(\corr((\calP,\{\id_{k}\}_{k \in K}))^\otimes=\calP^{\op,\coprod}\).
        The claim now follows from \cite[Lemma~2.2.11]{HeyerMann}, observing that 
	\[(\calC^{\calP},E^{\calP})=(\calC,E)\times (\calP,\{\id_{k}\}_{k \in K})\]
	in the category of geometric setups.
	\end{proof}

	A map $[j \leftarrow i]\in \calP$ defines a morphism $\corr(\calC,E) \times (\Delta^{1})^{\op} \ra \corr(\calC,E) \times \calP^{\op}$ along which we can restrict $\calD$ to obtain $\corr(\mathcal{C},E) \times (\Delta^{1})^{\op} \ra \widehat{\Cat}$, which will define for us a natural transformation 
	\begin{equation}{\label{eqn: naturaltransformationofsixfunctors}}
		\eta^{*}_{i \ra j}: \phantom{}^{j}\calD \ra \phantom{}^{i}\calD
	\end{equation}
	of 3-functor formalism, as desired.
	Further restricting to $\calC_E\times (\Delta^{1})^{\op}$ produces natural transformations \(\phantom{}^{j}\calD_!\to\phantom{}^{i}\calD_!\), that compare $\phantom{}^j f_!$ and $\phantom{}^i f_!$. 

	More precisely, given $[f: X \ra Y] \in \calC_{E}$ and $[j\leftarrow i]\in \calP$, we get a map $\Delta^1\times(\Delta^1)^\op\to \calC_E\times (\Delta^1)^\op$ which gives rise to the following commuting diagram after applying $\calD$ 
	\begin{equation*}
		\begin{tikzcd}
			\phantom{}^{j}\calD(X) \arrow[d, "\phantom{}^{j}f_!"] \arrow[r, "\eta_{i \ra j}^*"] & \phantom{}^{i}\calD(X) \arrow[d, "\phantom{}^{i}f_!"] \\
			\phantom{}^{j}\calD(Y) \arrow[r, "\eta_{i \ra j}^*"]                       & \phantom{}^{i}\calD(Y).                     
		\end{tikzcd}
	\end{equation*}

	In this way, by considering 3-functor formalisms on \(\calC^{\calP}\), we can encode natural transformations between $3$-functor formalisms. 
	
	The advantage of regarding natural transformations as 3-functor formalisms is that we can apply the constructions of \S \ref{section on six functors} to construct such natural transformations.
	Let us discuss how to apply \Cref{construction: 3 functors} to this specific setup.
	\begin{lemma}{\label{lemma: HowtoGetaSixfunctorformalsimontheinterval}}
		Fix a geometric setup \((\calC,E)\) and suppose that we are given a suitable decomposition of $E$ into $I$ and $P$ and a lex posetal category $\calP$ attached to a poset $K$.
		Suppose that we are given a functor \(\calD_0\from(\calC^\calP)^\op\to\CAlg(\widehat{\Cat})\)
		such that \(I\) and \(P\) satisfy the assumptions in \Cref{construction: 3 functors} individually with respect to the functors $\phantom{}^{k}\calD_0$ for all $k\in K$. 
		Fix $[f:c\to d]\in \calC$ and $[j\leftarrow i]\in \calP$ and assume that the following additional compatibilities hold:
		\begin{enumerate}
			\item If $f\in I$, then \(\eta_{i \ra j}^* \phantom{}^{j}f_\natural \cong \phantom{}^{i}f_\natural  \eta_{i \ra j}^*\) via the natural map.
			\item If $f\in P$, then \(\eta_{i \ra j}^*\phantom{}^{j}f_* \cong \phantom{}^{i}f_*\eta_{i \ra j}^*\) via the natural map.
		\end{enumerate}
		Then the geometric setup $(\calC^{\calP}, E^{\calP})$ from \Cref{defn: deltaextendedgeometricsetup} and the functor $\calD_0$ satisfies the axioms of \Cref{construction: 3 functors}, and therefore we get a 3-functor formalism 
		\[\calD:=\LiuZheng(\calC^\calP,E^\calP,P^\calP,I^\calP,\calD_{0}) \from\corr(\calC^\calP,E^{\calP})\to \CAlg(\widehat{\Cat}).\] 
	\end{lemma}
	\begin{proof}
	This is straightforward.
	\end{proof}

	Another advantage of this rephrasing is that the intuition one has for 3-functor formalisms carries over to natural transformations between them.
	For example, one can consider how the $\eta^{*}_{i \ra j}$ interact with the category of kernels.
	Recall from \cite[Theorem 4.2.4]{HeyerMann} that, on any geometric setup $(\calC,E)$, the construction of kernel categories upgrades to a lax symmetric monoidal functor 
	\[\calK_{\calD,(-)}:\Corr(\calC,E)^\otimes\to \on{Cat}_2^\times \text{ with } S\mapsto \calK_{\calD,S}.\]
	Moreover, on arrows of the form $\{S\xleftarrow{f} T\xrightarrow{\id} T\}$ the induced $2$-functor
	\[f^*:\calK_{\calD,S}\to \calK_{\calD,T}\] 
	acts on objects as $X\mapsto X\times_S T$ and on morphism categories as 
	\[\calD(X\times_S Y)\xrightarrow{f_{X \times_{S} Y}^*} \calD((X\times_S Y)\times_S T).\]

	In our category $\calC^{\calP}$, this gives the following.

	\begin{lemma}{\label{lemma: analytificationinducesmaponKernelCategories}}
		Fix a 3-functor formalism $\mathcal{D}: \Corr(\calC^{\calP},E^{\calP}) \ra \widehat{\Cat}$, fix $S \in \mathcal{C}$, and fix $[j \leftarrow i ]\in \calP$.
We have a $2$-functor
		\[ \eta_{i \ra j}^{*}: \mathcal{K}_{S,\phantom{}^{j}\calD} \ra \mathcal{K}_{S,\phantom{}^{i}\calD} \]
		on the $2$-categories of kernels. 
		On objects $[X\to S]\in (\calC_E)_{/S}$ it is the identity and on categories of morphisms it is given by $\eta^{*}_{i \ra j, X \times_{S} Y}:  \phantom{}^{j}\calD(X \times_{S} Y) \ra \phantom{}^{i}\calD(X \times_{S} Y)$.
	\end{lemma}
	\begin{proof}
		This follows from applying \cite[Proposition~4.2.1 (i)]{HeyerMann} to the natural transform $\eta_{i \ra j}^{*}: \phantom{}^{j}\calD \ra \phantom{}^{i}\calD$ of (\ref{eqn: naturaltransformationofsixfunctors}).
	\end{proof}

	\Cref{lemma: analytificationinducesmaponKernelCategories} gives the following. 
	\begin{corollary}{\label{cor: DeltathreefunctorsgivesTransformonKernels}}
		Fix a 3-functor formalism, $\mathcal{D}: \Corr(\calC^{\calP} ,E^{\calP}) \ra \widehat{\Cat}$, an object $f:[X \ra S] \in \mathcal{C}_{E}$, and $[i \leftarrow j] \in \calP$. 
		The natural transformation map $\eta^{*}_{i \ra j,X}$ induces a functor 
		\[ \eta^{*}_{i \ra j, X}: \Suave_{\phantom{}^{j}f,\phantom{}^{j}\calD}(X) \ra \Suave_{\phantom{}^{i}f,\phantom{}^{i}\calD}(X), \]
		satisfying the property that $\eta_{i \ra j,X}^{*} \SD_{\phantom{}^{j}f} \simeq \SD_{(\phantom{}^{i}f)}\eta_{i \ra j,X}^{*}$.
In particular, if $\phantom{}^{j}f$ is $\phantom{}^{j}\calD$-suave (resp. $\phantom{}^{j}\calD$-smooth) then $\phantom{}^{i}f$ is $\phantom{}^{i}\calD$-suave (resp. $\phantom{}^{i}\calD$-smooth).
		A similar statement holds for prim objects and prim duality.
	\end{corollary}
	\begin{proof}
		This follows from \Cref{lemma: analytificationinducesmaponKernelCategories}, as functors between $2$-categories preserve adjoints between $1$-morphisms. 
		Smoothness follows from the fact that $\eta^{*}_{i \ra j,X}$ is symmetric monoidal, and that if the suave dual of $\mathbf{1}$ is invertible, its image under $\eta^{*}_{i \ra j,X}$ will remain invertible.
	\end{proof}

	We turn our attention to investigate how the operation $(\mathcal{C},E) \mapsto (\mathcal{C}^{\calP},E^{\calP})$ interacts with the various extension procedures described in the previous section. 
	In particular, we spell out what happens when we do \Cref{cons: ExtensiontoRepresentableMaps}. 
		Let $(\calC,E)$ be a geometric setup, where $\calC$ is endowed with a Grothendieck topology $\tau$, and let $\calP$ be a lex posetal category attached to a partially ordered set $K$.
		Now, $\calC^\calP$ inherits a topology $\tau^\calP$ whose associated topos is simply $\calX^\calP$.
		Further, $E^\calP_{\rep}=(E_\rep)^\calP$ and we have a morphism of geometric setups
		\[(\calC^\calP,E^\calP)\hookrightarrow (\calX^\calP,E_\rep^\calP)\]
		which is obtained from 
		\[(\calC,E)\hookrightarrow (\calX,E_\rep)\]
		by considering the product with $(-)\times \calP$.
\begin{proposition}
	\label{path-topos}
		Suppose we are given a 3-functor formalism 
		\[\calD:\Corr(\calC^{\calP},E^{\calP}) \ra \widehat{\Cat}\]
		such that, for all $k \in K$, the induced 3-functor formalism $\phantom{}^k\calD$ is sheafy. 
		Let $\calD'$ denote the 3-functor formalism obtained by applying \Cref{cons: ExtensiontoRepresentableMaps} to $\calD$.
		Let $(\phantom{}^k\calD)'$ denote the 3-functor formalism obtained by applying \Cref{cons: ExtensiontoRepresentableMaps} to $\phantom{}^k\calD$.
		Finally, let $\phantom{}^k(\calD')$ denote the restriction of $\calD'$ along $\ins_{k}^{(\calX,E_\rep)}$. 
		Then, $\phantom{}^k(\calD')$ and $(\phantom{}^k\calD)'$ are canonically equivalent 3-functor formalisms on $(\calX,E_\rep)$.
		In particular, for all $[j \leftarrow i]\in \calP$ we obtain a natural transformation of 3-functor formalisms on $(\calX,E_\rep)$ 
		\[\eta_{i \ra j}^{*}:(\phantom{}^j\calD)'\ra (\phantom{}^i\calD)'.\] 
\end{proposition}
\begin{proof}
	This follows directly from the uniqueness part of \Cref{cons: ExtensiontoRepresentableMaps}, since $\phantom{}^k(\calD')$ and $(\phantom{}^k\calD)'$ are sheafy and restrict to $\phantom{}^k\calD$ on $\Corr(\calC,E)$ for all $k \in K$.	
\end{proof}

	Applying \Cref{shriek hull} to $\Corr(\calC^{\calP},E^{\calP})$ is a little more subtle.
	The reason is that, given a 3-functor formalism $\calD$ on $(\calC^\calP,E^\calP)$, a morphism $[f: X \ra Y] \in \calC$, and $k \in K$, whether the map $\iota_{k}(f): (X,k) \ra (Y,k)$ is of universal $\calD^{! (\text{or }*)}$-descent is not in general the same as $f$ being of universal $\phantom{}^{k}\calD^{! (\text{or }*)}$-descent. 
	Indeed, for $\ins_{k}(f)$ to be of universal $\calD^{!(\text{or }*)}$-descent, one needs to ask that $f$ is of universal $\phantom{}^j\calD^{!(\text{or }*)}$-descent for all $[k \leftarrow j] \in \calP$, since the pullback of $\iota_{k}(f)$ along $(Y,j) \ra (Y,k)$ is precisely $\iota_{j}(f)$.

	In what follows, we consider conditions on a map in $E_{{}^k\calD,!}$ (the $!$-able hull of $E_{{}^k\calD,!}$) that automatically ensure it will also lie in $E_{{}^j\calD,!}$.
	The point is that, although the $!$-able hull $E_{{}^k\calD,!}$ could contain some morphisms that are $!$-able for a random reason, there should be a better behaved subset of maps that are $!$-able for a good reason.
	Those maps that are $!$-able for a good reason will automatically be $!$-able for sheafy 6-functor formalisms ${}^j\calD$ that receives a map from ${}^k\calD$.
	In practice, the subset of maps in $E_{{}^k\calD,!}$ that we will be interested in always have a good reason to be $!$-able.

	We restrict to presentable 3-functor formalisms.
This already ensures that $\calD(X)$ is complete and cocomplete.
    
    \begin{definition}
	    \label{defi:suaveprim-cover}
        Let \(\calD\from\Corr(\calC,E)\to\PrL\) be a sheafy presentable 3-functor formalism for a Grothendieck topology \(\tau\) on $\calC$.
        Write \(\calX\) for the category of \(\tau\)-sheaves in anima on \(\calC\).
	Assume now that we are given a class \(\tilde{E}\supset E_{\rep}\) of arrows in \(\calX\) for which $(\calX,\tilde{E})$ is a geometric setup. 
	Suppose that \(\calD\) uniquely extends to a sheafy presentable 3-functor formalism \(\tilde{\calD}\) on \((\calX,\tilde{E})\).

        \begin{enumerate}
		\item We say that $f$ is a \emph{$(\tilde{\calD},\tilde{E},\tau)$-suave cover} if it is $\tilde{\calD}$-suave and surjective in $\calX$. 
		\item We say that $f$ is a \emph{$(\tilde{\calD},\tilde{E},\tau)$-prim cover} if it is $\tilde{\calD}$-prim and $f_* \mathbf{1}$ is descendable. 
		\item We say that $f$ is a $(\tilde{\calD},\tilde{E},\tau)$\emph{-smooth cover}, if it is a $(\tilde{\calD},\tilde{E},\tau)$-suave cover and $\tilde{\calD}$-smooth.
        \end{enumerate}
    \end{definition}

    \begin{lemma}\label{lem: suave prim smooth cover stable under base change}
	    The base change of a \((\tilde{\calD},\tilde{E},\tau)\)-suave (resp. smooth, resp. prim) cover is a \((\tilde{\calD},\tilde{E},\tau)\)-suave (resp. smooth, resp. prim) cover.
    \end{lemma}
    \begin{proof}
        By \cite[Lemma 4.4.8]{HeyerMann}, the base change remains suave (resp. prim).
	By \cite[Remark 4.5.2]{HeyerMann}, and the fact that $f^*$ preserves invertible objects, the base change of a smooth map remains smooth.
        Since $\tau$-surjectivity is stable under base change, the suave and smooth cases are settled.
        
        Finally, for the prim case, note that given a pullback diagram
        \begin{equation*}
            \begin{tikzcd}
X' \arrow[d, "f'"] \arrow[r, "g'"] & X \arrow[d, "f"] \\
Y' \arrow[r, "g"]                  & Y.               
\end{tikzcd}
        \end{equation*}
        By \cite[Lemma 4.5.13.(ii)]{HeyerMann} we have \(g^*f_*\textbf{1}\simeq f'_*\textbf{1}\), and since \(g^*\) is symmetric monoidal, we deduce that \(f'_*\textbf{1}\) stays descendable by \cite[Corollary 3.21]{MatthewGaloisGroup}.
    \end{proof}

    \begin{definition}
	    \label{a standard cover}
	    Let the setup be as in \Cref{defi:suaveprim-cover}.
	    We say that a map $g:X\to Y$ is a \textit{standard $!$-cover} if $g=g_1\circ\dots \circ g_m$, where each $g_i$ is either a $(\tilde{\calD},\tilde{E},\tau)$-prim cover or a $(\tilde{\calD},\tilde{E},\tau)$-suave cover.
    \end{definition}

    \begin{definition}
	    \label{standard-local}
	    Let the setup be as in \Cref{defi:suaveprim-cover} and \Cref{a standard cover}.
        \begin{enumerate}[(i)]
		\item We say that \(\tilde{E}\) is \textit{standard $!$-local on the source} if whenever \(\tilde{f}\from \tilde{X}\to\tilde{Y}\in\calX\) and there is a standard $!$-cover \(\tilde{g}\from\tilde{X}'\to\tilde{X}\) such that \(\tilde{g}, \tilde{f} \circ \tilde{g}\in\tilde{E}\), then \(\tilde{f}\in\tilde{E}\).
		\item We say that \(\tilde{E}\) is \textit{standard $!$-local on the target} if whenever \(\tilde{f}\from \tilde{X}\to\tilde{Y}\in\calX\) and there is a standard $!$-cover \(\tilde{g}\from\tilde{Y}'\to\tilde{Y}\) such that the base change \(\tilde{f}'\from \tilde{X}\times_{\tilde{Y}}\tilde{Y'}\to\tilde{Y'}\in\tilde{E}\), then \(\tilde{f}\in\tilde{E}\).
        \end{enumerate}
    \end{definition}

	\begin{definition}
		\label{lets define fine finally}
		Let the setup be as in \Cref{defi:suaveprim-cover} and \Cref{a standard cover}, except that we also ask that $\tilde{E}=E_{\on{rep}}$.
		We let $E_\fine$ (or $E^{\on{fine}}$) denote the set of maps $[f:X\to Y]\in \calX$ satisfying that there is a standard $!$-cover $g:U\to X$ of the form $g=g_1\circ\dots \circ g_m$ such that 
		\begin{enumerate}
			\item Each $g_i$ is either a $(\calD,E_{\on{rep}},\tau)$-prim cover or a $(\calD,E_{\on{rep}},\tau)$-suave cover. 
			\item Each $g_i\in E_\rep$ and $f\circ g_i\in E_\rep$. 
		\end{enumerate}
		We call this set of maps the \textit{standard fine} maps.
	\end{definition}

	For most practical purposes working with the set of standard fine maps suffices.

	\begin{lemma}
		\label{the fine maps are geometric setup}
		The pair $(\calX,E_\fine)$ is a geometric setup.	
	\end{lemma}
	\begin{proof}
		It follows easily from \Cref{lem: suave prim smooth cover stable under base change} that $E_\fine$ is stable under base change.
		Let $\calS$ denote the set of maps of the form $g=g_1\circ\dots \circ g_m$ where each $g_i$ is in $E_\rep$ and it is either a $(\calD,E_\rep,\tau)$-prim cover or a $(\calD,E_\rep,\tau)$-suave cover. 
		Fix $[g:X\to Y]$ and $[f:Y\to Z]$ both in $E_\fine$. 
	Fix maps $X'\to X$ and $Y'\to Y$ in $\calS$ such that both are an atlas exhibiting that $g,f\in E_\fine$.
		That $E_\fine$ is stable under composition follows from contemplating the following diagram with Cartesian squares
		\begin{center}
		\begin{tikzcd}
			X'' \arrow{r} \arrow{dd}{\in \calS}  & Y'_X \arrow{r} \arrow{d}{\in \calS}  & Y' \arrow{d}{\in \calS} \arrow{rd}{\in E_\rep} \\
		  & X \arrow{r} & Y \arrow{r} & Z \\
		  X' \arrow{ur}{\in \calS} \arrow{urrr}[swap]{\in E_\rep}.
		\end{tikzcd}
		\end{center}
        Here we note that, by assumption $X' \ra Y$ is representable, since $X' \ra X$ is an atlas for $X \ra Y$, and therefore $X'' \ra Y'$ as the base-change is also representable.
Moreover, since $Y' \ra Y$ is an atlas for $Y' \ra Z$ the map $Y' \ra Z$ is representable by assumption.
This in turn implies that $X'' \ra Z$ is representable as claimed.
In particular, the diagram constructs $X''\to X$ in $\calS$, and, since $X'' \to Z$ is representable, the map $X''\to X$ exhibits that $f\circ g\in E_\fine$.
	That $E_\fine$ is right cancellative similarly follows from contemplating a commutative diagram with Cartesian square of the form 
		\begin{center}
		\begin{tikzcd}
			X'' \arrow{r}{\in \calS} \ar{ddr}& 	X'\arrow{r}{\in \calS} \arrow{dd}  &X \arrow{dd}\arrow{dr} \\
			   & & & Z \\
			   & Y'\arrow{r}{\in \calS}  & Y \arrow{ur}
		\end{tikzcd}
		\end{center}
		Here, $Y'\to Y$ is an atlas exhibiting that $Y\in (\calX_{E_{\on{fine}}})_{/Z}$, and $X''\to X'$ is an atlas exhibiting that $X'\in (\calX_{E_{\on{fine}}})_{/Z}$ since it is the composition of two maps in $E_{\on{fine}}$.
		The cancellative property of $E_\rep$ will show that $X''\to Y'$, and consequently $X''\to Y$, are in $E_\rep$. 
		This shows that $X''\to X$ is an atlas exhibiting that $X\to Y$ is in $E_\fine$.
	\end{proof}

	\begin{definition}
		\label{define-Artin}
		Fix the setup as in \Cref{shriek hull}, and fix $[f:X\to Y]\in \calX_{E_{\calD,!}}$.
	We say $f$ is $(\calD,E,\tau)$\emph{-Artin} if there is a representable $(\calD,E_{\rep},\tau)$-smooth cover $[g:U\to X]\in \calX_{E_\rep}$ such that $f\circ g\in \calC_{E_\rep}$.  
	If $\calD$, $\tau$ and $E$ are clear from the context we simply say that $f$ is Artin, and we write $E^{\on{Art}}$ for the set of Artin maps. 
	We are mostly interested in objects $X\in \calX$ whose structure map $X\to \ast$ is in $E^\Art$. 
	We let $\calX^\Art$ the class of such objects.
\end{definition}

\begin{remark}
	The same proof as in \Cref{the fine maps are geometric setup} shows that $(\calX,E^{\on{Art}})$ is a geometric setup. 	
\end{remark}

For us, instead of studying morphisms of universal $\calD^!$-descent, it will be more useful to consider prim and suave covers and their composition (see \Cref{a standard cover}).
The following statement is a version of \Cref{shriek hull} where we replace the role of the $!$-covers for the stronger notion standard $!$-cover of \Cref{a standard cover}.

\begin{proposition}
		\label{suavecito-gives-fine}
	Let the setup be as in \Cref{shriek hull}. 
	There exists a set $E^{\on{std}}_{\calD,\tau,!}$ of maps $[f\from X\to Y]\in \calX$ satisfying the following conditions
	\begin{enumerate}[(i)]
				\item The morphism $\calC\to \calX$ induces a morphism of geometric setups \((\calC,E)\hookrightarrow (\calX,E^{\on{std}}_{\calD,\tau,!})\) and $\calD$ extends uniquely along this morphism to a sheafy $6$-functor formalism on $(\calX,E^{\on{std}}_{\calD,\tau,!})$. 
				\item $E^{\on{std}}_{\calD,\tau,!}$ is \(*\)-local on the target, i.e, for any map $V\to Y$ with in $V\in \calC$ such that the pullback of $f$ to $V$ lies in $E^{\on{std}}_{\calD,\tau,!}$, then $f\in E^{\on{std}}_{\calD,\tau,!}$.
				\item $E^{\on{std}}_{\calD,\tau,!}$ is standard $!$-local on the source  (see \Cref{standard-local}).
				\item $E^{\on{std}}_{\calD,\tau,!}$ is standard $!$-local on the target (see \Cref{standard-local}).
				\item $E^{\on{std}}_{\calD,\tau,!}$ is tame, i.e, if \([f:V\to Y]\in E^{\on{std}}_{\calD,\tau,!}\) with \(Y\in\calC\), it is !-locally on the source in \(E\).
				\item $E^{\on{std}}_{\calD,\tau,!}$ is \textit{standard tame}, i.e, if \([f:V\to Y]\in E^{\on{std}}_{\calD,\tau,!}\) with \(Y\in\calC\), it is standard !-locally on the source in \(E\).
				\item $E^{\on{std}}_{\calD,\tau,!}$ is minimal among the class of edges that satisfy $(i)-(v)$. 
		\end{enumerate}
    \end{proposition}

Before proving this proposition we make the following definition and remarks.

	\begin{definition}
		\label{define-standard-Shriek}
We refer to maps in $E^{\on{std}}_{\calD,\tau,!}$ as $(\calD,E,\tau)$\emph{-standard $!$-able} maps.
If $\calD$ and $\tau$ are clear from the context then we will simply say that $f$ is $E$-\emph{standard $!$-able} and write $E^{\on{std}}_!$ for $E^{\on{std}}_{\calD,\tau,!}$.
\end{definition}
\begin{remark}{\label{rem: standard!ableis!able}}
	We note that condition (v) of \Cref{suavecito-gives-fine} tautologically follows from condition (vi).
	Nevertheless, condition (v) is useful to formulate the minimality property that we are interested in.
	In general, the collection of classes of maps that satisfy (v) and (vi) might differ.
	Now, the $!$-able hull $E_{\calD,!}$ of \Cref{shriek hull} satisfy the axioms (i)-(v) of \Cref{suavecito-gives-fine}. 
	In particular, the minimality claim in \Cref{suavecito-gives-fine} gives us an inclusion $E^{\on{std}}_{\calD,\tau,!} \subset E_{\calD,!}$ allowing us to see that \Cref{suavecito-gives-fine} is a controlled refinement of the extension procedure in \Cref{shriek hull}.
This is because suave and prim covers are universal !-covers.
Indeed, this is the content of \Cref{suave-implies-fine} in the suave case (using surjectivity of \(f\) and the assumed sheafiness of $\calD$ to see that \(f^*\) is conservative) and \Cref{prim standard !-able} in the prim case. 
\end{remark}

\begin{remark}
\label{iv-is-automatic}
    Using \Cref{lem: suave prim smooth cover stable under base change}, one can easily see that a set of edges $E'$ in $\calX$ satisfying condition (i) and (iii) of \Cref{suavecito-gives-fine} automatically satisfies condition (iv). 
\end{remark}

    \begin{proof}[Proof of \Cref{suavecito-gives-fine}]
We give a construction of $E^{\on{std}}_{\calD,\tau,!}$ to show its existence.
The construction and proof will follow almost verbatim the construction in \cite[Theorem 3.4.11]{HeyerMann}, for this reason we do not give full details and only provide a sketch that highlights the differences.
With notation as in \cite[Theorem 3.4.11]{HeyerMann}, we let $A$ denote the class of all collections of edges $E'\subseteq \calX$ satisfying $(i)$ and $(v)$ in \Cref{suavecito-gives-fine}, so that \(\calD\) extends uniquely to \((\calX,E')\).
This class is stable under filtered unions and contains $E_\rep$.
As in \cite[Theorem 3.4.11]{HeyerMann}, one shows the following statement
\begin{itemize}
\item[(*)] For every $E'\in A$, there is a minimal way to enlarge it to $E'_!\in A$ such that $E'_!$ satisfies $(iii)$  of \Cref{suavecito-gives-fine}. 
	In particular, \(\calD\) extends uniquely to \((\calX,E'_!)\) by the assumption $E'_!\in A$ and, by \Cref{iv-is-automatic} $E'_!$ also automatically satisfies $(iv)$ of \Cref{suavecito-gives-fine}.
\end{itemize}
        To argue, it suffices to replace in the argument for \cite[Theorem 3.4.11]{HeyerMann}, every application of !-locality (i.e., the condition in \cite[Theorem 3.4.11 (iii)]{HeyerMann}) with our variant notion of standard $!$-locality.
	As in their proof, one uses \cite[Proposition 3.4.8 (ii)]{HeyerMann} (or \Cref{cons: ExtensiontoFineMaps}) to verify that there is a unique way of extending $\calD$ to $(\calX,E'_!)$.

As in \cite[Theorem 3.4.11]{HeyerMann} the next statement to show is the following.
\begin{itemize}
	\item[(**)] For every $E'\in A$ that satisfies $(iii)$ of \Cref{suavecito-gives-fine} (and consequently $(iv)$ using \Cref{iv-is-automatic}), there is a minimal way to enlarge it to $E'_*\in A$ such that $E'_*$ satisfies $(ii)$ of \Cref{suavecito-gives-fine}.
This extension has the property that \(\calD\) extends uniquely to \((\calX,E'_*)\).
\end{itemize}
This part of the argument follows verbatim.
Finally, one defines recursively $E_0=E_\rep$, $E_{2n+1}:=(E_{2n})_!$ (as in $(*)$) for $n\geq 0$, $E_{2n+2}:=(E_{2n+1})_*$ (as in $(**)$) and $E^{\on{std}}_{\calD,\tau,!}=\bigcup E_n$.
To show minimality of $E^{\on{std}}_{\calD,\tau,!}$ it suffices to show that any other $E'$ contains $E_n$ for all $n$, which can be shown by induction and using the minimality of the closure operations $(*)$ and $(**)$.
	\end{proof}

    \begin{remark}
        The class of edges $E^{\on{std}}_{\calD,\tau,!}$ is fairly abstract, for practical purposes it suffices to work with a much smaller class of edges in $\calX$.
	Note that Condition (i) and (ii) ensure that all edges in $E_\rep$ are standard !-able.
	Moreover, condition (iii) will then tell us that all standard fine maps (\Cref{lets define fine finally}) are standard !-able.
	In the main body of the article, knowing that $E_\fine\subseteq E^{\on{std}}_!$ will suffice, since all the edges we are interested in will lie in $(\calX,E_\fine)$.
    \end{remark}

We record the following variant which will also play a role later on. 

\begin{proposition}
		\label{suavecito-gives-fine take number 666}
	Let the setup be as in \Cref{shriek hull} and \Cref{suavecito-gives-fine}. 
	Let $F$ be a set of arrows in $\calX$ containing $E_\rep$ and such that $(\calX,F)$ is a geometric setup. 
	There exists a set $E^{\on{std}}_{\calD,\tau,F,!}$ of maps $f\from X\to Y$ satisfying the following conditions
	\begin{enumerate}[(i)]
				\item The morphism $\calC\to \calX$ induces a morphism of geometric setups \((\calC,E)\hookrightarrow (\calX,E^{\on{std}}_{\calD,\tau,F,!})\), and $\calD$ extends uniquely along this morphism to a sheafy $6$-functor formalism on $(\calX,E^{\on{std}}_{\calD,\tau,F,!})$. 
				\item $E^{\on{std}}_{\calD,\tau,F,!}$ is $F$-\(*\)-local on the target, i.e, given $[f:X\to Y]\in F$ such that for any map $V\to Y$ with $V\in \calC$ the pullback of $f$ to $V$ lies in $E^{\on{std}}_{\calD,\tau,F,!}$, then $f\in E^{\on{std}}_{\calD,\tau,F,!}$.
				\item $E^{\on{std}}_{\calD,\tau,F,!}$ is standard $!$-local on the source.
More precisely, if $[f:X\to Y]$ is an arrow in $F$ such that we may find a standard $!$-cover $g:U\to X$ with $f\circ g$ and $g$ in $E^{\on{std}}_{\calD,\tau,F,!}$, then $f$ is also in $E^{\on{std}}_{\calD,\tau,F,!}$.
				\item $E^{\on{std}}_{\calD,\tau,F,!}$ is both standard $!$-local local on the target.
				\item $E^{\on{std}}_{\calD,\tau,F,!}$ is tame, i.e, if \([f:V\to Y]\in E^{\on{std}}_{\calD,\tau,!}\) with \(Y\in\calC\), it is !-locally on the source in \(E\).
				\item $E^{\on{std}}_{\calD,\tau,F,!}\subseteq F$.
				\item $E^{\on{std}}_{\calD,\tau,F,!}$ is minimal among the class of edges that satisfy $(i)-(vi)$. 
		\end{enumerate}
    \end{proposition}
    \begin{proof}
	    A variant of the construction and proof as in \cite[Theorem 3.4.11]{HeyerMann} (or \Cref{suavecito-gives-fine}) shows this more general case. 
    Indeed, we let $E_0=E_\rep$, $E_{2n+1}:=(E_{2n})_{!_F}$ where $(E_{2n})_{!_F}$ is an $F$-variant of (*), $E_{2n+2}:=(E_{2n+1})_{*_F}$ where $(E_{2n+1})_{*_F}$ is an $F$-variant of (**), finally $E^{\on{std}}_{\calD,\tau,F,!}=\cup E_n$. 
    The proof of \cite[Theorem 3.4.11]{HeyerMann} shows that if $E'$ is a geometric setup in $A$ (with notation as in \Cref{suavecito-gives-fine}), then letting $E''$ be those maps that are standard $!$-locally on the source in $E'$ gives another geometric setup with $E'\in A$. 
    Moreover, if $E'$ also satisfied condition $(vi)$, then $E''\cap F$ is still in $A$ and satisfies $(vi)$. 
    Replacing the role of $E'$ by $E''\cap F$ and passing to unions gives rise to $(E')_{!_F}\in A$ which satisfies $(vi)$, $(iii)$ and $(iv)$.

    Similarly, the proof of \cite[Theorem 3.4.11]{HeyerMann} shows that if $E'$ is in $A$ and satisfies $(iii)$ and $(iv)$, then letting $E''$ be the set of arrows whose pullback to $\calC$ is in $E'$ gives a geometric setup in $A$.
    If $E'$ also satisfies $(vi)$ then $E''\cap F$ is in $A$ and satisfies $(ii)$ and $(vi)$.
    \end{proof}

    We now explain why standard !-able maps behave better with respect to analytification.
    We have the following lemma on prim and suave covers.
    \begin{lemma}\label{lem: suave and prim cover behave nice under analytification}
	    Let $(\calC,E)$ be a geometric setup endowed with a Grothendieck topology $\tau$.
	    Let \(\calP\) be a lex posetal category with underlying poset \(K\), fix \(q\in\calP\), and fix \[\calD^\calP\from\Corr(\calC^\calP,E^\calP)\to\PrL\] a sheafy presentable 3-functor formalism. 
        We write \({}^q\calD\defined\calD^\calP\circ\ins_q\) as a presentable 3-functor formalism on \(\Corr(\calX,E_\rep)\).
	Let $E_\rep\subseteq E'\subseteq \ins_q^{-1}(E_{\calD^\calP,!})$ be a collection of arrows in $\calX$ such that $(\calX,E')$ is a geometric setup and ${}^q\calD$ extends uniquely to $(\calX,E')$.
	Then, \([f:X\to S]\in E'\) is a \(({}^q\calD,E,\tau)\)-suave (resp.
smooth, resp. prim) cover if and only if \(\ins_q(f)\in E^\calP\) is a \((\calD^\calP,E^\calP,\tau^\calP)\)-suave (resp. smooth, resp. prim) cover.
    \end{lemma}
    \begin{proof}
        Observe that for any \(S\in\calX\), the functor \(\ins_q\) induces a 2-functor \(\ins_q\from\calK_{{}^q\calD,S}\to\calK_{\calD^\calP,\ins_q(S)}\) using \cite[Proposition 4.2.1 (ii)]{HeyerMann}.
	Using the explicit description of the kernel $2$-category (see \cite[Definition 4.1.3]{HeyerMann}), we obtain that this is even a fully faithful functor of 2-categories whose essential image is precisely $\ins_q([\calX_{E'}]_{/S})$.
	In other words, for all $X,Y\in [\calX_{E'}]_{/S}$ we have an equivalence 
	\begin{equation}
		\label{functor categories for the argument}
	\ins_q:\on{Fun}_{\calK_{{}^q\calD,S}}(X,Y)\simeq \on{Fun}_{\calK_{\calD^\calP,\ins_q(S)}}(\ins_q(X),\ins_q(Y))	
	\end{equation}
	of mapping categories.

	Since being suave (resp. prim) is defined (using adjunction see \cite[Definition 4.4.1]{HeyerMann}) in terms of the kernel category, it follows that $f\in E'$ is suave (resp. smooth, resp. prim) if and only if $\ins_q(f)$ is.
	Indeed, admitting a left or right adjoint is a property that can be verified after applying $\ins_q$ as in \eqref{functor categories for the argument}. 			
	Moreover, whether $f^!\textbf{1}$ is invertible, (resp. whether $f_*\textbf{1}$ is descendable) only depends on ${}^q\calD(X)\simeq \calD^\calP(\ins_q(X))$ (resp. ${}^q\calD(Y)\simeq \calD^\calP(\ins_q(Y))$), so the two notions agree.
	This already settles the prim case.
	To show the suave and smooth cases, observe that $\tau$ surjectivity of $f$ is also equivalent to $\tau^\calP$ surjectivity of $\ins_q(f)$.
    \end{proof}
    \begin{proposition}\label{prop: standard shriekable is ubiquitously shriekable}
        Let \(\calP\) be a lex posetal category with underlying poset \(K\), \(q\in\calP\), and \(\calD^\calP\from\Corr(\calC^\calP,E^\calP)\to\PrL\) a $\tau^\calP$-sheafy presentable 3-functor formalism. 
        We write \({}^q\calD\defined\calD^\calP\circ\ins_q\) which is a $\tau$-sheafy presentable 3-functor formalism on \(\Corr(\calC,E)\).
	    Then
        \begin{itemize}
            \item[a)] \(\ins_q^{-1}(E^{\calP,\on{std}}_{\calD^\calP,\tau^\calP,!})=E^{\on{std}}_{{}^q\calD,\tau,!}\), and
            \item[b)]  \(E^{\calP,\on{std}}_{\calD^\calP,\tau^\calP,!}=\bigcup_q \ins_q(E^{\on{std}}_{{}^q\calD,\tau,!})\).
        \end{itemize}
               In particular, $\ins_q(E^{\on{std}}_{{}^q\calD,\tau,!})\subseteq E^{\calP,\on{std}}_{\calD^\calP,\tau,!}$.
    \end{proposition}
    \begin{proof}
    This statement is shown by inspecting the proof of \Cref{suavecito-gives-fine}.
    Indeed, with notation as in \Cref{suavecito-gives-fine}, we show inductively that
\begin{itemize}
            \item[a)] \(\ins_q^{-1}(E^\calP_n)={}^qE_n\), and that
            \item[b)]  \(E^{\calP}_n=\bigcup_q \ins_q({}^qE_n)\).
        \end{itemize}
        Here $E^\calP_n$ and ${}^qE_n$ are defined recursively starting with 
        \begin{itemize}
            \item $E^\calP_0:=E^\calP_\rep$, ${}^qE_0:=E_\rep$,
            \item with $E^\calP_{2n+1}:=(E^\calP_{2n})_!$ and ${}^qE_{2n+1}:=({}^qE_{2n})_!$ (using $(*)$) for $n\geq 0$,
            \item with $E^\calP_{2n+2}:=(E^\calP_{2n+1})_*$ and ${}^qE_{2n+2}:=({}^qE_{2n+1})_*$ (using $(**)$) for $n\geq 0$.
        \end{itemize}

        The inductive argument consists in showing that the closure operations $(\ast)$ and $(\ast \ast)$ are compatible with $\ins_q$, and that they do not produce edges of the form $(X,i)\to (Y,j)$ with $i\neq j$. 
        That this is the case for $(\ast)$, follows from \Cref{lem: suave and prim cover behave nice under analytification}.
        The $(**)$ case follows from the explicit description of $(E^\calP_{2n+1})_*$ as the set of edges whose pullback to any element of $\calC^\calP$ lies in $E^\calP_{2n}$ (see proof of \cite[Theorem 3.4.11]{HeyerMann}).
        Indeed, it is easy to show that the base change of $(X,q)\to (Y,q)$ by any element in $\calC^\calP$ lies in $E^\calP_{2n}$ if and only if the similar statement is already true when we restrict to objects in $\calC\times \{q\}\subseteq \calC^\calP$.
        To wit, this follows from the canonical factorization
        \[(C,j)\to (C,q)\to (Y,q).\]

Finally, since $E^{\calP,\on{std}}_{\calD^\calP,\tau^\calP,!}=\bigcup_n E^{\calP}_n$ and $E^{\on{std}}_{{}^q\calD,\tau,!}=\bigcup_n {}^qE_n$, we can show the identities
\begin{itemize}
            \item[a)] \(\ins_q^{-1}(E^{\calP,\on{std}}_{\calD^\calP,\tau^\calP,!})=E^{\on{std}}_{{}^q\calD,\tau,!}\), and
            \item[b)]  \(E^{\calP,\on{std}}_{\calD^\calP,\tau^\calP,!}=\bigcup_q \ins_q(E^{\on{std}}_{{}^q\calD,\tau,!})\).
        \end{itemize}
by inspection.    
    \end{proof}

Summarizing, the previous discussion, we obtain the following.

\begin{proposition}
	\label{path-topos-sheafy}
	Let $(\calC,E)$ be a geometric setup, let $\calP$ be a lex posetal category attached to a poset $K$.
	Suppose that $\calC$ is endowed with a Grothendieck topology $\tau$.
		Suppose we are given a presentable 3-functor formalism 
		\[\calD:\Corr(\calC^{\calP},E^{\calP}) \ra \PrL\]
		such that, for all $k \in K$, the induced 3-functor formalism $\phantom{}^k\calD$ is sheafy with respect to $\tau$. 
		The following statements hold.
		\begin{enumerate}
			\item $\calD$ is sheafy with respect to $\tau^\calP$.
			\item If $[f:X\to Y]\in \calX$ is $(\phantom{}^k\calD,E,\tau)$-standard $!$-able, then $\ins_k(f)$ is $(\calD,\tau^\calP,E^\calP)$-standard $!$-able.
				In particular, $f$ is $(\phantom{}^j\calD,E,\tau)$-standard $!$-able for all $ j\geq k$.
            \item We have a morphism of geometric setups 
				\[\ins_k:(\calX,E^{\on{std}}_{^k\calD,\tau,!})\to (\calX^\calP,E^{\calP,\on{std}}_{\calD,\tau^\calP,!}).\]
			\item  	Let $\calD'$ denote the 6-functor formalism on $(\calX^\calP,E^{\calP,\on{std}}_{\calD,\tau^\calP,!})$ obtained from $\calD$ by unique extension.
				For $k\in K$, let $(^k\calD)'$ denote the 6-functor formalism on $(\calX,E^{\on{std}}_{^k\calD,\tau,!})$ obtained from $^k\calD$ by unique extension.
Then we have a unique equivalence $\calD'\circ \ins_k\simeq (^k\calD)'$ that restricts to the identity on \(\Corr(\calC,E)\).
			\item 	Given $[j\leftarrow i]\in \calP$ the natural map 
				\[\corr(\calX,E^{\on{std}}_{^j\calD,\tau,!})\times (\Delta^1)^\op\to \corr(\calX^\calP,E^{\calP,\on{std}}_{\calD,\tau^\calP,!})\] 
				induces a natural transformation of 6-functor formalisms on $\corr(\calX,E^{\on{std}}_{^j\calD,\tau,!})$  
				\[\eta^*_{i\to j}:(\phantom{}^j\calD)'\ra (\phantom{}^i\calD)'.\]

		\end{enumerate}
\end{proposition}

In other words, \Cref{path-topos-sheafy} states that given a natural transformation of $\tau$-sheafy 6-functor formalisms on a geometric setup $(\calC,E)$, one can always extend this natural transformation uniquely to the associated topos functorially on all maps that are standard $!$-able with respect to $E$.

\section{Analytification for stacks}
\label{s:Analyt for stacks section}

    \noindent
    \textbf{Soft preamble to the section:}
    In \cite[\S 27]{Sch17}, Scholze introduces the functor $\diamondsuit$ from the category of perfect schemes in characteristic $p$ to the category of small v-sheaves.
    This functor is well behaved and supports a functor $c^*$ of the form
    \[c_X^*:\calD_\et(X)\to \calD_\et(X^\diamondsuit),\]
    that compares \'etale sheaves on the scheme-theoretic setup, and \'etale sheaves on the analytic setup.
    As shown in \cite[Proposition 27.4]{Sch17}, $c^*$ commutes with lower $!$-functors.
    The purpose of this section is to construct $c^*$ (and the commutativity data for lower $!$-functors) for more general scheme-theoretic \'etale stacks. 
    Although analytification for stacks has been discussed in the literature before (see \cite[Appendix A]{AGLR22}), to the best of our knowledge, this is the first time that compatibility with lower-$!$-functors is thoroughly discussed beyond the case of schemes.  \\

    \noindent
    \textbf{Technical preamble to the section:}
    \begin{enumerate}
	    \item[\S  \ref{subsection algeb 6-ff}] Recalls various versions of 6-functor formalisms considered in the literature that deserve the name of ``\'etale cohomology'' in the scheme-theoretic context. 
		    This part is purely expository, but it serves the purpose of specifying the precise coefficient theory that we work with. 
	    \item[\S \ref{ss: analytic 6-ff go}] Recalls the theory of \'etale cohomology of diamonds. 
		    We discuss the precise $6$-functor formalism defined by Mann \cite{Mann2022NuclearSheaves}, but we also consider a variant which seems to work better with respect to analytification.
		    The construction of this variant itself is identical to Mann's, except that instead of working with the v-topology we work with the pro\'etale topology.  
	    \item[\S \ref{subsec: DiamondFunctor}] Provides a very general construction of Scholze's comparison functor $c^*$ together with coherent lower $!$-compatibilities in the context of scheme-theoretic \'etale stacks.  
		    We formulate two versions (\Cref{thm: DiamondAnalytificationUltimate} and \Cref{thm: DiamondAnalytificationUltimate2}), one is with respect to the pro\'etale topology and the other one is with respect to the more traditional v-topology. 
		    The first version has the feature that every relatively representable map of prestacks analytifies to a $!$-able map, but it departs from the convention in the literature of treating small v-stacks as the basic geometric object.
		    The second version has the feature that all considerations stay within the realm of small v-stacks, but it has the downside that we must impose restrictions on the type of scheme-theoretic \'etale stack that we are allowed to analytify.  
		    Indeed, we consider what we call \textit{resilient stacks} (\Cref{def: resiliant stacks}) and justify that in practice, all stacks we need to work with are resilient (\Cref{enough-resilients}). 
    \end{enumerate}

\subsection{Algebraic 6-functor formalisms}
\label{subsection algeb 6-ff}
We discuss four variants of 3-functor formalisms that deserve the name of \'etale cohomology of schemes (following \cite[Appendix to Lecture VII]{scholze6ff}).
Let us first fix our geometric setup to be $\calC=\Sch^{\qcqs}$ the category of qcqs schemes, and let $E=E_{\on{fin.exp.}}$ be the collection of separated maps of ``finite expansion'' which we recall below. 
	\begin{definition}[{\cite[Definitions~1.1,1.4]{HammacherExtensionof6functors}, \cite[Appendix to \S VII]{scholze6ff}}]
    We say that
	 \begin{enumerate}
		 \item A map of affine schemes $\Spec A\to \Spec B$ is of \emph{finite expansion} if we have a factorization $B\to B[x_1,\dots,x_n]\to A$ where the map $B[x_1,\dots x_n]\to A$ is integral.
		 \item A map of qcqs schemes $f:X\to Y$ is of \emph{finite expansion} if for every pair of open subsets $\Spec A\subseteq X$ and $\Spec B\subseteq Y$ such that $f(\Spec A)\subseteq \Spec B$ the induced map is of finite expansion. 
	 \end{enumerate}	
	\end{definition}
	The perfection morphism $X^\perf\to X$ is an important example of a morphism of finite expansion, since it is integral.
	We let $I_{\finexp}$ denote the open immersions and $P_{\finexp}$ the separated and universally closed maps (we will refer to these as proper maps\footnote{Strictly speaking $P_{\finexp}$ contains more maps than the standard literature considers as ``proper''.
Indeed, any integral map is in $P_{\finexp}$, but in standard references for algebraic geometry, proper maps are usually required to be of finite type, which we do not wish to enforce.}).
	Since open immersions are finitely presented, one easily sees that $I_{\finexp}\subseteq E_\finexp$.
	Moreover, it follows from \cite[Corollary 1.9]{HammacherExtensionof6functors} that $P_{\finexp}\subseteq E_{\on{fin.exp.}}$.

	\begin{proposition}[{\cite[Proposition~7.12]{scholze6ff}, \cite[Proposition~1.8, Theorem~1.17]{HammacherExtensionof6functors}}]
		Let $[X\to Y]\in\calC_E$.
Then the following statements hold.	
		\begin{enumerate}
			\item The map $f$ can be written as an integral map $[X\to X']$ composed with a separated map of finite presentation $[X'\to Y]$.
			\item If $f\in P_{\finexp}$, then it can be written as an integral map $[X\to X']$ composed with a proper map of finite presentation $[X'\to Y]$.
			\item The map $f$ can be written as an open immersion $[X\to \bar{X}]$ composed with $[\bar{X}\to Y]\in P_{\on{fin.exp.}}$.
		\end{enumerate}
	\end{proposition}

	From this, we can deduce the following.
	\begin{corollary}{\label{cor: finiteexpansiongeometricsetup}}
		The tuple $(\Sch^\qcqs,E_{\on{fin.exp.}},I_{\finexp},P_{\finexp})$ is a Nagata setup. %
	\end{corollary}

	We wish to use \Cref{construction: 3 functors}, which requires that we give as input a map
	\[\calD_0:\Sch^\qcqs\to \CAlg(\widehat{\Cat}).\]
This map will be provided from the general theory of sites. 
Fix $\Lambda$ an $\ell$-torsion ring. 
Recall that to any ordinary topos $\calX$ and any ordinary ring $\Lambda$, one can attach three versions of the derived category of sheaves of $\Lambda$-modules which come equipped with morphisms of commutative algebra objects in $\LinCat_\Lambda$,
	\begin{equation}
		\label{sequenceofthree}
		\fShv(\calX,\Lambda) \ra \calD(\calX,\Lambda) \ra \hat{\calD}(\calX,\Lambda).
	\end{equation}
	The first, $\fShv(\calX,\Lambda)$, is the category of ${\rm Mod}_\Lambda$-valued sheaves in the sense of Lurie, compare with \cite[Definition 6.2.2.6, Notation 6.3.5.16]{HTT}. 
	The second one $\calD(\calX,\Lambda)\subseteq \fShv(\calX,\Lambda)$ is the subcategory of those sheaves that are hypercomplete as explained before \cite[Lemma 6.5.2.9]{HTT} (this is equivalent to effective descent with respect to hypercovers by \cite[Theorem 6.5.3.12]{HTT}), which gives rise to a reflective localization map $\fShv(\calX,\Lambda) \to \calD(\calX,\Lambda)$ known as the hypercompletion, see \cite[Lemma 6.5.2.10]{HTT}. 
	For certain $\calX$ (like the \'etale site of scheme), this category can also be identified with the $\infty$-categorical enhancement of the unbounded derived category of sheaves of $\Lambda$-modules on $\calX$, (see \cite[\S~2.2, proof of Proposition 2.2.2.6]{GaitsgoryLurieWeilsConjectureI}, or \cite[proof of Proposition 7.1]{scholze6ff}).
	
	Finally, $\fShv(\calX,\Lambda)$ and $\calD(\calX,\Lambda)$ carry compatible standard t-structures (in cohomological notation: $(\fShv(\calX,\Lambda)^{\leq 0}, \fShv(\calX,\Lambda)^{\geq 0})$ and $(\calD(\calX,\Lambda)^{\leq 0}, \calD(\calX,\Lambda)^{\geq 0})$ respectively), and the third category of sheaves $\hat{\calD}(\calX,\Lambda)$ is defined as the left-completion with respect to either of these t-structures, (see \cite[\S3.3 and \S5]{BS15} and \cite[\S Lecture VII, \S Appendix to Lecture VII]{scholze6ff}). 

	Furthermore, whenever we have a geometric morphism of topoi $f:\calY\to \calX$, the pullback $f^*$ provides a symmetric monoidal functor for each of the sheaf theories considered above, and these pullbacks respect the morphisms (\ref{sequenceofthree}).
	This allow us to consider for each of the rules $\calD_0\in \{\fShv(-,\Lambda), \calD(-,\Lambda), \hat{\calD}(-,\Lambda)\}$ a map 
		\begin{align}
		\label{rule-forsymme}
		\calD_0:(\Sch^{\qcqs})^\op& \to \CAlg(\widehat{\Cat}) \\
			X& \mapsto \calD_0(X_\et,\Lambda) \nonumber.
		\end{align}

	Since the natural maps provide identifications (in cohomological notation)
	\begin{equation}{\label{eqn: identificationofboundedbelowcategories}}
    \fShv^{\geq n}(X_\et,\Lambda)\simeq \calD^{\geq n}(X_\et,\Lambda)\simeq \hat{\calD}^{\geq n}(X_\et,\Lambda),
    \end{equation}
we can unambiguously define $\calD_{\et}^{+}(X,\Lambda)$, the category of bounded below complexes, as the union over \(n\in\bbZ\).
Similarly, we can consider $\calD_\cons(X,\Lambda)\subseteq \calD_{\et}^{+}(X,\Lambda)$ the subcategory of bounded below complexes with constructible cohomology.
Alternatively, this category can be characterized as the smallest $\Lambda$-linear idempotent complete stable $\infty$-subcategory containing $j_!\Lambda$ for all maps $j:U\to X$ that are qcqs and \'etale. 
Since $f^*$ and $\otimes$ preserves constructibility, $\calD_\cons(-,\Lambda)$ also gives rise to a functor of the form described in \eqref{rule-forsymme}. 

\begin{convention}
	\label{conv: names for sheaf categories}
	When the coefficients are clear from the context, we will often simplify the notation as follows.
	\begin{enumerate}
		\item We let $\fShv(X):=\fShv(X_\et,\Lambda)$.
		\item We let $\calD(X_\et):=\calD(X_\et,\Lambda)$.
		\item We let $\calD_\et(X):=\hat{\calD}(X_\et,\Lambda)$.
		\item We let $\calD_\et^+(X):=\calD_\et^+(X,\Lambda)$.
		\item We let $\calD_\cons(X):=\calD_\cons(X,\Lambda)$.
	\end{enumerate}
    Let us recall how these functors upgrade to 3-functor formalisms.
\end{convention}

    \begin{proposition}[{\cite[Proposition~2.1, Theorem~7.15]{scholze6ff}}]
	    \label{scholzes-version}
	    The following statements hold.
	    \begin{enumerate}
		    \item The classes $I_{\finexp}$, $P_{\finexp} \subseteq E_{\on{fin.exp}}$ satisfy the hypothesis of \Cref{construction: 3 functors} with respect to $\{\fShv(-), \calD(-_\et), \calD_\et(-)\}$.
		    \item For $[f:X\to Y]\in I$, we have a commutative diagram
			   \begin{center}
			   \begin{tikzcd}
				   \fShv(X) \arrow{r} \arrow{d}{f_\natural}  & \calD(X_\et) \arrow{r} \arrow{d}{f_\natural}  & \calD_\et(X) \arrow{d}{f_\natural} \\
				   \fShv(Y)  \arrow{r} & \calD(Y_\et) \arrow{r} & \calD_\et(X).
			   \end{tikzcd}
			   \end{center}
		    \item For $[f:X\to Y]\in P$, we have a commutative diagram
			   \begin{center}
			   \begin{tikzcd}
				   \fShv(X) \arrow{r} \arrow{d}{f_*}  & \calD(X_\et) \arrow{r} \arrow{d}{f_*}  & \calD_\et(X) \arrow{d}{f_*} \\
				   \fShv(Y)  \arrow{r} & \calD(X_\et) \arrow{r} & \calD_\et(X).
			   \end{tikzcd}
			   \end{center}
 
		   \item Using \Cref{lemma: HowtoGetaSixfunctorformalsimontheinterval}, we obtain natural transformations 
			   \[\fShv(-)\Rightarrow \calD(-_\et)\Rightarrow \calD_\et(-)\]
			   of 6-functor formalisms on $(\Sch^{\qcqs},E_{\on{fin.exp}})$ with values in $\LinCat_\Lambda$.
	    \end{enumerate}
	        \end{proposition}

		Since we also wish to work with $\calD_{\cons}$, and the functor $f_*$, for $f$ a morphisms of finite expansion, does not always preserve constructibility, we will restrict the 3-functor formalisms considered above to a smaller geometric setup. 

		We consider the category $\calC=\PSch^\qcqs$ of perfect schemes that are qcqs over $\Spec k$ (we recall that $k$ denotes a fixed algebraically closed field in which $\ell$ is invertible). 
		Let $E_\pfp$ denote the set of separated and perfectly finitely presented morphisms inside $\PSch^{\qcqs}$. 
		The pair $(\PSch^{\qcqs},E_\pfp)$ defines a geometric setup and the natural inclusion 
		\[(\PSch^{\qcqs},E_\pfp)\to(\Sch^{\qcqs},E_{\on{fin.exp.}})\]
		is a map of geometric setups, along which we can precompose to obtain 6-functor formalisms as in \Cref{scholzes-version}. 
		Since it is well known (see \cite[Proposition 10.14]{Zhu25}) that, for perfectly finitely presented maps, the functor $f_!$ preserves constructibility, this leads to the following.

        \begin{proposition}[{\cite[Theorem~10.15]{Zhu25}}]{\label{prop: Dconsthreefunctorformalism}}
		The classes $I$ of open immersions and $P_\pfp$ of perfectly finitely presented proper morphisms are a suitable decomposition for $E_\pfp$.
        The functor $\Dcons(-)$ satisfies the assumptions of \Cref{construction: 3 functors} on the geometric setup $(\PSch^{\qcqs},E_\pfp)$, with respect to $I$ and $P$. 
	In particular, the functor extends along the embedding
        \[ \PSch^{\qcqs,\op} \hookrightarrow \Corr(\PSch^{\qcqs},E_\pfp) \]
        to a 3-functor formalism 
        \[ \Dcons(-): \Corr(\PSch^{\qcqs},E_\pfp) \rightarrow \LinCat_{\Lambda}^{\mathrm{sm}}.  \]
        \end{proposition}
        We note that this 3-functor formalism is not presentable nor a 6-functor formalism. 
	Indeed, the $*$-pullback functor does not always admit a right adjoint, since one rarely expect $f_*$ to preserve constructible complexes whenever $f$ is not perfectly finitely presented (see \cite[Remark~10.16]{Zhu25}).
    
	We summarize what we need in the following theorem.

        \begin{theorem}
		\label{collecting-formalisms pfprep}
		The functors $\Dcons(-)$, $\fShv(-)$, and $\Detale(-)$ upgrade to 3-functor formalisms on the geometric setup $(\PSch^{\qcqs},E_\pfp)$, where the latter two are presentable 3-functor formalisms. 
		We have natural transformations 
		\[ \Dcons(-) \Rightarrow \fShv(-) \Rightarrow \Detale(-) \]
        of 3-functor formalisms on $(\PSch^{\qcqs},E_\pfp)$. 
	Furthermore, the first functor is fully faithful and the second functor is a morphism in $\LinCat_\Lambda$.
        \end{theorem}

\subsubsection{Basic properties of algebraic 3-functor formalisms}
We now highlight some important structural properties of the 3-functor formalisms introduced above.
Recall that a functor $\calF$ is said to be finitary if it preserves filtered colimits.
Recall that a prestack $\calF$ on $\Sch$ that satisfies Zariski descent is completely determined by the covariant functor it induces on $\CAlg$.
We say that a prestack on $\Sch$ is a finitary $\tau$-sheaf if it satisfies $\tau$-descent for a topology refining the Zariski topology, and the induced functor on $\CAlg$ is finitary. We have the following.

	\begin{theorem}[{\cite[Theorem 2.2]{HS23}}]
		\label{thm-Hansen-Scholze}{\label{thm: Dconsarcsheaf}}
		The functor
		\[\calD_\cons(-):\Sch^{\qcqs,\op} \to \LinCat^\sm_\Lambda\]
		is a finitary arc-sheaf.
In particular, it is a finitary schematic v-sheaf.
	\end{theorem}

One advantage of working with the left-completed sheaf theory $\calD_\et(X)$ is that it has very strong descent properties.

	\begin{theorem}[{\cite[Theorem 5.7]{HS23}}]{\label{thm: schematicDetisavsheaf}}
		Fix $X\in \Sch^\qcqs$.
		\begin{enumerate}
			\item The map of topoi $(X_\proet)\to (X_\et)$ induces a fully faithful embedding
				\[\Detale(X) \subset \calD(X_{\proet}).\]
    \item The functors $\calD_{\et}(-)$ and $\calD_{\et}^{+}(-)$
		are sheaves of symmetric monoidal categories for the topology of universal submersions.
In particular, they are also v-sheaves.
		\end{enumerate}
	\end{theorem}

On the other hand, one does not have such nice descent properties when working with $\fShv$.
Indeed, as explained in \cite[Example 10.24]{Zhu25} even \'etale descent can fail when we work over general base fields $k$.
This is only one of the reasons where the assumption that $k$ is algebraically closed and of characteristic $p\neq \ell$ becomes important. 
In this situation, at least \'etale descent holds.

\begin{proposition}[{\cite[Proposition 10.25]{Zhu25}}]
	\label{Shv-is-sheafy-etale}
The 6-functor formalism	$\fShv(-)$ is sheafy for the \'etale topology on $\PSch$.
\end{proposition}

	The following statement explains that $\calD_\cons(X)$ and $\fShv(X)$ can be obtained from the other by applying a formal procedure. 

	\begin{proposition}
		\label{compact-generation-shv*}
		Fix $X\in \PSch^\qcqs$, the following statements hold.
		\begin{enumerate}
			\item The category $\fShv(X)$ is compactly generated.
			\item $\calD_{\cons}(X)\simeq \fShv(X)^\omega$, equivalently $\on{Ind}\calD_\cons(X)\simeq \fShv(X)$.   
		\end{enumerate}
	\end{proposition}
    \begin{proof}
    We first assume that $X$ is pfp over $\ast$; in this case, we have that $X$ has finite cohomological dimension and the result follows from \cite[Proposition~6.4.8]{BS15} and \Cref{lemma: BasicPropertiesofDconsDet} below.
In general, we may apply Noetherian approximation \cite[Tag~07SU]{StaProj} to obtain a presentation
    \[ X := \lim_{i \in I} X_{i}, \]
    where each $X_{i} \in \PSch^{\pfp}$ and the transition morphisms are affine.
Using \cite[Proposition~7.14]{scholze6ff}, we now have that 
    \begin{equation}{\label{eqn: colimitofpullbacksfShv}}
     \fShv(X) \simeq \colim_{i \in I} \fShv(X_{i}) 
    \end{equation}
    and 
    \begin{equation}{\label{eqn: colimitofpullbacksDcons}}
     \Dcons(X) \simeq \colim_{i \in I} \Dcons(X_{i}), 
    \end{equation}
    where the transition maps are given by $*$-pullback.
    
    Each of the transition morphisms in (\ref{eqn: colimitofpullbacksfShv}) preserves constructible sheaves and commutes with colimits; therefore, it follows from \cite[\S~1.9.2]{drinfeld2015compactgenerationcategorydmodules}, that $\fShv(X)$ is compactly generated and that we have an identification 
    \[ \fShv(X)^{\omega} \simeq \colim_{i \in I} \fShv(X_{i})^{\omega} \]
    However, the right-hand side, identifies with $\colim_{i \in I} \Dcons(X_{i})$ by the pfp case discussed above, and therefore the claim follows from (\ref{eqn: colimitofpullbacksDcons}).
    \end{proof}

Compact generation is certainly an appealing formal property of a category, since it implies that it is dualizable. 
Dualizability can be used to perform higher-algebraic manipulations, like defining and computing the categorical traces of an endofunctor. 
This advantage is one of the main reasons that \cite{Zhu25} uses $\fShv(X)$ and variants thereof as the chosen theory to work with.  \\

Having discussed the main advantages of working with one theory over the other, we record the following statement which explains that three of the variants give rise to the same category under some finiteness assumptions. 
\begin{theorem}[{\cite[Propositions~3.3.7]{BS15}, \cite[Theorem 7.13]{scholze6ff}}]{\label{lemma: BasicPropertiesofDconsDet}}
	    Let $X \in \Sch^{\qcqs}$, suppose that there exists $d\in \bbN$ such that, for every qcqs \'etale map $U\to X$, the (derived global sections) functor $\Gamma(U,-)$ has cohomological dimension bounded by $d$. 
	    Then the natural maps 
	    \[\fShv(X)\to \calD(X_\et)\to \calD_\et(X), \]
	    are equivalences.
	    In particular, since the \'etale site is invariant under universal homeomorphisms \cite[Tag~03SI]{StaProj}, for every $X$ of finite expansion over $k$ we have isomorphisms
	    \[\fShv(X)\simeq \calD(X_\et)\simeq \calD_\et(X). \]
    \end{theorem}

	\subsubsection{Extending the algebraic 3-functor formalisms}
	The category of prestacks on $\PSch^\qcqs$ and on $\PSchf$ clearly disagree. 
	As we said in \S \ref{ss: perfectPrestacks}, we find it convenient to treat every geometric object as an object of $\PreStk$ due to the presence of many different sheaf theories with slightly different descent properties.
	Since the 3-functor formalisms that we consider are sheafy for the \'etale topology and in particular the Zariski topology (\Cref{thm-Hansen-Scholze}, \Cref{thm: schematicDetisavsheaf} and \Cref{Shv-is-sheafy-etale}), there is no harm in restricting them to $\PSchf$.

More precisely, consider $\PSch^\qcqs$ as a site endowed with the \'etale topology. 
By applying \Cref{cons: ExtensiontoRepresentableMaps} and \Cref{path-topos}, we get natural transformations of 3-functor formalisms on $(\SchStk_\et, E^\rep_\pfp)$, where $E^\rep_\pfp$ denotes maps that are relatively representable in schemes perfectly finitely presented and separated 
\begin{equation}{\label{eqn: naturalTransformofSixFunctors1}}
\calD_\cons\Rightarrow \fShv \Rightarrow \calD_\et.
\end{equation}
We let $E^\rep_{\pfp,\PreStk}$ denote $(-)^{-1}_\et(E^\rep_\pfp)$ the preimage of $E^\rep_\pfp$ under \'etale sheafification.

We can restrict the 3-functor formalisms of \Cref{lemma: BasicPropertiesofDconsDet} to $(\PreStk, E^\rep_{\pfp,\PreStk})$ along the  morphism of geometric setups 
\[(\PreStk, E^\rep_{\pfp,\PreStk}) \ra (\SchStk_{\et},E^{\rep}_{\pfp}) \]
that \'etale sheafification induces. 
From now on, if the context is clear, we will abuse the notation and use $E^\rep_\pfp$ to denote what would be more correctly to denote by $E^\rep_{\pfp,\PreStk}$.

\begin{remark}
	\label{sheafification-observ}
	Although in our setup $\calD_\cons$ (respectively $\fShv$ and $\calD_\et$) takes input values a general prestack $X\in \PreStk$, 
	the functor naturally factors through arc-sheafification $\PreStk\to \SchStk_{\on{arc}}$ (respectively \'etale sheafification and v-sheafification), using \Cref{thm: Dconsarcsheaf} (respectively \Cref{Shv-is-sheafy-etale} and \Cref{thm: schematicDetisavsheaf}).
	Indeed, this follows from \cite[Proposition 1.3.1.7]{SAG} (see also the proof of \Cref{prop: descentproperties}). 
\end{remark}
\begin{remark}{\label{rem: KanextendingNaturalTransforms}}
Consider the natural transformations 
    \[ \calD_{\cons,0} \Rightarrow \fShv_{0} \Rightarrow \calD_{\et,0} \]
    obtained by restricting \eqref{eqn: naturalTransformofSixFunctors1} along $\PreStk^{\op} \hookrightarrow \Corr(\PreStk,E_{\pfp}^{\rep})$.
   By \Cref{lemma: rightKanextensionDoesWhatisExpected}, an alternative way of constructing it is by taking the right Kan-extension along the embedding $(\PSch^{\aff})^{\op} \hookrightarrow \PreStk^{\op}$ of the similar sequence 
    \[ \calD_{\cons,0} \Rightarrow \fShv_{0} \Rightarrow \calD_{\et,0} \]
    of functors defined on $(\PSch^{\aff})^{\op}$.
\end{remark}

In \S \ref{subsec: DiamondFunctor}, we will work with the analytic version of $\calD_\et(-)$. 
To avoid confusion, we will give the 3-functor formalism just constructed on $\PreStk$ a new name. 
\begin{definition}{\label{defn: DLambdaSch}}
	For $X\in \PreStk$, we let $\calD_\Lambda^\sch(X):=\calD_\et(X)$, and we let
	\[\calD^\sch_\Lambda:\Corr(\PreStk, E^\rep_\pfp)\to \LinCat_\Lambda\]
	denote the 6-functor formalism just constructed.
\end{definition}

We make some straightforward observations.

\begin{proposition}
	For any $X\in \PreStk$, the map $\calD_\cons(X)\to \calD^\sch_\Lambda(X)$ is fully faithful. 	
\end{proposition}
\begin{proof}
This follows from \Cref{rem: KanextendingNaturalTransforms}, and the fact that limits of fully faithful maps are always fully faithful. 
\end{proof}

       \begin{definition}
	       \label{schematic standard shirek}
       We let $E^{\on{std}}_\sch$ denote the class of maps in $\PreStk$ whose \'etale sheafification belongs to the set of $(\fShv,\et,E)$-standard $!$-able maps in the sense of \Cref{define-standard-Shriek}.       	
       \end{definition}

We apply \Cref{path-topos-sheafy} and \Cref{suavecito-gives-fine} to get the following 6-functor formalism which will suffice for our purposes.\footnote{Strictly speaking, we apply \Cref{path-topos-sheafy} and \Cref{suavecito-gives-fine} to $(\PSch^\qcqs,E_\pfp)$ which gives rise to a 6-functor formalism on $(\SchStk_\et,E^{\on{std}}_{\fShv,\et,!})$.
But we can always pre-compose with the \'etale sheafification map $\PreStk\to \SchStk_\et$.}
        \begin{proposition}
		\label{important analytification prop}
		The following statement hold.
		\begin{enumerate}
			\item The 6-functor formalisms $\fShv$ and $\calD^\sch_\Lambda$ defined on $(\PreStk,E^\rep_\pfp)$ extend uniquely to an \'etale-sheafy 6-functor formalisms defined on $(\PreStk,E^{\on{std}}_\sch)$ with values in $\LinCat_\Lambda$.
			\item  We have a natural transformation of 6-functor formalisms  
		defined on $\Corr(\PreStk,E_\sch^\sstd)$ 
		\[\fShv\Rightarrow \calD^\sch_\Lambda.\]
	\item 	The morphism $\fShv(X)\to \calD^\sch_\Lambda(X)$ 
		is compatible with the natural embeddings $\calD_\cons(X)\subseteq \fShv(X)$ and $\calD_\cons(X) \subseteq \calD^\sch_\Lambda(X)$.
		\end{enumerate}
        \end{proposition}

	\begin{remark}
		For general $\fShv$-$!$-able maps, it might not be true that $f_!:\fShv(X)\to \fShv(Y)$ preserves constructibility. 
		For example, we could consider the compactly supported cohomology of $B\bb{G}^\perf_{m,k} := [\Spec(k)/\mathbb{G}^\perf_{m,k}]$, which has infinite amplitude. 
	\end{remark}

We introduce the following shorthand for the relevant notion of Artin stack of \Cref{define-Artin} in the $\fShv$ formalism.
 \begin{definition}{\label{defn: SchematicArtinPreStks}}
\label{defn: schematicfine}
\label{defi relatively Artin}
We say that $X\in \PreStk$ is Artin, if the \'etale sheafification of the structure map $f:X\to \ast$ is $(\fShv,E^\rep_\pfp,\et)$-Artin in the sense of \Cref{define-Artin}. 
We denote by $\PreStk^\Art\subseteq \PreStk$ the subcategory of Artin prestacks.
\end{definition}

\subsection{Analytic 6-functor formalisms}
\label{ss: analytic 6-ff go}

Before the language of 6-functor formalisms and correspondence categories became widespread, one would talk about 6-functor formalisms in a non-technical sense, by explicitly constructing the $6$-operations and showing various desired compatibilities at the level of triangulated categories. 
The main goal of \cite{Sch17} is to construct and elucidate a version of \Cref{prop basic 6-functors analytic} below, without the higher-categorical considerations. 
We follow Mann's higher-categorical incarnation of this work.
\subsubsection{Usual \'etale cohomology of diamonds.}
Recall that, if $X$ is a locally spatial diamond, then one can attach a well-behaved small \'etale site (see \cite[Definition~14.1]{Sch17}).
	As in the algebraic case, one defines $\calD_\et(X)$ to be the left-completion of $\calD(X_\et) := \calD(X_\et,\Lambda)$, where $X_\et$ denotes the topos of \'etale sheaves. 
	Moreover, by \cite[Proposition 14.15]{Sch17} the natural map of sites 
	\[X_v\to X_\et\]
	induces a fully faithful embedding
	\[\calD_\et(X)\subseteq \calD(X_v).\]

	As shown in \cite[\S 14]{Sch17} (or \cite[Propositions 3.16, 3.20]{Mann2022NuclearSheaves}), the functor $\calD_\et(-)$ satisfies v-descent, and in particular extends to a functor on $\AnStk_v$.

As in the previous section, we now explain how this functor extends to a 6-functor formalism on $\AnStk_v$. 
    Following \cite{Sch17} and \cite{GHW22}, Mann has defined a 6-functor formalism of \'etale sheaves on $\AnStk_{v}$ (or rather the $1$-truncated objects $\vStk \subset \AnStk_{v}$, but the construction is the same) with the so-called fdcs maps (which stands for finite dimension, compactifiable and spatial).
	\begin{definition}[{\cite[Definition~5.4]{Mann2022NuclearSheaves}}]{\label{defn: fdcs}}
		We say a map $[f: X \ra Y] \in \AnStk_{v}$ is \emph{fdcs} if the following hold. 
		\begin{enumerate}
			\item It is representable in locally spatial diamonds. 
			\item It has locally finite transcendence degree in the sense of \cite[Definition~21.7]{Sch17}.
			\item For every point $x\in |X|$, there is an open immersion $U\subseteq X$ containing $x$ such that $f{|_U}:U\to Y$ is compactifiable in the sense of \cite[Definition~22.2]{Sch17}. 
		\end{enumerate}
        We write $E_\fdcs$ for the set of fdcs maps in $\AnStk_{v}$. 
	\end{definition}
        
\begin{theorem}[{\cite[Proposition 5.6]{Mann2022NuclearSheaves}}]
		\label{prop basic 6-functors analytic}
		The following statements hold.
		\begin{enumerate}
			\item The pair $(\AnStk_{v},E_\fdcs)$ defines a geometric setup. 
			\item $\calD_\et(-)$ extends to a sheafy presentable 3-functor formalism on $\Corr(\AnStk_{v},E_\fdcs)$ with values in $\LinCat_\Lambda$ which we will denote by $\calD_{\Lambda}^{\an}$.
			\item If $f$ is \'etale, then $f_!\simeq f_\natural$ is the left-adjoint to $f^*$. 
			\item If $f$ is proper, then $f_!\simeq f_*$ is the right-adjoint to $f^*$.
		\end{enumerate}
	\end{theorem}

    As before, one can perform \Cref{shriek hull} or its standard incarnation \Cref{suavecito-gives-fine}. 
    We let $E^\sstd_\an$ denote the set of edges in $\AnStk_v$ that are $(\calD^\an_\Lambda,v, E_\fdcs)$-standard $!$-able in the sense of \Cref{define-standard-Shriek}. 

    \begin{theorem}[{\cite[Theorem 5.11]{Mann2022NuclearSheaves}}]
	    \label{applied the standard shriekable already}
	    We have morphisms of geometric setups
	    \[(\AnStk_v,E_\fdcs) \to (\AnStk_{v},E^\sstd_\an)\to (\AnStk_{v},E_{\calD^\an,!})\]
	    The presentable 3-functor formalism $\calD_{\Lambda}^{\an}$, originally defined for $\Corr(\AnStk_v,E_\fdcs)$, extends uniquely to presentable 3-functor formalisms on $\Corr(\AnStk_{v},E^\sstd_\an)$ and $\Corr(\AnStk_{v},E_{\calD^\an,!})$, which we also denote by $\calD_{\Lambda}^{\an}$.
    \end{theorem}

    \begin{remark}{\label{rem: perfectoidprestacks}}
    As before, we will occasionally precompose $\calD_{\Lambda}^{\an}$ along the natural map of geometric setups 
    \[ (\AnPreStk,E^\sstd_\an) \ra (\AnStk_v,E^\sstd_\an), \]
    where we abuse notation by continuing to denote by $E^\sstd_\an$ its preimage under the v-sheafification functor. 
    \end{remark}
    As in \Cref{defn: SchematicArtinPreStks}, we consider the notion of Artin stacks. 
     \begin{definition}{\label{defn: AnalyticArtinPreStks}}
We say that a v-stack $X\in \AnStk_{v}$ (resp. a perfectoid prestack $X \in \AnPreStk$) is an Artin $v$-stack, if the structure map $f:X\to \ast$ is $(\calD_{\Lambda}^{\an},E_\fdcs,v)$-Artin. 
We denote by $\AnStk_{v}^\Art\subseteq \AnStk_{v}$ (resp. $\AnPreStk^{\Art} \subset \AnPreStk$) the subcategory of Artin v-stacks. 
\end{definition}

    \subsubsection{A pro\'etale version.}
    \label{sec: uop topology}
For technical reasons, we will need a variant of \Cref{applied the standard shriekable already} for the pro\'etale topology, which will interact better with analytification.
We start with the following.

\begin{proposition}
	\label{uop properties}
	The following statements hold.
	\begin{enumerate}
		\item If $X$ is a spatial diamond, then there is a surjective and universally open cover of the form $\Spa(R,R^+)\to X$ with $\Spa(R,R^+)$ a strictly totally disconnected space.
		\item If $X\to Y$ is quasiseparated map of v-sheaves, $Y$ and $\tilde{Y}$ are locally spatial diamonds and $\tilde{Y}\to Y$ is a v-surjective quasi-pro\'etale map such that $\tilde{X}:=X\times_Y \tilde{Y}$ is a locally spatial diamond, then $X$ is a locally spatial diamond. 
		\item If $X$ is a qs v-sheaf and $f:Y\to X$ is a qs universally open quasi-pro\'etale surjective map with $Y$ a locally spatial diamond, then $X$ is locally spatial diamond.
	\end{enumerate}
\end{proposition}
\begin{proof}
The first claim is \cite[Proposition 11.24]{Sch17}.	
The second claim follows from \cite[Proposition 13.4.(iv)]{Sch17}.
Indeed, any v-surjective quasi-pro\'etale map of locally spatial diamonds is surjective as a map of pro\'etale stacks.

To show the third claim, first take $U\subseteq Y$ with $U$ a spatial diamond. 
Then the map $U\to f(U)$ is qcqs universally open quasi-pro\'etale and surjective, and replacing the role of $X$ and $Y$ by $f(U)$ and $U$ we may reduce to the case in which $X$ is qcqs and $Y$ is a spatial diamond.
In this setup, one readily verifies that $X$ is a spatial v-sheaf as in \cite[Definition 12.12]{Sch17}.
Indeed, for any $U\subseteq Y$ quasicompact $f(U)\subseteq X$ is an open subset that is quasicompact over $X$, and since $f$ is an open map, the family $\{f(U)\}_{U\subseteq Y}$ as we range over the quasicompact open subsets of $Y$ is a basis of the topology for $X$.
One can then appeal to \cite[Theorem 12.18]{Sch17}, to show that $X$ is a spatial diamond.
\end{proof}

The main reason that the pro\'etale topology will be useful for us is the following statement.
It explains that to verify if a map is $E_\fdcss$ it suffices to do it pro\'etale locally.
Recall that $\AnStk_\uop$ denotes the category of small pro\'etale sheaves on $\Perf^\aff$ with values in anima (\ref{eqn: SubcategoriesofAnalyticStacks}).

\begin{proposition}
	\label{pro-etale locality}
Let $f:X\to Y$ be a map in $\AnStk_\uop$.
Suppose that $Y$, $\tilde{Y}\in \on{sLocSptl}$, and that there is a v-surjective quasi-pro\'etale map $\tilde{Y}\to Y$ such that the base change $\tilde{f}:X\times_Y \tilde{Y}\to \tilde{Y}$ is in $E_{\fdcss}$.
Then $X\in \on{sLocSptl}$ and $f\in E_{\fdcss}$.
\end{proposition}
\begin{proof}
	We note that $X$ is already a v-sheaf, indeed one can argue as in \Cref{enough-resilients}.(2) below, by noting that surjective quasi-pro\'etale maps are surjective as pro\'etale stacks and recalling that all diamonds are $v$-sheaves.
We want to show that $X\to Y$ is in $E_\fdcss$. 
	Being separated and locally compactifiable can be tested v-locally on the target (see \cite[Proposition 10.11]{Sch17}, \cite[Proposition 22.3.(iii)]{Sch17}).
In particular, we know that $X \ra Y$ is separated (since $\tilde{f}$ is), and it follows by \Cref{uop properties} that $X\in \on{sLocSptl}$.
Being representable in locally spatial diamonds holds for any morphism of locally spatial diamonds (see \cite[Proposition 13.4.(ii)]{Sch17}).
	Finally, we must show that it is  of locally finite transcendence degree.
By definition (see \cite[Definition 21.7]{Sch17}), one has to take a supremum of the transcendence degree $C(x)/C(y)$ along the set of commutative diagrams 
	\begin{center}
	\begin{tikzcd}
		\Spa(C(x),C(x)^+) \arrow{r} \arrow{d}  & \Spa(C(y),C(y)^+) \arrow{d} \\
	 X \arrow{r} & Y,
	\end{tikzcd}
	\end{center}
	where the two vertical maps are quasi-pro\'etale morphisms from geometric points.
	Since geometric points split every pro\'etale cover, any such commutative diagram lifts to $\tilde{f}:X\times_Y \tilde{Y}\to \tilde{Y}$, and since $\tilde{f}$ is locally of finite transcendence degree the map $f$ also is. 
\end{proof}

We will now have to dig into the precise construction of the $6$-functor formalism from \Cref{prop basic 6-functors analytic}. 
Recall that the proof of {\cite[Proposition 5.6]{Mann2022NuclearSheaves}} goes through first considering an auxiliary geometric setup 
\[(\calC,E)=(\on{sLocSptl},E_{\on{fdcss}}).\]
Here $\on{sLocSptl}$ is the category of locally spatial diamonds that are separated over $\ast$, and $E_{\on{fdcss}}$ is the set of morphisms that are fdcs in the sense of Definition \ref{defn: fdcs}.
Note that, by \cite[Proposition 13.4]{Sch17} and the cancellative property of separatedness, any map in $(\on{sLocSptl},E_{\on{fdcss}})$ is automatically representable in locally spatial diamonds and separated. 
This justifies fdcss as an acronym for finite dimensional, compactifiable, spatial and separated.

The proof of {\cite[Proposition 5.6]{Mann2022NuclearSheaves}} considers an auxiliary 6-functor formalism 
\[\calD_{\on{sLS}}:\Corr(\on{sLocSptl},E_{\on{fdcss}})\to \LinCat_\Lambda.\]
Using that locally spatial diamonds are $v$-sheaves \cite[Proposition~11.9]{Sch17}, one can endow $\on{sLocSptl}$ with the v-topology and the induced topos is still $\AnStk_v$. 
Moreover, $\calD_{\on{sLS}}$ is sheafy for the v-topology.
In this way, Mann applies \Cref{cons: ExtensiontoRepresentableMaps} to obtain a 6-functor formalism  
\[\calD_s:\Corr(\AnStk_v,E^{\on{rep}}_{\on{fdcss}})\to \LinCat_\Lambda.\]
Finally, Mann applies \Cref{cons: ExtensiontoFineMaps} to extend from $E^{\on{rep}}_{\on{fdcss}}$ to $E_{\on{fdcs}}$.
It follows from revisiting this proof that we have an equality
\[E^{\on{std}}_\an=E^{\on{std}}_{\calD_{\on{sLS}},v,!},\] where the right-hand side is as defined in \Cref{define-standard-Shriek}. 

We now move on to discuss the pro\'etale version of this theory.
We have a sequence of full subcategories
\[\on{Std}\subseteq \Perf^\aff\subseteq  \on{sLocSptl}\subseteq \AnStk_v\subseteq \AnStk_\uop \subseteq \AnPreStk.\]

We can endow $\on{sLocSptl}$ with the Grothendieck topology $\tau_\uop$ in which $X\to Y$ is surjective if and only if it is a surjection of pro\'etale stacks.

Using \Cref{uop properties}, we see that strictly totally disconnected perfectoid spaces are a basis for the topology $\tau_\uop$ on $\on{sLocSptl}$.
In particular, the category of small prestacks on $\on{sLocSpt}$ with values in anima that are sheaves with respect $\tau_\uop$ agrees with $\AnStk_\uop$.
Moreover, since 
\[\calD_{\on{sLS}}:\Corr(\on{sLocSptl},E_{\on{fdcss}})\to \LinCat_\Lambda,\]
is sheafy for the v-topology, it is also sheafy for the pro\'etale topology on $\on{sLocSptl}$.
This allow us to apply \Cref{shriek hull} and \Cref{suavecito-gives-fine} to obtain 6-functor formalisms
\[\calD^\uop_\Lambda:(\AnStk_\uop, E_{\calD_{\on{sLS}},!})\to \LinCat_\Lambda\]
and 
\[\calD^\uop_\Lambda:(\AnStk_\uop, E^{\on{std}}_{\calD_{\on{sLS}},\uop,!})\to \LinCat_\Lambda, \]
where the former is an extension of the latter (see \Cref{rem: standard!ableis!able}).
\begin{remark}{\label{rem: WHatArtinandFineMeansintheAnalyticContext}}
We note that each of the 6-functor formalisms $\calD^{\uop}_{\Lambda}$ and $\calD_{\Lambda}^{\an}$, will have their own intrinsic notion of Artin stacks as in \Cref{define-Artin}.
 We refer to notion of Artin stacks attached to the $\calD^{\uop}_{\Lambda}$ formalism as Artin pro\'etale-stacks and to the notion of Artin stacks attached to the $\calD_{\Lambda}^{\an}$ formalism as Artin $v$-stacks, as in \Cref{defn: AnalyticArtinPreStks}. We emphasize that in the $\calD^{\uop}_{\Lambda}$ case, the atlas is required to be surjective for the pro\'etale topology.
\end{remark}

\begin{remark}
	Similarly to \Cref{rem: WHatArtinandFineMeansintheAnalyticContext}, the 6-functor formalisms $\calD^\uop_\Lambda$ and $\calD^\an_\Lambda$ have their own intrinsic notion of standard $!$-covers (see \Cref{a standard cover}) and consequently of standard $!$-ability (see \Cref{define-standard-Shriek}). 
\end{remark}

For notational convenience, we make the following definition.
\begin{definition}
We let $E^\sstd_\uop$ denote the set of edges in $\AnStk_\uop$, with $E^\sstd_\uop:=E^{\on{std}}_{\calD_{\on{sLS}},\uop,!}$. 	
Similarly, we let $E_{\calD_\Lambda^\uop,!}:=E_{\calD_{\on{sLS}},!}$.
\end{definition}
As we have mentioned, for technical reasons we will naturally be lead to work with $\calD^\uop_\Lambda$. 
Nevertheless, we care mostly about $\calD^\an_\Lambda$. 
The following two statements shows that these two constructions can be compared.
\begin{proposition}
	\label{everything looks normal after rest}
	The following statement holds.
	\begin{enumerate}
		\item The natural inclusions $\on{sLocSptl}\subseteq \AnStk_v\subseteq \AnStk_\uop$ induce a commutative diagram of geometric setups 
			\begin{center}
			\begin{tikzcd}
				(\on{sLocSptl},E_\fdcss) \arrow{r} \arrow{d}  & (\AnStk_v,E^\sstd_\an) \arrow{d} \\
				(\AnStk_\uop,E^\sstd_\uop)
 \arrow{r} & (\AnStk_\uop,E_{\calD^\uop_\Lambda,!}).
			\end{tikzcd}
			\end{center}

		\item We have a unique equivalences of 6-functor formalisms 
			\[(\calD^\uop_\Lambda)_{\mid_{\Corr(\AnStk_v,E^{\on{std}}_\an)}}\simeq \calD^\an_\Lambda,\]
			extending the equivalence of functors 
			\[(\calD^\uop_{\Lambda})_{\mid_{\Corr(\on{sLocSptl},E_{\on{fdcss}})}}\simeq (\calD^\an_{\Lambda})_{\mid_{\Corr(\on{sLocSptl},E_{\on{fdcss}})}}.\] 
	\end{enumerate}
\end{proposition}
\begin{proof}
	Since $\calD_{\on{sLS}}$ is sheafy for the v-topology, $\calD^\uop_{\Lambda,0}$ factors through the v-sheafification map as
	\[\AnStk_\uop^{\op} \to \AnStk_v^{\op} \xrightarrow{\calD^\an_{\Lambda,0}}\LinCat_\Lambda.\]
	In particular, we obtain a unique equivalence 
	\begin{equation}{\label{eqn: AnStackProEtStackEquation}}
    (\calD^\uop_{\Lambda,0})_{\mid_{\AnStk_v^\op}}\simeq \calD^\an_{\Lambda,0}.
    \end{equation}
	By construction, the classes of edges $E_{\calD^\uop_\Lambda,!}$ and $E^{\on{std}}_\an$ were obtained by appealing to \Cref{shriek hull} and \Cref{suavecito-gives-fine} respectively.
	In particular, we must show that
	$(E_{\calD^\uop_\Lambda,!})_{\mid_{\AnStk_v}}$ satisfies conditions $(i)-(v)$ of \Cref{suavecito-gives-fine} with respect to the v-topology, in order to show the containment $E^{\on{std}}_\an \subseteq (E_{\calD^\uop_\Lambda,!})_{\mid_{\AnStk_v}}$. 
	Verifying conditions (i), (ii) and (v) is straightforward, and condition (iv) follows from (iii) using \Cref{iv-is-automatic}. 
	To show (iii), fix a map $f:Y\to Z\in \AnStk_v$.
	Recall that a standard $!$-cover will have the form $g:X\to Y\in \AnStk_v$ with $g=g_1\circ \dots g_n$ where each $g_i\in (E_{\calD^\uop_\Lambda,!})_{\mid_{\AnStk_v}}$ and it is either a prim or a suave v-cover. 
	Condition (iii) posits that if $f\circ g\in E_{\calD^\uop_\Lambda,!}$, then $f\in E_{\calD^\uop_\Lambda,!}$ must also hold.
	Our claim is that if $g$ is a standard $!$-cover for the v-topology, then it is a universal $(\calD^\uop_\Lambda)^!$-cover (see \Cref{defn: universal*and!covers}). 
	If this claim holds true, by $!$-locality of $E_{\calD^\uop_\Lambda,!}$ we will conclude that $f\in E_{\calD^\uop_\Lambda,!}$, as required.
	To see that $g$ is a universal $(\calD^\uop_\Lambda)^!$-cover we may use induction to reduce to the case $g=g_1$, so that $g$ is either a prim or a suave v-cover. 
	We note that the notion of prim v-cover agrees with the notion of prim pro\'etale cover (see \Cref{defi:suaveprim-cover} and note that the $*$-pushforward functors agree by passing to adjoints in (\ref{eqn: AnStackProEtStackEquation})), in this case the statement follows from \Cref{prim standard !-able}. 
	If $g$ is a suave v-cover, then $g:X\to Y$ is $\calD^\uop_\Lambda$-suave map of pro\'etale stacks that is v-surjective (but not necessarily pro\'etale surjective). 
	Since $\calD^\uop_\Lambda$ satisfies descent for the v-topology, then $g^*$ is conservative, and we may use \Cref{suave-implies-fine} to conclude that $g$ is a universal $(\calD^\uop_\Lambda)^!$-cover (but not necessarily a pro\'etale standard $!$-cover).

	 Point (2) follows from the uniqueness part of the extension procedure \Cref{suavecito-gives-fine}.
\end{proof}

Let $(\AnStk_v,E^\sstd_\uop)$ denote the geometric setup obtained from restricting $(\AnStk_\uop,E^\sstd_\uop)$ along the inclusion map $\AnStk_v\to \AnStk_\proet$.

\begin{proposition}
	\label{more-maps-of-geometric-setups}
	We have a maps of geometric setups 
	\[(\on{sLocSptl},E_\fdcss)\to (\AnStk_v,E^\sstd_\uop)\to (\AnStk_v,E^\sstd_\an),\]
	and there is a unique equivalence
	\[(\calD^\uop_\Lambda)_{\mid_{\corr(\AnStk_v,E^\sstd_\uop)}}\simeq (\calD^\an_\Lambda)_{\mid_{\corr(\AnStk_v,E^\sstd_\uop)}}\]
			extending the equivalence of functors 
			\[(\calD^\uop_{\Lambda})_{\mid_{\Corr(\on{sLocSptl},E_{\on{fdcss}})}}\simeq (\calD^\an_{\Lambda})_{\mid_{\Corr(\on{sLocSptl},E_{\on{fdcss}})}}.\] 
\end{proposition}
\begin{proof}
	The key point is that any standard $!$-cover for the pro\'etale topology is automatically a standard $!$-cover for the v-topology. 
	At a technical level, the statement follows from applying \Cref{suavecito-gives-fine take number 666} with $F=E^\sstd_\uop$. 
	One verifies explicitly that for such a choice of $F$, $E^{\on{std}}_{\calD_{\on{sLS}},v,F,!}=E^\sstd_\uop$. 
\end{proof}

To summarize the situation we combine \Cref{everything looks normal after rest} and \Cref{more-maps-of-geometric-setups} to obtain a commutative diagram of geometric setup, on which we have defined a 6-functor formalism that uniquely extends $\calD_{\on{sLS}}$.
\begin{center}
\begin{tikzcd}
	(\on{sLocSptl},E_\fdcss) \ar{rd} 
		\ar[bend right=20]{rdd}
		\ar[bend left=20]{rrd}
	& &  \\
			 &(\AnStk_v,E^\sstd_\uop)	 \arrow{r} \arrow{d}  & (\AnStk_v,E^\sstd_\an) \arrow{d} \\
			 &	(\AnStk_\uop,E^\sstd_\uop) \arrow{r} & (\AnStk_\uop,E_{\calD^\uop_\Lambda,!}) 
\end{tikzcd}
\end{center}

    We now turn our attention to studying how the algebraic and analytic 6-functor formalisms interact under the $\diamondsuit$-functor introduced in \S \ref{sec: the Diamond Functor}.

    \subsection{The $\diamondsuit$-functor and 6-functors}{\label{subsec: DiamondFunctor}}
    Recall (\cite[\S~27]{Sch17}) that, for a perfect scheme $X\in \PSch$, there are natural maps of sites $c_{X}: (X^{\dia})_{v} \rar (X_\proet)\to (X_\et)$, which induce a commutative diagram of fully faithful functors 
    \begin{center}
    \begin{tikzcd}
	    \calD_\et(X) \arrow{r}{c^*_X} \arrow{d}  & \calD_\et(X^\diamondsuit) \arrow{d} \\
	    \calD(X_\proet)\arrow{r}{c^*_X} & \calD(X^\diamondsuit_v).
    \end{tikzcd}
    \end{center}

    By \cite[Proposition~27.1]{Sch17}, this induces a natural transformation
    \[\calD^\sch_\Lambda(-)\Rightarrow \calD^\an_\Lambda(-^\diamondsuit)\]
    of functors of the form
    \[\calD_0:\PSch^\op\to \CAlg(\LinCat_\Lambda).\]

    Ideally, one would want to upgrade this to a natural transformation of 6-functor formalisms on a geomeric setup of the form $(\PreStk, E)$, where $E$ is a set of edges containing the representable ones and such that $E^\diamondsuit\subseteq E^\sstd_\an$. 

    The first observation is that 
    \[(\PreStk,\diamondsuit^{-1}(E^\sstd_\an))\]
    is a geometric setup, since $\diamondsuit$ is compatible with finite limits and composition.
    Precomposition gives a 6-functor formalism 
    \begin{align}
	    \Corr(\PreStk,\diamondsuit^{-1}(E^\sstd_\an))&\to \LinCat_\Lambda\\
     X&\mapsto \calD^\an_\Lambda(X^\diamondsuit)\\	
     [X\leftarrow Z \rightarrow Y] &\mapsto \calD^\an_\Lambda([X^\diamondsuit\leftarrow Z^\diamondsuit \rightarrow Y^\diamondsuit]),
    \end{align}
    which we denote by $\calD^\diamondsuit_\Lambda$.

    \begin{proposition}
	    \label{Comparing-underlying-cats}
	    The restriction of $\calD^\diamondsuit_\Lambda$ along \[\PreStk^\op\hookrightarrow \Corr(\PreStk,\diamondsuit^{-1}(E^\an_\fine))\] is the right Kan extension of 
	    \[\calD^\an_{\Lambda,0}(-^\diamondsuit):\PSch^\op\to \CAlg(\LinCat_\Lambda)\] 
	    along the inclusion $\PSch^\op \hookrightarrow \PreStk^\op$.
    \end{proposition}
    \begin{proof}
	    Recall (see the discussion after \Cref{defn: diamondfunctor}) that the functor $\diamondsuit$ is a left Kan-extension of the functor 
        \[ \diamondsuit: \PSchf \ra \AnStk_v \subset \AnPreStk \]
        \[ \Spec(A) \mapsto \Spec(A)^{\diamondsuit}, \]
	along $\PSchf \subset \PreStk$, and therefore it commutes with colimits in $\AnPreStk$ by \cite[Proposition~6.2.1.9]{SAG}.
Furthermore, the restriction of $\calD^\an_\Lambda$ may, in light of \Cref{prop basic 6-functors analytic}.(2), be viewed as a right Kan extension of $\Detale(-)$ along the embedding $(\Perff)^\op \hookrightarrow \AnStk^{\op}$.
This similarly gives that $\calD^\diamondsuit_\Lambda$ commutes with limits on $\PreStk^\op$.
	    These facts (see \cite[Lemma 5.1.5.5]{HTT}) identify $\calD^\diamondsuit_\Lambda$ as the right Kan extension along  
	    \[\CAlg^\perf\hookrightarrow \PreStk^\op\]
	    of the functor 
	    \[A\mapsto \calD^\an_\Lambda(\Spec A^\diamondsuit).\]
	    By \cite[Proposition 4.3.2.8]{HTT}, it suffices to show that 
	    \[X\mapsto \calD^\an_\Lambda(X^\diamondsuit)\]
	    is also the right Kan extension along the inclusion 
	    \[\CAlg^\perf\to \PSch^\op.\]

	    From the definition of Kan extensions (see \cite[Definition 4.3.2.2]{HTT}) we must show that for all $X\in \PSch$
\[\calD^\an_\Lambda(X^\diamondsuit) \simeq \varprojlim_{\Spec A\to X} \calD^\an_\Lambda(\Spec A^\diamondsuit)  \]
	    is an equivalence. 
	    Since $\calD^\an_\Lambda$ is v-sheafy and commutes with limits in $\AnStk_v^\op$, it suffices to show that 
	    \[\varinjlim_{\Spec A\to X} \Spec A^\diamondsuit\to X^\diamondsuit \]
	   is an equivalence in the category of v-sheaves. 
	    But this follows from the fact that the $\diamondsuit$-functor takes Zariski covers to v-covers. 
	    Indeed, as will be explained in \Cref{prop: hcoversgivediamondcovers} below, something stronger holds.
    \end{proof}

    Using \Cref{Comparing-underlying-cats}, one can already construct a natural transformation 
    \[\calD^\sch_{\Lambda,0}\Rightarrow \calD^\diamondsuit_{\Lambda,0}\]
    as functors from $\PreStk^\op$ with values in $\CAlg(\LinCat_\Lambda)$.
    To construct 
    $\calD^\sch_\Lambda\Rightarrow \calD^\diamondsuit_\Lambda$
    as a natural transformation of 6-functor formalisms, we will construct an intermediary 6-functor formalism in the spirit of \Cref{section-natural-transforms}.

	\begin{lemma}{\label{lemma: pfpanalytifiestofdcs}}
		Let $[f: X \ra Y ]\in \PSch$ be a separated map of perfect schemes and let $f^{\diamondsuit}: X^{\diamondsuit} \ra Y^{\diamondsuit}$ be the induced map in $\vShv$.
		\begin{enumerate}
			\item If $f$ is a closed immersion, then $f^\diamondsuit$ is a closed immersion. \label{lemma: pfpanalytifiestofdcs-item:closed}
			\item $f^\diamondsuit$ is separated and partially proper. \label{lemma: pfpanalytifiestofdcs-item:separated}
             \item If $f$ is \'etale (resp. an open immersion), then $f^{\diamondsuit}$ is \'etale (resp. an open immersion).
Moreover, the left adjoints to $f^*$ and $f^{\diamondsuit,*}$ (which we write as $f_{\natural}$ and $f^{\diamondsuit}_{\natural}$) satisfy $c_{Y}^{*}f_{\natural} \simeq f^{\diamondsuit}_{\natural}c_{X}^{*}$.  \label{lemma: pfpanalytifiestofdcs-item:etale}
	     \item  If $f\in E_\pfp$, then $f^{\diamondsuit}\in E^\rep_{\on{fdcss}}$. \label{lemma: pfpanalytifiestofdcs-item:fdcs}
			\item If $f$ is perfectly proper, then $f^{\diamondsuit}$ is proper and the right adjoints to $f^*$ and $f^{\diamondsuit,*}$ (which we write as $f_{*}$ and $f^{\diamondsuit}_{*}$) satisfy $c_{Y}^{*}f_{*} \simeq f^{\diamondsuit}_{*}c_{X}^{*}$.  \label{lemma: pfpanalytifiestofdcs-item:proper}
		\end{enumerate}
	\end{lemma}
	\begin{proof}
		Note that property \ref{lemma: pfpanalytifiestofdcs-item:closed} readily implies the separatedness part of property \ref{lemma: pfpanalytifiestofdcs-item:separated}.
		That $f^\diamondsuit:X^\diamondsuit\to Y^\diamondsuit$ is partially proper follows directly from the definition.
		Verifying that a map of $v$-sheaves is a closed immersion (resp. \'etale, resp. an open immersion) may be done $v$-locally on the target by \cite[Proposition~10.11]{Sch17}.
    In particular, it suffices to show that these properties hold for $X^\diamondsuit\times_{Y^\diamondsuit}\Spa(R,R^+)\to \Spa(R,R^+)$ as we range over a basis for the v-topology on $(\AnStk_v)_{/Y^\diamondsuit}$.
    Therefore, we may restrict our attention to those maps $\Spa(R,R^+)\to Y^\diamondsuit$ with $\Spa(R,R^+)$ a product of points as in \cite[Definition 2.14]{GIZ25}.
    We note that, for a product of points, $Y^\diamondsuit(R,R^+)=Y^{\diamondsuit_\pre}(R,R^+)$, and in particular we have a factorization 
    \[\Spa(R,R^+)\to (\Spec R)^\diamondsuit\to Y^\diamondsuit.\]
    Without loss of generality, we may assume $Y=\Spec R$.
    To prove property \ref{lemma: pfpanalytifiestofdcs-item:closed}, we observe that if $X\to \Spec R$ is a closed immersion corresponding to an ideal $I\subseteq R$, then $X^\diamondsuit\times_{\Spec(R)^\diamondsuit}\Spa(R,R^+)\to \Spa(R,R^+)$ is the Zariski closed immersion (as in \cite[Definition 5.7]{Sch17}) corresponding to the same ideal.

    Property \ref{lemma: pfpanalytifiestofdcs-item:etale} follows from observing that, since $R$ is a comb (see \cite[Proposition 2.19]{GIZ25}), every \'etale map over $S\to \Spec R$ admits a decomposition of the form $S=\bigcup_{i\in \calI}S_i\to \Spec R$ where $S_i=\Spec R[\frac{1}{r_i}]$ for some $r_i\in R$.  
    Since arbitrary unions of separated \'etale maps remain separated \'etale, it suffices to show that $X^\diamondsuit\to Y^\diamondsuit$ is an open immersion in the case $Y = \Spec A$, and $X=\Spec A[\frac{1}{a}]$ for some $a\in A$. 
    It follows from the definitions that, for any map $\Spa(R,R^+)\to (\Spec A)^\diamondsuit$ from a perfectoid space, the base change $X^\diamondsuit\times_{Y^\diamondsuit} \Spa(R,R^+)\subseteq \Spa(R,R^+)$ is the open subsheaf corresponding to the subset given by $\{x\in \Spa(R,R^+)\mid |f^{*}(a)|_x\neq 0\}$, which implies the desired claim.  
    The isomorphism on the second part of property \ref{lemma: pfpanalytifiestofdcs-item:etale} can be shown by passing to geometric points $\Spa(C,C^+)\to Y^\diamondsuit$, as pulling back along such geometric points defines a conservative family.
This reduces the comparison to the case in which $Y=\Spec C$ an algebraically closed field and $X=\coprod_{i\in \calI} \Spec C$, but here it is an explicit calculation (see also \cite[Proposition~27.4]{Sch17}).
        
		Let us prove property \ref{lemma: pfpanalytifiestofdcs-item:fdcs}, we have already shown that the map is separated by property \ref{lemma: pfpanalytifiestofdcs-item:separated}. 
        By \cite[Lemma~5.5 (1)]{Mann2022NuclearSheaves}, the property of being fdcs can be checked analytic locally on both source and target. 
	Therefore, we may assume that $X$ and $Y$ are affine.
In this case, we can factor $f$ as $X\xrightarrow{g} \bbA^{n,\perf}_Y\xrightarrow{\pi_Y} Y$ with $g$ a closed immersion.
		We note that $g^\diamondsuit$ is also a closed immersion from which it follows that it is fdcs. 
		That $\pi_{Y}$ is fdcs is evident since its base change by any affinoid perfectoid $\Spa(R,R^+)\to Y^\diamondsuit$ is representable by the perfectoid affine space $\bbA^n_R$ over $\Spa(R,R^+)$.

		Let us prove property \ref{lemma: pfpanalytifiestofdcs-item:proper}. 
	We have seen that $X^\diamondsuit\to Y^\diamondsuit$ is separated and partially proper, so it suffices to show that the map is quasicompact. 
		Quasicompactness is v-local by \cite[Proposition 10.11]{Sch17}, so we may assume $Y$ is affine.
        Since the map $X\to Y$ is pfp we may use Noetherian approximation to write $X\to Y$ as the base change of a pfp map $X_0\to Y_0$ with a pfp structure map $Y_0\to \ast$. 
	By Chow's lemma \cite[Tag 0200, 0201]{StaProj}, we may find a perfectly projective cover $Z_0\to X_0$ such that $Z_0\to Y_0$ is also a perfectly projective map.
        In particular, we are reduced to showing that the maps $Z_0^\diamondsuit\to X_0^\diamondsuit$ and $Z_0^\diamondsuit\to Y_0^\diamondsuit$ are quasicompact.
Indeed, in this case $Z_0^\diamondsuit \to X_0^\diamondsuit$ is automatically v-surjective (see \Cref{prop: hcoversgivediamondcovers} below), which also shows that $X_0^\diamondsuit\to Y_0^\diamondsuit$ is quasicompact. 
        Finally, dealing with perfectly projective maps splits into closed immersions and the $\bbP_k^{n,\mathrm{perf}}\to \ast$ case: the former is property \ref{lemma: pfpanalytifiestofdcs-item:closed} and the latter follows from \Cref{lm:proper_gives_cartesian_diam} and \Cref{case-of-thepoint} which we show below.
	The commutativity of pushforward with $c^*$ is the content of \cite[Proposition 27.4]{Sch17}.
	\end{proof}
    
	Consider the posetal category $\calP=\{\{\an\} \rightarrow \{\sch\}\}$ containing two objects $\{\sch\}$ and $\{\an\}$, and consider the geometric setup $(\PSch^\calP,E_\pfp^\calP)$. 
	Since $c^*_X$ is defined through a map of sites, we obtain a functor 
	\begin{align}
		\calD_\Lambda^{c^*}:(\PSch^\calP)^\op &\to \CAlg(\LinCat_\Lambda)\\
		\calD_\Lambda^{c^*}((X,\an))&:= \calD^\diamondsuit_\Lambda(X) \\
		\calD_\Lambda^{c^*}((X,\sch))&:= \calD^\sch_\Lambda(X) 
	\end{align}
	satisfying $\calD^{c^*}_\Lambda\circ \ins_{\an} \simeq \calD^\diamondsuit_\Lambda$ and $\calD^{c^*}_\Lambda\circ \ins_{\sch}\simeq \calD^\sch_\Lambda$.
	As in \Cref{remark-factoring-morphisms}, the functor $\calD^{c^*}_\Lambda$ is simply encoding the natural transformation $\calD^\sch_\Lambda\Rightarrow \calD^\diamondsuit_\Lambda$. 

	Recall that we have a suitable decomposition $P_\pfp^\calP\subseteq E_\pfp^\calP$ and $I_\pfp^\calP\subseteq E_\pfp^\calP$.
	Using \Cref{lemma: pfpanalytifiestofdcs}, we apply \Cref{construction: 3 functors} in the form of \Cref{lemma: HowtoGetaSixfunctorformalsimontheinterval} to obtain a 6-functor formalism
	\[\calD^{c^*}_\Lambda:(\PSch^\calP,E^\calP_\pfp)\to \LinCat_\Lambda.\]

\begin{proposition}
\label{comparing-intermediary-step}
The map $\diamondsuit:\PSch\to \vShv$ induces a morphism of geometric setups
\[(\PSch,E_\pfp)\to (\vShv,E^\rep_{\on{fdcss}}).\]
	Moreover, we have natural equivalences 
	\begin{enumerate}
		\item $\calD^{c^*}_\Lambda \circ \ins_\an \simeq \calD^\diamondsuit_\Lambda$
		\item $\calD^{c^*}_\Lambda \circ \ins_\sch \simeq \calD^\sch_\Lambda$
	\end{enumerate}
	of 6-functor formalisms on $(\PSch,E_\pfp)$.
In particular, we obtain a morphism of 6-functor formalisms
	\[\calD^\sch_\Lambda\Rightarrow \calD^\diamondsuit_\Lambda\]
	on $(\PSch,E_\pfp)$.
\end{proposition}
\begin{proof}
	The first statement was shown in \Cref{lemma: pfpanalytifiestofdcs}.  
Let us recall how these 6-functor formalisms got constructed. 
Observe that 
\begin{enumerate}
	\item $\calD^\sch_\Lambda\simeq \calL\calZ(\PSch,E_\pfp,P_\pfp,I,\calD^\sch_{0,\Lambda})$ by construction,
	\item $\calD^{c^*}_\Lambda\simeq \calL\calZ(\PSch^\calP,E^\calP_\pfp,P^\calP_\pfp,I^\calP,\calD^{c^*}_{0,\Lambda})$ also by construction,
	\item but $\calD^\diamondsuit_\Lambda$ was not obtained directly by the Liu--Zheng construction. 
\end{enumerate}
Using \Cref{cor: LiuZhengRestriction}, we see that 
\begin{align*}
	\calD^{c^*}_\Lambda\circ \ins_\sch &\simeq \calL\calZ(\PSch^\calP,E^\calP_\pfp,P^\calP_\pfp,I^\calP,\calD^{c^*}_{0,\Lambda}) \circ \ins_\sch \\
			    &\simeq \calL\calZ(\PSch,E_\pfp,P_\pfp,I,\calD^{c^*}_{0,\Lambda}\circ \ins_\sch)\\
			    &\simeq \calL\calZ(\PSch,E_\pfp,P_\pfp,I,\calD^{\sch}_{0,\Lambda})\\
			    &\simeq \calD^\sch_\Lambda 
\end{align*}
 Similarly, $\calD^{c^*}_\Lambda \circ \ins_\an \simeq \calL\calZ(\PSch,E_\pfp,P_\pfp,I,\calD^{\diamondsuit}_{0,\Lambda})$, and we wish to show that
\[\calD^\diamondsuit_\Lambda \simeq \calL\calZ(\PSch,E_\pfp,P_\pfp,I,\calD^{\diamondsuit}_{0,\Lambda}).\]
It is not at all obvious to show this by hand, however, the characterization of the Liu--Zheng construction given in \Cref{thm: LiuZhengUniquenesss} shows this. 
Indeed, the result follows by combining \Cref{thm: LiuZhengUniquenesss} with the following key facts that follow from \Cref{prop basic 6-functors analytic} and \Cref{lemma: pfpanalytifiestofdcs}. 
\begin{enumerate}
	\item Both theories take $P_\pfp$ to cohomologically proper maps.
	\item Both theories take $I$ to cohomologically \'etale maps. 
	\item By definition $\calD^{c^*}_{\Lambda,0} \circ \ins_\an \simeq \calD^\diamondsuit_{\Lambda,0}$. \qedhere
\end{enumerate}
\end{proof}

Using \Cref{prop: descentproperties}.(5) and \Cref{thm: schematicDetisavsheaf}, we apply \Cref{cons: ExtensiontoRepresentableMaps} in the form of \Cref{path-topos} to obtain an extended 6-functor formalism
\[\calD^{c^*}_\Lambda:(\SchStk_\et^\calP,(E^\rep_\pfp)^\calP)\to \LinCat_\Lambda.\]
As usual, one can run the descent machinery (\Cref{shriek hull}) to obtain a 6-functor formalism on the $!$-able hull of $\calD_\Lambda^{c^*}$ 
\[\calD^{c^*}_\Lambda:\corr(\SchStk_\et^\calP,(E^\calP_{\pfp})_{\calD_\Lambda^{c^*},!})\to \LinCat_\Lambda. \]
One is naturally lead to the following question.

\begin{question}
	\label{Dc* vs D diamond}
How do $\calD_\Lambda^{c^*}\circ \ins_{\on{an}}$ and $\calD^\diamondsuit_\Lambda$ compare as $6$-functor formalisms? 
\end{question}

As we have warned the reader throughout \S \ref{section-natural-transforms}, a key subtlety is that it is not automatic that $\diamondsuit(E_{\calD^\sch_\Lambda,!})\subseteq E_{\calD^\an_\Lambda,!}$. 
The intersection 
\[E'=E_{\calD^\sch_\Lambda,!}\cap \diamondsuit^{-1}(E_{\calD^\an_\Lambda,!}),\] 
makes $(\SchStk_\et,E')$ into a geometric setup on which both $\calD_\Lambda^{c^*}\circ \ins_{\on{an}}$ and $\calD^\diamondsuit_\Lambda$ are defined. 
Nevertheless, since $E'$ was not directly constructed by extension procedures, it is not very clear if the morphism of geometric setups 
\[(\PSch,E_\pfp)\to (\SchStk_\et,E')\] 
enjoys a ``uniqueness of the extension'' property, which makes it very subtle to understand to what extent $\calD_\Lambda^{c^*}\circ \ins_{\on{an}}$ and $\calD^\diamondsuit_\Lambda$ agree.
This will force us to restrict the $E'$ that we consider in order to ensure uniqueness of the extension. 
This also leads to the observation that $\diamondsuit$ might not be what we should consider for our purposes.
Actually, it is not even clear to us if $\diamondsuit(E^\rep_\pfp)\subseteq E_{\calD^\an_\Lambda,!}$ holds or not. 
For this reason, it becomes technically more convenient to make the following two changes.
\begin{enumerate}
	\item Instead of working with all of $E_{\calD^\sch_\Lambda,!}$, we work instead with the standard $!$-able maps $E^\sstd_\sch$ (as in \Cref{define-standard-Shriek} and \Cref{schematic standard shirek}). 
		Actually, as we will see, we work with the variant from \Cref{suavecito-gives-fine take number 666}.
	\item Instead of working with $\AnStk_v$ we work with $\AnStk_\uop$ (see \S \ref{sec: uop topology}).
\end{enumerate}

We consider the pro\'etale version of $\diamondsuit$,
\[\diamondsuit_\uop:\PreStk\to \AnStk_\uop.\]

It may be constructed by the formula $(-)_\uop\circ \diamondsuit_\pre$, where $(-)_\uop$ denotes sheafification for the pro\'etale topology. 
We may follow the construction of $\calD_\Lambda^{c^*}$ and $\calD_{\Lambda}^{\Diamond}$ , replacing the v-topology with the pro\'etale topology to obtain a 6-functor formalisms

\[\calD^{c^*_\uop}_\Lambda:\corr(\SchStk_\et^\calP,(E^\rep_\pfp)^\calP)\to \LinCat_\Lambda. \]
and 
\[\calD^{\diamondsuit_\uop}_\Lambda:\corr(\SchStk_\et,E^\rep_\pfp) \to \LinCat_\Lambda. \]

We have the following extension\footnote{Recall that the analytic sheafification of $\diamondsuit_{\pre}$ applied to perfect schemes is already a $v$ and hence a proetale sheaf, as in the discussion after \Cref{defn: diamondfunctor}.} of \Cref{comparing-intermediary-step} to stacks.

\begin{proposition}
\label{comparing-intermediary-step-number 2}
The map $\diamondsuit_\uop:\SchStk_\et\to \AnStk_\uop$ induces a morphism of geometric setups
\[(\SchStk_\et,E^\rep_\pfp)\to (\AnStk_\uop,E^\rep_{\on{fdcss}}).\]
	Moreover, we have natural and unique equivalences 
	\begin{enumerate}
		\item $\calD^{c_\uop^*}_\Lambda \circ \ins_\an \simeq \calD^{\diamondsuit_\uop}_\Lambda$
		\item $\calD^{c_\uop^*}_\Lambda \circ \ins_\sch \simeq \calD^\sch_\Lambda$
	\end{enumerate}
	of 6-functor formalisms on $(\SchStk_\et,E^\rep_\pfp)$, extending the 6-functor formalism from \Cref{comparing-intermediary-step}.
\end{proposition}
\begin{proof}
    We first note that one has a map of geometric setups
    \[(\PSch,E_\pfp)\to (\AnStk_\uop,E^{\on{rep}}_{\on{fdcss}}).\]
    that factors through the map $(\vShv,E^{\on{rep}}_{\on{fdcss}}) \ra (\AnStk_\uop,E^{\on{rep}}_{\on{fdcss}})$ of geometric setups, by combining  \Cref{comparing-intermediary-step}  with the observation that the analytic sheafification of $X^{\diamondsuit_{\pre}}$ is already a $v$-sheaf, as in the discussion after \Cref{defn: diamondfunctor}. 
    In particular, \Cref{comparing-intermediary-step} gives us the claim on $(\PSch,E_{\pfp})$.

	The second step is to show that $\diamondsuit_\uop(E^\rep_\pfp)\subseteq E^\rep_{\on{fdcss}}$.
	This part of the argument requires us to work with $\AnStk_\uop$ instead of $\AnStk_v$.  
	Fix $X,Y\in \SchStk_\et$ and a map $f:X\to Y$ in $E^\rep_\pfp$.
	Let $T\in \on{sLocSptl}$ be a separated locally spatial diamond with a map $T\to Y^{\diamondsuit_\uop}$ and let 
	\[S:=X^{\diamondsuit_\uop}\times_{Y^{\diamondsuit_\uop}} T.\]
	By construction $S$ is a pro\'etale stack, and we wish to show that it is a separated locally spatial diamond for which the map $S\to T$ is fdcss.
	Since $X^{\diamondsuit_\uop}\to Y^{\diamondsuit_\uop}$ is partially proper by construction, (local) compactifiability and separatedness of $S$ over $T$ is automatic.
	The condition on the finiteness of transcendence degrees may be tested on geometric points of $T$.
 Moreover, every $\Spa(C,C^+)$ point of $Y^{\diamondsuit_\proet}$ factors as
	\[\Spa(C,C^+)\to (\Spec C)^{\diamondsuit_\uop}\to Y^{\diamondsuit_\uop}, \]
	since every pro\'etale cover of $ \Spa(C,C^+)$ admits a section.
	In particular, the condition on finite transcendence degree reduces to the case where $f$ is a pfp map of perfect schemes which was explained above.

	We now show that $S$ is a locally spatial diamond, which, using \cite[Proposition~13.4 (ii)]{Sch17}, will tell us that the map $S \ra T$ is representable in locally spatial diamonds, as desired.

	This easily reduces to the case in which $T$ is itself spatial, and by \Cref{uop properties} we may even assume that $T=\Spa(R,R^+)$ for a strictly totally disconnected space. 
	By definition of $\diamondsuit_\uop$, there is a pro\'etale cover $T'\to T$ such that the induced map $T'\to Y^{\diamondsuit_\uop}$ factors as 
	\[T'\to (\Spec A)^{\diamondsuit_\uop}\to Y^{\diamondsuit_\uop}.\]
	By \Cref{uop properties}, we may show instead that $S'=S\times_T T'$ is locally spatial. 
	If we let $Z=X\times_Y \Spec A$, then 
	\[S'=Z^{\diamondsuit_\uop}\times_{(\Spec A)^{\diamondsuit_\uop}} T'.\]
	That $S'$ is a locally spatial diamond now follows from the pfp case explained above.
This finishes the second step, and finishes showing the first claim.

	The second claim can be easily reduced to a pfp map of perfect schemes using the uniqueness part of \Cref{cons: ExtensiontoFineMapsI}.   
	\end{proof}

	We can run descent again (\Cref{shriek hull}) to obtain a 6-functor formalism on the $!$-able hull of $\calD_\Lambda^{c_\uop^*}$ 
\[\calD^{c_\uop^*}_\Lambda:\corr(\SchStk_\et^\calP,(E^\calP_{\pfp})_{\calD_\Lambda^{c_\uop^*},!})\to \LinCat_\Lambda, \]
but we will still run into a similar question to \Cref{Dc* vs D diamond}.

It works better to use \Cref{path-topos-sheafy} to obtain a $6$-functor formalism
\[\calD^{c_\uop^*}_\Lambda:\Corr(\SchStk_\et^\calP,(E^\sstd_\sch)^\calP)\to \LinCat_\Lambda.\]
but even if we do this, it is still not clear if $\diamondsuit_\uop(E^\sstd_\sch)\subseteq E^\sstd_{\uop}$.
To resolve this, we define a class of $!$-\textit{analytifiable maps.}

\begin{definition}
	\label{!-analitifyiable}
	Let $F:=\diamondsuit^{-1}_\uop(E^\sstd_\uop)$ denote the set of maps in $\SchStk_\et$.
	We let $E^\an_\sch :=E^\sstd_{\calD_{\Lambda}^{\mathrm{sch}}, \et,F,!}$ denote the set of edges obtained from \Cref{suavecito-gives-fine take number 666}.
	If $f\in E^\an_\sch$ we call it $!$-analytifiable.
\end{definition}

By construction, $E^\rep_\pfp \subseteq E^\an_\sch\subseteq E^\sstd_\sch$, $\diamondsuit_\uop(E^\an_\sch)\subseteq E^\sstd_{\uop}$, and $(\SchStk_\et,E^\an_\sch)$ is a geometric setup.  
Crucially, we have the following statement.

\begin{proposition}\label{cstar-extension}
    There is a unique equivalence \[\calD^{c_\uop^*}_\Lambda\circ \ins_\an\overset{\alpha}{\simeq} \calD^{\diamondsuit_\uop}_\Lambda,\]
			of 6-functor formalisms on $\corr(\SchStk_\et,E_\sch^\an)$, extending the equivalence of \Cref{comparing-intermediary-step-number 2}.   
\end{proposition}
\begin{proof}
	The uniqueness of the extension is part of the statement in \Cref{suavecito-gives-fine take number 666}.
\end{proof}

\begin{remark}
	\label{fine contains analytifiable}
	Although $E^\an_\sch$ is fairly abstractly defined, it contains all the standard fine maps of \Cref{lets define fine finally}, and consequently all the Artin maps of \Cref{define-Artin}. 
	For most practical purposes, it suffices to consider these.
\end{remark}

As usual, we may precompose along the map $\PreStk\to \SchStk_\et$ and redefine $E^\an_\sch$ as the pullback of the previously defined version.
Let us summarize our findings.
The technical version of \Cref{thm: DiamondAnalytificationUltimate} is \Cref{cstar-extension}.  

	\begin{theorem}{\label{thm: DiamondAnalytificationUltimate}}
		Fix $X\in \PreStk$ and a map $[f:X_1\to X_2]\in \PreStk$.
		Then we have constructed an analytification functor 
		\[c^*_X:\calD^\sch_\Lambda(X)\to \calD_\Lambda^\uop(X^{\diamondsuit_\uop})\simeq \calD^\an_\Lambda(X^\diamondsuit),\]
		and a class of edges $E^\an_\sch$ with $E_\pfp\subseteq E_{\on{fine}}\subseteq E^\an_\sch\subseteq E^\sstd_\sch$ that satisfy the following. 
\begin{enumerate}
	\item $c^*_X:\calD^\sch_\Lambda(X)\to \calD^\uop_\Lambda(X^{\diamondsuit_\uop})\simeq \calD_\Lambda^\an(X^\diamondsuit)$ is a fully faithful morphism in $\LinCat_\Lambda$.
	\item $c^*_{X_1}\circ f^*\simeq f^{\diamondsuit_\uop,*}\circ c^*_{X_2}$, satisfying higher coherences.
	\item $ f_*\circ c_{X_1,*}\simeq c_{X_2,*} \circ f_*^{\diamondsuit_\uop} $, satisfying higher coherences.
	\item If $f\in E^\an_\sch$, then $f$ and $f^{\diamondsuit_\uop}$ are $!$-able and 
		\[ f^{\diamondsuit_\uop}_!\circ c^*_{X_1}\simeq c^*_{X_2} \circ f_!\] satisfying higher coherences.
	\item The functor $\diamondsuit_\uop$ together with the family of functors $\{c_Y^*\}_{Y\in (\PreStk)_{/X}}$ organize into a morphism of the $2$-categories of kernels 
		\[ \mathcal{K}_{X,\calD_{\Lambda}^{\mathrm{sch}}} \ra \mathcal{K}_{X^{\diamondsuit_\uop},\calD_{\Lambda}^{\uop}}.\] 

		\item If $f\in E^\an_\sch$, then $c^*_{X_1}(\Suave_{f,\calD^{\mathrm{sch}}_{\Lambda}}(X_1))\subseteq \Suave_{f^{\diamondsuit_\uop},\mathcal{D}^{\uop}_{\Lambda}}(X_1^{\Diamond_{\proet}})$. 
			Moreover, for $A\in \Suave_{f,\calD^{\mathrm{sch}}_{\Lambda}}(X_1)$ we have a canonical identification  
			\[c^*_{X_1}\bbD_{X_1/X_2,\calD^\sch_\Lambda}(A)\simeq \bbD_{X^{\diamondsuit_\uop}_1/X^{\diamondsuit_\uop}_2,\calD^\uop_\Lambda}c^*_{X_1}A.\]
		\item If $f\in E^\an_\sch$, and $f$ is $\calD_{\Lambda}^{\mathrm{sch}}$-suave (resp. $\calD_{\Lambda}^{\mathrm{sch}}$-smooth, resp. $\calD_{\Lambda}^{\mathrm{sch}}$-unipotent), then $f^{\diamondsuit_\uop}$ is $\calD^{\uop}_{\Lambda}$-suave (resp. $\calD^{\uop}_{\Lambda}$-smooth, resp. $\calD^{\uop}_{\Lambda}$-unipotent).
\end{enumerate}
	\end{theorem}
    \begin{proof}
    Everything follows from the above discussion, where we note that we implicitly used \Cref{cor: DeltathreefunctorsgivesTransformonKernels} for points (5)-(7).
    \end{proof}

An annoying feature of our setup is that we have deviated from the usual practices in the literature of working with $\AnStk_v$ in order to define analytification.
This motivates the following definition.
\begin{definition}
	\label{def: resiliant stacks}
	We say that $X\in \SchStk_\et$ is \textit{resilient} if the natural map $X^{\diamondsuit_\uop}\to X^\diamondsuit$ is an equivalence.
	Alternatively, $\diamondsuit_\uop^{-1}(\AnStk_v)\subseteq \SchStk_\et$ is the full subcategory of resilient \'etale stacks.
	We denote this category by $\resSchStk_{\et}$.
\end{definition}

\begin{proposition}
	\label{technical-theorem}
	There is a morphism of 6-functor formalisms
	\[\calD^\sch_\Lambda\Rightarrow \calD^\diamondsuit_\Lambda\]
	on $\corr(\resSchStk_\et,E^\an_\sch)$, extending the one from \Cref{comparing-intermediary-step}.
\end{proposition}
\begin{proof}
	In \Cref{cstar-extension}, we have explained that we have a morphism of 6-functor formalisms
	\[\calD^\sch_\Lambda\Rightarrow \calD^{\diamondsuit_\uop}_\Lambda\]
on $\corr(\SchStk_\et,E^\an_\sch)$, that uniquely extends the one from \Cref{comparing-intermediary-step}.
We may restrict this along the inclusion $(\resSchStk_\et,E^\an_\sch)\hookrightarrow (\SchStk_\et,E^\an_\sch)$ to obtain a morphism of 6-functor formalisms
	\[\calD^\sch_\Lambda\Rightarrow \calD^{\diamondsuit_\uop}_\Lambda\]
on $\corr(\resSchStk_\et,E^\an_\sch)$.
It then suffices to construct an equivalence of 6-functor formalisms 
\begin{equation}
	\label{eq:inmidle-of-proof-loreleiloreleiloreiiah}
	\calD^{\diamondsuit_\uop}_\Lambda\simeq \calD^{\diamondsuit}_\Lambda
\end{equation}
on $\corr(\resSchStk_\et,E^\an_\sch)$.
As we explained in \Cref{more-maps-of-geometric-setups}, we have an equivalence
\begin{equation}
	\label{eq:inmidle-of-proof-loreleiloreleiloreiiah2}
\calD^\uop_\Lambda\simeq \calD^\an_\Lambda
\end{equation}
of 6-functor formalisms on $(\AnStk_v,E^\sstd_\uop)$.   
But, by \Cref{def: resiliant stacks} and \Cref{!-analitifyiable}, we have a map of geometric setups
\[\diamondsuit:(\resSchStk,E^\sch_\an)\to (\AnStk_v,E^\sstd_\uop).\]
Moreover, both $\calD^{\diamondsuit_\uop}_\Lambda$ and $\calD^{\diamondsuit}_\Lambda$ are respectively obtained from $\calD^\uop_\Lambda$ and $\calD^\an_\Lambda$ by precomposing with this map of geometric setups.
In particular, the equivalence from (\ref{eq:inmidle-of-proof-loreleiloreleiloreiiah2}) can be used to construct the equivalence of (\ref{eq:inmidle-of-proof-loreleiloreleiloreiiah}).
\end{proof}

We can restate \Cref{thm: DiamondAnalytificationUltimate} in terms of $\diamondsuit$, as long as we restrict to resilient stacks.

	\begin{theorem}{\label{thm: DiamondAnalytificationUltimate2}}
		Fix $X\in \resSchStk_\et$ and a map $[f:X_1\to X_2]\in \resSchStk_\et$.
		Then we have constructed an analytification functor 
\[c^*_X:\calD^\sch_\Lambda(X)\to \calD_\Lambda^\an(X^\diamondsuit),\]
satisfying the following.
\begin{enumerate}
	\item $c^*_X:\calD^\sch_\Lambda(X)\to \calD_\Lambda^\an(X^\diamondsuit)$ is a fully faithful morphism in $\LinCat_\Lambda$.
	\item $c^*_{X_1}\circ f^*\simeq f^{\diamondsuit,*}\circ c^*_{X_2}$, satisfying higher coherences.
	\item $ f_*\circ c_{X_1,*}\simeq c_{X_2,*} \circ f_*^{\diamondsuit} $, satisfying higher coherences.
	\item If $f\in E^\an_\sch$, then $f$ and $f^\diamondsuit$ are $!$-able and 
		\[ f^\diamondsuit_!\circ c^*_{X_1}\simeq c^*_{X_2} \circ f_!\] satisfying higher coherences.
	\item The functor $\diamondsuit$ together with the family of functors $\{c_Y^*\}_{Y\in (\resSchStk_\et)_{/X}}$ organize into a morphism of the $2$-categories of kernels 
		\[ \mathcal{K}_{X,\calD_{\Lambda}^{\mathrm{sch}}} \ra \mathcal{K}_{X^\diamondsuit,\calD_{\Lambda}^{\an}}.\] 

		\item If $f\in E^\an_\sch$, then $c^*_{X_1}(\Suave_{f,\calD^{\mathrm{sch}}_{\Lambda}}(X_1))\subseteq \Suave_{f^{\diamondsuit},\mathcal{D}^{\an}_{\Lambda}}(X_1^{\Diamond})$. 
			Moreover, for $A\in \Suave_{f,\calD^{\mathrm{sch}}_{\Lambda}}(X_1)$ we have a canonical identification  
			\[c^*_{X_1}\bbD_{X_1/X_2,\calD^\sch_\Lambda}(A)\simeq \bbD_{X^\diamondsuit_1/X^\diamondsuit_2,\calD^\an_\Lambda}c^*_{X_1}A.\]
		\item If $f\in E^\an_\sch$, and $f$ is $\calD_{\Lambda}^{\mathrm{sch}}$-suave (resp. $\calD_{\Lambda}^{\mathrm{sch}}$-smooth, resp. $\calD_{\Lambda}^{\mathrm{sch}}$-unipotent), then $f^{\diamondsuit}$ is $\calD^{\an}_{\Lambda}$-suave (resp. $\calD^{\an}_{\Lambda}$-smooth, resp. $\calD^{\an}_{\Lambda}$-unipotent).
\end{enumerate}
	\end{theorem}
    \begin{proof}
	    This is a reformulation of \Cref{technical-theorem}, where we have used \Cref{cor: DeltathreefunctorsgivesTransformonKernels} for points (5)-(7).
    \end{proof}

In what follows, we show a lemma stating that many \'etale stacks that one cares about are resilient.
These examples will suffice for our purposes.
\begin{lemma}
	\label{enough-resilients}
The following statements hold.
\begin{enumerate}
	\item If $X\in \PSch$, then $X$ is resilient.
	\item If $[f:Y\to X] \in \SchStk_\et$, is a map such that $X\in \resSchStk_\et$ and for every strictly totally disconnected $\Spa(R,R^+)$ the base change $X^{\diamondsuit_\uop}\times_{Y^{\diamondsuit_\uop}} \Spa(R,R^+)$ is a v-sheaf of anima, then $Y\in \resSchStk_\et$.
\item If $[f:Y\to X] \in \SchStk_\et$, is a map such that $X\in \resSchStk_\et$ and for every strictly totally disconnected $\Spa(R,R^+)$ the base change $X^{\diamondsuit_\uop}\times_{Y^{\diamondsuit_\uop}} \Spa(R,R^+)$ is a diamond, then $Y\in \resSchStk_\et$.
	\item If $[f:Y\to X] \in \SchStk_\et$ such that $X\in \resSchStk_\et$ and such that for every $\Spec A \in \PSch^\aff$ the base change $Y\times_X \Spec A\in \resSchStk_\et$ is resilient, then $Y$ is resilient.
	\item If $Y\to X$ is in $E^\rep_\pfp$ and $X$ is resilient, then $Y$ is resilient.
	\item For any linear algebraic group $G$ over $k$, $\bbB_\et G$ is resilient.
    \item For $K$ a profinite group, we have that $\bbB_{\mathrm{fpqc}}\underline{K}$ is resilient.
\end{enumerate}
\end{lemma}
\begin{proof}
	Point (1) follows from the fact that, for a perfect scheme $X$, the sheafification of $X^{\diamondsuit_{\pre}}$ with respect to the analytic topology on $\Spa(R,R^+)$ is already a $v$-sheaf, as explained in the discussion proceeding \Cref{defn: diamondfunctor}. 
	In particular, we have that $X^{\diamondsuit_{\uop}} = X^{\diamondsuit}$. 

	For point (2), we contemplate the map $Y^{\diamondsuit_\uop}\to Y^{\diamondsuit}$ and we wish to show it is an equivalence. 
	We may regard both as pro\'etale sheaves. 
	In particular, it suffices to show that $Y^{\diamondsuit_\uop}(R,R^+)\to Y^{\diamondsuit}(R,R^+)$ on a basis for the pro\'etale topology.
	By \Cref{uop properties}, the strictly totally disconnected spaces form a basis for the pro\'etale topology.
	Note that since $X$ is resilient, we have a commutative diagram of pro\'etale sheaves
	\begin{center}
	\begin{tikzcd}
		Y^{\diamondsuit_\uop} \arrow{dr} \ar{r}  & Y^\diamondsuit \arrow{d} \\
		 & X^{\diamondsuit_\uop}
	\end{tikzcd}
	\end{center}
	We may instead show that
	\[Y^{\diamondsuit_\uop}\times_{X^\diamondsuit_\uop} \Spa(R,R^+)\to Y^\diamondsuit\times_{X^{\diamondsuit}} \Spa(R,R^+)\] is an equivalence on $\Spa(R,R^+)$-valued points.
	But by assumption, $Y^{\diamondsuit_\uop}\times_{X^\diamondsuit_\uop} \Spa(R,R^+)$ is a v-sheaf of anima, so the claim follows. 
	Point (3) follows from Point (2) and \cite[Proposition~11.9]{Sch17}.
	
	For point (4), we want to verify that the hypothesis of point (2) hold.
	In other words, 
	We want to show that for any strictly totally disconnected space and a map $\Spa(R,R^+)\to X^\diamondsuit$, $Y^{\diamondsuit_\uop}\times_{X^\diamondsuit}\Spa(R,R^+)$ is a v-sheaf. 
	The claim is clear if the map $\Spa(R,R^+)\to X^\diamondsuit$ factors through a map $f^\diamondsuit:(\Spec A)^\diamondsuit\to X^\diamondsuit$.
	But this always holds pro\'etale-locally.
	More precisely, let $\Spa(R_0,R_0^+)\to \Spa(R,R^+)$ be a pro\'etale cover for which the induced map $\Spa(R_0,R_0^+)$ factors through some affine $\Spec A^{\diamondsuit} \ra X^{\diamondsuit}$, as above.
	Let $\Spa(R_\bullet, R^+_\bullet)$ denote the \v{C}ech nerve.
	We consider the following commutative diagram in $\AnStk_\uop$ with Cartesian squares
	\begin{center}
	\begin{tikzcd}
		Y^{\diamondsuit_\uop}\times_{X^\diamondsuit} \Spa(R_\bullet,R^+_\bullet)	\arrow{r}{\theta_\bullet} \arrow{d}  & Y^\diamondsuit\times_{X^\diamondsuit} \Spa(R_\bullet,R^+_\bullet)\arrow{r} \arrow{d}  & \arrow{d} \Spa(R_\bullet,R^+_\bullet) \\
	 Y^{\diamondsuit_\uop} \arrow{r} & Y^{\diamondsuit} \arrow{r} & X^\diamondsuit 
	\end{tikzcd}
	\end{center}
	From the hypothesis it follows that for all $n\in \Delta$ 
	\[Y^{\diamondsuit_\uop}\times_{X^\diamondsuit} \Spa(R_n,R^+_n)	\xrightarrow{\theta_n}  Y^\diamondsuit\times_{X^\diamondsuit}\Spa(R_n,R^+_n)\]
	is an equivalence. 
	Since $\AnStk_\uop$ is a topos, colimits are universal and we can conclude that
	\[\colim_{n\in \Delta^\op}\theta_n:Y^{\diamondsuit_\uop}\times_{X^\diamondsuit} \Spa(R,R^+)	\rightarrow Y^\diamondsuit\times_{X^\diamondsuit}\Spa(R,R^+)\]
	is also an equivalence as we wanted to show.

	Point (5) follows from Point (4) and Point (1).	

	For Point (6), we choose a closed embedding  $G \ra \GL_{n}$, which induces a representable map
\[ \bbB_{\et}G \ra \bbB_{\et}\GL_{n}, \]
and, by point (2) this allows us to assume that $G = \GL_{n}$. 
To check that $(\bbB_{\et}\GL_{n})^{\diamondsuit_\uop}$ satisfies $v$-descent it suffices to do so on a basis for the pro\'etale topology on $\Perff$.  
In particular, we may do this on affinoid perfectoid spaces, where it reduces to $v$-descent for vector bundles on perfectoid spaces (\cite[Lemma~17.1.8, Theorem~5.2.8]{SW20}).

For point (7), we first observe that $(\bbB_{\rm fpqc}\ul{K})^{\diamondsuit_\uop} = [\ast/\ul{K}]_{\uop}$. 
Indeed, recall that $(\bbB_{\rm fpqc}\ul{K})\simeq (\bbB_{\profet}\ul{K})$ where $\profet$ denotes the pro-finite \'etale topology (see \cite[Remark 10.113]{Zhu25}).
In particular, their pro\'etale sheafification will also agree.
Now, for any strictly totally disconnected $\Spa(R,R^+)$, and any $\calT\in (\bbB_{\profet}\ul{K})(\Spec R)$, $\calT$ is isomorphic to the trivial torsor. 
Indeed, it suffices to find a section to $\calT\to \Spec R$.
This section may be constructed as in the proof of \cite[Proposition 2.20]{GIZ25}. 
Indeed, for all $x\in \pi_0(\Spec R)$ the corresponding connected closed subscheme is of the form $\Spec R_x$ with $R_x$ a strictly henselian local ring by \cite[Lemma 2.4.17]{kedlaya_liu_relative_p_adic_hodge_theory_foundations}.
In particular, each $\Spec R_x$ splits every finite \'etale map, and by Noetherian approximation this splitting spreads to an open neighborhood of $x\in \pi_0(\Spec R)$. 
Finally, we must show that $[\ast/\ul{K}]_\uop$ is already a v-sheaf, or equivalently that, for every $\Spa(R,R^+)\in \Perf^\aff$, the map of groupoids 
\[[\ast/\ul{K}]_\uop(R,R^+)\to [\ast/\ul{K}]_v(R,R^+)\] is an equivalence.
The map is clearly fully-faithful and the essential image consists of the $\ul{K}$-torsors for the v-topology that are pro\'etale locally trivial.
Let $\calT\in [\ast/\ul{K}]_v(R,R^+)$, we observe that the map $\calT\to \Spa(R,R^+)$ is a pro\'etale map that is also surjective.
Indeed, since $\Spa(R,R^{+})$ is strictly totally disconnected, pro\'etale and quasi-pro\'etale are the same (see \cite[Definition 10.1]{Sch17}), and quasi-pro\'etaleness of a map is v-local on the target (see \cite[Proposition 10.11.(v)]{Sch17}).
Since $\calT\to \Spa(R,R^+)$ is pro\'etale, and the pullback of any torsor to itself trivializes the torsor, we see that $\calT\in [\ast/\ul{K}]_\uop(R,R^+)$, as we wanted to show.
\end{proof}

    We finish this section with a statement that explains why the pro\'etale topology and the 6-functor formalism $\calD^\uop_\Lambda$ will be a technically convenient tool even if our main interest is studying $\calD^\an_\Lambda$.

\begin{proposition}
	\label{precise consequences of smoothness}
	Let $X_\uop, Y_\uop, Z_\uop\in \AnStk_\uop$ with maps $f_\uop:X_\uop\to Y_\uop$ and $g:Z_\uop\to Y_\uop$ such that $g_\uop\in E^\sstd_{\uop}$ is $\calD^\uop_\Lambda$-suave. 
	Let $X, Y, Z$ and $f, g$ denote the respective v-sheafification.
Consider a Cartesian diagram 
\begin{center}
\begin{tikzcd}
	W	\arrow{r}{f'} \arrow{d}{g'}  & Z \arrow{d}{g} \\
 X \arrow{r}{f} & Y
\end{tikzcd}
\end{center}
Then, for every $A\in \calD^\an_\Lambda(X)$, the natural map is an equivalence
\[g^*f_* A\xrightarrow{\simeq} f'_*g'^*A.\]
\end{proposition}
\begin{proof}
	Consider the pro\'etale version
\begin{center}
\begin{tikzcd}
	W_\uop	\arrow{r}{f'_\uop} \arrow{d}{g'}  & Z_\uop \arrow{d}{g_\uop} \\
 X_\uop \arrow{r}{f_\uop} & Y_\uop.
\end{tikzcd}
\end{center}
This diagram satisfies base-change against $*$-pushforward, by \cite[Lemma~4.5.13 (i)]{HeyerMann}.
But since the two formalisms agree for $*$-pullback (i.e., $\calD^\uop_{\Lambda,0}\simeq \calD^\an_{\Lambda,0}$) the base-change property also holds for $\calD^\an_\Lambda$. 
Indeed, the $*$-pushfoward can be recovered purely from $\calD^\uop_{\Lambda,0}\simeq \calD^\an_{\Lambda,0}$ by passing to right adjoints.
\end{proof}

	\section{Computing analytification}
    \label{sec: compute analytif}

    \noindent
    \textbf{Soft preamble to the section:}
    After a formal but lengthy dévisage, that $\pitch$ is an equivalence ultimately boils down to rather concrete computations.
    Some of this computations are closely related to the same ideas that placid geometry is founded on. 
    In a sense, this section is the technical backbone to justifying that one can can still do placid geometry in the analytic context after finding suitable hypothesis.  
    This section is technically intricate, but it is also a technical pillar for our argument.\\

    \noindent
    \textbf{Technical preamble to the section:}
    \begin{enumerate}
	    \item[\S \ref{ss: 6ff formalism computs}] Reaps the benefits of our work in \S \ref{s:Analyt for stacks section} and \S \ref{sec: analytificiation and 6-functor formalisms}.
		    Indeed, framing analytification as a 6-functor formalism gives plenty of commutativities for free, which we record in this section. 
	    \item[\S \ref{subsec: threeanalytificationfunctors}] Recalls the three analytification constructions $\diamond$, $\dagger$, $\diamondsuit$, and how they interact with one another under properness assumptions (see \Cref{lemma: geometricpropertiesofmapdaggertodiatodiamond,lm:proper_gives_cartesian_diam}). 
		    We also discuss how each of these functors induce a topology on the category of schemes, and relate these topologies with more classical ones.
		    For our purposes, the most relevant part is that all of these functors are compatible with the \'etale topology (see \Cref{prop: descentproperties}). 
    \item[\S \ref{Computing d}]  Explains why under finiteness assumptions $\diamond_\uop$ preserves suave maps. 
	    This does not follow from the 6-functor formalism, instead it is the observation that for a scheme $X$ perfectly finitely presented over $k$, the natural map $X^{\diamond}\subseteq X^\diamondsuit$ is an open immersion, and that suaveness still passes to the stacky setup.
	    This plays and important role when attempting to apply the techniques of \S \ref{ss: pro-unip base change}. 
    \item [\S \ref{overconv repl}] As it was observed in \cite{GIZ25}, the stack of analytic shtukas can be obtained from the schematic one via the $\dagger$-construction. 
	    The $\dagger$-construction does not preserve smooth or suave maps, which precludes the possibility of doing placid geometry directly on $\Sht^\an_{\calG}$.
	    In this subsection, we introduce a technique that we call overconvergent replacement, which captures at a technical level the observation that $X^\diamond$ and $X^\dagger$ only differ on higher rank points. 
	    In particular, computations with overconvergent sheaves on each of the setups should yield the same result.
	    The advantage is that by our work in \S \ref{Computing d}, the $X^\diamond$ version is more responsive to placid-geometrical techniques. 
    \item [\S \ref{ss: pro-unip base change}] Discusses pro-unipotent base change in the analytic setup. 
	    The main point is that, in the analytic setup, qcqs assumptions are more crucial than in the schematic one.
	    This section together with \S \ref{Computing d} and \S \ref{overconv repl} form the core techniques that play the role of placid geometry, but in the analytic setup.
    \end{enumerate}

	\subsection{6-functor formalism computations}
	\label{ss: 6ff formalism computs}
	Our first observation is that since analytification can be thought of as a 6-functor formalism, one can use formal properties of general 6-functor formalisms to perform some computations.
        \begin{proposition}{\label{lemma: constructibleissuave}}
		Suppose that $X\in \PSch^\pfp$, let $[f:X \ra \ast]$ denote the structure map, then any $A \in \calD_{\cons}(X)$ is $f$-suave.
	More precisely, the natural fully faithful maps 
	\[\calD_\cons(X)\xleftarrow{\simeq} \on{Suave}_{f,\calD_{\cons}}(X) \xrightarrow{\simeq} \on{Suave}_{f,\fShv}(X) \]
    are equivalences.
        \end{proposition}
        \begin{proof}
		The equivalence
		\[\calD_\cons(X)\xleftarrow{\simeq} \on{Suave}_{f,\calD_{\cons}}(X)\]
		is the $(i)\iff (ii)$ part of \cite[Theorem 4.4]{HS23}.
		The equivalence
		\[\on{Suave}_{f,\calD_{\cons}}(X) \xrightarrow{\simeq} \on{Suave}_{f,\fShv}(X)\]
		is the content of \cite[Proposition 3.4.(iii)]{HS23}.
		Strictly speaking, the reference works with finitely presented maps, but one can easily resolve this by choosing a finitely presented deperfection $X_0\to \Spec k$, and recalling that $X_\et\simeq (X_0)_\et$ canonically.
        \end{proof}
        We now want to extend this to a larger class of stacks. 
	Recall our notion of schematic Artin prestacks from \Cref{defn: SchematicArtinPreStks}.

	\begin{proposition}{\label{prop: allconstructiblesheavesSuaveunderembedding}}
		Let $X\in \PreStk^\Art$ be a schematic Artin pre-stack with structure map $f:X\to \ast$, and let $A\in \calD_\cons(X)$. 
		Then, the image of $A$ in $\fShv(X)$ (respectively in $\calD^\sch_\Lambda(X)$) is $f$-suave.
		Moreover, $\mathbb{D}_X(A)\in \calD_\cons(X)$.
	\end{proposition}
	\begin{proof}
		By definition of $\PreStk^\Art$ we may find an \'etale surjective and $\fShv$-suave map $[g:U\to X]\in E^\rep_\pfp$ such that $U\to \ast\in E_\pfp$.
		Since $\fShv$ is a sheaf for the \'etale topology by \Cref{Shv-is-sheafy-etale}, $g$ is a universal $\fShv^*$-cover, and by \cite[Lemma~4.5.8 (i)]{HeyerMann} it suffices to verify that $g^*A$ is $f\circ g$-suave.
		Since $g^*A\in \calD_\cons(X)$, by \Cref{lemma: constructibleissuave}, it is already $(\fShv,f\circ g)$-suave as we wanted to show.

		For the second part, $\mathbb{D}_X(A)$ is the suave dual of $A$ (see \Cref{prop: VerdierDualityandSuaveDuality}), so $\mathbb{D}_X(A)$ is also $f$-suave.
		With notation as above, $g^*\mathbb{D}_X(A)$ is $(f\circ g)$-suave, so by \Cref{lemma: constructibleissuave} $g^*\mathbb{D}_X(A)\in \calD_\cons(U)$.  
		Finally, since constructibility can be checked on \'etale covers (even arc-covers) by \Cref{thm: Dconsarcsheaf}, we conclude $\mathbb{D}_X(A)\in \calD_\cons(X)$.
	\end{proof}

	\begin{proposition}{\label{cor: lalggivescfine}}
		Let $X \in \PreStk^\Art$ with structure map $f:X\to \ast$.
The following is true. 
		\begin{enumerate}
			\item We have that $f\in E^\an_\sch$. 
			\item For $A \in \Dcons(X) \subset \calD_{\Lambda}^{\mathrm{sch}}(X)$, we have that $c_{X}^{*}(A)$ is $f^{\diamondsuit_\uop}$-suave.
		\end{enumerate}
	\end{proposition}
	\begin{proof}
		From the definitions it follows that $f\in E_\pfp^\fine$, and by \Cref{fine contains analytifiable} all standard fine maps are $!$-analytifiable.
		The second claim follows from combining \Cref{thm: DiamondAnalytificationUltimate}.(6) with \Cref{prop: allconstructiblesheavesSuaveunderembedding}.
	\end{proof}

	By \Cref{thm: DiamondAnalytificationUltimate}, we know that $c^*$ commutes in a higher-categorical way with $*$-pullback and $!$-pushforward, since this is part of what a natural transformation of 6-functor formalisms provides.
	On the other hand, from the 6-functor formalism structure and adjunction, one also gets a lax-comparison with respect to $*$-pushforward and $!$-pullback.  
	Nevertheless, one should not expect the maps arising from the lax-structure to be equivalences outside of particular circumstances.
	We record some instances where the equivalence holds, and which will be key for our purposes.

	\begin{proposition}\label{prop:analytification_commutes_!pull}
		Let $X$ and $Y$ be in $\PreStk^\Art$.
		Consider $f:X\to Y$ a map such that $f\in E^\an_\sch$.
		Then the following is true. 
		\begin{enumerate}
			\item  If $A\in \calD_\cons(Y)$, then the natural map
			\[ c_{X}^{*}f^{!} A \ra f^{\diamondsuit_\uop,!}c_{Y}^{*} A \]
			is an equivalence.
		\item If one of the following conditions hold 
			\begin{enumerate}
				\item $A\in \calD_\cons(X)$ and $f_!\bbD_{X}(A)\in \calD_\cons(Y)$,
				\item or $A\in \calD_\cons(X)$ and $f\in E^\rep_\pfp$, 
			\end{enumerate}
            then the map
				\[  c_{Y}^{*}f_{*} A  \ra  f^{{\diamondsuit_\uop}}_{*}c_{X}^{*} A\] 
				is an equivalence.
		\end{enumerate}
	\end{proposition}
	\begin{proof}
		First, we observe that since $A$ is constructible, then $f^*\bbD_Y A$ is also constructible (see \Cref{prop: allconstructiblesheavesSuaveunderembedding}). 
		In particular, by \Cref{cor: lalggivescfine}.(2), we deduce that both objects are suave over $\ast$.
	In particular, we may write 
	\[f^!A\simeq  \bbD_X f^*\bbD_Y A\simeq \SD_X f^* \SD_Y A\] using the suaveness of $A$ and of $f^*\bbD_Y A$.
		This shows that
		\begin{align*}
			c_{X}^{*}f^{!}A & \simeq c_{X}^{*} \SD_X f^*\SD_Y A\\
					& \simeq  \SD_{X^{\diamondsuit_\uop}} c_{X}^{*} f^*\SD_Y A\\
					& \simeq  \SD_{X^{\diamondsuit_\uop}} f^{{\diamondsuit_\uop},*} c_{Y}^{*} \SD_{Y} A\\
					& \simeq  \bbD_{X^{\diamondsuit_\uop}} f^{{\diamondsuit_\uop},*} \SD_{Y^{\diamondsuit_\uop}} c_{Y}^{*}  A\\
					& \simeq  f^{{\diamondsuit_\uop},!} \bbD_{Y^{\diamondsuit_\uop}}  \SD_{Y^{\diamondsuit_\uop}} c_{Y}^{*}  A\\
			 & \simeq  f^{{\diamondsuit_\uop},!} c_{Y}^{*}  A.
		\end{align*}
        Here we have repeatedly used the identification between suave duality and Verdier duality described in \Cref{prop: VerdierDualityandSuaveDuality}. 
	This finishes the proof of the first claim.

		Similarly, we note that in $(2)$ condition $(b)$ implies condition $(a)$, since $f_{!}$ preserves constructible sheaves when $f \in E^{\rep}_{\pfp}$, as in (\ref{eqn: naturalTransformofSixFunctors1}).
Since we assume that $A$ and $f_!\bbD_X A$ are constructible, and by \Cref{cor: lalggivescfine}.(2) suave over $*$, we can write 
\[f_*A\simeq \bbD_Y f_! \bbD_X A\simeq \SD_Y f_! \SD_X A.\]
We can compute that
		\begin{align*}
			c_{Y}^{*}f_{*}A & \simeq c_{Y}^{*} \SD_Y f_!\SD_X A\\
					& \simeq  \SD_{Y^{\diamondsuit_\uop}} c_{Y}^{*} f_!\SD_X A\\
					& \simeq  \SD_{Y^{\diamondsuit_\uop}} f^{{\diamondsuit_\uop}}_! c_{X}^{*} \SD_X A\\
					& \simeq  \bbD_{Y^{\diamondsuit_\uop}} f^{{\diamondsuit_\uop}}_! \SD_{X^{\diamondsuit_\uop}} c_{X}^{*}  A\\
			 & \simeq  f^{{\diamondsuit_\uop}}_* \bbD^{\diamondsuit_\uop}_X  \bbD^{\diamondsuit_\uop}_X c_{X}^{*}  A\\
			 & \simeq  f^{{\diamondsuit_\uop}}_* c_{Y}^{*}  A. \qedhere
		\end{align*}
	\end{proof}

	\subsection{On $\diamond$, $\dagger$, and $\diamondsuit$.}{\label{subsec: threeanalytificationfunctors}}
	We now introduce some important variants of analytification, moving beyond the functors $\diamondsuit$ and $\diamondsuit_{\pre}$ introduced in \S \ref{sec: the Diamond Functor}. 
	Let $X=\Spec A$ be an affine scheme over $\Spec k$. 
	When considering the analytification of $X$ it is natural to consider the following three perfectoid presheaves $X^{\diamond}$, $X^{\dagger}$, $X^\diamondsuit:(\Perff)^\op\to \Sets$
		\begin{align*}
			 X^{\diamond}:(R,R^{+}) &\mapsto X(\Spec{R^{+}})=\{f\mid f:A \to R^+\} \\
			 X^{\dagger}:(R,R^{+}) &\mapsto X(\Spec{R^{\circ}})=\{f\mid f:A \to R^\circ\} \\
			 X^\diamondsuit:(R,R^{+}) &\mapsto X(\Spec{R})=\{f\mid f:A \to R\}, 
		\end{align*}
		where $R^{\circ}$ denotes the ring of power bounded elements.

	It follows from \cite[Theorem 8.7]{Sch17} that all of the above functors are v-sheaves over $\Spd k$.
	Alternatively, recall that for a Huber pair $(A,A^+)$ one can define a perfectoid prestack $\Spd(A,A^{+})$ via the formula
	\begin{equation}{\label{eqn: SpdHuberpairs}}
		\Spd(A,A^+)[(R,R^+)]=\{f:(A,A^+)\to (R,R^+)\mid \text{with } f \text{ a map of Huber pairs}\}.
    \end{equation}
	Then $X^\diamond=\Spd(A,A)$ and $X^\diamondsuit=\Spd(A,\widetilde{k}^{\on{min}})$ where $A$ is endowed with the discrete topology and $\widetilde{k}^{\on{min}}$ denotes the integral closure of $k$ in $A$. 

	The functor $X^\dagger$ is not representable by an adic space, but it can be characterized as the smallest closed subsheaf of $X^\diamondsuit$ containing $X^\diamond$ (see \cite[Definition 2.2, Proposition 2.25]{Gle24} and also \cite[Definition 2.5]{GIZ25}).
	Alternatively, $X^\dagger$ is the canonical compactification of the structure map $X^\diamond \to \ast$ in the sense of \cite[Proposition 18.6]{Sch17}.
	We have maps of v-sheaves 
	\begin{equation}{\label{eqn: mapsofanalytifications}}
    X^\diamond\xrightarrow{a_X} X^\dagger \xrightarrow{b_X} X^\diamondsuit
    \end{equation}
	induced from the natural inclusions of rings $R^+\subseteq R^\circ\subseteq R$.
	We let $d_X:X^\diamond\to X^\diamondsuit$ denote the composition. 
For the reader's convenience, we discuss a key example.	
    \begin{example}{\label{example: Theaffineline}}
    Let $X = \bb{A}^{n,\mathrm{perf}}_{k}$ denote the perfection of the $n$-dimensional affine line over $k$. 
    The functors in \Cref{eqn: mapsofanalytifications} give rise to maps 
    \[ (\bb{A}^{n,\mathrm{perf}}_{k})^{\diamond} \xrightarrow{a_{\bb{A}^{n,\mathrm{perf}}_{k}}} (\bb{A}^{n,\mathrm{perf}}_{k})^{\dagger} \xrightarrow{b_{\bb{A}^{n,\mathrm{perf}}_{k}}} (\bb{A}^{n,\mathrm{perf}}_{k})^{\diamondsuit} \]
    of v-sheaves. 
    Consider the Huber pair $(k[t_{1}^{\onepinfty},\ldots,t_{n}^{\onepinfty}],k[t_{1}^{\onepinfty},\ldots,t_{n}^{\onepinfty}])$ equipped with the discrete topology and write $\bb{B}_{k}^{n} := \Spd(k[t_{1}^{\onepinfty},\ldots,t_{n}^{\onepinfty}],k[t_{1}^{\onepinfty},\ldots,t_{n}^{\onepinfty}])$ for the v-sheaf attached to it via \eqref{eqn: SpdHuberpairs}. 
    We see from the definition that 
    \[\bbB_k^n(R,R^+):=\{\text{Huber pair maps } f:(k[t_{1}^{\onepinfty},\ldots,t_{n}^{\onepinfty}],k[t_{1}^{\onepinfty},\ldots,t_{n}^{\onepinfty}])\to (R,R^+)\}.\]
    The definition of maps of Huber pairs forces each $t_i$ to map to $R^+$, and since $R^+$ is a perfect $k$-algebra, the whole map is determined by the image of the $t_i$. 
    We see that $\bbB_k^n$ represents the functor 
    \[\bbB^n_k:(R,R^+)\mapsto (R^+)^n, \]
    so in particular we have an identification $(\bb{A}^{n,\perf}_{k})^{\diamond} \simeq \bbB^{n}_{k}$. 
    Moreover, we see that $\bbB_k^n$ has the property that, for any map $\Spd(R,R^{+}) \ra \Spd(k)$, the base-change $\bb{B}_{k}^{n}\times_{\Spd k} \Spd(R,R^+)$ identifies with 
    \[(\bbB^{n,\on{perfd}}_R)^{\diamondsuit}= \Spd(R\langle T_{1}^{\onepinfty},\ldots,T_{n}^{\onepinfty} \rangle,R^{+}\langle T_{1}^{\onepinfty},\ldots,T_{n}^{\onepinfty} \rangle), \] 
    which is the diamond attached to the usual $n$-dimensional perfectoid closed unit ball $\bbB^{n,\on{perfd}}_R$ over $\Spd(R,R^{+})$. 
    Similarly, one can verify that $\Spd(k[t_{1}^{\onepinfty},\ldots,t_{n}^{\onepinfty}],k)$ represents the functor 
    \[(\bbA^{n,\mathrm{perf}}_k)^\diamondsuit:(R,R^+)\mapsto (R)^n,\]
    and that the pullback of $(\bb{A}^{n,\mathrm{perf}}_{k})^{\diamondsuit}\to \ast$ to $\Spd(R,R^{+})$ identifies with the diamond attached to the perfectoid $n$-dimensional affine space 
    \[\bbA^{n,\on{perfd}}_R := \bigcup_{m \geq 0} \Spa(R\langle \varpi^{m}T_{1}^{\onepinfty},\ldots,\varpi^{m}T_{n}^{\onepinfty} \rangle, R^{+}\langle\varpi^{m}T_{1}^{\onepinfty},\ldots,\varpi^{m}T_{n}^{\onepinfty}\rangle) \]
    over $(R,R^{+})$, where $\varpi \in R^{\circ}$ is a choice of pseudo-uniformizer. 
    In particular, we see that the natural map
    \[ d_{\bb{A}^{n,\mathrm{perf}}_{k}}: (\bb{A}^{n,\mathrm{perf}}_{k})^{\diamond} \hookrightarrow (\bb{A}^{n,\mathrm{perf}}_{k})^{\diamondsuit} \]
    is an open immersion. 

    Finally, consider the closure $\overline{\bbB}^{n,\mathrm{perfd}}_{R}$ of ${\bbB^{n,\mathrm{perfd}}_R}$ within $\bbA^{n,\on{perfd}}_R$.
    It is represented by the following Huber pair
    \[\overline{\bbB}^{n,\on{perfd}}_{R}= \Spa(R\langle T_{1}^{\onepinfty},\ldots,T_{n}^{\onepinfty} \rangle,R^{\on{min}}),\] 
    where $R^{\on{min}}\subseteq R\langle T_{1}^{\onepinfty},\ldots,T_{n}^{\onepinfty} \rangle$ is the smallest open bounded integrally closed subring containing $R^+$.
We have inclusions $(\bbB^{n,\on{perfd}}_R)^\diamondsuit\subseteq  (\overline{\bbB}_R^{n,\on{perfd}})^\diamondsuit\subseteq (\bbA^{n,\on{perfd}}_R)^\diamondsuit$ as diamond over $\Spd(R,R^+)$.
    On points, we have that
    \[(\overline{\bbB}^{n,\on{perfd}}_R)^\diamondsuit:(A,A^+)\mapsto (A^\circ)^n,\]
    for all perfectoid Huber pairs $(A,A^+)$ over $(R,R^+)$.
    The functor $(\bbA^n_k)^\dagger$ has formula
    \[(R,R^+)\mapsto (R^\circ)^n,\]
    and therefore we deduce that the base change of $(\bbA^n_k)^\dagger\to \ast$ along $\Spd(R,R^+)\to \ast$ agrees with $(\overline{\bbB}_R^{n,\on{perfd}})^{\diamondsuit}$.
For this reason, we also denote $(\bbA^{n,\perf}_{k})^{\dagger}$ by $\ol{\bb{B}}^{n}_{k}$.
    \end{example}

	\begin{lemma}{\label{lemma: geometricpropertiesofmapdaggertodiatodiamond}}
		Let $X=\Spec A$ with $A \in \CAlg^{\perf}$.
Then the following is true. 
		\begin{enumerate}
			\item The map  $d_X:X^{\diamond} \ra X^{\diamondsuit}$ is a qcqs pro-open monomorphism.
				If $X\in \PSch^\pfp$ then $d_X$ is an open immersion.
			\item The map $a_X:X^{\diamond} \ra X^{\dagger}$ is a qcqs pro-open monomorphism, and induces a bijection on rank one points. 
				If $X\in \PSch^\pfp$ then $a_X$ is an open immersion.
			\item The map $b_X:X^\dagger \to X^\diamondsuit$ is a closed immersion.
		\end{enumerate} 
	\end{lemma}
	\begin{proof}
		Let $X=\Spec A$. 
		For all $a\in A$, we let $X_a\subseteq X^\diamondsuit$ be the quasicompact open subsheaf with $X_a(R,R^+)\subseteq X(R,R^+)$ if the image of $a$ in $R$ lies in $R^+$. 
		More precisely, the choice $a\in A$ gives a map $X^\diamondsuit\to (\bbA^{1,\mathrm{perf}}_{k})^\diamondsuit$ and $X_a$, by definition, fits in the following Cartesian diagram
		\begin{center}
		\begin{tikzcd}
		 X_a \arrow{r} \arrow{d}  & \arrow{d} \bbB^1_{k} \\
		 X^\diamondsuit \arrow{r} & (\bbA^{1,\mathrm{perf}}_{k})^\diamondsuit,
		\end{tikzcd}
		\end{center}
		where $\bbB^1_{k}=(\bbA^{1,\mathrm{perf}}_{k})^\diamond$ is as in  \Cref{example: Theaffineline}.
As seen in this example, the map $\bbB^1_{k} \to (\bbA^{1,\mathrm{perf}}_{k})^\diamondsuit$ is a qcqs open map, since this is true after any base-change from $\Spd(k)$ to $\Spa(R,R^+)\in \Perff$.
Therefore, the same is true for $X_a\to X^\diamondsuit$.

	Clearly, we have $X^\diamond=\cap_{a\in A} X_a$. 
		For a finite set $I\subseteq A$, we let $X_I=\cap_{a\in I} X_a$, then $X^\diamond=\cap_{I\subseteq A} X_I$ is a presentation as a pro-open monomorphism.
		If $A$ is pfp over $\Spec k$ and $I\subseteq A$ is a finite set of generators of $A$ as a perfect $k$-algebra, then $X^\diamond=X_I$ and thus $X^\diamond\to X^\diamondsuit$ is an open immersion in this case.
		This finishes the proof of the first claim. 
		
		Since $X^\diamond=X^\diamond\times_{X^\diamondsuit} X^\dagger$ the properties of $a_X$ follow from those of $d_X$. 
		It follows directly from the definition that it induces a bijection of rank one points.
		This finishes the proof of the second claim.
		
		The map $b_X:X^\dagger \to X^\diamondsuit$ is clearly a monomorphism, so it is separated. 
		Now, $X^\dagger\to \Spd k$ is qcqs by \cite[Lemma 2.26]{Gle24}, and therefore $b_{X}$ is qcqs.
Moreover, $b_{X}$ is partially proper by definition.
In summary, we see that $b_{X}$ is a proper monomorphism; in other words, a closed immersion.
This finishes the proof of the third claim.
	\end{proof}
	
	One can extend the definition of $\diamond$, $\dagger$ and $\diamondsuit$ to prestacks $X\in \PreStk$ by considering functors 
	\begin{align*}
		\label{some functors pre}
		X^{\diamond_\pre}:(R,R^{+})& \mapsto X(\Spec{R^{+}}) \\
		X^{\dagger_\pre}:(R,R^{+})& \mapsto X(\Spec{R^\circ}) \\
		X^{\diamondsuit_\pre}:(R,R^{+})& \mapsto X(\Spec{R}). 
	\end{align*}
	These will give rise to perfectoid prestacks $X^{\diamond_\pre}$, $X^{\dagger_\pre}$ and $X^{\diamondsuit_\pre}$. 
	Finally, given $X\in \PreStk$ we let $X^\diamond$, $X^\dagger$ and $X^\diamondsuit$ denote the respective sheafifications.
	Formally, we simply took the left Kan extension along the Yoneda embedding of the previously defined functors as in the following diagram 

	\begin{center}
	\begin{tikzcd}
		\PSchf \arrow{r}{y}  \arrow[d,  "{\diamond, \dagger, \diamondsuit}" description]  & \PreStk \arrow[dl, dotted, "{\diamond, \dagger, \diamondsuit}"] \\
		\AnStk_{v}.  & 
	\end{tikzcd}
	\end{center}

	We can of course sheafify instead with respect to the pro\'etale topology which gives rise to spaces $X^{\diamond_\uop}$, $X^{\dagger_\uop}$ and $X^{\diamondsuit_\uop}$.
	This gives a diagram

	\begin{center}
	\begin{tikzcd}
		\PSchf \arrow{r}{y}  \arrow[d,  "{\diamond, \dagger, \diamondsuit}" description]  & \PreStk \arrow[dl, dotted, "{\diamond_\uop, \dagger_\uop, \diamondsuit_\uop}"] \\
		\AnStk_{\uop}.  & 
	\end{tikzcd}
	\end{center}

	Most of our reasoning will go through $\diamond$, $\dagger$ and $\diamondsuit$, but we will occasionally need $\diamond_\uop$ and $\diamondsuit_\uop$ to argue.
	Let us record the following useful statement comparing these functors when applied to proper maps.

	\begin{lemma}\label{lm:proper_gives_cartesian_diam}
		Let $f \colon X \to Y$ be a proper morphism in $(\PreStk,E^\rep_\pfp)$. 
		Then all the squares in the following commutative diagram 
		\begin{center}
			\begin{tikzcd}
				X^\diamond \arrow{r} \arrow{d}{f^\diamond} & X^\dagger \arrow{r} \arrow{d}{f^\dagger} & X^\dia \arrow{d}{f^\dia} \\
				Y^\diamond \arrow{r} & Y^\dagger \arrow{r} & Y^\dia
			\end{tikzcd}
		\end{center}
		are Cartesian.
	\end{lemma}
	\begin{proof}
		As the functors $\diamond$, $\dagger$, $\diamondsuit$ take values in v-sheaves, it suffices to replace them by the corresponding presheaf-functors $\diamond_{\rm pre}$, $\dagger_{\rm pre}$, $\diamondsuit_{\rm pre}$ and show that the diagram of presheaves is Cartesian on a basis for the v-topology. 
		We can test this on strict products of points $S = \Spa(R,R^+)$, as these form a basis for the v-topology, see \cite[Definition 2.14, Remark 2.15]{GIZ25}. 
		For the right hand square, this amounts to showing that any commuting square
		\begin{center}
			\begin{tikzcd}
				\Spec R \arrow{r} \arrow[hookrightarrow]{d} & X \arrow{d} \\
				\Spec R^\circ \arrow{r} & Y
			\end{tikzcd}
		\end{center}
		admits a unique commuting section $\Spec R^\circ \rightarrow X$. 
		Without loss of generality $Y=\Spec R^\circ$, so we may assume $X$ is a pfp proper scheme over $Y$, by the assumed representability of $f$. 
		By \cite[Proposition 2.19]{GIZ25} (and its proof), $\Spec R^\circ$ is a strict comb and $\Spec R$ is an open subset meeting each connected component. 
		Thus, componentwise, the claim follows from the valuative criterion of properness. 
		For $x \in \pi_0(\Spec R^\circ)$, let $R^\circ_x$ be the local ring of $\Spec R^\circ$ at the unique closed point $x_0$ of $x$. 
		It suffices to show that the unique section $s_x \colon \Spec R_x^\circ \rightarrow X$ extends to a section on a clopen neighboorhood of $x$. 
		Replacing $Y$, and then $X$, by (any) affine open containing the image of $x_0$, we may assume that $X = \Spec B$, $Y = \Spec A$ are affine. 
		Then $s_x$ corresponds to an $A$-algebra map $B \rightarrow R_x^\circ = \colim_{x\in V} R_V^\circ$, where the colimit is taken over all clopen $x \in V \subseteq \pi_0(\Spec R)$. 
		Moreover, replacing $A \rightarrow B$ by a deperfection, we may assume that $B$ is finitely presented over $A$. 
		It then follow that there is some clopen neighboorhood $x \in V \subseteq \pi_0(\Spec R^\circ)$, such that $B \rightarrow R_x^\circ$ factors through $B \rightarrow R_V^\circ$, i.e., $s_x$ extends to the preimage of $V$ in $\Spec R^\circ$, and we are done. 
		The proof for the left square is the same, with $\Spec R \rightarrow \Spec R^\circ$ replaced by $\Spec R^\circ \rightarrow \Spec R^+$.
	\end{proof}

	\begin{remark}
		\label{case-of-thepoint}
		The particular case of \Cref{lm:proper_gives_cartesian_diam} when $Y=\ast=\Spec k$ simply says that $X^\diamond=X^\dagger=X^\diamondsuit$, since by definition we are forcing $\ast^\diamond=\ast^\dagger=\ast^\diamondsuit=\Spd(k,k)$.
		Indeed, our functors are only defined on perfectoid spaces over $\Spd(k,k)$. 
		In particular, note that in this case $X^\diamond\to \ast$ is also qcqs, by combining \Cref{lm:proper_gives_cartesian_diam} with \Cref{prop: hcoversgivediamondcovers}. 
	\end{remark}

	Extracting the formal ingredients from the proof of \Cref{lm:proper_gives_cartesian_diam}, we arrive at the following generalization. 

	\begin{lemma}\label{lm:proper_gives_cartesian_diam2}
		Let $f \colon X \to Y$ be a morphism of prestacks. 
		Suppose that the following conditions hold.
		\begin{enumerate}
			\item $f$ is relatively representable in Zariski sheaves of sets.
			\item $f$ satisfies the valuative criterion of properness.
In other words, for any valuation ring $V$ with fraction field $K$ the natural map 
				\[X(V)\xrightarrow{\simeq} Y(V)\times_{Y(K)} X(K)\]
				is an equivalence.
			\item $f$ is a finitely presented map of functors.
In other words, for every ring $A$ and every $A$-algebra $B$ the natural map 
				\[ \colim X(B_i)\times_{Y(B_i)} Y(A)\xrightarrow{\simeq} X(B)\times_{Y(B)} Y(A),\]
				is an equivalence.
Here $B=\colim B_i$ is a presentation of $B$ as a filtered colimit of perfectly finitely presented $A$-algebras.
			
		\end{enumerate}
		Then all the squares in the following commutative diagram 
		\begin{center}
			\begin{tikzcd}
				X^\diamond \arrow{r} \arrow{d}{f^\diamond} & X^\dagger \arrow{r} \arrow{d}{f^\dagger} & X^\dia \arrow{d}{f^\dia} \\
				Y^\diamond \arrow{r} & Y^\dagger \arrow{r} & Y^\dia
			\end{tikzcd}
		\end{center}
		are Cartesian.
In particular, this holds when $f$ is relatively representable by perfectly finitely presented and proper algebraic spaces.	
	\end{lemma}

        We now discuss how the three functors interact with analytification. 
	\subsubsection{Three analytification functors}
	As discussed in \Cref{thm: DiamondAnalytificationUltimate} and \Cref{thm: DiamondAnalytificationUltimate2}, 
	given $X\in \PreStk$, we get a fully faithful analytification functor 
	\begin{equation} \label{c analytification} c_{X}^{*}: \calD^\sch_\Lambda(X) \rar 
 \calD_{\Lambda}^{\diamondsuit}(X) = \calD_{\Lambda}^{\an}(X^{\diamondsuit}). \end{equation} 

 In what follows, we study the following variants obtained by composition 
	\begin{equation} \label{b analytification} b^{*,\an}_{X} := b^{*}_{X} \circ c_{X}^{*}: \calD^\sch_\Lambda(X) \ra \calD^\an_\Lambda(X^{\dagger}) \end{equation} 
	and 
	\begin{equation} \label{d analytification} d^{*,\an}_{X} := d^{*}_{X} \circ c_{X}^{*}: \calD^\sch_\Lambda(X) \ra \calD^\an_\Lambda(X^{\diamond}). \end{equation} 
	The functor $b_X^{*,\an}$ will be key to our construction of the equivalence $\pitch$. 
	On the other hand, the functors $d_X^{*,\an}$ will prove to be very useful to perform computations.
The idea is roughly that we can use properties of $d_{X}^{*,\an}$ to deduce properties of $b_{X}^{*,\an}$ using that analytification takes values in overconvergent sheaves and that difference between $b_{X}^{*,\an}$ and $d_{X}^{*,\an}$ on overconvergent sheaves is immaterial.
    
	Just as we organized $\calD_\et(X^\diamondsuit)$ into a functor 
	\[\calD_\Lambda^\diamondsuit: \PreStk^{\op} \to \LinCat_\Lambda,\] it will be useful to do the same for the $\dagger$ and the $\diamond$ constructions. 

	\begin{definition}
	We define functors, 
	\[\calD^\diamond_\Lambda(-),\calD^\dagger_\Lambda(-),\calD^\diamondsuit_\Lambda(-): \PSch^{\aff,\op} \to \LinCat_\Lambda\]
	with formula	
	\[A\mapsto \calD_\Lambda^\an(\Spec(A)^\diamond),\text{ } \calD_\Lambda^\an(\Spec(A)^\dagger) \text{ and } \calD_\Lambda^\an(\Spec(A)^\diamondsuit)\] 
	respectively. 
	We define 
	\[\calD^\diamond_\Lambda(-),\calD^\dagger_\Lambda(-),\calD^\diamondsuit_\Lambda(-): \PreStk^\op\to \LinCat_\Lambda\]
	as the right Kan extension along the Yoneda embedding of the functors above.
	\end{definition}

	By definition, we have a tautological commutative diagram  
	\begin{center}
	\begin{tikzcd}
		(\PSchf)^\op \arrow{rr}{\diamond\, \dagger\, \diamondsuit} \arrow[dr, "\calD_\Lambda^\diamond\, \calD_\Lambda^\dagger\, \calD_\Lambda^\diamondsuit", swap] &  & \AnStk^{\op}_v \arrow{dl}{\calD_{\Lambda}^{\an}(-)} \\
					 & \LinCat_\Lambda & 
	\end{tikzcd}
	\end{center}

	The following statement extends this to prestacks. 
        \begin{proposition}{\label{prop: compatibilityofKanextensionswithDan}}
        We have a commutative diagram of the following form
        \begin{center}
	\begin{tikzcd}
		\PreStk^\op \arrow{rr}{\diamond, \dagger, \diamondsuit} \arrow[dr, "\calD_\Lambda^\diamond\, \calD_\Lambda^\dagger\, \calD_\Lambda^\diamondsuit", swap] &  & \AnStk_{v}^{\op} \arrow{dl}{\calD_{\Lambda}^{\an}} \\
					 & \LinCat_\Lambda. & 
	\end{tikzcd}
	\end{center}
        \end{proposition}
        \begin{proof}
        For $? \in \{\diamond,\dagger,\diamondsuit\}$, we need to show that there exists a natural equivalence of the form
        \[
        \calD^{?}_{\Lambda}(-) \simeq \calD_{\Lambda}^{\an}((-)^{?}).
        \]
	We observe that all of the functors $\calD^\an_\Lambda:\AnStk^\op_v\to \LinCat_\Lambda$, $\calD^?_\Lambda:\PreStk^\op \to \LinCat_\Lambda$ and $(-)^?:\PreStk^\op\to \AnStk_v^\op$ preserve limits.
	Indeed, $(-)^?$ and $\calD^?_\Lambda$ are obtained as a right Kan extension of the Yoneda embedding $(\PSchf)^\op\subseteq \PreStk^\op$, so \cite[Lemma 5.1.5.5.(1)]{HTT} applies.
	On the other hand, $\calD^\an_\Lambda$ satisfies v-descent and can be regarded as a right Kan extension along the Yoneda embedding $(\Perff)^\op\subseteq \AnPreStk^\op$ which happens to factor through v-sheafification $\AnPreStk^\op\to\AnStk_v^\op$ (see \Cref{rem: perfectoidprestacks}).
	In particular, \cite[Lemma 5.1.5.5.(1)]{HTT} applies again. 
	Since both $\calD^{?}_{\Lambda}(-)$ and $\calD_{\Lambda}^{\an}((-)^{?})$ commute with limits, by \cite[Theorem 5.1.5.6]{HTT}, it suffices to construct an equivalence of functors after restricting to $\PSchf$, but on this subcategory the functors agree by definition.
        \end{proof}
	A consequence of \Cref{prop: compatibilityofKanextensionswithDan} is that we get natural transformations of functors on $\PreStk^\op$ with values in $\CAlg(\LinCat_\Lambda)$ 
	\begin{equation}{\label{eqn: naturaltransformsofsheaftheories}}
    \calD^\sch_\Lambda \xRightarrow{c^*} \calD^\diamondsuit_\Lambda \xRightarrow{b^*} \calD^\dagger_\Lambda \xRightarrow{a^*} \calD^\diamond_\Lambda, 
    \end{equation}
    where the first natural transformation is as in \Cref{thm: DiamondAnalytificationUltimate}, and the other natural transformations come via pulling back along the natural maps $X^{\diamond} \overset{a_X}{\ra} X^{\dagger} \overset{b_X}{\ra} X^{\diamondsuit}$. 
    As above, we let $b^{*,\an}_{(-)}:=b_{(-)}^*\circ c_{(-)}^*$ and $d^{*,\an}_{(-)}=a_{(-)}^* \circ b_{(-)}^*\circ c_{(-)}^*$.

	\begin{proposition}
		\label{all-are-fully-fatihful}
		For all $X\in \PreStk$, the maps $c_X^*$, $b_X^{*,\an}$ and $d_X^{*,\an}$ are fully faithful.
	\end{proposition}
	\begin{proof}
		Since every $X\in \PreStk$ is a colimit of representables objects and limits of fully faithful functors are again fully faithful, we may assume that $X=\Spec A$ by using \Cref{prop: compatibilityofKanextensionswithDan}.  
		In this case, the claim for $c_X^*$ is \cite[Proposition 27.2]{Sch17}.
		The claim for $d_X^{*,\an}$ is \cite[Proposition 4.2]{GL22}.
		Finally, the claim for $b_X^{*,\an}$ follows from \cite[Lemma 4.1]{GL22}.
	\end{proof}
	
	\subsubsection{Descent properties}
	Recall that $\calD^{\on{sch}}_\Lambda(-)$ satisfies schematic v-descent. 
	This is not at all clear (and might actually fail) for $\calD^\diamondsuit_\Lambda(-)$. 
	Indeed, as explained in \cite[Example~A.3]{AGLR22} the functor $(-)^{\diamondsuit}$ does not preserve v-covers. 
	This motivates the following definition. 
	\begin{definition}
		\label{diamond topology}
		Given a family of maps of affine schemes $\{X_i\to Y\}_{i\in I}$, we say that it is a $\diamondsuit$-cover if the map 
		\[\coprod_{i\in I}X_i^\diamondsuit\to Y^\diamondsuit\] 
		is v-surjective.	
		The $\diamondsuit$-topology is the topology on $\PSchf$ that is generated by $\diamondsuit$-covers. 
	\end{definition}
	
	\begin{remark}
		Analogously to \Cref{diamond topology}, one can define the $\diamond$-topology and the $\dagger$-topology. 
		Nevertheless, by \Cref{preserves v-covers small d} and \Cref{arc covers2} we simply recover the schematic v-topology and the arc-topology respectively. 
	\end{remark}

	\begin{proposition}
		\label{preserves v-covers small d}
		Let $f:X\to Y$ be a map of perfect affine schemes.
The following are equivalent.
		\begin{enumerate}
			\item $f:X\to Y$ is a schematic v-cover.
			\item $|f^\diamond|:|X^\diamond|\to |Y^\diamond|$ is surjective.
			\item $f^\diamond:X^\diamond\to Y^\diamond$ is v-surjective.
		\end{enumerate}
	\end{proposition}
	\begin{proof}
		This follows from the proof of \cite[Proposition 3.7]{Gle24}.
	\end{proof}

	\begin{proposition}
		\label{arc covers2}		
		Let $f:X\to Y$ be a map of perfect affine schemes.
The following are equivalent.
		\begin{enumerate}
			\item $f:X\to Y$ is an arc-cover.
			\item $|f^\dagger|:|X^\dagger|\to |Y^\dagger|$ is surjective.
			\item $f^\dagger:X^\dagger\to Y^\dagger$ is v-surjective.
		\end{enumerate}
	\end{proposition}
	\begin{proof}
		By \cite[Lemma 2.26]{Gle24} and \cite[Lemma 12.11]{Sch17}, (2) and (3) are equivalent. 
        
		We now show the equivalence of (1) and (2).
Let $X=\Spec A$ and let $Y=\Spec B$.
We note that $X\to Y$ is an arc-cover in the sense of \cite[Definition 1.2]{BhattMathew_21} if and only if every map $\Spec(V)\to Y$ from a rank $1$ valuation ring lifts to $\Spec(W)\to X$ for some choice of rank $1$ valuation ring $W$ over $V$, such that $V\to W$ is an extension of valuation rings (i.e., a faithfully flat map $V\to W$).
		We may restrict our attention to absolutely integrally closed valuation rings, which we may further assume to be complete, since the completion $V\to \hat{V}$ is an integral extension of valuation rings whenever $V$ is of rank $1$. 
		Write $K_V$ for $\mathrm{Frac}(V)$ and $K_W$ for $\mathrm{Frac}(W)$. 
		Endowing $V$ and $W$ with their adic topologies, we get algebraically closed perfectoid fields $(K_V,V)$ and $(K_W,W)$ and we see that we have an identification
		\[Y^\dagger(K_V,V)=Y(V) \text{ and } X^\dagger(K_W,W)=X(W).\]

Let us show the implication $(1)\implies (2)$.
Since $X^\dagger\to Y^\dagger$ is partially proper, the map is specializing and it suffices to show surjectivity on rank $1$ points.
Any rank $1$ point $y\in \mid Y^\dagger\mid$ is represented by a map of the form $\Spa(C_1,O_{C_1})\to Y^\dagger$, which we may interpret as a map of the form $\Spec O_{C_1}\to Y$. 
By (1) and our considerations above, we may find an extension of complete rank $1$ valuation ring $O_{C_1}\to O_{C_2}$ and a map $\Spec O_{C_2}\to X$ lifting the map $\Spec O_{C_1}\to Y$.
This constructs a map of the form $\Spa(C_2,O_{C_2})\to X^\dagger$ lifting $\Spa(C_1,O_{C_1})\to Y^\dagger$.

Conversely, if $\mid X^\dagger\mid \to \mid Y^\dagger\mid$ is surjective, the corresponding map at the level of rank $1$ points is also surjective. 
A map of the form $\Spec V\to Y$ with $V$ complete and having algebraically closed field, produces a geometric point $\Spa(K_V,V)\to Y^\dagger$.
By the assumed surjectivity, we may find a geometric point $\Spa(K_W,W)\to X^\dagger$ lifting $\Spa(K_V,V)\to Y^\dagger$, which produces an extension of valuation rings $V\to W$ and a map $\Spec W\to X$ lifting $\Spec V\to Y$.
	\end{proof}

	\begin{proposition}
		If $f:X\to Y$ is a map of affine schemes and it is a $\diamondsuit$-cover, then it is an arc-cover.
	\end{proposition}
	\begin{proof}
		Let $X=\Spec A$ and $Y=\Spec B$.
		Without loss of generality $B=V$ is a rank one valuation ring with pseudo-uniformizer $\pi\in V$.
		Let $A_s$ and $A_\eta$ denote the perfection of $A/\pi$ and $A[\frac{1}{\pi}]$, respectively.
		Suppose that $X\to Y$ is not an arc-cover, then $\Spec A_\eta\subseteq \Spec A$ is a constructible open subset which is also stable under specialization, which implies that $\Spec A_\eta\subseteq \Spec A$ is closed in the Zariski topology. 
		This gives $\Spec A_s\coprod \Spec A_\eta=\Spec A$ and $f$ factors through a map $X\to \Spec V[\frac{1}{\pi}]\coprod \Spec (V/\pi)^\perf$.
		Let $Z=\Spec V[\frac{1}{\pi}]\coprod \Spec (V/\pi)^\perf$, it suffices to show that $g:Z\to Y$ is not a $\diamondsuit$-cover.
		Since $g$ is a monomorphism $g^\diamondsuit$ is also a monomorphism, if $g$ was a $\diamondsuit$ cover then $g^\diamondsuit$ would be an isomorphism, but $Z^\diamondsuit$ has two connected components while $Y^\diamondsuit$ has one.
	\end{proof}

    	We warn the reader that the $\diamondsuit$-topology is not finitary in contrast to the arc- and v-topologies on perfect schemes.
		For this reason, we consider the following variant of \Cref{diamond topology}.
		\begin{definition}
			We define the $\diamondsuit^{\on{fin}}$-topology on $\PSchf$ to be the finitary topology generated by maps $f:X\to Y$ in $\PSchf$ such that $f$ is a $\diamondsuit$-cover as in \Cref{diamond topology}. 
		\end{definition}

	As it turns out, there are  $\diamondsuit$-covers that are not schematic v-covers (see \Cref{A counterexample}). 
	This is intimately related to the failure of the $\diamondsuit$-topology to be finitary.
	Similarly, there exist maps of affine schemes which are schematic v-covers (or $\diamond$-covers) but are not $\diamondsuit$-covers, see \cite[Example~A.3]{AGLR22}.
To remedy this, we recall that a map $X\to Y$ in $\PSchf$ is an \emph{h-cover} if it is pfp and a v-cover then we have the following.
	\begin{proposition}{\label{prop: hcoversgivediamondcovers}}
		If $X\to Y \in \PSchf$ is an h-cover then it is also a $\diamondsuit$-cover. 
		In particular, \'etale covers are $\diamondsuit$-covers.
	\end{proposition}
	\begin{proof}
		This is the content of \cite[Lemma~A.2]{AGLR22} or \cite[Proposition 5.4]{Gle24}.
	\end{proof}

	For each of $\tau\in \{\diamondsuit, \diamondsuit^{\on{fin}}, \on{v}, \on{arc}, \on{h}, \et\}$ we consider the subcategory $\SchStk_{\tau} \subseteq \PreStk$ of $\tau$-sheaves.  
	We summarize the relationship between the different topologies described above with the following diagram in which each of the arrows are inclusions that admit a left-adjoint (sheafification).

	\begin{center}
	\begin{tikzcd}
 (\SchStk)_{\diamondsuit} \ar{r} & (\SchStk)_{\diamondsuit^{\on{fin}}} \ar{rd}&  &		\\
 (\SchStk)_{\on{arc}} \ar{ur}\ar{dr}	& & (\SchStk)_{h} \ar{r}  &(\SchStk)_{\et} \\
						& (\SchStk)_{\on{v}}\ar{ur}  & & 		\\
	\end{tikzcd}
	\end{center}

	These considerations and the fact that $\calD^\an_\Lambda$ satisfies analytic v-descent (Remark \ref{rem: perfectoidprestacks}) has the following consequence.

	\begin{proposition}{\label{prop: descentproperties}}
		Consider $\calD^\diamond_\Lambda$, $\calD^\dagger_\Lambda$, $\calD^\diamondsuit_\Lambda$ as functors on $\PreStk$ with values in $\LinCat_\Lambda$. 
		Then the following statements hold.
		\begin{enumerate}
			\item $\calD_\Lambda^\diamond$ is a v-sheaf, so it factors through $(\SchStk)_{\on{v}}$.
			\item $\calD_\Lambda^\dagger(-)$ is an arc-sheaf, so it factors through $(\SchStk)_{\on{arc}}$.
			\item $\calD^\diamondsuit_\Lambda(-)$ is a $\diamondsuit$-sheaf so it factors through $(\SchStk)_{\diamondsuit}$. 
			\item The 3 functors considered above are h-sheaves, so they all factor through $(\SchStk)_{\on{h}}$. 
			\item The 3 functors considered above are \'etale sheaves, so they all factor through $(\SchStk)_{\et}$. 
		\end{enumerate}
	\end{proposition}
	\begin{proof}
		For the restriction of the functors to $(\PSchf)^{\op} \subset \PreStk^{\op}$, the desired claim was discussed above. 
		To deduce the factorization of the map $\PreStk^\op\to \CAlg(\LinCat_\Lambda)$ through sheafification $\PreStk^\op\to \SchStk^\op_{?}$ if suffices to recall \cite[Proposition 1.3.1.7]{SAG}. 
		Indeed, by \cite[Proposition 5.5.4.20]{HTT} forming categories of limit preserving functors gives a fully faithful embedding
		\[\on{Fun}^R((\SchStk_{?})^\op,\CAlg(\LinCat_\Lambda))\to \on{Fun}^R(\PreStk^\op,\CAlg(\LinCat_\Lambda)).\]
		Moreover, by Yoneda \cite[Theorem 5.1.5.6]{HTT} we have an equivalence
		\[\on{Fun}^R(\PreStk^\op,\CAlg(\LinCat_\Lambda))\simeq \on{Fun}(\PSchf,\CAlg(\LinCat_\Lambda)).\]
		This gives overall a fully faithful embedding
		\[\on{Fun}^R((\SchStk_{?})^\op,\CAlg(\LinCat_\Lambda))\to \on{Fun}(\PSchf,\CAlg(\LinCat_\Lambda)),\]
		and \cite[Proposition 1.3.1.7]{SAG} identifies the essential image. 
	\end{proof}
	\begin{remark}
		At this point we should emphasize that if we are given a presentation of a quotient stack $Y=[X/G]\in \PreStk$, and we wish to compute $\calD_\Lambda^\diamondsuit(Y)$, it is cardinal to clarify if there was a sheafification process involved in the definition of $Y$. 
		Indeed, if $Y_v$ denotes the schematic v-sheafification of $Y$, then there is no reason to expect that $\calD^\diamondsuit_\Lambda(Y)$ should agree with $\calD^\diamondsuit_\Lambda(Y_v)$. 
	\end{remark}

	\subsection{Computing $d^*_X$ for Artin stacks} 
	\label{Computing d}
	In contrast with $\diamondsuit$ and $\diamondsuit_\uop$, the analytification functors $\diamond$ and $\dagger$ do not behave correctly with respect to $f_!$ and do not give rise to morphisms of 6-functor formalisms.
	In particular, it is no longer automatic that the functors $\diamond$ and $\dagger$ preserve suave maps.
	For example, the open immersion $\bbG^{\mathrm{perf}}_m\to \bbA^{1,\mathrm{perf}}$ becomes a closed immersion $\bbG^{\mathrm{perf},\dagger}_m\to \bbA^{1,\mathrm{perf},\dagger}$, and this latter morphism is certainly not $\calD^\an_\Lambda$-suave. 
	In what follows, we will show that, under finiteness assumptions, the $\diamond_\uop$-functor preserves suave maps.

	\begin{lemma}{\label{lemma: Diamondtodiamond}}
		Suppose that we have a map $[f: X \ra Y] \in (\PreStk,E^\rep_\pfp)$ and that it is $\calD^{\mathrm{sch}}_{\Lambda}$-suave, then $[f^{\diamond_\uop}:X^{\diamond_\uop} \ra Y^{\diamond_\uop}]$ is $\calD^{\uop}_{\Lambda}$-suave. 
	\end{lemma}
	\begin{proof}
		We have to show that $f^{\diamond_\uop}$ is $!$-able for $\calD^\uop_\Lambda$, and that it is $\calD^{\uop}_{\Lambda}$-suave, we start with a reduction step.
		We claim that we have $f^{\diamond_\uop}\in E^\rep_\fdcss$. 
		Since $E_\fdcss$ is pro\'etale local (see \Cref{pro-etale locality}), and pro\'etale locally every locally spatial diamond is a disjoint union of strictly totally disconnected spaces (see \Cref{uop properties}), we may assume $Y=\Spec R^+$ and that $X\to \Spec R^+$ is in $E_\pfp$.  
		Indeed, pro\'etale locally, any map $\Spa(R,R^+)\to Y^{\diamond_\uop}$ factors as
		\[\Spa(R,R^+)\to (\Spec R^+)^{\diamond_\uop} \to Y^{\diamond_\uop}.\] 

		By \cite[Lemma~4.5.7]{HeyerMann}, we may check suaveness locally. 
		Applying a similar reasoning as above, suaveness can also be verified in the case $Y=\Spec R^+$ for $\Spa(R,R^+)$ a strictly totally disconnected once we have shown it is $!$-able.
		This finishes the reduction step which allow us to assume that $Y$ is affine and that $X\to Y$ is in $E_\pfp$.
		In this case, $\diamondsuit$ agrees with $\diamondsuit_{\uop}$, and $\diamond$ agrees with $\diamond_\uop$. 

		The second reduction step goes as follows. 
		Fix an open cover $X=\cup_{i\in I} U_i$ by affine schemes, we obtain an open cover of v-sheaves $X^\diamond=\cup_{i\in I}U_i^\diamond$. 
		Note that suaveness, and $!$-ability, is suave-local on the source (see \cite[Lemma~4.5.8 (i), Remark 4.4.11, Definition 3.1.6]{HeyerMann} for the precise statement). 
        In this particular case, it suffices to show that $U_i^\diamond \to Y^\diamond$ is $\calD^\an_\Lambda$-suave.
		Replacing $X$ by $U_i$, we may now assume that both $X$ and $Y$ are affine.

		Finally, consider the following commutative diagram with Cartesian square
		\[ 
		\begin{tikzcd}
			X^{\diamond} \arrow[r] \ar{rd} & X^{\diamondsuit} \times_{Y^{\diamondsuit}} Y^{\diamond} \arrow[d] \arrow[r] & Y^{\diamond} \arrow[d]  \\
			& X^{\diamondsuit} \arrow[r] & Y^{\diamondsuit}.
		\end{tikzcd}
		\]
		By \Cref{thm: DiamondAnalytificationUltimate}.(7) applied to $f^\diamondsuit$ together with \cite[Lemma~4.5.9 (i)]{HeyerMann}, the map $X^{\diamondsuit} \times_{Y^{\diamondsuit}} Y^{\diamond} \ra Y^{\diamond}$ is $\calD^{\an}_{\Lambda}$-suave. 
		This reduces us to showing that $[X^{\diamond} \ra X^{\diamondsuit} \times_{Y^{\diamondsuit}} Y^{\diamond}]$ is $\calD^{\an}_{\Lambda}$-suave. 
		By \Cref{lemma: geometricpropertiesofmapdaggertodiatodiamond} and the fact that pro-open immersions are right cancellative, $X^\diamond\to X^{\diamondsuit} \times_{Y^{\diamondsuit}} Y^{\diamond}$ is a pro-open monomorphism, and we claim that in this case it is actually open and consequently $\calD^{\an}_{\Lambda}$-suave. 
		Indeed, after fixing a closed embedding $g:X\hookrightarrow Y\times \bbA_k^{n,\mathrm{perf}}$, from \Cref{lm:proper_gives_cartesian_diam} applied to $g$ we get a commutative diagram with Cartesian squares 
		\begin{center}
\begin{tikzcd}[row sep=small, column sep=small]
	X^\diamond \arrow[dr, "\simeq"] \arrow[rr] \arrow[dd] & & (Y\times \bbA_k^{n,\mathrm{perf}})^\diamond \arrow[dd] \arrow[dr, "\simeq"]  & &  \\
				     & X^\diamond \arrow[rr] \arrow[dd]         & & (Y\times \bbA_k^{n,\mathrm{perf}})^\diamond \arrow[dd]            & &   \\
  X^\diamondsuit\times_{Y^\diamondsuit} Y^\diamond \arrow[rr] \arrow[dr]           & & Y^\diamond\times (\bbA_k^{n,\mathrm{perf}})^\diamondsuit \arrow[dr] \arrow[rr]           & & Y^\diamond\arrow[dr] \\
				    & X^\diamondsuit \arrow[rr]                    & & (Y\times \bbA_k^{n,\mathrm{perf}})^\diamondsuit \arrow[rr]                      & &Y^\diamondsuit 
\end{tikzcd}
		\end{center}
		from which we can extract a Cartesian square
		\begin{center}
		\begin{tikzcd}
		 X^\diamond \arrow{r} \arrow{d}  & Y^\diamond \times (\bbA_k^{n,\mathrm{perf}})^\diamond \arrow{d} \\
		 X^{\diamondsuit} \times_{Y^{\diamondsuit}} Y^{\diamond} \arrow{r} & Y^\diamond \times (\bbA_k^{n,\mathrm{perf}})^\diamondsuit.
		\end{tikzcd}
		\end{center}
		We finish by observing that $(\bbA_k^{n,\mathrm{perf}})^\diamond\to (\bbA_k^{n,\mathrm{perf}})^\diamondsuit$ is an open immersion, as in \Cref{example: Theaffineline} and \Cref{lemma: geometricpropertiesofmapdaggertodiatodiamond}.(1).
	\end{proof}

\begin{proposition}
	\label{algebraic-2-diamonds}
Fix an Artin stack $[f: X \ra \ast] \in \PreStk^\Art$ and a constructible complex $A \in \calD_{\cons}(X)$.
The following statements hold.
\begin{enumerate}
	\item The natural map $d_{X}: X^{\diamond_\uop} \ra X^{\diamondsuit_\uop}$ is $\calD^{\uop}_{\Lambda}$-suave. 
	\item $c^*_X A$ is $f^{\diamondsuit_\uop}$-suave. 
	\item $d^{*,\an}_X A$ is $f^{\diamond_\uop}$-suave. 
\end{enumerate}
\end{proposition}
\begin{proof}

	For the first claim, choose a $\calD^{\mathrm{sch}}_{\Lambda}$-suave atlas $[U \ra X]\in E^\rep_\pfp$, then we have a commutative (but usually not Cartesian) diagram 
		\[ 
		\begin{tikzcd}
			U^{\diamond_\uop} \arrow[r] \arrow[d] & X^{\diamond_\uop} \arrow[d] \\
			U^{\diamondsuit_\uop} \arrow[r] & X^{\diamondsuit_\uop}.  
		\end{tikzcd}
		\]
		From \Cref{lemma: Diamondtodiamond}, we know that $U^{\diamond_\uop} \ra X^{\diamond_\uop}$ is a $\calD_{\Lambda}^{\uop}$-suave cover. 
		Using \cite[Lemma 4.5.8.(i)]{HeyerMann}, we reduce to showing that $U^{\diamond_\uop} \ra X^{\diamondsuit_\uop}$ is $\calD_{\Lambda}^{\uop}$-suave. 
		By \Cref{thm: DiamondAnalytificationUltimate}.(7), we know that $U^{\diamondsuit_\uop} \ra X^{\diamondsuit_\uop}$ is $\calD_{\Lambda}^{\uop}$-suave. 
		Since $U \ra \ast$ is in $E_\pfp$, it follows from \Cref{lemma: geometricpropertiesofmapdaggertodiatodiamond} (1) that $U^{\diamond} \ra U^{\diamondsuit}$ is an open immersion and $\calD_{\Lambda}^{\uop}$-suave. 
		The claim now follows since suave maps are stable under composition (\cite[Lemma~4.5.9 (i)]{HeyerMann}).
		The second claim is \Cref{cor: lalggivescfine}.
		The third claim follows from combining the first claim with the fact that pullback under suave maps preserves suave objects by \cite[Lemma~4.5.16 (i)]{HeyerMann}.
\end{proof}

	\subsection{Overconvergent replacement}
	\label{overconv repl}
	The functors $\calD^\dagger_\Lambda$ and $\calD^\diamond_\Lambda$ on $\PreStk$ are different, but as we explain below, they agree on overconvergent objects.

		Let $S\in \Perff$ be a geometric point i.e., $S=\Spa(C,C^+)$ with $C$ an algebraically closed non-Archimedean field and $C^+\subseteq C$ an open and bounded valuation subring. 
		Let $\pi:S\to \ast$ denote the structure map to the final object.  
		In this setup, the constant sheaf functor 
		\[\pi^*:{\rm Mod}_\Lambda\simeq \calD^\an_\Lambda(\ast)\to \calD^\an_\Lambda(S)\] 
		is fully faithful using \cite[Theorem~1.13]{Sch17}, and we say that a sheaf $A\in \calD^\an_\Lambda(S)$ is \emph{overconvergent} if it lies in the essential image of $\pi^*$.

	\begin{definition}
		Fix $X\in \AnStk_v$ and $A\in \calD^\an_\Lambda(X)$. 
		\begin{enumerate}
			\item We say that $A$ is overconvergent if $\bar{x}^*A$ is overconvergent in $\calD^\an_\Lambda(S)$ for every geometric point $\bar{x}:S\to X$ with $S=\Spa(C,C^+)$. 
			\item We let $\calD_\Lambda^\oc(X)\subseteq \calD^\an_\Lambda(X)$ denote the subcategory of overconvergent sheaves.
                \item  For $X \in \PreStk$ and $? \in \{\diamond,\dagger,\diamondsuit\}$, we write $\calD^{?,\oc}_{\Lambda}(X) \subset \calD^{?}_{\Lambda}(X)$ for the full subcategory of $\calD^{?}_{\Lambda}(X) \simeq \calD^{\an}_{\Lambda}(X^{?})$ (here we have used \Cref{prop: compatibilityofKanextensionswithDan}) whose objects are overconvergent.  
		\end{enumerate}
	\end{definition}
	\begin{remark}
		\label{ULA-is-overconvergent}
		Given $[f:X\to \ast]\in (\AnStk_{v},E_\fdcs)$ and a sheaf $A\in \calD^\an_\Lambda(X)$, the condition that $A$ is overconvergent is necessary for it to be $f$-suave (or in other words $f$-ULA in terminology of \cite{FS21}, see \cite[Definition IV.2.1, IV.2.4]{FS21}).
		In other words, for all such $f$, we have an inclusion
		\[\on{Suave}_{f,\calD^\an_\Lambda}(X)\subseteq \calD^\oc_\Lambda(X).\]
	\end{remark}
        \begin{remark}{\label{rem: naturaltransformationsofsheaftheoriesRespectOverconvergence}}
        We note that the overconvergent condition is tautologically stable under pullback along maps of prestacks. 
	In particular, for $X \in \PreStk$, the natural transformations \eqref{eqn: naturaltransformsofsheaftheories} induce natural transformations
        \[ \calD_{\Lambda}^{\diamondsuit,\oc}(X) \xRightarrow{b^*} \calD_{\Lambda}^{\dagger,\oc}(X) \xRightarrow{a^*} \calD_{\Lambda}^{\diamond,\oc}(X).  \]
        \end{remark}

	\begin{proposition}
\label{rem: overconvergentstableundercolimits2}
		For all $X\in \AnStk_v$, the inclusion $\calD^\oc_\Lambda(X)\subseteq \calD^\an_\Lambda(X)$ is stable under colimits.  	
	\end{proposition}
	\begin{proof}
	Since pullback commutes with colimits, we may assume $X=\Spa(C,C^+)$.
Let $\pi:X\to \ast$ denote the structure map.
Since $\pi^*$ also commutes with colimits, the claim follows.
	\end{proof}
	
	\begin{proposition}{\label{prop: overconvergentisvsheaf}}
		The following statements hold.
		\begin{enumerate}
\item The condition of being overconvergent can be verified v-locally. 
			\item The functor $\calD^{\oc}_\Lambda$ is a v-sheaf on $\Perff$.
			\item $\calD_\Lambda^{\diamond,\oc}$ is a schematic v-sheaf on $\PSchf$. 
			\item $\calD_\Lambda^{\dagger,\oc}$ is an arc-sheaf on $\PSchf$. 
		\end{enumerate}
	\end{proposition}
	\begin{proof}
		By \Cref{prop: descentproperties} (1)-(2), all of the claims reduce formally to the first claim, which we now justify. 
		Let $X=\Spa(C_1,C_1^+)$, $Y=\Spa(C_2,C_2^+)$, $f:Y\to X$ a v-cover and $A\in \calD_\et(X)$. 
		Since overconvergence is a condition on geometric points, one can easily reduce to showing that $f^*A$ is overconvergent if and only if $A$ is overconvergent. 
		This case follows directly from \cite[Theorem 1.13.(iii)]{Sch17}.
	\end{proof}

	\begin{proposition}
		For any $X\in \PreStk$, the analytification map $c^*_X:\calD^{\sch}_{\Lambda}(X) \to \calD_{\Lambda}^{\diamondsuit}(X)$ factors through $\calD_{\Lambda}^{\diamondsuit,\oc}(X)\subseteq \calD_{\Lambda}^{\diamondsuit}(X)$. 	
		In other words, we have a factorization of fully faithful maps
		\[\calD^\sch_\Lambda \xRightarrow{c^*} \calD^{\diamondsuit,\oc}_\Lambda\subseteq \calD^\diamondsuit_\Lambda.\]
	\end{proposition}
	\begin{proof}
		Since the functors $\calD^\sch_\Lambda$, $\calD^\oc_\Lambda$ and $\calD^\diamondsuit_\Lambda$ are all v-sheaves by \Cref{prop: overconvergentisvsheaf}, it suffices to show the factorization for affine schemes.
		Fix $X=\Spec R$ and $A\in \calD_\et(X)$, we wish to show that $c_X^*A\in \calD^\oc_\et(X^\diamondsuit)$.
		Observe that analytification and pullback are t-exact and commute with truncation functors. 
		Moreover, $\calD^\sch_\Lambda(X)$ and $\calD^\diamondsuit_\Lambda(X)$ are left-complete.
		Thus we may assume $A\in \calD^+_\et(X)$.
		Every object in $\calD^+_\et(X)$ is a colimit of objects in $\calD_{\on{cons}}(X)$ using \Cref{compact-generation-shv*} and the identification (\ref{eqn: identificationofboundedbelowcategories}), so by \Cref{rem: overconvergentstableundercolimits2} we may assume $A\in \calD_{\on{cons}}(X)$.
	Therefore, we may assume that $A=j_!\Lambda$ for a qcqs \'etale map $j:U\to X$.
	In this case $c^*_X A\simeq j_!^\diamondsuit \Lambda$, which one can easily verify is overconvergent, since the map $j^\diamondsuit:U^\diamondsuit\to X^\diamondsuit$ is partially proper. 
	\end{proof}

	\begin{lemma}
		\label{lemma open containing all rank 1-v2}
		If $j:X\hookrightarrow Y$ is a quasicompact monomorphism of  v-stacks, such that every rank $1$ geometric points of $Y$ lifts to $X$ (necessarily uniquely), then 
		\[j^*:\calD_\Lambda^\oc(Y)\to \calD_\Lambda^\oc(X)\]
		is an equivalence with inverse $j_*$.
	\end{lemma}
	\begin{proof}
		It suffices to show that for $A\in \calD_\Lambda^\oc(Y)$, the adjunction map $A\to j_*j^*A$ is an isomorphism.	
		This can be checked after pullback to geometric points of $Y$.
		Writing $A$ as the limit of its left truncations and since $j^*$ commutes with this Postinkov limit (being t-exact), we may assume that $A\in \calD^+_\et(Y)$.
		By quasicompact base change (see \cite[Corollary 16.10]{Sch17} and \cite[Proposition 4.5, Corollary 4.6]{HeyerMann}) and the assumption that the maximal generalizations defined by rank one points lift, we may assume that $Y=\Spa(C,C^+)$ and that $X=Y\times_{\ul{|Y|}}\ul{|Z|}$ for a pro-constructible generalizing subset of $\ul{|Y|}$ \cite[Corollary 10.6]{Sch17}.
		In other words, $X=\Spa(C,C'^+)$ for $C^+\subseteq C'^+\subseteq C$.
		In this case, $\calD_\Lambda^\oc(Y)$ is the essential image of ${\rm Mod}_\Lambda$ under $\pi^*$ for $\pi:Y\to \ast$ the structure map. 
		Since the composition 
		\[{\rm Mod}_\Lambda\simeq \calD^\an_\Lambda(\ast)\xrightarrow{\pi^*} \calD^\an_\Lambda(Y)\xrightarrow{j^*} \calD^\an_\Lambda(X) \]
		is fully faithful.
		The functor $j^*$ is fully faithful on $\calD^\oc_\Lambda(Y)$.
		Moreover, one can verify by hand, using \cite[Corollary 16.8]{Sch17}, that with this setup $j_*(\calD_\Lambda^\oc(X))\subseteq \calD_\Lambda^\oc(Y)$.
	\end{proof}
	
	\begin{remark}
		\label{argument affines-v2}	
		Combining \Cref{lemma open containing all rank 1-v2} and \Cref{lemma: geometricpropertiesofmapdaggertodiatodiamond} it follows that $a_X^*:\calD_{\Lambda}^\oc(X^\dagger)\to \calD_{\Lambda}^\oc(X^\diamond)$ is an equivalence whose inverse is $a_{X,*}$ for all affine schemes $X$.
		We will need a version of this for more general prestacks $X\in \PreStk$. 
		However, note that in this generality the maps $a_X:X^\diamond\to X^\dagger$ and $b_X:X^\dagger \to X^\diamondsuit$ are no longer monomorphisms, so one cannot proceed naively.
	\end{remark}

	\begin{proposition}
		\label{equivalence-aX}
		The map of functors on $\PreStk$ with values in $\LinCat_\Lambda$ 
		\[\calD^{\dagger,\oc}_\Lambda\xRightarrow{a^*} \calD^{\diamond,\oc}_\Lambda\]
		is an equivalence.
	\end{proposition}
	\begin{proof}
		By \Cref{prop: overconvergentisvsheaf}, both functors are schematic v-sheaves, so it suffices to show the equivalence for affine schemes. 
		This is the content of \Cref{argument affines-v2}.
	\end{proof}

	We will need a relative version of this statement.
	For this fix $X\in \PreStk$ and a map $X^\dagger \to X_0$ for $X_0\in \AnStk_v$, and fix a map $Y_0\to X_0$.
	We can construct two functors 
	\[\calD^{\diamond,\oc}_{\Lambda,Y_0}, \, \calD^{\dagger,\oc}_{\Lambda,Y_0}: \PreStk_{/X} \to \LinCat_\Lambda\]
		with formula 
		\[Z\mapsto \calD^{\an,\oc}_\Lambda(Z^\diamond\times_{X_0} Y_0) \text{ and } Z\mapsto \calD^{\an,\oc}_\Lambda(Z^\dagger \times_{X_0} Y_0).\] 

	\begin{proposition}
		\label{equivalence-aX-rel}
		With notation as above, the map of functors on $\PreStk_{/X}$ with values in $\LinCat_\Lambda$ 
		\[\calD^{\dagger,\oc}_{\Lambda,Y_0}\xRightarrow{a^*_{Y_0}} \calD^{\diamond,\oc}_{\Lambda,Y_0}\]
		is an equivalence.
	\end{proposition}
	\begin{proof}
	The same proof as in \Cref{equivalence-aX} applies. 	
	\end{proof}
	
	As established in Proposition \ref{equivalence-aX}, for general $X\in \PreStk$ we have an equivalence $a_X^*:\calD_\Lambda^\oc(X^\dagger)\to \calD_\Lambda^\oc(X^\diamond)$. 
	This says that for every $B\in \calD_\Lambda^\oc(X^\diamond)$ there is a unique (up to contractible choice) $A\in  \calD_\Lambda^\oc(X^\dagger)$ such that $a_X^*A\simeq B$. 
	Nevertheless, it is not at all clear that $A\simeq a_{X,*}B$. 
	Indeed, for this to hold one would have to know that $a_{X,*}B$ is overconvergent.
	As we explain below, in certain circumstances one can directly show that $A\simeq a_{X,*}B$.  
	
	\begin{lemma}
		\label{lemma: overconvergent replacement coorect adjoint}
		If $f:X\to Y$ is a qcqs map of v-stacks and $A\in \Detale^{+,\oc}(X)$, then $f_*A\in \Detale^{+,\oc}(Y)$.
		In particular, $((f^*)_{|\Detale^{+,\oc}(Y)},(f_*)_{|\Detale^{+,\oc}(X)})$ form an adjoint pair.
	\end{lemma}
	\begin{proof}
		By quasicompact base change (see \cite[Corollary 16.10]{Sch17} or \cite[Corollary 4.6.(ii)]{Mann2022NuclearSheaves}), we may assume that $Y=\Spa(C,C^+)$. 
		Let $Y_0=\Spa(C,O_C)$ and $X_0=Y_0\times_Y X$.
		Let $j_Y:Y_0\to Y$, $j_X:X_0\to X$ and $f_0:X_0\to Y_0$ denote the evident maps.
		We wish to show that 
		\[f_*A\to j_{Y,*}j_Y^*f_*A\]
		is an isomorphism.
		By \Cref{lemma open containing all rank 1-v2}, $A\simeq j_{X,*}j_X^*A$ and since $f\circ j_X=j_Y\circ f_0$, then $f_*A\simeq j_{Y,*}f_{0,*} j_X^*A$, and $f_{0,*} j_X^*A\simeq j_Y^*f_*A$ again by quasicompact base change.
		This finishes showing the claim.
	\end{proof}

	\begin{lemma}
		\label{qcqs:is-easy}
		Let $[f:X\to Y]\in \SchStk_{\et}$ be qcqs, then $f^\diamond$ and $f^\dagger$ are qcqs. 
	\end{lemma}
	\begin{proof}
		Since being qcqs is v-local on the target by \cite[Proposition 10.11 (o)]{Sch17}, we may assume that $Y=\Spec R$ .
		The argument for quasiseparatedness follows from that of quasicompactness applied to the diagonal.
		Since $f$ is quasicompact there is a v-surjective cover $\Spec A\to X$, so it suffices to show that $(\Spec A)^\dagger \to (\Spec R)^\dagger$ and $(\Spec A)^\diamond\to (\Spec R)^\diamond$ are quasicompact.
This follows from the proof of \cite[Lemma 2.26]{Gle24}.
	\end{proof}
	
	\begin{corollary}
		\label{overconvergent-replacement}
		Let $X\in \PreStk$ such that its v-sheafification is qcqs, then $X^\diamond \to X^\dagger$ is qcqs. 
		In particular, we have an equivalence of endofunctors on $\calD^{\dagger, \oc,+}_\Lambda(X)$ 
		\[\id\simeq a_{X,*}a_X^*,\]
		for any such $X$.
	\end{corollary}
	\begin{proof}
		Recall that in $(\SchStk)_{\on{v}}$ the final object is quasiseparated. 
		In particular, $X\in (\SchStk)_{\on{v}}$ is qcqs if and only if $X\to \ast$ is qcqs. 
		Applying \Cref{qcqs:is-easy} to the structure map $X\to \ast$, we obtain that $X^\diamond$ and $X^\dagger$ are both qcqs over $\ast$. 
		In particular, any map between them is also qcqs.
		The second claim follows from \Cref{qcqs:is-easy}, \Cref{lemma: overconvergent replacement coorect adjoint} and \Cref{equivalence-aX}.  
	\end{proof}
	
	As in \Cref{equivalence-aX-rel}, \Cref{overconvergent-replacement} has a relative version.

	\begin{corollary}
		\label{overconvergent-replacement-rel}
		With the setup as in \Cref{equivalence-aX-rel}. 
		Let $Z\in \PreStk_X$ and suppose that v-sheafification of $Z$ is qcqs, then we have an equivalence of endofunctors on $\calD^{\dagger, \oc,+}_{\Lambda, Y_0}(Z)$ 
	\[\id\simeq a_{Z\times_{X_0} Y_0,*}a_{Z\times_{X_0} Y_0}^*.\]
	\end{corollary}

	\subsection{Pro-unipotent base change}
	\label{ss: pro-unip base change}
	Later on we will work with maps of spaces that are infinite dimensional, but whose infinite dimensional part is essentially contractible. 
	This is the subject of placid geometry.
	In \cite[Proposition 10.32]{Zhu25}, we find analogues of smooth base change in that context. 
	In this section we exploit the techniques discussed in \S \ref{Computing d} and \S \ref{overconv repl} to prove an analogue in the context of v-sheaves and analytification.
 We first have the following basic result.

    \begin{lemma}
	    Let $f: X \ra Y $ be a $\calD_{\Lambda}^{\mathrm{an}}$-smooth map in $\AnStk_{v}$.
Then $f$ is unipotent if and only if it is unipotent after pulling back to every geometric point. 
    \end{lemma}
    \begin{proof}
	    Since $f$ is $\calD_{\Lambda}^{\mathrm{an}}$-smooth, it follows that we have a natural transformation 
    \[ f^{!}(-) \simeq f^{!}(\Lambda) \otimes f^{*}(-), \]
    and that $f^{!}(\Lambda)$ is invertible. 
    This shows that the fully faithfullness of $f^{*}$ is equivalent to the fully faithfulness of $f^{!}$, which is in turn equivalent to showing that the adjunction map
    \[ f_{!}f^{!}(A) \ra A \]
    is an isomorphism for all $A \in \calD_{\Lambda}^{\an}(Y)$. 
    This can be checked after pulling back to any geometric point $x:\Spa(C,C^{+}) \ra Y$ since this defines a conservative family.  
    Since $f$ is $\calD_{\Lambda}^{\mathrm{an}}$-smooth, and in particular $\calD_{\Lambda}^{\mathrm{an}}$-suave, it follows that $f^{!}$ commutes with $*$-pullback by \cite[Lemma~4.5.13 (i)]{HeyerMann}, and similarly $f_{!}$ commutes with $*$-pullback by proper base-change. 
    In other words, $f^*$ is fully faithful if and only if $f_{x,!}f^!_x x^*A\to x^*A$ is an isomorphism for every geometric point of $x: \Spa(C,C^{+}) \ra Y$ where $f_{x}$ denotes the pullback of $f$ along $x$, as desired.
    \end{proof}

In the algebraic context we have a stronger statement.

	\begin{lemma}{\label{lemma: CheckUnipotenceonPointsschematic}}
		Let $[f \colon X \to Y] \in (\PreStk,E^\rep_\pfp)$ be $\calD^{\mathrm{sch}}_{\Lambda}$-suave.
The following hold.
		\begin{enumerate}
			\item $f$ is $\calD_{\Lambda}^{\mathrm{sch}}$-smooth if it is so after pullback to any geometric point of $Y$.
			\item $f$ is $\calD_{\Lambda}^{\mathrm{sch}}$-unipotent if it is so after pullback to any geometric point of $Y$.
		\end{enumerate}

	\end{lemma}
	\begin{proof}
        As $f$ is $\calD_{\Lambda}^{\mathrm{sch}}$-suave, we need to show that the conditions that $f^{!}(\Lambda)$ is invertible and that $f^{*}$ is fully faithful hold if and only if they hold after pulling back to any geometric point. 
	Since the condition of being $\calD_{\Lambda}^{\mathrm{sch}}$-suave is stable under pullback by \cite[Lemma~4.5.9 (i)]{HeyerMann}, this can be checked after pulling back to any $Y \in \PSch^{\aff}$, which using the fact that upper $!$ and upper $*$ commute for $\calD_{\Lambda}^{\mathrm{sch}}$-suave maps by \cite[Lemma~4.5.13 (i)]{HeyerMann}, reduces us to the case that $[f: X \ra Y] \in (\PSch, E_{\pfp})$ and $Y\in \PSchf$. 
	More precisely, for $S \ra Y$ with  $S \in \PSch^{\aff}$ we write $f_{S}: X \times_{S} Y \ra S$ for the base-change.
We have a commutative diagram
        \[ 
        \begin{tikzcd}
        \calD_{\Lambda}^{\sch}(X) \arrow[r] & \lim_{S \ra Y} \calD_{\Lambda}^{\mathrm{sch}} (X \times_{Y} S) & \\
        \calD_{\Lambda}^{\sch}(Y) \arrow[u,"f^{?}"] \arrow[r] & \lim_{S \ra Y} \calD_{\Lambda}^{\sch}(S), \arrow[u,"\lim_{S \ra Y} f_{S}^{?}"] & 
        \end{tikzcd}
        \]
        for $? \in \{!,*\}$. 
	Here the horizontal arrows and the transition maps in the limit are given by $*$-pullbacks.
The commutativity for $?= \text{ !}$ and the well-definedness of $\lim_{S \ra Y} f_{S}^{!}$ follows from the $\calD^{\sch}_{\Lambda}$-suaveness. 
	Moreover, the horizontal arrows are equivalences by definition of $\calD_{\Lambda}^{\mathrm{sch}}$ as a right Kan extension.
        
       Using the diagram when $? = \text{ !}$, one can show invertibility of $f^{!}(\Lambda)$ by showing invertibility of $f_{S}^{!}(\Lambda)$ for varying $S$ using that the $*$-pullback functors are symmetric monoidal. 
       Using the diagram when $? = *$, one can show fully faihtfullness of $f^{*}$ by showing fully faithfullness of $f_{S}^{*}$ for varying $S$, since limits of fully faithful maps are fully faithful. 
     Therefore, we may assume that $Y\in \PSch^\aff$ and $f: X \ra Y \in E_\pfp$. 

        Now it follows from \cite[Lemma~10.44]{Zhu25} that the $\calD_{\Lambda}^{\mathrm{sch}}$-smoothness can be checked after pulling back to each geometric point.
Similarly, if $f: X \ra Y$ is $\calD_{\Lambda}^{\mathrm{sch}}$-smooth, and by \cite[Lemma~10.48]{Zhu25}, the unipotence is equivalent to showing that the map
        \[ f_{!}f^{!}(\Lambda) \ra \Lambda \]
        is an isomorphism.
This can be checked on geometric points, since $f_{!}$ and $f^{!}$ commute with $*$-pullback, as discussed above.  
	\end{proof}

Unipotent maps are particular instances of smooth maps.
In particular, $*$-pullback commutes with $*$-pushforward \cite[Lemma~4.5.13 (i)]{HeyerMann} in this case.
We show that, under qcqs assumptions, pro-unipotent maps still have this property.

	\begin{lemma}
		\label{pro-unip-base change}
		Let $[f:Z\to Y]$ be a map of analytic prestacks in $E^\rep_\fdcss$.
		Suppose that $Y$ is either an Artin pro\'etale-stack or an Artin v-stack (see \Cref{rem: WHatArtinandFineMeansintheAnalyticContext}).
		Let $g:X\to Y$ be a map such that $X=\varprojlim_{i\in I} X_i$, where each $[g_i:X_i\to Y]\in E^\rep_\fdcss$ is $\calD^\uop_\Lambda$-unipotent (resp. $\calD^\an_\Lambda$-unipotent, see \Cref{everything looks normal after rest}) and qcqs. 
		Consider the Cartesian diagram
		\begin{center}
			\begin{tikzcd}
				Z\times_Y X\arrow{r}{f'} \arrow{d}{g'}  & X \arrow{d}{g} \\
				Z \arrow{r}{f} & Y, 
			\end{tikzcd}
		\end{center}
		then the natural map $g^*f_*\to f'_*g'^*$ is an isomorphism. 
	\end{lemma}
	\begin{proof}
		We must show $g^*f_* A \to f'_*g'^* A$ for all $A\in \calD^\uop_\Lambda(Z)$ (resp. $A\in \calD^\an_\Lambda(Z)$). 
		We start with some reduction steps.
		Since $Y$ is Artin it admits a fdcss and $\calD^\uop_\Lambda$-smooth cover (resp. $\calD^\an_\Lambda$-smooth cover) by a locally spatial diamond. 
		Using smooth base change, we may replace $Y$ by a disjoint union of spatial diamonds. 
	After this replacement, we may assume that $X$, $Y$ and $X_i$ are spatial diamonds, and that $Z$ is a locally spatial diamond. 
	At this point we notice that all of the geometric objects involved, $X$, $Y$, $X_i$ and $Z$ are in $\on{sLocSptl}$, and that when we restrict to this category $\calD^\an_\Lambda$ and $\calD^\uop_\Lambda$ agree (see \Cref{everything looks normal after rest}).

		Using a Postnikov-tower argument, we may assume $A$ is bounded below. 
		For a geometric point $\overline{x}\to X$, we show that the map
		\[(g^*f_* A)_{\overline{x}} \to (f'_*g'^* A)_{\overline{x}}\]
    	on stalks is an isomorphism. 
		Both can be rewritten as colimits over qcqs \'etale neighborhoods $\overline{x}\to V\to X$. 
		So it suffices to show $\Gamma(V,g^*f_* A)\to \Gamma(V,f'_*g'^* A)$ is an isomorphism for all $V\to X$ \'etale.
		By \cite[Proposition~11.23]{Sch17}, we may assume $V=X\times_{X_i} V_i$ for some $i$ and for a qcqs \'etale map $V_i\to X_i$. 

		Since $V_i\to Y$ is smooth and satisfies smooth base change we may replace  $Y$ by $V_i$, $Z$ by $Z\times_Y V_i$, and $X$ by $V$ and the role of $X_j$ by $X_j\times_{X_i} V_i$.
	In other words, after reorganizing the notation we see that it suffices to show that
	\begin{equation}
\label{computation: key}
		\Gamma(X,g^*f_* A)\to \Gamma(X,f'_*g'^* A)
	\end{equation}
	is an isomorphism. 

	Let $Z_i=Z\times_{Y} X_i$ and let $g'_i$ (resp. $f'_i$) denote the pullback of $g_i$ (resp. $f$) along $f$ (resp. $g_i$).  
		Since $Y$ is spatial and the map $g$ is a limit of qcqs maps, we may apply \cite[Proposition 14.9]{Sch17} to the left-hand side of \eqref{computation: key} to compute: 

	\begin{align*}
		\allowdisplaybreaks
		\Gamma(X, g^*f_* A)& \simeq  \varinjlim_{i\in I} \Gamma(X_i, g_i^*f_*A)	 \\
				   & \simeq \varinjlim_{i\in I} \Gamma(X_i, f'_{i,*} g'^*_i A)	 \\
				   & \simeq \varinjlim_{i\in I} \Gamma(Z_i, g'^*_i A).
	\end{align*}

	On the other hand, we always have $\Gamma(X,f'_* g'^* A)\simeq \Gamma(Z, g'^* A)$.
	So it suffices to show that 
	\begin{equation}
		\label{more computations yay}
\varinjlim_{i\in I} \Gamma(Z_i, g'^*_{i} A)\to \Gamma(Z\times_Y X, g'^* A)
	\end{equation}
	is an isomorphism. 
We cannot apply \cite[Proposition 14.9]{Sch17} directly to the tower $Z_i$, since these are only locally spatial diamonds. 
Notice also that the reasoning above has not used the unipotence hypothesis which is certainly necessary.

We argue as follows.
Note that the unipotence assumption shows that the left-hand side of \eqref{more computations yay} reduces to $\Gamma(Z,A)$. 
In particular, it suffices to see that the natural map $\Gamma(Z,A)\to \Gamma(Z\times_Y X, g'^* A)$ is an isomorphism. 
Using smooth base change and hyperdescent for $\Gamma(-,-)$ in the analytic topology, this isomorphism can be shown on an open hypercover of the form $Z_\bullet\to Z$.
Choosing $Z_n$ of the form $\coprod_{i\in \calI_n} S_i$ with $S_i$ spatial, we may reduce to the case in which $Z=S_i$ is a spatial diamond.
In this case, \cite[Proposition 14.9]{Sch17} already applies.
	\end{proof}

	\section{Placid geometry}
	\label{Dictionary-with-Zhu}

    \noindent
    \textbf{Soft preamble to the section:}
	Many of the geometric objects that naturally arise when one studies the schematic local Langlands category are not finite dimensional, but they often admit a presentation as a highly organized (co)limit of finite dimensional pieces.  
	Experience with this type of geometric objects has lead mathematicians working in the geometric Langlands program to consider placid geometry and co-sheaf theories instead of the standard 6-functor formalisms of \'etale sheaves (see \cite{GaitsgoryLurieWeilsConjectureI} and \cite{perverse-sheaves-infinite-dimensional-stacks}).
	In our case, the schematic local Langlands category is obtained as the category of co-sheaves on a sind-placid stack \cite[Definition 10.157]{Zhu25} in Hemo--Zhu's terminology.

	In this section, we recall generalities of \'etale co-sheaves as discussed by Zhu, as well as the main class of geometric objects of schematic nature for which we study this co-sheaf theory. 
	More crucially, we learn one way in which we can analytify in this placid setup.
	\\

    \noindent
    \textbf{Technical preamble to the section:}
    \begin{itemize}
	    \item[\S \ref{Section: Sheaf and Cosheaf Theories}] This part is mostly expository.
		    Zhu's work considers the interplay between two sheaf theories, namely $\Shv^*$ and $\Shv^!$.
		    In this section, we run the bureaucracy to show that $\Shv^*$ (as defined by Zhu) agrees with the 6-functor formalism $\fShv$ that we work with throughout \S \ref{s:Analyt for stacks section}. 
		    We also explain the paradigm of co-sheaves $\Shv^!$ and collect some useful facts from \cite{Zhu25}. 
	    \item[\S \ref{ss: placid stacks}] This part is purely expository. 
		    We recall the formalism of placid stacks, ind-placid stacks and sind-placid stacks.
		    We also collect useful facts from \cite{Zhu25}.
	    \item[\S \ref{sec:cosheaves_and_analytification}] This subsection underlies the technical pillar to the construction of $\pitch$. 
		    Indeed, the sheaf theory that Zhu relies on, i.e., $\Shv^!$, is not very compatible with the analytification functor constructed in \S \ref{s:Analyt for stacks section}. 
		    To address this, we introduce the formalism of sind-$\dagger$-correspondences (see \Cref{sind-dagger-correspondence}). 
		    This formalism makes precise the idea that for a sind-placid stack, like the stack of isocrystals, one should not analytify directly. 
		    Instead, one should use a suitable sind-presentation in order to define analytification, and one should also justify that the presentation did not play a role.
		    This formalism in particular allows us to show that our functor $\pitch$ does not depend on the choice of an auxiliary parahoric.
	    \item[\S \ref{ss: some computations with placid}] Although we do not develop a general theory for the functoriality of sind-$\dagger$-correspondences or a general theory of placid analytic geometry, in this section we collect several useful statements along those lines which are enough for our purposes.   
		    For example, \Cref{base change formulas} reaps the benefits of \S \ref{Computing d}, \ref{overconv repl}, \ref{ss: pro-unip base change} and shows a general base change result that will in particular apply to the stack of analytic shtukas $\Sht^\an_\calG$.
		    These computations are key to showing that $\pitch$ is semi-orthogonal.
    \end{itemize}
	
	\subsection{Co-sheaf theory}{\label{Section: Sheaf and Cosheaf Theories}}
	Recall that in \Cref{subsection algeb 6-ff} we defined $\calD^\sch_\Lambda(-)$ on $\SchStk_v$ as the only hypercomplete v-sheaf
	\begin{equation}
		\label{new schematic equation}
	\calD_\Lambda^\sch(-):(\SchStk_v)^\op\to \LinCat_\Lambda
	\end{equation}
    whose value on schemes $Z$ is $\calD_\et(Z)$. 
	Here the transition maps in \eqref{new schematic equation} are given by $*$-pullback. 
	We can say colloquially that $\calD^\sch_\Lambda$ is ``$*$-glued''.
	Moreover, if $Z\in \PSch^\qcqs$ then $\calD_\et(Z) = \calD_{\Lambda}^{\sch}(Z)$ is the left-completion of the ind-completion of $\calD_\cons(Z)$, as explained in \S \ref{subsection algeb 6-ff}. 

	The theory of co-sheaves considered in \cite{Zhu25}, similarly to $\calD^\sch_\Lambda$, can also be obtained from $\calD_\cons(-)$ after performing categorical constructions.
	There are two main differences in Zhu's approach. 
	\begin{enumerate}
		\item In contrast to $\calD^\sch_\Lambda$, it follows directly from the definition that the value of Zhu's functor on qcqs schemes is a dualizable category (in fact, compactly generated, cf. \Cref{compact-generation-shv*}). 
		\item The co-sheaf theory of Zhu is, informally speaking, $!$-glued. 
			This gives it the colloquial name of co-sheaf theory. 
	\end{enumerate}
	
	A disadvantage of working with co-sheaves is that one has to construct and study de novo the six operations. 
	Fortunately for us, Zhu has given a thorough account of the theory (see \cite[\S 10]{Zhu25}), some of which we recall below.
	Following Zhu, we will consider our basic building blocks to be algebraic spaces instead of schemes when we discuss the co-sheaf theory.
	For the rest of the section, in contrast to the previous sections, when we say that a map is representable we mean that it is representable in algebraic spaces. 
	\begin{definition}
		\label{first defi Zhu setup}
		Consider the functors 
		\[\calF: (\AlgSp^{\pfp})^\op\to \CAlg(\LinCat_\Lambda)\]
		on perfectely finitely presented algebraic spaces with $\calF\in \{\Shv^*_c,\Shv^*,\Shv_c^!,\Shv^!\}$ given by the rules 
		\begin{align*}
			\Shv^*_c(X)&:=\calD_\cons(X) \\
		\Shv^*(X)&:=\Ind\calD_\cons(X)\\
		\Shv^!_c(X)&:=\calD^\op_\cons(X)\\
		\Shv^!(X)&:=\Ind\calD^\op_\cons(X)
		\end{align*}
		with transition maps given by $f^*$ and $\Ind(f^{*})$ in the former two theories, respectively, and $f^{*,\op}$ and $\Ind(f^{*,\op})$ in the latter two theories, respectively.
Here the superscript $(-)^{\op}$ denotes the opposite category.
	\end{definition}

	\begin{definition}
		\label{second defi zhu setup}
		For $\calF\in \{\Shv^*_c,\Shv^*,\Shv_c^!,\Shv^!\}$ we extend along $(\AlgSp^\pfp)^\op\subseteq (\AlgSp^\qcqs)^\op$ by taking its left Kan extension. 
		We keep the notation of $\calF$ for the functors obtained in this way.
		We will call the functors $\calF\in \{\Shv^*_c,\Shv^*\}$ the \emph{sheaf setup} and we will call the functors $\calF\in \{\Shv_c^!,\Shv^!\}$ the \emph{co-sheaf setup}.
	\end{definition}

	Our \Cref{first defi Zhu setup} and \Cref{second defi zhu setup} were obtained from Zhu's definition by gathering \cite[Remark 10.28, Remark 10.12, Equation (10.6)]{Zhu25}.

	\begin{remark}
		If $[f: X \ra Y] \in \AlgSp^\pfp$ one can use Verdier duality on $X$ and $Y$ to functorially identify $(\Shv_c^!(X),f^{*,\op})$ with $(\Shv^*_c(X),f^!)$. 
		As in the following commutative diagram
		\begin{center}
		\begin{tikzcd}
			\Shv_c^!(Y)	 \arrow{r}{\bbD_Y} \arrow{d}{f^{*,\op}}  & \Shv_c^*(Y) \arrow{d}{f^!}  \\
			\Shv_c^!(X) \arrow{r}{\bbD_X} & \Shv_c^*(X)
		\end{tikzcd}
		\end{center}
		This justifies \Cref{cosheaf convention} below. 
	\end{remark}
	
	\begin{convention}
		\label{cosheaf convention}
		If $[f:X\to Y]\in \AlgSp^\qcqs$ we let $\ul{f^!}:=f^{*,\op}$ denote the transition functors 
		\[\ul{f^!}:\Shv_c^!(Y)\to \Shv^!_c(X)\]
		and
		\[\ul{f^!}:\Shv^!(Y)\to \Shv^!(X).\]
		Since throughout the text we will be forced to mix sheaf and co-sheaf theories, we will underline all the functors that appear in a co-sheaf setup for clarity. 
		We leave the functors associated with sheaf setups and the 6-functor formalisms considered in \S \ref{sec: analytificiation and 6-functor formalisms} and \S \ref{s:Analyt for stacks section} unadorned. 
	\end{convention}

	\begin{definition}
		We extend the domain of $\calF$ of \Cref{first defi Zhu setup} along $\AlgSp^\qcqs\subseteq \PreStk$ by taking its right Kan extension.
	This gives us functors
	\[\calF: \PreStk^\op\to \CAlg(\LinCat^\sm_\Lambda) \text{ (resp. }\calF: \PreStk^\op\to \CAlg(\LinCat_\Lambda))\]
	in the cases $\calF\in \{\Shv^*_c, \Shv^!_c\}$ and $\calF\in \{\Shv^*, \Shv^!\}$ respectively.
	For \(\calF\in\{\Shv^*_c,\Shv^*\}\), the transition maps are denoted  by $f^*$.
In the case of \(\calF\in\{\Shv_c^!,\Shv^!\}\), they are denoted $\ul{f^!}$. 
	For any $X\in \PreStk$ we still have fully faithful embeddings $\Shv^!_c(X)\subseteq \Shv^!(X)$ and $\Shv^*_c(X)\subseteq \Shv^*(X)$. 
	Any object in the essential image is called \emph{constructible}.
    \end{definition}

	\begin{remark}
		We note that for general prestacks it might not hold that $\Shv^!(X)$ (resp. $\Shv^*(X)$) is compactly generated, and even if it is compactly generated it might not be true that $\Shv^!(X)^\omega\simeq \Shv^!_c(X)$ (resp. $\Shv^\ast(X)^\omega\simeq \Shv^\ast_c(X)$), see \cite[Example 10.134]{Zhu25}.
	\end{remark}

	\begin{remark}{\label{rem: ShvcisopShvcstar}}
		\label{opposites}
		Since $(-)^\op$ is an involution in $\LinCat^\sm_\Lambda$, it commutes with the formation of limits and colimits. 
		It follows that, for all $X\in \PreStk$, there is a functorial identification 
		\[\Shv^!_c(X)\simeq \Shv^*_c(X)^\op.\]
	\end{remark}

	We can derive the following formal consequences from \Cref{thm-Hansen-Scholze}. 
	\begin{proposition}{\label{prop: DconsmapsintoDet}}
		The following hold.
		\begin{enumerate}
			\item $\calD_\cons(-)\simeq \Shv_c^*(-)$ in $\Fun(\PreStk^\op,\CAlg(\LinCat_\Lambda^\sm))$. 
			\item $\calD_\cons(-)^\op\simeq \Shv_c^!(-)$ in $\Fun(\PSch^\qcqs,\CAlg(\LinCat_\Lambda^\sm))$. 
			\item Let $\calF\in \{\Shv_c^*(-), \Shv_c^!(-)\}$, then $\calF$ satisfies descent for the schematic v-topology on $\SchStk_v$. 
				In particular, $\calF:\PreStk\to \CAlg(\LinCat_\Lambda^\sm)$ factors through sheafification $\PreStk\to \SchStk_v$ composed with
			the restricted functor $\calF:\SchStk_v\to \LinCat^\sm$. 
		\end{enumerate}
	\end{proposition}
	\begin{proof}
		For the first claim, we observe that $\calD_\cons(-)\simeq \Shv_c^*(-)$ as presheaves on $\PSch^\qcqs$ with values in $\CAlg(\LinCat_\Lambda^\sm)$. 
	Indeed, this follows from the definition of $\Shv_c^*(X)$ as a left Kan extension and the finitary part of \Cref{thm-Hansen-Scholze}.
	The passage to $\Fun((\PreStk)^\op,\CAlg(\LinCat_\Lambda^\sm))$ is formal since both functors are obtained as right Kan extensions from $\PSchf$.
	Indeed, in the case of $\calD_\cons$ this is by definition.
In the case of $\Shv^*_c$, it follows from the fact that it satisfies \'etale descent so the Kan extensions from qcqs perfect algebraic spaces and from qcqs perfect schemes agree.
		The second claim follows from the first claim and \Cref{opposites}.
		Finally, the last claim follows from \cite[Proposition 1.3.1.7]{SAG} and \Cref{thm: Dconsarcsheaf}, by similar reasoning to the proof of \Cref{prop: descentproperties}.
	\end{proof}

	\begin{remark}
		\label{fast observations zhu setup}
		We note that the functor that Zhu considers, $\Shv^*$, can take objects in $\AlgSp^\qcqs$ as input, while a priori the functor $\fShv$ that we consider only takes as input values objects in $\PSch^\qcqs$.
		Since both functors satisfy \'etale descent it follows formally that $\Shv^*$ is the right Kan extension of $\fShv$ along the embedding $(\PSch^\qcqs)^\op\subseteq (\AlgSp^\qcqs)^\op$.
		\Cref{identification zhu vs standard} shows a stronger claim.
	\end{remark}

	\begin{proposition}
		\label{identification zhu vs standard}
		We have identifications $\Shv^*(-)\simeq \fShv(-)$ as functors on $\PreStk^\op$ with values in $\CAlg(\LinCat_\Lambda)$.
	\end{proposition}
	\begin{proof}
		On $\PSch^\qcqs$ both functors are the ind-completion of $\Shv_c^*(-)$ and $\calD_\cons(-)$ respectively.
		So the restrictions to functors on $(\PSch^\qcqs)^\op$ are equivalent by \Cref{prop: DconsmapsintoDet}.
	The passage to $\Fun((\PreStk)^\op,\CAlg(\LinCat_\Lambda))$ is formal since both functors are obtained as right Kan extensions from $\PSch^{\qcqs}$.
	\end{proof}

	Before we discuss the co-sheaf setup, we discuss a version of \Cref{identification zhu vs standard} for the correspondence category i.e., we compare $\Shv^*$ to the 6-functor formalism $\fShv$ constructed in \S \ref{subsection algeb 6-ff}. 
	We extracted the following \Cref{prop: algebraicspacesixfunctorformalism} from combining \cite[Remark 10.28, Theorem 10.15, Equation 10.26, Theorem 10.17, Proposition 10.25]{Zhu25}.
	We set some notation first, our geometric setup will be $\calC=\AlgSp^\qcqs$ and $E=E^\Alg_\pfp$ (i.e., separated maps of qcqs perfect algebraic spaces and perfectly finitely presented maps between them.
Note that we don't impose that the maps are relatively representable on schemes).
	We let $P^\Alg_\pfp\subseteq E^\Alg_\pfp$ denote the subclass of proper perfectly finitely presented maps and let $I^\Alg_\pfp\subseteq E^\Alg_\pfp$ denote the subclass of \'etale maps.

	\begin{proposition}{\label{prop: algebraicspacesixfunctorformalism}}
		The functor $\calF= \Shv^*_c$ (resp. $\calF=\Shv^*$) promotes to a 3-functor formalism
		\[\calF:\Corr(\AlgSp^\qcqs,E^\Alg_\pfp)\to \LinCat^\sm_\Lambda \text{ (resp. } \calF:\Corr(\AlgSp^\qcqs,E^\Alg_\pfp)\to \LinCat_\Lambda \text{)}\]
		satisfying
		\begin{enumerate}
			\item If $[f:X\to Y]\in E^\Alg_\pfp$ is \'etale, then $f_!\simeq f_\natural$ is left-adjoint to $f^{*}$.
			\item If $[f:X\to Y]\in E^\Alg_\pfp$ is proper, then $f_! \simeq f_{*}$ is right-adjoint to $f^{*}$.
			\item The restriction of $\calF$ to $\Corr(\AlgSp^\pfp,E_{\on{All}})$ is a 6-functor formalism.
			\item $\calF$ is sheafy for the \'etale topology.
		\end{enumerate}
	\end{proposition}

\begin{proposition}{\label{prop: comparisonofShv*andfShv}}
	The following statements hold. 
	\begin{enumerate}
		\item $(\AlgSp^\qcqs,E^\Alg_\pfp,P^\Alg_\pfp,I^\Alg_\pfp)$ is a Nagata setup, and the 6-functor formalism $\Shv^{*}$ is Nagata with respect to this setup, so that, 
        by \Cref{thm: LiuZhengUniquenesss}, we have an isomorphism 
        \[ \Shv^{*} \simeq \LiuZheng(\AlgSp^{\qcqs},E^{\Alg}_{\pfp},P^{\Alg}_{\pfp},I^{\Alg}_{\pfp},\Shv_0^{*}(-)). \]
            \item There is a unique equivalence between of $\fShv$ and $\Shv^*$ as 6-functor formalisms on $\Corr(\AlgSp^\qcqs,E_\pfp^\Alg)$, that extends the equivalence from \Cref{identification zhu vs standard}.
		    Here, $\fShv$ denotes the $6$-functor formalism of \Cref{important analytification prop} restricted along the inclusion 
		    \[\Corr(\AlgSp^\qcqs,E_\pfp^\Alg)\to \Corr(\PreStk,E^\sstd_\sch).\] 
		\item We have inclusion of classes of morphisms $E_\pfp\subseteq E^\Alg_\pfp\subseteq E^\an_\sch\subseteq  E^\sstd_\sch$. 
        \end{enumerate}
	In particular, given a map $[f:X\to Y]\in E^\Alg_\pfp$ the following statements hold.
        \begin{itemize}
		\item If $f\in P^\Alg_\pfp$ (i.e., it is proper pfp), then $f$ is $\fShv$-proper. 
		\item If $f\in I^\Alg_\pfp$ (i.e., it is an \'etale map of algebraic spaces), then $f$ is $\fShv$-\'etale.
	\end{itemize}
\end{proposition}
\begin{proof}
Part (1) follows from \Cref{prop: algebraicspacesixfunctorformalism}.
Indeed, the only claim that is not clear is that one can factor every pfp map as an open immersion and a pfp proper map. 
By writing $f: X \ra Y$ as the perfection of a finitely presented map and using that perfection is functorial, we reduce to showing that any finitely presented map of algebraic spaces factors as an open immersion and a proper finitely presented map.
This follows from \cite[Theorem~1.2.1]{ConradNagataCompactification}.
The fact that $\Shv^{*}$ is Nagata with respect to this Nagata setup now follows from \Cref{prop: algebraicspacesixfunctorformalism} (1)-(2).
For part (2), we start with the isomorphism 
\begin{equation}{\label{eqn: Shv*comparisonI}}
 \Shv^{*} \simeq \LiuZheng(\AlgSp^{\qcqs},E^{\Alg}_{\pfp},P_{\pfp}^{\Alg},I^{\Alg}_{\pfp},\Shv^{*}(-)), 
\end{equation}
given by part (1), and therefore, by \Cref{cor: LiuZhengRestriction}, it follows that we have an isomorphism
\[ \Shv^{*}|_{\Corr(\PSch^{\qcqs},E_{\pfp})} \simeq \LiuZheng(\PSch^{\on{qcqs}},E_{\pfp},I_{\pfp},P_{\pfp},\Shv^{*}(-)|_{\PSch^{\qcqs,\op}}). \]
Now, we recall, by definition, we have an identification
\[ \fShv|_{\Corr(\Sch^{\qcqs},E_{\finexp})} \simeq \LiuZheng(\Sch^{\qcqs},E_{\finexp},I_{\finexp},P_{\finexp},\fShv(-)), \]
where $I_{\finexp},P_{\finexp}$, and $E_{\finexp}$ are as in \Cref{cor: finiteexpansiongeometricsetup}. 
We note that we have proper inclusions $E_{\pfp} \subset E_{\finexp}$, $P_{\pfp} \subset P_{\finexp}$, and $I_{\pfp} \subset I_{\finexp}$, since perfection is an integral map. 
Therefore, we are again in a situation where we can apply \Cref{cor: LiuZhengRestriction}, and this tells us that we have an isomorphism
\begin{equation}{\label{eqn: Shv*comparisonII}}
 \fShv|_{\Corr(\PSch^{\qcqs},E_{\pfp})} \simeq \LiuZheng(\PSch^{\qcqs},E_{\pfp},I_{\pfp},P_{\pfp},\fShv(-)). 
\end{equation}

Since we have an equivalence $\fShv(-) \simeq \Shv^{*}(-)$ of functors on $\PSch^{\qcqs,\op}$ in light of \Cref{compact-generation-shv*}, by combining equations \ref{eqn: Shv*comparisonI} and \ref{eqn: Shv*comparisonII}, we can appeal to \Cref{thm: LiuZhengUniquenesss} to see that we have a unique equivalence
\[ \Shv^*|_{\Corr(\PSch^{\qcqs},E_{\pfp})} \simeq \fShv|_{\Corr(\PSch^{\qcqs},E_{\pfp})}. \]
We now claim that this upgrades to a unique equivalence 
\[ \Shv^* \simeq \fShv|_{\Corr(\AlgSp^{\qcqs},E^{\mathrm{Alg}}_{\pfp})}. \]
This follows from applying \Cref{cons: ExtensiontoRepresentableMaps} and \Cref{cons: ExtensiontoFineMapsI} once to obtain an equivalence
\[ \Shv^* |_{\Corr(\AlgSp^{\qcqs},E^\rep_{\pfp})}\simeq \fShv|_{\Corr(\AlgSp^{\qcqs},E^\rep_{\pfp})}, \]
and then applying \Cref{cons: ExtensiontoFineMaps}.
Indeed, every map of algebraic spaces in $E^{\mathrm{Alg}}_{\pfp}$ has the property that \'etale locally on the source it lies in $E_{\pfp}$ (i.e., it is relatively representable in pfp schemes). 
This \'etale local presentation also shows that $E^\Alg_\pfp\subseteq E^\fine_\pfp$, and consequently $E^\Alg_\pfp\subseteq E^\an_\sch$ (see \Cref{fine contains analytifiable}).
\end{proof}

\begin{proposition}
	\label{algebraic spaces are also resilient}
The following statements hold.	
\begin{enumerate}
	\item If $Y\in \AlgSp^\qcqs$, then $Y$ is resilient.
	\item If $f:Y\to X$ is a map of \'etale stacks such that the base change $Y\times_X \Spec A\to \Spec A$ is in $E^\Alg_\pfp$ for any $\Spec(A) \ra X$ and $X$ is resilient, then $Y$ is resilient. 
	\item If $Y,X\in \AlgSp^\qcqs$ and $Y\to X$ is a map in $E^\Alg_\pfp$, then $f^{\diamondsuit}: Y^{\diamondsuit}\to X^{\diamondsuit}$ is in $E^\rep_\fdcss$. 
\end{enumerate}
\end{proposition}
\begin{proof}
	For the first one, it suffices to see that $Y^\diamondsuit$ and $Y^{\diamondsuit_\uop}$ agree on strictly totally disconnected spaces since these are simultaneously a basis for the pro\'etale and v topologies.
	Let $W\to Y$ be an \'etale map in $E^{\on{rep}}_\pfp$, such that $W\in \PSch^\qcqs$.  
	Then, 
	\[Y^{\diamondsuit_\uop}=\on{co.eq.}(W^{\diamondsuit_\uop}\times_{Y^{\diamondsuit_\uop}} W^{\diamondsuit_\uop} \rightrightarrows W^{\diamondsuit_\uop}),\]
	and 
	\[Y^{\diamondsuit}=\on{co.eq.}(W^{\diamondsuit}\times_{Y^{\diamondsuit}} W^{\diamondsuit} \rightrightarrows W^{\diamondsuit}),\]
	where the colimits are taken in the category of pro\'etale-sheaves and v-sheaves respectively.
	We claim that if $\Spa(R,R^+)$ is a strictly totally disconnected space and $\Spa(R,R^+)\to Y^{\diamondsuit_\uop}$ (resp. $\Spa(R,R^+)\to Y^{\diamondsuit}$) is any point, it lifts to a point $\Spa(R,R^+)\to W^{\diamondsuit_\uop}$ (resp. $\Spa(R,R^+)\to W^{\diamondsuit}$).
	Since $W^\diamondsuit\simeq W^{\diamondsuit_\uop}$ and $W^{\diamondsuit}\times_{Y^{\diamondsuit}} W^{\diamondsuit}\simeq W^{\diamondsuit_\uop}\times_{Y^{\diamondsuit_\uop}} W^{\diamondsuit_\uop}$, this will show that $Y^{\diamondsuit}$ and $Y^{\diamondsuit_{\uop}}$ have the same  $(R,R^{+})$-points, as desired.  
	We prove the claim for the $\diamondsuit_\uop$ case, with the other one being analogous.
	Whenever $\Spa(R,R^+)\to Y^{\diamondsuit_\uop}$ factors through a map $f^{\diamondsuit_\uop}:(\Spec R)^{\diamondsuit_\uop}\to Y^{\diamondsuit_\uop}$, the claim follows from \Cref{lemma: pfpanalytifiestofdcs} since $(W\times_Y \Spec R)\to \Spec R$ is an \'etale map in $E_\pfp$, and $\Spa(R,R^+)$ splits every \'etale map (see \cite[Definition 7.15]{Sch17}).
	From the nature of sheafification, every map $\Spa(R,R^+)\to Y^{\diamondsuit_\uop}$ factors pro\'etale locally through $\Spec R$. 
	Then our claim follows from v (or pro\'etale) descent of \'etale maps (see \cite[Proposition 10.11]{Sch17}).

	The second point follows from the first one and from \Cref{enough-resilients}.(4). 	

	For the third point, we fix a map $\calS \to X^\diamondsuit$ where $\calS$ is a locally spatial diamond, we must show that $Y^\diamondsuit\times_{X^\diamondsuit}\calS$ is a locally spatial diamond, and that $Y^\diamondsuit\times_{X^\diamondsuit}\calS\to \calS$ is fdcss.
	We may use the identity $f^{\diamondsuit_\uop}=f^\diamondsuit$ given by Point (1), \Cref{uop properties} and \Cref{pro-etale locality} to reduce to the case that $\calS=\Spa(R,R^+)$ for $\Spa(R,R^+)$ a strictly totally disconnected space, $X=\Spec R$, $Y$ is an algebraic space perfectly finitely presented over $\Spec R$.

	It is not hard to see that $Y^\diamondsuit$ is separated over $X^\diamondsuit$, since $\diamondsuit$ preserves closed immersions.
	We fix an \'etale atlas $W\to Y$ with $W\to X$ in $E_\pfp$, as above.
	Since the map $W\to X$ is in $E_\pfp$, then by \Cref{lemma: pfpanalytifiestofdcs} $W^\diamondsuit\to X^\diamondsuit$ is in $E^{\on{rep}}_\fdcss$. 
	In particular, $W^\diamondsuit\times_{X^\diamondsuit} \Spa(R,R^+)$ is a locally spatial diamond. 
	Using \Cref{uop properties} and the \'etale map $W^\diamondsuit \to Y^\diamondsuit$ we can deduce that $Y^\diamondsuit\times_{X^\diamondsuit} \Spa(R,R^+)$ is also a locally spatial diamond. 
	From the $\diamondsuit$ construction, the map $Y^\diamondsuit\to X^\diamondsuit$ is automatically partially proper, hence compactifiable.
	The finiteness in transcendence degree of $Y^\diamondsuit\times_{X^\diamondsuit} \Spa(R,R^+)\to \Spa(R,R^+)$ can be deduced from the analogous statement for $W^\diamondsuit\times_{X^\diamondsuit} \Spa(R,R^+)\to \Spa(R,R^+)$ since the map $W^\diamondsuit  \to Y^\diamondsuit$ lifts all geometric points.
\end{proof}

	\subsubsection{Six operations on the co-sheaf setup}
	At a technical level, the co-sheaf setup is a 6-functor formalism, but one has to be careful since this can potentially lead to confusion (see \Cref{remark warn reader 3 functor}).
Indeed, as we will see, in Zhu's setup the ``usual'' 6-operations do not coincide with the 6-operations that a 6-functor formalism naturally provides.
We will clarify this point below.

The way Zhu upgrades $\Shv^!_c$ and $\Shv^!$ to 3-functor formalisms on $\Corr(\AlgSp^\qcqs,E_\pfp^\Alg)$ is by recalling that $\Shv^!_c=\calD_\cons^\op$ and that, for a map $[f:X\to Y]\in (\AlgSp^\qcqs, E_\pfp^\Alg)$, the functor $f_!$ preserves constructibility.
In this way Zhu defines $\ul{f_*}:= f_!^\op$ in $\Shv_c^!(X)\to \Shv_c^!(Y)$ and $\ul{f_*}:\Shv^!(X)\to \Shv^!(Y)$ is defined as the ind-extension.
Formally,
\[\Shv^!:\Corr(\AlgSp^\qcqs,E_\pfp^\Alg)\to \LinCat_\Lambda\]
is defined as $\Shv^!:=\iota\circ (-)^\vee \circ \Shv^*$ where $(-)^\vee$ denotes the self-equivalence
\begin{equation}
	\label{eq: dual-cat-defi}
(-)^\vee:\LinCat_\Lambda^{\on{dual}}\to \LinCat_\Lambda^{\on{dual}}
\end{equation}
discussed in \cite[Remark 7.24]{Zhu25} given by taking the dual categories, and $\iota$ denotes the (non-full) inclusion 
\[\iota:\LinCat_\Lambda^{\on{dual}} \to \LinCat_\Lambda \]
of the dualizable $\Lambda$-linear categories.

One can of course apply \Cref{cons: ExtensiontoRepresentableMaps} to obtain 3-functor formalisms
\begin{equation}{\label{eqn: Shv!candShvcformalism}}
	\Shv^!_c:\Corr(\PreStk,E^{\Alg,\rep}_\pfp)\to \LinCat^\sm_\Lambda \quad \text{ and } \quad \Shv^!:\Corr(\PreStk,E^{\Alg,\rep}_\pfp)\to \LinCat_\Lambda,
\end{equation}
where again $E^{\Alg,\rep}_\pfp$ denotes the maps of prestacks that are relatively representable in algebraic spaces, and in $E^\Alg_\pfp$.

	\begin{convention}
		\label{underline convention}
		As in \Cref{basic-obs}, we give names to the 3 basic operations furnished by the 3-functor formalism (\ref{eqn: Shv!candShvcformalism}).
		\begin{enumerate}
			\item For $f:X\to Y\in E^{\Alg,\rep}_{\pfp}$ we denote by $\ul{f_*}:\Shv^!(X)\to \Shv^!(Y)$, the image of $\{X\overset{\id}{\leftarrow}  X \overset{f}{\ra} Y\}$. 	
			\item For $f:X\to Y\in \PreStk$ we denote by $\ul{f^!}:\Shv^!(X)\to \Shv^!(Y)$, the image of $\{X\overset{f}{\leftarrow}  Y \overset{\id}{\ra} Y\}$.
			\item We denote by $- \boxtimes^! -:\Shv^!(X)\otimes_{\LinCat_\Lambda} \Shv^!(Y)\to \Shv^!(X\times Y)$ the map coming from lax-monoidality.
			\item We denote by $-\otimes^! -:\Shv^!(X)\otimes_{\LinCat_\Lambda} \Shv^!(X)\to \Shv^!(X)$ the composition of $\ul{\Delta^!}$ with $\boxtimes^!$.
		\end{enumerate}
	\end{convention}

    \begin{remark}
	    Let $X\to Y$ be a map in $(\PreStk,E^{\Alg,\rep}_\pfp)$.
		Under the identification $\calD_\cons(X)^\op\simeq \Shv^!_c(X)$ the operation $\ul{f_*}$ corresponds to $f_!^\op:\calD_\cons(X)^\op\to \calD_\cons(Y)^\op$.
		This justifies the notation used in \Cref{underline convention}.
	\end{remark}

	\begin{remark}
		\label{remark warn reader 3 functor}
		We warn the reader that although $\Shv^!$ automatically encodes $6$-operations (see \Cref{rem: presentable3functorsare6functors}), some of these operations should be treated as exotic operations. 
	Indeed, in this setup $\ul{f^!}$ and $\ul{f_*}$ always admit a right-adjoint when customarily one would expect the functors to admit a left-adjoint. 
	\end{remark}

	\begin{definition}
		For any map $f:X\to Y$ in $\PreStk$ we let $\ul{f_\flat}:\Shv^!(Y)\to \Shv^!(X)$ denote the right adjoint of $\ul{f^!}$.
	\end{definition}
	
	Although it is technically helpful to study $\ul{f_\flat}$ and the question of when this functor preserves constructibility as in \cite[\S 10.4.3]{Zhu25}, we are ultimately interested in defining and understanding the more familiar functors $\ul{f_!}$ and $\ul{f^*}$, which will be our notation for the left adjoints of $\ul{f^{!}}$ and $\ul{f_{*}}$. 
	However, these latter functors might not always exist in the co-sheaf setup.  
	Recall the following statement that isolates some situations where they do.
	\begin{proposition}[{\cite[Proposition 10.87]{Zhu25}}]
		\label{prop some adjoints exist}
		Let $[f:X\to Y]\in(\PreStk, E^{\Alg,\rep}_\pfp)$ be a morphism of prestacks.	
		\begin{enumerate}
			\item If $f$ is \'etale, then $\ul{f_*}$ is a right adjoint to $\ul{f^!}$ (i.e., $\ul{f^*}:=\ul{f^!}$).
			\item If $f$ is pfp proper, then $\ul{f_*}$ is a left adjoint to $\ul{f^!}$ (i.e., $\ul{f_!}:=\ul{f_*}$). 
			\item If $f$ is $\fShv$-suave, then $\ul{f_*}$ admits a left adjoint $\ul{f^*}$, which preserves constructibility. 
				\label{third case base change}
			In addition, for a pullback square of prestacks 
	\begin{equation}
		\label{sample pullback square}
		\begin{tikzcd}
			X_W \arrow{r}{f'} \arrow{d}{g'}  &W \arrow{d}{g} \\
			X \arrow{r}{f} &Y  
		\end{tikzcd}
	\end{equation}
			the natural map
			$\ul{(f')^*}\circ \ul{ g^!} \cong \ul{(g')^!}\circ \ul{ f^*}$
			is an isomorphism.

		\end{enumerate}

	\end{proposition}

	\begin{remark}{\label{rem: IndExtendingUnderDuality}}
		Let $[f:X\to Y]\in (\AlgSp^\qcqs,E^\Alg_\pfp)$ be a representable morphism that is $\fShv$-suave.
		Under the identification $\calD_\cons(X)^\op\simeq \Shv^!_c(X)$, the operation $\ul{f^*}$ induced on constructible sheaves by \Cref{prop some adjoints exist} (3) corresponds to $f^{!,\op}:\calD_\cons(Y)^\op\to \calD_\cons(X)^\op$.
		Moreover, since $\ul{f^*}$ is a left-adjoint to $\ul{f_*}$ it commutes with colimits, and in this case $\Shv^!(X)$ is the ind-extension of $\Shv^!_c(X)$, by the same reasoning as in the proof of \Cref{compact-generation-shv*}.
In particular, $\ul{f^*}$ is determined by its value on constructible sheaves. 
		For more general prestacks $X$ and $Y$, the functor $\ul{f^*}$ might not admit such a concrete description and its existence is instead deduced from \cite[Lemma 8.46]{Zhu25}.
	\end{remark}
    \begin{remark}
    We emphasize that, when formulating the notion of suave in \Cref{prop some adjoints exist} (3), we really do mean that it is suave with respect to the $6$-functor formalism $\fShv$ or equivalently $\Shv^{*}$ in light of \Cref{prop: comparisonofShv*andfShv} (which we implicitly used when citing \cite{Zhu25}).
Throughout this document we formulate notions of suaveness, unipotence, and smoothness with respect to the $*$-glued theory as in \cite{Zhu25}, and use it to deduce consequences for the $!$-glued theory.
    \end{remark}

	The following statement says that $\Shv^!(X)$ satisfies excision with respect to well-behaved closed immersions.
	Here one implicitly uses \Cref{prop some adjoints exist} to know that all of the functors in the statement below are defined.
	\begin{proposition}[{\cite[Lemma 10.90]{Zhu25}}]{\label{prop: excisionforPrestacks}}
		Let $j : U \to X$ be a quasicompact open embedding of prestacks with a pfp closed complement $i : Z \to X$.
Then:
		\begin{enumerate}
			\item $\ul{i^!}\circ \ul{j_*} \simeq 0$ and $\ul{j^!}\circ \ul{i_*} \simeq 0$.
			\item The functors $\ul{i_*}$ (resp. $\ul{j_*}$) are fully faithful, with essential image consisting of $\calF \in \Shv^!(X)$ with $\ul{j^!} \calF \simeq 0$ (resp. $\ul{i^!} \calF \simeq 0$).
			\item For $\calF \in \Shv^!(X)$ we have a functorial fiber sequence
			\[
			\ul{i_*} \circ \ul{i^!} \calF \to \calF \to \ul{j_*} \circ \ul{j^!} \calF.
			\]
			where the morphisms are the unit and co-unit of the corresponding adjunctions.
		\end{enumerate}
	\end{proposition}

	\subsection{Placid geometry}{\label{ss: placid stacks}}
	\subsubsection{Placid stacks}
	We can think of the category $\PreStk^\Art$ of \Cref{defn: SchematicArtinPreStks} as those geometric objects that are finite dimensional over $\Spec k$. 
	For our considerations, we will need to enlarge this category.  
	The category of very placid stacks of \Cref{definition very placid} (also \cite[Definition 10.114]{Zhu25}) is an appropriate enlargement and it plays a crucial role in this article for the following reasons.
	\begin{enumerate}
		\item Very placid stacks have well-behaved co-sheaf theory.
		\item The moduli spaces of shtukas are very placid.
	\end{enumerate}
	
	\begin{definition}[{\cite[Definition 10.49]{Zhu25}}]{\label{defn: proellsmoothandprounipotent}}
		Let $f:X\to Y$ be a map in $\AlgSp^{\qcqs}$. 
		\begin{enumerate}
			\label{defn: proellsmoothandprounipotent:1}
		\item We say that $f$ is \emph{cohomologically pro-smooth} (or \emph{coh.~pro-smooth}) if $X=\varprojlim_{i\in I}X_i$ with $X_i\in \AlgSp^\qcqs$ such that:
			\begin{itemize}
				\item each map $[f_i:X_i\to Y]$ is $\fShv$-smooth and in $E^\Alg_\pfp$,  
				\item each of the transition maps $X_i\to X_j$ is in $E^\rep_\pfp$, affine and $\fShv$-smooth,  
				\item and the limit is taken over a cofiltered diagram $I$. 
			\end{itemize}
			\label{defn: proellsmoothandprounipotent:2}
			\item We say that a coh.~pro-smooth map is \emph{strongly cohomologically pro-smooth} if the affine transition maps $X_i\to X_j$ in a presentation as above are all surjective.
			\label{defn: proellsmoothandprounipotent:3}
			\item We say $X\to Y$ is \emph{essentially cohomologically pro-smooth} (resp. \emph{essentially pro-\'etale}) if there is a factorization $X\to Z\to Y$ where $X\to Z$ is coh.~pro-smooth (resp. pro-\'etale) and $Z\to Y$ is in $E^\Alg_\pfp$.
		\end{enumerate}
	\end{definition}
	
	Recall the definition of a unipotent map \Cref{defn: Suavesmoothunipotent}, specifically we will consider \(\fShv\)-unipotent maps and their implications to the $\Shv^!$-theory.  
	
	\begin{definition}[{\cite[Definition 10.59]{Zhu25}}]{\label{defn: proellunipotent}}
		Let $f:X\to Y$ be a map in $\AlgSp^\qcqs$. 
		\begin{enumerate}
			\label{defn: proprounipotent:1}
			\item We say that $f$ is \emph{pro-unipotent} if we have a cofiltered presentation $X=\varprojlim_{i\in I}X_i$ with $X_i\in \AlgSp^\qcqs$ such that each map $f_i:X_i\to Y$ is $\fShv$-unipotent, and such that each transition map $X_i\to X_j$ is in $E^\rep_\pfp$, affine and $\fShv$-unipotent.
			\label{defn: prounipotent:2}
			\item We say $X\to Y$ is \emph{essentially pro-unipotent} if there is a factorization $X\to Z\to Y$ where $X\to Z$ is pro-unipotent and $Z\to Y$ is in $E^\Alg_\pfp$.
		\end{enumerate}
	\end{definition}

	We recall the definition of placid stacks in Zhu's context (which is a variant of the concept introduced in \cite{perverse-sheaves-infinite-dimensional-stacks}).
	Roughly speaking, placid stacks are infinite dimensional spaces which exhibit cohomological behavior that is similar to finite dimensional spaces. 
	Our basic building block is as follows.
	\begin{definition}[{\cite[Definition 10.61]{Zhu25}}]
		\label{standard-placidness}
		We say an algebraic space $X \in \AlgSp^{\qcqs}$ is \emph{standard placid} if the structure morphism $X \ra \Spec{k}$ is essentially cohomologically pro-smooth in the sense of \Cref{defn: proellsmoothandprounipotent}.
We let $\AlgSp^{\spl} \subset \AlgSp^{\qcqs}$ denote the full subcategory of standard placid spaces.
	\end{definition}

	\begin{definition}[{\cite[Definition 10.104]{Zhu25}}]{\label{defn: universalhomologicaldescent}}
		A morphism of perfect prestacks $f: X \ra Y$ is said to be of \emph{homological descent} if the canonical map
		\[ \Shv^!(Y) \ra \prolim_{n \in \Delta} \Shv^!(X_{n}) \]
		induced by $\underline{(-)}^{!}$ is an equivalence.
Here $X_{n}$ is the $n$th term of the \v{C}ech nerve of $X \ra Y$ and $\Delta$ denotes the simplex category. 
		Similarly, we say it is of \emph{universal homological descent} if this is true after any base-change $Y' \ra Y$.
	\end{definition}

	\begin{remark}
		We note that being of (universal) \emph{homological descent} is simply saying that $f:X\to Y$ is a (universal) $\calD^*$-cover in the sense of \Cref{defn: universal*and!covers} for $\calD=\Shv^!$. 	
	\end{remark}
	
	We extend \Cref{standard-placidness} to stacks as follows.
	\begin{definition}[{\cite[Definition 10.114]{Zhu25}}]
		\label{definition very placid}
		Let $X\in \SchStk_\et$ be a quasicompact \'etale stack.  
		\begin{enumerate}
			\item We say $X$ is a \emph{quasicompact quasi-placid \'etale stack} if there exists a surjective representable coh.~pro-smooth morphism $U \ra X$ with $U \in \AlgSp^{\spl}$. 
				We refer to the map $U\to X$ as a \textit{quasi-placid atlas} for $X$.  
			We denote by $\SchStk_\et^\qpl\subseteq \PreStk$ the full subcategory of quasicompact quasi-placid \'etale stacks.
		\item We say $X$ is a \emph{quasicompact placid \'etale stack} if $X\in \SchStk_\et^\qpl$ and there is a quasi-placid atlas that is a) strongly coh.~pro-smooth in the sense of \Cref{defn: proellsmoothandprounipotent} (3), and is b) of universal homological descent in the sense of \Cref{defn: universalhomologicaldescent}. 
			We call such quasi-placid atlas a \emph{placid atlas}. 
			We denote by $\SchStk_\et^\pl\subseteq \PreStk$ the full subcategory of quasicompact placid \'etale stacks.
			\item We say $X$ is a \emph{quasicompact very placid \'etale stack} if $X\in \SchStk_\et^\qpl$ and there is a quasi-placid atlas which is also essentially cohomologically pro-unipotent in the sense of \Cref{defn: proellunipotent}. 
			We call such a placid atlas a \emph{very placid atlas}.
			We denote by $\SchStk_\et^\vpl\subseteq \PreStk$ the full subcategory of very placid stacks.
		\end{enumerate}
	\end{definition}

	\begin{remark}
		\label{quasicompactness-assumption}
 To simplify the exposition, we only give \Cref{definition very placid} for quasicompact quasi-placid (resp. quasicompact placid, resp. quasicompact very placid) stacks.  
Zhu gives a more general definition (\cite[Definition 10.114]{Zhu25}). 	
In our work, we will only consider quasi-placid \'etale stacks that are quasicompact and, if the context is clear, we will drop the word ``quasicompact'' from the terminology. 
	\end{remark}
	
	\begin{proposition}[{\cite[Proposition 10.108]{Zhu25}}]
		Let $f:U\to X$ be a map of prestacks.
If $f$ is representable and essentially pro-unipotent, then $f$ is of universal homological descent.
In particular, if $X\in \SchStk^\vpl$ and $U\to X$ is a very placid atlas, then 
		\[ \Shv^!(X) \simeq \prolim_{n \in \Delta} \Shv^!(U_{n}), \]
		where $U_n$ denotes the $n$-fold fiber product of $U$ over $X$.
	\end{proposition}

	All of the operations discussed in \Cref{prop some adjoints exist} are with respect to maps $f:X\to Y$ that are representable and pfp (i.e., in $E^{\Alg,\rep}_\pfp$).  
	For our purposes we will need to discuss $\ul{f_!}$ (the left-adjoint of $\ul{f^!}$) when the map $f$ is not necessarily pfp. 
	For example, the atlas map $f:U\to X$ for $X$ a very placid stack.
	We have the following.

	\begin{proposition}[{\cite[Proposition 10.124.(1)]{Zhu25}}]
		\label{existence of the adjoint}
		Let $f : X \to Y$ be a morphism of quasicompact quasi-placid stacks. 
		If $f$ is representable, essentially pro-unipotent, then $\ul{f^!}:\Shv_c^!(Y)\to \Shv_c^!(X)$ admits a left adjoint when restricted to constructible subcategories $\ul{f_!} : \Shv^!_c(X) \to \Shv^!_c(Y)$.
	\end{proposition}

	We also have the following.

	\begin{proposition}[{\cite[Proposition 10.145]{Zhu25}}]
		\label{existence of the adjoint II}
		Let $f : X \to Y$ be a morphism of quasicompact placid stacks. 
		If $f$ is representable, essentially pro-unipotent, then $\ul{f^!}:\Shv^!(Y)\to \Shv^!(X)$ admits a left adjoint.
	\end{proposition}
	\begin{proof}
	\cite[Proposition 10.145]{Zhu25} only states the existence of the adjoint when $X$ and $Y$ are very placid, but the paragraph below it clarifies that the existence of the adjoint only requires $X$ and $Y$ to be placid.
	\end{proof}

	Almost by definition (see the discussion preceding \Cref{eq: dual-cat-defi}), when $X\in \AlgSp^\qcqs$ the inclusion map $\Shv_c^!(X)\hookrightarrow \Shv^!(X)$ induces an identification $\Shv^!_{c}(X)\simeq \Shv^!(X)^{\omega}$. 
	For general prestacks this is far from true;  nevertheless, we still have the following.
	\begin{lemma}[{\cite[Lemma~10.127, Proposition~10.144]{Zhu25}}]
		\label{quasi-placid have compact generated}
		Let $X$ be a quasicompact quasi-placid stack.
The following is true.
		\begin{enumerate}
			\item The essential image of the natural inclusion $\Shv^!(X)^\omega \hookrightarrow \Shv^!(X)$ lies in $\Shv^!_c(X)$. 
			\item Let $X$ be quasicompact very placid with corresponding very placid atlas $f: U \ra X$. 
			Then the category $\Shv^!(X)$ is compactly generated. 
			The functor $\ul{f_{!}}$ preserves compact objects and a family of compact generators is given by $\ul{f_{!}}(\mathcal{F})$ for $\mathcal{F} \in \Shv^!_{c}(U)$.
			\item In the situation of (2), if $f$ is cohomologically pro-unipotent (not just essentially pro-unipotent), then we have that $\Shv^!_{c}(X) \simeq \Shv^!(X)^{\omega}$. 
		\end{enumerate}
	\end{lemma}
	
	We also need the following statement on the existence of $\ul{f^*}$ (the left-adjoint of $\ul{f_*}$).	

	\begin{lemma}
		\label{placid:have-left-adjoint}
		Let $f : X \to Y\in \SchStk^\pl_\et$ be a morphism of quasicompact placid stacks. 
		Suppose that $f$ is in $E^{\Alg,\rep}_\pfp$, then $\ul{f_*}$ has a left-adjoint $\ul{f^*}$. 
	\end{lemma}
	\begin{proof}
		This follows from the proof of \cite[Proposition 10.145]{Zhu25} and the paragraph above it.   
	\end{proof}

	A useful consequence of this is the following.

	\begin{corollary}
		\label{somecompactgeneration}
	Let $f:X\to Y$ be a pfp closed immersion of quasicompact placid stacks.	
	If $\Shv^!(Y)$ is compactly generated, then $\Shv^!(X)$ is also compactly generated.
	\end{corollary}
	\begin{proof}
		If $\{a_i\}_{i\in I}$ is a family of compact generators for $\Shv^!(Y)$, then $\{\ul{f^*}a_i\}_{i\in I}\subseteq \Shv^!(X)$ is a family of compact generators for $\Shv^!(X)$. 
		Indeed, $\ul{f^*}$ preserves compact objects since $\ul{f_*}$ admits a further right adjoint, and since $\ul{f_*}$ is conservative (even fully faithful by \Cref{prop: excisionforPrestacks} (2)) the family $\{\ul{f^*}a_i\}_{i\in I}$ generates. 
	\end{proof}

	\subsubsection{Ind-placid and sind-placid stacks}
	In this section we discuss ind-placid and sind-placid stacks. 
	Intuitively speaking, these categories are obtained from the categories of quasi-placid stacks by formally adding filtered and sifted colimits along proper maps.
	As mentioned before, the schematic local Langlands category is the category of co-sheaves on the Kottwitz stack which is a sind-very placid stack. 
	In this sense, sind-placid stacks also play a key role in our considerations. 
	
	The following statement, which can be deduced from applying the universal property of the category of presheaves and \cite[Theorem~5.1.5.6]{HTT}, is the main reason one can control the category of co-sheaves for sind-placid stacks. 
	\begin{lemma}{\label{lemma: LurieColimits}}
		Let $X = \colim_{i \in I} X_{i}$ be a colimit in $\PreStk$ then we have an identification 
		\[ \Shv^!(X) \simeq \varprojlim_{\substack{i \in I^\op\\ \ul{f^!}}} \Shv^!(X_{i}) \]
		of categories where the transition maps on the right-hand side are formed with respect to $!$-pullback internal to the $\Shv^!(-)$-theory. 
	\end{lemma}

	Combining \Cref{lemma: LurieColimits} above with \cite[Theorem 5.5.3.18, Corollary 5.5.3.4]{HTT}, we deduce the following. 
	\begin{proposition}{\label{prop: Shvshriektheoryisnice}}
		Let $X = \colim_{i \in I} X_{i}$ be a colimit of prestacks with transition maps representable in pfp proper algebraic spaces, then we have an equivalence 
		\[ \Shv^!(X) \simeq \varinjlim_{\substack{i \in I \\ \ul{f_*}}} \Shv^!(X_{i}) \]
		in $\LinCat_\Lambda$, where the transition maps on the right-hand side are formed with respect to $*$-pushforward internal to the $\Shv^!(-)$-theory.
	\end{proposition}
	
	This motivates the following definition.

	\begin{definition}[{\cite[Definition 10.93]{Zhu25}}]
		\label{defn: indproper}
		A morphism of prestacks $f: X \ra Y$ is called \emph{ind-pfp proper} if, for all $S \in \AlgSp^{\qcqs}$, with a map $S \ra Y$ the pullback $f_{S}: X_{S} \ra S$ admits a presentation as a colimit $\colim_{i \in I} X_{i}$ with transition maps $X_{i} \ra X_{j}$ being pfp closed immersions of algebraic spaces such that $X_{i} \ra S$ is \emph{pfp proper} (i.e., $[X_i\ra S]\in P^\Alg_\pfp$).
	\end{definition}
    The left adjoint to \(\ul{f^!}\), for an ind-pfp proper map \(f\), behaves as expected.
    \begin{lemma}[{\cite[Lemma 10.100]{Zhu25}}]{\label{lemma: continuousrightadjointpropermpas}}
        Let \(f\from X\to Y\) be an ind-pfp proper morphism of prestacks.
        Then \(\ul{f_*}\) is the left adjoint to \(\ul{f^!}\).
    \end{lemma}
    We thus obtain an analogue of \Cref{prop: Shvshriektheoryisnice}, by the exact same line of reasoning.
    \begin{lemma}\label{lem: colimit ind-proper transitions}
        Let $X = \colim_{i \in I} X_{i}$ be a colimit of prestacks with ind-pfp proper transition maps, then we have an equivalence 
		\[ \Shv^!(X) \simeq \varinjlim_{\substack{i \in I \\ \ul{f_*}}} \Shv^!(X_{i}). 
        \]
    \end{lemma}
    Ind-proper surjective maps satisfy descent.
More precisely, we have the following.
	\begin{proposition}[{\cite[Proposition 10.106]{Zhu25}}]\label{ind-pfp proper}
		Let $f\from X\to Y$ be a ind-pfp proper and surjective morphism of prestacks. 
		Then $f$ is of universal homological descent in the sense of \Cref{defn: universalhomologicaldescent}.
	\end{proposition}

	\begin{remark}
	In \Cref{ind-pfp proper}, the precise meaning of the representable map $f$ being surjective is that, for any map $S\to Y$ for any qcqs scheme $S$, the base-change morphism $f_S:X\times_Y S \to S$ is surjective as a map of algebraic spaces (i.e a surjection at the level of topological spaces).
The argument of \cite[Proposition 10.106]{Zhu25} shows that such a surjective ind-pfp proper has the property that the base-change $f_S:X\times_Y S \to S$ is surjective for the h-topology (See \cite[Remark~10.107]{Zhu25}). 
	\end{remark}
	\begin{definition}[{\cite[Definition 10.153, 10.157]{Zhu25}}]
		Let $X\in \SchStk_{\et}$.
		\begin{enumerate}
			\item We say $X$ is \emph{quasicompact Ind-quasi-placid (resp. quasicompact Ind-placid, quasicompact Ind-very-placid)} if $X = \colim_{i \in I} X_{i}$ is a filtered colimit of quasicompact quasi-placid (resp. quasicompact placid, quasicompact very placid) stacks $X_{i}$ along representable pfp closed immersions. 
			We let 
			\[\IndStk^{\vpl}_{\et} \subset \IndStk^\pl_{\et} \subset \IndStk^{\qpl}_{\et} \subset \PreStk\] denote the corresponding full subcategories of objects in $\SchStk_{\et}$. 
			\item We say $X \in \SchStk_\et$ is a \emph{quasicompact sind-placid stack (resp. quasicompact sind-very placid)} if it admits an  \'etale-surjective ind-pfp proper morphism in the sense of \Cref{defn: indproper} from a quasicompact Ind-placid stack (resp. quasicompact Ind-very-placid stack) $V \ra X$. 
			We call such a $V \ra X$ a sind-placid (resp. sind-very placid) atlas.
			We let $\sStk_\et^\vpl\subset \sStk_\et^\pl \subseteq \SchStk_\et$ denote the categories of qc sind-very placid and qc sind-placid stacks.
		\end{enumerate}
	\end{definition}

	\begin{remark}
		We warn the reader that, despite the terminology, quasicompact Ind-quasi-placid and quasicompact sind-quasi-placid stacks are not quasicompact as \'etale stacks. 
		Since we only work with quasicompact Ind-quasi-placid and quasicompact sind-quasi-placid stacks, we will often drop the word ``quasicompact'' from the terminology.
	\end{remark}

	Let $X\in \IndStk^\qpl_\et$ be an Ind-quasi-placid stack $X$ with presentation $\varinjlim_{i \in I} X_{i}$, where each $X_i\in \SchStk^\qpl$ and each transition map $X_i\to X_j$ is a pfp closed immersion. 
	By \Cref{prop: Shvshriektheoryisnice}, we obtain a natural equivalence 
	\[ \Shv^!(X) \simeq \varinjlim_{\substack{i \in I \\ \ul{f_*}}}  \Shv^!(X_{i}). \]

	Similarly, by \Cref{ind-pfp proper} for $X\in \sStk^\pl$ and $V\to X$ a sind-placid atlas we have 
	\[ \Shv^!(X) \simeq \varprojlim_{\substack{n \in \Delta \\ \ul{f^!}}}  \Shv^!(V^{n}) \simeq \varinjlim_{\substack{n \in \Delta^\op \\ \ul{f_*}}}  \Shv^!(V^{n}),\]
	where $V^n$ denotes the $n$-fold product of $V$ over $X$.

We have the following analogue of \Cref{placid:have-left-adjoint} for ind-placid stacks.
\begin{proposition}
	\label{pass-to-bounded-piece-compacts}
	Let $X\in \IndStk^\pl_\et$ with ind-presentation $X \simeq \varinjlim_{i\in I} X_i$ with $X_i\in \SchStk^\pl_\et$.  	
	Let $f_i:X_i\to X$ denote the pfp closed immersion. 
	Then, $\ul{f_{i,*}}$ has a left-adjoint $\ul{f_i^*}$. 
	Moreover, $\Shv^!(X)$ is compactly generated if and only if every $\Shv^!(X_i)$ is compactly generated.
\end{proposition}
\begin{proof}
	The existence of the left-adjoint follows from \Cref{placid:have-left-adjoint} and from \cite[Proposition 5.1.8.(c)]{perverse-sheaves-infinite-dimensional-stacks}.
	If each $\Shv^!(X_i)$ is compactly generated, it follows from \Cref{prop:colimit_commutes_dual} below that $\Shv^!(X)$ is also compactly generated.
	Conversely, if $\Shv^!(X)$ is compactly generated we may argue as in \Cref{somecompactgeneration} to show that each $\Shv^!(X_i)$ is also compactly generated.
\end{proof}

	We also have the following more refined result on excision for these particular types of stacks.
	\begin{proposition}[{\cite[Proposition 10.177]{Zhu25}}]
		{\label{prop: semiorthogonaldecompositionsonsindplacid}}
		Let $Y$ be a sind-very placid stack and let $j : U \to Y$ be a qcqs open embedding with a closed complement $i : Z \to Y$.
Then $U$ and $Z$ are sind- very placid and $i$ is pfp.
In particular, by \Cref{prop: excisionforPrestacks}, we have that.
		\begin{enumerate}
			\item $\ul{i_!} = \ul{i_* },\ \ul{j^!}\circ \ul{i_*} \simeq 0,\ \ul{i^!}\circ  \ul{j_*} \simeq 0,\ \text{and}\ \ul{i^*}\circ  \ul{j_!} \simeq 0$.
			\item The functor $\ul{i_*}$ (resp. $\ul{j_*}$, resp. $\ul{j_!}$) is fully faithful, with essential image consisting of those $\calF \in \Shv^!(X)$ satisfying $\ul{j^!} \calF = 0$ (resp. $\ul{i^!} \calF = 0$, resp. $\ul{i^*} \calF = 0$).
			\item For every $\calF \in \Shv(Y)$, we have the following canonical fiber sequences:
			\[
			\ul{i_*}\circ \ul{i^!} \calF \to \calF \to \ul{j_*}\circ \ul{j^*} \calF,
			\]
			and
			\[
			\ul{j_!}\circ \ul{j^*} \calF \to \calF \to \ul{i_*}\circ \ul{i^*} \calF.
			\]
		\end{enumerate}
	\end{proposition}
	
	\subsection{Cosheaves and analytification}	\label{sec:cosheaves_and_analytification}
	We start this subsection, by offering some informal explanations and perspectives of the difficulties related to reconciling the cosheaf formalism of Zhu introduced in the previous sections and Scholze's analytification functor $c^*$ introduced in \S \ref{sec: analytificiation and 6-functor formalisms}.
	We do this to guide the reader in the discussion that follows.

	The first observation is that Scholze's functor is more directly related to $\Shv^*(X)$ and that unfortunately, for a general prestacks $X\in \PreStk$ it is hard to understand the relation between $\Shv^*(X)$ and $\Shv^!(X)$.
	When $X\in \AlgSp^\qcqs$, the categories are canonically equivalent through ind-extending Verdier duality on constructible sheaves. %
	When $X$ is quasi-placid the categories $\Shv^*_c(X)$ and $\Shv^!_c(X)$ are still equivalent, but the equivalence will usually rely on additional data (\cite[Lemma 10.129, Proposition 10.130]{Zhu25}).  
	For more general prestacks, it is hard to understand if there is a relation at all. 

	\begin{question}
		Does there exist an equivalence between $\Shv^*(\ICG)$ and $\Shv^!(\ICG)$? Or equivalently, are $\Shv^*(\ICG)$ and $\Shv^!(\ICG)$ dual to each other? 
	\end{question}

	Since the analytification functors $c^*$ are naturally defined for the $\Shv^*$-theory, two possibilities come to mind:
	\begin{itemize}
		\item Develop a cosheaf theory in the analytic setup over which an analytification functor $\ul{c^!}$ is naturally defined. 
		\item Study the precise relation between $\Shv^!(X)$ and $\Shv^*(X)$ for enough geometric objects.
	\end{itemize}

	The first option would take us too far afield, and would require us to overcome substantial technical problems.
    However, it would be interesting if such a theory became available. 
	We will take the second approach. 
	Let us recall the following statement.
	\begin{proposition}[{\cite[Proposition 10.148]{Zhu25}}]\label{very-placid-givesequiv}
		Suppose that $X=U/H$ where $U$ is a standard placid algebraic space and $H$ is an affine group scheme over $k$ fitting in an exact sequence of the form
		\[e\to H_0\to H\to H'\to e\]
		such that $H_0$ is coh.~pro-unipotent over $k$ and $H'$ is pfp over $k$.
		The following statements hold.
		\begin{enumerate}
			\item $X$ is very placid.
		\item The canonical identification $\Shv^!_c(X)\simeq (\Shv^*_c(X))^\op$ induces an identification 
			\[(\Shv^!(X))^\omega\simeq (\Shv^*(X)^\omega)^\op.\]
			\item There exist an identification 
				\[\id^\eta:\Shv^!(X)\simeq \Shv^*(X),\] 
				relying on the choice of a generalized constant sheaf $\Lambda_{X}^{\eta}\in \Shv^!_c(X)$ (see \cite[Definition 10.63, Definition 10.128]{Zhu25}) and the discussion on \cite[\S 10.4.2]{Zhu25}).
		\end{enumerate}
	\end{proposition}
Using \Cref{very-placid-givesequiv}, one can easily define an analytification functor 
\[\ul{c^!_\eta}:\Shv^!(X)\to \calD^{\diamondsuit_{\uop}}_{\Lambda}(X) \simeq \calD_{\Lambda}^{\an}(X^{\Diamond}) \]
by composing $\id^\eta:\Shv^!(X)\to \Shv^*(X)$ with the natural maps 
\[\Shv^*(X)\to \calD^\sch_\Lambda(X) \to \calD^\diamondsuit_\Lambda(X)\]
from \Cref{thm: DiamondAnalytificationUltimate}, \Cref{important analytification prop}, and \Cref{identification zhu vs standard}.

As the reader might have noticed, the analytification functor proposed $\ul{c^!_\eta}$ will bring hell and tempest to the poor soul that attempts to make it functorial, since it clearly depends on the choice of $\Lambda_X^\eta$.

An alternative approach, is to combine \Cref{very-placid-givesequiv} with \Cref{quasi-placid have compact generated} to observe that, whenever $X=U/H$ is very placid, then $\Shv^!(X)$ and $\Shv^*(X)$ are canonically dual. 
In other words,
\begin{equation}
	\label{duality-quasi-placid}
\Shv^*(X)\simeq \Shv^!(X)^\vee.
\end{equation}
We will take a version of this approach that will work for more general stacks. 
This will require us to study deeper the work of Zhu.
More precisely, we recall the ind-finitely generated sheaves of Zhu. 
We begin by recalling the following definition.
\begin{definition}{\cite[\S 10.5.4]{Zhu25}}
	A map $f:X\to Y\in \SchStk_{\et}^\qpl$ is in $\on{V}^!_c$ if the following conditions hold: 
	\begin{itemize}
		\item $f$ is $!$-able for the 3-functor formalism $\Shv^!$ (i.e., $\ul{f_*}$ is well-defined) in the sense of the extension procedure of \Cref{shriek hull} applied to (\ref{eqn: Shv!candShvcformalism}).
		\item For all $Y'\to Y\in \SchStk_{\et}^\qpl$ the pullback $X\times_Y Y'\in \SchStk_{\et}^\qpl$, and if $f':X\times_Y Y'\to Y'$ denotes the base change, then $\ul{f'_*}$ preserves constructible sheaves.
	\end{itemize}
Analogously, a map $f:X\to Y\in \SchStk_{\et}^\qpl$ is in $\on{V}^*_c$ if the following conditions hold: 
	\begin{itemize}
		\item $f$ is $!$-able for the 3-functor formalism $\Shv^*$ (i.e., $f_!$ is well-defined).
		\item For all $Y'\to Y\in \SchStk_{\et}^\qpl$ the pullback $X\times_Y Y'\in \SchStk_{\et}^\qpl$, and if $f':X\times_Y Y'\to Y'$ denotes the base change, then $f'_!$ preserves constructible sheaves.
	\end{itemize}

\end{definition}

\begin{remark}
	We note that $(\SchStk_{\et}^\qpl,\on{V}^!_c)$ (resp. $(\SchStk_{\et}^\qpl,\on{V}^*_c)$) is a geometric setup that might not have finite limits, since fiber products of quasi-placid stacks might not remain quasi-placid.
	This is an instance where the notion of geometric setup discussed in \Cref{defn: geometricsetup} does not suffice and one really needs the more general notion defined in \cite[Definition 2.1.1]{HeyerMann}. 
	We also note that $E^{\Alg,\rep}_\pfp\subseteq \on{V}^!_c\cap \on{V}^*_c$, and that $(\SchStk_\et, E^{\Alg,\rep}_\pfp)$ is already a geometric setup in the stricter sense of \Cref{defn: geometricsetup}.
	Working with this geometric setup will suffice for our purposes.
\end{remark}

Almost by definition, the co-sheaf $3$-functor formalism on $\PreStk$ and $\Shv^!$-$!$-able maps restricts to a $3$-functor formalism
\[\Shv_c^!:\Corr(\SchStk_{\et}^\qpl,\on{V}^!_c)\to \LinCat_\Lambda^\sm,\]
we can ind-extend it to obtain a presentable $3$-functor formalism
\begin{equation}
	\label{defining-ind-finitely-quasi-placid}
\IndShv^!:\Corr(\SchStk_{\et}^\qpl,\on{V}_c^!)\to \LinCat^\cg_\Lambda \hookrightarrow \LinCat_{\Lambda}.
\end{equation}
These are so-called ind-finitely generated sheaves (see \cite[Equation (10.55)]{Zhu25}).

Analogously, we get 
\begin{equation}
	\label{defining-ind-finitely-quasi-placid-star-v}
\IndShv^*:\Corr(\SchStk_{\et}^\qpl,\on{V}^*_c)\to \LinCat_\Lambda.
\end{equation}

For $X\in \SchStk_{\et}^\qpl$, the natural inclusion $\Shv_{c}^!(X)\subseteq \Shv^!(X)$ gives rise to a natural transformation of $6$-functor formalisms on $\Corr(\SchStk_{\et}^\qpl,\on{V}_c)$ (see \cite[Equation (10.56)]{Zhu25})
\[\Psi:\IndShv^!\Rightarrow \Shv^! \]
given by formally taking the colimit.
Similarly, we have a natural transformation 
\begin{equation}{\label{eqn: naturaltransformationonShv*formalism}}
 \IndShv^{*} \Rightarrow \Shv^{*} 
\end{equation}
given by the natural inclusion $\Shv_{c}(X) \subseteq \Shv^{*}(X)$ for $X \in \SchStk_{\et}^{\qpl}$ and formally taking colimits.

Assume now that $\Shv^!(X)$ is compactly generated (e.g if $X$ is very placid, by \Cref{quasi-placid have compact generated}).
Then it follows that $\Psi$ admits a fully faithful left-adjoint 
\begin{equation}
	\label{the-left-adjoint-useful}
 \Psi^L:\Shv^!(X)\to \IndShv^!(X), 
\end{equation}
which can be obtained by Ind-extending the inclusion $\Shv^!(X)^\omega\subseteq \Shv^!_c(X)$ granted by \Cref{quasi-placid have compact generated} (1).

\begin{construction}\label{construction:from_Shv!dual_to_Shv*}
Assume that $X\in \SchStk_{\et}^\qpl$ and $\Shv^!(X)$ is compactly generated. 
The map $\Psi^L$ above admits a conjugate map 
\[(\Psi^L)^o:\Shv^!(X)^\vee\to \IndShv^!(X)^\vee,\]
which can be obtained by ind-extending the inclusion $(\Shv^!(X)^\omega)^\op\subseteq (\Shv^!_c(X))^\op$.  %
Note that, by Ind-extending the identification in \Cref{rem: ShvcisopShvcstar} for $X\in \SchStk_{\et}^\qpl$, we have an equivalence 
\[\IndShv^!(X)^\vee\simeq \IndShv^*(X),\]
and that $\IndShv^*(-) \simeq \on{Ind}(\Shv^*_c)(-)$ as functors over $(\SchStk_{\et}^\qpl)^\op$.
In particular, this allows us to define the composition
\[
\Shv^!(X)^\vee \stackrel{(\Psi^L)^o}{\to} \IndShv^!(X)^\vee\simeq \IndShv^*(X) \to \Shv^*(X),
\]
where the last map is as in \eqref{eqn: naturaltransformationonShv*formalism}.
\end{construction}

The map from \Cref{construction:from_Shv!dual_to_Shv*} can now be composed with the natural transformation of \Cref{important analytification prop} (2) and the analytification map of \Cref{thm: DiamondAnalytificationUltimate}.
We summarize the situation as follows.

\begin{proposition}
	\label{analytification-cosheave-quasi-plac}
If $X\in \SchStk_{\et}^\qpl$ then there is a well-defined analytification map 
\[c^{*,\vee}_\fg:\IndShv^!(X)^\vee\to \calD_\Lambda^\diamondsuit(X).\]
Furthermore, if $\Shv^!(X)$ is compactly generated, then the composition of $(\Psi^L)^o$ with $c^{*,\vee}_\fg$ gives a well-defined analytification map
\[c^{*,\vee}:\Shv^!(X)^\vee\to \calD_\Lambda^\diamondsuit(X).\]
\end{proposition}

\begin{remark}
	Evidently, whenever $\Shv^!_c(X)=\Shv^!(X)^\omega$ and $\Shv^*_c(X)=\Shv^*(X)^\omega$ hold, then $\Shv^!(X)^\vee\simeq \Shv^*(X)$ and under this identification $c^{*,\vee}$ and $c^*$ agree. 
\end{remark}

In particular, we see that even when $\Shv^*(X)$ does not necessarily agree with $\Shv^!(X)^\vee$, one can still define an analytification functor for objects in $\Shv^!(X)^\vee$, whenever $X\in \SchStk_{\et}^\qpl$ and $\Shv^!(X)$ is compactly generated.
We set some notation.

\begin{definition}
	We let $\SchStk_{\et}^\qome\subseteq \SchStk_{\et}^\qpl$ denote the subcategory of quasicompact quasi-placid \'etale stacks $X$ for which $\Shv^!(X)$ is compactly generated.
\end{definition}

\begin{remark}
	\label{rmk:naturality-subtle}
We warn the reader that the naturality of the analytification functor
\[c^{*,\vee}:\Shv^!(X)^\vee\to \calD_\Lambda^\diamondsuit(X),\]
is subtle. 
Indeed, given a map of stacks $f:X\to Y\in \SchStk_{\et}^\qome$ it is not clear if $\ul{f^!}$ will preserve compact objects.  
If it does, one obtains a commutative diagram 
\begin{center}
\begin{tikzcd}
	\Shv^!(Y)^\vee \arrow{r}{(\Psi^L)^o} \arrow{d}{(\ul{f^!})^o}  & \IndShv^*(Y)^\vee \arrow{d}{{\on{Ind}}(f^*)} \\
	\Shv^!(X)^\vee \arrow{r}{(\Psi^L)^o} & \IndShv^*(X)^\vee.
\end{tikzcd}
\end{center}
\end{remark}

As we pointed out in \Cref{rmk:naturality-subtle}, it is better to analytify ind-finitely generated cosheaves. 
Indeed, for those we can formulate a version of \Cref{analytification-cosheave-quasi-plac} that is functorial in the category of correspondences. 
Recall that we have an involution 
\[(-)^\vee:\LinCat_\Lambda^\cg\to \LinCat_\Lambda^\cg,\]
given by taking the dual category, which may be identified with $\mathcal{C}^{\vee} = \Ind(\mathcal{C}^{\omega,\op})$, as in \cite[Section~1.5.3]{drinfeld2015compactgenerationcategorydmodules}.
Moreover, given any $3$-functor formalism we can consider the result of composing it with $(-)^\vee$.
This allow us to define a $3$-functor formalism 
\[(\IndShv^!)^\vee:\Corr(\SchStk_{\et}^\qpl,\on{V}^!_c)\to \LinCat_\Lambda^\cg.\]

\begin{remark}
	Given a map $f:X\to Y\in (\SchStk_{\et})^\qpl_{\on{V}^!_c}$ we obtain a functor \[(\IndShv^!)^\vee([X=X\to Y]):(\IndShv^!)^\vee(X)\to (\IndShv^!)^\vee(Y),\]
	which we denote by $T_{X \ra Y}$. 
	The functor $T_{X \ra Y}$ is the ind-extension of 
\[(\ul{f_*})^\op:\Shv_c^!(X)^\op\to \Shv^!_c(Y)^\op.\] 
One could hope that if $f\in V^!_c\cap V^*_c$, then one would be able to rewrite $T_{X \ra Y}$ simply as the ind-extension of 
\[f_!:\Shv^*_c(X)\to \Shv^*_c(Y) \]
via the natural equivalence $\IndShv^{!}(-)^{\vee} \simeq \Ind(\Shv^{!}_{c}(-)^{\op})$. 
When 
\[[f:X\to Y]\in (\SchStk_\et, {E^{\Alg,\rep}_\pfp})\] 
this holds, since in this case by definition $\ul{f_*}:=(f_!)^\op$ on constructible categories. 
Nevertheless, outside of this case, the identification becomes much more subtle and it is not at all clear in which generality this should hold.
Indeed, although for $f\in V^!_c\cap V^*_c$ the functor $\ul{f_*}$ is defined and preserves constructible categories, it is defined through extensions procedures and, in general, its behavior might be hard to control. 
\end{remark}

For us, the following will suffice, which follows from the above discussion.
\begin{proposition}
	\label{more-6-functor-comparisons}
	The following statements hold.
	\begin{enumerate}
		\item We have an equivalence of $6$-functor formalisms on $\corr(\SchStk_{\et}^\qpl,E^{\Alg,\rep}_\pfp)$ with values in $\LinCat^\cg_\Lambda$ 
	\[\IndShv^*(-)\simeq \IndShv^!(-)^\vee.\]
\item The functors $c^{*,\vee}_\fg$ fit in a natural transformation of $6$-functor formalisms defined over $\corr(\SchStk_{\et}^\qpl,E^{\Alg,\rep}_\pfp)$
	\[c^{*,\vee}_\fg:\IndShv^!(-)^\vee\Rightarrow \calD_\Lambda^{\diamondsuit_\uop}(-).\]
\item The functors $c^{*,\vee}$ fit in a natural transformation of functors 
	\[(\SchStk_{\et}^\qome)_{P^{\Alg,\rep}_\pfp}\to \LinCat_\Lambda^\cg \]
	\[c^{*,\vee}:\Shv^!(-)^\vee\Rightarrow \calD_\Lambda^{\diamondsuit_\uop}(-).\]
	\end{enumerate}
\end{proposition}
\begin{proof}
	The first two items were dealt with in the discussion preceding this statement, so we only provide the argument for the third.
	Note that we have a natural transformation of $6$-functor formalisms on $(\SchStk_\et,E^{\Alg,\rep}_\pfp)$, of the form
	\[\Psi: \IndShv^! \Rightarrow \Shv^!.\]
	When we restrict this to the subcategory $(\SchStk^\qome_\et)_{P^{\Alg,\rep}_\pfp}\hookrightarrow \Corr(\SchStk_\et,E^{\Alg,\rep}_\pfp)$, this natural transformation takes the form
	\begin{center}
	\begin{tikzcd}
		\IndShv^!(X) \arrow{r}{\Psi_X} \arrow{d}{\Ind(\ul{f_*})}  & \Shv^!(X) \arrow{d}{\ul{f_*}} \\
		\IndShv^!(Y) \arrow{r}{\Psi_Y} & \Shv^!(Y)
	\end{tikzcd}
	\end{center}
	on an arrow $f:X\to Y\in (\SchStk^\qome_\et)_{P^{\Alg,\rep}_\pfp}$, and this square is adjointable. 
	Indeed, as discussed above, the hypothesis $X, Y\in (\SchStk^\qome_\et)$ ensure that each $\Psi_X$, and $\Psi_Y$ admit fully faithful left-adjoint functors $\Psi^L_X$, and $\Psi^L_Y$ respectively, and the hypothesis $f\in P^{\Alg,\rep}_\pfp$ ensures that $\ul{f_*}$ preserves constructibility.
	Since passing to adjoints along adjointable squares is functorial, we get a natural transformation 
	\[\Psi^L:\Shv^!\Rightarrow \IndShv^!\]
	of functors on $(\SchStk^\qome_\et)_{P^{\Alg,\rep}_\pfp}$. 
	Passing to conjugate functors and composing with $c^{*,\vee}_\fg$ we obtain the desired natural transformation.
\end{proof}

\subsubsection{Finitely generated sheaves for sind stacks}

We now extend the discussion of analytification in the previous context to the context of ind-finitely generated sheaves on ind-placid and sind-placid stacks.
This will in particular give rise to the desired analytification operation in the relevant context of $\ICG$.
We have the following which essentially follows directly from the construciton of $\IndShv^!$.

\begin{definition}[{\cite[Definition 10.160]{Zhu25}}]
	\label{finitely-gen-ind-placid}
	Let $X$ be a quasicompact ind-quasi-placid stack. 
	An object $A\in \Shv^!(X)$ is called \emph{finitely generated} if it is of the form $\ul{i_*} B$ for some pfp closed immersion $i:Z\to X$ such that $Z\in \SchStk_\et^\qpl$ and $B\in \Shv^!_c(Z)$.
	We denote $\Shv^!_\fg(X)\subseteq \Shv^!(X)$ the full subcategory of finitely generated sheaves.
	We denote by $\IndShv^!(X)$ its ind-completion, and we let 
	\[\Psi:\IndShv^!(X)\to \Shv^!(X)\]
	be the ind-extension of the natural inclusion.
\end{definition}

\begin{proposition}[{\cite[Lemma 10.161]{Zhu25}}]\label{prop:ind_fin_gen_sheaves_on_ind_qpl_stack}
		Let $X$ be a quasicompact ind-quasi-placid stack with presentation $X\simeq \colim_{i\in I} X_i$, then the $*$-pushforward defines an equivalence in $\LinCat_\Lambda^\sm$ 
		\[ \varinjlim_{\substack{i \in I \\ \ul{f_*}}} \Shv^!_{c}(X_{i}) \simeq \Shv^!_{\fg}(X). \]
		In other words, we have an equivalence in $\LinCat_\Lambda^\cg$ 
		\[ \varinjlim_{\substack{i \in I \\ \on{Ind}(\ul{f_*})}} \IndShv^!(X_{i}) \simeq \IndShv^!(X). \]
	\end{proposition}

	We wish to have control of $\IndShv^!(X)^\vee$ for $X\in \IndStk_\et^\qpl$. 
	Recall the following statement.

	\begin{proposition}\label{prop:colimit_commutes_dual}
		If we have a filtered colimit $C=\varinjlim_{i\in I} C_i$ of compactly generated categories for which the transition maps preserve compact objects, then $C$ is compactly generated and $C^\vee\simeq \varinjlim_{i\in I} C^\vee_i$, where the transition morphisms  $F_{ij}^{o}: C_{i}^{\vee} \ra C_{j}^{\vee}$ are given by the conjugate functors to the transition functors $F_{ij}: C_{i} \ra C_{j}$.%
	\end{proposition}
    \begin{proof}
        This follows from \cite[Proposition  5.5.7.10,  5.5.7.11]{HTT} and \cite[Proposition 7.3.2]{study-in-DAG-1}.
    \end{proof}

\begin{proposition}\label{prop: indshvfg presentation}
	If $X\in \IndStk_\et^{\qpl}$ and we have a presentation $X=\varinjlim_{i\in I} X_i$ with each $X_i$ quasicompact quasi-placid and transition morphisms being pfp closed immersions, then the following is true
	\begin{enumerate}
		\item We have an identification 
			\[\IndShv^!(X)^\vee\simeq \varinjlim_{\substack{i \in I \\ \on{Ind}({f_!})}} \IndShv^*(X_{i}).\]
		\item If each $\Shv^!(X_i)$ is compactly generated, then $\Shv^!(X)$ is compactly generated, the transition morphisms $\ul{f_*}$ preserve compact objects and therefore gives rise to a conjugate map 
			\[(\Psi^L)^o:\Shv^!(X)^\vee\to \IndShv^!(X)^\vee. \]
	\end{enumerate}
\end{proposition}
\begin{proof}
This follows from \Cref{prop:ind_fin_gen_sheaves_on_ind_qpl_stack}, \Cref{prop:colimit_commutes_dual}, \Cref{more-6-functor-comparisons}(1) and the fact that in this setup the transition maps $\ul{f_*}$ admit a colimit-preserving right-adjoint (i.e., $\ul{f^!}$).
Indeed, the transition maps are pfp closed immersions.
\end{proof}

For sind-placid stacks there is also a theory of finitely generated sheaves together with its ind-extension \cite[\S 10.6.5]{Zhu25}.
\begin{definition}[{\cite[Equation (10.63) and the paragraph below]{Zhu25}}]
	For a sind-placid stack $X$, we let 
	\[\IndShv^!(X):=\varinjlim_{\substack{\on{Ind}(\ul{f_*})\\ V\to X}}\IndShv^!(V)\]
	as $V\to X$ ranges over the ind-pfp morphisms with $V$ placid.
\end{definition}
By \cite[Equation (10.64)]{Zhu25}, for all $X\in \sStk_\et^\pl$, we still have a functor 
\[\Psi:\IndShv^!(X)\to \Shv^!(X).\]
A difficulty of working with finitely generated sheaves on general sind-placid stacks is that one might not be able to access all finitely generated sheaves of $X$ via a single atlas $V\to X$.
Indeed, in general we only have a fully faithful embedding (see \cite[Proposition 10.181, Remark 10.182]{Zhu25})
\[|\IndShv^!(V^\bullet)|\hookrightarrow \IndShv^!(X).\]
For this reason, it is hard to construct analytification of ind-finitely generated cosheaves for general sind-placid stacks.
In contrast, the actual category of sheaves can still be accessed by an atlas. 
Indeed we have formulas (see \cite[Equation (10.61)]{Zhu25})
\[\Shv^!(X)\simeq \on{Tot}(\Shv(V^\bullet))\simeq |\Shv(V^\bullet)|.\]
Since we have only defined analytification for $\Shv^!(X)^\vee$ (and quasicompact quasi-placid stacks) when this category is compactly generated (i.e., when $X\in \SchStk^\qome_\et$), this restriction will naturally carry to the sind-placid stacks that we will be able to consider.

\subsubsection{Sind-$\dagger$-correspondences}
\label{sind-dagger-subsection}
The sheaf theory that Zhu works with in \cite{Zhu25} is designed to interact nicely with respect to ind and sind constructions.
This is not the case for the analytic sheaf theory of \cite{Sch17}. 
For this reason, to discuss analytification in the ind-placid and sind-placid context, it will work better to consider functors that one should think of as being the composition of an analytification with a correspondence.
The subtlety is that we want to consider correspondences that might not be $!$-able for all analytic sheaves, but become $!$-able once we restrict to sheaves that come from analytification.

\begin{definition}
	\label{dagger-correspondence}
	\begin{enumerate}
		\item A $\dagger$-\emph{correspondence} consists of a tuple $(X,Y,\alpha)$ with $Y\in \AnStk_v$, $X\in \resSchStk_{\et}^\qpl$ and $\alpha:X^\dagger\to Y$ a $\calD_{\Lambda}^\an$-$!$-able map. 
		\item We let $\Corr_\dagger^{\qpl}$ denote the category of $\dagger$-correspondences, which we can formally define as the pullback
        \begin{equation*}
            \begin{tikzcd}
		    {\Corr_\dagger^{\qpl}} \arrow[d] \arrow[r] & \Fun(\Delta^1,(\AnStk_v)_{E_{\calD^\an,!}}) \arrow[d, "s^*"] \\
		    (\resSchStk_{\et}^\qpl)_{E^{\Alg,\rep}_\pfp} \arrow[r, "(-)^\dagger"]      & \AnStk_v,                                 
\end{tikzcd}
        \end{equation*}
        where \(s^*\) is induced by the inclusion \(\{0\}\to\Delta^1\).
			Informally, a map of $\dagger$-correspondences $\Theta:(X_1,Y_1,\alpha_1)\to (X_2,Y_2,\alpha_2)$
			consists of maps $f:X_1\to X_2$, $g:Y_1\to Y_2$ together with commutation data for the following diagram 
            \begin{equation}{\label{eqn: commdiagforsinddagger}}
			\begin{tikzcd}
				X_1^\diamondsuit  \arrow{d}{f^\diamondsuit}  & X_1^\dagger \arrow{r}{\alpha_1} \ar{l}{b_{X_1}} \arrow{d}{f^\dagger}  & Y_1 \arrow{d}{g} \\
				X_2^\diamondsuit  & X_2^\dagger  \ar{l}{b_{X_2}} \arrow{r}{\alpha_2} &  Y_2.
			\end{tikzcd}
            \end{equation}
            Note that the left hand square automatically commutes (i.e., comes with canonical commutation data), where as for the right hand square one needs to provide additional data.
		\item With the notation as above, we say that a map of $\dagger$-correspondences $\Theta$ is in $P^\dagger$, if the map $f:X_1\to X_2$ is in $P^{\Alg,\rep}_\pfp$.
			We will write \(\Corr_{P^\dagger}^{\qome}\subset\Corr_\dagger^{\qpl}\) for the (non-full) subcategory of \(\dagger\)-correspondences \((X,Y,\alpha)\) such that \(\Shv^!(X)\) is compactly generated, and whose morphisms lie in $P^\dagger$. 
			Note that in this case the left hand square of the commutative diagram (\ref{eqn: commdiagforsinddagger}) will be Cartesian in light of \Cref{lm:proper_gives_cartesian_diam} and \Cref{lm:proper_gives_cartesian_diam2}.
\item Given $(X,Y,\alpha)\in \Corr_\dagger^{\qpl}$ we define the functor
	\[c^{*,\alpha}_\fg:\IndShv^!(X)^\vee\to \calD^\an_\Lambda(Y)\]
	as the composition $\alpha_!\circ b^*_X\circ c^{*,\vee}_\fg$.
If $\Shv^!(X)$ is compactly generated, then we define
	\[c^{*,\alpha}:\Shv^!(X)^\vee\to \calD^\an_\Lambda(Y)\]
	as the composition $c^{*,\alpha}:=c^{*,\alpha}_\fg\circ (\Psi^L)^o$.

\item  We note that these constructions are natural on $P^\dagger$ and give rise to functors

        \begin{equation*}
		c^{*,\alpha}_\fg\from\Corr_{P^\dagger}^{\qome}\to\Fun(\Delta^1,\LinCat_\Lambda).
        \end{equation*}
	and
        \begin{equation*}
		c^{*,\alpha}\from\Corr_{P^\dagger}^{\qome}\to\Fun(\Delta^1,\LinCat_\Lambda).
        \end{equation*}
	Indeed, this follows from proper base-change applied to the left hand square of \eqref{eqn: commdiagforsinddagger} (which is cartesian since $\Theta \in P^{\dagger}$) together with \Cref{lm:proper_gives_cartesian_diam}, \Cref{lm:proper_gives_cartesian_diam2} and \Cref{more-6-functor-comparisons}.

	\end{enumerate}
\end{definition}

\begin{remark}\label{rem:compatibility_dagger_correspondences}
	Informally, \Cref{dagger-correspondence}.(5) says that for any map $\alpha:X^\dagger \to Y$ we get a functor 
	\[\Shv^!(X)^\vee\xrightarrow{c^{*,\alpha}} \calD^\an_\Lambda(Y),\]
	and for every commutative diagram 
	\begin{center}
	\begin{tikzcd}
		X_1^\dagger \arrow{r}{\alpha_1} \arrow{d}{f^\dagger}  & Y_1 \arrow{d}{g} \\
		X_2^\dagger \arrow{r}{\alpha_2} & Y_2
	\end{tikzcd}
	\end{center}
with $f\in P^{\Alg,\rep}_\pfp$ and $\alpha_1,\alpha_2,g\in E^\an_{\calD^\an_\Lambda,!}$ 
	we also get a commutative diagram
	\begin{center}
	\begin{tikzcd}
		\Shv^!(X_1)^\vee \arrow{r}{c^{*,\alpha_1}} \arrow{d}{(\ul{f_*})^o}  & \calD^\an_\Lambda(Y_1) \arrow{d}{g_!} \\
		\Shv^!(X_2)^\vee \arrow{r}{c^{*,\alpha_2}} & \calD^\an_\Lambda(Y_2).
	\end{tikzcd}
	\end{center}
	
	An important special case of this is when $Y_1=Y_2$
\end{remark}

Before we define analytification for sind-placid stacks, we need the following definition.

\begin{definition}
Let $X\in \sStk^\pl_\et$. 
By a \emph{bounded piece} of $X$ we mean a map $g:V\to X$ that is ind-pfp proper with $V\in \SchStk^\pl_\et$.
\end{definition}

We wish to define analytification functors on $X$ by descending analytification functors defined on its bounded pieces.
The following statement will help us achieve this.

\begin{lemma}
	\label{easy-to-verify-dagger-diamond-thingy}
Suppose that $X\in \sStk^\pl_\et$ and that we have a commutative diagram
\begin{center}
\begin{tikzcd}
   & V_2 \arrow{d}{g_2} \\
	V_1 \ar{ru}{f} \arrow{r}{g_1} & X 
\end{tikzcd}
\end{center}
where each $g_1$ and $g_2$ are bounded pieces of $X$.
The following statements hold.
\begin{enumerate}
	\item The formula $V_i^\dagger \simeq V_i^\diamondsuit\times_{X^\diamondsuit} X^\dagger$ holds.
	\item $f\in P^{\Alg,\rep}_\pfp$.
	\item If $\Shv^!(V_1)$ and $\Shv^!(V_2)$ are compactly generated and $V_2\in \resSchStk_\et$, then $V_1\in \resSchStk_\et$ and there is a natural commutative diagram 
\begin{center}
\begin{tikzcd}
	\Shv^!(V_1)^\vee \arrow{r}{c^{*,\vee}} \arrow{d}{(\ul{f_*})^o}  & \calD_\Lambda^\diamondsuit(V_1) \arrow{d}{f^\diamondsuit_!} \ar{r}{b^*_{V_1}} & \calD_\Lambda^\dagger(V_1) \ar{d}{f^\dagger_!} \\
	\Shv^!(V_2)^\vee  \arrow{r}{c^{*,\vee}} & \calD_\Lambda^\diamondsuit(V_2) \ar{r}{b^*_{V_2}} & \calD_\Lambda^\dagger(V_2). 
\end{tikzcd}
\end{center}
\end{enumerate}
\end{lemma}
\begin{proof}
    For the first point, consider a totally disconnected perfectoid Huber pair \((R,R^+)\).
    It suffices to show that \(V_i(\Spec R^\circ)\cong V_i(\Spec R)\times_{X(R)}X(\Spec R^\circ)\) 
    since this will show the formula we want when we pass to the v-sheafification.
    If \(V_i\to X\) is pfp proper on the nose  (i.e., in $P^{\Alg,\rep}_\pfp$), the claim follows from \Cref{lm:proper_gives_cartesian_diam} and \Cref{lm:proper_gives_cartesian_diam2}.
    In general, one can write \(V_i\times_X \Spec R^\circ=\colim W_j\) where each \(W_j\to \Spec R^\circ\) is pfp proper.
    Any map \(\Spec R \to V_i\times_X \Spec R^\circ \) factors through a \(W_j\) for some $j$ by \cite[Lemma 10.154]{Zhu25}, so the existence of the unique lift $\Spec R^\circ\to V_i\times_X \Spec R^\circ$ extending \(\Spec R \to V_i\times_X \Spec R^\circ \) follows from the pfp proper case.

    For the second point, consider $W=V_1\times_X V_2$. 
    The first projection map exhibits $W$ as ind-pfp proper over $V_1$.
    By \cite[Lemma 10.155]{Zhu25}, we may write $W=\varinjlim_j W_{j,1}$ as filtered colimit along closed immersions where each $W_{j,1}\to V_1$ is pfp proper. 
    In particular, each $W_{j,1}\in \SchStk^{\pl}_\et$.
    We have a similar presentation, $W=\varinjlim_k W_{k,2}$ as a filtered colimit along closed immersions where each $W_{k,2}\to V_2$ is pfp proper and $W_{k,2}\in \SchStk^{\pl}_\et$.
    By \cite[Lemma 10.154]{Zhu25}, these families are cofinal, from which we deduce that for either family the two projection maps to $V_1$ and to $V_2$ are both pfp proper.
    By assumption, there is a section $V_1\to W$ which factors through $W_{k,2}$ for some $k$.
    This exhibits $V_1$ as a closed substack of $W_{k,2}$ which is pfp proper over $V_2$.
    Consequently, the map $f:V_1\to V_2$ is also pfp proper.

    For the last point, the first part follows from \Cref{algebraic spaces are also resilient} (2).
    Also, note that the map \((\Psi^L)^o:\Shv^!(X)^\vee\to\IndShv^!(X)^\vee\) is natural in those \(X\) such that \(\Shv^!(X)\) is compactly generated and those morphisms \(g\from X\to Y\in V_c^!\) for which $\ul{g_*}$ preserve compact objects.
    In our case, since $f$ is proper, $\ul{f_*}$ preserves compact objects, by \Cref{lemma: continuousrightadjointpropermpas}. 
    Indeed, it is left-adjoint to $\ul{f^!}$ which is colimit-preserving.
To prove our claim, it suffices to show that the diagrams 
\begin{center}
\begin{tikzcd}
	\IndShv^!(V_1)^\vee \arrow{r}{c_\fg^{*,\vee}} \arrow{d}{\Ind(\ul{f_*})^o}  & \calD_\Lambda^\diamondsuit(V_1) \arrow{d}{f^\diamondsuit_!}  & \calD_\Lambda^\diamondsuit(V_1) \ar{r}{b^*_{V_1}}\arrow{d}{f^\diamondsuit_!} & \calD_\Lambda^\dagger(V_1) \ar{d}{f^\dagger_!} \\  
	\IndShv^!(V_2)^\vee  \arrow{r}{c_\fg^{*,\vee}} & \calD_\Lambda^\diamondsuit(V_2) & \calD_\Lambda^\diamondsuit(V_2) \ar{r}{b^*_{V_2}} & \calD_\Lambda^\dagger(V_2). 
\end{tikzcd}
\end{center}
commute, but these follow from \Cref{more-6-functor-comparisons} and \Cref{lm:proper_gives_cartesian_diam2}, respectively, as in
\Cref{rem:compatibility_dagger_correspondences}. 
Here we have implicitly used Theorem \ref{thm: DiamondAnalytificationUltimate2} and the assumed resilience of $V_{1}$ and $V_{2}$ to identify the arrows with their $\calD_{\Lambda}^{\an}$-versions as opposed to their $\calD_{\Lambda}^{\uop}$-versions. 
\end{proof}
\begin{definition}
	\label{sind-dagger-correspondence}
	\begin{enumerate}
		\item	A \emph{sind}-$\dagger$-\emph{diagram} is a triple $(X,Y,\alpha)$ where 
			\begin{itemize}
		\item $X\in \sStk_\et^{\pl}$.
		\item $Y\in \AnStk_v$
		\item $\alpha:X^\dagger\to Y$ is a map in $\AnStk_v$. 
		\end{itemize}
\item Given a sind-$\dagger$-diagram and a bounded piece $g:V\to X$, we let $\alpha_V:V^\dagger\to Y$ be the map induced by $\alpha$. 
	In other, words $\alpha_V=\alpha \circ g^\dagger$ (see also \Cref{easy-to-verify-dagger-diamond-thingy}).

\item We say that a bounded piece $g:V\to X$ is $\alpha$\emph{-compatible} if all of the following statements hold.
	\begin{enumerate}
		\item $V\in \resSchStk_\et$ (i.e., $V$ is resilient).
		\item $V\in \SchStk^\qome_\et$ (i.e., $\Shv^!(V)$ is compactly generated)
		\item $\alpha_V$ is $\calD^\an_\Lambda$-$!$-able.
	\end{enumerate}
	For a sind-$\dagger$-diagram $(X,Y,\alpha)$, we let $\calC_\alpha\subseteq (\SchStk_\et)_{/X}$ denote the full subcategory spanned by maps $V\to X$ that are $\alpha$-compatible bounded pieces of $X$.
\item We say that a sind-$\dagger$-diagram $(X,Y,\alpha)$ is a \emph{sind-$\dagger$-correspondence} if the following natural map 
	\[\colim_{[V\to X]\in \calC_\alpha}V\to X\]
is an equivalence of \'etale stacks.
    \end{enumerate}
    \end{definition}

    In what follows, we will define an analytification map
	\[c^{*,\alpha}:\Shv^!(X)^\vee\to \calD^\an_\Lambda(Y)\]
	whenever $(X,Y,\alpha)$ is a sind-$\dagger$-correspondence.

Fix \((X,Y,\alpha)\) to be a sind-\(\dagger\)-correspondence.
By the previous constructions, for any \(\alpha\)-compatible piece \(V\to X\) there is a analytification functor 
\begin{equation}{\label{eqn: AnalytificationonAlphaBoundedPieces}} 
c_{V\to X}^{*,\alpha_V}\from\Shv^!(V)^\vee\to\calD^\an_{\Lambda}(Y), 
\end{equation} 
that has formula $c_{V\to X}^{*,\alpha_{V}}:=\alpha_{V,!}\circ b_V^*\circ c^{*,\vee}$.
    We want to extend this construction to \(X\).
    To make things functorial, we will define such an analytification functor for all \'etale stacks over \(X\), although for most objects the result of this construction will be of no interest to us.

    By definition, \(\calC_\alpha\) is the full subcategory of \((\SchStk_{\et})_{/X}\) spanned by \(\alpha\)-compatible pieces.
    By \Cref{easy-to-verify-dagger-diamond-thingy},
    we have a functor 
    \[\calC_\alpha\to \Corr_{P^\dagger}^{\qome}\]
    described informally as 
    \[V\mapsto [V^\diamondsuit \leftarrow V^\dagger \xrightarrow{\alpha_V} Y].\]
    We can compose this functor with the map $\Corr_{P^\dagger}^{\qome}\to \Fun(\Delta^1,\LinCat_\Lambda)$ of \Cref{dagger-correspondence}.(5).
    This construction gives a functor \(c^{*,\alpha}_{(-)}\from\calC_\alpha\to(\LinCat_\Lambda)_{/\calD_\Lambda^\an(Y)}\).
    We define \[c^{*,\alpha}_{(-)}\from(\SchStk_{\et})_{/X}\to(\LinCat_{\Lambda})_{/\calD_\Lambda^\an(Y)}\]
    as the left Kan extension of this functor along the inclusion $\calC_\alpha\subseteq (\SchStk_\et)_{/X}$.
    \begin{lemma}
	    If $(X,Y,\alpha)$ is a sind-$\dagger$-correspondence, then the evaluation 
	    \[c^{*,\alpha}_{(X\to X)}\in (\LinCat_{\Lambda})_{/\calD_\Lambda^\an(Y)}\] 
	    is an arrow of the form 
	    $\Shv^!(X)\to \calD^\an_{\Lambda}(Y)$.
    \end{lemma}
    \begin{proof}
	    A priori, $c^{*,\alpha}_{(X\to X)}$ takes the form $A\to \calD^\an_{\Lambda}(Y)$ and the content is in showing that $A\simeq  \Shv^!(X)$.
	    Recall the general fact that for a slice category $\calC_{/c}$ the forgetful functor $\calC_{/c}\to \calC$ reflects and preserves colimits. 
	    In particular, it suffices to compute $\varinjlim_{V\in \calC_\alpha} \Shv^!(V)\in \LinCat_\Lambda$.
	    By definition of sind-$\dagger$-correspondence, $\varinjlim_{V\in \calC_\alpha} V\simeq X$ in $\SchStk_\et$.
	    The claim now follows from \Cref{lemma: LurieColimits}, \Cref{prop: Shvshriektheoryisnice} and the fact that $\Shv^!$ is a sheaf for the \'etale topology (see \cite[Proposition 10.74.(2)]{Zhu25}).
    \end{proof}
   If the context is clear, we abbreviate the terminology and write 
    \begin{equation}\label{sind dagger analytifiction}
        c^{*,\alpha}\from\Shv^!(X)^\vee\to\calD^\an_\Lambda(Y)
    \end{equation}
    for the map attached to the sind-\(\dagger\)-correspondence \((X,Y,\alpha)\) that we just constructed.
    It can explicitly be described as the colimit of maps
	\[\varinjlim_{\substack{(\ul{f_*})^o\\ V\to X}}[ \Shv^!(V)^\vee\xrightarrow{c^{*,\alpha_V}} \calD^\an_\Lambda(Y)],\]
    where $V \to X$ ranges over a cofinal subset of the $\alpha$-compatible bounded pieces and the maps are defined as in \Cref{eqn: AnalytificationonAlphaBoundedPieces}. 

\begin{remark}
    If \(V_1\to X\) and \(V_2\to X\) are \(\alpha\)-compatible bounded pieces, then any map \(V_1\to V_2\) is proper by \Cref{easy-to-verify-dagger-diamond-thingy}(2).
    In particular, \(\underline{f_*}\) preserves compact objects, and its right adjoint is $\underline{f^!}$. 
    It follows from \cite[Lemma 5.5.7.6]{HTT} that \(\Shv^!(X)\) is compactly generated whenever $X$ forms part of a sind-$\dagger$-correspondence $(X,Y,\alpha)$. 
    Indeed, by assumption, $\Shv^!(V_1)$ is compactly generated whenever $V_1$ is $\alpha$-compatible, and as \cite[Lemma 5.5.7.6]{HTT} explains, limits of compactly generated categories along right adjoint functors still gives rise to compactly generated categories.
\end{remark}

By construction, analytification for a sind-$\dagger$-correspondence is computed as a colimit over all $\alpha$-compatible bounded pieces.
Often, this colimit can be simplified.

\begin{proposition}
	\label{two-analytif-quasi-placid}
	If $X\in \resSchStk_{\et}^\qome$ and $(X,Y,\alpha)$ is a $\dagger$-correspondence, then it is a sind-$\dagger$-correspondence and the two functors $c^{*,\alpha}$ defined in \eqref{sind dagger analytifiction} and \Cref{dagger-correspondence} agree. 
\end{proposition}
\begin{proof}
	Under the assumptions, \(\id\from X\to X\) is an \(\alpha\)-bounded piece, and the colimit computing (\ref{sind dagger analytifiction}) is taken over a diagram with a final object.
\end{proof}

\begin{remark}
For a given sind-$\dagger$-correspondence $(X,Y,\alpha)$, morally, 
\[``c^{*,\alpha}=\alpha_!\circ b^*_X\circ c^{*,\vee}"\]
and this is indeed true in the setup of \Cref{two-analytif-quasi-placid}.
Nevertheless, we are defining $c^{*,\alpha}$ even in situations where $\alpha_!$ might not exist.

Notice that even if $X\in \resSchStk_{\et}^\qome$ a sind-$\dagger$-correspondence $(X,Y,\alpha)$ might not come from $\dagger$-correspondence (i.e., $\alpha_!$ might not exist).
This could happen in circumstances where $V\to X$ is an $\alpha$-compatible bounded piece which is a surjection for the \'etale topology, but $V^\dagger \to X^\dagger$ is not a $\calD^\an$-$!$-cover.
However, this type of example will not really show up for us in practice. 
\end{remark}

\begin{proposition}\label{prop: analytification ind quasi-placid}
	Let $X\in \IndStk_{\et}^{\pl}$, with an ind-presentation $X\simeq \varinjlim X_i$, and let $(X,Y,\alpha)$ be a sind-$\dagger$-diagram. 
	Suppose that each $X_i$ is $\alpha$-compatible and $\alpha_{i} := \alpha_{X_{i}}$, then the triple $(X,Y,\alpha)$ is a sind-$\dagger$-correspondence and in this case, we have the following formula for analytification
	\[c^{*,\alpha}\simeq \varinjlim_{i\in I} c^{*,\alpha_i}.\]
\end{proposition}
\begin{proof}
	By construction, the functor \(c^{*,\alpha}_{(-)}\from(\SchStk_\et)_{/X}\to(\LinCat_{\Lambda})_{/\calD^\an_{\Lambda}(Y)}\) is a colimit over the $\alpha$-bounded pieces $V\in \calC_\alpha$.
	By \cite[Lemma 10.154]{Zhu25}, the family $\{X_i\}_{i\in I}$ is cofinal in $\calC_\alpha$. 
	In particular, \[\colim_{V\in \calC_\alpha} V \simeq \colim_{i\in I} X_i \simeq X,\]
	this shows that $(X,Y,\alpha)$ is a sind-$\dagger$-correspondence and that the analytification formula holds.
\end{proof}

\begin{proposition}
	\label{criterion-for-sind}
Let $(X,Y,\alpha)$ be sind-$\dagger$-diagram.
Let $f:V\to X$ be a ind-pfp proper atlas with $V\in \IndStk^{\pl}_\et$.
Suppose that $V$ admits an ind-presentation $V\simeq \varinjlim_{i\in I} V_i$ such that each $V_i$ is $\alpha$-compatible. 
Suppose that for every $\alpha$-compatible bounded piece $W\in \calC_\alpha$ we have that $\Shv^!(V\times_X W)$ is compactly generated. 
Then $(X,Y,\alpha)$ is a sind-$\dagger$-correspondence, each $(V^n,Y,\alpha_{V^n})$ is a sind-$\dagger$-correspondence and we have a colimit formula
	\[c^{*,\alpha}\simeq \varinjlim_{n\in \Delta^\op} c^{*,\alpha_{V^n}}.\]
	Here, $V^n$ denotes the $n$-fold fiber product of $V$ over $X$, and the colimit is taken in $(\LinCat_\Lambda)_{/\calD^\an_\Lambda(Y)}$.
\end{proposition}
\begin{proof}
	By \Cref{prop: analytification ind quasi-placid}, $(V,Y,\alpha_V)$ is a sind-$\dagger$-correspondence, let us check that $(V^n,Y,\alpha_{V^n})$ is also sind-$\dagger$-correspondences for all $n$.
	Note that every $V^n$ admits an ind-presentation of the form $V^n=\varinjlim_{j\in J_n} W^n_j$ where each $W^n_j\to V$ factors through a pfp proper map $W^n_j\to V_i$ for some $i$ (see \cite[Lemma 10.155]{Zhu25}).
	In particular, by \Cref{algebraic spaces are also resilient}.(2), $W^n_j\in \resSchStk_\et$. 
	We show, inductively, that each $\Shv^!(W^n_j)$ and $\Shv^!(V^n)$ are compactly generated for all $n\in \bbN$, and that the $W^n_j$ are $\alpha$-compatible bounded pieces.
	Indeed, for $n=1$ this is part of the hypothesis. 
	For the inductive step, we see that $V^{n+1}=\colim_{j\in J_n} W^n_j\times_X V$. 
	The inductive hypothesis states that each $W^n_j$ is $\alpha$-compatible and by our assumptions it then follows that $\Shv^!(W^n_j\times_X V)$ is compactly generated.  
	By \Cref{prop:colimit_commutes_dual}, $\Shv^!(V^{n+1})$ is compactly generated and applying \Cref{pass-to-bounded-piece-compacts}, we see that the categories $\Shv^!(W^{n+1}_j)$ are all compactly generated. 

	Since $W^n_j\to V_i$ is pfp, it follows from \cite[Lemma 10.119]{Zhu25} that $W^n_j\in \SchStk^\pl$ which shows that each $W^n_j$ is a bounded piece of $X$. 
By \Cref{thm: DiamondAnalytificationUltimate2}, \Cref{lm:proper_gives_cartesian_diam}, and \Cref{lm:proper_gives_cartesian_diam2}, the map $(W^n_j)^\dagger\to V_i^\dagger$ is $!$-able. 
	This, combined with the hypothesis that $V_i$ is $\alpha$-compatible, shows that $W^n_j$ is $\alpha_{V^n}$-compatible (or equivalently, $\alpha$-compatible).
	In particular, each $V^n$ is the colimit over its $\alpha_{V^n}$-compatible pieces.
	Consequently, we can conclude that each $(V^n,Y,\alpha_{V^n})$ is a sind-$\dagger$-correspondence. 
	
	Finally, let us show that $(X,Y,\alpha)$ is a sind-$\dagger$-correspondence. 
	Note that since $\SchStk_\et$ is a topos and since $V\to X$ is, by definition, a surjection for the \'etale topology, we have the identity in $\SchStk_\et$ 
	\begin{equation}
		\label{eq: cech-colimit}
	\colim_{n\in \Delta^\op}V^n\simeq X.
	\end{equation}
	This is not enough, we should also show that the identity 
	\begin{equation}
		\label{eq: correct-colimit}
\colim_{W\in \calC_\alpha}W \simeq X,
	\end{equation}
holds.
Let $\calC_\alpha^+$ denote the category obtained from $\calC_\alpha$ by formally adding a final object which we denote with $\ast$. 
There is a tautological functor $\calC^+_\alpha\to (\SchStk_\et)_X$ with 
\[W\mapsto W
	\text{ and }
\ast\mapsto X.\]
To show that \eqref{eq: correct-colimit} holds is equivalent to showing that the above functor is the left Kan extension of its restriction to $\calC_\alpha$, which is what we will show.
Let $\calU^V_X$ the sieve (in $\calC^+_\alpha$) over $\ast$ generated by $V$.
For all $W\in \calC_\alpha$, we let $\calU^V_W$ be the sieve over $W$ obtained from $\calU^V_X$ by pulling it back along $W\to X$. 
Notice that, since $\calC_\alpha$ does not necessarily have fiber products, it is not necessarily true that $\calU^V_W$ is generated by $V\times_X W$.
Nevertheless, we claim that $V\times_X W$ has an ind-presentation $V\times_X W\simeq \varinjlim_{i\in I} W_{V,i}$ such that $W_{V,i}\in \calC_\alpha$.
Indeed, using \cite[Lemma 10.155]{Zhu25} we may write $V\times_X W\simeq \varinjlim_{i\in I} W_{V,i}$ where $W_{V,i}\to W$ is pfp proper.
Arguing as in the above, we see that each $W_{V,i}\in \calC_\alpha$ since by assumption $\Shv^!(V\times_X W)$ is compactly generated.
Another use of \cite[Lemma 10.155]{Zhu25} allows us to conclude that the sieve $\calU^V_W$ is generated by the $W_{V,i}$. 
Since $\SchStk_\et$ is a topos, colimits in this category are universal. 
It follows from \eqref{eq: cech-colimit} that, for all $W\in \calC_\alpha$, we have
\begin{equation}
	\label{eq: colimit-sieve expression}
X\times_X W\simeq \colim_{n\in \Delta^\op} V^n\times_X W \simeq \colim_{n\in \Delta, (i_1,\dots,i_n)\in I^n} W_{V,i_0}\times_W \dots \times_W W_{V,i_n}.
\end{equation}
We also have that 
\begin{equation}
\colim_{U\in \calU^V_W} U \simeq \colim_{n\in \Delta, (i_1,\dots,i_n)\in I^n} W_{V,i_0}\times_X \dots \times_X W_{V,i_n}.
\end{equation}
Indeed, as in the proof of \cite[Lemma A.4.6]{HeyerMann}, one shows that the index category on the right is cofinal within the index category on the left.
Overall, this gives
\[W\simeq \colim_{U\in \calU^V_W} U .\]
This shows that the tautological functor $\calC_\alpha\to \SchStk_\et$ is the left Kan extension of its restriction to $\calU^V_X\subseteq \calC_\alpha$.
We claim that the tautological functor on $\calC^+_\alpha$ is the left Kan extension of its restriction $\calU^V_X$, which would show our original claim since Kan extensions compose. 
The only thing left to prove is the formula
\[X\simeq \colim_{U\in \calU^V_X} U.\]
By cofinality, this reduces to \eqref{eq: cech-colimit} which we know holds.
This finishes the proof that $(X,Y,\alpha)$ is a sind-$\dagger$-correspondence.

	To show the colimit formula, recall that $c^{*,\alpha}_{(-)}$ was defined as a left Kan extension. 
	In particular, 
	\[c^{*,\alpha}_{(X\to X)}\simeq \colim_{n\in \Delta^\op} c^{*,\alpha}_{(V^n\to X)}.\]
	Now, each term $c^{*,\alpha}_{(V^n\to X)}$ is of the form 
	\[c^{*,\alpha}_{(V^n\to X)}\simeq \varinjlim_{j\in J_n} c^{*,\alpha}_{(W^n_j\to X)}\simeq \varinjlim_{j\in J_n} \alpha_{W^n_j,!} \circ b^*_{W^n_j}\circ c^{*,\vee}\]
	which we may regroup as
	\[\varinjlim_{j\in J_n} \alpha_{W^n_j,!}\circ b^*_{W^n_j}\circ c^{*,\vee}\simeq \varinjlim_{j\in J_n} c^{*,\alpha_{V^n}}_{(W^n_j\to V^n)} \simeq c^{*,\alpha_{V^n}}\]
		by \Cref{prop: analytification ind quasi-placid}. 
\end{proof}

In the very placid case, our assumptions can be simplified. 
\begin{lemma}
	\label{compact generation helper lemma}
	Let $(X,Y,\alpha)$ a sind-$\dagger$-diagram, suppose that $X\in \sStk^\vpl_\et$ and that $V\to X$ is a sind-very placid atlas with $V\in \IndStk^\vpl_\et$. 
	If $W\to X$ is ind-pfp proper, then $\Shv^!(W\times_X V)$ is compactly generated.
	Consequently, for all $W\in \calC_\alpha$ we have that $\Shv^!(V\times_X W)$ is compactly generated. 
\end{lemma}
\begin{proof}
	Write $V\simeq \varinjlim_{i\in I} V_i$ with each $V_i\in \SchStk^\vpl_\et$.
	By \Cref{prop: indshvfg presentation}(2), it suffices to show that $\Shv^!(V_i\times_X W)$ is compactly generated.
	By \cite[Lemma 10.155]{Zhu25}, we have an ind-presentation $V_i\times_X W\simeq \varinjlim_{k\in K} W_{i,k}$ where each $W_{i,k}$ is pfp proper over $V_i$.
	Again, by \Cref{prop: indshvfg presentation}(2) it suffices to show that the $\Shv^!(W_{i,k})$ are compactly generated. 
	But by \cite[Lemma 10.119]{Zhu25}, the $W_{i,k}$ are quasicompact very placid stacks (i.e. $W_{i,k}\in \SchStk^\vpl$), so the claim follows by \Cref{quasi-placid have compact generated} (2).
\end{proof}

\begin{corollary}\label{criterion-for-sind-very}
Let $(X,Y,\alpha)$ be sind-$\dagger$-diagram.
Let $f:V\to X$ be a ind-pfp proper atlas with $V\in \IndStk^{\vpl}_\et$.
Suppose that $V$ admits an ind-presentation $V\simeq \varinjlim_{i\in I} V_i$ with each $V_i$ $\alpha$-compatible.
Then each $(X,Y,\alpha)$ is a sind-$\dagger$-correspondence, each $(V^n,Y,\alpha_{V^n})$ is a sind-$\dagger$-correspondence and we have a colimit formula
	\[c^{*,\alpha}\simeq \varinjlim_{n\in \Delta^\op} c^{*,\alpha_{V_n}}\]
	in $(\LinCat_\Lambda)_{/\calD^\an_\Lambda(Y)}$.
\end{corollary}
\begin{proof}
This follows from \Cref{compact generation helper lemma} and \Cref{criterion-for-sind}.  	
\end{proof}

Unfortunately, the classifying stacks $\bbB\ul{K}$ of locally profinite groups $K$ are not always sind-very placid and only sind placid \cite[Lemma 3.50]{Zhu25}. 
This will force us to consider cases of \Cref{criterion-for-sind} that are not covered by \Cref{criterion-for-sind-very}. 
Fortunately, $\bbB \ul{K}$ still has good enough ind-atlas to apply \Cref{criterion-for-sind}.

\begin{lemma}
	\label{sind-dagger-for-classifying-stacks}
	Let $(X,Y,\alpha)$ a sind-$\dagger$-diagram.
Suppose that $f:V\to X$ is an ind-pfp proper atlas and that $\ul{f_*}$ is universally conservative. 
	Then, for all $W\in \calC_\alpha$, we still have that $\Shv^!(V\times_X W)$ is compactly generated. 
	In particular, if $V$ additionally admits an ind-presentation $V\simeq \varinjlim_{i\in I} V_i$ such that each $V_i\in \calC_\alpha$, then $(X,Y,\alpha)$ is a sind-$\dagger$-correspondence, each $(V^n,Y,\alpha_{V^n})$ is a sind-$\dagger$-correspondence and we have a colimit formula
	\[c^{*,\alpha}\simeq \varinjlim_{n\in \Delta^\op} c^{*,\alpha_{V^n}},\]
	in $(\LinCat_\Lambda)_{/\calD^\an_\Lambda(Y)}$.

\end{lemma}
\begin{proof}
	The second part follows from the first part and \Cref{criterion-for-sind}. 
	For the first part we argue as above. 
	Namely, we use \cite[Lemma 10.155]{Zhu25} to write $V\times_X W$ as $\varinjlim_{i\in I} W^V_i$ where each $f_i:W^V_i\to W$ is pfp proper. 
	To show that $\Shv^!(V\times_X W)$ is compactly generated, it suffices to show that $\Shv^!(W^V_i)$ is compactly generated (\Cref{prop: indshvfg presentation}(2)). 
	By assumption, $\ul{f_{i,*}}$ is conservative, and by \Cref{placid:have-left-adjoint}, this functor admits a left adjoint $\ul{f^*_i}$.
	If $\{\calO_k\}_{k\in K}$ is a family of compact generators for $\Shv^!(W)$, then $\{\ul{f_i^*}\calO_k\}_{k\in K}$ is a family of compact generators for $\Shv^!(W^V_i)$.
\end{proof}

\subsubsection{An important special case}
There is an important class of sind-placid stacks for which the definition of analytification simplifies significantly.

We fix the setup.
Suppose that $X$ is a sind-very-placid stack, that $V\to X$ is an ind-pfp-proper atlas  and that $V=U/H$ with $U$ an ind-standard-placid algebraic space and $H$ an affine algebraic group as in \Cref{very-placid-givesequiv}. 

\begin{proposition}\label{prop:Shv!_of_sind_placid_stack}
	Let the setup be as above.
Then $\Shv^!(X)$ is compactly generated and there is a natural equivalence 
	\[\Shv^!(X)^\vee\simeq  \varinjlim_{\substack{W\to X \\ f_!}} \Shv^*(W)\]
	as $W\to X$ ranges over the bounded pieces of $X$ that factor as the composition of a closed immersion $W\to V^n$ and the natural map $V^n\to X$, and $f:W_1\to W_2$ denotes a transition map. 
\end{proposition}
\begin{proof}
	We have that 
	\[\Shv^!(X)\simeq \varprojlim_{\substack{W\to X \\ ({\ul{f^!}})}} \Shv^!(W),\]
    by \Cref{lemma: LurieColimits}
	so it follows from \cite[Proposition 1.8.3]{drinfeld2015compactgenerationcategorydmodules} and \Cref{quasi-placid have compact generated} (2) that $\Shv^!(X)$ is dualizable and that 
	\begin{equation}
		\label{easy-part-of-comput}
	\Shv^!(X)^\vee\simeq \varinjlim_{\substack{W\to X \\ ({\ul{f^!}})^\vee}} (\Shv^!(W))^\vee \simeq \varinjlim_{\substack{W\to X \\ (\ul{f_*})^o}} (\Shv^!(W)^\vee). 
	\end{equation}
	Here, the last identity follows from the identification $(\ul{f^!})^\vee\simeq (\ul{f_*})^o$ which follows from the adjunction $\ul{f_*}\dashv \ul{f^!}$ guaranteed by the properness of the transition maps.
    The hypothesis and \Cref{very-placid-givesequiv} give rise to the commutative diagram   
	\begin{center}
	\begin{tikzcd}
		\Shv^!(W)^\vee		\arrow{r}{\simeq } \arrow{d}{(\Psi^L)^o}  & \Shv^*(W)  \\
		\IndShv^!(W)^{\vee}  \arrow{r}{\simeq} & \IndShv^*(W) \arrow{u}{``\colim" \mapsto \colim}.
	\end{tikzcd}
	\end{center}
	This diagram is functorial in $f$ as long as $\ul{f_*}$ preserves compact objects (which is the case for proper morphisms).
	More precisely, if we have maps $W_1\to W_2$ over $X$, then we get commutative diagrams
	\begin{center}
	\begin{tikzcd}
		\Shv^!(W_1)^\vee \arrow{r}{(\ul{f_*})^o} \arrow{d}  & \Shv^!(W_2)^\vee \arrow{d} \\
		\IndShv^*(W_1)\arrow{r}{\on{Ind} f_!} \ar{d}{``\colim"\mapsto \colim} & \IndShv^*(W_2) \ar{d}{``\colim"\mapsto \colim} \\
		\Shv^*(W_1)\arrow{r}{f_!} & \Shv^*(W_2). 
	\end{tikzcd}
	\end{center}
	Indeed, this follows from the definition of the maps involved, \Cref{more-6-functor-comparisons} and the fact that $\ul{f_*}$ preserves compact objects.
	This allows us to rewrite \eqref{easy-part-of-comput} as 
	\[\varinjlim_{\substack{W\to X \\ (\ul{f_*})^o}} (\Shv^!(W)^\vee)\simeq \varinjlim_{\substack{W\to X \\ f_!}} (\Shv^*(W))\]
		as we wanted to show.
\end{proof}

\begin{remark}
	For a sind-$\dagger$-correspondence $(X,Y,\alpha)$ where $X$ is as in \Cref{prop:Shv!_of_sind_placid_stack}, the formula for $c^{*,\alpha}$ simplifies to be 
	\[\varinjlim_{\substack{W\to X \\ f_!}} \alpha_{W,!}\circ b^{\on{an},*}_W.\]
\end{remark}

\subsection{Some computations}
\label{ss: some computations with placid}
Carefully developing the functoriality of sind-$\dagger$-correspondences would also take us too far afield. 
Nevertheless, we will still need some basic functoriality to perform our computations.
Indeed, this will be particularly useful to study how the semi-orthogonal decomposition of $\ICG$ interacts with our functor.

\begin{lemma}
	\label{passing-to-pfp}
	Let $(X,Y,\alpha)$ be a sind-$\dagger$-correspondence.
Suppose that $V\to X$ is a sind-very placid atlas with $V\in \IndStk^\vpl$ and such that $V$ admits an ind-presentation of the form $V=\varinjlim V_i$ with $V_i\in \calC_\alpha$. 
	Suppose that $i:Z\to X$ (resp. $j:U\to X$) is a pfp closed immersion (resp. a pfp open immersion), then $(Z,Y,\alpha\circ i^\dagger)$ (resp. $(U,Y,\alpha\circ j^\dagger)$) is also a sind-$\dagger$-correspondence.
\end{lemma}
\begin{proof}
We use the criterion \Cref{criterion-for-sind-very} for a triple to be a sind-$\dagger$-correspondence to deduce the claim. 	
The proofs for $Z$ and $U$ are analogous, so we only provide the details for one of the arguments.
Let $\alpha^U:=\alpha\circ j^\dagger$.
Then $V\times_X U\to U$ is a sind-very placid atlas with ind-presentation $V\times_X U=\varinjlim V_i\times_X U$.
It suffices to show that $V_i\times_X U\in \calC_{\alpha^U}$.
The hypothesis imply that $V_i\times_X U\in \resSchStk^\vpl$ and consequently $\Shv^!(V_i\times_X U)$ is compactly generated, by Lemma \ref{quasi-placid have compact generated} (2). 
Let $f_i:V_i\times_X U\to U$ denote the evident map.
We must show that the morphism 
\[\alpha^U\circ f^\dagger_i:(V_i\times_X U)^\dagger\to Y\]
is $!$-able. 
This map fits in the following commutative diagram in $\AnStk_{v}$\begin{center}
\begin{tikzcd}
	(V_i\times_X U)^\dagger\arrow{r} \arrow{d}  & \arrow{d} V^\dagger_i \ar{r} & Y  \\
 U^\dagger \arrow{r} & X^\dagger 
\end{tikzcd}
\end{center}
 with Cartesian square. 
 This reduces us to showing that $U^\dagger \to X^\dagger$ is $!$-able.
However, one now observes that this is a proper monomorphism and consequently it is a closed immersion.
Indeed, it is clearly a monomorphism, it is partially proper by definition of $(-)^{\dagger}$ and it is qcqs in light of \Cref{qcqs:is-easy}, so proper by \cite[Proposition~18.3]{Sch17}. 
Since closed immersions are in $E^\rep_\fdcss$, they are $!$-able with respect to $\calD^\an_\Lambda$.
\end{proof}

\begin{lemma}
	\label{also-preserves-compact}
	Let $(X,Y,\alpha)$ be a sind-$\dagger$-correspondence. 
	Let $h:V\to X$ be an $\alpha$-compatible bounded piece.
	Let $f:W\to V$ be a pfp map and let $g=h\circ f$. 
	Suppose that $\Shv^!(W)$ is compactly generated and that $\ul{f_*}$ preserves compact objects, then $\ul{g_*}$ also preserves compact objects and we have a commutative diagram
	\begin{equation}
		\label{functoriality sind dagger 1}
	\begin{tikzcd}
		\Shv^!(W)^\vee \arrow{r}{c^{*,\vee}} \arrow{d}{(\ul{g_*})^o}  & \calD^\diamondsuit_\Lambda(W) \ar{r}{b^*_{W/X}} & \calD^\an_\Lambda(W^\diamondsuit\times_{X^\diamondsuit} X^\dagger \arrow{dl}{\alpha_{W,!}}) \\
		\Shv^!(X)^\vee \arrow{r}{c^{*,\alpha}} & \calD^\an_\Lambda(Y).
	\end{tikzcd}
	\end{equation}
	Here $b_{W/X}:W^\diamondsuit\times_{X^\diamondsuit} X^\dagger \to W^\diamondsuit$ is the natural projection.
\end{lemma}
\begin{proof}
Since $\ul{g_*}=\ul{h_*}\circ \ul{f_*}$, to show that $\ul{g_*}$ preserves compact objects, it suffices to show that $\ul{h_*}$ and $\ul{f_*}$ do.
This holds for $\ul{f_*}$ by hypothesis, that this holds for $\ul{h_*}$ follows from \Cref{lemma: continuousrightadjointpropermpas}.
Indeed, by assumption $h$ is ind-pfp proper, and any functor that admits a colimit-preserving right-adjoint necessarily preserves compact objects.

Now we construct the diagram \eqref{functoriality sind dagger 1}.
	Note that, by \Cref{algebraic spaces are also resilient}, the map $W^\diamondsuit\to V^\diamondsuit$ is in $E^\rep_\fdcss$ and any pullback of it is $!$-able. 
	Consequently, $\alpha_W$ is $!$-able.
	Also, since the map $f:W\to V$ is in $E^\Alg_\pfp$, $\ul{f_*}$ is well defined and preserves constructible objects.
	This gives rise to the following commutative diagram
	\begin{equation}
		\label{functoriality sind dagger 2}
	\begin{tikzcd}
		\IndShv^!(W)^\vee		\arrow{r}{c^{*,\vee}_\fg} \arrow{d}{(\on{Ind} \ul{f_*})^o}  & \ar{r}{b^*_{W/X}} \calD^\diamondsuit_\Lambda(W) \arrow{d}{f^\diamondsuit_!} & \calD^\an_\Lambda(W^\diamondsuit \times_{X^\diamondsuit} X^\dagger)  \ar{d} \\
		\IndShv^!(V)^\vee \arrow{r}{c^{*,\vee}_\fg} & \ar{r}{b^*_{V/X}} \calD^\diamondsuit_\Lambda(V) & \calD^\an_\Lambda(V^\diamondsuit \times_{X^\diamondsuit} X^\dagger),  
	\end{tikzcd}
	\end{equation}
    as in the discussion in \Cref{rem:compatibility_dagger_correspondences}.
	From our assumption on $\Shv^!(W)$ and $f$, we also have the following commutative diagram
	\begin{equation}
		\label{functoriality sind dagger 3}
	\begin{tikzcd}
		\Shv^!(W)^\vee \arrow{r}{\Psi^L} \arrow{d}{(\ul{f_*})^o}  & \IndShv^!(W)^\vee  \arrow{d}{\on{Ind}(f_*^\op)} \\
		\Shv^!(V)^\vee \arrow{r}{\Psi^L}  & \IndShv^!(V)^\vee.
	\end{tikzcd}
	\end{equation}
	The definition of $\alpha$-compatible piece, and of $c^{*,\alpha}$ provides the following commutative diagram
	\begin{equation}
		\label{functoriality sind dagger 4}
	\begin{tikzcd}
		\Shv^!(V)^\vee \arrow{dr}{(\ul{h_*})^o} \ar{r}{\Psi^L} & 		\IndShv^!(V)^\vee \arrow{r}{c_\fg^{*,\vee}}   & \calD^\diamondsuit_\Lambda(V) \ar{r}{b^*_{V/X}} & \calD^\an_\Lambda(V^\diamondsuit\times_{X^\diamondsuit} X^\dagger \arrow{dl}{\alpha_{V,!}}) \\
						  & \Shv^!(X)^\vee \arrow{r}{c^{*,\alpha}} & \calD^\an_\Lambda(Y)
	\end{tikzcd}
	\end{equation}

	The diagram of \eqref{functoriality sind dagger 1} is obtained by combining the diagrams in \eqref{functoriality sind dagger 2}, \eqref{functoriality sind dagger 3} and \eqref{functoriality sind dagger 4}.
\end{proof}

	The following two statements will be crucial to compute the analytification of sheaves on the stacks of shtukas, and it is one of the key inputs in the proof of \Cref{right orthogonality verified}.

	\begin{definition}
		\label{strong-essential-unipotence}
		We say that a map of \'etale stacks $Z\to Y$ is \textit{strongly essentially pro-unipotent} if there is a factorization $Z\to Y'\to Y$ such that $Y'\to Y$ is $\fShv$-suave and $Z\to Y'$ has a presentation $Z\simeq \varprojlim Z_i$ where each $Z_i\to Y$ is $\fShv$-unipotent and each transition map $Z_i\to Z_j$ is in $E^\rep_\pfp$, affine and $\fShv$-unipotent.
	\end{definition}

	\begin{lemma}
		\label{some base change lemma here}
	Let $g:Y\to Y_0$ be a map $\SchStk_\et$ such that there exists surjective map $f:Z\to Y$ of \'etale stacks satisfying that both $h=g\circ f$ and $f$ are strongly essentially pro-unipotent as in \Cref{strong-essential-unipotence}.
	Then $g^\diamond:Y^\diamond\to Y_0^\diamond$ satisfies base-change with respect to $*$-pushforward. 
	\end{lemma}
	\begin{proof}
	Fix $m:X\to Y_0^\diamond$ a map in $\AnStk_v$ and we form the commutative diagram with Cartesian squares 
	\begin{center}
	\begin{tikzcd}
		W_Z \arrow{r}{m''} \arrow{d}{f'}  & \arrow{d}{f^\diamond} Z^\diamond \\
		W \arrow{r}{m'} \arrow{d}{g'}  & \arrow{d}{g^\diamond} Y^\diamond \\
	 X	 \arrow{r}{m} & Y_0^\diamond,
	\end{tikzcd}
	\end{center}
    and set $h' = g' \circ f'$.
We wish to show that, for all $A\in \calD^\an_\Lambda(X)$, the base-change map
	\[g^{\diamond,*}m_*A \xrightarrow{\simeq} m'_*g'^*A\]
	is an isomorphism.
	Since the map $Z^\diamond\to Y^\diamond$ is v-surjective by \Cref{preserves v-covers small d}, $f^{\diamond,*}$ is conservative. 
	So it suffices to show the analogous base-change formulas 
	\[h^{\diamond,*}m_*A \xrightarrow{\simeq} m''_*h'^*A \text{ and } f^{\diamond,*}m'_*A \xrightarrow{\simeq} m''_*f'^*A.\]
	In other words, we are reduced to show base-change formulas for arbitrary strongly essentially pro-unipotent maps $Z\to Y_0$. 
	Choose a factorization $Z\to Y'\to Y_0$ into a $\fShv$-suave part and a pro-unipotent part as in \Cref{strong-essential-unipotence}.
	The $\calD^\uop_\Lambda$-version of these base change formulas follow from \Cref{lemma: Diamondtodiamond} for the smooth part (i.e., $Y'\to Y_0$), and from \Cref{qcqs:is-easy} and \Cref{pro-unip-base change} for the pro-unipotent part (i.e., $Z\to Y'$).
Arguing as in \Cref{precise consequences of smoothness}, we may use the $\calD^\uop_\Lambda$-version to conclude the $\calD^\an_\Lambda$-version. 
	\end{proof}

	\begin{proposition}
		\label{base change formulas}
		Let $f:X\to Y$ be a map in $\resSchStk_\et^{\vpl}$. 
		Suppose that $f$ comes via base change from a map $f_0:X_0\to Y_0$ in $(\PreStk^\Art,E^\rep_{\pfp})$, and that the projection map $g_Y:Y\to Y_0$ satisfies the hypothesis of \Cref{some base change lemma here}.  
		Let $g_X:X\to X_0$ denote the base change map.
		Suppose that $A=g_X^* A_0$ for some $A_0\in \Shv_{c}^{*}(X_{0})$. 
		Let $V=X^\diamondsuit \times_{Y^\diamondsuit} Y^\dagger$ and let $\varphi:V\to Y^\dagger$ and $\eta:V\to X^\diamondsuit$ denote the projection maps.
		Then 
		\[{b_Y^{\on{an},*}}  {f_*}A\simeq  {\varphi_*} {\eta}^* {c^*_X} A 
			\quad \text{ and } \quad 
		{b_Y^{\on{an},*}} {f_!}A\simeq  {\varphi_!} {\eta}^* {c^*_X} A \]
	\end{proposition}
	\begin{proof}
		Consider the following commutative diagram with Cartesian squares
		\begin{center}
		\begin{tikzcd}
			X^\diamondsuit\times_{Y^\diamondsuit} Y^\diamond  \ar[rrr, bend left, "\epsilon"] \ar{r}{\delta}\ar{d}{\phi}	& V  \arrow{r}{\eta} \arrow{d}{\varphi}  &  X^\diamondsuit \arrow{d}{f^\diamondsuit} \ar[r, "g_X^\diamondsuit", swap] &  X_0^\diamondsuit \ar{d}{f_0^\diamondsuit} \\
						   Y^\diamond \ar{r}{a_Y} \ar[rrr, bend right, "\calE"] & Y^\dagger \arrow{r}{b_Y} & Y^\diamondsuit \ar{r}{g_Y^\diamondsuit} &  Y_0^\diamondsuit.
		\end{tikzcd}
		\end{center}
		One of the equations follows from proper base change and the fact that the analytification functor $c^{*}$ is a natural transformation of $6$-functor formalisms. 
		Indeed,
			\begin{align}
			b_Y^{\on{an},*} {f_!}A & \simeq  b_Y^{*}c_Y^* {f_!}A  \\
						 & \simeq b_Y^{*}{f^\diamondsuit_!} c_X^*A \\
						 & \simeq \varphi_! \eta^* c_X^*A .
			\end{align}
			For the other, we argue as follows

\begin{align}
    b_Y^{\on{an},*} {f_*}A & \simeq  b_Y^{*}c_Y^* {f_*}A  \tag{1} \\
    & \simeq b_Y^{*} c_Y^* {f_*}g_X^*A_0 \tag{2} \\
    & \simeq b_Y^{*} c_Y^* g_Y^* {f_{0,*}}A_0 \tag{3} \\
    & \simeq b_Y^{*}  g^{\diamondsuit,*}_{Y} c_{Y_0}^* {f_{0,*}}A_0 \tag{4} \\
    & \simeq b_Y^{*}  g_{Y}^{\diamondsuit,*} {f^\diamondsuit_{0,*}} c_{X_0}^* A_0 \tag{5} \\
    & \simeq a_{Y,*}a_Y^* b_Y^{*}  g_{Y}^{\diamondsuit,*} {f^\diamondsuit_{0,*}} c_{X_0}^* A_0 \tag{6} \\
    & \simeq a_{Y,*} {\phi_*} \epsilon^* c_{X_0}^* A_0 \tag{7} \\
    & \simeq \varphi_* \delta_* \epsilon^* c_{X_0}^* A_0 \tag{8} \\
    & \simeq \varphi_* \delta_* (\delta^* \eta^* g_X^{\diamondsuit,*}) c_{X_0}^* A_0 \tag{9} \\
    & \simeq \varphi_* \eta^* g_X^{\diamondsuit,*} c_{X_0}^* A_0 \tag{10} \\
    & \simeq \varphi_* \eta^* c^*_X g_X^* A_0 \tag{11} \\
    & \simeq \varphi_* \eta^* c^*_X A \tag{12}
\end{align}
The steps $(1),(2),(4),(8),(9),(11),(12)$ are completely formal. 
Step $(3)$ follows from schematic pro-unipotent base change \cite[Proposition 10.124(1)]{Zhu25}. 
Strictly speaking, \cite[Proposition 10.124(1)]{Zhu25} is formulated for $\Shv^!_c$ and functors of the form $\ul{t_!}$ and $\ul{t^!}$.
Nevertheless, a similar (dual) statement with identical proof holds for $\Shv^*_c$ and functors of the form $t_*$ and $t^*$.
To show step $(3)$ we apply \cite[Proposition 10.124(1)]{Zhu25} twice, one time for the map $Z\to Y$ and another time for the map $Z\to Y_0$ where $Z\to Y$ is an auxiliary map as in \Cref{some base change lemma here}. 

Step $(5)$ follows from \Cref{prop:analytification_commutes_!pull}.  
Since $X\in \resSchStk^\vpl$ it is in particular qcqs.
Then step $(6)$ and $(10)$ follow from overconvergent replacement (i.e., \Cref{overconvergent-replacement} and \Cref{overconvergent-replacement-rel}).  
Finally, step $(7)$ follows from pro-unipotent base change.  
More precisely, we have a commutative diagram  
\begin{center}
\begin{tikzcd}
	Y^\diamond \arrow{r}{g_Y^\diamond} \arrow{d}{d_Y}  \ar{rd}{\calE} & Y_0^\diamond \arrow{d}{d_{Y_0}} \\
 Y^\diamondsuit \arrow{r} &Y_0^\diamondsuit,
\end{tikzcd}
\end{center}
and the pullback of $\calE$ along $f_0^\diamondsuit$ is precisely $\epsilon$.
By \Cref{algebraic-2-diamonds}.(1) $d_{Y_0}$ is $\calD^\uop_\Lambda$-suave and $d_{Y_{0}}^{*}$ satisfies base-change against $*$-pushforward, by \Cref{precise consequences of smoothness}.
By \Cref{some base change lemma here}, $g_Y^\diamond$ also satisfies base change against $*$-pushforward.
	\end{proof}

	\section{Semi-orthogonal decompositions}
    \label{sec: B(G)-filt}
    
    \noindent
    \textbf{Soft preamble to the section:}
        In this section, we formulate and prove an inductive criterion to detect when a functor between two categories that are endowed with a semi-orthogonal decomposition is an equivalence.
	We also show that any such equivalence automatically upgrades to a semi-orthogonal equivalence.
	Later in \S \ref{sec: TheEquivalence}, we will verify this criterion for the functor $\pitch$ constructed in \S \ref{sec:construction_functor} using the theory of kimberlites (see \S \ref{sec: kimberlitecalculations}). 
	Our discussion is lean, and relatively self-contained once \cite[Appendix~A.8]{LurieHigherAlgebra} is taken as a given. 
	Nevertheless, we refer the interested reader to \cite{ayala2023stratifiednoncommutativegeometry} where a more thorough treatment, with a different goal than ours, is given. 
\\

    \noindent
    \textbf{Technical preamble to the section:}
    \begin{enumerate}
	    \item[\S \ref{subsec: recollementsandSemiOrthogonalDecompositions}] Recalls recollement as discussed by Lurie in \cite[Appendix~A.8]{LurieHigherAlgebra}.
	    \item[\S \ref{semi: as a 2ff subs}] Defines our working notion of semi-orthogonal decomposition as particular instances of $2$-functor formalisms (as in \cite[Definition 2.15]{dauser2024uniquenesssixfunctorformalisms}).  
		    We find this language convenient since our semi-orthgonal decompositions all come from either a 6-functor formalism, or by passing to adjoints.
		    We also find it convenient to appeal to the Nagata perspective, discussed in \cite{dauser2024uniquenesssixfunctorformalisms}, to obtain uniqueness statements with minimal effort.
	    \item[\S \ref{ss: some basic examples of semi}] Provides the main examples of semi-orthogonal decomposition that we use in the body of the article.
	    \item[\S \ref{ss: here we prove the ind criterion}] Provides the inductive criterion that shows when a functor is an equivalence by exploiting the semi-orthogonal decomposition structure.
		    It is an easy inductive argument taking \cite[Proposition A.8.14]{LurieHigherAlgebra} as the base case.
	    \item[\S \ref{some non-sense about fully faithful maps1}] In this section we benefit from the Nagata perspective in order to show that an equivalence respecting the semi-orthogonal decomposition in a weak sense automatically promotes to a morphism in the category of semi-orthogonal decompositions. 
		    This part is purely formal, and it is a more sophisticated version of the fact that for any category $\calC$ and any object $X\in \calC$, the full subcategory $\calC^{\on{f.f.}}_{/X}$ of $\calC_{/X}$ whose objects consist of those fully faithful maps $Y\to X$ is a posetal category.
    \end{enumerate}

	\subsection{Recollement}
	{\label{subsec: recollementsandSemiOrthogonalDecompositions}}
        In this subsection, we review the theory of recollements on a presentable stable category, following mostly \cite[Appendix~A.8]{LurieHigherAlgebra}.     
       Before defining recollements, we introduce the following notation. 
       We recall that, since the category $\widehat{\Cat}$ admits all small limits by \cite[Corollary 4.2.4.8]{HTT}, we may form $\mathrm{Ker}(f)$ for any exact functor $f \colon \mathcal{C}\to \mathcal{D}$ of stable categories as the pullback of $f$ along $0 \to \mathcal{D}$, the inclusion of the $0$ object inside $\calD$. 
      
	\begin{definition} 
		\label{recollement}
	Let $\mathcal{X} \in \widehat{\Cat}$ be a stable category. 
	We say that a \emph{recollement} of $\mathcal{X}$ is a pair of stable categories $(\mathcal{U},\mathcal{Z})$ together with fully faithful functors $j_{*}: \mathcal{U} \hookrightarrow \mathcal{X}$ and $i_{*}: \mathcal{Z} \hookrightarrow \mathcal{X}$ such that $j_{*}$ and $i_{*}$ admit left adjoints $j^{*}$ and $i^{*}$, and $\mathrm{Ker}(j^*)=\mathrm{Im}(i_*)$. 
	When the context is clear, we say that $(\calZ,\calU)$ is a recollement of $\calX$, and omit $i_*$ and $j_*$ from the notation.
	\end{definition}

	\begin{remark}
		\label{remarks about the verdier quotient2}
		
    The equality $\mathrm{Ker}(j^*)=\mathrm{Im}(i_*)$ is to be interpreted as asking that the natural map $\mathrm{Ker}(j^*)\to \calX$ factors through the full subcategory $\mathrm{Im}(i^*)$ via an equivalence.
    \end{remark}

    \begin{remark}
		\label{remarks about the verdier quotient}
Recall that the Verdier quotient of a fully faithful functor $f\colon \calD \to \calC$ of small stable categories is defined as the Dwyer--Kan localization of \(\calC\) along the class of arrows whose cone lies in \(\calD\), compare with \cite[Theorem I.3.3]{nikolaus-scholze}.
    In the presentable case, the Verdier quotient is the defined as the cofiber in category of presentable stable categories, see \cite[Definition 5.4]{BGT13} and \cite[Theorem I.3.5]{nikolaus-scholze}.
    In particular, $j^*$ realizes $\mathcal{U}$ as the Verdier quotient of $\mathcal{X}$ by the full subcategory $\mathcal{Z}$ embedded along $i_{*}$, compare with \cite[Proposition 5.6]{BGT13}.
	\end{remark}
    
	\begin{remark}{\label{rem: LurieRecollement}}
	We note that $\mathcal{X} \in \widehat{\Cat}$ being a recollement in the above sense implies that it is a recollement in the sense of Lurie \cite[Definition~A.8.1]{LurieHigherAlgebra}.
In fact, since all categories are stable, the categories $\calX$, $\calU$, and $\calZ$ have all finite limits, and the left adjoints $j^{*}$ and $i^{*}$ automatically commute with them, which implies assumptions \cite[Definition~A.8.1 (a)-(c)]{LurieHigherAlgebra}. Moreover, the assumption $\mathrm{Ker}(j^{*})= \mathrm{Im}(i_{*})$ clearly implies \cite[Definition~A.8.1 (d)-(e)]{LurieHigherAlgebra}.
In light of this, we will freely cite results from \cite[Appendix~A.8]{LurieHigherAlgebra} throughout this section. 
	\end{remark}
    \begin{remark}{\label{rem: j!i!exists}}
    Given a recollement $(\calZ,\calU)$ of a stable $\calX$, it follows by \cite[Remark~A.8.5]{LurieHigherAlgebra} that the functor $i_*$ automatically admits a right adjoint denoted $i^!$. 
    Similarly, \cite[Corollary~A.8.13]{LurieHigherAlgebra} states that the functor $j^{*}$ has a fully faithful colimit-preserving left adjoint denoted $j_{!}$. 
    By Yoneda, we note that the identity $\mathrm{Ker}(j^{*}) = \mathrm{Im}(i_{*})$ implies that $\mathrm{Im}(j_{!}) = \mathrm{Ker}(i^{*})$ by passing to adjoints. 
    In particular, the map $i^{*}$ identifies $\calZ$ with the Verdier quotient $\calX/\mathrm{Im}(j_{!})$. 
    Moreover, for all $M\in \calX$, the recollement gives rise to a fiber sequence
    \[j_!j^*M\to M \to i_*i^*M.\]
    Indeed if $F$ denotes the fiber of $M\to i_*i^*M$, then $F\simeq j_!A$ for some $A$, from the identity $\mathrm{Im}(j_{!}) = \mathrm{Ker}(i^{*})$.
    Then it follows that $j^*F\simeq j^*j_!A\simeq A$, but also $j^*F\simeq j^*M$ since $\mathrm{Ker}(j^{*}) = \mathrm{Im}(i_{*})$.
    \end{remark}

    In particular, the previous remark yields a diagram of adjoints
    \[ 
		\begin{tikzcd}
			\mathcal{Z} \arrow[r,hook,"i_{*}"] & \mathcal{X} \arrow[r,"j^{*}"]  \arrow[l, bend right=49,"i^{*}",swap] \arrow[l, bend left=49,"i^{!}"]    & \mathcal{U} \arrow[l, bend left=49,"j_{*}"] \arrow[l, bend right=49,"j_{!}",swap] 
		\end{tikzcd}
    \]
    in $\widehat{\Cat}$, where $(\calZ,\calU)$ is a recollement of a stable $\calX$.

    Next, we will show how the datum of a recollement on $\calX$ can be reconstructed from the datum of the subcategory $i_{*}: \calZ \hookrightarrow \calX$ subject to certain axioms.
To this end, we need to define left and right orthogonals.
    
    \begin{definition}
    Let $i_{*}: \calZ \hookrightarrow \calX$ be a fully faithful embedding of stable categories. 
    The left orthogonal $\phantom{}^{\perp_{\calX}}\calZ$ is the full subcategory whose objects $A \in \calX$ satisfy $\Hom(A,i_{*}(B)) = 0$ for all $B \in \calZ$.
Similarly, the right orthogonal $\calZ^{\perp_{\calX}}$ is the full subcategory whose objects $A\in \calX$ satisfies $\Hom(i_{*}(B),A) = 0$ for all $B \in \calZ$, respectively.
    \end{definition}

    When it is clear from the context, we will often omit the subscript $\calX$.
The fully faithful embeddings $j_!$ and $j_*$ of $\calU$ into $\calZ$ can be reinterpreted via orthogonality as follows.

     \begin{lemma}{\label{lemma: leftandrightorthogonals}}
    Let $\calX$ be a stable category with a recollement $(\calZ,\calU)$.
Then the maps $j_{*}$ and $j_{!}$ define natural equivalences $\calU \xrightarrow{\simeq} \calZ^{\perp}$ and $\calU \xrightarrow{\simeq} \phantom{}^{\perp}\calZ$, respectively.
    \end{lemma}

\begin{proof}
   Note that, for all $A_{U} \in \calU$ and $A_{Z} \in \calZ$, we have  
    \[ \Hom(i_{*}(A_{Z}),j_{*}(A_{U})) = \Hom(j^{*}i_{*}A_{Z},A_{U}) = 0, \]
    so in particular $j_*\colon \calU\to \calX$ has essential image equal to $\calZ^{\perp}$ thanks to the identity $\mathrm{Ker}(j^{*}) = \mathrm{Im}(i_{*})$.
Similarly, the embedding $j_!\colon \calU\to \calX$ has essential image equal to $\phantom{}^{\perp}\calZ$ due to the identity $\mathrm{Im}(j_{!}) = \mathrm{Ker}(i^{*})$ described in \Cref{rem: j!i!exists}.
\end{proof}

    From this point forward, our categories will lie in $\LinCat$, i.e., they will be presentable and stable, but we will still consider diagrams in $\widehat{\Cat}$, i.e., not every functor that we consider between them is necessarily colimit-preserving and exact.
We now isolate our definition of closed subcategory.

    \begin{definition}
	    For $\calX\in \LinCat$, we say that a full subcategory $i_\ast\colon \calZ\to \calX$ (with $i_\ast$ a morphism in $\LinCat$) is closed if $i_\ast$ admits a left adjoint $i^\ast$.    
    \end{definition}

   \begin{lemma}
	   \label{automatic recollement}
       Let $i_\ast\colon \calZ \to \calX$ be a closed subcategory. 
       Then, the functor $i_\ast:\calZ\to \calX$ fits in a recollement $(\calZ,\calU)$ of $\calX$. 
       Moreover, we have natural identifications $\calU\simeq \calZ^\perp$ and $\calU\simeq \calX/\calZ$ where $\calX/\calZ$ denotes the Verdier quotient. 
   \end{lemma}

   \begin{proof}
	   The existence of $\calU$ (and $j_*$) and its identification with $\calZ^\perp$ follows directly from \cite[Proposition A.8.20]{HTT}.
	   Indeed, by hypothesis $i_*$ admits a left adjoint since $i_*$ is closed, and also $i_*$ admits a right adjoint by the adjoint functor theorem (see \cite[Corollary 5.5.2.9]{HTT}). 
We consider the natural fully faithful functor $j_\ast\colon  \calZ^\perp \to \calX $ (which is often not colimit-preserving). 
There is a natural Verdier quotient map $j^\ast \colon \calX \to \calX/\calZ$ with kernel given by $\calZ$. 
Now, \cite[Proposition 5.6]{BGT13} implies that the composition $j^\ast j_\ast:\calZ^\perp\to \calX/\calZ$ is an equivalence, i.e., the Verdier quotient admits a description via a Bousfield localization on $\calZ$-local morphisms (i.e., morphism with whose cofiber lies in $\calZ$). 
This is precisely the identification we were looking for.
   \end{proof}

    In summary, given a closed subcategory $\calZ \hookrightarrow \calX$ of a presentable stable $\calX$, we obtain a recollement $(\calZ,\calX/\calZ)$ of $\calX$ by identifying the right and left orthogonals of $\calZ$ with $\calX/\calZ$.

    \subsection{Semi-orthogonal decompositions as 2-functor formalisms}
    \label{semi: as a 2ff subs}
    The way we want to generalize recollements to more general posets is via Nagata 2-functor formalisms in the sense of \cite{dauser2024uniquenesssixfunctorformalisms}. 
    Recall that a 2-functor formalism is another name for a functor $\Corr(\calC,E)\to \widehat{\Cat}$, where $(\calC,E)$ is a geometric setup: they send a correspondence $[X\leftarrow Z \to Y]$ to $f_!g^*$ and encode base change formulas. 
    However, in contrast to 3-functor formalisms (see \Cref{defn: 3functorformalism}), a 2-functor formalism does not provide a lax symmetric monoidal structure. 
    In particular, this means that there is no projection formula, and concepts like cohomologically proper and cohomologically \'etale become more subtle since the stability of these properties is no longer automatically guaranteed. 

    Recall the notions of \textit{a suitable decomposition} and \textit{a Nagata setup}, that we discussed in \Cref{defn: SuitableDecomposition} and \Cref{defn: nagatasetups}.
    Just as one can define what it meant for a 3-functor formalism to be Nagata with respect to a Nagata setup, one has a similar concept in the 2-functor formalism context (see \cite[Definition 2.15]{dauser2024uniquenesssixfunctorformalisms}). \\

    Before we define semi-orthogonal decompositions, we set some notation related to posets. 
    Throughout this section we will assume that our posets are finite.
    Given a poset $K$, we endow it with the natural order topology with basis given by the open subsets $U_{x} := \{y \in P| y \geq x\}$.
We use $K^{\op}$ to denote the poset with its order reversed.
We write $\Closed_{K}$ for the category of closed subsets of $K$ with morphisms given by inclusion of closed subsets. 
    We say that a subset $V \subset K$ is \emph{convex} if it can be written as the difference $V = Z_{1} \setminus Z_{2}$ of two closed subsets $Z_{1},Z_{2} \subset K$.  Note that if $K$ is finite then, for each $k \in K$, the subset $\{k\} \subset K$ is convex. 
    Again, the collection of convex subsets $V\subset P$ forms a  category $\Convex_K $ with  morphisms given by inclusion.  Note that this carries a natural suitable decomposition $(P,I)$ (in the sense of Definition \ref{defn: SuitableDecomposition}) with $P$ (resp.~$I$) being closed (resp.~open) embeddings $V_1\to V_2$ for the subspace topologies.
    Moreover, $(\Convex_\calP,\on{All},P,I)$ is a Nagata setup (in the sense of \Cref{defn: nagatasetups}.(1)).

\begin{definition}
	\label{def: semi-orthogonal}
	Let $K$ be a poset. 
    A semi-orthogonal decomposition of $\calX \in \LinCat$ with respect to $K$ is a presentable 2-functor formalism \[\calD\colon \Corr(\Convex_K, \on{All})\to \LinCat\] such that 
    \begin{enumerate}
        \item $\calD(K)=\calX$.
        \item It is Nagata (in the sense of \cite[Definition 2.15]{dauser2024uniquenesssixfunctorformalisms}) with respect to the suitable decomposition $(P,I)$ defined above.
        \item For any $i\colon Z\to V $ in $P$, let $j \colon U:=V\setminus Z\to V$ be the complementary open immersion in $I$, then we have excision with respect to the pair $(j,i)$.
More precisely, for any object \(M\in\calD(V)\), the adjunction maps give rise to a fiber sequence
        \begin{equation*}
            j_!j^*M\to M\to i_*i^*M.
        \end{equation*}
    \end{enumerate}
Given such a semi-orthogonal $K$-decomposition $\calD$ on $\calX$, we write $\calX_{V} := \calD(V)$.    
\end{definition}
\begin{remark}{\label{rem: PassageToCOnvexSubsets}}
A convex subset $V \subset K$ is also a poset, and the corresponding topology agrees with the subspace topology on $V$.
We note that we obtain a natural morphism of geometric setups $(\Convex_{V},\mathrm{All}) \ra (\Convex_{K},\mathrm{All})$ such that pulling back $\calD$ along the induced map $\Corr(\Convex_{V},\mathrm{All}) \ra \Corr(\Convex_{K},\mathrm{All})$ gives rise to a semi-orthogonal decomposition $\calD_{V}$ with respect to the poset $V$ on $\calX_{V}$, by design. 
\end{remark}

As a first sanity check, we have the following observation.

    \begin{lemma}{\label{lemma: SemiOrthogonalisRecollement}}
        A semi-orthogonal decomposition $\calD$ of $\calX$ with respect to $\{0\leq 1\}$ is equivalent to the data of a recollement on $\calX$.
More precisely, the closed subcategory pinning down the recollement attached to $\calD$ via Lemma \ref{automatic recollement} is given by $\calZ = \calD(\{0\})$ with embedding $i_{*} := \calD(\{0\} \xleftarrow{\mathrm{id}} \{0\} \ra \{0,1\}): \calZ \ra \calX$.
    \end{lemma}
\begin{proof}
   Let $(\calZ,\calU)$ be a recollement for $\calX$.
We claim that this upgrades to a semi-orthogonal decomposition with respect to the poset $K:=\{ 0\leq 1\}$. Note that $\Convex_{K}=\{ \emptyset, 0, 1, \{0,1\}\}$.
First, we define $\calD_0\colon \Convex_{K} \to \LinCat$ via $\calD_0(\emptyset)=0$, $\calD_0(0)=\calZ$, $\calD_0(\{0,1\})=\calX$, and $\calD_0(1)=\calU $ with maps between them given by $j^*$ and $i^*$.
This is Nagata with respect to the open-closed suitable decomposition because $i$ and $j$ are monomorphisms such that $i_!=i_*$ and $j^!=j^*$ such that we have an identity 
   $\mathrm{Im}(j_{!}) = \mathrm{Ker}(i^{*})$ guaranteeing excision, as already explained in Remark \ref{rem: j!i!exists}. 
   
   Conversely, such a 2-functor formalism allows us to recover the categories $\calZ$ and $\calU$ by reversing the previous steps, together with the triple (by the Nagata condition) adjoint functors $(j_!,j^*,j_*)$ and $(i^*,i_*,i^!)$. 
   Excision implies that $\mathrm{Ker}(j^*)=\mathrm{Im}(i_*)$.
   \end{proof}
   In particular, we deduce the following consequence of this. 
   \begin{corollary}{\label{cor: lotsofrecollements}}
    Let $\calX \in \LinCat$ be a category equipped with a semi-orthogonal decomposition $\calD$ with respect to $K$.
For a convex subset $V$, we set $\calX_{V} := \calD(V)$.
For $i: Z \ra V$ a closed embedding with complementary open $U$, the category $\calX_{V}$ admits a recollement into $(\calX_{Z},\calX_{U})$ where $\calD(Z \xleftarrow{\mathrm{id}} Z \xrightarrow{i} V)$ defines the embedding $\calX_{Z} \hookrightarrow \calX_{V}$ of the closed subcategory characterizing the recollement via Lemma \ref{automatic recollement}.
   \end{corollary}
   \begin{proof}
    We consider the obvious functor $\Convex_{\{0 \leq 1\}} \ra \Convex_{K}$ sending $\{0\}$ to $Z$, $\{1\}$ to $U$,  $\{0,1\}$ to $V$, and $\emptyset$ to $\emptyset$.
This defines a map of geometric setups 
    \[ (\Convex_{\{0 \leq 1\}},\mathrm{All}) \ra (\Convex_{K},\mathrm{All})  \]
    and precomposing $\calD$ with the induced map $\Corr(\Convex_{\{0 \leq 1\}},\mathrm{All}) \ra \Corr(\Convex_{K},\mathrm{All})$ implies the desired claim, by combining with Lemma \ref{lemma: SemiOrthogonalisRecollement}.
   \end{proof}

\begin{remark}
	\label{rem: alternative defininition of semi-orthogonal}
    If we drop the first axiom from \Cref{def: semi-orthogonal}, one can form the category of semi-orthogonal decompositions with respect to $K$ as a full subcategory of the category of presentable 2-functor formalisms, that we will denote by $\on{Semi}_K$.
When the poset is clear from the context, we will omit it from the notation.
    We will write \(\ev_K\from\on{Semi}_K\to\LinCat\) for the forgetful functor given by evaluating on $K \in \Convex_{K}$.
\end{remark}

\subsection{An example}
\label{ss: some basic examples of semi}
Let's exhibit some examples of how semi-orthogonal decompositions will arise for us. 
Fix $X\in \AnStk_v$, a poset $K$, and a continuous map
\[f:|X|\to K.\]
Given $k\in K$, we let $T_k=f^{-1}(k)$.
We say that $X$ is weakly-stratified\footnote{We do not ask that the closure relationships among the $T_k$ agree with $K$.} by $K$ if for all $k\in K$, $T_k$ is a weakly generalizing subset of $X$ (as in \cite[Definition 2.3]{AGLR22}). 
In this circumstance, $T_k\simeq |X_k|$ for a unique (see \cite[Lemma~2.7]{AGLR22}) locally closed immersion $X_k\to X$, which can be defined by the formula 
\[X_k:=X\times_{\ul{|X|}}\ul{T_k}.\]
We set $X_V:=X\times_{|X|} \underline{f^{-1}(V)}$, where $V\subset K$ is an arbitrary convex subset.
Below, we discuss how to promote this to a semi-orthogonal decomposition.

\begin{proposition}
	\label{exceptional anlaytic}
   Let $X \in \AnStk_v$ be weakly-stratified by a poset $K$.
   The following hold:
   \begin{enumerate}
       \item There is a semi-orthogonal $K$-decomposition mapping $V$ to $\calD^\an_\Lambda(X_V)$ and $[U\leftarrow V\to W]$ to $i_{VW!}\circ i_{UV}^*$.
       \item Assume that for all convex subsets $V \subset W\subset K$, the map $i_{VW!}$ (resp.~$i_{VW}^*$) has an additional left adjoint $i_{VW}^\flat$ (resp.~$i_{VW\sharp}$).
Then there is a semi-orthogonal $K^\op$-decomposition mapping $V$ to $\calD^\an_\Lambda(X_V)$ and $[U\leftarrow V\to W]$ to $i_{VW\sharp}\circ i_{UV}^\flat$.
   \end{enumerate}
\end{proposition}

\begin{proof}
The preimage of any set along $\lvert X \rvert \to K$ is weakly generalizing in the sense of \cite{AGLR22}. 
In particular, it follows that $X_{V_1}\to X_{V_2}$ can be written as the composition of an open and a closed immersion. 
These classes of maps are both in $E^\rep_\fdcss$ and consequently $!$-shriekable for the 6-functor formalism $\calD_{\Lambda}^{\an}$.

It's not hard to construct a functor 
\[\Corr(\Convex_K,\on{All})\to \Corr(\AnStk_v,E^\rep_\fdcss)\]
with the rules
\[V\mapsto X_V\]
and 
\[ [U\leftarrow V\to W]\mapsto [X_U\leftarrow X_V \to X_W].\]
Indeed, this follows from functoriality of fiber products and of the construction $(-)\mapsto \ul{(-)}$ that takes topological space to a (not-necessarily small) v-sheaf. 
Composing with $\calD_\Lambda^\an$ gives us a presentable 2-functor formalism on $\Corr(\Convex_K)$. 
To show that it is Nagata, we observe that closed immersions are $\calD_\Lambda^\an$-proper and open immersions are $\calD_\Lambda^\an$-étale.
In order to prove excision, we follow \cite[Proposition 4.8.6]{HeyerMann}: it is enough to show that, for a pair of complementary immersions $(i,j)$, the family of functors $(i^*, j^*)$ is conservative. 
But this is clear, as they cover the underlying topological space.

For the second part, we use part (1) and pass to left adjoint functors to obtain a presentable $2$-functor formalism $\calD^{\on{ex}}:\Corr(\Convex_K)\to \LinCat$ sending $V$ to $\calD^\an_\Lambda(X_V)$ and the correspondence $[U\leftarrow V\to W]$ to $i_{VW\sharp}\circ i_{UV}^\flat$, passing to adjoints via the equivalence \(\LinCat\simeq\LinCat^R\), see \cite[Corollary 5.5.3.4]{HTT}.

We still need to verify that this is a semi-orthogonal decomposition with respect to the opposite category. 
We identify $\Convex_K$ with $\Convex_{K^\op}$, which exchanges closed and open immersions. 
If $V\to W$ is an open immersion in $K$, then $i_{VW}^!=i_{VW}^*$, so passing to left adjoints yields $i_{VW\sharp}=i_{VW!}$, which shows that elements in $I$ are $\calD^{\on{ex}}$-proper, as desired (because the semi-orthogonal decomposition is now indexed by $K^\op$). 
Similarly, if $V\to W$ is a closed immersion in $\Convex_{K}$, then $i_{VW!}=i_{VW\ast}$ yields $i_{VW}^*=i_{VW}^\flat$ by passing to left adjoints. 
This shows that closed immersons in $\Convex_{K}$ are $\calD^{\on{ex}}$-étale, as desired. 
Finally, the excision sequence can be deduced as well by part (1) and by passing to left adjoints.
Indeed, as above it suffices to show that $(i_{VW}^\flat,j_{W\setminus V,W}^\flat)$ is conservative, where $V\to W$ is a closed immersion in $\Convex_{K}$ with complementary open $W\setminus V$. 
If $i_{VW}^\flat A\simeq i_{VW}^* A\simeq 0$ then $A\simeq j_{W\setminus V,W!}B$, for some $B$ and $j_{W\setminus V,W}^\flat A \simeq B$, since $j_{W\setminus V,W!}$ is fully-faithful. 
\end{proof}

The second semi-orthogonal decomposition obtained by passing to left-adjoints will be very useful for our purposes. 

    \subsection{An inductive criterion}
    \label{ss: here we prove the ind criterion}
    In what follows, we provide an inductive criterion showing that if a functor $F:\calX\to \calY$ between two categories that admit a semi-orthogonal $K$-decomposition respects enough of the structure, then the equivalence can be tested on filtered pieces. 
We start by recalling the case of recollement.
    If $(\calX_0,\calX_1)$ is a recollement of $\calX$, we let $(\calX_1)_!\subseteq \calX$ and $(\calX_1)_*\subseteq \calX$ denote the images of $j_!$ and $j_*$ respectively. 
    For notational coherence, we also let $(\calX_0)_*=(\calX_0)_!$ denote the essential image of $i_*:\calX_0\to \calX$.
\begin{lemma}\label{lem_switch_persp_functors}
      Let $\calX$ and $\calY$ be presentable stable categories equipped with recollements $(\calX_0,\calX_1)$ and $(\calY_0,\calY_1)$. 
      Suppose that we have a functor $F\colon \calX\to \calY$ and a commutative diagram
      \begin{center}
      \begin{tikzcd}
	      \calX_0 \arrow{r}{F_0} \arrow{d}{i_*}  & \calY_0 \arrow{d}{i_*} \\
	      \calX \arrow{r}{F} & \calY.
      \end{tikzcd}
      \end{center}
      Then there is a unique (up to contractible choice) functor $\bar{F}:\calX_1\to \calY_1$ fitting in a commutative diagram

      \begin{center}
      \begin{tikzcd}
	      \calX \arrow{r}{F} \arrow{d}{j^*}  & \calY \arrow{d}{j^*} \\
	      \calX_1 \arrow{r}{\bar{F}} & \calY_1.
      \end{tikzcd}
      \end{center}

      In particular, if $F(\calX_{1!})\subset \calY_{1!} $, then we can construct a commutative square 
      \begin{center}
      \begin{tikzcd}
	      \calX_{1!} \arrow{r}{F} \arrow{d}{\simeq,j^*}  & \calY_{1!} \arrow{d}{\simeq,j^*} \\
	      \calX_1 \arrow{r}{\bar{F}} &\calY_1 
      \end{tikzcd}
      \end{center}
      Similarly, if $F(\calX_{1\ast})\subseteq \calY_{1\ast}$, then we can construct a commutative square 
      \begin{center}
      \begin{tikzcd}
	      \calX_{1*} \arrow{r}{F} \arrow{d}{\simeq,j^*}  & \calY_{1*} \arrow{d}{\simeq,j^*} \\
	      \calX_1 \arrow{r}{\bar{F}} & \calY_1
      \end{tikzcd}
      \end{center}
\end{lemma}

\begin{remark}
As we will discuss in \S \ref{some non-sense about fully faithful maps1}, the $F_0$ in \Cref{lem_switch_persp_functors}, if it exists, it is unique up to contractible choice.
\end{remark}

\begin{proof}
	The first statement follows directly from the universal property of Verdier quotients (see \cite[Definition 5.4]{BGT13}), and the fact that the pair $(\calX_1,j^*)$ is a model for the Verdier quotient (see \Cref{remarks about the verdier quotient}).
	Indeed, our assumptions imply that the essential image of $j^*\circ F$ applied to $\calX_0$ is $0\in \calY$.

	The second and third statement follows from concatenating the commutative diagram just obtained with 

      \begin{center}
      \begin{tikzcd}
	      \calX_{1} \arrow{r}{F} \arrow{d}{j_?}  & \calY_{1} \arrow{d}{j_?} \\
	      \calX \arrow{r}{{F}} & \calY
      \end{tikzcd}
      \end{center}
      with $?\in \{!,*\}$.

\end{proof}

\begin{proposition}[{\cite[Proposition A.8.14]{LurieHigherAlgebra}}]
	\label{case of recollement}
   Let $\calX$ and $\calY$ be presentable stable categories equipped with recollements $(\calX_0,\calX_1)$ and $(\calY_0,\calY_1)$.
Suppose $F\colon \calX\to \calY$ is an exact functor such that
   \begin{enumerate}
       \item $F(\calX_{i!})\subset\calY_{i!}$ for $i=0,1$ and the restricted functors are equivalences.
       \item $F(\calX_{i\ast})\subset \calY_{i\ast}$ for $i=0,1$.
   \end{enumerate}
   Then $F$ is an equivalence. 
\end{proposition}
\begin{proof}
   Condition (1) and (3) in \cite[Proposition A.8.14]{LurieHigherAlgebra} hold by hypothesis. 
   As for condition (2) in \cite[Proposition A.8.14]{LurieHigherAlgebra}, we have assumed that $F$ maps $\calX_{i\ast}$ to $\calY_{i\ast}$ for $i=0,1$ and we just need to verify that this happens via an equivalence. 
   For $i=0$, this is already given, since $\calX_{0!}=\calX_{0\ast}$, and in the $i=1$ case, we just remark that the restricted functor on $!$-embedded categories agrees with the one on $\ast$-embedded ones by \Cref{lem_switch_persp_functors}.
   In particular, since we are assuming one of them is an equivalence, the other one is also an equivalence.

   Finally, we must show condition (4) in \cite[Proposition A.8.14]{LurieHigherAlgebra}.
   This amounts to showing that the adjunction map 
   \[F(A)\to F(i^\calX_*i_\calX^*A)\]
   is isomorphic to the adjunction map
   \[F(A)\to i^\calY_*i_\calY^*F(A),\]
   for all $A$ of the form $A=j^\calX_*B$ and $B\in \calX_1$.
   But this is clear since, by hypothesis, $F(j^\calX_!B)\simeq j^\calY_!D$ for some $D$, and we have a fiber sequence 
   \[j^\calY_!D\to F(A)\to F(i^\calX_*i_\calX^*A).\]
   In particular, maps of the form $F(A)\to i^\calY_* T$ factor uniquely through $F(i^\calX_*i_\calX^*A)$, which is precisely the universal property of $i^\calY_*i_\calY^*F(A)$.  
\end{proof}

Now, we generalize \Cref{case of recollement} for a finite poset $K$.
Suppose that $\calX$ is a category in $\LinCat$ endowed with a semi-orthogonal $K$-decomposition.
Given two convex subsets $V\subseteq K$ and $W\subseteq K$ with $V \subset W$, we let $\calX_{VW!}\subseteq \calX_W$ (resp. $\calX_{VW*}\subseteq \calX_W$) denote the essential image of the functor $i_{VW!}$ (resp. $i_{VW*}$) provided by the $2$-functor formalism. 
If $W=K$ we simply write, $\calX_{V!}\subseteq \calX$ and $\calX_{V*}\subseteq \calX$.

\begin{proposition}
	\label{prop: equivalence criterion semi-orthogonal functor many criteria}
   Let $K$ be a finite poset.
Let $\calX$ and $\calY$ be presentable stable categories equipped with semi-orthogonal $K$-decompositions. 
   Suppose $F\colon \calX\to \calY$ is a functor in $\LinCat$ such that the following holds.
   \begin{enumerate}
       \item
       $F(\calX_{V!})\subset \calY_{V!}$ for all convex subsets $V\subset K$.
	  \item[(2)] For any open subset $U\subseteq K$ and any $k\in K$ that is a closed point of $U$, the map $F_U:\calX_U\to \calY_U$ (furnished by \Cref{lem_switch_persp_functors}) satisfies $F_U(\calX_{U\setminus \{k\} U*}) \subseteq \calY_{U\setminus \{k\} U*}$.
          \item[(3)] The functor $F$ maps $\calX_{k!}$ to $\calY_{k!}$ via an equivalence for all $k\in K$. 
     \end{enumerate}
   Then $F$ is an equivalence.
\end{proposition}
\begin{proof}
	By (1), \Cref{lem_switch_persp_functors}, and \Cref{cor: lotsofrecollements}, we may consider for every open $U\subseteq K$ with complement $Z=K\setminus U$ a functor
	\[F_U:\calX_U\to \calY_U.\]
	We prove by induction on the cardinality of $U\subseteq K$ that $F_U$ is an equivalence.
We note that it is a functor between categories with semi-orthogonal decompositions with respect to the poset $U$, by Remark \ref{rem: PassageToCOnvexSubsets}.
We first show that the first and third conditions on $F$ pass to $F_U$.
That $F_U(\calX_{VU!})\subseteq \calY_{VU!}$ whenever $V\subseteq U$ follows easily from the commutative diagram from \Cref{lem_switch_persp_functors}
      \begin{center}
      \begin{tikzcd}
	      \calX \arrow{r}{F} \arrow{d}{j_U^*}  & \calY \arrow{d}{j_U^*} \\
	      \calX_U \arrow{r}{F_U} & \calY_U.
      \end{tikzcd}
      \end{center}
      Indeed, the 2-functor formalisms on $\calX$ and $\calY$ give commutative diagrams
      \begin{center}
      \begin{tikzcd}
	      \calX_V \arrow{r}{i_{VK!}} \arrow{d}{=}  & \calX \arrow{d}{j_U^*} & \calY_V \arrow{r}{i_{VK!}} \arrow{d}{=}  & \calY \arrow{d}{j_U^*} \\
	      \calX_V \arrow{r}{i_{VU!}} & \calX_U & \calY_V \arrow{r}{i_{VU!}} & \calY_U.
      \end{tikzcd}
      \end{center}
      Similarly, from \Cref{lem_switch_persp_functors}, we have a commutative diagram
      \begin{center}
      \begin{tikzcd}
	      (\calX_k)_! \arrow{r} \arrow{d}{F,\simeq}  &(\calX_U)_! \arrow{d}{F} \arrow{r}{\simeq} & \calX_U \arrow{d}{F_U} \\
	      (\calY_k)_!\arrow{r} & (\calY_U)_!  \arrow{r}{\simeq} & \calY_U,
      \end{tikzcd}
      \end{center}
      which shows that the restriction of $F_U$ to $\calX_{kU!}$ induces an equivalence 
      \[F_U:\calX_{kU!}\to \calY_{kU!}.\]

      In the case where $k_{\on{max}}\in K$ is the maximal element, by hypothesis (3) (and \Cref{lem_switch_persp_functors} and \Cref{cor: lotsofrecollements}) we have an equivalence
      \[F_{k_{\on{max}}}:\calX_{k_{\on{max}}}\xrightarrow{\simeq} \calY_{k_{\on{max}}}.\]
	If $U\subseteq K$ is a general open and $k\in U$ is closed, then we have two recollement $(\calX_k,\calX_{U\setminus k})$ and $(\calY_k,\calY_{U\setminus k})$.
	By the inductive hypothesis, $F_{U\setminus k}$ is an equivalence, by $(3)$ $F_k$ is an equivalence, and by $(2)$ the functor $F_U$ satisfies the hypothesis of \Cref{case of recollement}. 
	This shows that $F_U$ is an equivalence. 
\end{proof}

\begin{proposition}
	\label{some auxiliary in-between}
	In the context of \Cref{prop: equivalence criterion semi-orthogonal functor many criteria}, it suffices to show that for $k\in K$ $F(\calX_{k!})\subseteq \calY_{k!}$ to conclude that $F(\calX_{V!})\subseteq \calY_{V!}$ for every convex subset $V\subseteq K$.
\end{proposition}
\begin{proof}
	Observe that $\calX_{V!}$ is generated under colimits by the $\calX_{k!}$ for $k\in V$. 
	The claim follows since $\calY_{V!}$ is stable under colimits for every convex subset $V$ and $F$ commutes with colimits (being a functor in $\LinCat$). 
\end{proof}

We can give a final reformulation of the inductive criterion.
For any $k\in K$, we let $U_{k}$ denote the open subset of elements greater or equal than $k$.
\begin{proposition}
   Let $K$ be a finite poset.
Let $\calX$ and $\calY$ be presentable stable categories equipped with semi-orthogonal $K$-decompositions. 
   Suppose $F\colon \calX\to \calY$ is a functor such that the following holds.
   \begin{enumerate}
	   \item The functor $F$ maps $\calX_{k!}$ to $\calY_{k!}$ via an equivalence for all $k\in K$. 
	   \item For all $k_1\leq k_2$ the functor 
		   \[(i_{k_1U_{k_1}}^\calY)^!\circ F_{U_{k_1 }}\circ i^\calX_{k_2 U_{k_1} *}\simeq 0,\]
		   (i.e., $(i_{k_1U_{k_1}}^\calY)^!\circ F_{U_{k_1}}\circ i^\calX_{k_2 U_{k_1}*}$ identically vanishes).
   \end{enumerate}
   Then $F$ is an equivalence.
\end{proposition}
\begin{proof}
This follows from \Cref{prop: equivalence criterion semi-orthogonal functor many criteria} and \Cref{some auxiliary in-between}.
Indeed, the condition 
\[(i_{k_1U_{k_1}}^\calY)^!\circ F_{U_{k_1 }}\circ i^\calX_{k_2 U_{k_1} *}\simeq 0\]
is, up to simple combinatorics, a rephrasing of condition $(2)$ in \Cref{prop: equivalence criterion semi-orthogonal functor many criteria}.     
\end{proof}

\subsection{From plain functors to semi-orthogonal maps}
\label{some non-sense about fully faithful maps1}
In this subsection we explain why equivalences of semi-orthogonal $K$-decomposed categories that respect the $K$-structure in a weak sense automatically promote uniquely (up to contractible choice) to an equivalence in the category of semi-orthogonal $K$-decompositions.
For this, we will need to understand evaluations from functor categories, since the forgetful functor \(\on{Semi}_K\to\LinCat^R\) is essentially an evaluation functor.

\begin{lemma}
	\label{lemma on Cartesian fibrations}
	Let \(\calI\) be a category with terminal object \(*\), let \(\calC\) be any category admitting finite limits.
    Then \[\ev_*\from\Fun(\calI,\calC)\to\calC\] 
    is a Cartesian fibration.
\end{lemma}
\begin{proof}
	To show it is a Cartesian fibration we must provide $\ev_*$-Cartesian lifts for every edge $x\to y$ in $\calC$, and every functor $F:\calI\to \calC$ with $F(\ast)=y$.
	Evaluation at the terminal object has a right adjoint given by the constant functor.
	More precisely, we have that
	\[\on{Map}(\ev_*(F),z)\simeq \on{Map}(F,\on{const}_z).\]
	By adjunction, for every $F\in \Fun(\calI,\calC)$ we have a map 
	\[F\to \on{const}_{F(\ast)}.\]
	If we are given a morphism $x\to F(\ast)$, then we can form the Cartesian diagram in $\Fun(\calI,\calC)$ 
	\begin{center}
	\begin{tikzcd}
		F_{x\to F(\ast)} \arrow{r} \arrow{d}  & F \arrow{d} \\
	\on{const}_{x} \arrow{r} & \on{const}_{F(\ast)}
	\end{tikzcd}
	\end{center}
	It is not hard to verify from its construction that $F_{x\to F(\ast)}\to F$ is a $\ev_*$-Cartesian lift of $x\to F(\ast)$.
\end{proof}

	Let $K$ be a finite poset.
Let $\calX$ be a category in $\LinCat$ endowed with a semi-orthogonal $K$-decompositions $\calX_{(-)}$.	 Let 
\[\calD_{\calX,*}:\Convex_K \to \LinCat^R\]
be the functor described by the rule 
\[\calD_{\calX,*}(V):=\calX_V\]
and 
\[\calD_{\calX,*}(V\to W):=[f_*:\calX_V\to \calX_W].\]
Formally, this functor is obtained by restricting $\calX_{(-)}$ along the inclusion $\Convex^\op_K \hookrightarrow \corr(\Convex_K,\on{All})$, and passing to right-adjoint functors.

\begin{proposition}
	\label{from * to a true functor}
	Let $K$ be a finite poset.
Let $\calX$ and $\calY$ be two categories in $\LinCat$ endowed with a semi-orthogonal $K$-decompositions.	
Let $F:\calX\to \calY$ be a functor in $\LinCat^R$.
The following statements are equivalent.
\begin{enumerate}
	\item For all $V\in \Convex_K$, we have $F(\calX_{V*})\subseteq \calY_{V*}$. 
	\item There exists a natural transformation $\hat{F}_{(-)}:\calD_{\calX,*}(-)\Rightarrow \calD_{\calY,*}(-)$ with $\hat{F}_{K}\simeq F$.
\end{enumerate}
Moreover, whenever these statements hold, $\hat{F}$ is unique (up to contractible choice).
\end{proposition}
\begin{proof}
	Suppose that we have such a natural transformation $\hat{F}$. 
	Then, for all $V\in \Convex_K$, we have a commutative square
	\begin{center}
	\begin{tikzcd}
		\calX_V \arrow{r}{\hat{F}_V} \arrow{d}{i_{V*}}  & \arrow{d}{i_{V*}} \calY_V \\
	\calX	\arrow{r}{F} & \calY,
	\end{tikzcd}
	\end{center}
	which readily shows that $F(\calX_{V*})\subseteq \calY_{V*}$.
	Conversely, assume that $F(\calX_{V*})\subseteq \calY_{V*}$ holds for all $V$.
	We want to find a lift of $F:\calX\to \calY$ along the map
	\[\ev_K:\Fun(\Convex_K,\LinCat^R)\to \LinCat^R,\]
	and also show that it is unique.
	By \Cref{lemma on Cartesian fibrations}, and the universal property of Cartesian edges, the anima of lifts of $F$ agrees with the anima of maps 
	\[\Fun(\calD_{\calX,*},[\calD_{\calY,*}\times_{\on{const}_\calY}\on{const}_\calX])^\simeq.\]
	After this replacement, we may assume that $\calX=\calY$, and we can now interpret both $\calD_{\calX,*}$ and $\calD_{\calY,*}$ as taking values in $\LinCat^R_{/\calX}$.
	Moreover, since $\calD_{\calX,*}$ and $\calD_{\calY,*}$ come from a semi-orthogonal decomposition, they take values in the full subcategory $(\LinCat^{R})^{\on{f.f.}}_{/\calX}$, where the superscript specifies that these are the objects whose structure map $\calZ\to \calX$ is fully faithful.
	Since the category $(\LinCat^{R})^{\on{f.f.}}_{/\calX}$ is a poset, maps from $\calD_{\calX,*}$ to $\calD_{\calY,*}$ are unique when they exist, and they exist if and only if for all $V\in \Convex_K$, the map 
	\[i_{V*}:\calX_V \to \calX\]
	factors through the map 
	\[i_*:\calY_V\times_{\calY}\calX \to \calX.\]
	Moreover, this happens if and only if $F\circ i^\calX_{V*}$ factors through $i^\calY_{V*}$ which is precisely the condition $F(\calX_{V*})\subseteq \calY_{V*}$.
\end{proof}

\begin{lemma}
	\label{passing from ! to * semiorthogonal decomp equivlence}
	Let $K$ be a finite poset.
Let $\calX$ and $\calY$ be two categories in $\LinCat$ endowed with a semi-orthogonal $K$-decompositions.	
Let $F:\calX\to \calY$ be an equivalence in $\LinCat$.
Suppose that, for all $V\in \Convex_K$, we have that $F(\calX_{V!})\subseteq \calY_{V!}$, then for all $V\in \Convex_K$ we also have that $F(\calX_{V*})=\calY_{V*}$ (i.e., the essential image of $F$ restricted to $\calX_{V*}$ is $\calY_{V*}$).
\end{lemma}
\begin{proof}
	It suffices to show this when $Z\in \Convex_{K}$ is closed, and when $U\in \Convex_{K}$ is open. 
	Indeed, the hypothesis will pass to the Verdier quotient functors. 
	If $Z$ is closed, then by hypothesis $F(\calX_{Z*})\subseteq \calY_{Z*}$.
	Let $A\in \calY_{Z*}$, and $B\in \calX$ with $F(B)\simeq A$.
	We claim that $B\in \calX_{Z*}$.
	Indeed, if $U$ is the open complement of $Z$ in $K$, and we denote by $i:Z\subseteq K$ and $j:U\subseteq K$ the inclusions, then we have a fiber sequence 
	\begin{equation}\label{eq: fiber sequence}
        j_!j^*B\to B\to i_*i^*B,
    \end{equation}
	from which we deduce a fiber sequence
	\[j_!T\to F(B)\to i_*S\]
	for some $S$ and $T$ by applying \(F\) to \cref{eq: fiber sequence} and observing that it preserves !-included subcategories.
	Since $F(B)\in \calY_{Z*}$ it follows that $j_!T\in \calY_{Z*}$ or in other words $j_!T\simeq 0$. 
	Since $F$ is an equivalence we must have $j_!j^*B\simeq 0$ from which we deduce that $B\in \calX_{Z*}$.

	Similarly, if $i:U\to K$ is open, then one can verify that $F(\calX_{U*})\subseteq \calY_{U*}$, by observing that 
	\[\on{Map}(F(i_*B),F(j_*A))\simeq 0.\]
	for all $B\in \calX_{K\setminus U*}$.
	That $F$ restricted to $\calX_{U*}$ essentially surjects onto $\calY_{U*}$ follows a similar proof to the one for closed immersions, using the sequence
    \begin{equation*}
        i_*i^!B\to B\to j_*j^*B 
    \end{equation*}
    instead.
\end{proof}

\begin{proposition}
	\label{thm: equivalence criterion semi-orthogonal functor few criteria}
   Let $K$ be a finite poset.
Let $\calX$ and $\calY$ be presentable stable categories equipped with semi-orthogonal $K$-decompositions. 
   Suppose $F\colon \calX\to \calY$ is an equivalence in $\LinCat$ such that the 
 $F(\calX_{V!})\subset \calY_{V!}$ for all convex subsets $V\subset K$.
Then it lifts uniquely (up to contractible choice) to an equivalence of semi-orthogonal $K$-decompositions.
\end{proposition}
\begin{proof}
Since $F(\calX_{V!})\subseteq \calY_{V!}$, we may use \Cref{passing from ! to * semiorthogonal decomp equivlence} and \Cref{from * to a true functor} to construct an equivalence (with inverse $F^{-1}$) of functors
\[F^{-1}:\calD_{\calY,*}\Rightarrow \calD_{\calX,*},\]
with notation as in \Cref{from * to a true functor}.  
Passing to left-adjoints, gives rise to an equivalence of functors
\[F:\calX\Rightarrow \calY\]
of functors defined on $\Convex_K^\op$ with values in $\LinCat$.
But, since both $\calX$ and $\calY$ are Nagata, by \cite[Theorem 3.3]{dauser2024uniquenesssixfunctorformalisms}, this functor promotes uniquely to an equivalence of $2$-functor formalisms as we wanted to show.
\end{proof}

	\section*{\textbf{Part II. The argument}}
	\section{Constructing the functor}{\label{sec:construction_functor}}
In this section we explain the construction of the functor, but before doing this we carefully explain the key simpler example of classifying stacks. 

For a group object $G$ over $k$, and a topology $\tau$ on $\PSch^{\aff}$ (resp. $\Perf^{\aff}$), we write $\bbB_\tau G \in \SchStk_{\tau}$ (resp. $[\ast/G]_\tau \in \AnStk_\tau$) for the classifying $\tau$-stack of $G$.
We omit $\tau$ from notation if it is clear from the context. 

\subsection{Sheaves on classifying stacks}\label{sec:classifying_stacks}
Let $H$ be a locally profinite group.
We write $\underline{H}$ for the scheme (resp. adic space) representing the functor $S \mapsto C^0(|S|,H)$ on $\PSch^{\rm op}$ (resp. on $\Perf^{\rm op}$) sending $S$ to the set of continuous functions on the underlying topological space of the scheme (resp. adic space) $S$.
In this subsection, we only use the ${\rm fpqc}$-topology on $\PSch$ and the v-topology on the analytic side.
By abuse of notation, we write $\bbB \underline{H}$ for $\bbB_{\rm fpqc} \underline{H}$, the fpqc quotient by $H$ in the category of perfect fpqc-stacks, and $[\ast/\underline{H}]$ for $[\ast/\underline{H}]_v$, the v-stack quotient by $H$ in the category of perfectoid v-stacks.
We apply the above (co)sheaf formalisms to the stacks $\bbB \underline{H}$ and $[\ast/\underline{H}]$ and show that the analytification maps constructed in \S\ref{Dictionary-with-Zhu} are equivalences in this situation. %
We assume that \cite[Assumption 3.48]{Zhu25} holds for $H$ throughout \S \ref{sec:classifying_stacks}. I.e., we assume that
\begin{itemize}
	 \item $H$ admits a $\Lambda$-valued left Haar measure $dh$, that is, an $H$-equivariant map 
		 \begin{equation}\label{eq:assumption_haar_measure}f\mapsto \int_H f dh \colon C_c^\infty(H,\Lambda) \to \Lambda,\end{equation} such that there exists a compact open subgroup $K \subseteq H$ with ${\rm vol}(K) := \int_K dh \in \Lambda^\times$.
\end{itemize}
We call an open compact subgroup $K \subseteq H$ \emph{good}, if $\int_K \mathbbm{1} dh \in \Lambda^\times$ for one (equivalently, any) left Haar measure on $H$.
Thus, our assumption on $H$ guarantees that $H$ contains a good subgroup $K \subseteq H$.
Moreover, note that if $p \in \Lambda^\times$ and $H$ contains an open compact pro-$p$-subgroup, then $H$ satisfies the above assumption.

    We denote by $\Rep(H)$ the derived ($\infty$-)category of the abelian category of smooth $H$-representations on $\Lambda$-modules. 
	It is shown in \cite[Lemma 9.14, Example 3.47]{Zhu25} that, under our assumptions, the category $\Rep(H)$ is left-complete and compactly generated. 
	If $H$ is profinite, let $\Rep_c(H) \subseteq \Rep(H)$ denote the full subcategory of those objects, whose underlying complex of $\Lambda$-modules is perfect.

For profinite groups we have the following.

\begin{proposition}[{\cite[Prop. 10.110, Cor. 10.111, Exmpl. 10.123, and p.95]{Zhu25}}]\label{lm:reps_of_compact_groups}
	 Assume that $H = K$ is profinite and good, then the following statements hold. 
	 \begin{enumerate}
		 \item $\bbB \underline{K}$ is a placid in the sense of \Cref{definition very placid} (2). (i.e., $\bbB \underline{K}\in \SchStk^\pl_\et$).
		 \item There is a canonical t-exact, symmetric monoidal equivalence
	 \[ \Shv_{(c)}^!(\bbB \underline{K}) \simeq \Rep_{(c)}(K) \]
	such that $!$-pullback along the natural map $\Spec k \to \bbB \underline{K}$ identifies with the forgetful functor from $\Rep_{(c)}(K)$ to $\Lambda$-modules.
\item $\Shv^!(\bbB \underline{K})$ is compactly generated and $\Shv^!(\bbB \underline{K})^\omega = \Shv^!_c(\bbB \underline{K})$.
		 \item There is a canonical t-exact, symmetric monoidal equivalence
	 \[ \Shv_{(c)}^*(\bbB \underline{K}) \simeq \Rep_{(c)}(K) \]
	such that $*$-pullback along the natural map $\Spec k \to \bbB \underline{K}$ identifies with the forgetful functor from $\Rep_{(c)}(K)$ to $\Lambda$-modules.
\item $\Shv^*(\bbB \underline{K})$ is compactly generated and $\Shv^*(\bbB \underline{K})^\omega = \Shv^*_c(\bbB \underline{K})$.
\end{enumerate}
\end{proposition}

\begin{remark}
We note that, even though $\Shv^!(\bbB \underline{K})$ is compactly generated when $K$ is profinite and good as in \Cref{quasi-placid have compact generated} (2), often $\bbB \underline{K}$ is in general placid and not very placid.
\end{remark}

\begin{remark}[Very placid case]
Suppose that $H = G(E)$ for a connected linear algebraic group $G$ over $E$, and $K = \calG(\calO_E)$ for a smooth $\calO_E$-model $\calG$ of $G$ with connected special fiber.
Then $\bbB \underline{K}$ is very placid, $\bbB \underline{H}$ is sind-very-placid and $\bbB \underline{K} \to \bbB \underline{H}$ is an (ind-)very placid atlas.
Indeed, $L^+\calG$ is connected and hence $\bbB \underline{K} \cong L^+\calG / \Ad_\varphi L^+\calG$ by the Lang--Steinberg theorem.
As $\ker(L^+\calG \to L^+_1 \calG)$ (the first congruence subgroup) is coh.~pro-unipotent, it follows from Proposition \ref{very-placid-givesequiv} that $\bbB \underline{K}$ is very placid.
Here our loop group notation is described in \S \ref{sec: loop gropus etc} below. 
\end{remark}

\begin{remark}
	\label{clarification Shv BK}
	Although the references provided in \Cref{lm:reps_of_compact_groups} only deal with the case of $\Shv^!(\bbB \underline{K})$ (i.e., $(2)$ and $(3)$) the analogous statements for $\Shv^*(\bbB \underline{K})$ (i.e., $(4)$ and $(5)$) have identical proofs. 
	Indeed, one first shows that the map $\Spec k\to \bbB \underline{K}$ satisfies descent for the $\Shv^*$ functor as in \cite[Proposition 10.110]{Zhu25}\footnote{One should use here that $f:\Spec k\to \bbB \underline{K}$ is cohomologically proper for the $\Shv^*$ functor as discussed in \cite[Remark 5.3.1]{HeyerMann}, and descendability of $f_!\Lambda$ is done following the same argument as in \cite[Proposition 10.110]{Zhu25}.}, then one uses the explicit descent data to compare $\Shv^*(\bbB \underline{K})$ to $\on{QCoh}(\bbB \underline{K}_\Lambda)$ as in \cite[Corollary 10.111 and Example 10.123]{Zhu25}, and one can use this description to identify the compact objects with the constructible ones.
     \end{remark}

Now let $H$ be arbitrary with a good compact open subgroup $K \subseteq H$. 
The fibers $H/K$ of the projection $\bbB \underline{K} \to \bbB \underline{H}$ are discrete and in Zhu's formalism we may treat this map as an ind-pfp finite morphism. 
In particular, it is ind-pfp proper.

\begin{proposition}[{\cite[Lemma 3.50, Proposition 3.51]{Zhu25}}]\label{prop:classifying_stack_is_sind_placid} Let $K \subseteq H$ be a good compact open subgroup. 
	Then the following statements hold.
\begin{enumerate}
\item The morphism $\bbB \underline{K} \to \bbB \underline{H}$ is ind-pfp finite. 
	In particular, $\bbB \underline{H}$ is sind-placid with sind-placid atlas $\bbB \underline{K} \to \bbB \underline{H}$. 
\item Let ${\rm Hk}^\bullet(\bbB \underline{K})$ denote the \v{C}ech nerve of $\bbB \underline{K} \to \bbB \underline{H}$.
Then ${\rm Hk}^n(\bbB \underline{K})$ is ind-placid for each $n\geq 0$, and
	\[ \Shv^!(\bbB \underline{H}) = \varinjlim_{\substack{n\in \Delta^\op \\ \underline{f_\ast}}}  \Shv^!(\Hk^n(\bbB \underline{K})). \]

\item The category $\Shv^!(\bbB \underline{H})$ is compactly generated and there is a $t$-exact symmetric monoidal equivalence
\[\Shv^!(\bbB \underline{H}) \cong {\rm Rep}(H). \]
\item Under the equivalence of (3), $\underline{f^!}$ and $\underline{f_\ast}$ along $f \colon \bbB \underline{K} \to \bbB \underline{H}$ identify with the 
forgetful functor and with ${\on{c-Ind}}_K^H$ compact induction from $K$ to $H$, respectively.
\end{enumerate}
\end{proposition}

Recall from \Cref{construction:from_Shv!dual_to_Shv*} that, for a placid stack $X$, we have natural functors
\[ \Shv^!(X)^\vee \stackrel{(\Psi^L)^o}{\longrightarrow} \IndShv^!(X)^\vee \simeq \IndShv^\ast(X) \to \Shv^\ast(X). \]

\begin{proposition}\label{prop:!dual_and_star_agree_for_classifying_stacks}
Let $K$ be a good subgroup of the locally profinite group $H$. 
Then, for $X \in \{\bbB \underline{K}, \bbB \underline{H}\}$, we have a sequence of equivalences 
\[\Shv^!(X)^\vee\xrightarrow{\simeq} \Shv^*(X)\xrightarrow{\simeq} \calD^\sch_\Lambda(X).\]
\end{proposition}
\begin{proof}
	We first consider the case $X=\bbB \underline{K}$ with $K$ good.
	By \Cref{lm:reps_of_compact_groups}.(3) and \Cref{lm:reps_of_compact_groups}.(5) (see also \Cref{clarification Shv BK}) we have identifications  
\[\Shv^!(\bbB \underline{K})\xleftarrow{\simeq} \IndShv^!(\bbB \underline{K}) \quad \text{ and } \quad \IndShv^*(\bbB \underline{K}) \xrightarrow{\simeq} \Shv^*(\bbB \underline{K}).\]
Overall, this gives an identification 
\[\Shv^!(\bbB \underline{K})^\vee\xrightarrow{\simeq} \Shv^*(\bbB \underline{K}).\]

The identification $\Shv^*(\bbB \underline{K})\xrightarrow{\simeq} \calD_\et(\bbB \underline{K})$ follows from descent.
More precisely, consider the \v{C}ech nerve $(X_n)_{n\geq 0}$ of the map $\Spec k \to \bbB \underline{K}$. 
Then each $X_n$ is a profinite set, and in particular, a qcqs scheme. 
By inspection, one can see that the natural maps $\Shv^*(X_n)\to \calD_\et(X_n) = \calD_{\Lambda}^{\sch}(X_{n})$ are isomorphisms. 
Indeed, one can use \Cref{lemma: BasicPropertiesofDconsDet} and \Cref{identification zhu vs standard} and that each $X_n$ is of finite expansion over $\Spec k$. 
Since $\Spec k\to \bbB \underline{K}$ satisfies descent with respect to $\Shv^*$ by \Cref{clarification Shv BK} and descent with respect to $\calD_\Lambda^{\sch}$ by Theorem \ref{thm: schematicDetisavsheaf}, we also get the desired identification
\[\Shv^*(\bbB \underline{K})\xrightarrow{\simeq} \calD_\Lambda^{\sch}(\bbB \underline{K}).\]
Let us deal with the case $X=\bbB \underline{H}$ where $H$ is a locally profinite group containing an open subgroup $K\subseteq H$ that is good.
In this case, we have to first construct a functor
\[\Shv^!(\bbB \underline{H})^\vee\to \Shv^*(\bbB \underline{H}).\]
Indeed, $\bbB \underline{H}$ is only sind-placid, so its category $\IndShv^!(\bbB \underline{H})$ is more intricately defined, we will not rely on this category. 

We consider the map $f:\bbB \underline{K}\to \bbB \underline{H}$ which is ind-finite.
Let $\{Y_n\}_{n\in \Delta}$ denote the \v{C}ech nerve of $f$. 
We have an identification 
\begin{equation}
	\label{colimit-etale-cover}
\Shv^!(\bbB \underline{H}) \simeq \varinjlim_{\substack{n \in \Delta^\op\\ \ul{f_*}}} \Shv^!(Y_n).
\end{equation}

Applying \cite[Proposition 1.8.3]{drinfeld2015compactgenerationcategorydmodules} to \eqref{colimit-etale-cover} we get an identification  
\begin{equation}
	\label{colimit-etale-cover-dual}
	\Shv^!(\bbB \underline{H})^\vee \simeq \varinjlim_{\substack{n \in \Delta^\op\\ (\ul{f_{\bullet,*}})^o}} \Shv^!(Y_n)^\vee,
\end{equation}
where we have implicitly used the identification $(\ul{f_{\bullet,*}})^o\simeq (\ul{f_\bullet^{!}})^\vee$.
In particular, there is a unique map $\Shv^!(\bbB \underline{H})^\vee\to \Shv^*(\bbB \underline{H})$ that lies over the following commutative square
\begin{center}
\begin{tikzcd}
	\varinjlim_{\substack{n\in \Delta^\op \\ (\ul{f_{\bullet,*}})^o}} \Shv^!(Y_n)^\vee \arrow{r}{\simeq} \arrow{d}{(\ul{f_*})^o,\simeq}  & \varinjlim_{\substack{n\in \Delta^\op \\ {f_!}}} \Shv^*(Y_n) \arrow{d}{f_{\bullet,!}}  \\
	\Shv^!(\bbB \underline{H})^\vee  \arrow{r} & \Shv^*(\bbB \underline{H}).
\end{tikzcd}
\end{center}
To show the map
\begin{equation}
	\label{moreo-etale-comp2}
	\Shv^!(\bbB \underline{H})^\vee  \to \Shv^*(\bbB \underline{H})
\end{equation}
	that we just constructed is an equivalence, it suffices to show that 
	\[\varinjlim_{\substack{n\in \Delta^\op \\ {f_{\bullet,!}}}} \Shv^*(Y_n) \to \Shv^*(\bbB \underline{H})\]
	is an equivalence, since the vertical left hand arrow is an equivalence by \eqref{easy-part-of-comput}.
	Equivalently, using \cite[Corollary 5.5.3.4]{HTT}, it suffices to show that 
	\begin{equation}
		\label{etale-on-star-form}
	\Shv^*(\bbB \underline{H}) \to \varprojlim_{\substack{n\in \Delta \\ {f_\bullet^!}}} \Shv^*(Y_n)
	\end{equation}
	is an equivalence.
	Since $f$ is cohomologically \'etale for the $\Shv^*$ formalism, we may replace \eqref{etale-on-star-form} for  
	\begin{equation}
		\label{etale-on-star-form-2}
	\Shv^*(\bbB \underline{H}) \to \varprojlim_{\substack{n\in \Delta \\ {f_\bullet^*}}} \Shv^*(Y_n).
	\end{equation}
	But this is an equivalence since $\Shv^*$ is an \'etale sheaf and $\bbB \underline{K}\to \bbB \underline{H}$ is surjective for the \'etale topology.
This finishes the proof that \eqref{moreo-etale-comp2} is an equivalence.

	To show the equivalence 
\begin{equation}
	\label{moreo-etale-comp}
	\Shv^*(\bbB \underline{K})  \to \calD_\Lambda^{\sch}(\bbB \underline{K})
\end{equation}
we may argue by descent as above.
Indeed, each of the $Y_n$ are (disjoint unions of) classifying stacks of good compact open subgroups (note that a open subgroup of a good subgroup is again good).
This reduces \eqref{moreo-etale-comp} to the case $H$ is profinite good, which we have dealt with above.
Alternatively, we can appeal to \cite[Proposition 5.3.10]{HeyerMann}.
\end{proof}

\subsubsection{Bernstein--Zelevinsky duality}{\label{ss: BZDuality}}

Now we discuss duality following \cite[\S 3.3.2]{Zhu25}.
Let $H$ be a locally profinite group satisfying assumption \eqref{eq:assumption_haar_measure}. 
Then Bernstein--Zelevinsky duality on $H$, which we will call cohomological duality
\[\bbD_{\coh,H} \colon \Rep(H)^{\op} \to \Rep(H)\] 
is induced by a Frobenius structure on $\Rep(H)$ over $\on{Mod}_\Lambda$ in the sense of \cite[Definition 4.6.5.1]{LurieHigherAlgebra}, which is given by the evaluation map
\begin{equation}\label{eq:BZ_on_RepH}
{\rm coh} \colon \Rep(H) \to \Modd_\Lambda, \quad V \mapsto (V \otimes\Delta_H^{-1})_H,
\end{equation}
where $(\cdot)_H$ is the functor of derived $H$-coinvariants and $\Delta_H$ is the modulus character, given by $\Delta_H(h) = \frac{{\rm vol}(hKh^{-1})}{{\rm vol}(K)}$ for a (any) good subgroup $K \subseteq H$, see \cite[Remark 3.49]{Zhu25}.
More precisely, this means that the composition 
\[ \Rep(H) \otimes_{{\rm Mod}_\Lambda} \Rep(H) \xrightarrow{- \otimes -} \Rep(H) \xrightarrow{\mathrm{coh}} {\rm Mod}_\Lambda \]
is a perfect pairing that defines for us an equivalence
\begin{equation}{\label{eqn: DualityPairingViaFrobeniusStructure}}
 \mathbb{D}_{\coh,H}: \Rep(H)^{\vee} \xrightarrow{\simeq} \Rep(H). 
\end{equation}
Or in other words, a self-duality on $\Rep(H)$.

Since $\Rep(H)$ is compactly generated, this means that given a compact object $A \in \Rep(H)^\omega$ then, we may describe $\bb{D}_{\coh,H}(A)$ as the unique compact object for which the formula
\[ \Hom_{\Rep(H)}(\mathbb{D}_{\coh,H}(A),B) = (A \otimes B \otimes \Delta_{H}^{-1})_{H}, \]
holds for all $B \in \Rep(H)$ (see the discussion in \cite[Remark 7.55]{Zhu25} to see why this identity comes from the Frobenius algebra structure).
We emphasize here that $\Hom$ is the external Hom (i.e., the usual bifunctor $\Rep(H)^\op\times \Rep(H)\to \Modd_\Lambda$).
If we substitute in $B = \mathcal{H}(H) := C^{\infty}_{c}(H,\Lambda)$ to be the compactly supported smooth Hecke algebra then this degenerates into the relationship 
\begin{equation}
	\label{degenerate-hom-relation}
 \Hom_{\Rep(H)}(\mathbb{D}_{\coh,H}(A),\mathcal{H}(H)) =  (A \otimes \mathcal{H}(H) \otimes \Delta_{H}^{-1})_{H} \simeq A, 
\end{equation}
where the last map is the isomorphism in \cite[Proof of Proposition 3.53, Equation 3.40]{Zhu20}. 
Using that $\mathbb{D}_{\coh,H}$ is involutive, we obtain an identification
\begin{equation}
	\label{degenerate-hom-relation-invol}
 \mathbb{D}_{\coh,H}(A) \simeq {\Hom}_{\Rep(H)}(A,\mathcal{H}(H)) 
\end{equation}
for all $A \in \Rep(H)^\omega$. 
A priori, \eqref{degenerate-hom-relation} and \eqref{degenerate-hom-relation-invol} are identities that only hold in $\Modd_\Lambda$, but we note that we can regard $\mathcal{H}(H)$ as a bi-module. 
After taking external Hom with respect to the left action, the complex still carries a remaining right action. 
This gives rise to the involutive equivalence considered first by Bernstein 
\begin{equation}
\label{concrete-duality}
\Hom_{\Rep(H)}(-,\mathcal{H}(H)):\Rep(H)^{\omega,\op}\xrightarrow{\simeq} \Rep(H)^\omega,
\end{equation}
see \cite[\S 2]{MiedaZelevinsky} and \cite[\S 2]{FarguesZelevinsky}.
The duality \eqref{eqn: DualityPairingViaFrobeniusStructure} can then be more concretely realized as the Ind-extension of the equivalence in \eqref{concrete-duality},
via the identification $\Rep(H)^{\vee} \simeq \Ind(\Rep(H)^{\omega,\op})$ of the dual category.

We now discuss how the equivalence $\mathbb{D}_{\coh}$ is incarnated in the geometry of classifying stacks.
Let $K \subseteq H$ be a good subgroup. 
Under the equivalence $\Shv^!(\bbB \underline{K})\simeq \Rep(K)$ from \Cref{lm:reps_of_compact_groups}, the tensor unit $\omega_K$ matches $\Delta_K$.
Moreover, $\omega_K$ is a generalized constant sheaf on $\bbB \underline{K}$ in the sense of \cite[Definition 10.128]{Zhu25}, which we will denote $\Lambda_K^{\rm coh}$. 
As in \cite[Page 98]{Zhu25}, this choice of generalized constant sheaf induces a duality 
\[ 
{\rm id}_{\coh}^K \colon  \Shv^!(\bbB \underline{K}) \stackrel{\sim}{\to} \Shv^!(\bbB \underline{K})^\vee, 
\]
which upon restricting to $\Shv^!_c(\bbB \underline{K})$ and identifying $\Shv^!_{c}(\bbB \underline{K})$ with $\Rep_c K$, becomes the usual duality sending a $K$-representation on a perfect $\Lambda$-module $V$ to its linear dual $V^\ast$. 
If $K$ is a good subgroup of the locally profinite group $H$ we may descend to a duality
\begin{equation}
	\label{descended-duality-from-above}
{\rm id}_{\coh}^H \colon  \Shv^!(\bbB \underline{H}) \stackrel{\sim}{\to} \Shv^!(\bbB \underline{H})^\vee, 
\end{equation}
 by applying these considerations to the \v{C}ech nerve of $\bbB \underline{K} \to \bbB \underline{H}$ as in \cite[Proposition 3.56]{Zhu25}.
 As explained in \cite[Proposition 3.56]{Zhu25} the duality of \eqref{descended-duality-from-above} corresponds to the duality of \eqref{eqn: DualityPairingViaFrobeniusStructure} under the identification of \Cref{prop:classifying_stack_is_sind_placid}.

\subsubsection{$\pitch$ for classifying stacks}

Let $H$ be as above and let $K$ be a good subgroup.
Recall from \cite[Proposition 2.20]{GIZ25} that the natural maps give identifications 
\[ (\bbB\ul{H})^\diamond \simeq (\bbB\ul{H})^\dagger\simeq (\bbB\ul{H})^\diamondsuit\simeq [\ast/\ul{H}].\]
Moreover, recall that when $H$ admits an embedding into $\on{GL}_n(E)$ for some $n$ and some non-Archimedean local field $E$ of residue characteristic $p$, then $[\ast/\ul{H}]$ is an Artin v-stack by \cite[Example IV.1.9]{FS21}.
From \Cref{prop: compatibilityofKanextensionswithDan} and the above isomorphisms of stacks, we 
deduce canonical equivalences 
\begin{equation}\label{eq:identification_categories_classifying_stack}
\calD_\Lambda^{\diamondsuit}(\bbB \underline{H}) \simeq \calD_\Lambda^\an((\bbB \underline{H})^\diamondsuit) \simeq \calD_\Lambda^\an([\ast/\underline{H}]).
\end{equation}

Whenever $K$ is profinite and good, $\bbB \underline{K}\in \SchStk_\et^\pl$ (see \Cref{lm:reps_of_compact_groups}) and $\Shv^{!}(\bbB\underline{K})$ is compactly generated.
In particular, we have an analytification functor
\[ c^{\ast,\vee}_{\bbB \underline{K}} \colon \Shv^!(\bbB \underline{K})^\vee \to D_{\Lambda}^\diamondsuit(\bbB \underline{K}) \simeq D_\Lambda^\an([\ast/\underline{K}]),\] 
as defined in \Cref{analytification-cosheave-quasi-plac}. 

\begin{lemma}\label{lm:analytification_isom_for_profinite_groups}
Let $K$ be profinite and good. 
The analytification functor 
\[ c^{\ast,\vee}_{\bbB \underline{K}} \colon \Shv^!(\bbB \underline{K})^\vee \to \calD_{\Lambda}^\diamondsuit(\bbB \underline{K}) \simeq \calD_\Lambda^\an([\ast/\underline{K}]),\] 
is an equivalence. 
Moreover, we have the commutative diagram of equivalences
\[ \begin{tikzcd} 
	& \Rep(K) \ar[dl,  "\on{Zhu}"'] \ar{dr}{\on{FS}} & \\
	\Shv^!(\bbB \underline{K}) \arrow{r}{\id^K_\coh} & \Shv^!(\bbB \underline{K})^\vee \arrow[r, "c^{\vee,\ast}_{\bbB \underline{K}}"]  & \calD^\an_\Lambda([\ast/\underline{K}]) 
		\end{tikzcd} 
\]
where the left diagonal map is the equivalence of \Cref{lm:reps_of_compact_groups} and the right diagonal map is the equivalence constructed in \cite[Theorem V.1.1]{FS21}.
\end{lemma}
\begin{proof} 
	We have a sequence of maps 
	\begin{equation}
		\label{sequence of comparisons sheaves on BK}
	\Shv^!(\bbB \underline{K})\xrightarrow{\id^K_{\coh}} \Shv^!(\bbB \underline{K})^\vee\to \Shv^*(\bbB \underline{K})\to  \calD^\sch_\Lambda(\bbB \underline{K}) \xrightarrow{c^*}  \calD^\diamondsuit_\Lambda(\bbB \underline{K})=\calD^\an_\Lambda([\ast/\underline{K}])
	\end{equation}
	and all but $c^{*}$ have been shown to be equivalences (see \Cref{prop:!dual_and_star_agree_for_classifying_stacks}).
	Note that $c^{\ast,\vee}_{\bbB \underline{K}}$ is the composition of these maps.
	Consider the fpqc-cover (resp. v-cover) $f:\Spec \overline\bbF_q \to \bbB \underline{K}$ (resp. $f:\ast \to [\ast/\underline{K}]$) and let $\{X_n\}_{n\in \Delta}$ be the \v{C}ech nerve. 
The maps in \eqref{sequence of comparisons sheaves on BK} fit into a diagram

\tikzcdset{scale cd/.style={every label/.append style={scale=#1},
    cells={nodes={scale=#1}}}}
    
\begin{center}
\begin{tikzcd}[column sep=large, row sep=large,scale cd=0.8]
\Shv^{!}(\mathbb{B}\underline{K}) \arrow{r} \arrow{d}
  & \Shv^{!}(\mathbb{B}\underline{K})^{\vee} \arrow{r} \arrow{d}
  & \Shv^{*}(\mathbb{B}\underline{K}) \arrow{r} \arrow{d} 
  & \mathcal{D}^{\mathrm{sch}}_{\Lambda}(\mathbb{B}\underline{K}) \arrow{r} \arrow{d}
  & \mathcal{D}^{\diamondsuit}_{\Lambda}(\mathbb{B}\underline{K}) \arrow{d} \\
\displaystyle
\lim_{\substack{n\in\Delta\\ \underline{f^!_\bullet}}}\Shv^{!}(X_n) \arrow{r}
  & \displaystyle \lim_{\substack{n\in\Delta\\ (\underline{f^!_\bullet})^o}}\Shv^{!}(X_n)^{\vee} \arrow{r}
  & \displaystyle \lim_{\substack{n\in\Delta\\ f_\bullet^*}}\Shv^{*}(X_n) \arrow{r}
  & \displaystyle \lim_{\substack{n\in\Delta\\ f_\bullet^*}} \mathcal{D}^{\mathrm{sch}}_{\Lambda}(X_n) \arrow{r}
  & \displaystyle \lim_{\substack{n\in\Delta\\ f_\bullet^*}} \mathcal{D}^{\diamondsuit}_{\Lambda}(X_n),
\end{tikzcd}
\end{center}
where here the vertical maps are all equivalences, as in the proof of \Cref{prop:!dual_and_star_agree_for_classifying_stacks}, where for $\calD_{\Lambda}^{\Diamond}(-)$ this follows from the fact that $\calD_{\Lambda}^{\an}(-)$ satisfies $v$-descent.
One easily verifies that $\mathcal{D}^{\mathrm{sch}}_{\Lambda}(X_n) \ra \calD_{\Lambda}^{\an}(X_{n})$ is an equivalence, since $X_{n}$ is a profinite set, which shows that $c^{*}$ is an equivalence giving the first part of the claim.

For the second part, if we write $\Rep (K)=\lim_{\substack{n\in\Delta\\ f_\bullet^*}}\on{QCoh}(X_{n,\Lambda})$, then the identifications $\Rep(K) \ra \Shv^{!}(\bbB\underline{K})$ and $\Rep(K) \ra \calD_{\Lambda}^{\an}([\ast/\underline{K}])$ of Zhu and Fargues-Scholze respectively, are induced by taking limits of diagrams of the form
\begin{center}
\begin{tikzcd}
 & \on{QCoh}(X_{n,\Lambda}) \ar{rd} \ar{dl} &  \\
	\Shv^!(X_n)	\arrow{rr} & & \calD^\diamondsuit_\Lambda(X_n)
\end{tikzcd}
\end{center}
for all $n$, and it suffices to show that these diagrams commute.

We claim that this is true for any space of the form $S=\Spec \bar{\bbF}_p\times T$ where $T$ is a profinite set. 
Writing $S=\varprojlim S_i$ with each $S_i=\Spec \bar{\bbF}_p\times T_i$ with $T_i$ finite we have formulas 
\begin{center}
\begin{tikzcd}
 & \on{QCoh}(S_\Lambda)= \varinjlim \on{QCoh}(\Spec \Lambda^{T_i}) \ar{rd} \ar{dl} &  \\
	\Shv^!(S) = \varinjlim	\Shv^!(S_i)	\arrow{rr} & & \calD^\diamondsuit_\Lambda(S)= \varinjlim \calD^\diamondsuit_\Lambda(S_i)
\end{tikzcd}
\end{center}
which can be further rewritten as 
\begin{center}
\begin{tikzcd}
 & \varinjlim \on{QCoh}(\Spec \Lambda)^{\oplus |T_{i}|} \ar{rd} \ar{dl} &  \\
	\varinjlim \Shv^!(\Spec \bar{\bbF}_p)^{\oplus |T_i|}		\arrow{rr} & & \varinjlim \calD^\diamondsuit_\Lambda(\ast)^{\oplus |T_{i}|}.
\end{tikzcd}
\end{center}
The whole computation reduces to the commutativity of the diagram 
\begin{center}
\begin{tikzcd}
 &  \on{QCoh}(\Spec \Lambda) \ar{rd} \ar{dl} &  \\
	\Shv^!(\Spec \overline{\bbF}_p))	\arrow{rr} & & \calD^\diamondsuit_\Lambda(\ast).
\end{tikzcd}
\end{center}
Since the maps are ${\rm Mod}_\Lambda$-linear, it suffices to show that each of the maps preserves the tensor unit, but this can be done by a direct computation.
\end{proof}

We move on to study the case where $H$ is a locally profinite group admitting a compact open subgroup $K\subseteq H$ that is good.

\begin{lemma}
	\label{sind-dagger-classifying}
	The triple $(\bbB \underline{H}, [\ast/\underline{H}], \alpha_H:=\id_{[\ast/\underline{H}]})$ is a sind-$\dagger$-correspondence in the sense of \Cref{sind-dagger-correspondence}.(4).
\end{lemma}
\begin{proof}
Indeed, let $\alpha_{H} =\id_{[\ast/\underline{H}]}$, this is indeed a map 
\[\id_{[\ast/\underline{H}]}:(\bbB \underline{H})^\dagger \xrightarrow{\simeq} [\ast/\underline{H}],\]
that gives rise to a sind-$\dagger$-diagram.
Let $K\subseteq H$ be an open compact good subgroup of $H$.
Then $f:\bbB \underline{K}\to \bbB \underline{H}$ is a sind-placid atlas (see \Cref{prop:classifying_stack_is_sind_placid} (1)), with $\bbB \underline{K} \in \resSchStk^\pl$ (see \Cref{enough-resilients}). 
Note that $(\bbB \underline{K})^\dagger\times_{(\bbB \underline{H})^\diamondsuit} (\bbB \underline{H})^\dagger\simeq [\ast/\underline{K}]$.
Moreover, 
\[[\ast/\underline{K}] \to [\ast/\underline{H}]\]
is \'etale and in particular $\calD^\an_\Lambda$-$!$-able.
This shows that $\bbB \underline{K}\to \bbB \underline{H}$ is an $\alpha_H$-compatible bounded piece as in \Cref{sind-dagger-correspondence}.(3). 
By \Cref{sind-dagger-for-classifying-stacks}, it suffices to show that 
\[\ul{f_*}:\Shv^!(\bbB \underline{K})\to \Shv^!(\bbB \underline{H})\]
	is universally conservative to conclude that $(\bbB \underline{H}, [\ast/\underline{H}], \alpha_H:=\id_{[\ast/\underline{H}]})$ is a sind-$\dagger$-correspondence. 
	Let $X\to \bbB \underline{H}$ be a map, let $Y=X\times_{\bbB \underline{H}}\bbB \underline{K}$, and let $f_X:Y\to X$ be the base change map.
We wish to show that $\ul{f_{X,*}}$ is conservative. 
	We may easily reduce to the case in which $X\in \PSch^\qcqs$, so that $Y\in \PSch$ and $Y\to X$ is ind-finite.
	If $A\in \Shv^!(Y)$ and $A\neq 0$, then there exists a geometric point $\overline{y}:\Spec C \to Y$ with $\ul{\overline{y}^!}A\neq 0$ in $\Shv^!(\Spec C)\simeq {\rm Mod}_\Lambda$.
	Let $\overline{x}:\Spec C\to X$ be the induced geometric point.
	Then $\ul{x^!} \circ \ul{f_{X,*}}A\neq 0$ since one can show, using proper base change, that $\ul{y^!}A$ is a direct summand of $\ul{x^!} \circ \ul{f_{X,*}}A$. 
	In particular, $\ul{f_{X,*}}A\neq 0$.
\end{proof}

As one should expect, the analytification functor constructed in \Cref{sind-dagger-classifying} using a sind-$\dagger$-correspondence admits an easier description. 
    
	\begin{theorem}
		\label{analitifying cInd}
		Consider the sind-$\dagger$-correspondence $(\bbB \underline{H}, [\ast/\underline{H}], \alpha_H:=\id_{[\ast/\underline{H}]})$ as above.
		The following statements hold.
		\begin{enumerate}
			\item The analytification functor 
				\[c^{*,\alpha_H}\colon \Shv^!(\bbB \underline{H})^\vee \to \calD^\diamondsuit_\Lambda(\bbB \underline{H})\]
				factors as $\Shv^!(\bbB \underline{H})^\vee\to \Shv^*(\bbB \underline{H}) \xrightarrow{c^*}\calD_\Lambda^\diamondsuit(\bbB \underline{H})$ and in particular it is an equivalence.
				Here, the map $\Shv^!(\bbB \underline{H})^\vee\to \Shv^*(\bbB \underline{H})$ is the equivalence constructed in \Cref{prop:!dual_and_star_agree_for_classifying_stacks}. 
			\item If $K\subseteq H$ is compact and good, then we have a commutative diagram 
				\begin{center}
				\begin{tikzcd}
					\Shv^!(\bbB \underline{K})^\vee \arrow{r}{c^{*,\alpha_K}} \arrow{d}{(\ul{f_*})^o}  & \calD^\an_\Lambda([\ast/\underline{K}]) \arrow{d}{f_!} \\
					\Shv^!(\bbB \underline{H})^\vee \arrow{r}{c^{*,\alpha_H}} & \calD^\an_\Lambda([\ast/\underline{H}]).
				\end{tikzcd}
				\end{center}
			\item After precomposing with the identification $\id^K_\BZ$ and $\id^H_\BZ$ that come from duality, we get a commutative diagram 
				\begin{center}
				\begin{tikzcd}
					\Rep(K)\ar{r}{\simeq}\ar{d}{\on{c-Ind}_K^H}	& \Shv^!(\bbB \underline{K}) \arrow{r}{\simeq} \arrow{d}{(\ul{f_*})}  & \calD^\an_\Lambda([\ast/\underline{K}]) \arrow{d}{f_!} \\
					\Rep(H)\ar{r}{\simeq}		& \Shv^!(\bbB \underline{H}) \arrow{r}{\simeq} & \calD^\an_\Lambda([\ast/\underline{H}]).
				\end{tikzcd}
				\end{center}
			\item Finally, we have a commutative diagram of equivalences
\[ \begin{tikzcd} 
	& \Rep(H) \ar[dl,  "\on{Zhu}"'] \ar{dr}{\on{FS}} & \\
	\Shv^!(\bbB \underline{H}) \arrow{r}{\id^H_\coh} & \Shv^!(\bbB \underline{H})^\vee \arrow[r, "c^{*,\alpha_H}"]  & \calD^\an_\Lambda([\ast/\underline{H}]).
		\end{tikzcd} 
\]
		\end{enumerate}
	\end{theorem}
	
\begin{proof} 
	Let us explicate the \v{C}ech nerve appearing in \Cref{prop:classifying_stack_is_sind_placid}.
As in \cite[Page 95]{Zhu25} we have
\begin{equation}\label{eq:Hecke_of_classifying_stack}
	\Hk^n(\bbB \underline{K}) = [\underline{K} \backslash \underline{H} \times^{\underline{K}} \underline{H} \times^{\underline{K}} \dots \times^{\underline{K}} \underline{H}/\underline{K}]_{\rm fpqc}^{\sch} \cong [\underline{H} \backslash (\underline{H}/\underline{K})^{n+1}]_{\rm fpqc}^{\sch} = \coprod_{\underline{h} \in H\backslash (H/K)^{n+1}} \bbB (\underline{H_{\underline h}}),
\end{equation}
where $H$ acts diagonally on $(H/K)^{n+1}$ via left multiplication, and where the stabilizer is $H_{\underline h} = H_{h_0} \cap  \dots \cap H_{h_n}$ for $\underline{h} = (h_0,\dots,h_{n}) \in (H/K)^{n+1}$.  
Note that $H_{\underline h}$ is good since it is a subgroup of $h_0K h^{-1}_0$, and $K$ is good. 
To make the notation compatible with \Cref{prop:!dual_and_star_agree_for_classifying_stacks}, we let $Y_n=\Hk^n(\bbB \underline{K})$, and we let $\alpha_n:Y_n^\dagger \to [\ast/\underline{H}]$, be the map induced by $\alpha_H$ and the map $Y_n^\dagger \to (\bbB \underline{H})^\dagger$.  
On the analytic side the \v{C}ech nerve $\Hk^\bullet([\ast/\underline{K}])$ of the v-cover $[\ast/\underline{K}] \to [\ast/\underline{H}]$ admits a completely analogous explicit description, and we have $Y_n^\diamondsuit = Y_n^\dagger = \Hk^n([\ast/\underline{K}])$, by \cite[Proposition 2.20]{GIZ25}. 

By the second part of \Cref{sind-dagger-for-classifying-stacks}, justified by \Cref{sind-dagger-classifying} and its proof, each triple $(Y_n,[\ast/\underline{H}],\alpha_n)$ is a sind-$\dagger$-correspondence and 
\[c^{*,\alpha_H}=\colim_{n\in \Delta^\op} c^{*,\alpha_n}.\]
Since $\Shv^!(\bbB \underline{H})^\vee=\colim_{\substack{n\in\Delta^\op\\ (\ul{f_*})^o}} \Shv^!(Y_n)^\vee$ and $\Shv^*(\bbB \underline{H})=\colim_{\substack{n\in\Delta^\op\\ f_!}} \Shv^*(Y_n)$ (see \Cref{prop:!dual_and_star_agree_for_classifying_stacks}), it suffices to factor each $c^{*,\alpha_n}$ through the natural map $\Shv^!(Y_n)^\vee\to \Shv^*(Y_n)$.
In other words, we want to show that the following diagram
\begin{center}
\begin{tikzcd}
	\Shv^!(Y_n)^\vee \arrow{r} \arrow{d}{c^{*,\alpha_n}}  & \Shv^*(Y_n) \arrow{d}{c^*} \\
	\calD^\an_\Lambda([\ast/\underline{H}])  & \calD^\an_\Lambda(\Hk^n([\ast/\underline{K}])) \arrow{l}{f_!}
\end{tikzcd}
\end{center}
is commutative.
To this end, we unfold ${c^{*,\alpha_n}}$ using \Cref{prop: analytification ind quasi-placid}. 
Let 
\[ \calI_n = \{\mu \colon \mu \subseteq H\backslash (H/K)^{n+1} \text{ finite subset} \}\]
For each $\mu \in \calI_n$ we have the closed substack  
\[ 
\Hk_n(\bbB \underline{K})_{\leq \mu} = \coprod_{\underline{h} \in \mu} \bbB \underline{H_{\underline h}} \subseteq \Hk_n(\bbB \underline{K})
\]
defined in terms of the right hand side of \eqref{eq:Hecke_of_classifying_stack}. 
This gives an ind-presentation (resp.\ a colimit formula) 
\[Y_n=\colim_{\mu\in \calI_n} \Hk_n(\bbB \underline{K})_{\leq \mu} \text{ (resp.\ } \Shv^!(Y_n)^\vee\simeq\colim_{\mu\in \calI_n} \Shv^!(\Hk_n(\bbB \underline{K})_{\leq \mu})^\vee\text{)}  .\]
By \Cref{prop: analytification ind quasi-placid}, to produce the factorization of $c^{*,\alpha_n}$ it suffices to show that for all $\mu\in \calI_n$ the map 
\[\Shv^!(\Hk_n(\bbB \underline{K})_{\leq \mu})^\vee\to \calD^\an_\Lambda([\ast/\underline{H}])\]
factors through $\Shv^*(Y_n)$.
But this is the case, since we have the commutative diagram 
\begin{center}
\begin{tikzcd}

	\Shv^!(\Hk_n(\bbB \underline{K})_{\leq \mu})^\vee \arrow{rr} \arrow{d}  & &  \Shv^*(Y_n)^\vee \arrow{d} \\
	\Shv^*(\Hk_n(\bbB \underline{K})_{\leq \mu})\arrow{rr} \arrow{d}  & &  \Shv^*(Y_n) \arrow{d} \\
	\calD^\an_\Lambda(\Hk_n([\ast/\underline{K}])_{\leq \mu}) \arrow{rr} \arrow{dr}{\alpha_{(n,\mu),!}} & &  \calD^\an_\Lambda(\Hk_n([\ast/\underline{K}]))\ar{dl}{\alpha_{n,!}} \\
								 & \calD^\an_\Lambda([\ast/\underline{H}]) & 
\end{tikzcd}
\end{center}
Now (2) follows from the defintion of $c^{\ast,\alpha_H}$ and the functoriality in \Cref{rem:compatibility_dagger_correspondences} and (3) follows from (2) and \Cref{prop:classifying_stack_is_sind_placid}(4).
Finally, (4) follows from (3) and \Cref{lm:analytification_isom_for_profinite_groups} as we range over deeper good subgroups $K\subseteq H$ which produce compact generators for $\Rep(H)$.
\end{proof}

\subsection{Loop groups and Grassmannians}
\label{sec: loop gropus etc}
For $R\in \CAlg^{\perf}$ let $W(R)$ be the ring of $p$-typical Witt vectors of $R$. 
We put $\bbW(R) = W(R) \otimes_{\bbZ_p} O_E$ if ${\rm char}(E) = 0$, resp. $\bbW(R) = R[\![\pi]\!]$ otherwise. 
We let $\varphi:\bbW(R)\to \bbW(R)$ denote the automorphism of $\bbW(R)$ induced from $\on{Frob}$.

Let $G$ be a reductive group over $E$ and let $\calG$ be a smooth affine model of $G$ over $O_E$. 
We define the \emph{jet group} (or \emph{positive loop group}) of $\calG$ resp. the \emph{truncated jet group} of $\calG$, resp. the \emph{loop group} of $G$ as prestacks with formula
	\begin{align*}
		L^+\mathcal{G} &\colon \Spec R \mapsto \mathcal{G}(\bbW(R)) \\
		L^+_n\mathcal{G} &\colon \Spec R \mapsto \mathcal{G}(\bbW(R)/\pi^n\bbW(R)) \\
		LG &\colon \Spec R \mapsto G(\bbW(R)[\frac{1}{\pi}])
	\end{align*}
	For $n \geq 1$, $L^+_n\calG$ is a pfp perfect affine group scheme over $k$.
For $1\leq m \leq n$ the transition maps $L^+_n\calG \rar L^+_m\calG$ are unipotent and $L^+\mathcal{G} = \lim_n L^+_n\calG$ is a perfect affine group scheme.
For $n \geq 1$ let
	\[L^+\calG^{(n)} = \ker(L^+\calG \rar L^+_n\calG).\]

	Assume that $\calG$ is a parahoric model of $G$.
The quotient $\Gr_{\calG}^{\rm sch} = LG/L^+\calG$ is the \emph{affine Grassmannian} of $\calG$ and the stack quotient $\Hk^\sch_\calG:=L^+\calG \backslash LG/L^+\calG$ is the \emph{local Hecke stack}. 
	If $\calI$ is an Iwahori model of $G$ with $\calI(O_E)$ contained in $\calG(O_E)$, and $\widetilde W$ denotes the Iwahori-Weyl group of $G$ with respect to $\calI$, then the $L^+\calI$-orbits in $\Gr_{\calI}^{\rm sch}$ are parametized by $\widetilde W$. 
	Similarly, the $L^+\calG$-orbits in $\Gr_{\calG}^{\rm sch}$ are parametrized by $W_\calG \backslash \widetilde W/W_\calG$. 
	In geometric terms, this gives rise to an identification 
	\[|\Hk^\sch_\calG|=W_\calG \backslash \widetilde W/W_\calG.\]
	There is an evident quotient map $\Gr_{\calG}^{\rm sch}\to \Hk^\sch_\calG$.
	For each $w \in W_\calG \backslash \widetilde W/W_\calG$, we have a stratum $\Hk^\sch_{\calG,w}$ and its pullback gives rise to an orbit $\Gr^{\rm sch}_{\calG,w}$ which is the perfection of a smooth affine scheme. 
	We let $\Hk^\sch_{\calG,\leq w}$ and $\Gr^{\rm sch}_{\calG,\leq w}$ denote the closures.
	We thus have the inclusions
	\[
		\Hk^\sch_{\calG,w} \subseteq \Hk^\sch_{\calG,\leq w} \subseteq \Hk^\sch_{\calG} \text{ and }
	\Gr^{\rm sch}_{\calG,w} \subseteq \Gr^{\rm sch}_{\calG,\leq w} \subseteq \Gr^{\rm sch}_{\calG}.
	\]
	Let $LG_w$, $LG_{\leq w}$ denote the preimages of $\Gr^{\rm sch}_{\calG,w}$, $\Gr^{\rm sch}_{\calG,\leq w}$ in $LG$.
Then one has the presentations $LG_{\leq w} = \lim_n \Gr^{{\rm sch},(n)}_{\calG,\leq w}$ and $LG_w = \lim_n \Gr^{{\rm sch},(n)}_{\calG,w}$, where
	\[ \Gr^{{\rm sch},(n)}_{\calG,\leq w} = LG_{\leq w}/L^+\calG^{(n)} \quad \text{ (resp.} \quad \Gr^{{\rm sch},(n)}_{\calG,w} = LG_{w}/L^+\calG^{(n)}\text{)}\] are again perfections of normal projective schemes (resp. affine schemes), and transition maps are unipotent for $n\geq 1$.
From these presentations we have the following.

	\begin{lemma}[{\cite[\S 3.1.3]{Zhu25}}]\label{lemma:LG_ind_placid} $LG_{\leq w}$, $LG_w$ are standard placid over $k$. 
		The loop group $LG$ is an ind-placid scheme.
		The Hecke stack $\Hk^\sch_\calG$ is ind-very-placid.
	\end{lemma}

\subsection{Isocrystals and shtukas}\label{sec:isoc_and_shtukas}

Let $\calG$ be a parahoric model of $G$.
We consider  
	\begin{align*}
		\ICG &= [LG/\mathrm{Ad}_\varphi LG]_{\rm \acute{e}t} \\
		\Sht_\calG^{\rm sch} &= [LG/\mathrm{Ad}_\varphi L^+\calG]_{\rm \acute{e}t},
	\end{align*}
	where by the quotient we mean the \'etale sheafification of the prestack quotients. 
	The first is called the \emph{stack of isocrystals} and the second is the \emph{stack of schematic shtukas}. 
	Note that for $R \in {\rm CAlg}^{\rm perf}$, $\ICG(R)$ is the groupoid of pairs $(\calE,\varphi)$ where $\calE$ is a $G$-torsor on $\bbW(R)$, trivial \'etale locally on $\Spec R$, and $\varphi$ is an isomorphism of $\calE$ with its Frobenius pullback.

	The natural map
	\begin{equation}
		\label{Leibniz-map}	
	{\rm Nt} = {\rm Nt}_{\calG} \colon \Sht_\calG^{\rm sch} \rar \ICG.
	\end{equation}
	is called the \emph{Newton map}.

\begin{lemma}[{\cite[Lemma 3.28]{Zhu25}}] \label{lm:Nt is ind-pfp proper}
	The map ${\rm Nt}$ is ind-pfp proper (as in \Cref{defn: indproper}) and all of its geometric fibers are isomorphic to $\Gr_\calG$.
    \end{lemma}
    \begin{remark}
	    We emphasize that \Cref{lm:Nt is ind-pfp proper} only claims ind-representability of ${\rm Nt}$ in $E^{\Alg,\rep}_\pfp$, and that it doesn't claim its ind-representability in $E^\rep_\pfp$.
    \end{remark}

	There is also an evident map $\Sht_\calG\to \Hk^\sch_\calG$ which allow us to pullback the stratification on $\Hk^\sch_\calG$. 
	Let $w \in W_\calG \backslash \widetilde W / W_\calG=|\Hk^\sch_\calG|$.
We have the substacks
	\begin{align*}
		\Sht^{\rm sch}_{\calG,w} = LG_w / \Ad_\varphi L^+\calG \\
		\Sht^{\rm sch}_{\calG,\leq w} = LG_{\leq w}/ \Ad_\varphi L^+\calG
	\end{align*}
	of $\Sht^{\rm sch}_\calG$.

Recall from \cite[\S3.1.4-3.1.5]{Zhu25}, that the \v{C}ech nerve of the Newton map gives rise to the simplicial stack $\{{\rm Hk}^n(\Sht_{\calG}^{\sch})\}_{n\geq 0}$, that $\ICG = \colim_{n \in \Delta^{\rm op}} {\rm Hk}^n(\Sht_{\calG}^{\sch})$, and that each ${\rm Hk}^n(\Sht_{\calG}^{\sch})$ is ind-very placid.
We summarize what we need in the following theorem.

	\begin{theorem}[{\cite[Example 3.13, Lemma 3.28]{Zhu25}}]
		\label{schematic-geometry-ofshtukas}
		The following statements hold.
		\begin{enumerate}
			\item The stacks $\Sht^{\rm sch}_{\calG,\leq w}$ and $\Sht^{\rm sch}_{\calG,w}$ are very placid with very placid atlas $LG_{\leq w}\to \Sht^{\rm sch}_{\calG,\leq w}$ and $LG_{w}\to \Sht^{\rm sch}_{\calG,w}$ respectively.
In other words, $\Sht^{\rm sch}_{\calG,\leq w}, \Sht^{\rm sch}_{\calG,w}\in \SchStk^\vpl$.
			\item The stack $\Sht^{\rm sch}_\calG$ is ind-very-placid with ind-presentation 
				\[\Sht^{\rm sch}_\calG=\colim_{w\in W_\calG \backslash \widetilde W / W_\calG} \Sht^{\rm sch}_{\calG,\leq w}.\]
				In other words, $\Sht^{\rm sch}_\calG\in \IndStk^\vpl$. 
			\item The stack $\ICG$ is sind-very-placid with sind-very-placid atlas
				\[ 
	{\rm Nt} = {\rm Nt}_{\calG} \colon \Sht_\calG^{\rm sch} \rar \ICG.
\]
In other words, $\ICG\in \sStk^\vpl$.
		\end{enumerate}
	\end{theorem}
	\subsection{Analytic stacks of local shtukas}
Let us recall the perfectoid v-stack of local shtukas. 
Given $S=\Spa(R,R^+)$, we can form a space $\calY_S:=\Spa \bbW(R^+)\setminus V([\varpi])$ as in \cite[\S II.1.1]{FS21}.
\begin{definition}
We let $\Sht^\an_\calG\in \AnStk_v$ denote the moduli stack with formula 
\[\Sht_\calG^\an(R,R^+):=\{(\calE,\Phi)\}\]
where $\calE$ is a $\calG$-bundle on $\calY_S$ and 
\[\Phi:\varphi^*\calE\dasharrow \calE\]
is an isomorphism defined over $\calY_S\setminus V(\pi)$ that is meromorphic along $V(\pi)\subseteq \calY_S$ (see \cite[Definition 5.3.5]{SW20}).
\end{definition}

Recall that $S=V(\pi)\subseteq \calY_S$ and that completion along this closed immersion gives rise to a map of locally ringed spaces  
\[(S,\bbW(\calO_S)) \to (\calY_S,\calO_{\calY_S}).\]
Restriction along this map gives rise to morphism
\[\Sht_\calG^\an\to (\Sht^\sch_\calG)^\diamondsuit\to  (\Hk^\sch_\calG)^\diamondsuit,\]
which on points (and before sheafification), can be described as
\[(\calE,\Phi) \mapsto (\calE_1=\varphi^*\calE_{|\bbW(R)},\calE_2=\calE_{|\bbW(R)},\Phi).\] 
We may pull back along the closed immersion $ (\Hk^\sch_{\calG,\leq w})^\diamondsuit\subseteq (\Hk^\sch_\calG)^\diamondsuit$ to obtain closed substacks $\Sht^\an_{\calG,\leq w}\subseteq \Sht_\calG^\an$.

It is natural to wonder what the actual relationship is between these stacks of analytic shtukas and the previously introduced scheme-theoretic ones.
In fact, there is an isomorphism
\begin{equation}\label{meromorphic comparison_shtuka}
   (\Sht^\sch_\calG)^\dagger \simeq\Sht^\an_\calG 
\end{equation}
preserving the corresponding closed substacks given as the union of those Bruhat strata that are bounded by $w$ (see \cite[Theorem 7.13]{GIZ25}).

	\subsection{Construction of the functor}
	\label{ss: constructin the functors}
	Recall that in \S \ref{sind-dagger-subsection} we constructed analytification functors for sind-$\dagger$-correspondences.     
	In this subsection, we explain how to fit $\ICG$ and $\Bun_G$ in a sind-$\dagger$-correspondence.      
	In other words, we need to provide a well-behaved map 
	\begin{equation}
	    \sigma:\ICG^\dagger\to \Bun_G.
	\end{equation}
    In rough terms, this map is obtained, before sheafification, by pulling back $\varphi$-$G$-torsors along the map of locally ringed spaces 
    \[\calY_{(0,\infty), S} \to \Spec \bbW(\calO_S^\circ)[1/\pi],\]
    where the left-hand side is the usual Frobenius cover of the adic Fargues-Fontaine curve attached to $S$.
	There is an alternative way of thinking about this map. 
	In \cite{GIZ25} two of us, together with Zillinger, introduced a stack $\Bun_G^\mer$ fitting in a correspondence of the form (see \cite[Theorem 1.1]{GIZ25})
	\begin{equation}
		\label{BungMer-triangle}
	\begin{tikzcd}
		\Bun_G^\mer \arrow{r}{\sigma} \arrow{d}{\gamma}  & \Bun_G  \\
	\ICG^\diamondsuit.  & 
	\end{tikzcd}
	\end{equation}
	Moreover, the following identifications were constructed in loc.~cit. (see \cite[Theorems 1.4, 7.13]{GIZ25})
\begin{equation}
	\label{meromorphic comparison_bun}
	\ICG^\dagger\simeq \Bun_G^\mer
\end{equation}
together with a Cartesian diagram (see \cite[Theorem 6.6]{GIZ25})
	\begin{equation}
		\label{important Cartesian}
	\begin{tikzcd}
		\Sht_\calG^\an	\simeq (\Sht_\calG^\sch)^\dagger \arrow{r}{\on{Nt}^\dagger} \arrow{d}{b_{\Sht_\calG^\sch}}  & \ICG^\dagger\simeq \Bun_G^\mer \ar{d}{\gamma=b_{\ICG}} \\
		(\Sht_\calG^\sch)^\diamondsuit \ar{r}{\on{Nt}^\diamondsuit}  & \ICG^\diamondsuit. 
	\end{tikzcd}
	\end{equation}

	\begin{remark}
	For the purpose of this article one can take the identities of \eqref{meromorphic comparison_shtuka}, resp.~\eqref{meromorphic comparison_bun} as the definition of $\Sht^\an_\calG$ and of $\Bun_G^\mer$, but we want to point out that there are more natural presentations of these objects. 
	\end{remark}

	We want to use \Cref{criterion-for-sind} (or more precisely \Cref{criterion-for-sind-very})
  to show that the sind-$\dagger$-diagram $(\ICG,\Bun_G,\sigma)$ is a sind-$\dagger$-correspondence.
  Since $\Sht^\sch_\calG\to \ICG$ is a sind-very-placid atlas (see \Cref{schematic-geometry-ofshtukas}), it suffices to show, by \Cref{criterion-for-sind-very}, that the bounded pieces $\Sht^\sch_{\calG,\leq w}$ of the ind-presentation $\Sht^\sch_\calG\simeq \varinjlim \Sht^\sch_{\calG,\leq w}$ are $\sigma$-compatible in the sense of \Cref{sind-dagger-correspondence}.(3).
That $\Shv^{!}(\Sht^\sch_{\calG,\leq w})$ is compactly generated follows from \Cref{quasi-placid have compact generated} (2), since $\Sht^\sch_{\calG,\leq w}$ is quasi-compact very placid.
We have to check that if $\sigma_{\leq w}$ denotes the following composition 
	\[\sigma_{\leq w}:(\Sht^\sch_{\calG,\leq w})^\dagger \to (\Sht^\sch_\calG)^\dagger \xrightarrow{\on{Nt}^\dagger} \ICG^\dagger \xrightarrow{\sigma} \Bun_G,\]
    then, the map $\sigma_{\leq w}$ is $\calD_{\Lambda}^{\an}$-!-able.
This was proved by one of us (G.).
Indeed, the following statement is a variant of the main theorem of \cite[Theorem 1.1]{Gle23}. 
	\begin{theorem}{\cite[Theorem 1.1]{Gle23}}
			{\label{analytic Newton map is fdcs}}
		The map 
		\[\sigma_{\leq w}:(\Sht^\sch_{\calG,\leq w})^\dagger \to \Bun_G\]
		is fdcss.
In particular, it is $\calD^\an_\Lambda$-$!$-able. 
	\end{theorem}
    To show that $\Sht^\sch_{\calG,\leq w}$ is $\sigma$-compatible, all that remains to show is the following.
	\begin{lemma}
		\label{resilience lemma ofr sht}
	For all $w\in \widetilde W$, $\Sht^\sch_{\calG,\leq w}$ is resilient.
	\end{lemma}
	\begin{proof}
	By \cite[Corollary 17.1.9]{SW20}, $(\Sht^\sch_{\calG,\leq w})^{\diamondsuit_\pre}$ is already a v-sheaf, and consequently the map 
		\[(\Sht^\sch_{\calG,\leq w})^{\diamondsuit_\uop}\xrightarrow{\simeq}  (\Sht^\sch_{\calG,\leq w})^{\diamondsuit}\]
		is already an equivalence.
	\end{proof}
    We can summarize the above discussion as follows.

\begin{corollary}{\label{cor: ConstructionofMainSindDaggerCorrespondence}} 
	$(\ICG,\Bun_G,\sigma)$ is a sind-$\dagger$-correspondence. 
\end{corollary}

We are finally able to define the comparison functor. 
Its definition involves the canonical self-duality $\id_\BZ \colon \Shv^!(\ICG) \to \Shv^!(\ICG)^\vee$ of $\Shv^!(\ICG)$ constructed by Zhu in \cite[Theorem 1.3]{Zhu25}.
We will review $\id_\BZ$ in detail later on (see \Cref{prop: BasicPropertiesofBZDualityonIsoc}).

	\begin{definition}
		\label{Defining-the-functor}
		We let $\pitch^\vee:=c^{*,\sigma}$ denote the comparison functor 
		\[\pitch^\vee:\Shv^!(\ICG)^\vee\to \calD^\an_\Lambda(\Bun_G)\]
		defined by the sind-$\dagger$-correspondence $(\ICG,\Bun_G,\sigma)$.
		We let 
		\[\pitch:\Shv^!(\ICG)\to \calD^\an_\Lambda(\Bun_G)\]
		denote the composition of Zhu's duality identification 
		\[\id_\BZ:\Shv^!(\ICG)\to \Shv^!(\ICG)^\vee\] 
		followed by $\pitch^\vee$.
	\end{definition}

	\begin{remark}
		In summer of 2023, one of us (G.) gave a talk in Bonn predicting the existence of the equivalence $\pitch$ satisfying some key compatibilities, and explaining the evidence for the existence of this functor that the work \cite{GIZ25} provided.   	
		Partially inspired by these key compatibilities discussed, Hansen constructed (independently of the existence of $\pitch$) an exceptional t-structure on $\calD^\an_\Lambda(\Bun_G)$ which would match the perverse t-structure defined by Zhu on $\Shv^!(\ICG)$.
		Hansen calls it the hadal\footnote{The etymology of this adjective stems from Hades (the underworld), and was first coined by marine biologists to describe the deepest oceanic trenches.} t-structure (see Remark \ref{rem: HadalTStructure} and Theorem \ref{thm: texactnessstatement}, for the precise realization of this vision in terms of our functor).
		Although the first author's original intention was to write $\Psi$ on the blackboard (due to a vague connection to nearby-cycles), Hansen commented privately that it ended up looking like a trident\footnote{We note however that Hades carries a bident, and it is Neptune who carries a trident.
Though, despite our best efforts, the symbol $\pitch$ really resembles more a pitchfork than a trident, so perhaps we cannot escape the satanic nature of its construction.} ($\pitch$).
		To not overload the usage of $\Psi$, and to allude to the hadal t-structure of Hansen, we have opted to use $\pitch$ to denote our functor.
	\end{remark}

	\section{Local Langlands categories}

	In order to show that the functor $\pitch$ constructed in \Cref{Defining-the-functor} is an equivalence, one first studies the the local Langlands categories $\Shv^{!}(\ICG)$ and $\calD_{\Lambda}^{\an}(\Bun_{G})$ individually.
    In this section, we recall the theory of sheaves on the stack $\ICG$ of $G$-isocrystals and on the stack $\Bun_G$ of $G$-bundles on the Fargues--Fontaine curve, as well as their formal properties.
Both these objects are geometrizations of a certain explicit combinatorial object known as the Kottwitz set of $G$.
Before discussing the Kottwitz set, we recall our group theoretic notation.

    Recall that, for us, $G$ is a quasi-split reductive group over $E$ and that we have fixed a rationally defined maximal torus and Borel subgroups $T \subset B\subset G$ as in the introduction.
    We let $A\subseteq T$ be the maximal split subtorus. 
    Recall that we have fixed an Iwahori model $\calI$ of $G$ defined over $O_E$.
    Moreover, we will assume that $\calT(O_{\breve{E}})$ (the Ner\'on model) is contained in $\calI({O_{\breve{E}}})$.
    In terms of the building, $\calI$ corresponds to a $\varphi$-invariant alcove in the apartment defined by $A$ (or $T$).

    \subsection{The Kottwitz set}
    In this subsection, we discuss the combinatorial properties of the set $B(G):=G(\breve{E})/\Ad_\varphi(G(\breve{E}))$ of $\varphi$-conjugacy classes in $G(\breve{E})$, originally introduced by Kottwitz \cite{Kot85}.
The first main invariant of elements in the set $B(G)$ is its Newton point $\nu_G(b) \in X_\ast(T_{\ol{E}})_{\bbQ}^{+,\Gamma_E} $, see \cite[\S 4]{Kot85}, where $\Gamma_E$ denotes the absolute Galois group of $E$ acting on rationalized dominant geometric coweights $X_\ast(T_{\ol{E}})_\bbQ^+$. %

The second crucial combinatorial invariant of $B(G)$ is the Kottwitz point $\kappa_G(b) \in \pi_1(G)_{\Gamma_E}$.
Indeed, following Kottwitz, one defines a group homomorphism $G(\breve{E}) \to \pi_1(G)_{I_E}$ towards inertia coinvariants\footnote{Note that $\pi_1(G)_{I_E}$ can be shown to identify with the geometric connected components of $LG$.}, and one can pass to a $\varphi$-conjugacy classes to obtain the map (see \cite[\S7]{KottwitzII} and \cite[\P 2.a.2]{pappas_rapoport_twisted_loop_groups}).
Alternatively, one has an isomorphism $\pi_{1}(G)_{\Gamma_E} \simeq X^{*}(Z(\hat{G})^{{\Gamma_E}})$, where $Z(\hat{G})$ is the center of the Langlands dual group of $G$, and one can view $\kappa_{G}$ as a map $B(G) \ra X^{*}(Z(\hat{G})^{{\Gamma_E}})$, as in \cite[\S~5]{Kot85}.

Following \cite{RR_96}, one can use $\nu_G$ and $\kappa_G$ to define a partial order $\leq$ on $B(G)$ as follows: we say that $b_1\leq b_2$ if and only if $\nu_G(b_1)\leq \nu_G(b_2)$ for the usual rational Bruhat order (i.e., the difference $\nu_G(b_2)-\nu_G(b_1)$ is a non-negative rational linear combination of positive simple coroots of $G$) and $\kappa_G(b_1)=\kappa_G(b_2)$. 
This is clearly transitive and its anti-symmetry is a consequence of the injectivity of $(\nu_G,\kappa_G)$, see \cite[\P 4.13]{KottwitzII}. 
A key finiteness property of this order is the following:

\begin{lemma}\label{lem_finite_B(G)_order}
	The poset $(B(G),\leq)$ is down-finite, i.e., for all $b\in B(G)$ the set $\{b'\in B(G)\mid b'\leq b\}$ is finite.
\end{lemma}

\begin{proof}
	This is \cite[Proposition 2.4.(iii)]{RR_96}.
\end{proof}

The partial order on $B(G)$ allow us to define a natural topology on the set $B(G)$. 
At the same time, there is also the opposite partial order to which we will also attach a topological space. 
For later reference, we record our conventions concerning $B(G)$ regarded as a topological space.

\begin{definition}{\label{defn: topologiesonB(G)}}
We endow $B(G)$ with the topology in which open subsets are up-closed for the partial order $\leq$, i.e., they are arbitrary unions of the sets of the form $U_{b} := \{x \in B(G)\colon b \leq x\}$ where $b \in B(G)$. 
Moreover, we denote by $B(G)^{\op}$ the topological space whose underlying set is the Kottwitz set of $G$, but whose open subsets are down-closed, i.e., they are arbitrary unions of sets that have the form $\{x \in B(G)\colon b \geq x \}$ for $b \in B(G)$.
We recall that the minimal elements in this partial ordering are denoted by $B(G)_{\mathrm{basic}} \subset B(G)$ and are referred to as the basic elements. 
\end{definition}

We have the following reductive groups attached to elements $b \in B(G)$.

\begin{definition}[\cite{KottwitzII}, \S3.3 and Appendix A; \cite{RZ96}, 1.12]{\label{defn: sigmacentralizer}}
	Let $b \in B(G)$, and fix a representative $\dot{b}\in G(\breve{E})$.
	We define the smooth affine connected $E$-group $G_{\dot{b}}$ as representing the functor:
\begin{equation}
	G_{\dot{b}}(R):=\{ g\in G(\breve{E}\otimes_E R) \colon g^{-1}\dot{b}\varphi(g)=\dot{b}\}
\end{equation}
for any $E$-algebra $R$. 
\end{definition}

Note that if $\dot{b}_1=\gamma^{-1} \dot{b}_2 \varphi(\gamma)$, then $G_{\dot{b}_1}\simeq G_{\dot{b}_2}$ via $\gamma$-conjugation. 
For this reason, we will at times abuse notation and denote $G_b$ a representative of this isomorphism class of reductive groups. 
Also recall, from \cite[Cororllary 1.14]{RZ96} (or \cite[Proposition III.4.2]{FS21}), that $G_b$ is always an inner form of a Levi subgroup of the quasi-split inner form of $G$.
The locally pro-finite groups $G_b(E)$ naturally arise as the automorphism of the $G$-isocrystals associated with $b$. \\

We now recall an alternative way of understanding the set $B(G)$ coming from the combinatorics of affine Weyl groups through the so-called $\sigma$-straight $\sigma$-conjugacy classes in the Iwahori-Weyl group $\widetilde W = N_G(T)(\breve E)/\calT(O_{\breve E})$ of $G$ (see \cite{He2014GeometricHomological}).
Recall that $\widetilde{W}$ fits into a short exact sequence
    \begin{equation}{\label{eqn: RelativeWeylGroup}}
    1 \ra X_{*}(T_{\ol{E}})_{I_{F}} \ra \widetilde{W} \ra W_{0} \ra 1, 
    \end{equation}
    where $W_{0} := N_{G}(T)(\Breve{E})/T(\Breve{E})$ denotes the relative Weyl group of $G_{\Breve{E}}$ (see \cite[\S~3.1.1]{Zhu25}). 
	Recall that $\widetilde W$ is a quasi-Coxeter group, equipped with a length function $\ell \colon \widetilde W \rar \bbZ_{\geq 0}$, see \cite{HaiRap} and (for example) \cite[\S4.4]{HeNie_minlen}.
	Since we have fixed an alcove in the reduced Bruhat--Tits building of $G$ over $E$ (the one corresponding to $\calI$), the stabilizer in $\widetilde W$ of that alcove is naturally identified with the set of length zero elements. 
    Let $\sigma$ denote the automorphism of $\widetilde W$ induced by $\varphi$.
    More precisely, 
An element $w\in \widetilde W$ is called \emph{$\sigma$-straight} if, for all $n \geq 0$, one has $\ell(w \sigma(w) \dots \sigma^{n-1}(w)) = n\ell(w)$, see \cite[\S2.4]{He2014GeometricHomological}. 
	A $\sigma$-conjugacy class in $\widetilde W$ is called \emph{$\sigma$-straight} if it contains a $\sigma$-straight element.
Let $B(\widetilde W)$ denote the set of $\sigma$-conjugacy classes in $\widetilde W$ and let $B(\widetilde W)_{\rm str}$ be the subset of straight classes.
The inclusion $N_G(T)(\breve E) \subseteq G(\breve E)$ induces a map $B(\widetilde W) \rar B(G)$ and the composition
	\[
	B(\widetilde W)_{\rm str} \hookrightarrow B(\widetilde W) \rar B(G)
	\]
	is a bijection, by \cite[Theorem~3.7]{He2014GeometricHomological}. 
	This gives rise to a combinatorial parametrization  
    \begin{equation}{\label{eqn: Relationshipbetweenstraightelements}}
     B(\widetilde{W})_{\rm str} \xrightarrow{\simeq} B(G), 
    \end{equation}
    which allows us to define the following. 
    \begin{definition}{\label{defn: sigmastraightattatchedtob}}
	    For each $b \in B(G)$ we define the following.
    \begin{enumerate}
    \item We let $w_{b} \in B(\widetilde{W})_{\mathrm{str}}$ denote the preimage of $b \in B(G)$ under the isomorphism (\ref{eqn: Relationshipbetweenstraightelements}). 
    \item For any fixed lift $\dot w_b \in N_G(T)(\breve E)$ of $w_b$, we consider
	    \[I_{\dot w_b} := \{g \in L^{+}\calI(\overline{\bbF}_p)=\calI(\calO_{\breve{E}}) \colon g^{-1}\dot w_b \varphi(g) = \dot w_b\}. \]
	    As it turns out, this group is an Iwahori subgroup of $G_{\dot w_b}(E)\simeq G_{\dot{b}}(E)$, where $G_b$ is the algebraic group considered in \Cref{defn: sigmacentralizer} (see \cite[Remark 3.19]{Zhu25}). 
	We will occasionally abuse notation and denote this simply by $I_{b}$.
    \end{enumerate} 
    \end{definition}

	\subsection{The analytic LLC category}{\label{ss: theanalyticcategory}}
	In this section, we review the geometric structure of the perfectoid $v$-stack $\Bun_{G}$, and the structure of the category $\calD_{\Lambda}^{\an}(\Bun_{G})$. 
    
	Let $\Bun_{G}$ be the $v$-stack of $G$-bundles on $X_{T}$, the relative Fargues-Fontaine curve over $T$ (see \cite[Definition III.0.1]{FS21}). 
	We have the following basic facts on its geometric structure.
	\begin{theorem}{\label{thm: BunGGeometricFacts}}
		The following is true. 
		\begin{enumerate}
			\item The $v$-stack $\Bun_{G}$ is an $\ell$-cohomologically smooth Artin $v$-stack of pure $\ell$-dimension $0$.
			\item There is a natural homeomorphism 
			\[ |\Bun_{G}| \simeq B(G)^\op, \]
			where $|\Bun_{G}|$ denotes the underlying topological space of $\Bun_{G}$ and the right-hand side is the topology described in Definition \ref{defn: topologiesonB(G)}.
			\item For all $b \in B(G)$, the locally closed substacks $j_{b}: \Bun_{G}^{b} \hookrightarrow \Bun_{G}$ (referred to as the Harder-Narasimhan (abbv. HN) stratification) furnished by point (2) are isomorphic to $[\ast/\tilde{G}_{b}]$. 
				Here $\tilde{G}_{b}$ is a group v-sheaf admitting a split surjection
			\[ \tilde{G}_{b} \ra \ul{G_{b}(E)}. \]
			Moreover, the kernel of this surjection is a spatial kimberlite which is an iterated fibration of absolute positive Banach-Colmez spaces.  
			Here $G_{b}$ denotes the $\varphi$-centralizer of the element $b \in B(G)$, as in \Cref{defn: sigmacentralizer}.
			\item The map $p_{b}: [\ast/\tilde{G}_{b}] \ra [\ast/\ul{G_{b}(E)}]$ admits a section $s_{b}: [\ast/\ul{G_{b}(E)}] \ra [\ast/\tilde{G}_{b}]$ whose fibers are, again, iterated fibrations of positive Banach-Colmez spaces. 
		\end{enumerate}
	\end{theorem}
	\begin{proof}
		The first point is precisely \cite[Theorem~IV.1.19]{FS21}. 
		The second point is the main result of \cite{Viehamann_Newton_BdR}, the third and fourth points follows from \cite[Proposition III.5.1, Proposition~III.5.3]{FS21} (see also \cite[Theorem 5.4]{Gle23}).
	\end{proof}
         For any subset $S \subset B(G)$, we write $j_{S}: \Bun_{G}^{S} \hookrightarrow \Bun_{G}$ for the substack parametrizing $G$-bundles whose pullback to each geometric point has isomorphism class lying in the subset $S \subset B(G)$ via the isomorphism of sets in Theorem \ref{thm: BunGGeometricFacts} (2).
	Given a pair of subsets $S_1\subseteq S_2 \subseteq B(G)$, we get a map of stacks that we denote by
	\[j_{S_1,S_2}:\Bun_{G}^{S_{1}} \to \Bun_{G}^{S_{2}} .\]
	If $S_2=B(G)$ we omit it from the notation.
Moreover, for $b \in B(G)$, we use the notation 
\[
j_b \colon \Bun_{G}^{b} \to \Bun_{G} \quad j_{\leq b} \colon \Bun_{G}^{\leq b} \to \Bun_{G}, \quad j_{< b} \colon \Bun_{G}^{< b} \to \Bun_{G}
\]
to abbreviate $S_1=\{b\}$ resp. $S_1=\{b'\mid b'\leq b\}$, resp. $S_1=\{b'\mid b'< b\}$. 
	 The category of sheaves on these strata is built from the following categories of smooth representations.
	\begin{proposition}{\label{prop: semiorthogonaldecomp}}
		For all $b \in B(G)$, we have equivalences of categories 
		\[\calD_{\Lambda}^{\an}(\Bun_{G}^{b}) \simeq \calD_{\Lambda}^{\an}([\ast/\ul{G_{b}(E)}]) \simeq \Rep(G_{b}(E)).\]
		The first equivalence is induced by the pullback $p_{b}^{*}$ along the natural map $p_b: \Bun_G^b \simeq [\ast/\tilde{G}_{b}] \ra [\ast/\underline{G_{b}(E)}]$, and the second equivalence is the one appearing in \Cref{analitifying cInd}.(4) (see \cite[Theorem V.1.1]{FS21}). 
	\end{proposition}
	\begin{proof}
		This is precisely \cite[Proposition~VII.7.1]{FS21} and \cite[Theorem V.1.1]{FS21}. 
	\end{proof}	   
	\begin{remark}
		\label{absurdities}
    The notation $\tilde{G}_{b}$ and the identification $\Bun_{G}^{b} \simeq [\ast/\tilde{G}_{b}]$ described above is a slight abuse, as both $\tilde{G}_{b}$ and the identification $\Bun_{G}^{b} \simeq [\ast/\tilde{G}_{b}]$ really depend on a choice of representative $\dot{b}$ of the $\varphi$-conjugacy class $b \in B(G)$, as mentioned in Definition \ref{defn: sigmacentralizer}. 
    
    More precisely, if we fix a representative $\dot{b}\in G(\breve{E})$ of $b \in B(G)$, then this determines for us an isocrystal with $G$-structure and a map $\Spd \bar{\bbF}_p\to \Bun_G$ whose isomorphism class is uniquely determined by $b \in B(G)$. 
    While we will want to suppress the choice of representative from our notation in what follows, there will be indeed certain ambiguities that will be introduced if one completely ignores that a choice is implicitly involved. 
    Therefore, for certain arguments, it works better to fix once and for all a choice of representative $\dot{b}\in G(\breve{E})$ of $b\in B(G)$.
    For most abstract purposes, any representative works equally well, but for concrete computations it becomes important to choose $\dot{b}$ wisely.
    In \S \ref{ss: renormalizing the choice of b}, we discuss some conventions that makes the choice of $\dot{b}$ convenient for computational purposes.

    Once $\dot{b}$ is fixed, the identification 
    \[\calD^\an_\Lambda([\ast/\ul{G_{\dot{b}}(E)}])\simeq^{\eta_{\on{FS}}} \Rep (G_{\dot{b}}(E))\]
	    is canonical.
	    For this reason, we will often omit it from the notation, although implicitly it plays a role in our considerations.
	\end{remark}
    We now turn our attention to describing $\calD_{\Lambda}^{\an}(\Bun_{G})^{\omega}$ (i.e., the subcategory of compact objects inside $\calD_{\Lambda}^{\an}(\Bun_{G})$). 
			For any $b\in B(G)$ and a pro-$p$ compact open subgroup $K\subseteq G_b(E)$, we consider
			\begin{equation}
				\label{rb-of-bung1}
			\calR^{\Bun_G}_{b,K}:=j_{b!}p_{b}^{*}\on{c-Ind}_{K}^{G_b(E)}(\Lambda).
			\end{equation}

		Recall the following statement.

	\begin{theorem}[{\cite[Proposition~VII.7.4]{FS21}}]
{\label{thm: compactenerationofBunG}}
		The category $\calD_{\Lambda}^{\an}(\Bun_{G})$ is compactly generated.
An object $A \in \calD_{\Lambda}^{\an}(\Bun_{G})$ lies in $\calD_{\Lambda}^{\an}(\Bun_{G})^{\omega}$ if and only if  the following is true. 
		\begin{enumerate} 
			\item The object $A$ is supported on only finitely many Harder-Narasimhan strata (i.e., $j_b^*A\neq 0$ for finitely many $b\in B(G)$). %
			\item For all $b \in B(G)$, the restriction $j_{b}^{*}(A) \in \calD_{\Lambda}^{\an}(\Bun_{G}^{b})^\omega$.
		\end{enumerate}

	\end{theorem}
    	\begin{remark}
			\label{the-bung-b-case}
			Recall that if $H$ is a locally profinite group containing a good compact open subgroup $K\subseteq H$, then $\Rep (H)$ is compactly generated by \ref{prop:classifying_stack_is_sind_placid} (3).
			Moreover, recall that $\{\on{c-Ind}_{K}^{H}(\Lambda)\}_{K\subseteq H}$ is a family of compact generators when we range $K\subseteq H$ along the compact open good subgroups of $H$ (see \Cref{prop:classifying_stack_is_sind_placid} and \Cref{lm:reps_of_compact_groups}). 
			For the $p$-adic reductive group $H$ over $E$ (for example $H=G_b$), any pro-$p$ subgroup $K\subseteq H(E)$ is good.
			It follows from \Cref{prop: semiorthogonaldecomp} that $\{j_b^*\RKB^{\Bun_G}\}_{K\subseteq G_b(E)}$, as we range over $K\subseteq G_b(E)$ open compact pro-$p$ subgroups, is a family of compact generators for $\calD^\an_\Lambda(\Bun_G^b)$.
		\end{remark}

	Combining \Cref{the-bung-b-case} with \Cref{thm: compactenerationofBunG}, one obtains the following statement thanks to excision.

	\begin{proposition}
		\label{right-compact-generators-BunG}
		The family $\{\RKB^{\Bun_G}\}^{b\in B(G)}_{K\subseteq G_b(E)}$ ranging over $b\in B(G)$ and $K\subseteq G_b(E)$ open compact pro-$p$ subgroups, is a family of compact generators for $\calD^\an_\Lambda(\Bun_G)$.
				Similarly, the family $\{\RKB^{\Bun_G}\}^{b\in V}_{ K\subseteq G_b(E)}$ is a set of compact generators for $\calD_{\Lambda}^{\an}(\Bun_{G}^{V})$ and any convex subset $V \subset B(G)^{\op}$. 
	\end{proposition}

    We now describe how these compact generators behave under duality. 
    To do this, we recall that, for each $b \in B(G)$, there is an $\ell$-cohomologically smooth chart 
	\[ \sigma_{b}: \mathcal{M}_{b} \ra \Bun_{G}. \]
	It comes with a canonical map $\gamma_b:\calM_b\to \ICG^\diamondsuit_b$.
	Moreover, after fixing a choice of representative as in Convention \ref{conv: DualCpx}, we may rewrite it as $\gamma_{b}: \mathcal{M}_{b} \ra [\ast/\ul{G_{b}(E)}]$ (see \cite[Proposition~V.3.5, Theorem~V.3.7]{FS21}).\footnote{Note that in the reference provided $\sigma_b$ is denoted by $\pi_b$ and $\gamma_b$ by $q_b$.} 
	Recall from \cite[Proposition V.3.5, Proposition V.2.1]{FS21} that the map $\gamma_b$ is unipotent, it follows that 
	\[\gamma_b^!\Lambda\simeq \gamma_b^*\omega_{\calM_b}\]
	for a unique line bundle $\omega_{\calM_b}\in \calD^\an_\Lambda([\ast/\ul{G_b(E)}])$ with formula
	\[\omega_{\calM_b}\simeq \gamma_{b*}\gamma_b^!\Lambda.\]

	After fixing a representative $\dot{b}$ of $b$ (see \Cref{absurdities}), we may abuse the notation to obtain functors
\[s^*_b:\calD^\an_\Lambda(\Bun_G^b)\to \Rep(G_{b}(E))\] 
and 
\[\gamma_b^*:\Rep(G_{b}(E))\to \calD^\an_\Lambda(\calM_b).\]
	We have the following fact.

	\begin{proposition}[{\cite[Proposition~VII.7.2]{FS21}}]{\label{prop: jbnatural}}
		The functor 
		\[s^*_b j_{b}^{*}: \calD_{\Lambda}^{\an}(\Bun_{G}) \ra \Rep(G_{b}(E)) \]
		admits a left adjoint functor with formula 
		\[  \sigma_{b\natural}\gamma_{b}^{*}: \Rep(G_{b}(E)) \ra  \calD_{\Lambda}^{\an}(\Bun_{G}), \]
		where $\sigma_{b\natural}$ is the left adjoint to $\sigma_{b}^{*}$. 
		We can rewrite this as 
		\[\sigma_{b\natural}\gamma_{b}^{*}\simeq \sigma_{b!}\gamma_{b}^{*}(- \otimes \omega_{\calM_b}).\]
    Moreover, after fixing $\dot{b}$ as in \Cref{conv: DualCpx}, we may further write it as 
		\[\sigma_{b\natural}\gamma_{b}^{*}\simeq \sigma_{b!}\gamma_{b}^{*}(- \otimes \delta_{b})[2\langle 2\rho_{G}, \nu_{b} \rangle] \]
		where $\delta_b$ is the character $G_b(E)\to \Lambda^\times$ of \cite[Definition 3.14]{HIDualCpx} defined with respect to the datum fixed in \Cref{conv: DualCpx}.
	\end{proposition}
	\begin{proof}
		The first statement is \cite[Proposition~VII.7.2]{FS21}. 
		For the second statement, recall the formula for $\ell$-cohomologically smooth maps
		\[\sigma_\natural(-)\simeq \sigma_!(-\otimes \sigma^!\Lambda).\]
		We need to show that we have an identification
		\[ \sigma_{b}^{!}(\Lambda) \simeq \gamma_b^*(\omega_{\calM_b})\simeq \gamma_{b}^{*}(\delta_{b})[2\langle 2\rho_{G}, \nu_{b}\rangle]. \]
		However, by \cite[Proposition~1.1]{HIDualCpx} the dualizing complex on $\Bun_{G}$ is the constant sheaf, from this it follows that 
		\[\sigma_b^!\Lambda\simeq \gamma^!_b\Lambda\simeq \gamma_b^*\omega_{\calM_b}.\]

		That we can write the dualizing complex on $\mathcal{M}_{b}$ explicitly as $\gamma_{b}^{*}(\delta_{b})[2\langle 2\rho_{G}, \nu_{b}\rangle]$ follows from \cite[Corollary~1.6]{HIDualCpx}.
	\end{proof}
        \begin{remark}{\label{rem: jbsharprestricted}}
		Using \Cref{prop: jbnatural}, we can define a functor $j_{b \sharp}:\calD^\an_\Lambda(\Bun_G^b)\to \calD^\an_\Lambda(\Bun_G)$ that is left adjoint to $j_b^*$ and has formula 
		\[j_{b \sharp}=\sigma_{b\natural}\gamma_{b}^{*}s_b^*\]
		or equivalently, it satisfies
		\[j_{b \sharp}p_{b}^{*} \simeq \sigma_{b\natural}\gamma_{b}^{*}. \]
		From the above description of $j_{b\sharp}$, it follows that this functor naturally factors through the fully faithful functor 
        \[ j_{\leq b!}: \calD_{\Lambda}^{\an}(\Bun_{G}^{\leq b}) \hookrightarrow \calD_{\Lambda}^{\an}(\Bun_{G}). \]
	Indeed, we note that the image of the map $\sigma_{b}$ lies inside $\Bun_{G}^{\leq b}$, as easily follows from simple Harder--Narasimhan polygon considerations (see \cite[Theorem V.3.7]{FS21} and \cite[Lemma~8.3]{HamGeomES}).
        \end{remark}

	Similarly to \eqref{rb-of-bung1}, we define

			\begin{equation}
				\label{rb-of-bung}
				\calL^{\Bun_G}_{b,K}:=j_{b\sharp}p_{b}^{*}\on{c-Ind}_{K}^{G_b(E)}(\Lambda) \simeq \sigma_{b\natural}\gamma_b^*\on{c-Ind}_{K}^{G_b(E)}(\Lambda).
			\end{equation}

	The existence of the left adjoint $j_{b\sharp}$ naturally leads to a dual description of the compact generators described in \Cref{right-compact-generators-BunG}.
	We have the following statement.
	\begin{proposition}{\label{prop: compactgeneratorssecondpass}}
The family $\{\LKB^{\Bun_G}\}^{b\in B(G)}_{K\subseteq G_b(E)}$, as we range over $b\in B(G)$ and $K\subseteq G_b(E)$ open compact pro-$p$ subgroups, is a family of compact generators for $\calD^\an_\Lambda(\Bun_G)$.
	\end{proposition}
	\begin{proof}
		This is simply a reformulation of the second part of \cite[Proposition~VII.7.4]{FS21}.
	\end{proof}

    Let us recall Fargues--Scholze's Bernstein--Zelevinsky duality functor (\cite[\S~V.5]{FS21})
    \[\bbD_{\BZ}:\calD^\an_\Lambda(\Bun_G)\xrightarrow{\simeq} \calD^\an_\Lambda(\Bun_G)^\vee.\] 
	Recall that, for $A \in \calD_{\Lambda}^{\an}(\Bun_{G})^{\omega}$, by \cite[Proposition~VII.7.6]{FS21}, there exists a unique compact object $\bbD_{\BZ}(A) \in \calD_{\Lambda}^{\an}(\Bun_{G})^{\omega}$ characterized by the property that 
	\[ \ExtRHom(\bbD_{\BZ}(A),B) \simeq \pi_{\natural}(A \otimes B)\simeq \pi_{!}(A \otimes B), \]
	for all $B \in \calD^\an_\Lambda(\Bun_{G})$, where $\pi: \Bun_{G} \ra \ast$ is the natural projection and \(\pi_{\natural}\simeq \pi_!\) is the left adjoint of \(\pi^*\) (equivalently $\pi^!$). 
	\begin{theorem}{\label{thm: naturaland!exchangedunderBZ}}
		The mapping $A \mapsto \bbD_{\BZ}(A)$ described above restricts to an equivalence 
		\[ \bbD_{\BZ}: \calD_{\Lambda}^{\an}(\Bun_{G})^{\omega} \ra \calD_{\Lambda}^{\an}(\Bun_{G})^{\omega,\op}, \]
		which satisfies the following.
		\begin{enumerate}
			\item The equivalence $\bbD_{\BZ}$ is involutive (i.e., $\bbD^{\op}_{\BZ}\circ \bbD_{\BZ}\simeq \on{id}$).    
             \item If $A$ is supported on a quasicompact open subset $U \subset \Bun_{G}$ then so is $\bb{D}_{\BZ}A$ (cf. \Cref{rem: jbsharprestricted}). 
		     In particular, there is a unique equivalence
		\[ \bbD_{\BZ,U}: \calD_{\Lambda}^{\an}(\Bun^U_{G})^{\omega} \ra \calD_{\Lambda}^{\an}(\Bun^U_{G})^{\omega,\op}, \]
		fitting in a commutative diagram
		\begin{center}
		\begin{tikzcd}
			\calD_{\Lambda}^{\an}(\Bun^U_{G})^{\omega} \arrow{r}{\bbD_{\BZ,U}} \arrow{d}{j_{U!}}  & \calD_{\Lambda}^{\an}(\Bun^U_{G})^{\omega,\op} \arrow{d}{j^{\op}_{U!}} \\
			\calD_{\Lambda}^{\an}(\Bun_{G})^{\omega}\arrow{r}{\bbD_\BZ} & \calD_{\Lambda}^{\an}(\Bun_{G})^{\omega,\op}. 
		\end{tikzcd}
		\end{center}
     \end{enumerate}
     \begin{proof}
     	This is \cite[Proposition~VII.7.6]{FS21}. 
     \end{proof}
     \end{theorem}
     In \S \ref{subsec: ApplicationsofSemiorthogonal}, we will interpret \Cref{thm: naturaland!exchangedunderBZ}.(2) as saying that $\bbD_{\BZ}$ promotes to an equivalence in the category of semi-orthogonal decompositions as in \Cref{def: semi-orthogonal} and \Cref{rem: alternative defininition of semi-orthogonal}. 
     Formally from the theory of semi-orthogonal decompositions, we obtain for all $b\in B(G)$ a natural equivalence 
        \[ \bb{D}_{\BZ,b}: \calD_{\Lambda}^{\an}(\Bun_{G}^{b}) \xrightarrow{\simeq} \calD_{\Lambda}^{\an}(\Bun_{G}^{b})^{\vee},\]
        sitting in a natural commutative diagram 
        \begin{equation}{\label{eqn: semiorthogonaldecomposition}}
        \begin{tikzcd}
        \calD_{\Lambda}^{\an}(\Bun_{G}^{b})^\omega \arrow[d,"j_{b!}"] \arrow[r,"\bb{D}_{\BZ,b}"] & \calD_{\Lambda}^{\an}(\Bun_{G}^{b})^{\omega,\op} \arrow[d,"j_{b\sharp}^{\op}"]  \\
        \calD_{\Lambda}^{\an}(\Bun_{G})^\omega \arrow[r,"\bb{D}_{\BZ}"]   & \calD_{\Lambda}^{\an}(\Bun_{G})^{\omega,\op}   
        \end{tikzcd}
        \end{equation}
        (see \Cref{cor: SelfDualizableBunG} for details).

	The functor $\bb{D}_{\BZ,b}$ is intimately related to the cohomological duality functor
	\[ \bbD_{\coh,G_{b}(E)}: \Rep(G_{b}(E)) \ra \Rep(G_{b}(E))^{\vee} \]
	\[ A \mapsto \ExtRHom_{G_{b}(E)}(A,\mathcal{H}(G_{b})) \]
    described in \Cref{ss: BZDuality}. 
    Indeed, the duality $\bbD_{\BZ}$ on $\Bun_{G}$ comes from the Frobenius algebra structure (\cite[Definition~4.6.5.1]{LurieHigherAlgebra}) given by 
    \[\pi_!=\pi_{\natural}: \calD_{\Lambda}^{\an}(\Bun_{G}) \ra {\rm Mod}_\Lambda,\] 
    together with the standard symmetric monoidal structure on $\calD_{\Lambda}^{\an}(\Bun_{G})$ (i.e., $\otimes$).
    Similarly, $\bbD_{\BZ,b}$ is induced from $\otimes$ and 
    \[\pi_{b!}:\calD_{\Lambda}^{\an}(\Bun^b_{G}) \ra {\rm Mod}_\Lambda.\] 
    Where $\pi_b=\pi\circ j_b$ is the structure map $\Bun_G^b\to \ast$.
    
    For basic strata, and under the standard identification 
    \[\Rep(G_b(E)) \simeq \calD_{\Lambda}^{\an}([\ast/\ul{G_b(E)}]),\]
    the Frobenius algebra induced by $(\pi_{b!},\otimes)$ naturally identifies with the one discussed in \S \ref{ss: BZDuality}.
    In contrast, when $b$ is not basic, when we compare $\bbD_{\BZ,b}$ and $\bbD_{\coh,G_{b}(E)}$ under a choice of identification $\calD_{\Lambda}^{\an}(\Bun_{G}^{b}) \simeq \Rep(G_{b}(E))$, these functor differ by a line bundle twist. 
    We explain this line bundle now.

    \begin{proposition}{\label{prop: dualcpxonClassifyingStack}}
	    Let $\pi_b:\Bun_G^b\to \ast$ denote the structure morphism. 
    We have an identification
    \[p_{b*} \pi_b^!\Lambda\simeq s_b^*\pi_b^!\Lambda \simeq \omega_{\calM_b}^{-1}\]
    or equivalently 
    \[\pi_b^!\Lambda\simeq p_b^*(\omega^{-1}_{\calM_b}).\]
    Moreover, after fixing $\dot{b}$ as in \Cref{conv: DualCpx}, we may further write it as 
    \[\pi_b^!\Lambda\simeq p_{b}^{*}(\delta_{b}^{-1})[-2d_{b}],\]
    where $d_{b} := \langle 2\rho_{G},\nu_{b} \rangle$, and $\delta_{b}: G_{b}(E) \ra \Lambda^{\times}$ is the character described in \cite[Definition~3.14]{HIDualCpx} defined with respect to the datum fixed in \Cref{conv: DualCpx}.
    \end{proposition}
    \begin{proof}
	    This is \cite[Corollary~1.7]{HIDualCpx}\footnote{Compare \Cref{rem: ComparisonWithDualizingComplexCalc} with \cite[Definition~3.14]{HIDualCpx} and \cite[Theorem~1.6]{HIDualCpx}).}. 
    For the convenience of the reader we recall a sketch of the argument.
	    Recall that from the Harder-Narasimhan formalism we have a Cartesian diagram 
	    \begin{center}
	    \begin{tikzcd}
		    {[\ast/\ul{G_b(E)}]} \arrow{r}{s_b} \arrow{d}{i_b}  & \Bun^b_G \arrow{d}{j_b} \\
		    \calM_b \arrow{r}{\sigma_b} & \Bun_G.
	    \end{tikzcd}
	    \end{center}
	    By \cite[Proposition~1.1]{HIDualCpx}, the dualizing complex on $\Bun_G$ is the constant sheaf.  
	    Since the map $\calM_b\to \Bun_G$ is $\calD^\an_\Lambda$-smooth, it follows from smooth base change that
	    $s_b^*\pi_b^!\Lambda\simeq i^!_b\Lambda$.
	    By definition of $\omega_{\calM_b}$, the formula 
	    \[\gamma_b^!(\omega^{-1}_{\calM_b})\simeq \gamma_b^*\omega^{-1}_{\calM_b}\otimes \gamma_b^!\Lambda \simeq \gamma_b^*\omega^{-1}_{\calM_b}\otimes \gamma_b^*\omega_{\calM_b} \simeq \Lambda\]
	    holds.
	   Since, $i_b$ is a section to $\gamma_b:\calM_b\to [\ast/\ul{G_b(E)}]$, we obtain the desired formula
	   \[s_b^*\pi_b^!\Lambda\simeq i_b^!\Lambda\simeq\omega^{-1}_{\calM_b}.\]
	   For the second part, we can appeal to the formula for $\omega_{\calM_b}$ in terms of $\delta_b$ appearing in \Cref{prop: jbnatural}.
    \end{proof}

    Under the identification $\Rep(G_b(E))\simeq \calD^\an_\Lambda([\ast/\ul{G_b{(E)}}])$, the Frobenius structure corresponding to $\bbD_{\coh,G_b(E)}$ is
    \[\Gamma_c([\ast/\ul{G_b{(E)}}],-):\calD^\an_\Lambda([\ast/\ul{G_b{(E)}}]) \ra {\rm Mod}_\Lambda.\]
    In contrast, after transfer of structure via $s_b^*$, the Frobenius structure inducing $\bbD_{\BZ,b}$ is
    \[\Gamma_c(\Bun_G^b,p_b^*(-))\simeq \Gamma_c([\ast/\ul{G_b{(E)}}],p_{b!}p^*_b-).\]
    In other words, the two structures differ by precomposing the following automorphism 
    \[(-)\otimes p_{b!}\Lambda:\calD^\an_\Lambda([\ast/\ul{G_b{(E)}}])\to \calD^\an_\Lambda([\ast/\ul{G_b{(E)}}]).\]
    This formally leads to the formula
    \[p_{b*}\bbD_{\BZ,b}p_b^*\simeq \bbD_{\coh,G_b(E)}(-)\otimes (p_{b!}\Lambda)^{-1}.\]
    Moreover, the identities
    \[\Lambda\simeq p_b^!p_{b!}\Lambda\simeq p^*_bp_{b!}\Lambda\otimes p_b^!\Lambda,\]
   together with \Cref{prop: dualcpxonClassifyingStack} give the formula  
   \[p_{b!}\Lambda\simeq \omega_{\calM_b}.\]
   And if $\dot{b}$ is chosen as in \Cref{conv: sigmastraightelements}, then we also have
   \[p_{b!}\Lambda\simeq \delta_b[2d_b].\]
   Let us summarize this discussion.
   \begin{proposition}
	   \label{computation of duality Db}
   We have equivalences
   \[p_{b*}\bbD_{\BZ,b}p_b^*\simeq \bbD_{\coh,G_b(E)}(-)\otimes \omega^{-1}_{\calM_b}\]
   and
   \[\bbD_{\BZ,b}p_b^*\simeq p_b^*\bbD_{\coh,G_b(E)}(-)\otimes p_b^!\Lambda \simeq p_b^*\bbD_{\coh,G_b(E)}(-)\otimes \pi_b^!\Lambda.\]
    Moreover, after fixing $\dot{b}$ as in \Cref{conv: DualCpx}, we may further write it as 
   \[p_{b*}\bbD_{\BZ,b}p_b^*\simeq \bbD_{\coh,G_b(E)}(-)\otimes \delta_b^{-1}[-2d_b].\]
   \end{proposition}

	\subsection{The schematic LLC category}
	In this subsection we recall the geometric properties of $\ICG$ and the structure of Zhu's local Langlands category $\Shv^!(\ICG)$. 
	First, recall that isomorphism classes of $G$-isocrystals over an algebraically closed extension of $\bbF_q$ are parametrized by the Kottwitz set $B(G)$, see \cite[\S3]{Kot85}. 
	
	We set some notation.
	Given a subset $S\subseteq B(G)$, we let $\ICG_S$ denote the substacks parametrizing those pairs $(\calE,\varphi)$, whose restriction to any geometric point corresponds to $b \in S$.
	As in the previous section, given a pair of subsets $S_1\subseteq S_2 \subseteq B(G)$, we get a map of stacks that we denote by
	\[i_{S_1,S_2}:\ICG_{S_1}\to \ICG_{S_2}.\]
	If $S_2=B(G)$ we omit it from the notation.
	Moreover, for $b \in B(G)$, we use the notation 
\[
i_b \colon \ICG_b \to \ICG, \quad i_{\leq b} \colon \ICG_{\leq b} \to \ICG, \quad i_{< b} \colon \ICG_{< b} \to \ICG
\]
to abbreviate $S_1=\{b\}$ resp. $S_1=\{b'\mid b'\leq b\}$, resp. $S_1=\{b'\mid b'< b\}$. 
We have the natural inclusion $i_{b,\leq b} \colon \ICG_{b} \to \ICG_{\leq b}$ with $i_b = i_{\leq b} \circ i_{b,\leq b}$.
The following theorem has a very rich history, Zhu obtains it by combining several references (see the discussion after \cite[Theorem 3.31]{Zhu25} and Footnote 3 above \cite[Theorem I.2.1]{FS21} for brief historical remarks).

\begin{theorem}{ \label{thm:isocG_geometry}} Fix an element $b \in B(G)$, then the following hold
	\begin{itemize}
		\item[(1)] $i_{\leq b}$ is a pfp closed embedding.
		\item[(2)] $i_{b,\leq b}$ is an affine open embedding.
		\item[(3)] The closure of $\ICG_b$ in $\ICG$ is $\ICG_{\leq b}$.
		\item[(4)] We have $\ICG_b \cong \bbB_{\proet}\underline{G_b(E)}$.
        \item[(5)] The natural map is a homeomorphism 
        \[ B(G) \simeq |\ICG|. \]
        Here $|\ICG|$ denotes the underlying topological space of $\ICG$, and $B(G)$ is the topological space described in Definition \ref{defn: topologiesonB(G)}. 
		\end{itemize}
\end{theorem}

\begin{remark}
{\label{conv: identifyinghseavesonIsocwithReps}}
As in \S\ref{ss: theanalyticcategory}, the identification $\bbB_{\proet}\underline{G_b(E)} \cong \ICG_{b}$ of \Cref{thm:isocG_geometry}(4) depends on a choice of representative $\dot{b} \in G(\Breve{E})$ of the element $b \in B(G)$, as it is really the element $\dot{b}$ that gives rise to an isocrystal whose isomorphism class is determined by $b$ and in turn a map $\ast \ra \ICG$ factoring through $\ICG_{b} \hookrightarrow \ICG$. 
As in the previous section, for most abstract purposes any choice $\dot{b}$ of $b$ would be equally useful, but as we already remarked, for computational purposes it works better to choose $\dot{b}$ wisely as in \S\ref{ss: renormalizing the choice of b}.

Once a representative $\dot{b}$ for $b \in B(G)$ has been fixed we obtain an equivalence $\ICG_{b} \cong \bb{B}_{\proet}\underline{G_{b}(E)}$.
Moreover, using \Cref{prop:classifying_stack_is_sind_placid} (3), we get identifications
\[\Rep(G_{b}(E)) \simeq \Shv^{!}(\bb{B}_{\proet}\underline{G_{b}(E)})\simeq \Shv^!(\ICG_b).\]

We will often omit these identifications from the notation. 
In situations where it is important to emphasize the precise equivalence we denote them by 
\[ \eta^{\mathrm{Zhu}}: \Rep({G_{b}(E)}) \xrightarrow{\simeq} \Shv^{!}(\ICG_{b}) \text{ or }  \eta^{\mathrm{Zhu}}: \Shv^!(\bb{B}\ul{G_b(E)})\xrightarrow{\simeq} \Shv^{!}(\ICG_{b}).\]
\end{remark}

Now we consider the category of sheaves on $\ICG$.  

\begin{theorem}{\label{thm: CategoricalPropertiesofIsoc}}
\begin{itemize}
	\item[(1)] The functor $\ul{i_{b\ast}}$ (resp. $\ul{i_{b!}}$) preserves compact objects and identifies $\Shv^!(\ICG_{b})$ with the full subcategory of $\Shv^!(\ICG)$ of all objects $F$ for which $\ul{i_{b'}^!}F \simeq 0$ (resp. $\ul{i_{b'}^\ast} F \simeq 0$) for all $b' \neq b$ in $B(G)$.
	\item[(2)] The category $\Shv^!(\ICG)$ is compactly generated and an object $F \in \Shv^!(\ICG)$ is compact if and only if $\ul{i_b^!} F \in \Shv^!(\ICG_{b}) \cong \Rep(G_b(E))$ is compact and zero for almost all $b$.
	\end{itemize}
\end{theorem}
\begin{proof}
	This follows from \cite[Propositions 3.66, 3.68, 3.69, 3.70]{Zhu25}.
\end{proof}
    As in \S\ref{ss: theanalyticcategory}, we see using \Cref{thm: CategoricalPropertiesofIsoc} that a natural set of compact generators of the category $\Shv^{!}(\ICG)$ is given by $\underline{i_{b!}}\eta^{\Zhu}(\on{c-Ind}_{K}^{G_b(E)}(\Lambda))$ for $K \subset G(E)$ a varying pro-$p$ compact open subgroup.
Similarly, we have a set of compact generators $\underline{i_{b*}}\eta^{\Zhu}(\on{c-Ind}_{K}^{G_b(E)}(\Lambda))$, and these two set of generators are related to each other under duality.

    In \S \ref{sec: TheEquivalence}, we will need a better way of computing $\ul{i_{b*}}$ and $\ul{i_{b!}}$.
  Recall that we have the category $\Shv^!(\Sht^{\rm sch}_\calG)$ of sheaves on the stack of shtukas for a parahoric $\mathcal{G}/\mathcal{O}_{E}$, as introduced in \S \ref{sec:isoc_and_shtukas}.
It is also compactly generated and the Newton map $\Nt: \Sht_{\calG}^{\mathrm{sch}} \ra \ICG$ induces  adjoint functors 
    \[{\ul{\rm Nt}_\ast} \colon \Shv^!(\Sht^{\rm sch}_\calG) \rightleftharpoons \Shv^!(\ICG) \colon {\ul{\rm Nt}^!}, \]
    by \cite[Lemma~10.100]{Zhu25}, since the map $\mathrm{Nt}$ is ind-pfp proper \cite[Lemma~3.39]{Zhu25}.
For $b \in B(G)$, we denote by $w_b \in B(\widetilde W)_{\rm str}$ its unique $\sigma$-straight representative, as in Definition \ref{defn: sigmastraightattatchedtob} (1). 
Recall that we have fixed an Iwahori group scheme of $\mathcal{I}$ of  $G$, we let $\mathcal{G} = \mathcal{I}$.
Then, as in Definition \ref{defn: sigmastraightattatchedtob} (2), for any lift $\dot w_b \in G(\breve E)$ of $w_b$ we obtain an Iwahori subgroup $I_{\dot w_{b}} \subset G_{b}(E)$. 
The choice of $\dot{w_b}$ gives rise to a map $\dot w_b:\Spec \bar{\bbF}_p \rar \Sht^{\rm sch}_{\calG,w_b}$, which induces an isomorphism $\bbB_{\text{prof\'et}}\underline{I_{\dot w_b}} \cong \Sht^{\rm sch}_{\calG, w_b}$.
Moreover, by \cite[Lemma~3.39]{Zhu25}, the restriction of ${\rm Nt} \colon \Sht^{\rm sch}_{\calG} \rar \ICG$ to $\Sht^{\rm sch}_{\calG,w_b}$ induces an ind-finite surjective morphism
	\[
	\bbB_{\text{prof\'et}}\underline{I_{\dot w_b}} \cong \Sht^{\rm sch}_{\calG,w_b} \rar \ICG_{b} \cong \bbB_{\proet}\underline{G_b(E)}.
	\] 
    Similarly, we may consider $\Sht_{\calG,\leq w_{b}}^{\sch}$ which is closed inside $\Sht_{\calG}^{\sch}$. 
    It follows from \cite[Lemma 3.38]{Zhu25} that the induced map $\Nt_{\leq w_{b}}: \Sht_{\calG,\leq w_{b}}^{\sch} \ra  \Sht_{\calG}^{\sch}$ factors through the inclusion $i_{\leq b}: \ICG_{\leq b} \hookrightarrow \ICG$.  
    We record this discussion for future use. 
    \begin{lemma}{\label{lemma: compactgeneratordiagram}}
    For $b \in B(G)$ with associated $\sigma$-straight element $w_{b}$, we have the following commutative diagram
    \begin{equation}{\label{eqn: IwahoriUniformization}}
\begin{tikzcd}
	\Sht^{\sch}_{\calG,w_{b}} \arrow[d,"\Nt_{w_{b}}"] \arrow[r,"i_{w_{b},\leq w_b}"] 
& \Sht_{\calG,\leq w_{b}} \arrow[r,"i_{\leq w_{b}}"]  \arrow[d,"\Nt_{\leq w_{b}}"]  
& \Sht_{\calG}^{\sch} \arrow[d,"\Nt"] \\
\ICG_{b} \arrow[r,"i_{b,\leq b}"] 
& \ICG_{\leq b} \arrow[r,"i_{\leq b}"] 
& \ICG
\arrow[from=1-1, to=2-2, phantom, "(A)" description]
\end{tikzcd}
    \end{equation}
    of perfect stacks, in which the square $(A)$ is Cartesian. 
    This satisfies the following properties. 
    \begin{enumerate}
    \item If we fix a choice of representative $\dot w_{b}$ of $w_{b}$ in $G(\Breve{E})$ then the natural map $\dot w_{b} \ra \Sht_{\calG,w_{b}}^{\sch}$ induces an isomorphism
    \[ \bbB_{\profet}\underline{I_{\dot w_{b}}} \simeq \Sht^{\sch}_{\calG,w_{b}} \]
    such that, under the natural identification $\ICG_{b} \simeq \bbB_{\proet}\underline{G_{b}(E)}$ induced by \Cref{thm:isocG_geometry} (4), the left vertical arrow of (\ref{eqn: IwahoriUniformization}) identifies with the map of classifying stacks $\bbB_{\profet}\underline{I_{\dot w_{b}}} \ra \bbB_{\proet}\underline{G_{b}(E)}$ induced by the inclusion $I_{\dot w_{b}} \subset G_{b}(E)$. 
    \item The map $\Nt_{w_{b}}$ is ind-pfp finite and surjective.
More precisely, under the identifications of (1) and the identifications of \Cref{lm:reps_of_compact_groups}, the Cartesian square in diagram (\ref{eqn: IwahoriUniformization}) induces commutative diagrams
    \begin{equation}
	\label{this-rando1}
		\begin{tikzcd}
			\Rep(I_{\dot w_b}) \arrow{r}{\ul{i_{w_b*}}} \arrow{d}{\on{c -Ind}_{I_{\dot w_{b}}}^{G_{b}(E)}} & \Shv^!(\Sht^{\rm sch}_{\calG}) \arrow{d}{\ul{\mathrm{Nt}_{*}}} \\
			\Rep(G_b(E)) \arrow{r}{\ul{i_{b*}}} & \Shv(\ICG)
		\end{tikzcd}
	\end{equation}
    and 
    \begin{equation}
	\label{this-rando2}
		\begin{tikzcd}
			\Rep(I_{\dot w_b}) \arrow{r}{\ul{i_{w_b!}}} \arrow{d}{\on{c -Ind}_{I_{\dot w_{b}}}^{G_{b}(E)}} & \Shv^!(\Sht^{\rm sch}_{\calG}) \arrow{d}{{\rm \ul{Nt_\ast}}} \\
			\Rep(G_b(E)) \arrow{r}{\ul{i_{b!}}} & \Shv(\ICG)
		\end{tikzcd}
	\end{equation}
    after applying $\Shv^{!}(-)$. 
    \end{enumerate}
    \end{lemma}
    \begin{proof}
	    Let us deal with the first statement.
	    Most of it was explained in the paragraph above, that the square (A) is Cartesian follows from \cite[Proposition~3.29]{Zhu25}.
	    For the second statement, the claim that $\Nt_{w_{b}}$ is ind-pfp finite and surjective in Part (2) follows from \cite[Lemma~3.39]{Zhu25}, as explained above. 
	    To see that the induced commutative square \eqref{this-rando1} holds, we pass to the lower-$\ast$ functors on \eqref{eqn: IwahoriUniformization} and apply \Cref{prop:classifying_stack_is_sind_placid} (4).
	    To see that \eqref{this-rando2} holds, we observe that $\ul{i_{w,\leq w,*}}\simeq \ul{i_{w,\leq w,!}}$ and that $\ul{\Nt_{w_{b},*}}\simeq \ul{\Nt_{w_{b},!}}$.
	    Indeed, these observations follow from \cite[Proposition~10.73 (1)]{Zhu25}, \cite[Lemma~10.100]{Zhu25} and \cite[Lemma~3.8]{Zhu25}.
    \end{proof}

    To formulate \Cref{lemma: compactgeneratordiagram} a choice of $\dot{w}_b$ was involved.
    For many purposes, any choice works equally well, we spell our preferred choice of $\dot{w}_b$ in \Cref{conv: sigmastraightelements} below. 

    Just as we have seen in the analytic context and in the case of the classifying stack in \S \ref{ss: BZDuality}, $\Shv^{!}(\ICG)$ carries a Frobenius structure inducing a self-duality which is an involutive anti-equivalence on the subcategory of compact objects. 
    Zhu constructs the duality on $\Shv^{!}(\ICG)$ by bootstrapping Verdier duality on pfp schemes along a sind-presentation of $\ICG$. 
    To explain this, we let $\mathcal{G}=\calI$ be our chosen Iwahori model for $G$ and we may form the simplicial resolution along the Newton map from \eqref{Leibniz-map}
    \begin{equation}{\label{eqn: ColimitPresentationIsocSection}}
     \ICG = \colim_{n \in \Delta^{\op}} \Hk^{n}(\Sht_{\mathcal{G}}^{\sch}).
    \end{equation}
     Here all the $\Hk^{n}(\Sht_{\mathcal{G}}^{\mathrm{sch}})$ are ind-very-placid stacks, and the colimit formula holds in $\SchStk_\et$. 
     Moreover, 
     these admit ind-presentations of the form
     \[\Hk^{n}(\Sht_{\mathcal{G}}^{\sch}) := \colim_{w_{\bullet}} (\Sht_{\leq w_\bullet,\mathcal{G}}^{\sch}),\] where each $(\Sht_{\leq w_{\bullet},\mathcal{G}}^{\sch})\subseteq \Hk^{n}(\Sht_{\mathcal{G}}^{\mathrm{sch}})$ is a very placid stack, and the transition maps are closed immersions. 
     Here $w_{\bullet}$ ranges over appropriate tuples $w_{\bullet}=(w_0,\dots, w_n)$ such that each of the entries is an element in the affine Weyl group i.e., $w_i\in \widetilde{W}$  (see \cite[\S3.4.2]{Zhu25} for details). 
     Using this presentation of $\Hk^{n}(\Sht_{\mathcal{G}}^{\sch})$, Zhu defines a compatible system of generalized constant sheaves $\Lambda^{\eta}$, in the sense of \cite[Definition~10.165]{Zhu25}, which is denoted as 
     \[\Lambda^\eta=\Lambda_{\Sht_{\leq w_\bullet}}^{\on{can}}.\]
We refer the reader to \cite[\S3.4.2, Page~107-109]{Zhu25} for the precise definition of the family $\Lambda^\eta$.
After fixing $\Lambda^\eta$, the $\eta$-global sections functor defines a Frobenius structure 
    \[ R\Gamma_{\eta,\mathrm{f.g}}(\Hk^{n}(\Sht_{\mathcal{G}}^{\sch}), -): \IndShv^{!}(\Hk^{n}(\Sht_{\mathcal{G}}^{\sch})) \ra {\rm Mod}_\Lambda \]
    from the Ind-completion of finitely generated sheaves to the category of $\Lambda$-modules (as in \Cref{finitely-gen-ind-placid}). 
    Moreover, precomposing with the functor 
    \[\Psi^L:\Shv^{!}(\Hk^{n}(\Sht_{\mathcal{G}}^{\sch})) \ra \IndShv^{!}(\Hk^{n}(\Sht_{\mathcal{G}}^{\sch}))\]
    of \eqref{the-left-adjoint-useful} defines a Frobenius structure 
    \[ R\Gamma_{\eta}(\Hk^{n}(\Sht_{\mathcal{G}}^{\sch}), -): \Shv^{!}(\Hk^{n}(\Sht_{\mathcal{G}}^{\mathrm{sch}})) \ra {\rm Mod}_\Lambda. \]
    The compatibilities of this family of generalized constant sheaves gives rise to a simplicial functor, which overall give rise to a Frobenius structure in $\Shv^!(\ICG)$,
    \begin{equation}{\label{eqn: FrobeniusStructureonIsoc}}
     R\Gamma_{\can}(\ICG,-): \Shv^{!}(\ICG) \ra {\rm Mod}_\Lambda 
     \end{equation}
    via the presentation $\Shv^{!}(\ICG_{G}) \simeq \colim_{n \in \Delta^{\op}} \Shv^{!}(\Hk^{n}(\Sht_{\mathcal{G}}^{\sch})$ induced by (\ref{eqn: ColimitPresentationIsocSection}) and \Cref{prop: Shvshriektheoryisnice}. 
\begin{definition}{\label{defn: idBZDualityFunctor}}
We let 
\[\id_\BZ:\Shv^!(\ICG)\simeq \Shv^!(\ICG)^\vee\]
denote the duality that the Frobenius structure $R\Gamma_{\can}(\ICG,-)$ defines. 
\end{definition}
   As in \S \ref{ss: theanalyticcategory}, this is intimately related with the cohomological duality functor $\bb{D}_{\coh,G_{b}(E)}: \Rep(G_{b}(E)) \xrightarrow{\simeq} \Rep(G_{b}(E))^{\vee}$ discussed in \Cref{ss: BZDuality} under an identification $\Rep(G_{b}(E)) \simeq \Shv^{!}(\ICG_{b})$, as in \ref{conv: identifyinghseavesonIsocwithReps}. 
   Now we have the following.
    \begin{proposition}{\label{prop: BasicPropertiesofBZDualityonIsoc}}
   The duality
    \[ \id_{\BZ}: \Shv^{!}(\ICG) \xrightarrow{\simeq} \Shv^{!}(\ICG)^{\vee} \]
        enjoys the following properties: 
    \begin{enumerate}
    \item The functor restricts to an involutive equivalence 
    \[ \id_\BZ^\omega: \Shv^{!}(\ICG)^{\omega} \xrightarrow{\simeq} \Shv^{!}(\ICG)^{\omega,\op} \]
    on subcategories of compact objects.
    \item For any parahoric group scheme $\calG$ of $G$ we have a commutative diagram 
	   \[ 
	    \begin{tikzcd}
		    \Shv^!(\Sht_\calG)^\omega     \arrow{r}{\ul{\Nt_*}}\arrow{d}{\id_{\Sht_\calG}^\omega}  & \Shv^!(\ICG)^\omega \arrow{d}{\id_\BZ^\omega} \\
		    \Shv^!(\Sht_\calG)^{\omega,\op}   \arrow{r}{(\ul{\Nt_*})^\op} & \Shv^!(\ICG)^{\omega, \op}.
	    \end{tikzcd}
    \]

    \item If $A$ is supported on a finite closed subset $Z \subset B(G)$ (i.e., $\ul{i^!_b}A=0$ if $b\notin Z$), then so is $\id_\BZ(A)$. 
    \end{enumerate}
    \end{proposition}
    \begin{proof}
    The first part is formal, and the second part is precisely \cite[Proposition~3.82]{Zhu25}. 
    The third part follows the proof of \cite[Proposition~3.84]{Zhu25}, which we sketch in what follows. 
    The full subcategory 
    \[\ul{i_{Z*}}: \Shv^{!}(\ICG_{Z}) \ra \Shv^{!}(\ICG) \]
    has compact generators of the form $\ul{\Nt_*}\calF$, with $\calF\in \Shv^!(\Sht^\sch_\calG)^\omega$ of the form $\calF\simeq \ul{i_{w_b*}}A$ where 
    \[i_{w_b}:\Sht^\sch_{\calG,{w_b}}\to \Sht^\sch_\calG\]
    denotes the immersion, and $b\in Z$.
    Then 
    \[\id^\omega_\BZ( \ul{\Nt_*}\calF)\simeq \ul{\Nt_*}\ul{i_{w_b!}}B\]
    for some $B\in \Shv^!(\Sht^\sch_{\calG,w_b})^\omega$.
    Since the image of $\Sht^\sch_{\calG,\leq w_b}$ under $\Nt$ is $\ICG_{\leq b}$, and by assumption $Z$ is closed, the claim follows. 
    \end{proof}

    As in \S \ref{ss: theanalyticcategory}, we will interpret \Cref{prop: BasicPropertiesofBZDualityonIsoc}.(3) as saying that Zhu's canonical duality $\id_\BZ$ is an equivalence in the category of semi-orthogonal decompositions (in the sense of \Cref{def: semi-orthogonal} and \Cref{rem: alternative defininition of semi-orthogonal}). 
    As a consequence of this interpretation, we obtain a natural duality 
    \[\id_{\BZ,b}:\Shv^!(\ICG_b)\xrightarrow{\simeq}\Shv^!(\ICG_b)^\vee,\]
    fitting in a commutative diagram analogous to \eqref{eqn: semiorthogonaldecomposition},

        \begin{equation}{\label{eqn: semiorthogonaldecompositionIsoc}}
        \begin{tikzcd}
		\Shv^!(\ICG_b)^\omega \arrow[d,"\ul{i_{b*}}"] \arrow[r,"\id_{\BZ,b}"] & \Shv^!(\ICG_b)^{\omega,\op} \arrow[d,"(\ul{i_{b!}})^{\op}"]  \\
		\Shv^!(\ICG)^\omega \arrow[r,"\id_{\BZ}"]   & \Shv^!(\ICG)^{\omega,\op}   
        \end{tikzcd}
        \end{equation}
     (see Corollary \ref{cor: SelfDualizableIsoc}, for details).
    It follows that $\id_{\BZ,b}$ comes from a Frobenius structure on the symmetric monoidal category
    \[(\Shv^!(\ICG_b),\otimes^!)\simeq (\Rep(G_b(E)),\otimes),\]
    for all $b \in B(G)$.
In particular, by the relationship between the cohomological duality $\bb{D}_{\coh,G_{b}(E)}$ and Frobenius structure on $\Shv^{!}(\bb{B}\underline{G_{b}(E)})$ explained in \S \ref{ss: BZDuality}, we have an identification 
    \[ \id_{\BZ,b}\eta^{\mathrm{Zhu}}(-) \simeq (\eta^{\mathrm{Zhu}} \bb{D}_{\coh,G_{b}(E)}(-) )\otimes \Lambda^{\mathrm{can}}_{b}, \]
    of functors 
    \[  \Rep(G_{b}(E))^{\omega} \rightarrow \Shv^{!}(\ICG)^{\omega,\op}, \]
    where the first identification is as in \Cref{conv: identifyinghseavesonIsocwithReps}, and $\Lambda^{\mathrm{can}}_{b} \in \Shv^{!}(\ICG_{b})$ an invertible sheaf.
    Moreover, from the proof of \cite[Proposition 3.84]{Zhu25} one knows that $\Lambda^{\mathrm{can}}_{b}$ is concentrated in degrees in $2d_{b}$, where $d_{b} := \langle 2\rho_{G},\nu_{b} \rangle$. 
    Under the identification $\Shv^{!}(\ICG_{b}) \simeq^{\eta^{\mathrm{Zhu}}} \Rep(G_{b}(E))$ discussed in \ref{conv: identifyinghseavesonIsocwithReps}, it follows that $\Lambda^{\mathrm{can}}_{b}\simeq\eta^{\on{Zhu}} \delta^{-1}_{b,\on{Zhu}}[-2d_b]$ for some character 
    \[\delta_{b,\mathrm{Zhu}}: G_{b}(E) \ra \Lambda^{\times}.\] 
    Unfortunately, the precise form of this character has not been computed in the current version of \cite{Zhu25}. 
    However, after private communication with Zhu, the following assumption seems justified.
    \begin{Assumption/Conjecture}{\label{assumpconj: linebundletwist}}
	    After fixing $\dot{b}$ as in \Cref{conv: DualCpx}, we have an equality 
    \[ \delta_{b,\mathrm{Zhu}} = \delta_{b}, \]
    of characters $G_{b}(E) \ra \Lambda^{\times}$, where $\delta_{b}$ is the character in \cite[Definition~3.14]{HIDualCpx} defined with respect to the datum fixed in \ref{conv: DualCpx}.
    \end{Assumption/Conjecture}

    For our purposes, we will need the following weaker version of \Cref{assumpconj: linebundletwist}. 

    \begin{proposition}{\label{prop: deltabZhuisweaklyunramified}}
		    The character $\delta_{b,\mathrm{Zhu}}$ is weakly unramified in the sense of \cite[\S~3.3.1]{Hai14}, that is, it factors through the Kottwitz map $\kappa_{G_b} \colon G_b(E) \to \pi_1(G_b)_{\Gamma_E}$. 	
    \end{proposition}
    \begin{proof}
    We first claim that it suffices to show that the character $\delta_{b,\Zhu}$ has trivial restriction to an Iwahori $I_{w_{b}} \subset G_{b}(E)$.
Indeed, by \cite[Proposition 11.5.4]{KP21}, $\ker(\kappa_{G_b})=G_b(E)^0$, where $G_b(E)^0$ is as in \cite[Definition 2.6.23]{KP21}. However, by \cite[Theorem 7.5.3]{KP21} and the fact that parabolic subgroups in a Tits system generate the ambient group, it follows that the parahoric subgroups of $G_b(E)$ generate $G_b(E)^0$.
Moreover, any parahoric is clearly generated by all Iwahori subgroups contained in it.
Thus the claim follows from the fact that all Iwahori subgroups of $G_b(E)$ are conjugate.

	    Let $\calG=\calI$ be an Iwahori model of $G$.
	    Let $w_b$ denote the $\sigma$-straight element lifting $b$.
	    Recall that $\Sht_{\calG,\leq w_b}^{\sch}$ has $\bbB I_{w_b}\simeq\Sht_{\calG,w_b}^{\sch}$ as an open substack, and that $I_{w_b}\subseteq G_b(E)$ is an Iwahori subgroup.
To show that $\delta_{b,\mathrm{Zhu}}$ has trivial restriction to the Iwahori $I_{w_b}$, in turn translates to finding an isomorphism
	    \[\on{c-Ind}_{I_{w_b}}^{G_{b}(E)}(-)\otimes \delta_{b,\mathrm{Zhu}}\simeq  \on{c-Ind}_{I_{w_b}}^{G_{b}(E)}(- \otimes \delta_{b,\Zhu}|_{I_{w_{b}}}) \simeq \on{c-Ind}_{I_{w_b}}^{G_{b}(E)}(-). \]
	    Indeed, 
	    \[\on{c-Ind}_{I_{w_b}}^{G_{b}(E)}(\mathbbm{1})\simeq\on{c-Ind}_{I_{w_b}}^{G_{b}(E)}(\delta_{b,\Zhu}|_{I_{w_{b}}})\]
	    already implies that $\delta_{b,\Zhu}|_{I_{w_{b}}}\simeq \mathbbm{1}$ by adjunction, and since we always have a non-zero map 
	    \[\on{c-Ind}_{I_{w_b}}^{G_{b}(E)}(\mathbbm{1})\to \mathbbm{1}.\]
Using \Cref{lemma: compactgeneratordiagram} (2), we reduce to constructing a natural equivalence
	    \[\ul{\Nt_{w_b,*}}(-) \otimes^! \Lambda^{\can}_{b} \simeq \ul{\Nt_{w_b,*}}(-\otimes^{!} \underline{\mathrm{Nt}}^{!}_{w_{b}}(\Lambda^{\can}_{b})) \simeq  \ul{\Nt_{w_b,*}}(-) [-2d_b].\]
	    By the projection formula, it suffices to find an identification of line bundles 
	    \[\ul{\mathrm{Nt}^{!}}_{w_{b}}(\Lambda_b) \simeq \omega_{w_{b}}[-2d_{b}],\]
	    where $\omega_{w_{b}}$ denotes the tensor unit on $\bbB I_{w_b}\simeq\Sht_{\calG,w_b}^{\sch}$.

From formal yoga on Frobenius algebras, and since the Frobenius algebra on $\ICG$ is by definition compatible with the Frobenius algebra on $\Sht_\calG$,
the line bundle $\ul{\mathrm{Nt}^{!}_{w_{b}}}(\Lambda^{\can}_{{b}})$ identifies with the sheaf $\ul{i_{w_b,\leq w_b}^!}\Lambda_{\Sht_{\calG,\leq w_b}^{\sch}}^{\on{can}}$ where $i_{w_b,\leq w_b}$ is the open immersion $\Sht_{\mathcal{G},w_{b}}^{\mathrm{sch}} \hookrightarrow \Sht_{\mathcal{G},\leq w_{b}}^{\mathrm{sch}}$.
The precise description of the sheaf $\Lambda_{\Sht_{\calG,\leq w_b}^{\sch}}^{\on{can}}$  is given in \cite[\S3.4.2, Page~107-109]{Zhu25} and we recall this now.
We set $\mathcal{F}\ell_{\mathcal{G}}^{w_{b}} := [LG_{w_{b}}/L^{+}\mathcal{G}]$ and $\mathcal{F}\ell_{\mathcal{G}}^{\leq w_{b}} := [LG_{\leq w_{b}}/L^{+}\mathcal{G}]$ to be Schubert cell (resp.
Schubert variety) attached to $w_{b}$ in the affine flag variety, where notation is as in \S \ref{sec: loop gropus etc}.
    
    Then $\Lambda_{\Sht_{\calG,w_b}^{\sch}}^{\on{can}}$ is obtained as $\ul{p^!}\Lambda_{[L^{m}\calG\backslash \Fl_\calG^{\leq w_b}]}$ for the map 
    \[\Sht_{\calG,\leq w_b}^{\sch}\xrightarrow{p} [L^{m}\calG\backslash \Fl_\calG^{\leq w_b}]\]
    for a suitably truncated version of the loop group $L^{m}\calG$ with $m$ sufficiently large with respect to $w_b$ so that the kernel of $L^+\calG\to L^{m}\calG$ acts trivially on $\Fl_\calG^{\leq w_b}$.
We claim that, on the open substack $[L^{+}_{m}\calG\backslash \Fl_\calG^{w_b}]$, we have an identification 
    \[\omega_{[L^{+}_{m}\calG\backslash \Fl_\calG^{w_b}]} \simeq \Lambda_{[L^{+}_{m}\calG\backslash \Fl_\calG^{w_b}]}[2\ell(w_b)], \]
    and we recall that by \cite[\S2.4]{He2014GeometricHomological} $\ell(w_b)=d_b$.
    This claim implies our desired statement.
    Indeed, the identification of $\Lambda_{[L^{+}_{m}\calG\backslash \Fl_\calG^{w_b}]}[2\ell(w_b)]$ with the $\otimes^!$-tensor unit produces, after applying $\ul{(p\circ i_{w,\leq w})^!}$, an identification of $\ul{\mathrm{Nt}^{!}}_{w_{b}}(\Lambda_b)[2d_b]$ with the $\otimes^!$-tensor unit as we wanted to construct.

    Finally, to prove our claim we note that $[L^{+}_{m}\calG\backslash \Fl_\calG^{w_b}]$ is a perfect Artin stack with pfp atlas $\Fl_\calG^{w_b} \ra [L^{+}_{m}\calG\backslash \Fl_\calG^{w_b}]$ given by quotienting out by a connected group scheme.
Thus, by \cite[Lemma~A.1.2.]{Zhu16}, the category of local systems on $[L^{+}_{m}\calG\backslash \Fl_\calG^{w_b}]$ is a full subcategory of the category of local systems on $\Fl_\calG^{w_b}$ via $*$-pullback.
In particular, the claim is reduced to showing that $\Fl_\calG^{w_b}$ is perfeclty finitely presented smooth of dimension $\ell(w_{b})$ over $k$ which follows from \cite[Remark 2.1.22]{Zhu16}. 
    \end{proof}

   \subsection{Renormalized functors.}
   For computational purposes, it is useful to renormalize the functors $j_{b!}$ and $j_{b\sharp}$ in the analytic setup, or $\ul{i_{b*}}$ and $\ul{i_{b!}}$ in the schematic setup.
   In the analytic setup, this is done as follows. 
	Recall that $\pi_b:\Bun_G^b\to \ast$ denotes the structure morphism.
    We can then consider a sheaf $\ICoh_{\Bun_{G}^{b}} := \sqrt{\pi_b^!\Lambda}\in \calD_{\Lambda}^{\an}(\Bun_{G}^{b})$, if it exists. 
It follows from \Cref{prop: dualcpxonClassifyingStack} and the explicit formula for $\delta_b$ computed in \cite{HIDualCpx}, that $\ICoh_{\Bun_{G}^{b}}$ exists if and only if $\Lambda$ contains a square root of $q$.
After fixing $\sqrt{q} \in \Lambda$, which we do for the rest of the section, any choice of $\ICoh_{\Bun_{G}^{b}}$ is necessarily given by the formula
\[\ICoh_{\Bun_{G}^{b}}:= p_{b}^{*}(\delta_{b}^{-\frac{1}{2}})[-d_{b}].\]
By construction, $\ICoh_{\Bun_{G}^{b}}^{\otimes 2} \simeq \pi_b^!\Lambda$ and this sheaf is Verdier self-dual in $\calD^\an_\Lambda(\Bun_G^b)$. 
      This motivates the following.
	\begin{definition}\label{def:modulus_character}
    For $b \in B(G)$ and $? \in \{!,\sharp\}$, we consider the functor
		\[j^{\on{ren}}_{b?}:\calD^\an_\Lambda(\Bun_{G}^{b})\to \calD^\an_\Lambda(\Bun_G)\]
			with the formula $j_{b?}^{\ren}(-) := j_{b?}(- \otimes \ICoh_{\Bun_{G}^{b}})$. 
	\end{definition}

     \begin{theorem}
	     \label{some other label needed to be here }
	     \begin{enumerate}
			\item For all $b \in B(G)$ and $B \in \calD_\Lambda^{\an}(\Bun_{G}^{b})^{\omega}$, we have natural identifications
			\[ \bbD_{\BZ}j_{b!}B \simeq j_{b\sharp}\bbD_{\BZ,b}B \]
            and
	    \[ \bbD_{\BZ}j_{b\sharp}B \simeq  j_{b!}\bbD_{\BZ,b}B. \]

			\item For all $b \in B(G)$ and $A \in \Rep(G_{b}(E))^{\omega}$, we have natural identifications
			\[ \bbD_{\BZ}j_{b!}^{\ren}p_{b}^{*}A \simeq j_{b\sharp}^{\ren}p_{b}^{*}\bb{D}_{\mathrm{coh},G_{b}(E)}A \]
            and
			\[ \bbD_{\BZ}j_{b\sharp}^{\ren}p_{b}^{*}A \simeq  j_{b!}^{\ren}p_{b}^{*}\bb{D}_{\mathrm{coh},G_{b}(E)}A. \]

		\item For all $b \in B(G)$ and $K \subseteq G_b(E)$ compact open pro-$p$-group, there exists (non-natural) identifications
			\[ \bbD_{\BZ}\LKB^{\Bun_G}[-d_b] \simeq \RKB^{\Bun_G}[-d_b] \]
            and
	    \[ \bbD_{\BZ}\RKB^{\Bun_G}[-d_b] \simeq  \LKB^{\Bun_G}[-d_b]. \]
		\end{enumerate}
	\end{theorem}

	\begin{proof}
		The first claim follows from diagram \eqref{eqn: semiorthogonaldecomposition}.
		The second claim follows from the first and from \Cref{computation of duality Db}. 
		Indeed,
		\begin{align*}
			\bbD_{\BZ}j_{b\sharp}^{\ren}p_{b}^{*}A & \simeq \bbD_{\BZ}j_{b\sharp}(p_{b}^{*}A\otimes \ICoh_{\Bun_{G}^{b}})\\
							       & \simeq j_{b!} \bbD_{\BZ,b}(p_{b}^{*}A\otimes \ICoh_{\Bun_{G}^{b}})\\
							       & \simeq j_{b!}( \bbD_{\BZ,b}(p_{b}^{*}A)\otimes \ICoh^{-1}_{\Bun_{G}^{b}})\\
							       & \simeq j_{b!}( (p_b^*\bb{D}_{\mathrm{coh},G_{b}(E)}(A)\otimes \pi_b^!\Lambda)\otimes \ICoh^{-1}_{\Bun_{G}^{b}})\\
							       & \simeq j_{b!}(p_b^*\bb{D}_{\mathrm{coh},G_{b}(E)}(A)\otimes \ICoh_{\Bun_{G}^{b}})\\
							       & \simeq j_{b!}^{\ren}(p_{b}^{*} \bb{D}_{\mathrm{coh},G_{b}(E)}A).
		\end{align*}

		The third claim follows from the first and the from the (non-canonical) formulas
		\[\LKB^{\Bun_G}[-d_b]\simeq j_{b\sharp}^{\ren}p_{b}^{*}\on{c-Ind}_{K}^{G}(\Lambda) \text{ and } \RKB^{\Bun_G}[-d_b]\simeq j_{b!}^{\ren}p_{b}^{*}\on{c-Ind}_{K}^{G}(\Lambda).\]
		which are a consequence of $(\delta_{b})_{\mid_K}$ being the trivial character, since $\delta_{b}$ is weakly-unramified.
	\end{proof}

   We now discuss renormalization in the schematic setup. 
    We now make the following definition. 
    \begin{definition}{\label{defn: RenormalizedPushforwardsIsoc}}
	    Suppose that $\Lambda$ contains all square roots of $\delta_{b,\on{Zhu}}(G_b(E))\subseteq \Lambda^\times$.\footnote{Under assumption \ref{assumpconj: linebundletwist}, it suffices that $\Lambda$ contains a square root of $q$.}
	    We fix a square root  $\delta_{b,\mathrm{Zhu}}^{\frac{1}{2}} = \delta_{b}^{\frac{1}{2}}$. 
	    This allows us to define the following.
    \begin{enumerate}
	    \item We define $\ICoh_{\ICG_{b}} \in \Shv^{!}(\ICG_{b})$ to be $\eta^{\on{Zhu}}\delta_{b,\mathrm{Zhu}}^{-1/2}[-d_{b}]$ where $\eta^{\on{Zhu}}$ is the identification from \Cref{conv: identifyinghseavesonIsocwithReps}.   
    \item For $? \in \{!,*\}$, we consider
    \[ \ul{i_{b?}^{\ren}}(-): \Shv^{!}(\ICG_{b}) \ra \Shv^{!}(\ICG), \] 
    via the formula $\ul{i_{b?}}(- \otimes \ICoh_{\ICG_{b}})$.
    \end{enumerate}
    \end{definition}
    As in \Cref{some other label needed to be here }, this formally leads to the following commutative diagram deduced from \eqref{eqn: semiorthogonaldecompositionIsoc},
    \begin{equation}{\label{eqn: SemiorthogonalBZdualityIsoc}}
    \begin{tikzcd}
	    \Rep(G_b(E)) \arrow[rr,"\bb{D}_{\coh,G_{b}(E)}"] \arrow[d,"\ul{i^\ren_{b!}}\eta^{\on{Zhu}}"] & & \Rep(G_b(E))^{\vee}  \arrow[d,"\ul{i^{\ren,o}_{b*}}\eta^{\on{Zhu},o}"]   \\
	    \Shv^{!}(\ICG) \arrow[rr,"\id_{\BZ}"] & & \Shv^{!}(\ICG)^{\vee}.  
    \end{tikzcd}
    \end{equation}

    \subsection{Conventions on $\dot{b}$}
    \label{ss: renormalizing the choice of b}
    In order to make the character  
    \[\delta_{b}:G_b(E)\to \Lambda^\times\]
explicit, we use the following convention for the representatives $\dot{b}$ of an element $b \in B(G)$.
    \begin{convention}{\cite{HIDualCpx}}{\label{conv: DualCpx}}
    We recall that we are assuming the group $G$ is quasi-split. In the introduction, we recall that we have fixed the following data.
    \begin{itemize} 
    \item A Borel pair $(T,B)$ of $G$.
    \end{itemize}
   Consider an element $b \in B(G)$, we consider the standard Levi factor $M_{b}$ with respect to the Borel pair given by the centralizer of $\nu_{b}$.
We recall that there is a (unique) reduction $b_{M} \in B(M_{b})_{\mathrm{basic}}$ of $b \in B(G)$ along the map $B(M_{b}) \ra B(G)$ such that the slope homomorphism applied to $b_{M}$ is $G$-dominant. In particular, we may arrange that $\dot{b}$ lies in $M_{b}(\Breve{E})$, after conjugating by an appropriate element of the Weyl group $N_{G}(T)(\Breve{E})/T(\Breve{E})$. We fix such a representative of $\dot{b}$ with respect to the data fixed above.
    \end{convention}
    \begin{remark}{\label{rem: ComparisonWithDualizingComplexCalc}}
    The identification $[\ast/\tilde{G}_{\dot{b}}] \simeq \Bun_{G}^{b}$ fixed by such a representative $\dot{b}$ will by construction satisfy the following compatibilities.
We write $P_{b}$ for the standard parabolic with respect to the fixed Borel $B$ with Levi factor given by $M_{b}$ and let $P_{b}^{-}$ denote its opposite.
The natural maps $M_{b} \leftarrow P_{b}^{-} \ra G$ induce a diagram 
    \[ \begin{tikzcd}
     \Bun_{P_{b}^{-}} \arrow[r,"\mf{p}_{P_{b}^{-}}"] \arrow[d,"\mf{q}_{P_{b}^{-}}"] & \Bun_{G} \\
    \Bun_{M_{b}} & & 
    \end{tikzcd} \]
    of $v$-stacks.
We let $\dot{b}_{M} \in M_{b}(\Breve{E})$ be the associated element representing $b_{M}$ for a choice of $\dot{b}$ as in Convention \ref{conv: DualCpx}. This determines an isomorphism 
    \[ [\ast/\underline{G_{\dot{b}}(E)}] \simeq \Bun_{M_{b}}^{b_{M}}. \]  
    We write $\Bun_{P_{b}^-}^{b_{M}} \hookrightarrow \Bun_{P_{b}^{-}}$ for the open substack given by the pullback of the open immersion $\Bun_{M_{b}}^{b_{M}} \hookrightarrow \Bun_{M_{b}}$.
The splitting of the Harder-Narasimhan filtration then determines an isomorphism 
    \[  \Bun_{P_{b}^-}^{b_{M}} \simeq [\ast/\tilde{G}_{\dot{b}}], \]
    which in turn induces an identification 
    \begin{equation}{\label{eqn: identificationofclassifyingstacks}}
     [\ast/\tilde{G}_{\dot{b}}] \simeq \Bun_{P_{b}^{-}}^{b_{M}} \xrightarrow{\simeq} \Bun_{G}^{b}, 
    \end{equation}
    via the morphism $\mf{p}_{P_{b}^{-}}$.
This in turn gives the isomorphism determined by the representative $\dot{b}$ of $b$, as in \cite[Example~V.3.4]{FS21}.
We note that the appearance of the opposite parabolic comes from the fact that the HN-slopes of the $G$-bundle $\mathcal{E}_{b}$ attached to $b$ are normalized to be the negatives of the associated isocrystal slopes.
\end{remark}

Now we specify our convention for our choice of representative for the sigma straight element.
\begin{convention}{\label{conv: sigmastraightelements}}
   For each $b \in B(G)$ and a choice of Iwahori $\mathcal{I} \subset G$ with associated $\sigma$-straight representative $w_{b} \in B(\widetilde{W})_{\mathrm{str}}$, we note that, after replacing $\dot{w}_{b}$, by a conjugate by an element in the relative Weyl group of $G_{\Breve{E}}$ as in the short exact sequence (\ref{eqn: RelativeWeylGroup}), we may arrange that the representative $\dot{w}_{b} \in N_{G}(T)(\Breve{E}) \subset G(\Breve{E})$  of $w_{b}$ is of the form described in Convention \ref{conv: DualCpx}.
We will fix such a choice and normalize the isomorphisms $\Sht_{\mathcal{G},w_{b}} \simeq \bb{B}_{\profet}\underline{I_{\dot{w}_{b}}}$ and $\ICG_{b} \simeq \bb{B}_{\proet}\underline{G_{b}(E)}$ in \Cref{lemma: compactgeneratordiagram} (1) with respect to this choice.
    \end{convention}

    \section{Semi-orthogonal decompositions in practice}
    {\label{subsec: ApplicationsofSemiorthogonal}}
    We now want to describe a semi-orthogonal decomposition, on $\Shv^{!}(\ICG)$ and $\calD_{\Lambda}^{\an}(\Bun_{G})$, with respect to the partially ordered set $B(G)$. 
    Recall that we write simply $B(G)$ for the set equipped with the topology generated by the opens $U_{b} := \{x \in B(G)| x \geq b\}$ for $b \in B(G)$, and we write $B(G)^{\op}$ for the set $B(G)$ equipped with the natural topology given by the opens $\{x \in B(G)| b \geq x \}$ for $b \in B(G)$ (see \Cref{defn: topologiesonB(G)}).
We recall that we have homeomorphisms 
   \begin{equation}
	   \label{topology B(G)}
     |\ICG| \simeq B(G)  
     \end{equation}
    and 
   \begin{equation}
	   \label{topology BunG}
    |\Bun_{G}| \simeq B(G)^{\op}, 
     \end{equation}
    by \Cref{thm: BunGGeometricFacts} (2) and \Cref{thm:isocG_geometry} (5). 
    Recall that, for any subset $S \subset B(G)$, we have substacks $i_{S}: \ICG_{S} \hookrightarrow \ICG$ and $j_{S}: \Bun_{G}^{S} \hookrightarrow \Bun_{G}$, and that if we have inclusion of subsets $S \subset S'\subseteq B(G)$, we obtain monomorphisms $i_{SS'}: \ICG_{S} \rightarrow \ICG_{S'}$ and $j_{SS'}: \Bun_{G}^{S} \rightarrow \Bun_{G}^{S'}$. 

    In particular, by \eqref{topology B(G)}, whenever we have an inclusion $V' \subset V \subset B(G)$ of convex subsets, we obtain a natural locally closed immersion $i_{V'V}: \ICG_{V'} \ra \ICG_{V}$.
This is an open (resp.~closed) immersion if and only if $V'\subset V$ is an open (resp.~closed) subset.
    Similarly, by \eqref{topology BunG}, for the same inclusion of convex subsets $V' \subset V \subset B(G)$ we obtain an locally closed immersion $j_{V'V}: \Bun_{G}^{V'} \hookrightarrow \Bun_{G}^{V}$.
On the other hand $j_{V'V}$ is an open (resp.~closed) immersion if and only if $i_{V'V}$ is a closed (resp.~open) immersion.
Indeed, $V'\subset V \subset B(G)^\op$ carries the opposite topology. 

    \subsubsection{Semi-orthogonal decompositions on $\Shv^{!}(\ICG)$}
    \label{subsection semi on isoc}
    We want to construct a semi-orthogonal decomposition on the category $\Shv^{!}(\ICG)$ with respect to the partially ordered set $B(G)$. 
    To avoid technicalities, we only consider a semi-orthogonal decomposition of the category $\Shv^!(\ICG_X)$ with respect to the partially ordered set $B(G)_X$ when $X\subseteq B(G)$ is a finite convex subset.

    Consider the map
    \[\Corr(\Convex_{B(G)_X},\on{All})\to \Corr(\PreStk,E_\pfp).\]
    with formula
    \[V\mapsto \ICG_V,\]
    and 
    \[[V\leftarrow W \rightarrow U]\mapsto [\ICG_V\leftarrow \ICG_W \rightarrow \ICG_U], \]
    which is well-defined by \Cref{thm:isocG_geometry} and \Cref{lem_finite_B(G)_order}.
    We may compose it with 
    \[
    \Shv^!:\Corr(\PreStk,E_\pfp)\to \LinCat_\Lambda\]
    to obtain a functor 
    \begin{equation}
	    \label{eq: the-pre-semi}
	    \bbS^{\on{ex}}_{\ICG_X}:\Corr(\Convex_{B(G)_X},\on{All})\to \LinCat_\Lambda.
    \end{equation}
    This functor encodes on a correspondence $V\xleftarrow{i_{WV}} W\xrightarrow{i_{WU}} U$ the functor $\ul{i_{WU,*}}\circ \ul{i_{WV}^!}$. 
    \begin{remark}
	    We observe that $\bbS^{\on{ex}}$ is a semi-orthogonal decomposition on $\Shv^!(\ICG_X)$ with respect to $X^\op\subseteq B(G)^\op$.
    We call this the \textit{exceptional semi-orthogonal decomposition} for this reason, as the natural topology on $\ICG$ is given by $B(G)$.
    Observe that the right adjoint functors of $\ul{i_*}$ and $\ul{i^!}$ are $\ul{i^\sharp}$ and $\ul{i_\flat}$ respectively.
    In particular, if $i^U:U\to \ICG_X$ and $i^Z:Z\to \ICG_X$ are complementary open and closed immersions in $\SchStk_\et$ the natural triangles arising from this semi-orthogonal decomposition are 
    \[\ul{i^Z_*}\ul{i^{Z,!}}A \to A \to \ul{i^U_*}\ul{i^{U,*}}A\]
    and 
    \[\ul{i^U_*}\ul{i^{U,\sharp}}A\to A \to \ul{i^Z_\flat}\ul{i^{Z,!}}A. \]
    \end{remark}

    We will find it more convenient to work with the semi-orthogonal decomposition with respect to $X\subseteq B(G)$, which is the natural topology $|\ICG|$.
    For this, we may pass to left adjoints (which exist by \Cref{prop: semiorthogonaldecompositionsonsindplacid}).
    Using that $\Corr(\Convex_{B(G)_X},\on{All})^\op\simeq \Corr(\Convex_{B(G)_X},\on{All})$ we obtain a functor 
    \[\bbS_{\ICG_X}:\Corr(\Convex_{B(G)_X},\on{All})\to \LinCat_\Lambda,\]
    which encodes the rule that on a correspondence $V\xleftarrow{i_{WV}} W\xrightarrow{i_{WU}} U$ produces the functor $\ul{i_{WU,!}}\circ \ul{i_{WV}^*}$. 
    By \Cref{prop: semiorthogonaldecompositionsonsindplacid}, this functor is a semi-orthogonal decompositions of $\Shv^!(\ICG_X)$.

    \begin{remark}
	    As we mentioned $\bbS_{\ICG_X}$ is Nagata with respect to $X\subseteq B(G)$.  
	    In particular, if $i^U:U\to \ICG_X$ and $i^Z:Z\to \ICG_X$ are complementary open and closed immersions in $\SchStk_\et$ the natural triangles are the standard ones
	    \[\ul{i^Z_*}\ul{i^{Z,!}}A \to A \to \ul{i^U_*}\ul{i^{U,*}}A\]
    and 
    \[\ul{i^U_!}\ul{i^{U,*}}A\to A \to \ul{i^Z_*}\ul{i^{Z,*}}A. \]
    Informally, $\bbS_{\ICG_X}$ captures the following data. 
    For every finite convex set $V\subseteq X$ we obtain a category $\Shv^!(\ICG_V)$.
    Moreover, every inclusion of finite convex subsets $V\subseteq V'$ we obtain pairs of adjoint functors 
   \[\ul{i_{VV'!}}\dashv \ul{i^!_{VV'}} \text{ and } \ul{i^*_{VV'}}\dashv \ul{i_{VV'*}},\]
   that satisfy the usual base change formulas and the usual identifications $\ul{i^!_{VV'}}\simeq \ul{i^*_{VV'}}$ in the case of open immersions, and $\ul{i_{VV'*}}\simeq \ul{i_{VV'!}}$ in the case of closed immersions.  
    \end{remark}

   From \Cref{thm: CategoricalPropertiesofIsoc}, it follows that $\ul{i_{VV'*}}$ and $\ul{i_{VV'}^!}$ preserve compact objects, and we may pass to conjugate functors to define a semi-orthogonal $B(G)_X$-decomposition on $\Shv^!(\ICG_X)^\vee$ of the form 
    \[\bbS^\vee_{\ICG_X}:\Corr(\Convex_{B(G)_X},\on{All})\to \LinCat_\Lambda.\]
    This functor encodes on a correspondence $V\xleftarrow{j} W\xrightarrow{i} U$ the functor $(\ul{i_*})^o\circ (\ul{j^!})^o$. 

    It follows from \Cref{prop: BasicPropertiesofBZDualityonIsoc} and \Cref{thm: equivalence criterion semi-orthogonal functor few criteria} that the duality operation promotes uniquely to an equivalence of presentable 2-functor formalisms
    \[\id_{\BZ,X}:\bbS_{\ICG_X}\xrightarrow{\simeq}\bbS^\vee_{\ICG_X}.\]

    We summarize the above discussion with the following.
    \begin{corollary}
	    \label{cor: SelfDualizableIsoc}
	    Fix $X\subseteq B(G)$ a finite convex subset.
The following statements hold. 
        \begin{enumerate}
		\item The category $\Shv^{!}(\ICG_{X})$ is compactly generated with compact generators as in \Cref{thm: CategoricalPropertiesofIsoc} (2). 
		Moreover, the duality $\id_{\BZ}$ induces a duality
        \[ \id_{\BZ,X}: \Shv^{!}(\ICG_{X}) \xrightarrow{\simeq} \Shv^{!}(\ICG_{X})^{\vee} \]
        which is an involutive contravariant equivalence on the subcategory of compact objects.
\item For any inclusion of finite convex subsets $V\subseteq V'$ of $B(G)$, we then have a natural equivalence
	\[ \id_{\BZ,V'}\ul{i_{VV'!}} \simeq (\ul{i_{VV'*}})^{o}\id_{\BZ,V}\]
\text{ and }
\[\id_{\BZ,V}\ul{i^*_{VV'}} \simeq (\ul{i^!_{VV'}})^{o}\id_{\BZ,V'}.\]
        \end{enumerate}
        \end{corollary}

    \subsubsection{Semi-orthogonal decompositions on $\calD^{\an}_{\Lambda}(\Bun_{G})$}\label{sec:semiorthogonalDecOfDBunG}
    We first turn to the classical semi-orthogonal decomposition on $\calD_{\Lambda}^{\an}(\Bun_{G})$.  
    As in the previous section, we restrict to study $\calD_{\Lambda}^\an(\Bun_G^X)$ for a fixed convex subset $X\subseteq B(G)$.

    The natural topology of $\Bun_G^X$ endows through \Cref{exceptional anlaytic}.(1) with a semi-orthogonal $X^\op$-decomposition (i.e., Nagata with respect to $X^\op$), which we denote 
    \[\bbS_{\Bun^X_G}:\Corr(\Convex_{X^\op},\on{All})\to \LinCat_\Lambda.\]
    This functor encodes the rule
    \[V\mapsto \calD^\an_\Lambda(\Bun_G^V)\]
    and
    \[[V\xleftarrow{j_{WV}} W \xrightarrow{j_{WU}} U]\mapsto j_{WU!} \circ  j^*_{WV}.\]

    \begin{remark}
	    If $j^Z:Z\to \Bun_G^X$ and $j^U:U\to \Bun_G^X$ denote complementary closed an open immersions in $\AnStk_v$, then the excision triangles with respect to $\bbS_{\Bun_G^X}$ are 
	    \[{j^Z_*}{j^{Z,!}}A \to A \to {j^U_*}{j^{U,*}}A\]
    and 
    \[{j^U_!}{j^{U,*}}A\to A \to {j^Z_*}{j^{Z,*}}A. \]
    \end{remark}

    Now we wish to apply \Cref{exceptional anlaytic}.(2) to the above setup, in order to obtain a semi-orthogonal decomposition with respect to $X$. 
    For this, we justify that $j_{WU!}$ and $j^\sharp_{WV}$ admit left adjoint functors.
    \begin{proposition}
	    \label{exteince of left adjoints}
	    Let $W\subseteq V\subseteq B(G)$ denote two finite convex subsets then $j_{WV!}$ admits a left adjoint $j^\flat_{WV}$ and $j^*_{WV}$ admits a left adjoint $j_{WV\sharp}$.
    \end{proposition}
    \begin{proof}
	    We only argue that $j^*_{WV}$ admits a left adjoint since the argument for the other functor is analogous. 
	    We divide our analysis to the cases in which $W\subseteq V$ is closed or open.
	    In the first case, the map of stacks $j_{WV}:\Bun_G^W\to \Bun_G^V$ is open (by the reversal of topologies) and in particular $j_{WV\sharp}:=j_{WV!}$ is already a left adjoint to $j_{WV}^*$.
	    In the second case, $j_{WV}$ is a closed immersion and we construct $j_{WV\sharp}$ as the ind-extension of the functor which on compacts takes the form 
	    \[j_{WV\sharp}(A):= \bbD_{\BZ, V} \circ j^{\op}_{WV!}\circ \bbD_{\BZ, W}(A).\]
    We compute on compacts 
    \begin{align*}
    \mathrm{Hom}_{\calD_{\Lambda}^{\an}(\Bun_{G}^{V})^\omega}(\bb{D}_{\BZ,V}j_{WV!}^{\op}\mathbb{D}_{\BZ,W}A,B) \simeq &  \mathrm{Hom}_{\calD^{\an}_{\Lambda}(\Bun_{G}^{V})^{\omega,\op}}(j_{WV!}^{\op}\mathbb{D}_{\BZ,W}A,\bb{D}_{\BZ,V}B)   \\
    \simeq & \mathrm{Hom}_{\calD^{\an}_{\Lambda}(\Bun_{G}^{W})^{\omega,\op}}(j^{\op}_{WV*}\mathbb{D}_{\BZ,W}A,\bb{D}_{\BZ,V}B) \\ 
    \simeq & \mathrm{Hom}_{\calD^{\an}_{\Lambda}(\Bun_{G}^{W})^{\omega,\op}}(\mathbb{D}_{\BZ,W}A,j^{*,\op}_{WV}\bb{D}_{\BZ,V}B) \\ 
     \simeq & \mathrm{Hom}_{\calD^{\an}_{\Lambda}(\Bun_{G}^{W})}(A,\mathbb{D}_{\BZ,W}j^{*,\op}_{WV}\bb{D}_{\BZ,V}B) \\
     \simeq & \mathrm{Hom}_{\calD^{\an}_{\Lambda}(\Bun_{G}^{W})}(A,j^*_{WV}B), 
    \end{align*}
    The identification $\mathbb{D}_{\BZ,W}j^{*,\op}_{WV}\bb{D}_{\BZ,V}B \simeq j^*_{WV}B$ follows from Yoneda, the Frobenius algebra interpretation of $\bb{D}_{\BZ,V}$ and $\bbD_{\BZ,W}$, and from the projection formula. 
    This upgrade to arbitrary objects by Ind-extending.
    \end{proof}
\begin{remark}\label{can extend exceptional adjoints to all of BunG}
    The formation of \(j^\flat_{WV}\) and \(j_{WV\sharp}\) is compatible with varying \(V\), thus we also obtain a left adjoint $j_W^\flat$ for $j_{W!}$ and a left adjoint $j_{W\sharp}$ for $j_W^*$.
\end{remark}
    As we have discussed, using \Cref{exceptional anlaytic}.(2) and \Cref{exteince of left adjoints} we obtain an exceptional semi-orthogonal decomposition of $\calD^\an_\Lambda(\Bun_G^X)$ of the form 
    \[\bbS^{\on{ex}}_{\Bun_G^X}:\Corr(\Convex_X,\on{All})\to \LinCat_\Lambda.\]
    This encodes, for a correspondence $V\xleftarrow{j_{WV}} W \xrightarrow{j_{WU}} U$ the functor $j_{WU\sharp}\circ j^\flat_{WV}$.
        As before, for every pair of finite convex subsets $V\subseteq V'\subseteq X$ we automatically get pairs of adjoint functors
    \[j_{VV'\sharp}\dashv j^*_{VV'} \text{ and } j^\flat_{VV'} \dashv j_{VV'!},\]
    such that if $V\to W$ is closed, resp.~open, with the subspace topology of $B(G)$ (i.e., open, resp.~closed, with respect to the subspace topology of $B(G)^\op$), then $j_{VV'\sharp}\simeq j_{VV'!}$, resp.~$j^*_{VV'}\simeq j^\flat_{VV'}$.   

    \begin{remark}
	    If $j_{ZX}:\Bun_G^Z\to \Bun_G^X$ and $j_{UX}:\Bun_G^U\to \Bun_G^X$ denote complementary closed and open immersions in $\AnStk_v$, then the excision triangles with respect to $\bbS^{\on{ex}}_{\Bun_G^X}$ are 
	    \[{j_{ZX\sharp}{j_{ZX}^{*}}A \to A \to {j_{UX!}}j_{UX}^{\flat}}A\]
    and 
    \[j_{UX!}j_{UX}^{*}A\to A \to j_{UZ\sharp}j_{UZ}^{*}A. \]
    \end{remark}

    Reasoning as in the previous subsection (using \Cref{thm: compactenerationofBunG}), we can obtain from $\bbS_{\Bun^X_G}$ a semi-orthogonal $X$-decomposition on the dual category $\calD_{\Lambda}^\an(\Bun_G^X)^\vee$ 
    \[\bbS^\vee_{\Bun^X_G}:\Corr(\Convex_{X},\on{All})\to \LinCat_\Lambda.\]
    This encodes, for a correspondence $V\xleftarrow{j_{WV}} W \xrightarrow{j_{WU}} U$ the functor $j^o_{WU!}\circ j^{*,o}_{WV}$.
    Using \Cref{thm: naturaland!exchangedunderBZ} and \Cref{thm: equivalence criterion semi-orthogonal functor few criteria}, we get a duality identification of presentable 2-functor formalisms
    \[\bbD_{\BZ,Z}:\bbS^{\on{ex}}_{\Bun_G^X}\xrightarrow{\simeq} \bbS^\vee_{\Bun^X_G}.\]
    We summarize this with the following statement

        \begin{corollary}
	\label{cor: SelfDualizableBunG}
	    Let $X\subseteq B(G)$ denote a finite convex subset. 
	    The following statements hold. 
        \begin{enumerate}
		\item The categories $\calD^\an_\Lambda(\Bun^X_G)$ is compactly generated with compact generators as in \Cref{thm: compactenerationofBunG}. 
		Moreover, the duality $\bbD_{\BZ}$ induces a duality
        \[ \bbD_{\BZ,X}: \calD^\an_\Lambda(\Bun^X_G)  \xrightarrow{\simeq} \calD^\an_\Lambda(\Bun^X_G)^\vee\]
        which is an involutive contravariant equivalence on the subcategory of compact objects.
\item For any inclusion of finite convex subsets $V\subseteq V'$ we have identities

	\[\bbD_{\BZ,V'}j_{VV'\sharp}\simeq (j_{VV'!})^o\bbD_{\BZ,V}\]
	and 
	\[\bbD_{\BZ,V}j^\flat_{VV'}\simeq (j^*_{VV'})^o\bbD_{\BZ,V'}.\]
        \end{enumerate}
        \end{corollary}
For the convenience of the reader, we explain in more detail a specific part of this semi-orthogonal decomposition.
    \begin{example}
    If we consider the open immersion $j_{< b}: \Bun_{G}^{< b} \hookrightarrow \Bun_{G}^{\leq b}$ with its closed complement $j_{b}: \Bun_{G}^{b} \hookrightarrow \Bun_{G}^{\leq b}$ we see that the functor $j_{b}^{*}: \calD_{\Lambda}^{\an}(\Bun_{G}^{\leq b}) \ra \calD_{\Lambda}^{\an}(\Bun_{G}^{b})$ admits an exceptional left adjoint. 
    This is precisely the functor $j_{b\sharp}$ described in \Cref{prop: jbnatural}, which will satisfy the identity that $\bb{D}_{\BZ,\leq b}j_{b!}^{o}\bb{D}_{\BZ,b} \simeq j_{b\sharp}$ (\Cref{thm: naturaland!exchangedunderBZ}). 
    Analogously, we have the identity $j_{<b}^\flat\simeq \bb{D}_{\BZ,< b}j_{< b}^{*}\bb{D}_{\BZ,\leq b}$ for the left adjoint of $j_{< b!}$. 
    In particular, in general we note that, for an inclusion of closed subsets $Z' \subset Z \subset B(G)$ with $Z$ finite and $U = Z \setminus Z'$, the full diagram of adjoints attached to the recollement given by the closed subcategory $j_{Z'Z!}: \calD^{\an}_{\Lambda}(\Bun_{G}^{Z'}) \ra \calD^{\an}_{\Lambda}(\Bun_{G}^{Z})$ is given by
    \[ 
		\begin{tikzcd}
			\calD_{\Lambda}^{\an}(\Bun_{G}^{Z'}) \arrow[r,hook,"j_{Z'Z!}",swap] & \calD_{\Lambda}^{\an}(\Bun_{G}^{Z}) \arrow[l, bend left=49,"j_{Z'Z}^{!}"] \arrow[l, bend right=49,"j^\flat_{Z'Z}",swap] \arrow[rr,"j_{UZ}^{*}"]  & & \calD^{\an}_{\Lambda}(\Bun_{G}^{U}), \arrow[ll,"j_{UZ\sharp}",bend right = 49,swap,hook] \arrow[ll,"j_{UZ!}",bend left = 49,hook] 
		\end{tikzcd}
    \]
    which will match up,  under our desired equivalence $\pitch: \Shv^{!}(\ICG) \xrightarrow{\simeq} \calD_{\Lambda}^{\an}(\Bun_{G})$,  with the standard recollement 
    \[ 
		\begin{tikzcd}
			\Shv^{!}(\ICG_{Z'}) \arrow[r,hook,"\ul{i_{ZZ'*}}"] & \Shv^{!}(\ICG_{Z}) \arrow[r,"\ul{i_{UZ}^{*}}"]  \arrow[l, bend right=49,"\ul{i_{ZZ'}^{*}}",swap] \arrow[l, bend left=49,"\ul{i_{ZZ'}^{!}}"]    & \Shv^{!}(\ICG_{U}) \arrow[l, bend left=49,"\ul{i_{UZ*}}"] \arrow[l, bend right=49,"\ul{i_{UZ!}}",swap] 
		\end{tikzcd}
    \]
    coming from the semi-orthogonal decomposition described in \S \ref{subsection semi on isoc}.
    \end{example}

    \subsubsection{$B(G)$-semi-orthogonal functors}
    We can now specialize \Cref{thm: equivalence criterion semi-orthogonal functor few criteria} and \Cref{prop: equivalence criterion semi-orthogonal functor many criteria}  to the context we will be interested in. 
    For any finite convex subset $V\subseteq B(G)$, we consider the functors
    \[\ul{i_{V!}}:\Shv^!(\ICG_V)\to \Shv^!(\ICG).\]
    and
    \[\ul{i_{V*}}:\Shv^!(\ICG_V)\to \Shv^!(\ICG).\]
    We let $(\ICG_V)_!$ and $(\ICG_V)_*$ denote respectively the essential image of these functors.
    If $V\subseteq W$ for another finite convex subset of $B(G)$, we denote by $(\ICG_{VW})_!$ and $(\ICG_{VW})_*$ the essential images of the functors $\ul{i_{VW!}}$ and $\ul{i_{VW*}}$. 
    We employ the same notational conventions for $(\Bun_G^V)_\sharp$ and $(\Bun_G^V)_!$ and $(\Bun_G^{VW})_\sharp$ and $(\Bun_G^{VW})_!$.

    \begin{proposition}
	    \label{prop: altexplicatedsemiorthogonaldecompositioncriterion2}
    Let $F: \Shv^{!}(\ICG) \ra \calD_{\Lambda}^{\an}(\Bun_{G})$ be a functor in $\LinCat_\Lambda$ satisfying the following:
    \begin{enumerate}
	    \item For any finite convex subset $V\subseteq B(G)$, we have that 
		    \[F((\ICG_V)_!)\subseteq (\Bun_G^V)_\sharp\] as subcategories of $\calD^\an_\Lambda(\Bun_G)$.
	    \item For any finite convex subset $V$, fix a presentation $V=Z\setminus Z'$ as the difference of two finite closed subsets $Z'\subset Z$.
If $b\in V$ is a closed point and $W=V\setminus \{b\}$ its open complement, we have that $F_V((\ICG_{WV})_*)\subseteq (\Bun_G^{WV})_!.$  
		    Here, $F_V$ denotes the functor induced by $F$ through the Verdier quotient identification as in \Cref{lem_switch_persp_functors}.  
	    \item For any $b\in B(G)$, the induced map 
		    \[F_b:(\ICG_b)_!\to (\Bun_G^b)_\sharp\] is an equivalence.  
    \end{enumerate}
    Then $F$ is an equivalence.
    Moreover, for every finite closed $Z\subseteq B(G)$, the functor $F$ promotes uniquely to an equivalence of presentable 2-functor formalisms 
    \[F_Z:\bbS_{\ICG_Z}\xrightarrow{\simeq} \bbS^{\on{ex}}_{\Bun_G^Z}.\]
    Writing \(F_Z\) for the restriction of \(F\) to \(Z\), we have \(F=\colim_{Z\subset\ICG}F_Z\), equivalently \(F=\lim_{Z\subset\ICG}F_Z\) and it promotes to an equivalence of presentable 2-functor formalisms
    \begin{equation*}
        F:\bbS_{\ICG}\xto{\simeq}\bbS^{\on{ex}}_{\Bun_G}.
    \end{equation*}
    \end{proposition}
    
    \begin{proof} 
    Write $B(G) = \cup_{Z \subset B(G)} Z$, where $Z$ ranges over finite closed subsets of $B(G)$ using \Cref{lem_finite_B(G)_order} as before.
This allows us to write 
    \[ \ICG = \colim_{Z \subset B(G)} \ICG_{Z}, \]
    where the transition maps are given by closed immersions, and 
    \[ \Bun_{G} = \colim_{Z \subset B(G)} \Bun_{G}^{Z}, \]
    where the transition maps are given by open immersions. 

    We see that by \Cref{prop: Shvshriektheoryisnice} 
    \[ \Shv^{!}(\ICG) \simeq \colim_{Z \subset B(G)} \Shv^{!}(\ICG_{Z}), \]
    where the transition maps are given by $\ul{i_{ZZ'*}}$ along the closed immersions $\ICG_Z\to \ICG_{Z'}$. 

    Similarly, 
    \[ \calD_{\Lambda}^{\an}(\Bun_{G}) \simeq \colim_{Z \subset B(G)} \calD_{\Lambda}^{\an}(\Bun_{G}^{Z}), \]
    where the transition maps are given by $j_{ZZ'!}$ along the open immersions $\Bun_G^Z\to \Bun_G^{Z'}$.
    Indeed, by construction of $\calD_{\Lambda}^{\an}(\Bun_{G})$, we have an isomorphism 
    \[ \calD_{\Lambda}^{\an}(\Bun_{G}) \simeq \lim_{Z\subseteq B(G)} \calD_{\Lambda}^{\an}(\Bun_{G}^{Z}), \]
    where the limit is given by $j_{ZZ'}^*: \calD^{\an}_{\Lambda}(\Bun_{G}^{Z'}) \ra \calD^{\an}_{\Lambda}(\Bun_{G}^{Z})$ for the open immersions $\Bun_G^{Z}\to \Bun_G^{Z'}$. 
    By \cite[Lemma~7.17~(i)]{scholze6ff}, we can rewrite this as a colimit over the left adjoint functors 
    \[ \colim_{Z \subset B(G)} \calD_{\Lambda}^{\an}(\Bun_{G}^{Z}) \]
    given by $!$-pushforward. 

    The above considerations and the first hypothesis allow us to write 
    \[F\simeq \colim_{Z\subseteq B(G)} F_Z.\]
    In particular, it suffices to show $F_Z$ is an equivalence.  
    But this follows directly from \Cref{prop: equivalence criterion semi-orthogonal functor many criteria}, then \Cref{thm: equivalence criterion semi-orthogonal functor few criteria} shows that \(F_Z\) is an equivalence of presentable 2-functor formalisms.
    
    For the last point, the colimit presentation is true by construction.
    By the formula for colimits in \(\LinCat_\Lambda\), we also obtain that \(F=\lim_{Z\subset B(G)}F_Z\).
    This presentation allows us to show that \(F\) is even an equivalence of presentable 2-functor formalisms.
    To wit, observe that \(\Fun(\Corr(\Convex_{B(G)},\all),\LinCat_\Lambda)\) is identified with the limit
        \begin{equation*}
		\lim_{Z\subset B(G)}\Fun(\Corr(\Convex_{B(G)_Z},\on{All}),\LinCat_\Lambda)
        \end{equation*}
	as \(\colim_{Z\subset B(G)}\Corr(\Convex_{B(G)_Z},\on{All})=\Corr(\Convex_{B(G)},\all)\), using that this is just exhausting \(\Corr(\Convex_{B(G)},\all)\) by full subcategories.
        The collection of \(F_Z\) for \(Z\subset B(G)\) defines an arrow in this limit and thus we see that \(F\) promotes to an equivalence of presentable 2-functor formalisms.
    \end{proof}

In practice, the inclusions in the first condition of \Cref{prop: altexplicatedsemiorthogonaldecompositioncriterion2} will be checked on individual strata, as follows:

 \begin{lemma}{\label{lemma: ReducetoSingleb}}
	    \label{extend filtered functor from singleton to arbitrary convex}
	    In the context of \Cref{prop: altexplicatedsemiorthogonaldecompositioncriterion2}, if $F:\Shv^!(\ICG)\to \calD^\an_\Lambda(\Bun_G)$ is a functor in $\LinCat_\La$ %
        such that inclusion \begin{equation}\label{eqn: filt_funct_left_orth}F((\ICG_V)_!)\subseteq (\Bun_G^V)_\sharp\end{equation} 
      holds when $V=\{b\}$ is a singleton, then it also holds more generally when $V$ is an arbitrary convex subset of $B(G)$.
    \end{lemma}
    \begin{proof}
	    For any convex $V\subseteq B(G)$, the category $(\ICG_V)_!$ is generated under colimits by objects in the full subcategories $\{(\ICG_b)_!\}_{b\in V}$.
	    In particular, if $F((\ICG_b)_!)\subseteq (\Bun^b_G)_\sharp$ for all $b\in V$, it follows that $F(\ICG_V)_!\subseteq (\Bun_G^V)_\sharp$, since $(\Bun^V_G)_\sharp$ contains the stratawise $\sharp$-included category $(\Bun^b_G)_\sharp$.
	    This finishes the proof that the first condition holds.
    \end{proof}

	\section{Kimberlites and henselianity}
	\label{sec: kimberlitecalculations}
	The purpose of this section is to recall the theory of kimberlites as developed in \cite{Gle24,Gle23}.
	This theory is an attempt to answer the question: \emph{What is the analogue of a formal scheme in the context of v-sheaves?}
	The upshot is that spatial kimberlites (just like formal schemes) satisfy a henselian property along their non-analytic locus. 
	This will play a key role in our computations.
	We approach the question in steps, 
	\begin{align*}
		\{\on{Small\, v-sheaves}\}& \supseteq \{\on{Specializing\,v-sheaves}\} \\
						& \supseteq \{\on{Prekimberlites}\} \\
						& \supseteq \{\on{Valuative\,Prekimberlites}\} \\
						& \supseteq \{\on{Kimberlites}\} \\
						& \supseteq \{\on{Locally\,Spatial\,Kimberlites}\}. 
	\end{align*}
	Each of these categories is a full subcategory of the category of small v-sheaves \cite[Definition 12.1]{Sch17} obtained by adding axioms at each stage.

	We recall some of the terminology following \cite{Gle24}. 
	Note that the analytification functor $\diamond$ described in \S \ref{subsec: threeanalytificationfunctors} admits a right adjoint functor 
	\begin{equation}\red:\AnStk_v \to \SchStk_v,\end{equation}
	compare with \cite[Definition 3.12]{Gle24}.	We write $X^\Red=(X^\red)^\diamond$ for simplicity, which carries a canonical counit map $X^\Red\to X$.
A map $Y\to X$ is {\it formally adic} in the sense of \cite[Definition 3.20]{Gle24} if it preserves reductions under pullback, i.e., the diagram 
	\begin{equation}	\begin{tikzcd}
			Y^\Red \arrow{r} \arrow{d}  & \arrow{d} Y \\
			X^\Red \arrow{r} & X 
		\end{tikzcd}
	\end{equation}
    is Cartesian.
This allows us to define {\it formally separated} v-sheaves as in \cite[Definition 3.27]{Gle24} by requiring their diagonal to be a formally adic closed immersion.
Since kimberlites are meant to be formal integral models of diamonds, their building blocks should be $\mathrm{Spd}(R^+)$ instead of $\mathrm{Spa}(R,R^+)$. 
	Following \cite[Definition 4.6]{Gle24}, we say that a map $f:\mathrm{Spa}(R,R^+)\to X$ of v-sheaves is formalizable (equiv., $X$ \textit{formalizes} $f$) if there exists a dashed arrow completing the commutative diagram below
	\begin{equation}
		\begin{tikzcd}
			\Spa(R,R^+)\arrow{r}{f} \arrow{d}  & X  \\
			\Spd R^+.\arrow[dashed]{ru}  
		\end{tikzcd}
	\end{equation}
	Any such arrow is called a formalization of $f$ and it is unique if $X$ is formally separated by \cite[Proposition 4.9]{Gle24}. 
	We say $X$ is \textit{v-formalizing} if, for any $f$ as above, there is a v-cover $g:\Spa(S,S^+)\to \Spa (R,R^+)$ such that $X$ formalizes $g\circ f$. 
	If $X$ is moreover formally separated then it is called \textit{specializing}, compare with \cite[Definition 4.11]{Gle24}. 
	Specializing v-sheaves carry a topological specialization map \cite[Proposition 4.14]{Gle24}
	\[\on{sp}:|X|\to |X^\red|\]
	and a v-sheaf theoretic specialization map (or Heuer specialization map) \cite[\S 4.4]{Gle24}
	\[\on{Sp}:X\to X^{\on{H}}.\]
	Here $X^{\on{H}}:=(X^\red)^{\diamond/\circ}$, where, for $X\in \SchStk_v$, we let $X^{\diamond/\circ}$ denote the v-sheafification of the functor 
	\[X^{(\diamond/\circ)_{\on{pre}}}(R,R^+)=\{\Spec (R^+/R^{\circ \circ})\to X\}.\]

	We may finally introduce the notion of a kimberlite, as in \cite[Definitions 4.15, 4.30 and 4.35]{Gle24}.

	\begin{definition}	 \label{definitionofkimberlite}
		A specializing v-sheaf $X$ is a \textit{prekimberlite} if $X^\red\in \SchStk_v$ is representable by a scheme, and $X^\Red\to X$ is a closed immersion.
		We set $X^{\mathrm{an}}:=X\setminus X^\Red$ and call it the \textit{analytic locus} of $X$. 
		A prekimberlite $X$ is \textit{valuative} if $\on{Sp}:X\to X^{\on{H}}$ is partially proper. 
		Finally, a valuative prekimberlite $X$ is called a {\it kimberlite} if, for all affine open subsets $\Spec A\subseteq X^\red$, the open $U:=\on{sp}^{-1}(\Spec A)\cap X^\an$ is a spatial diamond.
	\end{definition}

	We have the following key concepts \cite[Definition 4.18, 4.38]{Gle24}.

	\begin{definition}
		\label{defi:tubular}
		Given a prekimberlite $X$ and a locally closed immersion $U\hookrightarrow X^\red$, we let the \emph{formal neighborhood of $X$ along $U$} denote the v-sheaf
		\[\hat{X}_{/U}:=X\times_{X^H} U^{\diamond/\circ}.\]
		With the setup as above, we let the \emph{tubular neighborhood of $X$ along $U$} denote the v-sheaf
		\[{X}^\circledcirc_{/U}:=\hat{X}_{/U}\cap X^\an.\]
	\end{definition}

	\begin{remark}
		Reinterpreting \Cref{definitionofkimberlite} in terms of \Cref{defi:tubular}, a valuative prekimberlite $X$ is a kimberlite if and only if $X^\circledcirc_{/U}$ is a spatial diamond for all $U\subseteq X^\red$ open affine subsets.
	\end{remark}

	The specialization map for kimberlites is particularly nice \cite[Theorem 4.40]{Gle24}. 
	Unfortunately, it is not clear if the property of being a kimberlite is stable under some natural constructions like passage to formal neighborhoods under closed immersions (i.e., $\hat{X}_{/Z}$ for $Z\subseteq X^\red$ a closed immersion).
	For this reason, it is better to work with the stronger notion of \emph{locally spatial kimberlites} introduced in \cite{Gle23}.

	\begin{definition}
		\label{affinespatialdefi}
		We say that an affine kimberlite $X$ is \textit{spatial} if the following conditions hold.
		\begin{enumerate}
			\item $X$ formalizes geometric points. 
			\item There is a qcqs formally adic v-cover $f:Y\to X$, for $Y=\Spd(B,B)$ where $B$ is an $I$-adic ring for $I\subseteq B$ a finitely generated ideal.
		\end{enumerate}
		We say that a kimberlite is locally spatial, if $\hat{X}_{/U}$ is spatial for all $U\subseteq X^\red$ open affine subsets. 
	\end{definition}
	
	\begin{example}
		\label{basic-example-of-kimberlite}
		It is not hard to see that, for any $I$-adic ring $B$ with $I\subseteq B$ a finitely generated ideal, the v-sheaf $\Spd(B,B)$ is an affine spatial kimberlite. 
		Moreover, the category of spatial kimberlites is stable under finite limits (see \cite[Proposition 4.16]{Gle23}).
		In particular, $(\Spd \bbZ_p)^n$ is an example of a spatial kimberlite.  
	\end{example}

	Besides \Cref{basic-example-of-kimberlite} there are plenty of examples (see \cite[\S 5]{Gle23}) of locally spatial kimberlites that do not come directly from a formal scheme. 
	Important instances of this are the local models $\calM_{\calG,\mu}$ studied in \cite{AGLR22} (see \cite[Theorem 5.7]{Gle23}), and, as we will discuss later, the integral models of moduli spaces of local shtukas are also locally spatial kimberlites (see \cite[Theorem 7.1.(3)]{Gle23}).

For our purposes, the most important property of locally spatial kimberlites is that they are henselian along their non-analytic locus. 	
Let us make this precise.
	
	Let $X$ be a spatial kimberlite let $Y=X^{\on{an}}$, $Z=X^\Red$ and denote by $j:Y\to X$ and $i:Z\to X$ the complementary open and closed immersions.
	We let $\pi_X:X\to \ast$ denote the structure morphism to $\ast$, similarly for $\pi_Y$ and $\pi_Z$.
	One of us (G.) proved the following statement that shows that $X$ satisfies two henselianity properties along $Z\to X$.
    \begin{theorem}[{\cite[Theorem 1.4]{Gle23}}]
		\label{nearbycyclestheorem}
		Let the notation be as above, and assume that the following hold
		\begin{enumerate}
			\item $[\pi^\red_X:X^\red\to \ast]\in \PSch_{P_\pfp}$ (i.e., it is locally of perfectly finite presentation over $\Spec k$ and proper). 
		\item $\on{dim.trg}\pi_Y<\infty$ (see \cite[Definition 21.7]{Sch17}).
	\end{enumerate}
	Then, for all $\calF\in \calD_{\acute{e}t}(Y)$ the following vanishing statements hold:
		\begin{enumerate}
			\item $\Gamma(X,j_!\calF)=0$. 
			\item $\Gamma_c(X,j_*\calF)=0$.
		\end{enumerate}
	\end{theorem}

	We wish to apply a version of \Cref{nearbycyclestheorem} to the map
	\[\Sht^\an_{\calG,\leq w,b}:=\Sht^\an_{\calG,\leq w}\times_{\Bun_G}\Bun^b_G\to \Bun_G^b.\]
	
	For this, we will need the following version whose proof is almost identical, but with the additional complexity of carrying a $K$-action.
	
	\begin{theorem}
		\label{hyperbolic localization for some quotients}
		Let $X$ be a spatial kimberlite, let $Y=X^\an$ and $Z=X^\Red$, with immersions $j:Y\to X$ and $i:Z\to X$ respectively.
		Let $K$ be a profinite group of pro-order coprime to $\ell$ acting on $X$ such that the action on $Y$ is free.
		Let $\pi_{X,K}:X/\underline{K}\to [\ast/\underline{K}]$ denote the structure map and denote the natural immersions 
		$j_K:Y/\underline{K}\to X/\underline{K}$ and $i_K:Z/\underline{K}\to X/\underline{K}$ respectively.
		Suppose the following hold.
		\begin{enumerate}
			\item The reduced special fiber $X^\red$ is pfp proper over $\ast$.
			\item One can find a free action of $\underline{K}$ on a spatial diamond $\frakD$ such that the structure map $\frakD\to \ast$ is fdcs and $\calD_{\Lambda}^{\an}$-smooth.
		\end{enumerate}
		Then, for all $A\in \calD_\et(Y/\underline{K})$, we have the following vanishing statement
		\[\pi_{X,K,!}j_{K,*}A\simeq 0.\]
	\end{theorem}
    \begin{remark}{\label{rem: trickeries-with-quotients well-defined}}
    Let $K$ be a profinite group that admits a closed embedding into $\GL_{n}(E)$ for some $n$.
    Without loss of generality, we may assume that if factors through $\GL_n(O_E)$ by fixing a lattice stabilized by $K$.
We consider a space that we denote $\GL_{n,E}^{\diamond}$, which is given by sending a perfectoid Huber pair $(R,R^{+}) \in \Perff$ to the isomorphism classes of pairs $((R^{\sharp},R^{\sharp+}),i)$ of untilts of $(R,R^{+})$ over $E$ together with a $R^{\sharp+}$-point of $\GL_{n,O_E}$.
It is not hard to verify that $\GL_{n,E}^{\diamond}$ is a spatial diamond and that its structure map to $\Spd E$ is fdcss and $\calD_{\Lambda}^{\an}$-smooth.
Moreover, the action of $\underline{\GL_{n}(O_E)}$ (and consequently of \ul{K}) on $\GL_{n,E}^{\diamond}$ is free.
In particular, we see that assumption (2) of \Cref{hyperbolic localization for some quotients} is satisfied in this case (see \cite[Example~IV.1.9 (iv)]{FS21}).
     \end{remark}

	\begin{remark}
		 Consider the following Cartesian diagram.
		\begin{center}
			\begin{tikzcd}
				X \arrow{r} \arrow{d}  & \ast \arrow{d} \\
				X/\underline{K}\arrow{r} & {[\ast/\underline{K}]} 
			\end{tikzcd}
		\end{center}
		If the map $\ast \to [\ast/\underline{K}]$ was smooth \Cref{hyperbolic localization for some quotients} could easily be reduced to \Cref{nearbycyclestheorem} using smooth and proper base change.
		Since this map is not smooth, we have to work around this issue. 
			\end{remark}
	
	\begin{proof}
		Let $X_K$, $Y_K$ and $Z_K$ denote the quotients of $X$, $Y$ and $Z$ by $K$ respectively.
		We wish to show that $\pi_{X,K,!}j_{K,*}:\calD^\an_\Lambda(Y)\to \calD^\an_\Lambda([\ast/\underline{K}])$ is the $0$-functor. 
		
		Let $S=\frakD/\underline{K}$ and $s:S\to [\ast/\underline{K}]$ be the structure map.
		Observe that $S$ is a locally spatial diamond $\calD_{\Lambda}^{\an}$-smooth over $\ast$ by \cite[Proposition~24.2]{Sch17}.
As $\frakD$ is qcqs, $S = \frakD/K$ is also qcqs.
Being qcqs and locally spatial, $S$ is spatial. 
        
		Furthermore, $s$ is also $\calD_{\Lambda}^{\an}$-smooth and $s^\ast$ is conservative. 
		It suffices to show that $s^*\pi_{X,K,!}j_{K,*}$ is the $0$-functor. 
		Let $p_{Y_K}:Y_K\times_{[\ast/\underline{K}]}S\to Y_K$ denote the base change of $s$ along $Y_K\to [\ast/\underline{K}]$.
		By smooth and proper base change, it suffices to show that $\pi_{s,K,!}\circ j_{s,K,*}p^*_{Y_K}\cong 0$ for the natural maps $j_{s,K}:Y_K\times_{[\ast/\underline{K}]} S\subseteq X_K\times_{[\ast/\underline{K}]} S$ and $\pi_{s,K}:X\times_{[\ast/\underline{K}]} S\to S$.
		We also consider the maps $j_s:Y\times \frakD\to X\times \frakD$ and $\pi_s:X\times \frakD\to \frakD$.
		Note that, by \Cref{nearbycyclestheorem}, we also have that $\pi_{s,!}\circ j_{s,*}p^*_{Y_K}\simeq 0$.
		
		By hypothesis, $Y$ is a spatial diamond since $X$ is a spatial kimberlite, and since the action of $\underline{K}$ on $Y$ is free, $Y_K$ is also a spatial diamond by the same reasoning as for $S$ above.
		Also, $Y\times \frakD$ is a locally spatial diamond and, since $Y_K\times_{[\ast/\underline{K}]} S$ is a quotient of $Y\times \frakD$ by  $\underline{K}$, we also get that $Y_K\times_{[\ast/\underline{K}]} S$ is a locally spatial diamond.
		Note that $Y\times \frakD\to Y_K\times_{[\ast/\underline{K}]} S$ is universally open.
		
		Let $T\to Y$ and $Q\to \frakD$ be universally open quasi-pro\'etale covers with $T=\Spa( R,R^+)$ and $Q=\Spa (L,L^+)$ totally disconnected perfectoid spaces. 
		Being spatial diamonds $Y$, $\frakD$, they are qcqs; also $T$ and $Q$ are qcqs.
As any map between qcqs objects is qcqs, the maps $T\to Y$ and $Q\to \frakD$ are automatically qcqs.
        Fix a $\varpi_R\in R$ and $\varpi_L\in L$ pseudo-uniformizers. 
		There is a continuous map $\kappa:|(\Spd R^+\times \Spd L^+)^{\on{an}}|\to [0,\infty]$ meassuring the relative value of $\varpi_L$ against $\varpi_R$ and for every interval $I\subseteq [0,\infty]$ we let $U_{T,Q,I}$ denote the open subspace of $(\Spd R^+\times \Spd L^+)^{\on{an}}$, associated with the interior of $\kappa^{-1}(I)$.
		For example, \[T\times \Spd L^+=U_{T,Q,[0,\infty)}=(\Spd R^+\times \Spd L^+)_{\varpi_R\neq 0}.\] 
		Now, for all $I \subseteq (0,\infty)$,
		we let $U^{{K}}_I\subseteq Y_K\times_{[\ast/\underline{K}]} S$ denote the image of $U_{T,Q,I}\to Y_K\times_{[\ast/\underline{K}]} S$.
		We let $U_I\subseteq Y\times \frakD$ denote the preimage of $U^{{K}}_I$ in $Y\times \frakD$.
		We note that if $I\subseteq (0,\infty)$ is compact then $U_I$ and $U^K_I$ are spatial diamonds.

		We claim that $|U^K_{[a,\infty)}|\cup |Z\times_{[\ast/\underline{K}]} S|\subseteq |X\times_{[\ast/\underline{K}]} S|$ defines a quasicompact open subset and in particular it gives rise to a spatial diamond that we denote $U^K_{[a,\infty]}\subseteq X\times_{[\ast/\underline{K}]} S$ with preimage a spatial diamond $U_{[a,\infty]}\subseteq Y\times \frakD$.
		We argue as follows, we have a map of spatial diamonds
		\[U_{T,Q,[0,\infty]}=(\Spd R^+\times \Spd L^+)^{\on{an}} \to (X\times \Spd L^+)^{\on{an}}\]
		coming from the formally adic map of spatial kimberlites $\Spd R^+\times \Spd L^+ \to X\times \Spd L^+$. 
		We restrict this map to obtain a sequence of qcqs maps of locally spatial diamonds
		\[U_{T,Q,(0,\infty]}=(\Spd R^+\times Q) \to X\times Q\to X\times \frakD \to X_K\times_{[\ast/\underline{K}]} S\]

		Let $f$ denote this composition.
		In general, $f:|U_{T,Q,(0,\infty]}| \to |X_K\times_{[\ast/\underline{K}]} S|$ is not necessarily an open map, but it is still a quotient map since it is qcqs surjective.
		The locus $|U^K_{[a,\infty)}|\cup |Z_K\times_{[\ast/\underline{K}]} S|$ corresponds precisely to the image of $U_{T,Q,[a,\infty]}$ which is quasicompact. 
		Moreover, we have that $f^{-1}(f(U_{T,Q,[a,\infty]}))=|U_{T,Q,[a,\infty]}|\cup f^{-1}(f(U_{T,Q,[a,\infty)})$ which is open. 
		This implies that $U^K_{[a,\infty]}\subseteq X_K\times_{[\ast/\underline{K}]} S$ is open quasicompact.
		Since $U_{[a,\infty]}\subseteq X\times \frakD$ is simply the preimage of $U^K_{[a,\infty]}$ it is also open quasicompact.

		In what follows we show that $\pi_{s,K,!}\circ j_{s,K,*}$ can be expressed in terms of the $U^K_{[a,\infty]}$ and $U_{[a,\infty]}$ to ease notation we only do one case the other one being analogous. 
		Let $\overline{U}^K_{[a,\infty]}$ denote the closure of $U^K_{[a,\infty]}$ inside of $U^K_{(0,\infty]}=X_K\times_{[\ast/\underline{K}]}S$.
		We define several maps.
		
		Let $k_{a,\infty}:U^K_{[a,\infty]}\to \overline{U}^K_{[a,\infty]}$ and $\overline{k}_{a,\infty}:\overline{U}^K_{[a,\infty]}\to U^K_{(0,\infty]}$ denote the open and closed immersions respectively.
		We let $k_{a}:U^K_{[a,\infty)}\to \overline{U}^K_{[a,\infty)}$ and $\overline{k}_{a}:\overline{U}^K_{[a,\infty)}\to U^K_{(0,\infty)}$ denote the open and closed immersions respectively obtained by restriction along ${U}^K_{(0,\infty)}\subseteq {U}^K_{(0,\infty]}$.
		We also let $j_a:U^K_{[a,\infty)}\to U^K_{[a,\infty]}$ and $\overline{j}_a:\overline{U}^K_{[a,\infty)}\to \overline{U}^K_{[a,\infty]}$ denote the open immersion.
		For the convenience of the reader we organize it in a commutative diagrams.
		
		\begin{center}
			\begin{tikzcd}
				U^K_{[a,\infty]} \arrow{r}{k_{a,\infty}}   & \overline{U}^K_{[a,\infty]} \arrow{r}{\overline{k}_{a,\infty}} & U^K_{(0,\infty]}\ar{r}{\pi_{s,K}}=X_K\times_{[\ast/\underline{K}]}S & S \\
				U^K_{[a,\infty)} \arrow{r}{k_a} \arrow{u}{j_a} &  \overline{U}^K_{[a,\infty)} \ar{r}{\overline{k}_a}\ar{u}{\overline{j}_a} & U^K_{(0,\infty)}=Y_K\times_{[\ast/\underline{K}]}S  \ar{u}{j_{s,K}}
			\end{tikzcd}
		\end{center}

		For any $A\in \calD^\an_\Lambda(X_K\times_{[\ast/\underline{K}]} S)$, we may compute $\pi_{s,K,!}A \in \calD^\an_\Lambda(S)$ as 
		\begin{equation}
			\label{anotherhenseliancomputation}
			\varinjlim \limits_{a \to 0}\pi_{s,K,*}\overline{k}_{a,\infty,*}k_{a,\infty,!}A_{|_{U^K_{[a,\infty]}}}.
		\end{equation}
		
		Indeed, this follows from the fact that $X_K\times_{[\ast/\underline{K}]} S\to S$ is a map of locally spatial diamond, from the fact that the map $\pi_{s,K}\circ \overline{k}_{a,\infty}:\overline{U}_{[a,\infty]}\to S$ is proper, and from the fact that the family $U^K_{[a,\infty]}$ is cofinal among quasicompact open subset of $U^K_{(0,\infty]}=X_K\times_{[\ast/\underline{K}]} S$ \cite[Definition 22.13, Definition 22.4]{Sch17}.

		When $A=j_{s,K,*}B$ we may rewrite this as 
		\begin{align}
			\label{henselian-computation-beg}
			\varinjlim \limits_{a \to 0}\pi_{s,K,*}\overline{k}_{a,\infty,*}k_{a,\infty,!} A_{|_{U^K_{[a,\infty]}}} 
			&\cong \varinjlim \limits_{a \to 0}\pi_{s,K,*}\overline{k}_{a,\infty,*}k_{a,\infty,!}j_{a,*}B_{|_{U^K_{[a,\infty)}}} \\
			&\cong \varinjlim \limits_{a \to 0}\pi_{s,K,*}\overline{k}_{a,\infty,*} \overline{j}_{a,*}k_{a,!}B_{|_{U^K_{[a,\infty)}}} \\
			\label{henselian-computation-end}
			&\cong \varinjlim \limits_{a \to 0}\pi_{s,K,*} \overline{j}_{s,K,*}\overline{k}_{a,*} k_{a,!}B_{|_{U^K_{[a,\infty)}}} 
		\end{align}
		Here the only subtle step of the computation is to justify 
		\[k_{a,\infty,!}j_{a,*}B_{|_{U^K_{[a,\infty)}}}\cong   \overline{j}_{a,*}k_{a,!}B_{|_{U^K_{[a,\infty)}}}.\] We argue as follows.
		Note that $i_{a,s}^! k_{a,\infty,!}j_{a,*}B_{|_{U^K_{[a,\infty)}}}\cong 0$ for the inclusion $i_{a,s}:Z_K\times_{[\ast/\underline{K}]} S\to \overline{U}^K_{[a,\infty]}$.
		Indeed, we have a factorization $i_{a,s}:Z_K\times_{[\ast/\underline{K}]} S\to U^K_{[a,\infty]}\to \overline{U}^K_{[a,\infty]}$.
		Applying excision to $k_{a,\infty,!}j_{a,*}B_{|_{U_{[a,\infty)}}}$, we obtain the desired isomorphism. 
		
		We let $\beta^K=\pi_{s,K}j_{s,K}$ and $t^K_a=\overline{k}_a\circ k_a$
		\[\beta^K_a:U^K_{(0,\infty)}=Y_K\times_{[\ast/\underline{K}]}S \to S\]
		\[t^K_a:U^K_{[a,\infty)}\to Y_K\times_{[\ast/\underline{K}]}S\]
		and define a functor 
		\[\beta^K_{!,+}:\calD^\an_\Lambda(Y_K\times_{[\ast/\underline{K}]}S)\to \calD^\an_\Lambda(S)\]
		with formula
		\[\beta^K_{!,+}B=\varinjlim_a \beta^K_*t^K_{a,!}B\]
		in analogy to the functors of constructed in \cite[Definition IV.5.2]{FS21}.
		Note that our setup is slightly different so $\beta^K_{!,+}$ is not literally one of the functors considered in loc. cit.
		The computation \eqref{henselian-computation-beg}-\eqref{henselian-computation-end} and \Cref{anotherhenseliancomputation} show that,
		for all $B\in \calD^\an_\Lambda(Y_K\times_{[\ast/\underline{K}]}S)$, we have an isomorphism
		\[\pi_{s,K,!}j_{s,K,*}B\simeq \beta^K_{!,+}B.\]
		
		Let us go back to studying $\pi_s:X\times \frakD\to \frakD$. 
		We define maps 
		\[\beta:Y\times \frakD \to \frakD\]
		and maps
		\[t_a:U_{[a,\infty)}\to Y\times \frakD\]
		by pulling back $\beta^K$ and $t^K_a$ along $\frakD\to S.$
		We consider the functor 
		
		\[\beta_{!,+}:\calD^\an_\Lambda(Y\times \frakD)\to \calD^\an_\Lambda(\frakD)\]
		with formula
		\[\beta_{!,+}B=\varinjlim_a \beta_*t_{a,!}B.\]
		This agrees with the functor constructed in \cite[Definition IV.5.2]{FS21}.
		In a completely analogous way with identical proof we obtain a formula 
		\[\pi_{s,!}j_{s,*}B\simeq \beta_{!,+}B.\]
		Moreover, since $Y$ is a spatial diamond that is fdcs and partially proper over $\ast$ 
		by \cite[Theorem IV.5.3]{FS21} $\beta_{!,+}B$ vanishes whenever $B=p_Y^*A$ for $p_Y:Y\times \frakD\to Y$ the projection map.
		In other words, $\pi_{s,!}j_{s,*}p_Y^*A\simeq 0$ for any $A\in \calD^\an_\Lambda(Y)$. 
		In contrast, our goal is to show that $\pi_{s,K,!}j_{s,K,*}p^{*}_{Y_K} A\simeq 0$ for and any $A\in \calD^\an_\Lambda(Y_K)$. 
		
		We consider the pair of Cartesian diagrams
		\begin{center}
			\begin{tikzcd}
				Y\times \frakD \arrow{r}{\gamma_Y} \arrow{d}{p_Y}  & Y_K\times_{[\ast/\underline{K}]} S \arrow{d}{p_{Y_K}} \\
				Y \arrow{r}{q_Y} & Y_K 
			\end{tikzcd}
		\end{center}
		and 
		\begin{center}
			\begin{tikzcd}
				U_{[a,\infty)} \ar{r}{t_a} \ar{d}{\gamma_a}& Y\times \frakD \arrow{r}{\beta} \arrow{d}{\gamma_Y}  & \frakD \arrow{d}{\gamma} \\
				U^K_{[a,\infty)} \ar{r}{t^K_a}	& Y_K\times_{[\ast/\underline{K}]}S \arrow{r}{\beta^K} & S 
			\end{tikzcd}
		\end{center}
		By quasicompact and proper base change it follows that
		\[\gamma^*\beta^K_{+,!}\simeq \beta_{+,!}\gamma_Y^*.\]
		Moreover, the functor $\gamma^*$ is conservative.
		We can finish the proof by computing 
		\begin{align*}
			\gamma^*\pi_{s,K,!}j_{s,K,*}p^{*}_{Y_K} A &\simeq \gamma^* \beta^K_{!,+} p^{*}_{Y_K} A \\
			& \simeq \beta_{!,+}\gamma_Y^* p^{*}_{Y_K} A \\
			& \simeq \beta_{!,+}p_Y^* q_Y^* A \\
			&\simeq 0 
		\end{align*}
	\end{proof}

	In the next section, we will exploit the fact that $\Sht_{\calG,\leq w}(b):=\Sht^\an_{\calG,\leq w,b}\times_{\Bun_G^b}\ast$ is a locally spatial kimberlite to study $\pitch$.
    Here we recall that the map
    \[\Sht^\an_{\calG,\leq w}\to \Bun_G\]
    is given as the composition of the map $\Nt^\dagger$ from \eqref{important Cartesian} and from $\sigma$ from \eqref{BungMer-triangle}. 
	As we discussed above, locally spatial kimberlites are special instances of prekimberlites (Definition \ref{defi:tubular}), and in particular we can define a reduced locus 
	\[(\Sht_{\calG,\leq w}(b))^\Red\subseteq \Sht_{\calG,\leq w}(b), \]
	which is a closed immersion 
	\[(\Sht_{\calG,\leq w}(b))^\Red\simeq [(\Sht_{\calG,\leq w}(b))^\red]^\diamond.\]
	Now, $(\Sht_{\calG,\leq w}(b))^\red$ has a more familiar expression.  
	Indeed, recall from \cite[Proposition 2.61]{Gle21} that the following holds.

	\begin{proposition}
		\label{adlv-formula}
		There is a canonical identification 
		\[(\Sht_{\calG,\leq w}(b))^\red\simeq X_{\leq w}(b)\]
		where $X_{\leq w}(b)$ denotes the affine Deligne--Lusztig variety (i.e $X_{\leq w}(b):= \Sht_{\calG,\leq w} \times_{\ICG} \ast$, where $\ast \ra \ICG$ is given by $b$).
		In particular, $(\Sht_{\calG,\leq w}(b))^\red$ is a locally perfectly finitely presented scheme over $\Spec k$, and each of its irreducible components are pfp proper over $\Spec k$ (see \cite[Theorem 1.2]{MR4148693} and \cite[Proposition 2.61]{Gle21}). 
	\end{proposition}

	\begin{remark}
	Strictly speaking the reference \cite[Proposition 2.61]{Gle21} has a different setup and only shows the formula $(\Sht_{\calG,\leq \mu}(b))^\red\simeq X_{\leq \mu}(b)$ where the bound ranges over dominant cocharacters.
	Nevertheless, the argument goes through in our setup with only superficial modifications.
	\end{remark}

	Recall the generic Newton polygon map $\gamma$ of \eqref{BungMer-triangle} and \eqref{important Cartesian}. 
	As it turns out, one can use $\gamma$ to understand the closed subsheaf
	\[(\Sht_{\calG,\leq w}(b))^\Red\subseteq \Sht_{\calG,\leq w}(b).\]
	Indeed, we have the following.
	\begin{proposition}
		\label{meromorphicity gives closed immersion}
		The reduced locus 
		\[(\Sht_{\calG,\leq w}(b))^\Red\]
		is always contained in the pullback of
        \[(\gamma \circ \on{Nt}^\dagger)^{-1}(\ICG^\diamondsuit_b)\]
		along $* \ra \Bun_{G}^{b}$, as closed subsheaves of $\Sht_{\calG,\leq w}(b)$. 
	\end{proposition}
	\begin{proof}
		Recall that $(R,R^+)$ points of $\Sht_{\calG,\leq w}(b)$ are equivalence classes of data $(\calE,\Phi, \rho)$.
		Here $\Phi:\varphi^*\calE\to \calE$ is a shtuka defined over $\calY^R_{[0,\infty)}$ with leg at $\pi=0$ and meromorphy bounded by $w$, and $\rho$ is a $\varphi$-equivariant isomorphism $\calE\simeq \calE_b$ defined over $\calY^R_{(0,\infty)}$ (see \cite[Definition 2.45]{Gle21}). 
		The proof of \cite[Proposition 2.61]{Gle21} shows that $(\Sht_{\calG,\leq w}(b))^\Red$ is the locus in which $\rho$ is meromorphic along $\pi=0$.
		Since $\rho$ is meromorphic, one can consider the restriction of $\rho$ to a formal neighborhood of $\pi=0$ which gives an isomorphism 
		\[(\calE_{\bbW(R)},\Phi_{\bbW(R)[\frac{1}{\pi}]})\simeq_{\rho} (\calE_{b,\bbW(R)},\Phi_{b,\bbW(R)}[\frac{1}{\pi}]).\]
		In other words, $\rho_{\bbW(R)[\frac{1}{\pi}]}$ exhibits that $\gamma \circ \on{Nt}^\dagger(\calE,\Phi)$ is isomorphic to $(\calE_b,\Phi_b)$ in $\ICG^\diamondsuit(R,R^+)$, or equivalently, $\gamma \circ \on{Nt}^\dagger(\calE,\Phi)\in \ICG_b^\diamondsuit(R,R^+)$ as we wanted to see.
	\end{proof}

	For $\beta \in B(G)$, we let 
	\[\Sht^{\gamma=\beta}_{\calG,\leq w} \subseteq \Sht_{\calG,\leq w}\]
	denote the locus of the form $(\gamma \circ \on{Nt}^\dagger)^{-1}(\ICG^\diamondsuit_\beta)$, and after pullback along $* \ra \Bun_{G}^{b}$  we denote it as
	\[\Sht^{\gamma=\beta}_{\calG,\leq w}(b) \subseteq \Sht_{\calG,\leq w}(b).\]
Here the maps are as in (\ref{important Cartesian}).
	The following consequence will be crucial for our study of $\pitch$ in the next section.

	\begin{corollary}
		\label{Sht gamma = beta contained in analytic locus}
		Let $b,\beta\in B(G)$ with $b\neq \beta$.
		Then $\Sht^{\gamma=\beta}_{\calG,\leq w}(b)$ is contained in the analytic locus of $\Sht_{\calG,\leq w}(b)$.
In other words,
		\[\Sht^{\gamma=\beta}_{\calG,\leq w}(b)\subseteq (\Sht_{\calG,\leq w}(b))^\an.\]
	\end{corollary}
	\begin{proof}
		This follows from \Cref{meromorphicity gives closed immersion}. 
		Indeed, it suffices to show that \[\Sht^{\gamma=\beta}_{\calG,\leq w}(b)\cap \Sht_{\calG,\leq w}(b)^\Red\] is empty.
		But since $b\neq \beta$, $\Sht^{\gamma=\beta}_{\calG,\leq w}(b)\cap \Sht^{\gamma=b}_{\calG,\leq w}(b)$ is already empty.
	\end{proof}

	Since $\Sht_{\calG,\leq w}(b)$ is only a locally spatial kimberlite and not a spatial kimberlite, one cannot use \Cref{hyperbolic localization for some quotients} directly.   
	To address this, we formulate the following statement.

	\begin{proposition}
		\label{colimit presentation of Sht G leq w with required properties}
		Let $K\subseteq G_b(E)$, be a pro-$p$ open compact subgroup.
		Then $\Sht_{\calG,\leq w}(b)$ admits a colimit presentation of the form
		\[\Sht_{\calG,\leq w}(b)\simeq \colim_{i\in I} X_i\]
		such that the following statements hold:
		\begin{enumerate}
			\item The index category is cofiltered.
			\item Each $X_i$ is a spatial kimberlite with proper reduced special fiber.
			\item The maps $X_i\to \Sht_{\calG,\leq w}(b)$ and the transitions maps $X_i\to X_j$ are open immersions.
			\item Each $X_i$ is $K$-stable.
		\end{enumerate}
	\end{proposition}
	\begin{proof}
		By \Cref{adlv-formula}, $\Sht_{\calG,\leq w}(b)^\red$ is the affine Deligne--Lusztig variety $X_{\leq w}(b)$. 
		Let $\calC$ denote the set of $K$-orbits of irreducible components of $X_{\leq w}(b)$.
		Our index set $I$ is simply the set of finite subsets of $\calC$, and for every $i\in I$ of the form $i=\{c_1,\dots,c_n\}$ we form the closed subset
		\[Z_i=\cup_{j=1}^n c_i\subseteq X_{\leq w}(b).\]
		This is a perfectly finitely presented proper scheme over $\Spec k$.
		We let $X_i:=\widehat{\Sht_{\calG,\leq w}(b)}_{/Z_i}$.
		We now verify all of our claims.
		The indexing category is cofiltered by construction. 
		That each $X_i$ is a spatial kimberlite follows from \cite[Proposition 4.12]{Gle23}, the definition \cite[Definition 4.4]{Gle23} and the fact that $Z_i$ is qcqs since it is a pfp proper scheme over $\Spec k$.
		That $X_i\to \Sht_{\calG,\leq w}(b)$ is an open immersion follows from \cite[Proposition 4.22]{Gle24}.
		That each $X_i$ is $K$-stable follows from the functoriality of the specialization map \cite[Proposition 4.14]{Gle24}.
	\end{proof}

	\section{Proof of the equivalence}
	\label{sec: TheEquivalence}
	We are finally ready to show that the functor $\pitch$ of \S \ref{sec:construction_functor} is an equivalence.
	We will do this by verifying $F=\pitch$ satisfies the assumptions of \Cref{prop: altexplicatedsemiorthogonaldecompositioncriterion2} with respect to the semi-orthogonal decomposition $\bb{S}_{\ICG}$ on $\Shv^{!}(\ICG)$ and  $\bb{S}_{\Bun_{G}}^{\mathrm{ex}}$ on $\calD_{\Lambda}^{\an}(\Bun_{G})$ defined in \S \ref{subsection semi on isoc} and \S \ref{sec:semiorthogonalDecOfDBunG}, respectively.
    
	The first step of the proof is to show that condition \Cref{prop: altexplicatedsemiorthogonaldecompositioncriterion2} (1) holds for the functor $\pitch$, see \Cref{extend filtered functor from singleton to arbitrary convex,left orthogonality is shown}. 
	The calculation itself relies on the behavior of the Fargues--Scholze charts $\calM_b$ and the locally closed substack $\Sht_{\calI,w_b}$ of the shtuka stack for $\sigma$-straight elements $w_b$.
	The calculation is also key to show that $\pitch$ restricts to an equivalence strata by strata, see \Cref{pitch equivalence crit (3),diagram spell pitch out}. 
    
    The next step will be to consider analogues 
    \[\pitch_V:\Shv^!(\ICG_V)\to \calD^\an_\Lambda(\Bun_G^V),\] 
    of our functor $\pitch$ that range over finite convex subsets $V\subseteq B(G)$. 
    We use the geometry of $\Bun_G^\mer$ to study the behavior of these functors. 
    Moreover, if $V=Z_1\setminus Z_2$ with $Z_1,Z_2\subseteq B(G)$ for two finite closed subsets, then the functor $\pitch_V$ will give a computable expression for the Verdier quotient map
    \[\Shv^!(\ICG_{Z_1})_!/\Shv^!(\ICG_{Z_2})_! \to  \calD^\an_\Lambda(\Bun_G^{Z_1})_\sharp/\calD^\an_\Lambda(\Bun_G^{Z_1})_\sharp.\] 
    Using our criterion \Cref{prop: altexplicatedsemiorthogonaldecompositioncriterion2}, we are left to show that condition \Cref{prop: altexplicatedsemiorthogonaldecompositioncriterion2} (2) holds for $V\setminus \{b\}$ in $V$ for a closed point $b\in V$ and an arbitrary finite convex $V\subseteq B(G)$. 
    This last orthogonality, is obtained using the henselianity for the locally spatial kimberlite $\Sht_{\calG,\leq w_b}(b)$ to control the removed locus, see \Cref{right orthogonality verified}.

	\subsection{Left semi-orthogonality and $\Bun_G^\mer$}
	By construction $\pitch$ (or more precisely $\pitch^\vee$) arises from the sind-$\dagger$-correspondence of the form 
	\begin{center}
	\begin{tikzcd}
		\ICG^{\dagger} \simeq \Bun_G^\mer \arrow{r}{\sigma} \arrow{d}{\gamma}  & \Bun_G \\
	 \ICG^\diamondsuit, & 
	\end{tikzcd}
	\end{center}
    described in \Cref{cor: ConstructionofMainSindDaggerCorrespondence}.

	A key geometric input is the following result that two of us proved with Zillinger.
	Let $\calM$ denote the moduli of filtered bundles with $G$-structure introduced in \cite[Definition V.3.2]{FS21}.  
	We colloquially refer to it as the Fargues--Scholze chart of $\Bun_G$, since there is a fdcs smooth surjective map \cite[Theorem V.3.7]{FS21}
	\[\pi:\calM\to \Bun_G.\]
	Recall that $\calM$ admits a decomposition into connected components indexed by $b\in B(G)$, and that we have a map
	\[\bigsqcup_{b\in B(G)} \gamma_b:\bigsqcup_{b\in B(G)} \calM_b\to \bigsqcup_{b\in B(G)} [\ast/\ul{G_b(E)}]\] 
	that sends a filtered bundle with $G$-structure to the corresponding semi-stable graded bundle with $G_b$-structure (see discussion below \cite[Definition V.3.2]{FS21}).
The corresponding diagram 
    \[ 
    \begin{tikzcd}
    \mathcal{M}_{b} \arrow[r,"\sigma_{b}"] \arrow[d,"\gamma_{b}"] & \Bun_{G} & \\
    \text{[}\ast/\underline{G_{b}(E)}\text{]} & &
    \end{tikzcd}
    \]
    is  precisely the diagram described in \S\ref{ss: theanalyticcategory}.
One has an isomorphism $\ICG_b^\diamondsuit\simeq [\ast/\ul{G_b(E)}]$, by \cite[Proposition~2.20]{GIZ25} and  \Cref{thm:isocG_geometry}.
This allows us to make sense of the following. 
	\begin{theorem}
		[{\cite[Theorem 1.1]{GIZ25}}]\label{thm: GIZ}	\label{cartesian isoc bungmer}
		We have a commutative diagram  
	\begin{center}
	\begin{tikzcd}
		\calM \arrow{r} \arrow{d}{q} \arrow[bend left]{rr}{\pi} & \Bun_G^\mer \arrow{d}{\gamma} \ar{r}{\sigma} & \Bun_G\,.  \\
	 \bigsqcup_{b\in B(G)} \ICG_b^\diamondsuit \arrow{r} & \ICG^\diamondsuit,
	\end{tikzcd}
	\end{center}
    where the square is Cartesian.
	\end{theorem}

	This has the following consequence for our functor.
	\begin{proposition}
		\label{left orthogonality is shown}
		The following statements hold.
		\begin{enumerate}
			\item For all $b\in B(G)$
				\[\pitch((\ICG_b)_!)\subseteq (\Bun_G^b)_\sharp.\]

			\item For all $Z\subseteq B(G)$, a finite closed subset,
				\[		\pitch(\Shv^!(\ICG_Z)_{*})\subseteq \calD^\an_\Lambda(\Bun_G^Z)_{!} \]
				holds.
	\item $\pitch$ satisfies condition (1) of \Cref{prop: altexplicatedsemiorthogonaldecompositioncriterion2} with respect to the semi-orthogonal decompositions $\bb{S}_{\ICG}$ and $\bb{S}_{\Bun_{G}}^{\mathrm{ex}}$. 
			
		\end{enumerate}
	\end{proposition}
	\begin{proof}
		The second and third claim follow from \Cref{extend filtered functor from singleton to arbitrary convex} and from the first claim.
		In turn, the first claim ultimately follows from the definitions, from \Cref{cartesian isoc bungmer} and from the Cartesian square \eqref{important Cartesian}, we expand the argument below. 
		For the rest of the proof, we fix $b\in B(G)$, we put ourselves in the situation of \Cref{lemma: compactgeneratordiagram}. 
		In particular, we recall there is a $\sigma$-straight element $w_{b}$ attached to $b \in B(G)$ and we have open and closed immersions
		\[\Sht^\sch_{\calG,w_b}\to \Sht^\sch_{\calG,\leq w_b} \to\Sht^\sch_{\calG},\]
		together with the Newton map 
        \[ \Nt: \Sht^{\sch}_{\calG} \ra \ICG. \]
        By \Cref{lemma: compactgeneratordiagram}, the map
		\[\Sht^\sch_{\calG,w_b}\to \ICG\]
		factors through $\ICG_b$ and we have an identification $\Sht^\sch_{\calG,w_b}\simeq \bbB_{\profet} \underline{I_{w_b}}$ over $\bbB_{\proet} \underline{G_b(E)}$, where we recall $I_{w_b}\subseteq G_b(E)$ is the Iwahori subgroup attached to some lift $\dot{w}_b$ of $w_b$, as in \Cref{defn: sigmastraightattatchedtob}. 
		In particular, the category $\Shv^!(\ICG_b)$ is generated under colimits by objects of the form $\ul{\on{Nt}_{w_b,*}}A$ with $A\in \Shv^!(\Sht^\sch_{\calG,w_b})$.
		For our purposes, it suffices to show that for all $A\in \Shv^!(\Sht^\sch_{\calG,w_b})$ 
		\[\pitch(\ul{i_{b,!}}\ul{\on{Nt}_{w_b,*}}A)\in (\Bun_G^b)_\sharp.\]
		Indeed, $\pitch\circ \ul{i_{b,!}}$ commutes with colimits, and $(\Bun_G^b)_\sharp$ is stable under colimits. 

		Recall that $\pitch=\pitch^\vee\circ \id_\BZ$.
		We first observe that 
		\begin{equation}{\label{eqn: BZCalc}}
        \id_\BZ(\ul{i_{b,!}}\ul{\on{Nt}_{w_b,*}}A)\simeq (\ul{i_{b,*}})^o (\ul{\on{Nt}_{w_b,*}})^o B
        \end{equation}
		for some $B\in \Shv^!(\Sht^\sch_{\calG,w_b})^\vee$.
		Indeed, this follows from the commutative diagram
		\begin{center}
		\begin{tikzcd}
			\Shv^!(\Sht^\sch_{\calG,w_b}) \arrow{r}{{\ul{\on{Nt}_{w_b,*}}}} \arrow{d}{\id_{\on{BZ},w_b}}  & \Shv^!(\ICG_b) \arrow{r}{\ul{i_{b!}}} \arrow{d}{\id_{\BZ,b}}  & \Shv^!(\ICG) \arrow{d}{{\id}_\BZ} \\
			\Shv^!(\Sht^\sch_{\calG,w_b})^\vee  \arrow{r}{{(\ul{\on{Nt}_{w_b,*}})^o}} & 	\Shv^!(\ICG_b)^\vee  \arrow{r}{(\ul{i_{b,*}})^o} & \Shv^!(\ICG)^\vee.
		\end{tikzcd}
		\end{center}
        Here we recall that the commutativity of the left square follows from  \Cref{lemma: compactgeneratordiagram} (cf. \Cref{prop: BasicPropertiesofBZDualityonIsoc} (2)), and the commutativity of the right most square follows from \eqref{eqn: semiorthogonaldecompositionIsoc}.
        We note that $\Sht^\sch_{\calG,w_b}$ is a placid stack, and that the map 
		\[\Sht^\sch_{\calG,w_b}\xrightarrow{ i_b\circ\on{Nt}_{w_b}} \ICG\]
		factors through $\Sht^\sch_{\calG,\leq w_b}$ which is a $\sigma$-compatible bounded piece of $\ICG$ (see \Cref{sind-dagger-correspondence} and \S\ref{ss: constructin the functors}). 
		Indeed, that $\Sht^\sch_{\calG,\leq w_b}$ is $\sigma$-compatible is the content of \Cref{analytic Newton map is fdcs} and \Cref{resilience lemma ofr sht}. 
		From this, and since $\Sht^\sch_{\calG,w_b}\to \Sht^\sch_{\calG,\leq w_b}$ is pfp representable (see \Cref{thm:isocG_geometry} and \Cref{lemma: compactgeneratordiagram}), we may use \Cref{also-preserves-compact} to write 
		\begin{align*}
		\pitch(\ul{i_{b,!}}\ul{\on{Nt}_{w_b,*}}A) := c^{*,\sigma}\id_\BZ(\ul{i_{b,!}}\ul{\on{Nt}_{w_b,*}}A) &\simeq^{(\ref{eqn: BZCalc})} c^{*,\sigma} \ul{(i_b\circ \on{Nt}_{w_b})}_*^o B \\ &\simeq (\sigma\circ \on{Nt}^\dagger_{\leq w_b})_! \circ k_{b,!} \circ  b^*_{\Sht^\sch_{\calG,w_b}/\ICG} c^{*,\vee}(B).
        \end{align*}
		Here $c^{*,\vee}$ is the analytification functor from \Cref{analytification-cosheave-quasi-plac}, and the rest of the maps fit in the commutative diagram with Cartesian squares
		\begin{equation}
			\label{diagram with a bunch of arrows}
		\begin{tikzcd}
			\text{[}\widetilde{\calM}_b/I_{w_b}\text{]} \arrow{r}{k_{b}} \arrow{d}{b_{\Sht^\sch_{\calG,w_b}/\ICG}}  & \arrow[rr, bend left=30, "\sigma_{\leq w_b}"] \Sht^\an_{\calG,\leq w_b} \arrow{r}{\on{Nt}_{\leq w_b}^\dagger} \arrow{d}{b_{\Sht^\sch_{\leq w_b}}}    & \Bun_G^\mer \arrow{d}{\gamma} \ar{r}{\sigma} & \Bun_G \\
			(\Sht^\sch_{\calG,w_b})^\diamondsuit \arrow{r} & (\Sht^\sch_{\calG,\leq w_b})^\diamondsuit \arrow{r}  & \ICG^\diamondsuit, 
		\end{tikzcd}
		\end{equation}
        obtained from \Cref{thm: GIZ} and \eqref{important Cartesian}.
		The map $[\widetilde{\calM}_b/I_{w_b}] \to \Bun_G$ factors through $\calM_b\to \Bun_G$ and this factorization fits in the following commutative diagram with Cartesian squares 
\begin{center}
\begin{tikzcd}
	\text{[}\widetilde{\calM}_b/I_{w_b}\text{]}  \arrow{r} \arrow{d}  & \calM_b \arrow[bend left]{rr}{\sigma_b} \arrow{r} \arrow{d}{\gamma_b}  & \Bun_G^\mer \arrow{d}{\gamma} \ar{r}{\sigma} & \Bun_G \\
	(\Sht^\sch_{\calG,w_b})^\diamondsuit \arrow{r}{\on{Nt}_{w_b}^\diamondsuit} & \ICG^\diamondsuit_b \arrow{r}{i_b^\diamondsuit} & \ICG^\diamondsuit. 
\end{tikzcd}
\end{center}
By proper base change, and since $\sigma_b$ is $\calD^\an_\Lambda$-$!$-able, we may rewrite 
\[(\sigma\circ \on{Nt}^\dagger_{\leq w_b})_! \circ k_{b,!} \circ  b^*_{\Sht^\sch_{\calG,w_b}/\ICG} c^{*,\vee}B \simeq \sigma_{b!}\gamma_b^* \on{Nt}_{w_b,!}^\diamondsuit c^{*,\vee}B.\]
If we let $D=\on{Nt}_{w_b,!}^\diamondsuit c^{*,\vee}B$, then we have  
\[\sigma_{b!}\gamma_b^* D\in (\Bun_G^b)_\sharp.\]
by \Cref{prop: jbnatural} (see also \Cref{rem: jbsharprestricted}).
\end{proof}

	\subsection{Understanding $\pitch_U$ and $\pitch_b$}
	Since we wish to apply \Cref{prop: altexplicatedsemiorthogonaldecompositioncriterion2} to $\pitch$, we will need to be able to describe, for any finite convex subsets $V\subseteq B(G)$ (including, in particular, single points $b \in B(G)$), all of the functors
	\[
    \pitch_V:\Shv^!(\ICG_V)\to \calD^\an_\Lambda(\Bun^V_G)\]
    that a map of semi-orthogonal decompositions specify.
Here we note that these functors can be defined through \Cref{lem_switch_persp_functors}.
More precisely, one can use the Verdier quotient presentation
\[\Shv^!(\ICG_V)\simeq \Shv^!(\ICG_Z)/\Shv^!(\ICG_{Z'})\]
to construct the maps, where $Z$ and $Z'$ are any two finite closed subsets $Z'\subset Z\subset B(G)$ with $V=Z\setminus Z'$.
In this subsection, we study a different presentation of the $\pitch_V$ via sind-$\dagger$-correspondences.
This other presentation will allow us to exploit the geometry of kimberlites to finish showing that $\pitch$ satisfies the hypothesis in \Cref{prop: altexplicatedsemiorthogonaldecompositioncriterion2}.

	In what follows, we fix $Z\subseteq B(G)$ a closed subset and $U\subseteq Z$ an open subset.
	We get maps 
	\[\ICG_U\to \ICG_Z\to \ICG\]
	which are an open and pfp closed immersion respectively by \Cref{thm:isocG_geometry}.
	It follows from \Cref{passing-to-pfp} and \Cref{cor: ConstructionofMainSindDaggerCorrespondence}, that the triples 
	\[(\ICG_U,\Bun_G,\sigma\circ i_U^\dagger) \text{ and } (\ICG_Z,\Bun_G,\sigma\circ i_Z^\dagger)\]
	are sind-$\dagger$-correspondences. 

	\begin{lemma}
		\label{correct factorizations}
		With the setup as above, let $X\in \{Z, U\}$, and consider the map 
		\[\sigma\circ i_X^\dagger:\ICG_X^\dagger\to \Bun_G.\]
		Then there is a factorization of the form 
		\begin{center}
		\begin{tikzcd}
			\ICG^\dagger_X \arrow{r}{\sigma_X} \arrow{d}{i_X^\dagger}  &\Bun_G^X \arrow{d}{j_X} \\
		 \Bun_G^\mer \ar{r}{\sigma} &\Bun_G. 
		\end{tikzcd}
		\end{center}
	\end{lemma}
	\begin{proof}
		It suffices to argue on topological spaces, i.e., we wish to show that $|\sigma\circ i_X^\dagger|(|\ICG^\dagger_X|)\subseteq |\Bun_G^X|$. 
		This in turn can be done at the level of geometric points. 
		Fix a map $x:\Spa(C,C^+)\to \ICG^\dagger$. 
		Let $k=O_C/C^{\circ \circ}$.
		By \cite[Proposition 7.4(4)]{GIZ25}, this corresponds to a map 
		\[\tilde{x}:\Spec O_C\to \ICG.\]
		Indeed, $\ICG^{\dagger_{\on{pre}}}(\Spa(C,C^+))$ is by definition $\ICG(\Spec O_C)$.
		Moreover, $\tilde{x}$ lifts to $\ICG_X^\dagger$ if and only if the map $\Spec O_C\to \ICG$ lifts to $\ICG_X$, and this happens if and only if the image of $\Spec C$ and of $\Spec k$ in $\ICG$ lie in $X$, but $\sigma(x)\in |\Bun_G|$ matches the image of $\Spec k$ in $|\ICG|$ under the (non-continuous) bijection 
		\[|\Bun_G| \simeq |\ICG|. \qedhere \]
	\end{proof}

	\begin{lemma}
		\label{Z to U cartesian}
	With the setup as above, we have a Cartesian diagram 	
	\begin{center}
	\begin{tikzcd}
	\ICG^\dagger_U \arrow{r} \arrow{d}  & \Bun_G^U \arrow{d} \\
	\ICG^\dagger_Z \arrow{r} & \Bun_G^Z.
	\end{tikzcd}
	\end{center}
    of $v$-stacks.
	\end{lemma}
	\begin{proof}
		As in the proof of \Cref{correct factorizations}, we may argue on geometric points. 
		In this case, it suffices to know that if we have a factorization $\Spec O_C\to \ICG_Z\subseteq \ICG$ and the induced closed point $\Spec k$ factors through $\ICG_U$, then all of $\Spec O_C$ factors through $\ICG_U$.
However, this is clear, since $U$ is open in $Z$. 
	\end{proof}

	The above considerations give rise to maps of sind-$\dagger$-correspondences

	\begin{center}{\label{eqn: SindDaggerCorrespondences}}
	\begin{tikzcd}
	\ICG_U^\diamondsuit  \arrow{d}  &\ICG_U^\dagger \arrow{l} \arrow{r}{\sigma_U} \arrow{d} \arrow[draw=none]{dr}[description]{\text{(A)}}  & \Bun_G^U \arrow{d} \\
	\ICG_Z^\diamondsuit  \arrow{d} \arrow[draw=none]{dr}[description]{\text{(B)}} & \arrow{l} \ICG_Z^\dagger \arrow{r}{\sigma_Z} \arrow{d}  & \Bun_G^Z \arrow{d} \\
	\ICG^\diamondsuit  & \Bun_G^\mer \arrow{l}{\gamma} \arrow{r}{\sigma} & \Bun_G  
	\end{tikzcd}
	\end{center}
	in which the squares $(A)$ and $(B)$ are Cartesian.
Indeed, the fact that $(A)$ is Cartesian is \Cref{Z to U cartesian}. 
That $(B)$ is a Cartesian follows from \Cref{lm:proper_gives_cartesian_diam}.\footnote{Note that the other squares are typically never Cartesian: indeed, the map $Y^\dagger \to Y^\diamondsuit$ does not pullback under open immersions, since an $O_C$-point of $\ICG$ with induced $k$-point in $Z$ might induce a $C$-point that misses $Z$.}

	Recall that we have an identification $\id_\BZ:\Shv^!(\ICG)\simeq \Shv^!(\ICG)^\vee$ respecting the semi-orthogonal decomposition on both categories, and that we have commutative diagrams 
	\begin{center}
\begin{tikzcd}
	\Shv^!(\ICG_U) \arrow{r}{\id_{BZ,U}} \arrow{d}{\ul{i_{U!}}}  & \Shv^!(\ICG_U)^\vee  \arrow{d}{(\ul{i_{U*}})^o} & 
	\Shv^!(\ICG_U) \arrow{r}{\id_{BZ,U}} \arrow{d}{\ul{i_{U*}}}  & \Shv^!(\ICG_U)^\vee  \arrow{d}{(\ul{i_{U!}})^o} \\
 \Shv^!(\ICG) \arrow{r}{\id_\BZ} &  \Shv^!(\ICG)^\vee  & 
 \Shv^!(\ICG) \arrow{r}{\id_\BZ} &  \Shv^!(\ICG)^\vee,   
\end{tikzcd}
\end{center}
as discussed in \Cref{cor: SelfDualizableIsoc}. 

	\begin{proposition}
		\label{truncating pitch}
		For all $Z\subseteq B(G)$ a finite closed subset, and, all $U\subseteq Z$ open, we have an identification of functors
		\[\pitch_U\simeq c^{*,\sigma_U}\circ ({\id}_{\BZ,U}), \]
        where the right-hand side is the analytification map (\ref{sind dagger analytifiction}) attached to the sind-$\dagger$- correspondence $(\ICG_U,\Bun^U_G,\sigma_U)$ described above.
	\end{proposition}
	\begin{proof}
		By \Cref{lem_switch_persp_functors}, we may use the Verdier quotient presentation of $\pitch_U$ to write 
		\[\pitch_U\simeq j^*_U\pitch \ul{i_{U!}}.\]
		From this presentation, we see that $\pitch_U=j_U^*c^{*,\sigma}(\ul{i_{U*}})^o(\id_{\BZ,U})$.
		The $\sigma$-compatible bounded piece $\Sht^{\sch}_{\calG,\leq {w_\bullet}}\to \ICG$ defines a $\sigma_U$-compatible bounded piece 
		\[\on{Nt}_{\leq {w_\bullet},U}:\Sht^{\sch}_{\calG,\leq {w_\bullet},U}\to \ICG_U.\] 
		Fix $B=\ul{\on{Nt}_{\leq {w_\bullet},U,*}}A$ for $A\in \Shv^!(\Sht^{\sch}_{\calG,\leq {w_\bullet},U})$.
		We also let $D_A$ denote the unique object in $\Shv^!(\Sht^{\sch}_{\calG,\leq {w_\bullet},U})$ such that $(\id_{\BZ,U}) B\simeq (\ul{\on{Nt}_{\leq {w_\bullet},U,*}})^o D_A$.
		The construction of $D_A$ is functorial in $A$ since its the dual in $\Shv^!(\Sht^{\sch}_{\calG,\leq {w_\bullet},U})$ is induced by the natural Frobenius algebra structure on $\Shv^!(\Sht^{\sch}_{\calG,\leq {w_\bullet},U})$ which is compatible with $\id_{\BZ,U}$.
		We write $i_{\leq {w_\bullet}, U}$ for the inclusion 
		\[i_{\leq {w_\bullet},U}:\Sht^{\sch}_{\calG,\leq {w_\bullet},U}\to \Sht^{\sch}_{\calG,\leq {w_\bullet}}.\]
		Now, we may compute $\pitch_U(B)$ as follows 
		\begin{center}
			\begin{align*}
			\pitch_U(B)&\simeq j_U^*c^{*,\sigma}(\ul{i_{U*}})^o(\id_{\BZ,U}) B\\ \nonumber 
				   &\simeq j_U^*c^{*,\sigma}(\ul{i_{U*}})^o (\ul{\on{Nt}_{\leq {w_\bullet},U,*}})^o D_A\\ \nonumber 
			&\simeq j_U^*c^{*,\sigma}(\ul{\on{Nt}_{\leq {w_\bullet},*}})^o(\ul{i_{\leq {w_\bullet},U*}})^o D_A  \\ \nonumber 
				 & \simeq j_U^*{\sigma_{\leq {w_\bullet},!}}b_{\Sht^{\sch}_{\calG,\leq {w_\bullet}}}^* c^{*,\vee}(\ul{i_{\leq {w_\bullet},U*}})^o D_A \\
				 & \simeq j_U^*({\sigma_{\leq {w_\bullet},!}})b_{\Sht^{\sch}_{\calG,\leq {w_\bullet}}}^*  ({i^\diamondsuit_{\leq {w_\bullet},U!}}) c^{*,\vee} D_A  \nonumber
			\end{align*}
		\end{center}
		We get a diagram 
	\begin{center}
	\begin{tikzcd}
		(\Sht^\sch_{\calG,{\leq w_\bullet}, U})^\diamondsuit    \arrow{d} \arrow[dd, bend right=70, "{i^\diamondsuit_{{\leq w_\bullet},U}}" description] & \arrow[draw=none]{dr}[description]{\text{(A)}} (\Sht^\sch_{\calG,{\leq w_\bullet}, U})^\dagger  \arrow{l}{b_{\Sht^{\sch}_{\calG,{\leq w_\bullet},U}}^*} \arrow{r}{\sigma_{{\leq w_\bullet},U}} \arrow{d}  & \Bun_G^U \arrow{d}{j_{UZ}} \arrow[dd, bend left=70, "{j_U}" description] \\
		\arrow[draw=none]{dr}[description]{\text{(B)}} 	(\Sht^\sch_{\calG,{\leq w_\bullet}, Z})^\diamondsuit  \arrow{d}  &(\Sht^\sch_{\calG,{\leq w_\bullet}, Z})^\dagger  \arrow{l} \arrow{r}{\sigma_{{\leq w_\bullet},Z}} \arrow{d}  & \Bun_G^Z \arrow{d}{j_Z} \\
	 (\Sht^\sch_{\calG,{\leq w_\bullet}})^\diamondsuit    &(\Sht^\an_{\calG,{\leq w_\bullet}})  \arrow{l}{b_{\Sht^{\sch}_{\calG,{\leq w_\bullet}}}^*} \arrow{r}{\sigma_{{\leq w_\bullet}}}   & \Bun_G  
	\end{tikzcd}
	\end{center}
	where the square $(A)$ and the square $(B)$ are Cartesian, as in (\ref{eqn: SindDaggerCorrespondences}). 
	An elementary, but lengthy, diagram chase using that all vertical arrows are locally closed immersions gives the formula
	\[j_U^*({\sigma_{{\leq w_\bullet},!}})b_{\Sht^{\sch}_{\calG,{\leq w_\bullet}}}^*  i_{{\leq w_\bullet},U,!}^\diamondsuit\simeq \sigma_{{\leq w_\bullet},U,!} b_{\Sht^{\sch}_{\calG,{\leq w_\bullet},U}}^* \]
	We continue computing 
		\begin{center}
			\begin{align*}
			\pitch_U(\ul{\on{Nt}_{{\leq w_\bullet},U,*}}A)&\simeq j_U^*({\sigma_{{\leq w_\bullet},!}})b_{\Sht^{\sch}_{\calG,{\leq w_\bullet}}}^*  ({i^\diamondsuit_{{\leq w_\bullet},U!}}) c^{*,\vee}D_A \\
&\simeq \sigma_{{\leq w_\bullet},U,!} b_{\Sht^{\sch}_{\calG,{\leq w_\bullet},U}}^*c^{*,\vee} D_A \\
& \simeq c^{*,\sigma_U} (\ul{\on{Nt}_{{\leq w_\bullet},U,*}})^o D_A 
			\end{align*}
		\end{center}
		This gives the following commutative diagram 
		\begin{center}
		\begin{tikzcd}
			\Shv^!(\Sht^{\sch}_{\calG,{\leq w_\bullet},U}) \arrow[dd, bend right=90, "\pitch_U\circ \ul{\on{Nt}_{{\leq w_\bullet},U,*}}"'] \arrow{r}{A\mapsto D_A} \arrow{d}{\ul{\on{Nt}_{{\leq w_\bullet},U,*}}}  & \Shv^!(\Sht^{\sch}_{\calG,{\leq w_\bullet},U})^\vee  \arrow{d}{(\ul{\on{Nt}_{{\leq w_\bullet},U,*}})^o} \\
			\Shv^!(\ICG_U) \arrow{r}{\on{id}_{\BZ,U}} &  \Shv^!(\ICG_U)^\vee \ar{dl}{c^{*,\sigma_U}} \\
								    \calD^\an_\Lambda(\Bun^U_G)
		\end{tikzcd}
		\end{center}
		Passing to the colimit as we vary among the bounded pieces of $\ICG_U$ (i.e., we vary $w_\bullet$) we obtain the desired commutative diagram.   
		\begin{center}
		\begin{tikzcd}
			\Shv^!(\ICG_U) \arrow{r}{\id_{\BZ,U}} \arrow{d}{\pitch_U}  & \Shv^!(\ICG_U)^\vee \arrow{dl}{c^{*,\sigma_U}}  \\
			\calD^\an_\Lambda(\Bun^U_G). & 
		\end{tikzcd}
		\end{center}

	\end{proof}

	We can easily specialize the discussion to $\pitch_b$ by taking $Z=\{\beta \in B(G)\mid \beta\leq b\}$ and by taking $U=\{b\}$.

	\begin{corollary}
		\label{pitch equivalence crit (3)}
		For all $b\in B(G)$, the functor 
		\[\pitch_b:\Shv^!(\ICG_b)\to \calD^\an_\Lambda(\Bun_G^b)\]
		is an equivalence.
In particular, the last condition in \Cref{prop: altexplicatedsemiorthogonaldecompositioncriterion2} holds. 
	\end{corollary}
	\begin{proof}
	This follows from \Cref{truncating pitch} and \Cref{analitifying cInd} (1). 
	\end{proof}
    \begin{remark}\label{diagram spell pitch out}
	It will be convenient to spell out precisely the functor $\pitch_b$.
	Recall that we have fixed standard identifications

	\begin{center}
	\begin{tikzcd}
	   & \Shv^!(\ICG_b)^\vee \ar{dr}{c^{*,\id}} \ar{rr}  & &  \calD^\an_\Lambda(\Bun_G^b)   \\
		\Shv^!(\ICG_b)  \ar{ru}{\bbD_{\on{coh},G_{b}(E)}}  \ar{rr}{c^{*,\id} \circ \bbD_{\on{coh},G_{b}(E)}} &   & \calD^\an([\ast/\ul{G_b(E)}]) \ar{ur}{s_{b*}} &   \\
												      & \Rep (G_b(E)) \ar{ul}{\on{Zhu}} \ar{ur}{\on{FS}} & &  .  \\
												  
	\end{tikzcd}
	\end{center}

Here all the identifications are normalized with respect to the same choice of representative $\dot{b} \in G(\Breve{E})$ of $b$, which we may choose to be as in \Cref{conv: DualCpx}.
We note that we have abused the notation and we have written $\bb{D}_{\coh,G_{b}(E)}$ for the equivalence $\Shv^{!}(\ICG_{b}) \xrightarrow{\simeq} \Shv^{!}(\ICG_{b})^{\vee}$ induced by $\eta^{\mathrm{Zhu}}$ and $\bb{D}_{\coh,G_{b}(E)}: \Rep(G_{b}(E)) \xrightarrow{\simeq} \Rep(G_{b}(E))^{\vee}$. 
In the following diagram, we fix these standard identifications which we will simply denote with $\mathrm{std}$.
With all the identifications fixed, the functor $\pitch_b$ can be computed as the only map making the following diagram commutative
	\begin{center}
	\begin{tikzcd}
		\Rep (G_b(E)) \arrow{rr}{\otimes \delta_b^\mathrm{Zhu}[2d_b]} \arrow{d}{\mathrm{std}}  &    & \Rep (G_b(E)) \arrow{d}{\mathrm{std}} \arrow{r}{\otimes \delta_b^{-1}[-2d_b]}   & \Rep (G_b(E)) \arrow{d}{\mathrm{std}} \\
		\Shv^!(\ICG_b) \arrow{r}{\id_{\BZ,b}} \ar[rrr, "\pitch_b", bend right] & \Shv^!(\ICG_{b})^\vee  \arrow{r}{c^{*,\id}} & \calD^\an_\Lambda([\ast/G_b(E)]) \arrow{r}{s_{b!}} & \calD^\an_\Lambda(\Bun_G^b),
	\end{tikzcd}
	\end{center}
    where the commutativity of the left-hand square follows from the discussion following (\ref{eqn: semiorthogonaldecompositionIsoc}) and \Cref{analitifying cInd} (4), and the commutativity of the right-hand square follows from \cite[Corollary~1.8]{HIDualCpx}, where the notation is as in \ref{conv: DualCpx}
    This tells us that the composition 
    \[ \Rep(G_{b}(E)) \xrightarrow{\eta^{\mathrm{Zhu}}} \Shv^{!}(\ICG_{b}) \xrightarrow{\pitch_{b}} \calD_{\Lambda}^{\an}(\Bun_{G}^{b}) \xrightarrow{s_{b}^{*}} \Rep(G_{b}(E)) \]
    is equivalent to the functor 
    \[(-) \otimes \delta_{b}^{-1} \otimes \delta_{b}^{\mathrm{Zhu}}[2d_{b} - 2d_{b}] \simeq (-) \otimes \delta_{b}^{\mathrm{Zhu}} \otimes \delta^{-1}_{b}.\] 
    We set $\chi_{b} := \delta_{b}^{\mathrm{Zhu}} \otimes \delta^{-1}_{b}$. 
    In particular, we see that if we choose the identifications $\Shv^{!}(\ICG_{b}) \simeq \Rep(G_{b}(E))$ and $\calD_{\Lambda}^{\an}(\Bun_{G}^{b}) \simeq \Rep(G_{b}(E))$ with respect to the same choice of representative $\dot{b}$ of $b \in B(G)$, it that follows $\pitch_{b}\simeq \id_{\Rep(G_{b}(E))}$ if and only if $\chi_b\simeq \mathbbm{1}$ (i.e., \Cref{assumpconj: linebundletwist} holds).
	\end{remark}

	\subsection{Henselianity and right-orthogonality of $\pitch$}
	In this section, we finish showing that $\pitch$ is an equivalence.

	\begin{proposition}
		\label{right orthogonality verified}
		Let $Z\subseteq B(G)$ be a finite closed subset, let $U\subseteq Z$ be an open subset and let $b\in U$ be a closed point.
		Then 
		\[\pitch_U(\Shv^{!}(\ICG_{U\setminus b})_*) \subset \calD^\an_\Lambda(\Bun_G^{U\setminus b})_!\]
		in $\calD^\an_\Lambda(\Bun^U_G)$.
		In particular, we deduce that $\pitch$ satisfies condition of \Cref{prop: altexplicatedsemiorthogonaldecompositioncriterion2} (2) with respect to the semi-orthogonal decomposition $\bb{S}_{\ICG}$ on $\Shv^{!}(\ICG)$ and the semi-orthogonal decomposition $\bb{S}_{\Bun_{G}}^{\mathrm{ex}}$ on $\calD_{\Lambda}^{\an}(\Bun_{G})$.
	\end{proposition}
	\begin{proof}
    For the last assertion, note that we have identified $\pitch_U$ with the Verdier quotient map induced by $\pitch$. 
    
    By definition, $\Shv^{!}(\ICG_{U\setminus b})_*$ is generated under colimits by objects of the form $\ul{i_{\beta U,*}}A$ with $\beta\in U \setminus \{b\}$. 
		Since $\pitch$ is colimit-preserving, it suffices to show that 
		\[\pitch_U(\ul{i_{\beta U,*}}A)\in (\calD^\an_\Lambda(\Bun^{U\setminus b}_G))_!.\]
		For now, we fix $B=\id_{\BZ,\beta}(A)\in \Shv^!(\ICG_\beta)^\vee$. 
		We have that 
		\begin{align}
			\label{computation pre orthogonality2}
			\pitch_U\ul{i_{\beta U,*}}A &      \simeq c^{*,\sigma_U} \id_{\BZ,U} 
            \ul{i_{\beta U,*}}  A \\ \nonumber
											& \simeq c^{*,\sigma_U} (\ul{i_{{\beta U},!}})^o B 
		\end{align}
		by \Cref{cor: SelfDualizableIsoc}, and we wish to show that $c^{*,\sigma_U} (\ul{i_{{\beta U},!}})^o B\in (\calD^\an_\Lambda(\Bun^{U\setminus b}_G))_!$ for all $B\in \Shv^!(\ICG_\beta)^\vee$.

		As in \Cref{defn: sigmastraightattatchedtob}, for $\beta\in B(G)$ we can attach a straight element $w_\beta\in B(\widetilde{W})$ and an Iwahori subgroup $I_\beta\subseteq G_\beta$. 
		Moreover, recall from \Cref{lemma: compactgeneratordiagram} and \Cref{analitifying cInd} (2), that if $K\subseteq I_\beta$ is a pro-$p$ subgroup, then we have an isomorphism
		\[(\ul{i_{\beta U,!}})^o \on{c-Ind}_{K}^{G_\beta}\Lambda\simeq (\ul{\Nt_{\leq w_b,U,*}})^o\circ (\ul{i_{w_{\beta U,!}}})^o\on{c-Ind}_{K}^{I_\beta}\Lambda.\]
		Recall the maps $\sigma_{\leq w_\beta}$ and $k_{\beta,K}$ from \eqref{diagram with a bunch of arrows}.
		We will compute as follows and justify the steps below,
		\begin{align}
			\label{computation pre orthogonality}
		c^{*,\sigma_U} (\ul{i_{{\beta U},!}})^o \on{c-Ind}_{K}^{G_\beta}\Lambda & \simeq c^{*,\sigma_U} (\ul{\Nt_{\leq w_b,U,*}})^o\circ (\ul{i_{w_{\beta} U,!}})^o \on{c-Ind}_{K}^{I_{\beta}}\Lambda \\
									     & \simeq  {\sigma_{\leq w_b,U,!}}\circ b_{\Sht^\sch_{\leq w_b,U}}^* c_{\Sht^\sch_{\leq w_b,U}}^{*,\vee} (\ul{i_{w_{\beta} U,!}})^o \on{c-Ind}_{K}^{I_{\beta}}\Lambda \\
									     & \simeq  {\sigma_{\leq w_b,U,!}}\circ b_{\Sht^\sch_{\leq w_b,U}}^* c_{\Sht^\sch_{\leq w_b,U}}^{*} {i_{w_{\beta} U,*}} \on{c-Ind}_{K}^{I_{\beta}}\Lambda \\
									     & \simeq  {\sigma_{\leq w_b,U,!}}{k_{{\beta},I_{w_\beta},*}} \on{c-Ind}_{K}^{I_{\beta}}\Lambda.
		\end{align}
		The first step was discussed above.  
		The second step is the definition of $c^{*,\sigma_U}$ in \eqref{sind dagger analytifiction} since $\Sht^\sch_{\leq w_b,U}$ is a $\sigma_U$-compatible bounded piece.
		The third step comes from the equivalence 
		\[\Shv^!(\Sht^\sch_{\leq w_b,U})^\vee\xrightarrow{\simeq} \Shv^*(\Sht^\sch_{\leq w_b,U})\]
		(see \Cref{very-placid-givesequiv}), the definition of $c^{*,\vee}$ (see \Cref{analytification-cosheave-quasi-plac}) and the by-hand verification that $\on{c-Ind}_{K}^{I_{\beta}}\Lambda$ is a compact in $\Shv^!(\bbB \underline{I_\beta})$. 
		Indeed, by definition of the conjugate functor, $(\ul{i_{w_{\beta},!}})^o$ may be computed by $(\ul{i_{w_{\beta},!}})^\op$ on compact objects.
		Moreover, we have an adjoint pair 
		\[ (\ul{i^{!}_{w_{\beta}}})^\op \dashv (\ul{i_{w_{\beta},!}})^\op \]
		of functors between the categories $\Shv^!(\Sht^\sch_{\calG,w_b})^{\omega,\op}$ and $\Shv^!(\Sht^\sch_{\calG,\leq w_b})^{\omega}$.
		Indeed, both $\ul{i^{!}_{w_{\beta}}}$ and $\ul{i_{w_{\beta},!}}$ admit continuous right adjoints (since $i_{w_{b}}$ is an open immersion Theorem \ref{thm:isocG_geometry} (2)), so they preserve compact objects. 
		Under the identifications (see \Cref{very-placid-givesequiv})  
		\[\Shv^!(\Sht^\sch_{\calG,w_b})^{\omega,\op}\simeq \Shv^*(\Sht^\sch_{\calG,w_b})^{\omega} \text{ and } \Shv^!(\Sht^\sch_{\calG,\leq w_b})^{\omega,\op}\simeq \Shv^*(\Sht^\sch_{\calG,\leq w_b})^{\omega},\]
		the functors $(\ul{i^{!}_{w_{\beta}}})^\op$ and $i^{*}_{w_\beta}$ are intertwined, so their right adjoints are also intertwined. 
		This finishes justifying the third step.
		The fourth step follows from \Cref{base change formulas}. 

		In what follows, we show that 
		\[{\sigma_{\leq w_{\beta},U,!}}{k_{{\beta},I_{\beta},*}} \on{c-Ind}_{K}^{I_{\beta}}\Lambda\in (\calD^\an_\Lambda(\Bun^{U\setminus b}_G))_!\] 
		for all $\beta\in U\setminus \{b\}$. 
		Equivalently, since $\Bun_G^b\subseteq \Bun_G^U$ is open it suffices to show that
		\[{j^*_b\sigma_{\leq w_{\beta},U,!}}{k_{{\beta},I_{\beta},*}} \on{c-Ind}_{K}^{I_{\beta}}\Lambda\simeq 0.\] 
	Our argument will show that  
	\[ {j^*_b\sigma_{\leq w_{\beta},U,!}}{k_{{\beta},I_{\beta},*}} D\simeq 0\] 
	holds, for arbitrary $D$. 
In particular, it will hold for $D=\on{c-Ind}_{K}^{I_{w_\beta}}\Lambda$.
		Consider the following diagram with Cartesian squares
		\begin{center}
			\begin{tikzcd}
				{\Sht_\calG^{\gamma=w_\beta}}(b) \arrow[d] \arrow{r} &   {\text{[}\Sht_{\calG,U}^{\gamma=w_\beta}}(b)/\ul{K_{b}}\text{]}\ar{d}{k^{b}_{\beta, I_\beta}} \ar{rr}                       &                        & {\text{[}\widetilde{\calM_{U,\beta}}/\ul{I_\beta}\text{]}} \arrow{d}{k_{\beta,I_\beta}} \\
				{\Sht_{\calG,\leq w_\beta}(b)} \arrow[d] \arrow[r]   & {\text{[}\Sht_{\calG,\leq w_\beta}(b)/\ul{K_{b}}\text{]}} \arrow{d}{\pi_{K_b}} \arrow{r} & {\Sht^{b}_{\calG,\leq w_\beta}}\arrow[r] \arrow[d] & {\Sht^\an_{\calG,\leq w_\beta, U}} \arrow{d}{\sigma_{\leq w_\beta, U}} \\
				{\ast} \arrow[r]             & {[\ast/\ul{K_{b}}]} \arrow[r]   \arrow[rr,"j_{K_{b}}"', bend right]        & {\Bun_G^{b}} \arrow[r, "j_{b}"']           & {\Bun^U_G}.         
			\end{tikzcd}
		\end{center}
		Here $K_b\subseteq G_b(E)$ is any auxiliary pro-$p$ compact open subgroup.
		Observe that the map $\bbB\ul{K_{b}}\to \Bun_G^U$ is $\calD_{\Lambda}^{\an}$-smooth, and that pullback along the map $[\ast/\underline{K_{b}}] \ra \Bun_{G}^{b}$ is conservative. 
		For all $K_{b}\subseteq G_{b}$ and all $D$ we have that
		\[j_{K_{b}}^*\sigma_{\leq w_b,U,!}k_{\beta, I_\beta, *}D \simeq \pi_{K_b,!}k^{b}_{\beta, I_\beta,*} E\]
		for some $E$ by smooth and proper base change. 
		Recall from \Cref{Sht gamma = beta contained in analytic locus} that 
		\[{\Sht_{\calG,U}^{\gamma=w_\beta}}(b)\subseteq {\Sht_{\calG,\leq w_\beta,U}(b)}^\an\] since $\beta\in U\setminus \{b\}$, so we may apply \Cref{hyperbolic localization for some quotients} to conclude  
		\[\pi_{K_b,!}k^{b}_{\beta,I_\beta,*} E\simeq 0. \]

	Indeed although, as explained in \S \ref{sec: kimberlitecalculations}, $\Sht_{\calG,\leq w_\beta,U}(b)$ is not a spatial kimberlite with proper reduced locus, it is a locally spatial kimberlite whose reduced locus is an ind-proper scheme which is locally pfp (see \Cref{colimit presentation of Sht G leq w with required properties}).  
		Writing
		\[\Sht_{\calG,\leq w_\beta,U}(b)\simeq \colim X_i\]
		as in \Cref{colimit presentation of Sht G leq w with required properties}, we can apply \Cref{hyperbolic localization for some quotients} to each piece of cohomology corresponding to $\widehat{\Sht_{\calG,\leq w_b}(b')}_{/S}$, with $S$ an suitable union of irreducible components of the affine Deligne--Lusztig variety (i.e., we apply \Cref{hyperbolic localization for some quotients} to each $X_i$). 
		Since the total compactly supported cohomology is the colimit as $S$ increases, and term by term on $S$ the complex vanishes we get a total vanishing. 
	\end{proof}

	\begin{theorem}
		\label{thm: MainTheoremPartial}
		The following statements hold.	
		\begin{enumerate}
			\item The functor $\pitch$ is an equivalence. 
			\item For every finite convex subset $V\subseteq B(G)$ we have an equivalence \[\pitch_V:\Shv^!(\ICG_V)\to \calD^\an_\Lambda(\Bun^V_G).\]
			\item For every pair of convex subsets $V_1\subseteq V_2$ we get intertwining formulas 
				\begin{enumerate}
                    \item $\pitch_{V_2} \ul{i^*_{V_1V_2}}\simeq j^\flat_{V_1V_2}\pitch_{V_1}$
					\item $\pitch_{V_2}\ul{i_{V_1V_2*}}\simeq j_{V_1V_2!}\pitch_{V_1}$
					\item $\pitch_{V_2} \ul{i^\sharp_{V_1V_2}} \simeq j^!_{V_1V_2}\pitch_{V_1}$.
					\item $\pitch_{V_2}\ul{i_{V_1V_2!}}\simeq j_{V_1V_2\sharp}\pitch_{V_1}$
					\item $\pitch_{V_2} \ul{i^!_{V_1V_2}} \simeq j^*_{V_1V_2} \pitch_{V_1}$ 
					\item $\pitch_{V_2}\ul{i_{V_1V_2\flat}}\simeq j_{V_1V_2*}\pitch_{V_1}$.
				\end{enumerate}
			\item For every $b\in B(G)$, we fix a choice of representative $\dot{b}$ of $b$ as in Convention \ref{conv: DualCpx} and consider the associated identifications as in \ref{conv: DualCpx} and \ref{conv: identifyinghseavesonIsocwithReps}, we let $\chi_{b} := \delta_{b,\Zhu} \otimes \delta^{-1}_{b}$ be as in \Cref{diagram spell pitch out}, then we have that:
				\begin{enumerate}
					\item[a)] For every $A\in \Rep(G_b(E))$ the identity $\pitch(\ul{i_{b,*}}\eta^{\mathrm{Zhu}}(A)) \simeq j_{b,!}p_{b}^{*}(A \otimes \chi_{b})$ holds.  
					\item[b)] For every $A\in \Rep(G_b(E))$ the identity $\pitch(\ul{i_{b,!}}\eta^{\mathrm{Zhu}}(A) ) \simeq j_{b,\sharp}p_{b}^{*}(A \otimes \chi_{b})$ holds.  
				\end{enumerate}
				Moreover, these identifications are functorial in $A$.
                \item There is a unique equivalence of presentable 2-functor formalisms 
	\[\pitch_\bbS:\bbS_{\ICG}\xrightarrow{\simeq}\bbS^{\mathrm{ex}}_{\Bun_G}.\]
Moreover, the equivalence $\pitch$ can be rewritten as
	\[\pitch\simeq \colim_{Z\subseteq B(G)} \pitch_Z\]
	where $Z\subseteq B(G)$ ranges over finite closed subsets of $B(G)$ and $\pitch_{Z}$.        
                \item There is an equivalence of presentable 2-functor formalisms 
			\[\pitch_{\bbS^{\on{ex}}}:\bbS^{\on{ex}}_{\ICG}\xrightarrow{\simeq}\bbS_{\Bun_G}.\]
		\end{enumerate}
	\end{theorem}

	\begin{proof}
		The fact that $\pitch$ is an equivalence was proven in several steps by means of \Cref{prop: altexplicatedsemiorthogonaldecompositioncriterion2}. 
		Indeed, by construction $\pitch$ is a functor in $\LinCat_\Lambda$, and the first condition in \Cref{prop: altexplicatedsemiorthogonaldecompositioncriterion2} was shown in \Cref{left orthogonality is shown}.
		The third condition was shown in \Cref{pitch equivalence crit (3)}. 
		Finally, the second condition of \Cref{prop: altexplicatedsemiorthogonaldecompositioncriterion2} is the content of \Cref{right orthogonality verified}.
		This finishes the proof of the first claim.
        
		For the fifth claim, note that \Cref{prop: altexplicatedsemiorthogonaldecompositioncriterion2} already provides an equivalence of 2-functor formalisms
		\[\pitch: \bbS_{\ICG}\xrightarrow{\simeq} \bbS^{\on{ex}}_{\Bun_G}.\]
		Evaluating it on $V$, i.e., $\pitch_V:\bbS_{\ICG}(V)\xrightarrow{\simeq}\bbS^{\on{ex}}_{\Bun_G}(V)$ gives rise to the equivalence
\[\pitch_V:\Shv^!(\ICG_V)\to \calD^\an_\Lambda(\Bun^V_G),\]
and shows the second claim.
The last claim follows from the fifth claim by passing to right adjoint functors. 

The intertwining in the third claim follow directly from fifth claim by virtue of having an equivalence of 2-functor formalisms.
Indeed, the fifth claim itself provides the intertwining of $(a)$ and $(d)$, and passing (once or twice) to right adjoints provides the rest of the intertwining formulas.

		The fourth claim follows from the third. 
		Indeed, for each $b\in B(G)$ we get formulas 
		\[\pitch(\ul{i_{b,*}}\eta^{\mathrm{Zhu}}(A)))\simeq j_{b,!}\pitch_b( \eta^{\mathrm{Zhu}}(A)) \text{ and } \pitch(\ul{i_{b,!}}\eta^{\mathrm{Zhu}}(A)))\simeq j_{b,\sharp}\pitch_b( \eta^{\mathrm{Zhu}}(A)), \]
		and the formula from the statement then follows from the identity
		\[\pitch_b(\eta^{\mathrm{Zhu}}(A))\simeq p_{b}^{*}(A \otimes \chi_{b})\]
		deduced in \Cref{diagram spell pitch out}. 
	\end{proof}

	\begin{remark}
    \label{convenient-table of ajd}
	For the convenience of the reader, we make the following table of adjunctions that are intertwined by $\pitch$.
\begin{table}[h]
\centering
\begin{minipage}{0.45\textwidth}
\centering
\begin{tabular}{ccccc}
$j^\flat$ & $\dashv$ &$j_!$ & $\dashv$ &$j^!$ \\
$\ul{i^*}$ & $\dashv$ &$\ul{i_*}$ &$\dashv$ & $\ul{i^\sharp}$ \\
\end{tabular}
\end{minipage}
\hfill
\begin{minipage}{0.45\textwidth}
\centering
\begin{tabular}{ccccc}
$j_\sharp $& $\dashv$ &$j^* $ & $\dashv$ &$j_* $\\
$\ul{i_!}$ & $\dashv$ &$\ul{i^!} $& $\dashv$ &$\ul{i_\flat} $\\
\end{tabular}
\end{minipage}
\end{table}

	\end{remark}

        \section{Properties and applications of the equivalence}\label{sec:properties_and_applications}
	Throughout this section, we will fix representatives for $b \in B(G)$ and implicitly identify $\Shv^{!}(\ICG_{b},\Lambda) \simeq_{\eta^{\on{Zhu}}} \Rep(G_{b}(E))$ and $\calD_{\Lambda}^{\an}(\Bun_{G}^{b}) \simeq_{s_b^*} \Rep(G_{b}(E))$ with the category of smooth represenations of $G_{b}(E)$ on $\Lambda$-modules, as in \Cref{conv: identifyinghseavesonIsocwithReps} and \Cref{conv: DualCpx}.
	Moreover, we will from now on omit $\eta^{\on{Zhu}}$ and $s_b^*$ from the notation.
We now combine \Cref{thm: MainTheoremPartial} with the results of \cite{Zhu25} on the tame categorical local Langlands conjecture.
We first review the precise statement of this result.
        \subsection{The Tame Categorical Local Langlands Conjecture}
        We now discuss the tame version of the CLLC proven by Zhu \cite{Zhu25}.  
It relates Ind-coherent sheaves on the stack of tame Langlands parameters to the subcategory of $\Shv^{!}(\ICG)$ consisting of tame sheaves. 
We first describe both sides of this equivalence, and then compare them with analogous constructions for $\calD_{\Lambda}^{\an}(\Bun_{G})$. 
Since some of what we discuss holds for arbitrary $\bbZ_\ell$-algebras, we will depart from our usual convention of assuming that $\Lambda/\bb{Z}_{\ell}$ is a torsion algebra and also allow for the case where $\Lambda = \mathcal{O}_{L},L$ for $L/\bb{Q}_{\ell}$ an algebraic extension throughout this subsection. 
Nevertheless, all of our new results will only be formulated and shown in the torsion case since they depend on \Cref{thm: MainTheoremPartial}. 
We first discuss the spectral side.
        \subsubsection{The Stack of Langlands Parameters}{\label{ss: stackofLanglandsparameters}}
        As before, $E$ is a non-Archimedean local field with Weil group $W_{E}$. 
	Given a $\bb{Z}_{\ell}$-algebra $A$, we can endow it with a topology (or the structure of a condensed ring) by writing $A = \colim_{A' \subset A} A'$, where $A'$ is a finitely generated $\mathbb{Z}_{\ell}$-module equipped with its $\ell$-adic topology. 
	We follow \cite[\S VIII]{FS21}, and consider the moduli space, denoted $\mathcal{Z}^{1}(W_{E},\hat{G})$, whose $A$-points are the continuous $1$-cocycles $W_{E} \rightarrow \hat{G}(A)$ with respect to the natural action of $W_{E}$ on $\hat{G}(A)$, where we regard $\hat{G}/\mathbb{Z}_{\ell}$, the split reductive group with root datum dual to the absolute root datum of $G$ as a reductive group over $\bb{Z}_{\ell}$. 
	
	This defines a scheme considered in \cite{DHKM20,Zhu20,FS21} over $\bb{Z}_{\ell}$ which, by \cite[Theorem~I.8.1,Theorem~VIII.0.1]{FS21}, can be written as a union of open and closed affine subschemes $\mathcal{Z}^{1}(W_{E}/P,\hat{G})$ as $P$ runs through the (sufficiently deep) finite index subgroups of the wild inertia of $W_{E}$. 
	For a fixed $P$, $\mathcal{Z}^{1}(W_{E}/P,\hat{G})$ is a flat local complete intersection over $\mathbb{Z}_{\ell}$ of dimension $\dim(G)$ (see \cite[Theorem VIII.1.3]{FS21}). 
	Concretely, $\mathcal{Z}^{1}(W_{E}/P,\hat{G})$ is the affine subscheme parametrizes the cocyles $W_{E} \ra \hat{G}(A)$ which factor through $W_{E}/P$. 
	We may then consider the stack quotient $[\mathcal{Z}^{1}(W_{E},\hat{G})/\hat{G}]$, where $\hat{G}$ acts on such cocycles via conjugation. 
	We denote this stack quotient by $\Par_{G,E}$ and refer to it as the stack of Langlands parameters. 
        
        We now wish to further decompose these stacks into connected components.
To do this, we consider the map 
        \begin{tikzcd}{\label{eqn: tikzcdmapfromcoarsequotienttothestackquotient}}
         q_{E}: \Par_{G,E} \ra \Par_{G,E}^{\coarse},
        \end{tikzcd}
        from the stacky quotient to the coarse quotient of $\mathcal{Z}^{1}(W_{E},\hat{G})$ by $\hat{G}$. 
	More precisely, 
	\[\Par_{G,E}^{\coarse} := \mathrm{Spf}(H^{0}(R\Gamma(\Par_{G},\mathcal{O}_{\Par_{G,E}})))\] 
	is the formal scheme attached to the ring of global functions (where we regard $H^{0}(R\Gamma(\Par_{G},\mathcal{O}_{\Par_{G}}))$ as a pro-object, as in the discussion proceeding \cite[Theorem~2.3]{Zhu25}). 
	
	The $\Lambda$-valued points of $\Par_{G,E}^{\coarse}$ for $\Lambda$ an algebraically closed field are given by $\hat{G}$-conjugacy classes of continuous semisimple (in the sense of \cite[Definition~VIII.3.1]{FS21})
    maps $W_{E} \ra \phantom{}^{L}G(\Lambda)$ (see \cite[Proposition~VIII.3.2]{FS21}), and we may think of the map $q_{E}$ as the semi-simplification map in terms of this moduli interpretation.

         We write $\Breve{E}$ for the completion of the maximal unramified extension of $E$. 
	 We have an analogous stack $\Par_{G,\Breve{E}}$ which parameterizes $\hat{G}$-conjugacy classes of continuous cocycles from the inertia subgroup $I_{E}$ of $W_{E}$ to $\hat{G}$, and an analogous coarse quotient 
	 \[\Par_{G,\Breve{E}}^{\coarse} := \mathrm{Spf}(H^{0}(R\Gamma(\Par_{G},\mathcal{O}_{\Par_{G,\Breve{E}}}))),\] 
	 we refer the reader to \cite[\S 2.1.3]{Zhu25}.

	 If we fix a representative of Frobenius $\sigma \in W_{E}$, then we can check that conjugation by $\sigma$ defines an automorphism 
	 
        \[ \phi: \Par_{G,\Breve{E}} \ra \Par_{G,\Breve{E}},  \]
        as in \cite[Equation 2.19]{Zhu25}.  We write $(\Par_{G,\Breve{E}}^{\coarse})^{\phi} \subset \Par_{G}^{\coarse}$ for the (classical) $\phi$-fixed points; in other words, the formal scheme  $\mathrm{Spf}(H^{0}(R\Gamma(\Par_{G},\mathcal{O}_{\Par_{G,\Breve{E}}}))^{\phi})$. 
	By \cite[Lemma~2.10]{Zhu25}, every connected component of $(\Par_{G,\Breve{E}}^{\coarse})^{\phi}$ is finite over $\bb{Z}_{\ell}$. 
        
	When $\Lambda/\bb{Z}_{\ell}$ is an algebraically closed field then we call a $\Lambda$-point of $(\Par_{G,\Breve{E}}^{\coarse})^{\phi}$ an \emph{inertial type} as in \cite[Definition 2.12]{Zhu25}. 

	By \cite[Lemma~2.13]{Zhu25}, an inertial type is the same as the datum of a semi-simple homomorphism $\breve{\tau}^{\mathrm{ss}}: I_{F} \ra \phantom{}^{L}G(\Lambda)$ with finite image, that can be extended to a homomorphism $\tau:W_F\to \phantom{}^{L}G(\Lambda)$. 
	When $\Lambda = \ol{\mathbb{Q}}_{\ell}$, these will be the natural $L$-parameters attached to the families of smooth irreducible representations living in the Bernstein blocks of $\Rep(G_{b}(E),\Lambda)$ under any reasonable form of the semi-simplified local Langlands, justifying the terminology (However, we note that of course multiple Bernstein blocks could give rise to the same inertial type in this sense, due to the existence of a non-singleton $L$-packets). 
	
	We now consider the composition
        \begin{equation}{\label{eqn: maptocoarsemodulispace}} 
        \Par_{G,E} \xrightarrow{\mathrm{res}} \Par_{G,\Breve{E}} \xrightarrow{q_{\Breve{E}}} \Par_{G,\Breve{E}}^{\mathrm{coarse}}, 
        \end{equation}
        where the first map is given by restriction to inertia and $q_{\Breve{E}}$ is defined analogously to $q_{E}$. 
	Since $\Par_{G,E}$ is given by the $\phi$-fixed points of $\Par_{G,\Breve{E}}$ as explained above, it is easy to see that the map (\ref{eqn: maptocoarsemodulispace}) factors through the subspace $(\Par_{G,\Breve{E}}^{\coarse})^{\phi} \subset \Par_{G,\Breve{E}}^{\coarse}$. 
	This allows us to define the following. 
        \begin{definition}{\label{defn: tameinertialtype}}
		For an inertial type $\zeta\in (\Par_{G,\Breve{E}}^{\mathrm{coarse}})^\phi(\Lambda)$ defined over an algebraically closed field $\Lambda$, we define $\Par_{G,E}^{\hat{\zeta}} \subset \Par_{G,E} \otimes_{\bb{Z}_{\ell}} \Lambda$ to be the base-change of $\Par_{G,E}$ along the morphism (\ref{eqn: maptocoarsemodulispace}) to the formal neighborhood of the closed point attached to $\zeta$ inside $(\Par_{G,\Breve{E},\Lambda}^{\coarse})^{\phi}$. 
        \end{definition}    
		By \cite[Lemma~2.15]{Zhu25}, this is a finite union of connected components of $\Par_{G,E} \otimes \Lambda$; in particular, it is a classical algebraic stack over $\Spec(\Lambda)$. 
        One easily sees that this gives a decomposition 
        \begin{equation}{\label{eqn: decompositionintoconnectedcomponentsFullStack}}
         \Par_{G,E} \otimes \Lambda := \bigsqcup_{\zeta} \Par_{G,E}^{\hat{\zeta}} 
        \end{equation}
        for an algebraically closed field $\Lambda$, where $\zeta$ runs over inertial types.
We now focus on a specific part of $\Par_{G,E}$; namely, the tame locus.

        \subsubsection{The Tame and Unipotent Stack of Langlands Parameters}{\label{ss: TameUnipotentStackofLanglandsParameters}}
        We assume for the rest of this subsection that $G$ is a connected reductive group which is tame; in other words, its splitting field, denoted $\tilde{E}/E$, is a tame extension.
We consider the inclusions 
        \[ P_{E} \subset I_{E} \subset W_{E} \]
        given by the inclusion of the wild inertia and inertia subgroup.
We consider the open and closed (by the discussion at the beginning of \S\ref{ss: stackofLanglandsparameters}) substack 
        \begin{equation}{\label{eqn: tamestackofparameters}}
         \Par_{G,E}^{\tame} := [\mathcal{Z}^{1}(W_{E}/P_{E},\hat{G})/\hat{G}] \subset \Par_{G}, 
        \end{equation}
        defined by parameters that factor through $W_{E}^{t} := W_{E}/P_{E}$, the tame Weil group.
Here we note that the moduli space of cocycles $\mathcal{Z}^{1}(W_{E}/P_{E},\hat{G})$ is well-defined by virtue of the fact that the action of $W_{E}$ on $\hat{G}$ factors through the quotient map to the tame Weil group by our assumption on $G$.
        
        As in \S \ref{ss: stackofLanglandsparameters}, we can analogously define spaces $\Par_{G,\Breve{E}}^{\tame}$ such that $\Par_{G,E}^{\tame}$ is the $\phi$-fixed points of $\Par_{G,\Breve{E}}^{\tame}$.
We then define $\Par_{G,\Breve{E}}^{\tame,\coarse}$ by looking at the ring of global functions on $\Par_{G,\Breve{E}}^{\tame}$ and let  $(\Par_{G,\Breve{E}}^{\tame,\coarse})^{\phi}$ be its fixed points.
As before, we have a map 
        \[ \Par_{G,E}^{\tame} \ra \Par_{G,\Breve{E}}^{\tame} \ra \Par_{G,\Breve{E}}^{\tame,\coarse},\]
        which factors through $(\Par_{G,\Breve{E}}^{\tame,\coarse})^{\phi} \subset \Par_{G,\Breve{E}}^{\tame,\coarse}$.
For $\Lambda$ an algebraically closed field, we say an inertial type $\zeta$, as defined in \Cref{ss: stackofLanglandsparameters}, is tame if it comes from a $\Lambda$-point of $(\Par_{G,\Breve{E}}^{\tame,\coarse})^{\phi}$.
We then have a map
        \[ \Par_{G,E}^{\tame} \ra \Par_{G,\Breve{E}}^{\tame} \ra \Par_{G,\Breve{E}}^{\tame,\coarse}, \]  
        which factors through $(\Par_{G,\Breve{E}}^{\tame,\coarse})^{\phi} \subset \Par_{G,\Breve{E}}^{\tame,\coarse}$, the subspace given by the $\phi$ fixed points.
This leads to a decomposition
        \begin{equation}{\label{eqn: decompositionintoconnectedcomponents}}
         \Par_{G,E}^{\tame} \otimes \Lambda := \bigsqcup_{\zeta} \Par_{G,E}^{\hat{\zeta}} 
        \end{equation}
        where $\zeta$ runs over tame inertial types, which is the base-change of the decomposition (\ref{eqn: decompositionintoconnectedcomponentsFullStack}) along the natural map $\Par_{G,E}^{\tame} \otimes \Lambda \hookrightarrow \Par_{G,E} \otimes \Lambda$.
        
        We let $\tau$ be a choice of generator of the tame inertia $I_{E}^{\mathrm{t}} := I_{E}/P_{E}$, and write $\ol{\tau}$ for the image of $\tau$ under the natural homomorphism $W_{E} \ra W_{E}/W_{\tilde{E}}$, where we recall that $\tilde{E}$ is the splitting field of $G$.
In particular, we note the element $\ol{\tau}$ acts on $\hat{G}$ and its maximal torus $\hat{T}$ via the Galois action on the absolute root datum of $G$. We let $\hat{S} := \hat{T}/(1 - \ol{\tau})\hat{T}$ denote the torus given by coinvariants with respect to this action.
Set $W_{G}$ to be the absolute Weyl group of $G$, and let $W_{0} := W_{G}^{\ol{\tau}}$ denote its $\ol{\tau}$-invariants.
This will act on the torus $\hat{S}$.
Now, we 
        note, by \cite[Lemma~2.36]{Zhu25}, that there is a correspondence between the set of tame inertial types and pairs $(w,\theta)$, where $\theta: I_{E}^{\mathrm{t}} \ra \hat{S}(\Lambda)$ is a homomorphism with finite image up to $W_{0}$-conjugacy, and $w \in W_{0}$ is some element such that $w\sigma(\theta) = \theta^{q}$, where $q$ is the order of the residue field of $E$.

        We consider now the special case of the tame inertial type where $\zeta = \mathrm{triv}$ corresponds to the trivial homomorphism $\theta: I_{E}^{\mathrm{t}} \ra \hat{S}(\Lambda)$ and $\ol{\tau} = 1$. \emph{In particular, we assume now that $G$ splits over an unramified extension of $E$.} Then, for every algebraically closed field $\Lambda$, the element $\zeta$ defines a set of connected components of $\Par_{G,E} \otimes \Lambda$.
However, as explained in \cite[\S~2.2.2]{Zhu25}, this actually comes from a closed substack 
        \begin{equation}{\label{eqn: unipotentstack}}
         \Par_{G,E}^{\widehat{\mathrm{unip}}} \hookrightarrow \Par_{G,E}, 
        \end{equation}
        over $\bb{Z}_{\ell}$, which is an ind-algebraic stack whose base-change to any algebraically closed field is just a regular algebraic stack.
We refer to this as the \emph{stack of unipotent Langalnds parameters}.
        \begin{remark}{\label{rem: CgroupvsLgroup}}
        We note that, in Zhu's work \cite{Zhu25}, in order to make the tame local Langlands correspondence extracted from his categorical equivalence independent of a choice of isomorphism $\ol{\bb{Q}}_{\ell} \simeq \bb{C}$, he works with the stack of Langlands parameters defined with respect to the $C$-group, while here we have chosen to work with the $L$-group instead in order to make it more compatible with the constructions of \cite{FS21}.
However, if we consider the extension $\bb{Z}_{\ell}[\sqrt{q}]/\bb{Z}_{\ell}$, where we have adjoined a square root of the residue characteristic $q$, then we will obtain that the base-change of $\Par_{G}$, $\Par_{G}^{\widehat{\unip}}$, and $\Par_{G}^{\tame}$ to  $\bb{Z}_{\ell}[\sqrt{q}]$, and under this identification we will obtain the stacks considered in \cite{Zhu25} (see \cite[Remark~2.4]{Zhu25}).
        \end{remark} 
        \subsubsection{Coherent Sheaves on the Stack of Parameters}
        We now fix our notation for categories  of coherent sheaves on the stack of Langlands parameters introduced in \Cref{ss: stackofLanglandsparameters} and \Cref{ss: TameUnipotentStackofLanglandsParameters}.
        \begin{definition}{\cite[Section~9.1,Definition 9.19]{Zhu25}}
        For $X$ an ind-algebraic stack, we consider the following categories.
        \begin{enumerate}
		\item We write $\Perf(X)$ (resp. $\QCoh(X)$) for the idempotent complete stable $\infty$-category of perfect complexes with quasicompact support on $X$ (resp. quasi-coherent sheaves on $X$).
        \item We write $\Coh(X)$ for the idempotent complete stable $\infty$-category category of quasi-coherent sheaves with bounded coherent cohomology and quasicompact support.
        \item We write  $\IndPerf(X)$ (resp. $\IndCoh(X)$) for the Ind-category of $\Perf(X)$ (resp. $\Coh(X)$).
We note that since $\Perf(X)$ (resp. $\Coh(X)$) are  idempotent complete, we have that $\IndPerf(X)^{\omega} = \Perf(X)$ (resp. $\IndCoh(X)^{\omega} = \Coh(X)$) form a set of compact generators.
        \end{enumerate}
        \end{definition}
        For $X$ an ind-algebraic stack, we recall that there is a natural action 
        \begin{equation}{\label{eqn: PerfActsonCoh}}
         \Perf(X) \times \Coh(X) \ra \Coh(X) 
        \end{equation}
        given by tensor product.
Ind-extending this functor defines for us a natural map
        \begin{equation}{\label{eqn: TensorProductofSheaves}}
         \IndPerf(X) \times \IndCoh(X) \ra \IndCoh(X). 
        \end{equation}
        
        We note that the decomposition in \Cref{eqn: decompositionintoconnectedcomponentsFullStack} (resp. \Cref{eqn: decompositionintoconnectedcomponents}) induces the following decomposition on these categories of coherent sheaves
        \begin{equation}{\label{eqn: FullDecompositionofIndCoherentSheaves}}
         \IndCoh(\Par_{G,E} \otimes \Lambda) \simeq \bigoplus_{\zeta} \IndCoh(\Par^{\zeta}_{G,E}) 
         \end{equation}
        (resp. 
        \begin{equation}{\label{eqn: TameDecompositionofIndCoherentSheaves}} 
        \IndCoh(\Par_{G,E}^{\tame} \otimes \Lambda) \simeq \bigoplus_{\zeta} \IndCoh(\Par_{G,E}^{\zeta}),
        \end{equation}
        for $\Lambda$ an algebraically closed field. 
        
        Here the direct sum runs over all inertial types (resp. tame inertial types).
Similarly, we recall that $Z(\hat{G})^{{\Gamma_E}}$ will act trivially on the space $\mathcal{Z}^{1}(W_{E},\hat{G})$.
In particular, the quotient $\Par_{G,E}$ of $\mathcal{Z}^{1}(W_{E},\hat{G})$ by $\hat{G}$ naturally has the structure of a $Z(\hat{G})^{{\Gamma_E}}$-gerbe over $[\mathcal{Z}^{1}(W_{E},\hat{G})/(\hat{G}/Z(\hat{G})^{\Gamma_E})]$.
As explained in \cite[Remark~2.74]{Zhu25} 
        and \cite[\S 3.2]{Zhu20}, this induces a decomposition 
        \begin{equation}{\label{eqn: ConnectedComponentsDecompositionSpectralSide}}
         \IndCoh(\Par_{G,E}) \simeq \bigoplus_{\nu \in \bb{X}^{*}(Z(\hat{G})^{{\Gamma_E}})} \IndCoh(\Par_{G,E})^{\nu}, 
        \end{equation}
        running over cocharacters (i.e., irreducible representations) of $Z(\hat{G})^{{\Gamma_E}}$.
This will refine the decompositions (\ref{eqn: FullDecompositionofIndCoherentSheaves}) and (\ref{eqn: TameDecompositionofIndCoherentSheaves}) after base-changing to an algebraically closed field, and we use the analogous notations for further projection to this direct summand on this subcategory.

        We will now study the automorphic incarnation of the direct sum decomposition (\ref{eqn: TameDecompositionofIndCoherentSheaves}).
        \subsubsection{The Schematic Tame Local Langlands Category}{\label{ss: TameLocalLanglandsCategory}}
        We will now introduce some variants of the schematic category $\Shv^{!}(\ICG)$\footnote{When $\Lambda$ is not a torsion algebra the construction of this in \cite{Zhu25} is very similar to \S \ref{Section: Sheaf and Cosheaf Theories}}.
We first introduce their representation theoretic analogues. 
        \begin{definition}[{\cite{Zhu25}}]{\label{defn: tamesubcategory}}
        For $H/E$ a connected reductive group, we define the following.
        \begin{enumerate}
        \item We write $\Rep^{\tame}(H(E),\Lambda) \subset \Rep(H(E),\Lambda)$ for the presentable stable category generated by $\on{c-Ind}_{\mathcal{G}(\mathcal{O}_{E})^{u}}^{H(E)}(\Lambda)$ under colimits, where $\mathcal{G}(\mathcal{O}_{E})^{u}$ is the pro-\(p\) radical of a parahoric subgroup $\mathcal{G}(\mathcal{O}_{E}) \subset H(E)$. 
        \item We define $\Rep_{\fg}(H(E),\Lambda) \subset \Rep(H(E),\Lambda)$ to be the idempotent complete  full subcategory generated by $\on{c-Ind}_{K}^{H}(\Lambda)$ for $K \subset H(E)$ any open compact subgroup. 
        \item We define $\Rep(H(E),\Lambda)^{\Adm} \subset \Rep(H(E),\Lambda)$ to be the idempotent complete (see \cite[Lemma~7.36]{Zhu25}) stable category of objects $A$ such that the invariants $A^{K} := \mathrm{RHom}(\on{c-Ind}_{K}^{H(E)}(\Lambda),A)$ is a perfect complex of $\Lambda$-modules for $K \subset H(E)$ any pro-$p$ compact open subgroup.
        \item We write $\Rep^{\widehat{\unip}}(H(E),\Lambda) \subset \Rep(H(E),\Lambda)$ for the presentable $\Lambda$-linear stable category generated by $\on{c-Ind}_{\mathcal{G}(\mathcal{O}_{E})}^{H(E)}(\pi)$, where $\pi$ is a representation of the Levi factor $L_{\mathcal{G}}$ of a parahoric $\mathcal{G}(\mathcal{O}_{E}) \subset H(E)$ which is a unipotent representation in the sense of \cite[Definition~4.86]{Zhu25}.
We write $\Rep_{\fg}^{\widehat{\unip}}(H(E),\Lambda) \subset \Rep_{\fg}(H(E),\Lambda)$ for the full idempotent complete subcategory generated by the compact inductions $\on{c-Ind}_{\mathcal{G}(\mathcal{O}_{E})}^{H(E)}(\pi)$ inside the category of finitely generated $H(E)$-representations.
        \end{enumerate}
        \end{definition} 
        
        We will now use these representation categories to define some variants of the schematic category $\Shv^{!}(\ICG)$.
We first start with the one coming from finitely generated representations, which we recall has an independent geometric meaning.

        For $H/E$ a connected reductive group, we recall (\Cref{prop:classifying_stack_is_sind_placid} (1)) that $\bbB \underline{H(E)}$ is a sind-placid stack, and therefore we have a natural map 
        \begin{equation}{\label{eqn: mapfromFinitelyGeneratedonClassifyingStack}}
        \Shv^{!}_{\fg}(\bbB \underline{H(E)}) \ra \Shv^{!}(\bbB \underline{H(E)}), 
        \end{equation}
        which a priori may not be fully faithful.
By \Cref{prop:classifying_stack_is_sind_placid} (3), we have an identification $\Shv^{!}(\bbB \underline{H(E)}) \simeq \Rep(H(E),\Lambda)$ with the category of smooth representations.
It then follows, by \cite[Proposition~3.57]{Zhu25}, that the map (\ref{eqn: mapfromFinitelyGeneratedonClassifyingStack}) is actually fully faithful and its essential image can be identified with $\Rep_{\fg}(H(E),\Lambda)$, as defined above. 
        
        Similarly, we can use this to see that the analogous map $\Shv^{!}_{\fg}(\ICG) \rightarrow \Shv^{!}(\ICG)$ for $\ICG$ is also fully faithful with essential image given by the smallest full stable idempotent complete subcategory generated by the $i_{b*}(A_{b})$ for $A_{b} \in \Rep_{\fg}(G_{b}(E),\Lambda)$ for $b \in B(G)$ varying (\cite[Proposition~3.96]{Zhu25}).
        
It follows by definition that we have a fully faithful embedding $\Rep(G_{b}(E),\Lambda)^{\omega} \hookrightarrow \Rep_{\fg}(G_{b}(E),\Lambda)$ and in turn we obtain a fully faithful embedding 
        \begin{equation}{\label{eqn: inclusionofcompactsintofinitelygenerated}}
         \Shv^{!}(\ICG)^{\omega} \hookrightarrow \Shv^{!}_{\fg}(\ICG) 
        \end{equation}
        of idempotent complete categories, which is in fact an equivalence if $\Lambda = \ol{\mathbb{Q}}_{\ell}$ by \cite[Corollary~3.58]{Zhu25}.
By Ind-extension, this gives rise to a fully faithful functor 
        \begin{equation}{\label{eqn: inclusionintoIndFinitelyGeneratedSheaves}}
         \Shv^{!}(\ICG) \hookrightarrow \IndShv^{!}(\ICG). 
        \end{equation}
        We now define the remaining important variants of $\Shv^{!}(\ICG)$. 
        \begin{definition}
        \begin{enumerate}
        \item We let $\Shv^{!,\tame}(\ICG) \subset \Shv^{!}(\ICG)$ be the full subcategory of objects $A \in \Shv^{!}(\ICG)$ satisfying the condition that $\ul{i_{b}^{!}}(A) \in \Rep^{\tame}(G_{b}(E),\Lambda)$ for all $b \in B(G)$. 
        \item We let $\Shv^{!}(\ICG)^{\Adm} \subset \Shv^{!}(\ICG)$ be the idempotent complete full subcategory of objects $A \in \Shv^{!}(\ICG)$ which satisfy the condition that $\ul{i_{b}^{!}}A \in \Rep(G_{b}(E),\Lambda)^{\Adm}$\footnote{We note that the notion of admissible objects makes sense in any dualizable category (see \cite[\S~7.2.3]{Zhu25}).
However, the definition here agrees with this abstract categorical definition in this particular case, by \cite[Theorem~1.3 (1)]{Zhu25}.}.
         \item We write $\Shv^{!,\widehat{\unip}}(\ICG)$ (resp. $\Rep_{\fg}^{\widehat{\unip}}(\ICG,\Lambda)$) for the full subcategory of objects $A \in \Shv^{!}(\ICG)$ satisfying the condition that $\ul{i_{b}^{!}}(A) \in \Rep^{\widehat{\unip}}(G_{b}(E),\Lambda)$ (resp. $\ul{i_{b}^{!}}(A) \in \Rep_{\fg}^{\widehat{\unip}}(G_{b}(E),\Lambda)$) for all $b \in B(G)$.
We write $\IndShv^{!,\widehat{\unip}}(\ICG) \subset \IndShv^{!}(\ICG)$ for its Ind-completion. 
        \end{enumerate}
        \end{definition}
        
        In light of \eqref{eqn: TameDecompositionofIndCoherentSheaves}, we expect that $\Shv^{!,\tame}(\ICG)$ should decompose as a direct sum over tame inertial types $\zeta = (w,\theta)$ if $\Lambda$ is an algebraically closed field. In order to make this expectation precise, we now describe the category that should correspond to the direct summand $\IndCoh(\Par_{G,E}^{\hat{\zeta}})$ under categorical Langlands.
        
        \emph{For the rest of the subsection, we assume that $G$ is tamely ramified  and that $\Lambda$ is algebraically closed.} We fix a choice of Iwahori $\mathcal{G}$ of $G$ and consider an element $w \in \widetilde{W}$ in the Iwahori-Weyl group.

We recall that this defines for us a locally closed subspace $i_{w}: \Sht_{\calG,w}^{\sch} \rightarrow \Sht_{\mathcal{G}}^{\sch}$ which we can postcompose with the Newton map $\Nt: \Sht_{\calG}^{\sch} \ra \ICG$ to obtain a map $\Nt_{w}: \Sht^{\sch}_{\mathcal{G},w} \ra \ICG$.
        
        We write $\mathcal{S}_{k}$ for the Levi factor of the Iwahori $\mathcal{G}$, which is a torus over $k$.
This is the dual torus of the torus $\hat{S} := \hat{T}/(1 - \ol{\tau})\hat{T}$ considered in \S \ref{ss: TameUnipotentStackofLanglandsParameters}.
We recall that this has an action by $W_{0} = W_G^{\ol{\tau}}$.
        We consider the image $\ol{w}$ of $w \in \widetilde{W}$ inside $W_{0}$ given by the reduction map.
Then we have an induced map
        \[ \phi_{w}: \mathcal{S}_{k} \ra \mathcal{S}_{k} \]
        \[ s \mapsto s^{-1}{\rm Ad}(\ol{w})(\sigma(s)), \]
        where $\sigma$ denotes the Frobenius on $k$. 
        We set $\mathcal{S}_{k}^{\ol{w}\sigma}$ to be the kernel of $\phi_{w}$, as in \cite[Equation~(4.37)]{Zhu25}. There is a natural map
        \[ \pr_{w}^{\sigma}: \Sht^{\sch}_{\mathcal{G},w} \ra \bb{B}\ul{S_{k}^{\ol{w}\sigma}}  \]
        towards the classifying stack attached to this finite group.
        In particular, if we are given a character $\theta: S_{k}^{\ol{w}\sigma} \ra \Lambda^\times$ of this finite group then we may regard it as a sheaf on $\bb{B}\ul{S_{k}^{\ol{w}\sigma}}$, 
This allows us, for $? \in \{!,*\}$, to define sheaves 
        \[ R_{w,\theta}^{?} := \ul{\Nt_{w*}}\ul{\pr_{w}^{\sigma?}}(\theta[-\ell(w)]) \in \Shv^{!}(\ICG), \]
        where $\ell: \widetilde{W} \ra \mathbb{N}_{\geq 0}$ denotes the natural length function.
In the case where $\theta$ is the trivial character we simply denote this by $R_{w}^{?}$.
     
        We can say a pair of elements $(w,\theta)$ and $(w',\theta')$ for $w,w' \in \widetilde{W}$ and $\theta$ (resp. $\theta'$) a pair of characters of $S_{k}^{\ol{w}\sigma}$ (resp. $S_{k}^{\ol{w}'\sigma}$) are \emph{geometrically conjugate} if the reductions $(\ol{w},\theta)$ and $(\ol{w}',\theta')$ are geometrically conjugate in the sense of Deligne-Lusztig.
We write $[(w,\theta)]$ for the equivalence class defined by taking $(w,\theta)$ up to geometric conjugacy.
It follows from \cite[Lemma~4.61]{Zhu25} that there is a natural correspondence between the conjugacy classes $[(w,\theta)]$ and tame inertial types $\zeta$ in the sense defined in \S \ref{ss: TameUnipotentStackofLanglandsParameters}.
This allows us to define the following. 
        \begin{definition}{\label{defn: categoriesassociatedtotypes}}
        For a tame inertial type $\zeta$ with associated geometric conjugacy class $[(w,\theta)]$ as described above, we define the subcategory
        \[ \Shv^{!,\hat{\zeta}}(\ICG) \subset \Shv^{!,\tame}(\ICG) \]
        to be the presentable stable full subcategory of $\Shv^{!,\tame}(\ICG)$ generated by the $R^{*}_{w',\theta'}$, for all $(w',\theta')$ lying in the equivalence class $[(w,\theta)]$.
        \end{definition}
        We now have the following consequence of \cite[Proposition~4.60]{Zhu25}, which provides the geometric mirror to the decomposition (\ref{eqn: TameDecompositionofIndCoherentSheaves}).  
        \begin{proposition}{\cite[Proposition~4.60]{Zhu25}}{\label{prop: directsumofIsocTame}}
        For $\Lambda$ an algebraically closed field, there is a decomposition
        \[ \bigoplus_{\zeta} \Shv^{!,\hat{\zeta}}(\ICG) \simeq \Shv^{!,\tame}(\ICG), \]
        where the direct sum ranges over tame inertial types.
        \end{proposition}
        \begin{remark}{\label{rem: DeligneLusztigTheory}}
        We note that, by the same argument as in \cite[Lemma~4.57]{Zhu25}, the stalk of $i_{b}^{!}R^{*}_{w,\theta}$ for all $b \in B(G)$ can be computed in terms of the $G_{b}(E)$-representation $C_{*}^{\mathrm{BM}}(X_{w}(b),\theta[-\ell(w)])$, the Borel-Moore homology of the affine Deligne-Lusztig variety attached to $b$ and $w$ with coefficients in the local system determined by $\theta$ (e.g $X_{w}(b) := \Sht^{\sch}_{\calG,w} \times_{\ICG} \ast$, where the fiber product is formed using $\ast \ra \ICG_{b} \ra \ICG$ and the Newton map $\Nt_{w}: \Sht_{\calG,w}^{\sch} \ra \ICG$).    
        \end{remark}
        This concludes our discussion of the schematic tame local Langlands category.
We can now state Zhu's version of the tame categorical local Langlands correspondence. 
        \subsubsection{Statement of the Schematic Tame Categorical Equivalence}{\label{ss: SchematicEquivalence}}
        With all the basic notations explained, we recall the salient aspects of the main theorem  of \cite{Zhu25}.
First, we recall that the usual statement of the tame categorical Langlands equivalence. As in all instances of the Langlands conjecture, this depends on some choice of auxiliary data.
Namely, for us the relevant notion will be the following.
Recall that we denote by $\bbF_q$ the residue field of $E$.
        \begin{definition}{\label{defn: IwahoriWhittaer}}
         We fix a choice $\psi: \bbF_q \ra \Lambda^\times$ of non-degenerate additive character, and a choice of Iwahori $\mathcal{I}/\mathcal{O}_{E}$.
We denote its pro-$p$ radical by $\mathcal{I}^{u}$ and write $U(\bbF_q)$ for its mod $\pi$-reduction.
Consider the composition
        \[ \psi_{e}: \mathcal{I}^{u} \ra U(\bbF_q) \xrightarrow{e} \bbF_q \xrightarrow{\psi} \Lambda^\times, \]
        where the first map is the map given by mod $\pi$-reduction and the second map is given by $U(\bbF_q) \ra U/[U,U](\bbF_q) \simeq \bigoplus_{i = 1}^{n} \bbF_q \xrightarrow{\mathrm{sum}} \bbF_q$.
With respect to the pair $(\mathcal{I},\psi)$, we define the Iwahori-Whittaker sheaf 
        \[ \on{c-Ind}_{\mathcal{I}^{u}}^{G(E)}(\psi_{e}) =: \IW_{\psi}. \]
        \end{definition}
        We fix such a choice from now on. Similarly, \emph{we fix a choice of square root of} $q$ \emph{in our coefficients} $\Lambda$ \emph{from now on} and implicitly pass between the $C$-group and $L$-group with respect to this choice, as in Remark \ref{rem: CgroupvsLgroup}. We have the following.
        \begin{theorem}{\label{thm: tamecategoricalLanglandsEquivalence}}{\cite[Theorem~1.6,1.7]{Zhu25}}
        Assume that $G/E$ is an unramified connected reductive group.
        \begin{enumerate}
        \item If $\Lambda = \ol{\mathbb{Q}}_{\ell}$, then there is a natural equivalence of categories 
        \[ \mathbb{L}_{G}^{\mathrm{tame}}:  \Shv^{!,\tame}(\ICG) \xrightarrow{\simeq} \IndCoh(\Par_{G,E}^{\tame} \otimes \Lambda), \]
        which, for each tame inertial type $\zeta$, restricts to an equivalence
        \[ \mathbb{L}_{G}^{\hat{\zeta}}:  \Shv^{!,\hat{\zeta}}(\ICG) \xrightarrow{\simeq} \IndCoh(\Par_{G,E}^{\hat{\zeta}}). \]
        The equivalence $\mathbb{L}_{G}^{\mathrm{tame}}$ is Ind-extended from an equivalence 
        \[ \mathbb{L}_{G}^{\mathrm{tame},\omega}:  \Shv^{!,\tame}(\ICG)^{\omega} \xrightarrow{\simeq} \Coh(\Par_{G,E}^{\tame} \otimes \Lambda), \]
        on compact generators.
We have that 
        \begin{equation}{\label{eqn: IwahoriWhittaker}}
         \bb{L}_{G}^{\tame}(\IW_{\psi}) \simeq \mathcal{O}_{\Par_{G,E}^{\tame} \otimes \Lambda}.
        \end{equation}
        \item If $\Lambda = \ol{\mathbb{F}}_{\ell}$ and $\ell$ satisfies \Cref{assump: modularcoefficients} below then there is a fully faithful embedding 
        \[ \mathbb{L}_{G,\fg}^{\widehat{\mathrm{unip}}}: \IndShv^{!,\widehat{\unip}}(\ICG) \hookrightarrow \IndCoh(\Par_{G,E}^{\widehat{\unip}} \otimes \Lambda), \]
        whose image is stable under the action of $\IndPerf(\Par_{G,E}^{\widehat{\unip}} \otimes \Lambda)$ described in (\ref{eqn: TensorProductofSheaves}).
Moreover, this is Ind-extended from a fully faithful embedding 
        \[ \bb{L}_{G,\fg}^{\widehat{\mathrm{unip}},\omega}: \Shv^{!,\widehat{\unip}}_{\fg}(\ICG) \hookrightarrow \Coh(\Par_{G,E}^{\widehat{\unip}} \otimes \Lambda) \]
        with image stable under the action of $\Perf(\Par_{G,E}^{\widehat{\unip}} \otimes \Lambda)$ described in (\ref{eqn: PerfActsonCoh}).
The functor $\mathbb{L}_{G,\fg}^{\widehat{\mathrm{unip}}}$ restricts to a fully faithful embedding
        \[ \mathbb{L}_{G}^{\widehat{\mathrm{unip}}}: \Shv^{!,\widehat{\unip}}(\ICG) \hookrightarrow \IndCoh(\Par_{G,E}^{\widehat{\unip}} \otimes\Lambda), \]
        whose image is stable under the action of $\Perf(\Par_{G,E}^{\widehat{\unip}} \otimes \Lambda)$. 
        \item Both the functor $\bb{L}_{G}^{\tame}$ and $\bb{L}_{G,\fg}^{\widehat{\unip}}$ match up the decomposition of (\ref{eqn: ConnectedComponentsDecompositionSpectralSide}) indexed by $\nu \in \bb{X}^{*}(Z(\hat{G})^{{\Gamma_E}})$ with the decomposition induced by $\ICG := \bigsqcup_{\nu \in \bb{X}_{*}(Z(\hat{G})^{{\Gamma_E}})} \ICG^{(\nu)}$ into the connected components up to a minus sign, where we recall that $\pi_{0}(|\ICG|) \simeq^{\kappa} \bb{X}^{*}(Z(\hat{G})^{{\Gamma_E}})$ via  \Cref{thm:isocG_geometry} (5) and \cite[\S~5]{Kottwitz}.
        \end{enumerate}
        \end{theorem}
        Here the assumption needed in the $\ell$-modular case is as follows. 
        \begin{assumption}{\label{assump: modularcoefficients}}
        We assume that the prime $\ell$ is greater than the Coxeter number of any simple factor of $G$ and $\ell \neq 19$ (resp. $\ell \neq 31$) if $G$ has a simple factor of type $E_{7}$ (resp. $E_{8}$).
        \end{assumption}
        
        We now transport this result to the analytic category via our functor $\pitch$.
To do this, we first introduce the analogous sheaf categories on $\Bun_{G}$.
        \subsubsection{The Analytic Tame Local Langlands Category}\label{sec:analytic_tame_LLC}
        We consider the category $\Dlis(\Bun_{G},\Lambda)$ of lisse-\'etale sheaves on $\Bun_{G}$, as defined in \cite[Definition~VII.6.1]{FS21}.
We recall that this is (non-obviously with non-torsion coefficients) equipped with a semi-orthogonal decomposition (see \cite[\S~3]{ImaConv} for details) such that the graded pieces identify with $\Dlis(\Bun_{G}^{b},\Lambda) \simeq \Rep(G_{b}(E),\Lambda)$, the unbounded derived category of smooth representations.
        
        In particular, for the inclusion $j_{b}: \Bun_{G}^{b} \hookrightarrow \Bun_{G}$, we have the full six operations, as well as an exceptional left adjoint $j_{b\sharp}$ to the $*$-pullback functor $j_{b}^{*}$ (by \cite[Proposition~VII.7.2]{FS21}).
Moreover, when $\Lambda$ is a torsion ring, we have an equivalence $\Dlis(\Bun_{G},\Lambda) \simeq \calD_{\Lambda}^{\an}(\Bun_{G})$ by \cite[Proposition~VII.6.6]{FS21} combined with \cite[Proposition~V.3.5]{FS21} and \cite[Theorem~V.3.7]{FS21}, and in this case it recovers the semi-orthogonal decomposition discussed in \S \ref{ss: theanalyticcategory}.
This allows us to define the following.
        \begin{definition}{\label{defn: analyticTameUnipotentCategory}}
        \begin{enumerate}
        \item We define $\Dlis^{\tame}(\Bun_{G},\Lambda) \subset \Dlis(\Bun_{G},\Lambda)$ to be the full subcategory of objects $A \in \calD(\Bun_{G},\Lambda)$ such that $j_{b}^{*}A \in \Rep^{\tame}(G_{b}(E),\Lambda)$. 
        \item We define $\calD_{\lisse,\fg}(\Bun_{G},\Lambda) \subset \calD_{\lisse}(\Bun_{G},\Lambda)$ to be the idempotent complete subcategory of objects $A$ such that $j_{b}^{*}(A) \in \Rep_{\fg}(G_{b}(E),\Lambda)$, and we write $\Ind\calD_{\lisse,\fg}(\Bun_{G},\Lambda)$ for its associated Ind-completion.
        \item We define $\calD^{\ULA}_{\lisse}(\Bun_{G},\Lambda) \subset \calD_{\lisse}(\Bun_{G},\Lambda)$ to be the idempotent complete subcategory of objects $A$ such that $j_{b}^{*}(A) \in \Rep(G_{b}(E),\Lambda)^{\Adm}$\footnote{The superscript $\mathrm{ULA}$ stands for universally locally acyclic in the sense of being ULA with respect to the structure morphism $\Bun_{G} \ra \ast$ (\cite[Definition~VII.7.7]{FS21}).
By Proposition \cite[Proposition~VII.7.9]{FS21}, this agrees with the definition given here.}. 
        \item For $G/E$ tamely ramified, if $\Lambda$ is an algebraically closed field and $\zeta$ a tame inertial type with associated geometric conjugacy class $[(w,\theta)]$ for $w \in \widetilde{W}$ and $\theta: \mathcal{S}_{k}^{\ol{w}\sigma} \ra \Lambda^{\times}$ as described in \Cref{defn: tameinertialtype}, we define $\Dlis^{\hat{\zeta}}(\Bun_{G},\Lambda)\subset \Dlis^{\tame}(\Bun_{G},\Lambda)$ to be the presentable stable subcategory generated by $\pitch(R^{*}_{w',\theta'})$ for all $(w',\theta') \in [(w,\theta)]$.
        \item We define $\Dlis^{\widehat{\unip}}(\Bun_{G},\Lambda) \subset \Dlis(\Bun_{G},\Lambda)$ (resp. $\calD^{\widehat{\unip}}_{\lis,\fg}(\Bun_{G},\Lambda)$) to be the full subcategory of objects satisfying the condition that $j_{b}^{*}(A) \in \calD^{\widehat{\unip}}(G_{b}(E),\Lambda)$ (resp. $j_{b}^{*}(A) \in \Rep_{\fg}^{\widehat{\unip}}(G_{b}(E),\Lambda)$).
We set $\Ind\calD_{\lis,\fg}^{\widehat{\unip}}(\Bun_{G},\Lambda) \subset \Ind\calD_{\lis,\fg}(\Bun_{G},\Lambda)$ to be its Ind-completion.
        \end{enumerate}
        \end{definition}
        We now deduce the following consequence of our main Theorem.
        \begin{corollary}{\label{cor: variantsoftheequivalence}}
        For $\Lambda/\bb{Z}_{\ell}$ a torsion ring, the following is true.
        \begin{enumerate}
        \item The functor $\pitch$ restricts to a $\Lambda$-linear equivalence \[ \pitch^{\tame}:  \Shv^{!,\mathrm{tame}}(\ICG,\Lambda)\xrightarrow{\simeq} \calD_{\lis}^{\mathrm{tame}}(\Bun_{G},\Lambda).\] 
        \item For $\Lambda$ an algebraically closed field and $\zeta$ a tame inertial type, the restriction $\pitch^{\tame}$ gives rise to a $\Lambda$-linear equivalence
        \[  \pitch^{\hat{\zeta}}: \Shv^{!,\hat{\zeta}}(\ICG,\Lambda) \xrightarrow{\simeq} \Dlis^{\hat{\zeta}}(\Bun_{G},\Lambda) \]
        such that $\pitch^{\tame} := \bigoplus_{\zeta} \pitch^{\hat{\zeta}}$.
In particular, by combining with Proposition \ref{prop: directsumofIsocTame}, we have a direct sum decomposition
        \begin{equation}\label{eqn: directsumdecompositionofDlisTame}
         \Dlis^{\tame}(\Bun_{G},\Lambda) \simeq \bigoplus_{\zeta} \Dlis^{\hat{\zeta}}(\Bun_{G},\Lambda). 
        \end{equation}
        \item  The functor $\pitch$ restricts to a $\Lambda$-linear equivalence 
        \[ \pitch_{\fg}: \Shv^{!}_{\fg}(\ICG,\Lambda) \xrightarrow{\simeq} \calD_{\lisse,\fg}(\Bun_{G},\Lambda),  \]
        which, by Ind-extension, gives rise to a $\Lambda$-linear equivalence
        \[ \Ind(\pitch_{\fg}): \IndShv^{!}(\ICG,\Lambda) \xrightarrow{\simeq}  \Ind\calD_{\lisse,\fg}(\Bun_{G},\Lambda). \]
        \item The functor $\pitch$ (resp. $\Ind(\pitch_{\fg})$) restricts to equivalences
        \[ \pitch^{\widehat{\unip}}: \Shv^{!,\widehat{\unip}}(\ICG,\Lambda) \xrightarrow{\simeq} \calD_{\lisse}^{\widehat{\unip}}(\ICG,\Lambda) \]
        (resp. 
        \[ \Ind(\pitch^{\widehat{\unip}}_{\fg}): \IndShv^{!,\widehat{\unip}}(\ICG,\Lambda) \xrightarrow{\simeq} \Ind\calD_{\lis,\fg}^{\widehat{\unip}}(\Bun_{G},\Lambda)) \]
        \item The functor $\pitch$ restricts to a $\Lambda$-linear equivalence
        \[ \pitch^{\omega}: \Shv^{!}(\ICG,\Lambda)^{\omega} \xrightarrow{\simeq} \calD_{\lis}(\Bun_{G},\Lambda)^{\omega} \]
        and 
        \[ \pitch^{\Adm}: \Shv^{!}(\ICG,\Lambda)^{\Adm} \xrightarrow{\simeq} \calD_{\lis}(\Bun_{G},\Lambda)^{\ULA}. \]
        \end{enumerate}
        \end{corollary}
        \begin{proof}
        Parts (1), (3) and (5) follow from the definitions and Theorem \ref{thm: MainTheoremPartial}. Parts (2) and (4) follow by the same logic.
However, as we do not know that Assumption \ref{assumpconj: linebundletwist} is true, we need to take account for the fact that the tensor action on $\Rep(G_{b}(E),\Lambda)$ by the character $\chi_{b} := \delta_{b} \otimes \delta_{b,\Zhu}^{-1}$ appearing in Theorem \ref{thm: MainTheoremPartial} (4) may not preserve the different subcategories used to define $\Shv^{!,\hat{\zeta}}(\ICG)$, $\Shv^{!,\widehat{\unip}}(\ICG)$, and their analogues on $\Bun_{G}$. However, this cannot happen in light of the following Lemma.
     \begin{lemma}{\label{lemma: RestrictiontoIwahoriTrivial}}
    Let $b \in B(G)$ and let $I_{b} \subset G_{b}(E)$ be an Iwahori subgroup.
Then $\chi_{b} := \delta_{b} \otimes \delta_{b,\Zhu}^{-1}$ is weakly unramified in the sense of \cite[\S~3.3.1]{Hai14}.
    \end{lemma}
    \begin{proof}
    This follows by combining \Cref{prop: deltabZhuisweaklyunramified} together with the fact that the character $\delta_{b}$ defined in \cite{HIDualCpx} satisfies the stronger condition of being unramified in the sense of \cite[\S~3.3.1]{Hai14}.
    \end{proof}
        \end{proof}
        \begin{remark}{\label{rem: DiscussionofSplitting}}
        We note that the existence of the direct sum decomposition (\ref{eqn: directsumdecompositionofDlisTame}) is already very interesting.
Indeed, we note that it follows by combining \Cref{thm: MainTheoremPartial} and \Cref{rem: DeligneLusztigTheory} that we have that 
        \[ j_{b}^{*}\pitch(R_{w,\theta}^{*}) \simeq C_{*}^{\mathrm{BM}}(X_{w}(b),\theta[-\ell(w)]). \]
        In particular, if one takes $w=w_b$, the $\sigma$-straight element attached to $b$, then, if $\calG$ is the Iwahori, we have an identification $\Sht_{\calG,w} \simeq \bb{B}\ul{I_{w_{b}}}$, as in \Cref{lemma: compactgeneratordiagram} (1). This implies that we have that $R_{w,\theta}^{*} \simeq \ul{i_{b*}}(\on{cInd}_{I_{w_{b}}}^{G_{b}(E)}(\theta))$ and in turn that $\pitch(R_{w,\theta}^{*}) \simeq j_{b!}(\on{cInd}_{I_{w_{b}}}^{G_{b}(E)}(\theta))$, by Theorem \ref{thm: MainTheoremPartial} and \Cref{lemma: RestrictiontoIwahoriTrivial}, where here we have inflated the character $\theta$ along the mod $\pi$-reduction map $I_{w_{b}} \ra S_{k}^{\ol{w_{b}}\sigma}$. We now look at two distinct geometric conjugacy classes $[(w,\theta)]$ and $[(w',\theta')]$ such that $[(w,\theta)] \neq [(w',\theta')]$ and consider two elements $b,b' \in B(G)$.
Suppose that $\pi_{b}$ (resp. $\pi_{b'}$) are smooth irreducible representations of $G_{b}(E)$ (resp. $G_{b'}(E)$) which occur as subquotients of a cohomology sheaf of $C_{*}^{\mathrm{BM}}(X_{w}(b),\theta)$ (resp. $C_{*}^{\mathrm{BM}}(X_{w'}(b'),\theta')$).
Then the direct sum decomposition (\ref{eqn: directsumdecompositionofDlisTame}) tells us that we have 
        \[ \RHom(j_{b!}(\pi_{b}),j_{b'!}(\pi_{b'})) = 0. \]
        In particular, we have that 
        \[\RHom(j_{b!}(\on{cInd}_{I_{w_{b}}}^{G_{b}(E)}(\theta)),j_{b'!}(\on{cInd}_{I_{w_{b'}}}^{G_{b'}(E)}(\theta'))) = 0. \]
        if $[(w_{b},\theta)] \neq [(w_{b'},\theta')]$. This would for example be an immediate consequence of showing that the Fargues-Scholze parameter of $j_{b!}(\pi_{b})$ and $j_{b'!}(\pi_{b'})$ are distinct (which would follow from a tame variant of \Cref{conj: IndPerfLinearity} below); however, this result shows this consequence immediately and can be made more explicit by combining with the calculations of the Borel-Moore homology of $X_{w'}(b')$ and $X_{w}(b)$ (cf. \cite[\S5]{He2014GeometricHomological}, resp. \cite{Ivanov_adlv,Zbarsky_09} for some special cases) and is completely unconditional.
       
        \end{remark}
        We now combine these variants of $\pitch$ with the schematic categorical Langlands equivalence described in \Cref{ss: SchematicEquivalence}.
        \subsection{The Analytic Tame Categorical Langlands Equivalence}
        We now want to combine some of our discussion with Theorem \ref{thm: tamecategoricalLanglandsEquivalence}.
In particular, we have the following immediate consequence of \Cref{cor: variantsoftheequivalence} and \Cref{thm: tamecategoricalLanglandsEquivalence}. 
        \begin{theorem}{\label{thm: AnalyticEquivalence}}
        For $G/E$ unramified, if $\Lambda = \ol{\bb{F}}_{\ell}$ and $\ell$ satisfies \Cref{assump: modularcoefficients} then we set $\bb{L}_{G,\fg}^{\an,\widehat{\unip}} := \bb{L}_{G,\fg}^{\widehat{\unip}} \circ \Ind(\pitch)_{\fg}^{-1}$, where $\Ind(\pitch_{\fg})$ is the equivalence of Corollary \ref{cor: variantsoftheequivalence} (4).
This induces a fully faithful embedding
        \begin{equation}{\label{eqn: embeddingIndFinitelyGeneratedSheaves}}
         \bb{L}_{G,\fg}^{\an,\widehat{\unip}}: \Ind\calD^{\widehat{\unip}}_{\lisse,\fg}(\Bun_{G},\Lambda) \hookrightarrow \IndCoh(\Par_{G,E}^{\widehat{\unip}} \otimes \Lambda)   
        \end{equation}
        with image stable under the action of $\Ind\Perf(\Par_{G,E}^{\widehat{\unip}} \otimes \Lambda)$.
More specifically, there is a fully faithful embedding  
        \[ \calD^{\widehat{\unip}}_{\lisse,\fg}(\Bun_{G},\Lambda) \hookrightarrow \Coh(\Par_{G,E}^{\widehat{\unip}} \otimes \Lambda)   \]
        with image stable under the action of $\Perf(\Par_{G,E}^{\widehat{\unip}} \otimes \Lambda)$ and the functor (\ref{eqn: embeddingIndFinitelyGeneratedSheaves}) is obtained by Ind-extending.
Moreover, the restriction of $\bb{L}_{G,\fg}^{\an,\widehat{\unip}}$ to the full subcategory $\calD^{\widehat{\unip}}_{\lisse}(\Bun_{G},\Lambda) \subset \Ind\calD^{\widehat{\unip}}_{\lisse,\fg}(\Bun_{G},\Lambda)$ is stable under the action of $\Perf(\Par_{G,E}^{\widehat{\unip}} \otimes \Lambda)$.
We denote this restriction by $\bb{L}_{G}^{\an,\widehat{\unip}}$. 
        \end{theorem}
        
        \emph{We assume that $\ell$ satisfies \Cref{assump: modularcoefficients} for the rest of the subsection.}
We will now draw the readers attention to one of the most important pieces of structure that comes from \Cref{thm: AnalyticEquivalence}.
Namely, we note that the tensor product action of $\Perf(\Par_{G,E}^{\widehat{\unip}})$ on the image of $\calD_{\lis}^{\widehat{\unip}}(\Bun_{G},\Lambda)$ under $\bb{L}_{G}^{\an,\widehat{\unip}}$ induces a spectral action.
Namely, a monoidal functor
        \begin{align}{\label{eqn: ZhuSpectralAction}}
        \Perf(\Par_{G,E}^{\widehat{\unip}} \otimes \ol{\bb{F}}_{\ell}) & \ra \End(\Dlis^{\widehat{\unip}}(\Bun_{G},\ol{\bb{F}}_{\ell})) \\ 
        A & \mapsto (- \mapsto A \star_{\mathrm{Zhu}} -). \nonumber
        \end{align}
        Concretely, given $V \in \Rep(\hat{G})$ an algebraic representation of the dual group of $G$, the underlying $\hat{G}$-representation gives rise to a vector bundle on $B\hat{G}$ which pulls back along the natural map $\Par_{G,E}^{\widehat{\unip}} \ra B\hat{G}$ to a vector bundle $C_{V}^{\widehat{\unip}}$ on $\Par_{G,E}^{\widehat{\unip}}$, and the endomorphism 
$C_{V}^{\widehat{\unip}} \star -$ on $\Dlis^{\widehat{\unip}}(\Bun_{G},\ol{\bb{F}}_{\ell})$ coming from the spectral action is given by the convolution action of the central sheaf attached to $V$ on $\Shv^{!}(\mathcal{I}\backslash LG/\mathcal{I})$ under the affine Deligne-Lusztig induction functor $\mathrm{Ch}^{\unip}_{\phantom{}^{L}G,\phi}$ after passing through the equivalence $\pitch$ (see the proof of \cite[Proposition~3.12]{YangZhuTorsion} for details).  
        
        On the other hand, assuming that $\ell \nmid \lvert\pi_{0}(Z(G))\rvert$, the category $\Dlis(\Bun_{G},\ol{\bb{F}}_{\ell})$ is also equipped with an action $\Perf(\Par_{G,E})$, of perfect complexes on the full stack of $L$-parameters by \cite[Theorem~I.10.1]{FS21}.
Concretely, a vector bundle $V \in \Rep(\hat{G})$ defines an object of $\Perf(B\hat{G})$, which, by pulling back along the natural morphism $\Par_{G,E} \otimes \ol{\bb{F}}_{\ell} \ra B\hat{G}$ to the classifying stack of $\hat{G}$ over $\ol{\bb{F}}_{\ell}$, defines a vector bundle $C_{V} \in \Perf(\Par_{G,E} \otimes \ol{\bb{F}}_{\ell})$. Then the action of $C_{V} \in \Perf(\Par_{G,E})$ on $\Dlis(\Bun_{G},\ol{\bb{F}}_{\ell})$ is given by the endofunctor 
        \begin{equation}{\label{eqn: HeckeOperatorofV}}
        T_{V}: \Dlis(\Bun_{G},\ol{\bb{F}}_{\ell}) \ra \Dlis(\Bun_{G},\ol{\bb{F}}_{\ell}) 
        \end{equation}
        defined by geometric Satake for the $\mathbb{B}_{\dR}^{+}$-grassmannian and the Hecke correspondence attached to $V$, as described in \cite[Chapter~IX]{FS21}.
This Hecke action upgrades using \cite[Theorem I.10.1.]{FS21} to a monoidal functor 
        \begin{align}{\label{eqn: FSSpectralAction}}
        \Perf(\Par_{G,E} \otimes \ol{\bb{F}}_{\ell}) & \ra \End(\Dlis(\Bun_{G},\ol{\bb{F}}_{\ell})) \\ 
        A & \mapsto (- \mapsto A \star_{\mathrm{FS}} -), \nonumber
        \end{align}
        which sends the canonical vector bundle $C_{V}$ defined by pulling back along $\Par_{G,E} \otimes \ol{\bb{F}}_{\ell} \ra B\hat{G}$ to the endofunctor $T_{V}$.
We now have what we consider to be one of the most important conjectures concerning our functor $\pitch$.
        \begin{conjecture}{\label{conj: IndPerfLinearity}}
        If $\ell$ satisfies \Cref{assump: modularcoefficients} and $\ell \nmid \lvert\pi_{0}(Z(G))\rvert$ and $G/E$ is unramified then the following should hold.
        \begin{enumerate}
        \item The spectral action functor of Fargues-Scholze defined by (\ref{eqn: FSSpectralAction}) preserves the full subcategory (in fact, direct summand by (\ref{eqn: directsumdecompositionofDlisTame})) $\Dlis^{\widehat{\unip}}(\Bun_{G},\ol{\bb{F}}_{\ell}) \subset \Dlis(\Bun_{G},\ol{\bb{F}}_{\ell})$. 
        \item The monoidal functor induced by (\ref{eqn: FSSpectralAction}) and point (1) 
        \[ \Perf(\Par_{G,E} \otimes \ol{\bb{F}}_{\ell}) \ra \End(\Dlis^{\widehat{\unip}}(\Bun_{G},\ol{\bb{F}}_{\ell})) \]
        factors through the projection map $\Perf(\Par_{G,E} \otimes \ol{\bb{F}}_{\ell}) \ra \Perf(\Par_{G,E}^{\widehat{\unip}} \otimes \ol{\bb{F}}_{\ell})$ to the direct summand, as in (\ref{eqn: FullDecompositionofIndCoherentSheaves}). 
        \item The functor 
        \[ \Perf(\Par_{G,E}^{\widehat{\unip}} \otimes \ol{\bb{F}}_{\ell}) \ra \End(\Dlis^{\widehat{\unip}}(\Bun_{G},\ol{\bb{F}}_{\ell})) \]
        induced by point (1) and (2) is naturally equivalent to the spectral action functor of Zhu described in (\ref{eqn: ZhuSpectralAction}) as monoidal functors. In other words, the functor $\pitch^{\widehat{\unip}}$ is linear with respect to the spectral actions of (\ref{eqn: ZhuSpectralAction}) and points (1)-(2).
        \end{enumerate}
        \end{conjecture}
        \begin{remark}
        Point (1) is essentially saying that the Hecke operators on $\Dlis(\Bun_{G},\ol{\bb{F}}_{\ell})$ preserve the unipotent subcategory. Given (1), it should be possible to verify (2). Indeed, it essentially reduces to showing that any unipotent representation is in fact given by the semi-simplification of a unipotent $L$-parameter (via the relationship between the spectral action of Fargues-Scholze and the excursion algebra (see \cite[Theorem~5.2.1]{ZouCatLanglandsforTorii})). For the $*$-stalks of the generators $\underline{i_{b*}}(\mathrm{cInd}_{I_{b}}^{G_{b}(E)}(\Lambda))$ where $G_{b}$ is quasi-split, this  can be routinely verified using compatability of the Fargues-Scholze correspondence with parabolic induction and local Langlands for tori \cite[Theorem~I.9.6. (i),(viii)]{FS21}. However, one also needs to account for the generators $\underline{i_{b*}}(\mathrm{cInd}_{I_{b}}^{G_{b}(E)}(\Lambda))$ where $G_{b}$ is not quasi-split as well as the case of the generators coming from arbitrary parahorics of $G_{b}$ and general unipotent representations on their Levi factors. Using (1), this should be approachable by combining \cite[Theorem~1.0.2]{KW} with analogous arguments using compatability of Fargues-Scholze with parabolic induction for general parabolics.
For a more detailed discussion of how to approach (1) and (2) using geometric Eisenstein series, see Remark \ref{rem: PreservationofUnipotentSubcategoryUnderHeckearguments}.
Point (3) tells us that there is a precise relationship between the convolution action of central sheaves described above under the affine Deligne-Lusztig induction functor $\mathrm{Ch}^{\unip}_{\phantom{}^{L}G,\phi}$ and the Hecke operators of Fargues-Scholze (see \cite[\S~3.1.4]{YangZhuTorsion} and the proof of \cite[Proposition~3.12]{YangZhuTorsion}) which is the deepest part of the conjecture.
        \end{remark}
        \begin{remark}
        When $\Lambda = \ol{\bb{Q}}_{\ell}$, one should of course be able to formulate analogues for the tame subcategory, as well as the categories attached to each tame inertial type.
However, this is not within the scope of this paper, as the functor $\pitch$ has not been properly defined with rational coefficients. 
We leave it as a (very hard) exercise to the reader to determine the precise formulation of Conjecture \ref{conj: IndPerfLinearity} in this level of generality.
        \end{remark}
         \begin{remark}{\label{rem: QCohLinearityIndFinitelyGenerated}}
        With $\ell$-modular coefficients, it would be interesting to consider a variant of this statement where one works with $\calD_{\lis,\fg}^{\widehat{\unip}}(\Bun_{G},\Lambda)$ together with its spectral action induced by Theorem \ref{thm: AnalyticEquivalence}.
However, here there is an additional subtlety.
Namely, we do not know that the spectral action of Fargues-Scholze preserves the full subcategory $\calD_{\lis,\fg}(\Bun_{G},\Lambda) \subset \calD_{\lis}(\Bun_{G},\Lambda)$ of finitely generated sheaves beyond the case where one has an equality $\calD_{\lis,\fg}(\Bun_{G},\Lambda) = \calD_{\lis}(\Bun_{G},\Lambda)^{\omega}$ (e.g if $\ell$ is banal with respect to $G_{b}$ for all $b \in B(G)$) in which case it follows from \cite[Theorem~IX.2.2]{FS21}. In particular, we have the following precise question, which seems interesting in modular characteristic.  
        \begin{question}{\label{question: HeckeOperatorsPreserveFinitelyGenerated}}
        For $V \in \Rep_{\Lambda}(\phantom{}^{L}G^{I})$, does the associated Hecke operator 
        \[ T_{V}: \calD_{\lis}(\Bun_{G},\Lambda) \ra \calD_{\lis}(\Bun_{G},\Lambda) \]
        on $\Bun_{G}$ constructed in \cite[\S~IX.2]{FS21}, where we have forgotten the Weil group action, preserve the full subcategory $\calD_{\lis,\fg}(\Bun_{G},\Lambda) \subset \calD_{\lis}(\Bun_{G},\Lambda)$ of finitely-generated sheaves, as in Definition \ref{defn: analyticTameUnipotentCategory} (2)? 
        \end{question}
        \end{remark}

        We will want to discuss an application of this conjecture to the cohomology of local Shimura varieties at the end of this paper, by combining it with the work of Yang-Zhu \cite{YangZhuTorsion} on the perverse $t$-exactness of the spectral action of (\ref{eqn: ZhuSpectralAction}).
To explain this properly, we will need to localize the statement of Conjecture \ref{conj: IndPerfLinearity} around different parameters.
For the rest of the section, \emph{we assume that $G/E$  is unramified}, and consider the natural map
        \[ q_{E}^{\widehat{\unip}}: \Par_{G,E}^{\widehat{\unip}} \ra \Par_{G,E}^{\coarse,\widehat{\unip}} \]
        given by restricting the map (\ref{eqn: maptocoarsemodulispace}) to the unipotent locus, where $\Par_{G,E}^{\coarse,\widehat{\unip}} = \Spec(\mathcal{Z}_{\hat{G}}^{\widehat{\unip}})$ for $\mathcal{Z}_{\hat{G}}^{\widehat{\unip}} := H^{0}(\Par_{G,E}^{\widehat{\unip}},\mathcal{O}_{\Par_{G,E}^{\widehat{\unip}}})$.
The closed $\Lambda$-points of $\Par_{G,E}^{\coarse,\widehat{\unip}}$ correspond to conjugacy classes of unramified $L$-parameters $\phi$. We abusively write $\phi$ to denote the associated closed point of $\Par_{G,E}^{\coarse,\widehat{\unip}}$.

We want now pass to the formal completion around the closed substack $(q_{E}^{\widehat{\unip}})^{-1}(\phi) \hookrightarrow \Par_{G,E}^{\widehat{\unip}}$, which we denote by $V_{\hat{\phi}}$.
We recall, as described in \cite[\S~3.2.3]{YangZhuTorsion}, that there is a fully faithful embedding
        \[ i_{\hat{\phi}}: \IndCoh(V_{\hat{\phi}}) \hookrightarrow \IndCoh(\Par_{G,E}^{\widehat{\unip}} \otimes \ol{\bb{F}}_{\ell}),  \]
        whose essential image we denote by $\IndCoh(\Par_{G,E}^{\widehat{\unip}} \otimes \ol{\bb{F}}_{\ell})_{\hat{\phi}}$.
By projection formula, the full subcategory $\IndCoh(\Par_{G,E}^{\widehat{\unip}} \otimes \ol{\bb{F}}_{\ell})_{\hat{\phi}}$ is stable under the action of $\Perf(\Par_{G,E}^{\widehat{\unip}} \otimes \ol{\bb{F}}_{\ell})$. 

        Similarly, on the geometric side, we may consider $\Shv^{!,\widehat{\unip}}(\ICG,\Lambda)$ and the action by $\Perf(\Par_{G,E}^{\widehat{\unip}})$  induced by the embedding $\bb{L}_{G}^{\unip}$ of Theorem \ref{thm: tamecategoricalLanglandsEquivalence}.
The action of $\mathcal{Z}_{\hat{G}}^{\widehat{\unip}}$ via endomorphisms of the category $\Perf(\Par_{G,E}^{\widehat{\unip}})$ by scaling of global sections induces an action of $\QCoh(\mathcal{Z}_{\hat{G}}^{\widehat{\unip}})$ on the category $\Shv^{!,\widehat{\unip}}(\ICG,\Lambda)$. We can then form the Lurie tensor product 
        \begin{equation}{\label{eqn: hatphicompletedsheavesonIsoc}}
         \Shv^{!,\widehat{\unip}}(\ICG,\ol{\bb{F}}_{\ell})_{\hat{\phi}} :=  \Shv^{!,\widehat{\unip}}(\ICG,\ol{\bb{F}}_{\ell}) \otimes_{\QCoh(\Spec(\mathcal{Z}_{\hat{G}}^{\widehat{\unip}}))} \QCoh(\Spec(\mathcal{Z}_{\hat{G}}^{\widehat{\unip}}))_{\hat{\mf{\phi}}}, 
        \end{equation}
        where $\QCoh(\Spec(\mathcal{Z}_{\hat{G}}^{\widehat{\unip}}))_{\hat{\phi}} \hookrightarrow \QCoh(\Spec(\mathcal{Z}_{\hat{G}}^{\widehat{\unip}}))$ is the full subcategory corresponding to the formal completion $\Spec(\mathcal{Z}_{\hat{G}}^{\widehat{\unip}})$ around the closed point corresponding to $\phi$.
The natural map $\Shv^{!,\widehat{\unip}}(\ICG,\ol{\bb{F}}_{\ell})_{\hat{\phi}} \ra \Shv^{!,\widehat{\unip}}(\ICG,\ol{\bb{F}}_{\ell})$ defines the inclusion of a full subcategory (see \cite[\S~3.2.3]{YangZhuTorsion} for more details), which is stable under the action of $\Perf(\Par_{G,E}^{\widehat{\unip}})$ by \cite[Lemma~3.24]{YangZhuTorsion}.
Moreover, the natural map 
        \[ i_{\hat{\phi}}^{!}: \IndCoh(\Par_{G,E}^{\widehat{\unip}} \otimes \ol{\mathbb{F}}_{\ell}) \ra \IndCoh(V_{\hat{\phi}}) \]
        induces a right adjoint functor
        \begin{equation}{\label{eqn: phicompletedlocalization}}
         (-)_{\hat{\phi}}: \Shv^{!,\widehat{\unip}}(\ICG,\ol{\bb{F}}_{\ell}) \ra \Shv^{!,\widehat{\unip}}(\ICG,\ol{\bb{F}}_{\ell})_{\hat{\phi}}.
        \end{equation}
        which we refer to as the $\phi$-completed localization functor.
        
        We now have the following variant of the schematic unipotent equivalence.
        \begin{theorem}{\label{thm: ParameterZhuResult}}{\cite[Proposition~3.12]{YangZhuTorsion}}
        For $G/E$ unramified, $\phi$ a conjugacy class of unramified $L$-parameters, and $\ell$ satisfying Assumption \ref{assump: modularcoefficients}, the functor $\bb{L}_{G}^{\widehat{\unip}}$ induces a fully faithful embedding
        \[ \bb{L}_{G,\hat{\phi}}^{\widehat{\unip}}: \Shv^{!,\widehat{\unip}}(\ICG,\ol{\bb{F}}_{\ell})_{\hat{\phi}} \hookrightarrow \IndCoh(V_{\hat{\phi}}) \]
        with image stable under the action of $\Perf(\Par_{G,E}^{\widehat{\unip}} \otimes \ol{\bb{F}}_{\ell})$. 
        \end{theorem}
        We now would like to transport this to $\calD_{\lisse}(\Bun_{G},\ol{\bb{F}}_{\ell})$.
We assume that $\ell \nmid \lvert\pi_{0}(Z(G))\rvert$. Then the spectral action of $\Perf(\Par_{G,E})$ (\ref{eqn: FSSpectralAction}) induces an action of $\mathcal{Z}_{\hat{G}} := H^{0}(\Par_{G,E},\calO)$ on $\calD_{\lisse}(\Bun_{G},\ol{\bb{F}}_{\ell})$ and in turn an action of $\QCoh(\Spec(\mathcal{Z}_{\hat{G}}))$. 

        Given any semisimple $L$-parameter $\phi$, we may form the analogue of (\ref{eqn: hatphicompletedsheavesonIsoc}); namely, we have that 
        \[ \calD_{\lisse}(\Bun_{G},\ol{\bb{F}}_{\ell})_{\widehat{\phi}} := \calD_{\lisse}(\Bun_{G},\ol{\bb{F}}_{\ell}) \otimes_{\QCoh(\Spec(\mathcal{Z}_{\hat{G}}))} \QCoh(\Spec(\mathcal{Z}_{\hat{G}}))_{\hat{\mf{\phi}}},  \]
        which will be a full subcategory of $\calD_{\lisse}(\Bun_{G},\ol{\bb{F}}_{\ell})$, is easily checked to be preserved under the action of $\Perf(\Par_{G,E} \otimes \ol{\bb{F}}_{\ell})$ described in (\ref{eqn: FSSpectralAction}). As before, $\QCoh(\Spec(\mathcal{Z}_{\hat{G}}))_{\hat{\mf{\phi}}} \hookrightarrow \QCoh(\Spec(\mathcal{Z}_{\hat{G}}))$ is the full sub-category corresponding to the formal completion around the closed point corresponding to $\phi$.
We refer to this subcategory as the set of $\hat{\phi}$-completed sheaves.
For the canonical vector bundles $C_{V}$ attached to $V \in \Rep(\hat{G})$, we denote the induced operator by $T_{V,\hat{\phi}}$ and refer to it as the $\hat{\phi}$-completed Hecke operator.
        \begin{remark}{\label{rem: philocalvsphicomplete}}
        As is explained in \cite[Remark~3.28]{YangZhuTorsion}, there is also the category $\calD_{\lisse}(\Bun_{G},\ol{\bb{F}}_{\ell})_{\phi}$, which was constructed in \cite[Appendix~A]{HamannLeeTorsion}.
This corresponds to simply the localization (i.e., the stalk) over the closed point in $\Spec(\calZ_{\hat{G}})$ corresponding to $\phi$ instead of passing to the formal completion.
There is a natural map
        \[ \calD_{\lisse}(\Bun_{G},\ol{\bb{F}}_{\ell})_{\phi} \ra \calD_{\lisse}(\Bun_{G},\ol{\bb{F}}_{\ell})_{\hat{\phi}} \]
        corresponding to formal completion, which will induce an equivalence after passing to the subcategories of ULA objects.
        \end{remark}
        We now have the following formal consequence of \Cref{conj: IndPerfLinearity}.
        \begin{proposition}{\label{prop: parameterwisepitch}}
        Let $G/E$ be unramified. Assume \Cref{conj: IndPerfLinearity} and Assumption \ref{assump: modularcoefficients} on $\ell$. Then, for all conjugacy classes of unramified $L$-parameters $\phi$, the functor $\pitch^{\widehat{\unip}}$ restricts to an equivalence of categories 
        \[  \pitch_{\hat{\phi}}: \Shv^{!,\widehat{\unip}}(\ICG,\ol{\bb{F}}_{\ell})_{\hat{\phi}} \xrightarrow{\simeq} \calD_{\lis}^{\widehat{\unip}}(\Bun_{G},\ol{\bb{F}}_{\ell})_{\hat{\phi}}, \]
        which is linear (as in \Cref{conj: IndPerfLinearity} (3)) with respect to the spectral action of $\Perf(\Par^{\widehat{\unip}}_{G,E} \otimes \ol{\bb{F}}_{\ell})$ given by Theorem \ref{thm: ParameterZhuResult} of Zhu and to the spectral action of $\Perf(\Par_{G,E}^{\widehat{\unip}} \otimes \ol{\bb{F}}_{\ell})$ of Fargues-Scholze given by Conjecture \ref{conj: IndPerfLinearity} (1)-(2) and tensoring by $\otimes_{\QCoh(\Spec(\mathcal{Z}_{\hat{G}}))} \QCoh(\Spec(\mathcal{Z}_{\hat{G}}))_{\hat{\mf{\phi}}}$
        \end{proposition}

        We now explain another unconditional application of our main theorem pertaining to these localizations.
        \subsection{Splitting of the semi-orthogonal decomposition for generic parameters}
        We assume that $\Lambda = \ol{\bb{F}}_{\ell}$ for the rest of this subsection.
We consider an unramified semi-simple $L$-parameter $\phi: W_E \ra \phantom{}^{L}G(\Lambda)$ and the localization
        \[ \Shv^{!,\widehat{\unip}}(\ICG,\Lambda)_{\hat{\phi}} \hookrightarrow \Shv^{!,\widehat{\unip}}(\ICG,\Lambda) \hookrightarrow \Shv(\ICG,\Lambda) \]
        introduced in the previous section. 
        
        As explained in \cite[\S~3.2.3]{YangZhuTorsion}, the semi-orthogonal decomposition on $\Shv^{!}(\ICG,\Lambda)$ (in the sense of \S \ref{subsection semi on isoc}) restricts to a semi-orthogonal decomposition on the full subcategory $\Shv^{!,\widehat{\unip}}(\ICG,\Lambda)_{\hat{\phi}}$ with graded denoted by $\Rep^{\widehat{\unip}}(G_{b}(E),\Lambda)_{\hat{\phi}}$.
We  write 
        \[ i_{b,\hat{\phi}!},i_{b,\hat{\phi}*}: \Rep^{\widehat{\unip}}(G_{b}(E),\Lambda)_{\hat{\phi}} \ra \Shv^{!,\widehat{\unip}}(\ICG,\Lambda)_{\hat{\phi}}\]
        for the associated $!$- and $*$-pushforward, and 
        \[ i_{b,\hat{\phi}}^{!},i_{b,\hat{\phi}}^{*}: \Shv^{!,\widehat{\unip}}(\ICG,\Lambda) \ra \Rep^{\widehat{\unip}}(G_{b}(E),\Lambda)_{\hat{\phi}} \]
        for their right adjoint functors.
These functors are the restrictions of the functors $i_{b!},i_{b*},i_{b}^{!},i_{b}^{*}$, respectively.
        
        We have the set of unramified elements $B(G)_{\un} := \mathrm{Im}(B(T) \ra B(G))$.
In this section, we will be interested in the following result of Yang-Zhu.
        \begin{proposition}{\cite[Corollary~3.32,Proposition~3.36]{YangZhuTorsion}}{\label{prop: splitSemiOrthogonalIsoc}}
        Let $G/E$ be unramified and $\ell$ be a prime satisfying Assumption \ref{assump: modularcoefficients}.
Consider an unramified $L$-parameter  $\phi: W_E \ra \phantom{}^{L}G(\ol{\bb{F}}_{\ell})$, which is of Langlands-Shahidi type in the terminology of \cite[Definition~6.2]{HamannLeeTorsion} or equivalently strongly generic in the sense of \cite[Definition~1.1]{YangZhuTorsion}.
Then the following is true. 
        \begin{enumerate}
        \item For all $b \in B(G)_{\mathrm{un}} \setminus B(G)$, the functors 
        \[ i_{b,\hat{\phi}}^{!},i_{b,\hat{\phi}}^{*}: \Shv^{!,\widehat{\unip}}(\ICG,\Lambda)_{\hat{\phi}} \ra \Rep^{\widehat{\unip}}(G_{b}(E),\Lambda)_{\hat{\phi}} \]
        vanish.
        \item For all $b \in B(G)_{\mathrm{un}}$, the natural transformations 
        \[ i_{b,\hat{\phi}}^{*} \ra i_{b,\hat{\phi}}^{!} \]
        and
        \[ i_{b,\hat{\phi}!} \ra i_{b,\hat{\phi}*} \]
        are equivalences.
        \item In particular, by (1) and (2), the semi-orthogonal decomposition on $\Shv^{!,\widehat{\unip}}(\ICG,\ol{\bb{F}}_{\ell})_{\hat{\phi}}$ splits, and we have an induced direct sum decomposition
        \[ \Shv^{!,\widehat{\unip}}(\ICG,\Lambda)_{\hat{\phi}} \simeq \bigoplus_{b \in B(G)_{\un}} \Rep^{\widehat{\unip}}(G_{b}(E),\Lambda)_{\hat{\phi}}. \] 
        \end{enumerate}
        \end{proposition}
        We now define $\calD^{\widehat{\unip}}_{\lisse}(\Bun_{G},\Lambda)'_{\hat{\phi}} := \pitch(\Shv^{!,\widehat{\unip}}(\ICG,\Lambda)_{\hat{\phi}})$.
\begin{remark}{\label{rem: primesuperscript}}
We note we have used the superscript $(-)'$, as we do not know that this agrees with $\calD^{\widehat{\unip}}_{\lisse}(\Bun_{G},\Lambda)_{\hat{\phi}}$, as defined in the previous section using the Fargues-Scholze spectral action. Indeed, this would follow from  \Cref{conj: IndPerfLinearity} (cf. \Cref{prop: parameterwisepitch}).
\end{remark}
The category $\calD^{\widehat{\unip}}_{\lisse}(\Bun_{G},\Lambda)'_{\hat{\phi}}$ is also a full subcategory, and, using Theorem \ref{thm: MainTheoremPartial}, we deduce that the semi-orthogonal decomposition described in \S \ref{sec:semiorthogonalDecOfDBunG} restricts to a semi-orthogonal decomposition on $\calD^{\widehat{\unip}}_{\lisse}(\Bun_{G},\Lambda)'_{\hat{\phi}}$ with graded isomorphic to $\Rep^{\widehat{\unip}}(G_{b}(E),\Lambda)_{\hat{\phi}} \otimes \chi_{b}  \hookrightarrow \Rep^{\widehat{\unip}}(G_{b}(E),\Lambda) \otimes \chi_{b} = \Rep^{\widehat{\unip}}(G_{b}(E),\Lambda)$, where we recall that this is not the natural semi-orthogonal decomposition coming from $\Bun_{G}$ as an object in $\AnStk_{v}$ equipped with the $\calD_{\Lambda}^{\an}$ $6$-functor formalism.
Here for the last equality we have used Lemma \ref{lemma: RestrictiontoIwahoriTrivial}.
In particular, we have functors
        \[ j_{b,\hat{\phi}\sharp},j_{b,\hat{\phi}!}: \Rep^{\widehat{\unip}}(G_{b}(E),\Lambda)_{\hat{\phi}} \otimes \chi_{b}  \ra \Dlis^{\widehat{\unip}}(\Bun_{G},\Lambda)'_{\hat{\phi}}\]
        with right adjoint functors
        \[ j_{b,\hat{\phi}}^{*},j_{b,\hat{\phi}}^{\flat}: \calD_{\lis}^{\widehat{\unip}}(\Bun_{G},\Lambda)_{\hat{\phi}}' \ra \Rep^{\widehat{\unip}}(G_{b}(E),\Lambda)_{\hat{\phi}} \otimes \chi_{b}.  \]
These are the restrictions of the functors $j_{b\sharp}$, $j_{b!}$, $j_{b}^{*}$, and $j_{b}^{\flat}$ to the subcategory $\calD^{\widehat{\unip}}_{\lisse}(\Bun_{G},\ol{\bb{F}}_{\ell})'_{\hat{\phi}}$.
In particular, we may deduce the following consequence of \Cref{thm: MainTheoremPartial} and \Cref{prop: splitSemiOrthogonalIsoc}. 
        \begin{corollary}{\label{cor: SplitSemiOrthogonalDecomposition}}
        Let $G$ be unramified and consider an unramified $L$-parameter  $\phi: W_E \ra \phantom{}^{L}G(\Lambda)$, which is of Langlands-Shahidi type in the terminology of \cite[Definition~6.2]{HamannLeeTorsion} or equivalently strongly generic in these sense of \cite[Definition~1.1]{YangZhuTorsion} and $\ell$ satisfies assumption \ref{assump: modularcoefficients}.
Then the following is true. 
        \begin{enumerate}
        \item For all $b \in B(G)\setminus B(G)_{\mathrm{un}} $, the functors 
        \[ j_{b,\hat{\phi}}^{\flat},j_{b,\hat{\phi}}^{*}: \calD_{\lis}^{\widehat{\unip}}(\Bun_{G},\Lambda)_{\hat{\phi}}' \ra \Rep(G_{b}(E),\Lambda)_{\hat{\phi}} \otimes \chi_{b} \]
        vanish.
        \item For all $b \in B(G)_{\mathrm{un}}$, the natural transformations 
        \[ j_{b,\hat{\phi}}^{\flat} \ra j_{b,\hat{\phi}}^{*} \]
        and
        \[ j_{b,\hat{\phi}\sharp} \ra j_{b,\hat{\phi}!} \]
        are equivalences. 
        \item In particular, by (1) and (2), the exceptional semi-orthogonal decomposition on $\Dlis(\Bun_{G},\Lambda)_{\hat{\phi}}'$ induced by \S \ref{sec:semiorthogonalDecOfDBunG} splits, and we have an induced direct sum decomposition
        \[ \calD_{\lisse}(\Bun_{G},\Lambda)_{\hat{\phi}}' \simeq \bigoplus_{b \in B(G)_{\un}} \Rep^{\widehat{\unip}}(G_{b}(E),\Lambda)_{\hat{\phi}} \otimes \chi_{b}. \] 
        \end{enumerate}
        \end{corollary} 
        This is an analogue of \cite[Corollary~4.29]{HamannLeeTorsion}, using Remark \ref{rem: philocalvsphicomplete}.
However, in order to make this precise, one actually needs to describe the structure of the categories, which, since $\Rep^{\widehat{\unip}}(G_{b}(E),\Lambda)_{\hat{\phi}}$ was defined in terms of the embedding provided by Theorem \ref{thm: tamecategoricalLanglandsEquivalence}, is in turn defined in terms of taking categorical traces.
However, now one is faced with the problem of explicating these categories $\Rep^{\widehat{\unip}}(G_{b}(E),\Lambda)_{\hat{\phi}} \otimes \chi_{b}$.
To this aim, we recall that, for all $b \in B(G)_{\mathrm{un}}$, there is an associated Borel $B_{b} \subset J_{b}$ as in \cite[\S~4.2.1]{HamannLeeTorsion}, which we assume is standard with respect to our fixed choice of Borel. 
        \begin{lemma}{(Proof of \cite[Corollary~3.35]{YangZhuTorsion})}
        Let $\phi: W_E \ra \phantom{}^{L}G(\Lambda)$ be an unramified parameter induced from a toral parameter $\phi_{T}: W_E \ra \phantom{}^{L}T(\Lambda)$.
For fixed $b \in B(G)_{\mathrm{un}}$, suppose that $\pi_{b} \in \Rep^{\widehat{\unip}}(G_{b}(E),\Lambda)_{\hat{\phi}}$ is a smooth irreducible representation.
Then it is an unramified principal series.
More precisely, it is a subquotient of the normalized parabolic induction $i_{B_{b}}^{J_{b}}(|\cdot|^{s})$, where $|\cdot|^{s}: T(\bb{Q}_{p}) \ra \Lambda^{\times}$ denotes an unramified character.
        \end{lemma}
        In particular, this tells us that such a $\pi_{b} \in \Rep^{\widehat{\unip}}(G_{b}(E),\Lambda)_{\hat{\phi}}$ should be given by a subquotient of certain parabolic induction $i_{B_{b}}^{G_{b}}(|\cdot|^{s})$ for some $s$.
The specific unramified twist appearing here should of course be completely specified by the parameter $\hat{\phi}$.
However, it is unclear (to us) which precise unramified twist should be occurring here, at least using results in the existing literature.
In particular, we have the following question/assumption, which is motivated by the description of the analogous $\phi$-localized category in \cite[Proposition~4.6]{HamannLeeTorsion}. 
        \begin{assumption/question}{\label{assumption/conj: UnramifiedTwistsaatIwahoriLevel}}
For $b \in B(G)$, we let $M_{b} \subset G$ be the standard Levi determined by the centralizer of the slope homomorphism and set $W_{b} := W_{G}/W_{M_{b}}$ to be the associated cosets of Weyl group. We fix a choice of representatives of minimal length in $W_{b}$.
For all $b \in B(G)$, do we have that 
        \[ \Red_{b}(\phi)_{w} := i_{B_{b}}^{J_{b}}(\chi^{w}) \otimes \delta_{b}^{-\frac{1}{2}}[-d_{b}] \in \Rep^{\widehat{\unip}}(G_{b}(E))_{\hat{\phi}} \otimes \chi_{b}? \]
        Here $\chi$ denotes the character attached to the toral parameter $\phi_{T}$ inducing $\phi$ under the geometric normalization of class field theory. We recall that $d_{b} := \langle 2\rho_{G},\nu_{b} \rangle$, where $\rho_{G}$ is the half sum of all positive roots.
        \end{assumption/question}
        \begin{remark}
        We note that, for the stratum corresponding to the trivial element, this should reduce to compatibility between the usual action of the spherical Hecke algebra on $\on{c-Ind}_{K}^{G}(\Lambda)$ and the Hecke operators on $\Shv^{!,\widehat{\unip}}(\ICG,\Lambda)$ which are induced by the endofunctors of $\Shv^{!,\widehat{\unip}}(\ICG,\Lambda)$ coming from convolution and geometric Satake via taking categorical trace (See \cite[Propositions~4.42,4.43]{YangZhuTorsion})
Similarly, for general $b \in B(G)$, one has to show a similar compatibility between the action  of Hecke operators on the neutral stratum of $\Shv^{!,\widehat{\unip}}(\mathfrak{B}(J_{b}),\Lambda)$ coming from geometric Satake and the one induced on the strata of $\Shv^{!,\widehat{\unip}}(\ICG,\Lambda)$. This should ultimately be done via some geometric constant terms on the affine Grassmannian.
We recall that these geometric constant terms naturally carry Tate twists and these should give rise to the modulus character twists appearing in Assumption \ref{assumption/conj: UnramifiedTwistsaatIwahoriLevel}, as in \cite[\S~IX.7.1]{FS21}. 
        
        More conceptually, this would follow for example from Proposition \ref{prop: parameterwisepitch} and the fact that we can compute the Fargues-Scholze-parameter of the smooth irreducible representations ocuring in $j_{b!}(i_{B_{b}}^{J_{b}}(\chi^{w})[-d_{b}] \otimes \delta_{b}^{-\frac{1}{2}})$ directly and show that it always agrees with $\phi$. 
        \end{remark}
        Admitting this for now, we deduce the following corollary. 
        \begin{corollary}{\label{cor: semiorthogonaldecompositionsplits}}
        Assume that Question \ref{assumption/conj: UnramifiedTwistsaatIwahoriLevel} has a positive answer and that Assumption \ref{assumpconj: linebundletwist} is true. Then for all $b \in B(G)_{\mathrm{un}}$ the natural transformations 
        \[ j_{b\sharp} \ra j_{b!}, \]
        and 
        \[ j_{b}^{\flat} \ra j_{b}^{*} \]
        are isomorphisms when evaluated on the representations $\Red_{b}(\phi)_{w}$ for all $w \in W_{b}$.
        \end{corollary}
        \begin{proof}
        This follows from \Cref{cor: SplitSemiOrthogonalDecomposition} and \Cref{assumption/conj: UnramifiedTwistsaatIwahoriLevel}, where we note that $\chi_{b}$ is trivial in light of Assumption \ref{assumpconj: linebundletwist}. 
        \end{proof}
        \begin{remark}{\label{rem: splittingofSemiorthogonalDecomposition}}
        We note that Corollary \ref{cor: semiorthogonaldecompositionsplits} is very much related to \cite[Theorem~8.1]{HamGeomES}, by using the identification 
        \[ \mathrm{nEis}_{B!}^{b_{T}}(\chi) := j_{b\sharp}(i_{\ol{B}_{b}}^{J_{b}}(\chi^{w}) \otimes \delta_{b}^{-1/2}[-d_{b}]) = j^{\ren}_{b\sharp}(i_{\ol{B}_{b}}^{J_{b}}(\chi^{w})). \]
        for $b_{T}$ a dominant reduction of $b \in B(G)_{\mathrm{un}}$ to $B(T) \simeq \bb{X}^{*}(\hat{T})$, which follows from combining \cite[Corollary~2.2.5 (5)]{HHS} with \cite[Theorem~4.28]{HIDualCpx} and arguing similarly to the beginning of the proof of \cite[Theorem~4.28]{HIDualCpx}.
Here $\mathrm{nEis}_{B}$ denotes the normalized geometric Eisenstein functor (see \Cref{ss: EisensteinFunctors} below).
Indeed, this result generalizes part of \cite[Theorem~8.1]{HamGeomES} to hold for a general unramified $G$ and an unramified toral parameter $\phi_{T}$ under the assumptions on $\ell$.
However, \cite[Theorem~8.1]{HamGeomES} also allows for general toral parameters $\phi_{T}$ of Langlands-Shahidi type (L.S.-type), but is contingent on compatibility of Fargues-Scholze with endoscopic classifications.
        \end{remark}
        As one might expect from the previous remark, our functor $\pitch$ should explain a connection between the tame categorical Langlands of Zhu and the geometric Eisenstein series functors studied in \cite{HamGeomES,HIDualCpx,HHS}.
We formulate this compatability precisely now and check it in some examples.
        \subsection{Compatibility of $\pitch$ with Eisenstein Functors}{\label{ss: EisensteinFunctors}}
        We assume $G/E$ is unramified throughout this section.
We now formulate a compatibility of the functor $\bb{L}_{G}^{\an,\widehat{\unip}}$ introduced in Theorem \ref{thm: tamecategoricalLanglandsEquivalence} with the geometric Eisenstein functors introduced in \cite{HamGeomES,HHS,HIDualCpx}.
We will assume that $\Lambda = \ol{\bb{F}}_{\ell}$ and also implicitly pass through the identification $\calD_{\lis}(X,\Lambda) \simeq \calD_{\Lambda}^{\an}(X)$ for suitably nice Artin $v$-stacks using \cite[Proposition~VII.6.6.]{FS21}. We fix a parabolic $P \subset G$ with Levi factor $M$ and induced diagram $M \leftarrow P \rightarrow G$.
This allows to define the normalized Eisenstein functor 
        \begin{equation}{\label{eqn: GeometricEisensteinFunctor}}
         \mathrm{nEis}_{P}(-): \Dlis(\Bun_{M},\Lambda) \ra \Dlis(\Bun_{G},\Lambda), 
        \end{equation}
        as in \cite[Definition~2.1.7]{HHS}, which is defined via the correspondence $\Bun_{M} \leftarrow \Bun_{P} \rightarrow \Bun_{G}$ and tensoring via a Verdier self-dual invertible sheaf, denoted $\mathrm{IC}_{\Bun_{P}}$, defined in terms of a square root\footnote{We fix this choice of square root to be compatible with the one fixed in \S \ref{ss: SchematicEquivalence}} of the modulus character of $P$.
Analogously, on the spectral side, we have the correspondence 
        \[ \Par_{M} \xleftarrow{\mf{q}^{\spec}} \Par_{P} \xrightarrow{\mf{p}^{\spec}} \Par_{G}, \]
        and we get a well-defined functor 
        \begin{equation}{\label{eqn: SpectralEisensteinFunctor}}
        \mathrm{nEis}_{P}^{\spec}(-) := \mf{p}_{\spec*}\mf{q}^{\spec*}(-): \Coh(\Par_{M}) \ra \Coh(\Par_{G}) 
        \end{equation}
        by \cite[Proposition~2.3.9]{Zhu20}\footnote{Strictly speaking, this result is for $\mf{q}^{\spec!}$ and not $\mf{q}^{*}$.
However, this does not make a meaningful difference by virtue of the fact that the map $\mf{q}$ is quasi-smooth (see the remarks after \cite[Conjecture~1.4.7]{HansenBeijingNotes})}.
We may then Ind-extend this functor to get a well defined functor 
        \[ \mathrm{nEis}_{P}^{\spec}(-): \IndCoh(\Par_{M}) \ra \IndCoh(\Par_{G}). \]
        We expect the functors (\ref{eqn: GeometricEisensteinFunctor}) and (\ref{eqn: SpectralEisensteinFunctor}) to be compatible under the conjectural categorical equivalence (see \cite[Conjecture~1.4.7]{HansenBeijingNotes}).
We note however that the functors $\mathrm{nEis}_{P}(-)$ are currently only describable in terms of the category $\Dlis(\Bun_{G},\Lambda)$ and not in terms of the category $\Shv(\ICG,\Lambda)$.
However, we can resolve this issue using Theorem \ref{thm: AnalyticEquivalence}.
In particular, we have the following.
        \begin{conjecture}{\label{conj: CompatabilityofOurFunctorWithEisensteinSeries}}
        Let $G/E$ be unramified and suppose that $\Lambda = \ol{\bb{F}}_{\ell}$ for $\ell$ satisfying assumption \ref{assump: modularcoefficients}.
Let $P \subset G$ be a parabolic with Levi factor $M$. 
        \begin{enumerate}
        \item We have an inclusion $\mathrm{nEis}_{P}(\calD_{\lis}^{\widehat{\unip}}(\Bun_{M},\Lambda)) \subset \calD_{\lis}^{\widehat{\unip}}(\Bun_{G},\Lambda)$.
        \item Assuming Part (1), the induced diagram
        \[
        \begin{tikzcd} 
         \calD_{\lisse}^{\widehat{\unip}}(\Bun_{G}, \Lambda) \arrow[r,"\bb{L}_{G}^{\an,\widehat{\unip}}"] &  \IndCoh(\Par^{\widehat{\unip}}_{G,E} \otimes \Lambda) & \\
        \calD_{\lisse}^{\widehat{\unip}}(\Bun_{M},\Lambda) \arrow[r,"\bb{L}_{M}^{\an,\widehat{\unip}}"] \arrow[u,"\mathrm{nEis}_{\ol{P}}"] &  \IndCoh(\Par_{M}^{\widehat{\unip}} \otimes \Lambda). \arrow[u,"\mathrm{Eis}_{P}^{\spec}"] &
        \end{tikzcd}
        \]
        commutes.
Here $\ol{P}$ denotes the opposite parabolic (cf. \cite[Conjecture~3.2 (ii)]{hellmann2021derived}).
        \end{enumerate}
        \end{conjecture}
        \begin{remark}{\label{rem: Eisensteinfunctorspreserveuniptoentsubcategroy}}
        The first part of the conjecture is not difficult to verify, so the real content is the second part.
Indeed, by virtue of the fact that the functor $\mathrm{nEis}_{P}$ commute with colimits (since it is a left adjoint), one may reduce to checking the claim on a set of compact generators of $\Dlis^{\widehat{\unip}}(\Bun_{M},\Lambda)$, which will be the form of $i_{b!}(A_{b})$, where $A_{b} \in \Rep^{\widehat{\unip}}(G_{b}(E),\Lambda)$ is a compact object (see the discussion after \cite[Definition 3.3]{YangZhuTorsion}). From here, one can argue as in the proof of \cite[Claim~3.14]{HHS}, where it reduces to the analogous claim for parabolic induction, and here it is easily verified (e.g using \cite[Corollary~3.6]{DatFinitudepourlesRepresentationsLissesdeGroupesPadiques}). 
        \end{remark}
        \begin{remark}{\label{rem: IndFinitelyGeneratedEisensteinSeries}}
        Analogous to Question \ref{question: HeckeOperatorsPreserveFinitelyGenerated}, one can ask whether the Eisenstein functors have the property that $\mathrm{nEis}_{P}(\calD_{\lis,\fg}(\Bun_{M},\Lambda)) \subset \calD_{\lis,\fg}(\Bun_{G},\Lambda)$.
This would then allow one to define a functor $\mathrm{nEis}_{P,\fg}: \Ind\calD_{\lis,\fg}(\Bun_{M},\Lambda) \ra \Ind\calD_{\lis,\fg}(\Bun_{G},\Lambda)$ by Ind-extending, and one can formulate a similar compatability with $\bb{L}_{G,\fg}^{\an,\widehat{\unip}}$. Unlike the case of Hecke operators, this is true.
In particular, one can apply the exact argument of \cite[Claim~3.14]{HHS} to formally reduce this to the claim that parabolic induction preserves finitely generated representations, which is true in any coefficient system over $\bb{Z}[\frac{1}{p}]$ by \cite[Corollary~1.5]{DHKMFinitenessforHeckeAlgebras}. We have chosen not to formulate things in this extra level of generality for simplicity.
        \end{remark}
        \begin{remark}{\label{rem: PreservationofUnipotentSubcategoryUnderHeckearguments}}
        We note that if one combines \Cref{conj: CompatabilityofOurFunctorWithEisensteinSeries} (1) (which should be fairly routine as described in Remark \ref{rem: Eisensteinfunctorspreserveuniptoentsubcategroy}) with \cite[Conjecture 1.5.2]{HansenBeijingNotes} that one can reduce checking \ref{conj: IndPerfLinearity} (1)-(2) for the Eisenstein part of $\calD_{\lis}^{\widehat{\unip}}(\Bun_{G},\ol{\bb{F}}_{\ell})$ (as in the semi-orthogonal decomposition of \cite[Theorem~1.3.2 (2)]{HHS}) to the Levi subgroup $M$.
In particular, this kind of analysis reduces one to checking \Cref{conj: IndPerfLinearity} (1)-(2) for the cuspidal part of $\calD_{\lis}^{\widehat{\unip}}(\Bun_{G},\ol{\bb{F}}_{\ell})$ (which in particular should be trivial, as there are no unipotent representations with supercuspidal $L$-parameter). The conjecture  \cite[Conjecture 1.5.2]{HansenBeijingNotes} is a theorem when $M = T$ is the maximal torus by \cite[Theorem~1.7]{HamGeomES}, and the general case will appear in forthcoming work of the second author (H.) with Hansen and Scholze.
Similar reduction steps should also apply to Conjecture \ref{conj: IndPerfLinearity} (3). In particular, one should be able to reduce it to the cuspidal part of $\calD_{\lis}^{\widehat{\unip}}(\Bun_{G},\ol{\bb{F}}_{\ell})$; however, this is much more involved.
        \end{remark}
        We now check this conjecture in some easy special cases.
To this aim, we recall that both functors admit a grading indexed by $\bb{X}^{*}(Z(\hat{G})^{{\Gamma_E}}) \simeq B(G)_{\basic}$. In particular, the Eisenstein functor $\mathrm{nEis}_{P}$ admits a decomposition into functors $\mathrm{nEis}_{P}^{\nu}$ indexed by cocharacters $\nu \in \bb{X}^{*}(Z(\hat{M})^{{\Gamma_E}})$ coming from taking the preimage of the connected component in $\Bun_{M}$ corresponding to $\nu \in \bb{X}^{*}(Z(\hat{M})^{{\Gamma_E}}) \simeq^{\kappa} B(G)_{\basic}$ via \Cref{thm: BunGGeometricFacts} (2). 
        Similarly, the spectral Eisenstein functors admit a grading 
        \[ \mathrm{nEis}_{P}^{\spec}(-) = \bigoplus_{\nu \in \bb{X}^{*}(Z(\hat{M})^{{\Gamma_E}})} \mathrm{nEis}_{P}^{\nu,\spec}(-) \]
        induced by the decomposition (\ref{eqn: ConnectedComponentsDecompositionSpectralSide}).
These two gradings should match under the categorical conjecture, as in Theorem \ref{thm: tamecategoricalLanglandsEquivalence} (4).
In particular, we now have the following slight refinement of Conjecture \ref{conj: CompatabilityofOurFunctorWithEisensteinSeries}, which is equivalent to Conjecture \ref{conj: CompatabilityofOurFunctorWithEisensteinSeries} using Theorem \ref{thm: tamecategoricalLanglandsEquivalence} (4) and the fact that our functor $\pitch$ clearly matches up the decomposition into connected components of $\Bun_{G}$ and $\ICG$ in an obvious way is as follows.
         \begin{conjecture}{\label{conj: CompatabilityofOurFunctorWithEisensteinSeriesConnectedComponents}}
        Suppose that $G$ is unramified, $\Lambda = \ol{\bb{F}}_{\ell}$, and $\ell$ satisfies assumption \ref{assump: modularcoefficients}.
Then assuming part (1) of \ref{conj: CompatabilityofOurFunctorWithEisensteinSeries} and combining it with \ref{thm: tamecategoricalLanglandsEquivalence} (4), the induced diagram
        \[
        \begin{tikzcd} 
         \Dlis^{\widehat{\unip}}(\Bun_{G}, \Lambda) \arrow[r,"\bb{L}_{G}^{\an,\widehat{\unip}}"] &  \IndCoh(\Par_{G,E} \otimes \Lambda) & \\
        \Dlis^{\widehat{\unip}}(\Bun_{M},\Lambda) \arrow[r,"\bb{L}_{M}^{\an,\widehat{\unip}}"] \arrow[u,"\mathrm{nEis}_{\ol{P}}^{-\nu}"] &  \IndCoh(\Par_{M,E} \otimes \Lambda) \arrow[u,"\mathrm{Eis}_{P}^{\nu,\mathrm{spec}}"].  &
        \end{tikzcd}
        \]
        commutes, for all $\nu \in \bb{X}^{*}(Z(\hat{M}))$. 
        \end{conjecture}
        We verify this conjecture in some simple cases.
To this aim, we recall that we have a natural map $\Par_{T,E} \ra B\hat{T}$ and pullback induces a natural map 
        \[ \Rep(\hat{T}) \ra \Perf(\Par_{T,E}), \]
        as in the definition of the canonical vector bundles in \S \ref{sec:analytic_tame_LLC}.
Therefore, for each $\lambda \in \bb{X}^{*}(\hat{T})$, we obtain a bundle which we denote by $\mathcal{O}_{\Par_{\hat{T}}}(\lambda)$, which will lie in $\IndCoh(\Par_{T,E})^{\lambda}$ under the direct sum decomposition described in (\ref{eqn: ConnectedComponentsDecompositionSpectralSide}).
For a general sheaf  $\mathcal{L} \in \IndCoh(\Par_{T,E})$, we write $\mathcal{L}(\lambda) := \mathcal{L} \otimes_{\mathcal{O}_{\Par_{\hat{T},E}}} \mathcal{O}_{\Par_{\hat{T},E}}(\lambda)$.
For each $\lambda \in \bb{X}^{*}(\hat{T})$, we write $b_{\lambda}^{T} \in B(T)$ for the associated element under the isomorphism $B(T) \simeq^{\kappa_{T}} \bb{X}^{*}(\hat{T}^{{\Gamma_E}})$ induced by the $\kappa$-invariant.
We now have the following.
        \begin{proposition}{\label{prop: EisensteinCompatability}}
        Suppose that $\ell$ does not divide the pro-order of any parahoric subgroup $K \subset G_{b}(E)$ for all $b \in B(G)$ (so that $\IndShv^{!,\widehat{\unip}}(\ICG,\Lambda) = \Shv^{!,\widehat{\unip}}(\ICG,\Lambda)$ (cf.
Remark \ref{rem: IndFinitelyGeneratedEisensteinSeries})) and that $G$ is unramified over $E$. Conjecture \ref{conj: CompatabilityofOurFunctorWithEisensteinSeries} (2) is true when $M = T$, and we evaluate on the compact objects $i_{b_{-\lambda}^{T}!}(\on{c-Ind}_{T(\mathcal{O}_{E})}^{T(E)}(\Lambda))$ of $\Shv^{!,\widehat{\unip}}(\frakB(T),\Lambda)$ for $\lambda$ either dominant or anti-dominant with respect to the Borel $B$.
        In particular, we have an isomorphism 
        \[ \bb{L}_{G}^{\an,\widehat{\unip}}(\mathrm{nEis}_{\ol{B}}(i_{b_{-\lambda}^{T}!}(\on{c-Ind}_{T(\mathcal{O}_{E})}^{T(E)}(\Lambda))))) \simeq \Eis_{B}^{\spec}\bb{L}_{T}^{\an,\widehat{\unip}}(\mathcal{O}(\lambda))  \]
        in $\IndCoh(\Par_{G,E}^{\widehat{\unip}}) \subset \IndCoh(\Par_{G,E})$ for $\lambda$ dominant or anti-dominant with respect to $B$.
        \end{proposition}
        \begin{proof}
        We assume without loss of generality that the Borel $B$ is the one we fixed.
        For an element $\lambda \in \bb{X}^{*}(\hat{T})$, we write $\ol{\lambda}$ for the induced element in $\bb{X}^{*}(\hat{T}^{{\Gamma_E}})$.
We let $b_{\lambda,T}$ be the corresponding element in $B(T) \simeq^{\kappa_{T}} \bb{X}^{*}(\hat{T}^{{\Gamma_E}})$ associated with $\ol{\lambda}$, and write $b_{\lambda} \in B(G)_{\mathrm{un}}$ for the image under the map $B(T) \ra B(G)$.
We write $I_{b_{\lambda}} \subset G_{b_{\lambda}}(E)$ for the Iwahori attached to our fixed choice of $\mathcal{I}/\mathcal{O}_{E}$, as in Convention \ref{conv: sigmastraightelements}. We let $w_{0} \in W_{G}$ denote the element of longest length and assume that $\lambda \in \bb{X}^{*}(\hat{T})$ is dominant with respect to the Borel $B$. We have the following, by \cite[Corollary~5.7]{Zhu25}.
        \begin{enumerate}
        \item There is an isomorphism
        \[ \bb{L}_{T}^{\widehat{\unip}}i_{b_{-\lambda}^{T}!}(\on{c-Ind}_{T(\mathcal{O}_{E})}^{T(E)}(\Lambda)) \simeq \mathcal{O}_{\Par_{\widehat{T}}^{\widehat{\unip}}}(\lambda), \]
        recalling the minus sign appearing in Theorem \ref{thm: tamecategoricalLanglandsEquivalence} (4).
        \item There is an isomorphism
        \[ \bb{L}_{G}^{\widehat{\unip}}(i_{b_{-\lambda}!}(\on{c-Ind}_{I_{b_{-\lambda}}}^{G_{b_{-\lambda}}(E)}(\Lambda)[-d_{b}])) \simeq \Eis_{B}^{\spec}(\mathcal{O}_{\Par_{\widehat{T}}^{\widehat{\unip}}}(\lambda)) = \Eis_{B}^{\spec,\lambda}(\mathcal{O}_{\Par_{\widehat{T}}^{\widehat{\unip}}}(\lambda)).   \]
        \item There is an isomorphism
        \[ \bb{L}_{G}^{\widehat{\unip}}(i_{b_{-\lambda}*}(\on{c-Ind}_{\ol{I}_{b_{-\lambda}}}^{G_{b_{-\lambda}}(E)}(\Lambda)[-d_{b}])) \simeq \Eis_{B}^{\spec}(\mathcal{O}_{\Par_{\widehat{T}}^{\widehat{\unip}}}(w_{0}(\lambda))) = \Eis_{B}^{\spec,w_{0}(\lambda)}(\mathcal{O}_{\Par_{\widehat{T}}^{\widehat{\unip}}}(w_{0}(\lambda))). \]
        Here $\ol{I}_{b}$ is the opposite Iwahori, and we recall that we have normalized our identifications of $\Rep(G_{b}(E),\Lambda) \simeq \Shv^{!}(\ICG_{b},\Lambda)$ with respect to the fixed choice of Borel $B$. 
        \end{enumerate}
        We also know, by Theorem \ref{thm: MainTheoremPartial} and Lemma \ref{lemma: RestrictiontoIwahoriTrivial}, that the following is true for all $\lambda \in \bb{X}^{*}(T)$.
        \begin{enumerate}
        \item There is an isomorphism 
        \[ \pitch_{T}(i_{b^{T}_{\lambda}!}(\on{c-Ind}_{T(\mathcal{O}_{E})}^{T(E)}(\Lambda))) = j_{b^{T}_{\lambda!}}(\on{c-Ind}_{T(\mathcal{O}_{E})}^{T(E)}(\Lambda)), \]
        \item There is an isomorphism
        \[ \pitch_{G}(i_{b_{\lambda}!}(\on{c-Ind}_{I_{b_{\lambda}}}^{G_{b_{\lambda}}(E)}(\Lambda)))[-d_{b}]) \simeq j_{b_{\lambda}\sharp}(\on{c-Ind}_{I_{b_{\lambda}}}^{G_{b_{\lambda}}(E)}(\Lambda)[-d_{b}]). \]
        \item There is an isomorphism
        \[ \pitch_{G}(i_{b_{\lambda}*}(\on{c-Ind}_{\ol{I}_{b_{\lambda}}}^{G_{b_{\lambda}}(E)}(\Lambda)))[-d_{b}]) \simeq j_{b_{\lambda}!}(\on{c-Ind}_{\ol{I}_{b_{\lambda}}}^{G_{b_{\lambda}}(E)}(\Lambda)[-d_{b}]). \] 
        \end{enumerate}
        In particular, for $\lambda$ dominant with respect to $B$, this reduces us to showing that we have an isomorphism
        \[ \Eis_{\ol{B}}(i_{b^{T}_{-\lambda}!}(\on{c-Ind}_{T(\mathcal{O}_{E})}^{T(E)}(\Lambda))) \simeq j_{b_{-\lambda}\sharp}(\on{c-Ind}_{I_{b_{-\lambda}}}^{G_{b_{-\lambda}}(E)}(\Lambda)[-d_{b}]), \]
        and an isomorphism:
          \[ \mathrm{nEis}_{\ol{B}}(i_{b^{T}_{-w_{0}(\lambda)}!}(\on{c-Ind}_{T(\mathcal{O}_{E})}^{T(E)}(\Lambda))) = j_{b_{-\lambda}!}(\on{c-Ind}_{\ol{I}_{b_{-\lambda}}}^{G_{b_{-\lambda}}(E)}(\Lambda)[-d_{b}]). \] 
        We note that the slopes of $b_{-\lambda}$ are dominant with respect to $\ol{B}$ while the slopes of $b_{-w_{0}(\lambda)}$ are anti-dominant with respect to $\ol{B}$.
To explicate the Eisenstein functors, we recall, for $b \in B(G)_{\mathrm{un}}$ unramified, the fixed Borel $B \subset G$ determines a Borel $B_{b} \subset G_{b}$, and it follows by \cite[Theorem~4.28]{HIDualCpx}, that we have a natural equivalence 
        \[ \Eis_{\ol{B}}i_{b^{T}_{-w_{0}(\lambda)}!}(-) \simeq j_{b_{-\lambda}!}(i_{\ol{B}_{b_{-\lambda}}}^{G_{b_{-\lambda}}}(-) \otimes \delta_{b_{-\lambda}}^{-\frac{1}{2}}[-d_{b}]) = j_{b_{-\lambda}!}^{\ren}i_{\ol{B}_{b_{-\lambda}}}^{G_{b_{-\lambda}}}(-). \]
        Moreover, it follows, by \cite[Corollary~2.2.5 (iv)]{HHS}, \cite[Theorem~4.28]{HIDualCpx}, and \cite[Lemma~4.19]{HIDualCpx}, that we have an isomorphism  
        \[ \mathrm{nEis}_{\ol{B}}i_{b^{T}_{-\lambda}!}(-) \simeq j_{b_{-\lambda}\sharp}(i_{B_{b_{-\lambda}}}^{G_{b_{-\lambda}}}(-) \otimes \delta_{b_{-\lambda}}^{-\frac{1}{2}}[-d_{b}]) = j_{b_{-\lambda}\sharp}^{\ren}i_{B_{b_{-\lambda}}}^{G_{b_{-\lambda}}}(-), \]
       
Now note that we may absorb the modulus character twist in the normalized parabolic induction into $\on{c-Ind}_{T(\mathcal{O}_{E})}^{T(E)}(\Lambda)$ and that we may absorb the modulus character twist by $\delta_{b_{-\lambda}}^{-1/2}$ into $\on{c-Ind}_{I_{b}}^{G_{b}(E)}(\Lambda)$.
In particular, we see that the claim is reduced to the identity 
        \[ \Ind_{B_{b}}^{G_{b}} \circ \on{c-Ind}_{T(\mathcal{O}_{E})}^{T(E)}(\Lambda) \simeq \on{c-Ind}_{I_{b}}^{G_{b}(E)}(\Lambda),  \]
        and similarly for the opposites, which is a straightforward calculation (see \cite[Lemma~1.6.1]{IwahoriHeckeAlgebras} for rational coefficients and \cite[Corollaire~3.6]{DatFinitudepourlesRepresentationsLissesdeGroupesPadiques} for an even more general claim with modular coefficients). 
        \end{proof}

        Now we explain how to combine Proposition \ref{prop: parameterwisepitch} with some $t$-exactness properties of the functor $\pitch$ to deduce some interesting consequences for the cohomology of local Shimura varieties.
We discuss this now. 
        \subsection{$t$-exactness of $\pitch$}
        The categories $\calD_{\lisse}(\Bun_{G},\Lambda)$ and $\Shv(\ICG,\Lambda)$ related under our functor $\pitch$, carry a wide variety of $t$-structures.
We recall the definitions of these $t$-structures in both contexts and then show that our functor is actually $t$-exact with respect to these $t$-structures.
        \subsubsection{$t$-structures in the schematic context}
        We recall that, since $\Shv^{!}(\ICG,\Lambda)$ is stratified by the gerbes $\bb{B}\underline{G_{b}(E)}$ for $b \in B(G)$, there is an obvious candidate for a perverse $t$-structure given by insisting that $*$ and $!$-pullbacks along the map $i_{b}: \Shv^{!}(\ICG_{b},\Lambda) \hookrightarrow \Shv^{!}(\ICG,\Lambda)$ sit in appropriate degrees determined by Verdier duality, where the role of Verdier duality is played by the canonical duality functor $\id_{\BZ}$. 
        
        However, we need to be a bit careful.
On the stack $\bb{B}\ul{G_{b}(E)}$ Verdier duality in the $\Shv^{!}$-will behave as if $\bb{B}\ul{G_{b}(E)}$ is of dimension $0$. This might lead one to believe that $\bb{B}\ul{G_{b}(E)}$ should be treated as a smooth space of dimension $0$ when defining a perverse $t$-structure. However, this is too naive.
Indeed, $\id_{\BZ}$ will intertwine $i_{b}^{*}$ and $i_{b}^{!}$ only up to a shift by $2d_{b}$ (as in \ref{eqn: SemiorthogonalBZdualityIsoc}), where we recall that we set $d_{b} := \langle 2\rho_{G},\nu_{b} \rangle$.
With this in mind, the obvious definition becomes the following.
        \begin{definition}{\label{defn: IsocPerverseTStructure}}
        We consider the pair of full subcategories $\Shv^{!}(\ICG,\Lambda)^{p,\leq 0}$ (resp. $\Shv^{!}(\ICG,\Lambda)^{p,\geq 0}$) of $\Shv^{!}(\ICG,\Lambda)$.
Here  $\Shv^{!}(\ICG,\Lambda)^{p,\leq 0}$ is the full subcategory generated under small colimits by $i_{b!}(\on{c-Ind}_{K}^{G_{b}(E)}(\Lambda)[n -d_{b}])$ for $K$ a compact open pro-$p$ subgroup and $n \in \bb{N}_{\geq 0}$, and $\Shv^{!}(\ICG,\Lambda)^{p,\leq 0}$ is the full subcategory of sheaves satisfying $i_{b}^{!}(A) \in \Rep^{\geq d_{b}}(\ICG_{b},\Lambda)$) for all $b \in B(G)$.
By \cite[Proposition~3.105]{Zhu25} specialized to the case of $\chi = 2\rho_{G}$, these give rise to a well-defined $t$-structure on $\Shv^{!}(\ICG,\Lambda)$ which we refer to as the perverse $t$-structure.
        \end{definition}
        On the other hand, the exceptional pullback $i_{b}^{\sharp}$ which is the right adjoint of $i_{b*}$, which ends up being intertwined with $i_{b}^{!}$ under what is known as admissible duality on the set of objects $\Shv^{!}(\ICG,\Lambda)_{\Adm}$ again up to a shift by $d_{b}$.
Using this, one can define the following exotic kind of $t$-structure. 
        \begin{definition}{\label{defn: IsocExoticTStructure}}
        We consider the pair of full subcategories $\Shv^{!}(\ICG,\Lambda)^{e,\leq 0}$ (resp. $\Shv^{!}(\ICG,\Lambda)^{e,\geq 0}$) of $\Shv^{!}(\ICG,\Lambda)$ satisfying that $i_{b}^{!}(A) \in \Shv^{!}(\ICG_{b},\Lambda)^{\leq d_{b}}$ (resp. $i_{b}^{\sharp}(A) \in \Shv^{!}(\ICG_{b},\Lambda)^{\geq d_{b}}$) for all $b \in B(G)$.
By \cite[Proposition~3.110]{Zhu25} specialized to the case of $\chi = 2\rho_{G}$, these give rise to a well-defined $t$-structure on $\Shv^{!}(\ICG,\Lambda)$, which we refer to as the exotic $t$-structure.
        \end{definition}
        We now turn our attention to the mirrors of these $t$-structures on $\calD_{\lisse}(\Bun_{G},\Lambda)$. 
        \subsubsection{$t$-structures in the analytic context}
        As before, since the moduli stack $\Bun_{G}$ is stratified by the gerbes $[\ast/\tilde{G}_{b}]$, there is a natural candidate for a perverse $t$-structure on this derived category, where we note that unlike the previous case of $\ICG$, it is completely clear that $[\ast/\tilde{G}_{b}]$ is $\ell$-cohomologically smooth of dimension $-d_{b}$ (due to the appearance of the unipotent part of $\tilde{G}_{b}$) and therefore such a definition should involve a shift by $d_{b}$.
However, unlike the case of $\ICG$, Verdier duality on the Artin $v$-stack $\Bun_{G}$ naturally incarnates as smooth duality under the equivalence $\Dlis(\Bun_{G}^{b},\Lambda) \simeq \Rep(G_{b}(E),\Lambda)$ not BZ-duality as was the case for $\ICG$.
As we will see, this is due to the fact that the natural perverse $t$-structure on $\Dlis(\Bun_{G},\Lambda)$ will match up with the exotic $t$-structure on $\Shv(\ICG,\Lambda)$ introduced in the previous section. 
        \begin{definition}{\label{defn: BunGPerverseTStructure}}
        We consider the full subcategory $\calD_{\lisse}^{p, \leq 0}(\Bun_{G},\Lambda)$ (resp. $\calD_{\lisse}^{p, \geq 0}(\Bun_{G},\Lambda)$) defined by the condition that $j_{b}^{*}(A) \in \calD_{\lisse}^{\leq d_b}(\Bun_G^b,\Lambda)$ (resp. $j_{b}^{!}(A) \in \calD_{\lisse}^{\geq d_b}(\Bun_G^b,\Lambda)$).
This gives rise to a well-defined $t$-structure by \cite[Proposition~8.1.5]{MingjiaPolI}, which we refer to as the perverse $t$-structure on $\Dlis(\Bun_{G},\Lambda)$.
        \end{definition}
        Now we are tasked with defining what matches up with the actual perverse $t$-structure on $\Shv^{!}(\ICG,\Lambda)$.
Analogous to the case of $\ICG$, this should be given by some exotic $t$-structure defined in terms of certain exceptional pullback functors.
Indeed, we have the following. 
        \begin{definition}{\label{defn: HadaltStructure}}
        We consider the full subcategories $\calD_{\lisse}^{e, \leq 0}(\Bun_{G},\Lambda)$ (resp. $\calD_{\lisse}^{e, \geq 0}(\Bun_{G},\Lambda)$).
Here $\calD_{\lisse}^{e, \leq 0}(\Bun_{G},\Lambda)$ is the full subcategory is generated by $j_{b\sharp}(\on{c-Ind}_{K}^{G(E)}(\Lambda)[n - d_{b}])$ for $K \subset G_{b}(E)$ a compact open pro-$p$ subgroup and $n \geq \bb{N}_{\geq 0}$, and $\calD_{\lisse}^{e, \geq 0}(\Bun_{G},\Lambda)$ is the full subcategory of objects $A \in \calD_{\lisse}(\Bun_{G},\Lambda)$ such that $j_{b}^{*}(A) \in \calD_{\lisse}^{\leq d_{b}}(\Bun_{G}^{b},\Lambda)$ for all $b \in B(G)$.
This gives rise to a well-defined $t$-structure by the exact same argument as \cite[Proposition~3.105]{Zhu25}, and we refer to it as the exotic $t$-structure on $\calD_{\Lambda}^{\an}(\Bun_{G})$.
        \end{definition}
        \begin{remark}{\label{rem: HadalTStructure}}
        The  $t$-structure introduced here is very much related to the hadal $t$-structure defined in \cite[Theorem~1.2.3]{HansenBeijingNotes}.
However, there actually a somewhat stronger claim is shown.
Indeed, in \textit{loc.cit} a $t$-structure is only constructed in the case that $\Lambda = \ol{\bb{Q}}_{\ell}$ on the subcategory of compact objects $\calD_{\lisse}(\Bun_{G},\Lambda)^{\omega} \subset \calD_{\lisse}(\Bun_{G},\Lambda)$.
This is a stronger statement since one needs to check that compact objects on individual HN-strata are preserved under standard truncation in order to see that this is well-defined, which requires the assumption that $\Lambda = \ol{\bb{Q}}_{\ell}$ together with a non-trivial result of Bernstein (see \cite[Remark~1.1.3]{HansenBeijingNotes}).
One can easily see however that the restriction of $\calD_{\lisse}^{e, \geq 0}(\Bun_{G},\Lambda)$ (resp. $\calD_{\lisse}^{e, \leq 0}(\Bun_{G},\Lambda)$) to the full subcategory of compact objects agrees with Hansen's construction when $\Lambda = \ol{\bb{Q}}_{\ell}$ (However, as just mentioned it is not clear that this restriction is still a $t$-structure). 
        \end{remark}
        \subsection{$t$-exactness of the equivalence}
        We now have the following $t$-exactness properties of the functor $\pitch$. 
        \begin{theorem}{\label{thm: texactnessstatement}}
        For $\Lambda/\bb{Z}_{\ell}$ a torsion ring, the following is true.
        \begin{enumerate}
        \item The functor $\pitch$ is $t$-exact for the perverse $t$-structure on $\Shv^{!}(\ICG,\Lambda)$ defined in \ref{defn: IsocPerverseTStructure}, and the exotic $t$-structure on $\calD_{\lisse}(\Bun_{G},\Lambda)$ defined in \ref{defn: HadaltStructure}.
        \item The functor $\pitch$ is $t$-exact for the exotic $t$-structure on $\Shv^{!}(\ICG,\Lambda)$ defined in \ref{defn: IsocExoticTStructure} and the perverse $t$-structure on $\calD_{\lisse}(\Bun_{G},\Lambda)$ as defined in \ref{defn: BunGPerverseTStructure}.
        \end{enumerate}
        \end{theorem}
        \begin{proof}
        This follows immediately from Theorem \ref{thm: MainTheoremPartial} (3)-(4).
        \end{proof}
        We now conclude our applications section by showcasing the power of Conjecture \ref{conj: IndPerfLinearity}, by combining Theorem \ref{thm: texactnessstatement} with the work of Yang-Zhu \cite{YangZhuTorsion}.
To explain this, we write $\mathrm{Tilt}(\hat{G}) \subset \Rep(\hat{G})$ for the subcategory of tilting modules (see \cite[\S~9.1]{HamGeomES} for a discussion of this).
For $V \in \mathrm{Tilt}(\hat{G})$, we consider the attached Hecke operator 
        \[ T_{V}: \Dlis(\Bun_{G},\Lambda) \ra \Dlis(\Bun_{G},\Lambda), \]
        as in (\ref{eqn: HeckeOperatorofV}).
If $\mu \in X_{*}(T_{\ol{E}})^{+}$ is a geometric dominant cocharacter of $G$ and $V = \mathcal{T}_{\mu} \in \mathrm{Tilt}(\hat{G})$ is the associated highest weight tilting module then we will denote this operator by $T_{\mu}$ We recall that this corresponds to the spectral action of the canonical vector bundle $C_{V}$ attached to $V$, and given a parameter $\phi$, the spectral action and in turn the Hecke operator preserves the full subcategory $\Dlis(\Bun_{G},\Lambda)_{\hat{\phi}} \hookrightarrow \Dlis(\Bun_{G},\Lambda)$ of $\phi$-complete sheaves.
By a completely analogous argument to \cite[Lemma~3.26,3.17]{YangZhuTorsion}, the perverse $t$-structure $(\Dlis^{p, \leq 0}(\Bun_{G},\Lambda),\Dlis^{p, \geq 0}(\Bun_{G},\Lambda))$ restricts to a well-defined perverse $t$-structure $(\Dlis^{p, \leq 0}(\Bun_{G},\Lambda)_{\hat{\phi}},\Dlis^{p, \geq 0}(\Bun_{G},\Lambda)_{\hat{\phi}})$ on $\Dlis(\Bun_{G},\Lambda)_{\hat{\phi}}$.
We now have the following.
        \begin{theorem}{\label{thm: PerverseTExactnessforHeckeOperatorsonBunG}}
        For $G/E$ unramified, we let $\hat{\phi}$ be an unramified semisimple $L$-parameter which is of weakly Langlands-Shahidi type (in the sense of \cite[Definition~6.2]{HamannLeeTorsion}) or equivalently generic (in the sense of \cite[Definition~1.1]{YangZhuTorsion}).
Assume that $\ell \nmid \lvert\pi_{0}(Z(G))\rvert$ satisfies assumption \ref{assump: modularcoefficients} and that \cref{conj: IndPerfLinearity} is true then the $\phi$-completed Hecke operator 
        \[ T_{V,\hat{\phi}}: \Dlis(\Bun_{G},\ol{\bb{F}}_{\ell})_{\hat{\phi}} \ra \Dlis(\Bun_{G},\ol{\bb{F}}_{\ell})_{\hat{\phi}} \]
        induced by $T_{V}$ on the category of $\phi$-completed sheaves is perverse $t$-exact.
        \end{theorem}
        \begin{proof}
        This is a formal consequence of \Cref{thm: texactnessstatement} (2) and  \Cref{prop: parameterwisepitch}.
        \end{proof}
        \begin{remark}
        We note that, in light of \Cref{rem: philocalvsphicomplete} this recovers (modulo some differences in the assumptions on $\ell$) \cite[Corollary~4.29]{HamannLeeTorsion} after restricting to the full subcategory of ULA objects.
        \end{remark}
        It is natural to ask what the point of this claim is given that the analogous claim also holds in the $\ICG$ setting by the work of \cite{YangZhuTorsion}.
One reason to care is that, unlike the spectral action induced by the embedding $\bb{L}_{G}^{\widehat{\unip}}$ which involves a categorical trace construction and is difficult to make explicit, the Hecke operator $T_{V}$ is readily related to the cohomology of local Shtuka spaces and Shimura varieties.
We explain this now.

    We recall that a local shtuka datum is a triple $(G,b,\mu)$ for $\mu$ a geometric dominant cocharacter of $G/E$ and $b \in B(G,\mu)$ an element of the $\mu$-admissible locus of the Kottwtiz set of $G$.
We let $E_{\mu}/E$ be the reflex field of $\mu$.
The triple $(G,b,\mu)$ defines a diamond
    \[ \Sht(G,b,\mu)_{\infty} \rightarrow \Spd(\Breve{E}_{\mu}) \]
    parameterizing modifications $\mathcal{E}_{b} \rightarrow \mathcal{E}_{G}^{0}$ with meromorphy bounded by $\mu$ on the Fargues-Fontaine curve $X$, where $\Breve{E}_{\mu} := \Breve{\bb{Q}}_{p}E_{\mu}$.
Here $\mathcal{E}_{b}$ is the $G$-bundle on $X$ attached to $\mathcal{E}$ and $\mathcal{E}_{G}^{0}$ is the trivial $G$-bundle.
The space carries an action of $G(\mathbb{Q}_{p}) \times G_{b}(\mathbb{Q}_{p})$ and a (non-effective) descent datum from $\Breve{E_{\mu}}$ down to $E_{\mu}$.
This allows us to consider the tower of quotients
    \[ [\Sht(G,b,\mu)_{\infty}/\underline{K}] =: \Sht(G,b,\mu)_{K} \]
    for varying open compact subgroups $K \subset G(\mathbb{Q}_{p})$.
We write $\mathcal{S}_{\mu}$ for the $\Lambda$-valued sheaf attached to the highest weight tilting module $\mathcal{T}_{\mu}$ of $W_{E} \ltimes \hat{G}$ as defined in \cite[\S~10.1]{HamGeomES}.
This is given by pulling back the sheaf on $\mathrm{Hck}_{G,E}$ defined by $\mathcal{T}_{\mu}$ and using geometric Satake for $\bb{B}_{\dR}^{+}$-grassmannian along the natural map $\Sht(G,b,\mu) \ra \mathrm{Hck}_{G,E}$.
In particular, the sheaf $\mathcal{S}_{\mu}$ is equivariant with respect to the actions of $G(\mathbb{Q}_{p})$ and $G_{b}(\mathbb{Q}_{p})$ by construction.
Letting $\Sht(G,b,\mu)_{K,\mathbb{C}_{p}}$ denote the base-change of these spaces to $\mathbb{C}_{p}$, we can now define the complex
\[ R\Gamma_{c}(G,b,\mu) := \colim_{K \rightarrow \{1\}} R\Gamma_{c}(\Sht(G,b,\mu)_{K,\mathbb{C}_{p}},\mathcal{S}_{\mu}) \]
of $G(\mathbb{Q}_{p}) \times J_{b}(\mathbb{Q}_{p}) \times W_{E}$-modules.
We now want to disentangle the $G(\mathbb{Q}_{p})$ and $G_{b}(\mathbb{Q}_{p})$ action.
To do this, for $\pi$ (resp. $\rho$) a smooth irreducible representation of $G(\mathbb{Q}_{p})$ (resp. $G_{b}(\mathbb{Q}_{p})$) on $\Lambda$-modules, we define the $\pi$ (resp. $\rho$)-isotypic part.
I.e we define the complexes 
\[ R\Gamma^{\flat}_{c}(G,b,\mu)[\rho] := \mathrm{Hom}_{G_{b}(\bb{Q}_{p})}(R\Gamma_{c}(G,b,\mu),\rho), \]
where the above denotes external (derived) Hom in the category of smooth $G_{b}(\bb{Q}_{p})$-representations.
It is complex of $G(\bb{Q}_{p})$-representations with a continuous action of $W_{E_{\mu}}$. 

In particular, we spell out one concrete consequence of \Cref{thm: PerverseTExactnessforHeckeOperatorsonBunG}, which seems to difficult to deduce directly on the geometry of $\ICG$ without consideration of the functor $\pitch$ (or at the very least nearby cycles) and Conjecture \ref{conj: IndPerfLinearity}.
\begin{corollary}
	\label{shtukas corollary}
Assume Conjecture \ref{conj: IndPerfLinearity} and that $\ell \nmid \lvert\pi_{0}(Z(G))\rvert$ and it satisfies \cref{assump: modularcoefficients}.
We let $\phi$ be an unramified $L$-parameter of weakly Langlands-Shahidi type, and let $\rho \in \Rep(G_{b}(E),\ol{\bb{F}}_{\ell})$ be a smooth irreducible representation such that its Fargues-Scholze parameter $\phi_{\rho}^{\mathrm{FS}}$ agrees with $\phi$ under the twisted embedding $\phantom{}^{L}G_{b} \ra \phantom{}^{L}G$ (as defined in \cite[\S~IX.7.1]{FS21}).
Then, for all geometric dominant cocharacters $\mu$, the complex 
\[ R\Gamma_{c}^{\flat}(G,b,\mu)[\rho] \in \Rep(G(E),\Lambda) \]
is concentrated in degrees $\geq \langle 2\rho_{G},\nu_{b} \rangle$. Moreover, if $\phi$ is of Langlands-Shahidi type then this complex is concentrated in degrees $\langle 2\rho_{G},\nu_{b} \rangle$.
\end{corollary}
\begin{proof}
We write $T_{\mu,\hat{\phi}}$ for $\hat{\phi}$-completed Hecke operator attached to $T_{\mu}$.
We note, by \cite[Lemma~4.7]{HIDualCpx}, \cite[Lemma~4.2 (1)]{HamannLeeTorsion} and Remark \ref{rem: philocalvsphicomplete}, that the assumption on $\rho$ guarantees that $j_{b*}(\rho[d_{b}]) \in \Dlis^{p,\geq 0}(\Bun_{G},\ol{\bb{F}}_{\ell})_{\hat{\phi}}$ lies in the $\hat{\phi}$-completed subcategory.
It follows that $T_{\mu}j_{b*}(\rho[d_{d}]) \simeq T_{\mu,\hat{\phi}}j_{b*}(\rho[d_{b}]) \in \Dlis^{p,\geq 0}(\Bun_{G},\ol{\bb{F}}_{\ell})_{\hat{\phi}}$ by \Cref{thm: PerverseTExactnessforHeckeOperatorsonBunG}, which implies that $j_{1}^{*}T_{\mu}j_{b*}(\rho[d_{b}]) \in \Rep^{\geq 0}(G(E),\Lambda)$.
By invoking \cite[Lemma~10.1]{HamGeomES} and moving over the shifts, this implies the desired claim. Similarly, if $\phi$ is of Langlands-Shahidi type then we claim that 
\[ j_{b*}(\rho) \simeq j_{b!}(\rho), \]
so that $j_{b*}(\rho[d_{b}])$ is actually a perverse sheaf. Indeed, it follows by Corollary \ref{cor: semiorthogonaldecompositionsplits}, Conjecture \ref{conj: IndPerfLinearity}, and Remark \ref{rem: primesuperscript} that we can rewrite $j_{b!}(\rho)$ as $j_{b\natural}(\rho)$. However, we easily see that $j_{b'}^{*}j_{b\natural}(\rho)$ is trivial for any $b' \geq b$ (see Remark \ref{rem: jbsharprestricted}), which implies that the natural map $j_{b!}(\rho) \simeq j_{b*}(\rho)$ is an isomorphism. By the exact same reasoning as above, we conclude that the complex is concentrated in degrees $\langle 2\rho_{G},\nu_{b} \rangle$.
\end{proof}
\begin{remark}
One can of course deduce many different variations of this Corollary by using the various comparisons of isotypic components and Hecke operators spelled out in \cite[Lemma~11.1]{HamGeomES}, as well as consider variations for the cohomology of the space of shtukas parametrizing modifications $\mathcal{E}_{b} \dashrightarrow \mathcal{E}_{b'}$ of meromorphy $\leq \mu$ for any $b,b' \in B(G)$.
We note that this and its variations prove special cases of an $\ell$-modular version of \cite[Conjecture~6.21]{KoshikawaShin}, where we note that the quantities $\langle \phi_{\pi_{b}}|_{\mathbb{G}_{m}^{D}},- \rangle$ appearing there will all be trivial by the weakly Langlands-Shahidi type assumption (see \cite[Remark~5.2]{KoshikawaShin}). It is also consistent with the results of \cite[Section~10]{HamGeomES}.
\end{remark}
\begin{remark}
We recall that the Hecke operator $T_{\mu}$ actually factors through the forgetful functor $\Dlis(\Bun_{G},\Lambda)^{BW_{E_{\mu}}} \ra \Dlis(\Bun_{G},\Lambda)$ from $E_{\mu}$-equivariant objects.
Similarly, the canonical vector bundle $C_{\mathcal{T}_{\mu}}$ attached to the highest weight tilting module $\mathcal{T}_{\mu}$ upgrades to a perfect complex with $E_{\mu}$-equivariant structure.
Under the spectral action of (\ref{eqn: FSSpectralAction}), these two structures are identified; in particular, using Conjecture \ref{conj: IndPerfLinearity} one would be able to compute the $W_{E_{\mu}}$ action on the complexes $R\Gamma_{c}^{\flat}(G,b,\mu)[\rho]$ in terms of the categorical conjecture.
For a flavor of how this explicitly should look like in the  Langlands-Shahidi type case, see for example \cite[Conjecture~11.18]{HamGeomES}. 
\end{remark}

	\appendix
	\section{A counterexample}
	
	\label{A counterexample}
	\begin{example}
		The purpose of this example is to exhibit a $\diamondsuit$-cover which is not a v-cover.
		For this we study a specific type of valuation ring.
		Let $\calP$ be a total order and let $\Gamma=\oplus_{p \in \calP}\gamma^\bbZ_p$. 
		We will write elements $\gamma\in \Gamma$ using multiplicative notation.
		We define a total order on $\Gamma$ as follows.
		Given two elements of the form $\prod_{i=1}^n \gamma^{a_i}_{p_i}$ and $\prod_{i=1}^n \gamma^{b_i}_{p_i}$ with $p_1<p_2<\dots<p_n$, then $\prod_{i=1}^n \gamma^{a_i}_{p_i}\leq \prod_{i=1}^n \gamma^{b_i}_{p_i}$ if and only if $(-a_1,\dots,-a_n)\leq (-b_1\dots,-b_n)$ in $\bbZ^n$ with its usual lexicographic order. 
		For example, $\gamma^2_p<\gamma_p$ since $-2<-1$ and if $p_1<p_2<p_3$ in $\calP$ then $\gamma_{p_1}\cdot \gamma^2_{p_3}<\gamma_{p_2}\cdot \gamma_{p_3}^2$, since $(-1,0,-2)<(0,-1,-2)$.
		Let 
		\[K=\{f:\Gamma\to k\mid \on{supp}(f) \text{ is well ordered}\}\] 
		be the field of Hahn series with coefficients in $k$ and values in $\Gamma$. 
		We write expressions of the form
		\[f=\sum_{\gamma\in \Gamma} f_{\gamma} \gamma\]
		for elements of $K$.
		We have a valuation map $|\cdot|:K\to \{0\}\cup \Gamma$ which sends $f$ to the smallest element of the support of $f$, which is well defined since the support is well-ordered.
		We let $V\subseteq K$ denote the valuation ring defined by $|\cdot|$. 
		
		For every element $p\in \calP$ the element $\gamma_p\in V$ generates a prime ideal $\frakp^+_p\subseteq V$. 
		Concretely, the prime ideal $\frakp_p^{+}$ consists of those elements in $g\in K$ such that $|g|\leq |\gamma_p|$. 
		We also consider $\frakp^-_p$ to denote the largest ideal in $V$ not containing $\gamma_p$.
		Observe that the functions
		\[\frakp_{(-)}:\calP\to \Spec V \text{ with } p\mapsto \frakp^+_p\text{ and } p\mapsto \frakp^-_p\]
		are order preserving in the sense that if $p_1<p_2$ then $\frakp^+_{p_1}\subseteq \frakp^+_{p_2}$.
		Indeed, if $p_1<p_2$ implies that $\gamma_{p_1}<\gamma_{p_2}$ so $\langle \gamma_{p_1}\rangle \subseteq  \langle \gamma_{p_2}\rangle$.
		
		Suppose that $M\in \calP$ is a maximum. 
		Then $\frakp_{M}\subseteq V$ is the maximal ideal. 
		Indeed, if $f\in V$ and it is not a unit then $|f|\in \Gamma$ has the form $|f|=\gamma_{p_1}^{a_1}\cdot(\prod_{i=2}^{n}\gamma_{p_i}^{a_i})$ with $a_1\geq 1$ and $p_1<p_i$.
		Since ${M}\geq p_1$, then $\gamma_{p_1}^{a_1}\cdot(\prod_{i=2}^{n}\gamma_{p_i}^{a_i})\leq \gamma_{M}$ and $f\in \frakp_{M}$. 
		Analogously, if $m\in \calP$ is a minimum then $\frakp_{p_m}^-=0$ and $\frakp_{p_m}^+$ is the unique smallest non-empty prime ideal of $V$.
		
		From now on, we assume that $\calP$ has both a maximum and a minimum which we denote by $M$ and $m$, respectively. 
		If we assume that $\calP\setminus \{M,m\}\neq \emptyset$ then the map 
		\[f:X=\Spec V/\frakp^+_m\coprod \Spec V_{(\frakp^-_M)}\to \Spec V=Y\]
		is an arc-cover that is not a v-cover.
		We claim that if we impose certain hypothesis on $\calP$ then this map is also a $\diamondsuit$-cover.
		
		Let us spell out what it means to be a $\diamondsuit$-cover in this context.
		Since products of points are basis for the v-topology and $f^\diamondsuit$ is partially proper by construction of $(-)^{\diamondsuit}$ it suffices to show that every map of the form $\Spa(R,R^\circ)\to Y^\diamondsuit$ lifts to $X^\diamondsuit$ where $R^\circ=\prod_{i\in I}O_{C_i}$ with $O_{C_i}\subseteq C_i$ the ring of integers of a non-archimedean field $C_i$.
		A map as above is specified by a ring map $g:V\to R$, which further induces ring maps $g_i:V\to C_i$.  
		The collection $\{g_i\}_{i\in I}$ determines $g$, but not every collection arises from a map $g:V\to R$.
		We can partition the set $I=I_0\coprod I_{\neq 0}$ according to whether $\on{ker}(g_i)\subseteq V$ is $0$ or whether $\frakp^+_m\subseteq \on{ker}(g_i)$, respectively.
		This partition gives a topological splitting $R=R_{I_0}\times R_{I_{\neq 0}}$ and maps $g:V\to R$ are in bijection with pairs of maps $(g_I,g_{I_{\neq 0}}):V\to R_{I_0}\times R_{I_{\neq 0}}$.
		By construction $g_{I_{\neq 0}}$ factors through the quotient $V/\frakp^+_m$, and in particular lifts to $X^\diamondsuit$. 
		It suffices then to consider maps $g:V\to R$ such that $\on{ker}(g_i)=0$ for all $i\in I$.
		
		It is not in general true that such maps $g:V\to R$ always lift to $X^\diamondsuit$. 
		Our claim is that for specific choices of $\calP$ the lift always exist (necessarily through $\Spec V_{(\frakp^-_M)}^\diamondsuit\to Y^\diamondsuit$).  
		To summarize the situation, our context is as follows. 
		We are given a ring map $g:V\to R$ with $(R,R^\circ)$ a product of rank $1$ points of index $I$, with the additional property that $g_i:V\to C_i$ factors through $\on{Frac}(V)\to C_i$, and we must show that for specific choices of $\calP$ the image of $\gamma_M$ in $R$ is automatically a unit.
		
		If the image of $\gamma_M$ is not a unit then there is a connected component $\alpha$ of $\Spa(R,R^\circ)$, for which the map $V\to C_\alpha$ factors through the residue field $V\to k$.
We fix such an $\alpha$.  
		In particular, for every element $f\in \frakp_M^+$ and every $n\in \bbN$ the set $\{|f|\leq |\varpi|^n\}\subseteq \Spa(R,R^+)$ is a quasicompact non-empty open subset for $\varpi \in R$ a uniformizing element.
		We will show, assuming this happens and under an additional assumption on $\calP$ (in particular, that the cardinality is sufficiently large), that there exists an element $v\in V$ such that the family of elements $g_i(v)\in C_i$ is not uniformly bounded,  contradicting that the family $g_i$ arose from a map $g:V\to R$. 
		
		For an element $t\in\calP$, we consider 
		\[I^t_{n\leq}=\{i\in I |  |g_i(\gamma_t)|_{C_i}\leq |\varpi_i^n|_{C_i}\}.\]
		This gives rise to an open and closed subset of $\Spa(R,R^+)$ which we denote by $U^t_{n\leq}$.
		
		Observe that \[\alpha\in \bigcap_{t\in \calP} \bigcap_{n\in \bbN} U^t_{n\leq},\]
		and that each finite intersection of the form $J=\bigcap^n_{j=1} I_{n_i\leq}^{p_i}$ gives rise to an infinite set $J\subseteq I$. 
		
		Suppose that $s<t$ in $\calP$.
We let $\on{Err}(s,t)\subseteq \bbN$ denote the subset of natural numbers $n\in \bbN$ for which $I^t_{n\leq}$ contains an element $i$ for which $|g_i(\gamma^2_t)|_{C_i}\leq |g_i(\gamma_s)|_{C_i}$.
		We claim that this set is finite.
		Indeed, if this set was infinite then there is an infinite sequence of elements 
		\[\{i_0,\dots, i_j\dots \}_{j\in \bbN} \subseteq I, \]
		and an increasing sequence of numbers $n_j$ such that 
		\[|g_{i_j}(\gamma^2_t)|_{C_{i_j}}\leq |g_{i_j}(\gamma_s)|_{C_{i_j}} \text{ and } |g_{i_j}(\gamma_t)|_{C_{i_j}} \leq |\varpi_{i_j}^{n_j}|.\] 
		In particular, the element $\gamma_s\cdot \gamma_t^{-3}\in V$ and $|g_{i_j}(\gamma_s\cdot \gamma_t^{-3})|_{C_{i_j}}\geq |\varpi^{-n_j}_{i_j}|_{C_{i_j}}$.
This shows that if $\on{Err}(s,t)$ is infinite then the family $\{g_i(\gamma_s\cdot \gamma_t^{-2})\}_{i\in I}$ is not uniformly bounded and it cannot arise from a map $V\to R$ in which case we would be done.
        
		We let $N_{(s,t)}$ be the smallest number such that, for all $i\in I^t_{n\leq}$, the inequality $|g_i(\gamma_s)|_{C_i}<|g_i(\gamma^2_t)|_{C_i}$ holds.
		
		We define a partial order $\calP_n\subseteq \calP$ for all $n\in \bbN$.
		The elements of $\calP_n$ are those elements of $t\in \calP$ such that $N_{(m,t)}, N_{(t,M)}\leq n$. 
		Given $s,t\in \calP_n$ we say that $s<_n t$ in $\calP_n$ if $s<t$ in $\calP$ and $N_{(s,t)}\leq n$.
		Equivalently, $s<_n t$ if and only if $s<t$ in $\calP$ and, for all $i\in I^M_{n\leq }$, 
		we have 
		\[|g_i(\gamma_m)|_{C_i}<|g_i(\gamma^2_s)|_{C_i} <|g_i(\gamma^4_t)|_{C_i}<|g_i(\gamma^8_M)|_{C_i}\leq \varpi^{8n}.\]
		We have that 
		\[\calP=\cup_{n\in \bbN} \calP_n.\]
		
		Moreover, every totally ordered subset of $\calP_n$ is finite. 
		Even more is true.
For all $n\in \bbN$, we can find a function $f_n:\calP\to \bbN\cup \{-1\}$ such that $f_n(p)=-1$ if $p\notin \calP_n$, $f_n(p)=0$ if $p=M$ and $f_n(p)=k$ if the largest chain starting in $p$ and ending in $M$ has length $k$. 
		We claim that these two properties already impose restrictions on $\calP$.
		Indeed, consider the function 
		\[\Theta:\calP\to (\bbN\cup \{-1\})^\bbN \text{ with } p\mapsto f_{(-)}(p).\]   
		We claim that $\Theta$ is injective. 
		Indeed, if $n\geq N_{(s,t)}, N_{(m,t)}, N_{(m,s)}, N_{(t,M)}, N_{(s,M)}$ then $f_n(s)\neq f_n(t)$ since $s<_n t$. 
		It follows that if we pick $\calP$ with a cardinality larger than $2^{\aleph_0}$ then $\mathrm{Err}(s,t)$ must be infinite and by the above reasoning the family $\{g_{i}(\gamma_{s}\gamma_{t}^{-2})\}_{i \in I}$ would not be uniformly bounded.
As explained above, this then implies that $V$ admits a $\diamondsuit$-cover that is not a v-cover.
	\end{example}

\printbibliography

\end{document}